\documentclass[12pt]{article}
\usepackage{amsmath,latexsym,amssymb}
\usepackage[OT2, OT1]{fontenc}
\usepackage[russian, english]{babel}
\begin{document}
\headsep 0.5 true cm

\begin{center}
{\Large\bf Geometric realizations of representations for
           $\text{PSL}(2, \mathbb{F}_p)$ and Galois representations
           arising from defining ideals of modular curves}
\vskip 1.0 cm
{\large\bf Lei Yang}
\end{center}
\vskip 1.5 cm

\begin{center}
{\large\bf Abstract}
\end{center}

\vskip 0.5 cm

  We construct a geometric realization of representations for
$\text{PSL}(2, \mathbb{F}_p)$ by the defining ideals of rational
models $\mathcal{L}(X(p))$ of modular curves $X(p)$ over $\mathbb{Q}$,
which gives rise to a Rosetta stone for geometric representations
of $\text{PSL}(2, \mathbb{F}_p)$. The defining ideal of a modular
curve, i.e., an anabelian counterpart of the Eisenstein ideal, is
the anabelianization of the Jacobian of this modular curve
and is a reification of the fundamental group $\pi_1$. We show that
there exists a correspondence among the defining ideals of modular
curves over $\mathbb{Q}$, reducible $\mathbb{Q}(\zeta_p)$-rational
representations
$\pi_p: \text{PSL}(2, \mathbb{F}_p) \rightarrow \text{Aut}(\mathcal{L}(X(p)))$
of $\text{PSL}(2, \mathbb{F}_p)$, and $\mathbb{Q}(\zeta_p)$-rational Galois
representations $\rho_p: \text{Gal}(\overline{\mathbb{Q}}/\mathbb{Q})
\rightarrow \text{Aut}(\mathcal{L}(X(p)))$ as well as their modular and
surjective realization. It is an anabelian counterpart of the global
Langlands correspondence for $\text{GL}(2, \mathbb{Q})$ by the \'{e}tale
cohomology of modular curves as well as an anabelian counterpart of Artin's
conjecture, Serre's modularity conjecture, the Fontaine-Mazur conjecture
and the $p$-adic local Langlands correspondence for $\text{GL}(2, \mathbb{Q}_p)$.
It is an ideal theoretic (i.e. nonlinear) counterpart of Grothendieck's
section conjecture and an ideal theoretic (i.e. nonlinear) reification
of ``arithmetic theory of $\pi_1$'' expected by Weil for modular curves.
It is also an anabelian counterpart of the theory of Kubert-Lang and Mazur-Wiles
on the cuspidal divisor class groups and the Eisenstein ideals of modular
curves.

\vskip 10.0 cm

\begin{center}
{\large\bf Contents}
\end{center}

$$\aligned
  &\text{1. Introduction}\\
  &\text{2. A Rosetta stone for geometric representations
            of $\text{PSL}(2, \mathbb{F}_p)$}\\
  &\quad \text{2.1. Representations of $\text{PSL}(2, \mathbb{F}_p)$
            and modular curves}\\
  &\quad \text{2.2. Analytic theory and de Rham realization of
            representations for}\\
  &\quad \quad \quad \text{$\text{PSL}(2, \mathbb{F}_p)$: Riemann
            surfaces}\\
  &\quad \text{2.3. Algebraic theory and $\ell$-adic realization
            of representations for}\\
  &\quad \quad \quad \text{$\text{PSL}(2, \mathbb{F}_p)$:
            Drinfeld-Deligne-Lusztig curves}\\
  &\quad \text{2.4. Arithmetic theory and $\mathbb{Q}(\zeta_p)$-rational
            realization of}\\
  &\quad \quad \quad \text{representations for $\text{PSL}(2, \mathbb{F}_p)$:
            $\mathcal{L}(X(p))$}\\
  &\quad \text{2.5. A comparison between Hecke's decomposition formulas and our}\\
  &\quad \quad \quad \text{decomposition formulas, i.e., motives versus anabelian geometry}\\
  &\text{3. The division of elliptic functions, the transformation of
            elliptic functions}\\
  &\quad \text{and their Galois groups}\\
  &\quad \text{3.1. The division equations of elliptic functions, their
            Galois groups}\\
  &\quad \quad \quad \text{and their relation with transformation equations}\\
  &\quad \text{3.2. The theory of the special transformation equations}\\
  &\quad \text{3.3. The transformation theory of elliptic functions and modular}\\
  &\quad \quad \quad \text{equations}\\
  &\quad \text{3.4. The Galois groups of the special transformation equations}\\
  &\quad \text{3.5. Galois' theorem}\\
  &\quad \text{3.6. The resolvent of the seventh degree}\\
  &\quad \text{3.7. The two resolvents of the eleventh degree}\\
  &\text{4. Galois representations arising from $\mathcal{L}(X(p))$ and their modularity}\\
  &\quad \text{4.1. The case of $p=7$ and Klein quartic curve}\\
  &\quad \text{4.2. The case of $p=11$ and Klein cubic threefold}\\
  &\quad \text{4.3. The case of $p=13$ and our quartic fourfold}\\
  &\text{5. An anabelian counterpart of the global Langlands correspondence}\\
  &\quad \text{for $\text{GL}(2, \mathbb{Q})$ by the cohomology of modular curves}\\
  &\text{6. An anabelian counterpart of the Shimura-Taniyama-Weil conjecture}\\
\endaligned$$
$$\aligned
  &\text{7. An anabelian counterpart of Artin's conjecture, Serre's modularity}\\
  &\quad \text{conjecture and the Fontaine-Mazur conjecture}\\
  &\quad \text{7.1. Relation with the Artin conjecture and Langlands' principle of}\\
  &\quad \quad \quad \text{functoriality}\\
  &\quad \text{7.2. Relation with Deligne's theorem ($k \geq 2$) and the Deligne-Serre}\\
  &\quad \quad \quad \text{theorem ($k=1$)}\\
  &\quad \text{7.3. Relation with Serre's modularity conjecture and Serre weights}\\
  &\quad \text{7.4. Relation with the Fontaine-Mazur-Langlands conjecture}\\
  &\text{8. An anabelian counterpart of the $p$-adic local Langlands correspondence}\\
  &\quad \text{for $\text{GL}(2, \mathbb{Q}_p)$}\\
  &\text{9. Comparison with the finite Langlands correspondence (the Macdonald}\\
  &\quad \text{correspondence)}\\
  &\text{10. An ideal theoretic counterpart of Grothendieck's section conjecture}\\
  &\quad \text{and an ideal theoretic reification of ``arithmetic theory of $\pi_1$''}\\
  &\text{11. An anabelian counterpart of the theory of Kubert-Lang and Mazur-}\\
  &\quad \text{Wiles on the cuspidal divisor class groups and the Eisenstein ideals}\\
  &\quad \text{11.1. Cyclotomic fields and Stickelberger ideals: the Herbrand-Ribet}\\
  &\quad \quad \quad \text{theorem}\\
  &\quad \text{11.2. Cuspidal divisor class groups and Eisenstein ideals of modular}\\
  &\quad \quad \quad \text{curves: the theory of Kubert-Lang and Mazur-Wiles}\\
  &\quad \text{11.3. Modular curves and cyclotomic fields: Sharifi's conjecture}\\
  &\quad \text{11.4. Division fields of elliptic curves: a non-abelian counterpart of}\\
  &\quad \quad \quad \text{the Herbrand-Ribet theorem}\\
  &\quad \text{11.5. Defining ideals of modular curves: an anabelian counterpart of the}\\
  &\quad \quad \quad \text{eigenspace decomposition and the cuspidal divisor class formula}\\
  &\text{12. Other consequences and relations with various topics}\\
  &\quad \text{12.1. Relation with Langlands' work}\\
  &\quad \text{12.2. Relation with motives and automorphic representations}\\
  &\quad \text{12.3. Relation with Kottwitz's conjecture and Arthur's conjecture}\\
  &\quad \text{12.4. Relation with the construction of interesting algebraic cycles}\\
\endaligned$$
$$\aligned
  &\quad \text{12.5. Relation with the realization of Shimura varieties}\\
  &\quad \text{12.6. Relation with the action of $\text{Gal}(\overline{\mathbb{Q}}/\mathbb{Q})$
                on Grothendieck's}\\
  &\quad \quad \quad \text{dessins d'enfants}\\
\endaligned$$

\begin{center}
{\large\bf 1. Introduction}
\end{center}

  We construct a geometric realization of representations for
$\text{PSL}(2, \mathbb{F}_p)$ by the defining ideals of rational
models $\mathcal{L}(X(p))$ of modular curves $X(p)$ over $\mathbb{Q}$.
In fact, the modular curve $X(p)$ is defined by a collection of
quartics which we can write down explicitly. More precisely, Klein
showed that $X(p)$ is isomorphic to the locus
$\mathcal{L}(X(p)) \subset \mathbb{CP}^{\frac{p-3}{2}}$ given by his
$z$-curves, i.e., the Klein $z$-curves (see Theorem 2.4.1, Theorem 2.4.4
and Theorem 2.4.3 for more details). In particular, Klein's $z$-curves
are nonsingular (see \cite{Ve}, Theorem 10.6). This leads us to study
the locus $\mathcal{L}(X(p))$, especially its defining ideal
$I(\mathcal{L}(X(p)))$ in the projective space $\mathbb{P}^{\frac{p-3}{2}}$.
Hence, for the irreducible representations of $\text{PSL}(2, \mathbb{F}_p)$,
whose geometric realizations can be formulated in three different scenarios
in the framework of Weil's Rosetta stone: number fields, curves over
$\mathbb{F}_q$ and Riemann surfaces.

$$\begin{matrix}
 &\text{Riemann surfaces} &\text{Curves over $\mathbb{F}_q$} &\text{Number fields}\\
 &\text{analytic theory} &\text{algebraic theory} &\text{arithmetic theory}\\
 &\text{de Rham cohomology} &\text{$\ell$-adic cohomology} &\text{defining ideals}\\
 &\text{(motives)} &\text{(motives)} &\text{(ideals)}\\
 &\text{linear structure} & \text{linear structure} & \text{nonlinear structure}\\
 &\text{de Rham realization} &\text{$\ell$-adic realization} &\text{rational realization}\\
 &\text{complex} &\text{$\ell$-adic} &\text{$\mathbb{Q}(\zeta_p)$-rational}\\
 &\text{representations} &\text{representations} &\text{representations}\\
 &(\ell=\infty) & (\ell \neq p) & (\ell=p)\\
 &\text{modular curves/$\mathbb{C}$} &\text{modular curves/$\mathbb{Z}_p$}
 &\text{modular curves/$\mathbb{Q}$}\\
 &\text{Hecke's level theory} &\text{Drinfeld's level theory} &\text{our level theory}\\
 &\dim H^1(X(p), \mathcal{O}_{X(p)}) & \dim H_c^1(\text{DL}_q, \overline{\mathbb{Q}}_{\ell})
 &\dim I(\mathcal{L}(X(p)))\\
 &\text{(genus of $X(p)$)} &\text{($\text{DL}_q$: $x y^q-x^q y=1$)} &\text{(dimension of ideals)}\\
 &\text{Jacobian $J(p)$: $H^1$} & \text{Jacobian: $H^1$} & \text{locus: $\mathcal{L}(X(p))$}\\
 &\text{(reification of $H_1$)} &\text{(reification of $H_1$)} & \text{(reification of $\pi_1$)}\\
 &\text{abelian}  &\text{abelian} &\text{anabelian}\\
 &\text{algebraic geometry} &\text{algebraic geometry} &\text{algebraic geometry}
\end{matrix}\eqno{(1.1)}$$

  In particular, we show that there exists a correspondence among
the defining ideals of modular curves over $\mathbb{Q}$, reducible
$\mathbb{Q}(\zeta_p)$-rational representations
$\pi_p: \text{PSL}(2, \mathbb{F}_p) \rightarrow \text{Aut}(\mathcal{L}(X(p)))$
of $\text{PSL}(2, \mathbb{F}_p)$, and $\mathbb{Q}(\zeta_p)$-rational Galois
representations $\rho_p: \text{Gal}(\overline{\mathbb{Q}}/\mathbb{Q})
\rightarrow \text{Aut}(\mathcal{L}(X(p)))$ as well as their modular and
surjective realization.

\textbf{Theorem 1.1. (Main Theorem 1).} {\it There is the following
decomposition as}
$\text{Gal}(\overline{\mathbb{Q}}/\mathbb{Q}) \times \text{PSL}(2,
\mathbb{F}_p)$-{\it representations, that is, there is an isomorphism}
$$V_{I(\mathcal{L}(X(p)))}=\bigoplus_{j} m_{p, j} V_{p, j},\eqno{(1.2)}$$
{\it where $j$ denotes the $j$-dimensional irreducible representations
of} $\text{PSL}(2, \mathbb{F}_p)$ {\it appearing in} (1.2), {\it and
$m_{p, j}$ denotes its multiplicities. Correspondingly, the defining ideal
$I(\mathcal{L}(X(p)))$ has the following decomposition as the intersection
of invariant ideals under the action of}
$\text{Gal}(\overline{\mathbb{Q}}/\mathbb{Q}) \times \text{PSL}(2, \mathbb{F}_p)$:
$$I(\mathcal{L}(X(p)))=\bigcap_{j} I_{p, j},\eqno{(1.3)}$$
{\it where the ideals $I_{p, j}$ correspond to the representations $V_{p, j}$.
This decomposition} (1.2) {\it gives rise to the following two kinds of
representations:}
$$\pi_{p, j}: \text{PSL}(2, \mathbb{F}_p) \longrightarrow \text{Aut}(V_{p, j})
  \eqno{(1.4)}$$
{\it is an irreducible $\mathbb{Q}(\zeta_p)$-rational representation
of} $\text{PSL}(2, \mathbb{F}_p)$, {\it and}
$$\rho_{p, j}: \text{Gal}(\overline{\mathbb{Q}}/\mathbb{Q}) \longrightarrow
  \text{Aut}(V_{p, j})\eqno{(1.5)}$$
{\it is an irreducible $\mathbb{Q}(\zeta_p)$-rational representation
of} $\text{Gal}(\overline{\mathbb{Q}}/\mathbb{Q})$. {\it On the other
hand},
$$\pi_p: \text{PSL}(2, \mathbb{F}_p) \longrightarrow \text{Aut}(\mathcal{L}(X(p)))\eqno{(1.6)}$$
{\it is a reducible $\mathbb{Q}(\zeta_p)$-rational representation of}
$\text{PSL}(2, \mathbb{F}_p)$, {\it and}
$$\rho_p: \text{Gal}(\overline{\mathbb{Q}}/\mathbb{Q}) \longrightarrow
  \text{Aut}(\mathcal{L}(X(p)))\eqno{(1.7)}$$
{\it is a reducible $\mathbb{Q}(\zeta_p)$-rational representation of}
$\text{Gal}(\overline{\mathbb{Q}}/\mathbb{Q})$. {\it Locally},
$\pi_{p, j}$ {\it mat-ches} $\rho_{p, j}$ {\it by the following
one-to-one correspondence}
$$\pi_{p, j} \longleftrightarrow \rho_{p, j}.\eqno{(1.8)}$$
{\it Globally}, $\pi_p$ {\it matches} $\rho_p$ {\it by the following
one-to-one correspondence}
$$\pi_p \longleftrightarrow \rho_p.\eqno{(1.9)}$$
{\it Local-to-globally, the set $\{ \pi_{p, j} \}$ correspond to the
set $\{ \rho_{p, j} \}$, where $\{ \pi_{p, j} \}$ are joint irreducible
$\mathbb{Q}(\zeta_p)$-rational representations of}
$\text{PSL}(2, \mathbb{F}_p)$, {\it and} $\{ \rho_{p, j} \}$ {\it are
joint irreducible} $\mathbb{Q}(\zeta_p)$-{\it rational representations
of} $\text{Gal}(\overline{\mathbb{Q}}/\mathbb{Q})$. {\it In particular,
the Galois representations $\rho_{p, j}$ are realized by the algebraic
cycles corresponding to the ideals $I_{p, j}$.}

\textbf{Corollary 1.2.} {\it The explicit decomposition formulas and
the dimension formulas of the defining ideals $I(\mathcal{L}(X(p)))$
for the special cases $p=7$, $11$ and $13$$:$

$(1)$ $p=7$, $V_{I(\mathcal{L}(X(7)))}=\mathbf{1}$ is the trivial representation,
$\dim V_{I(\mathcal{L}(X(7)))}=1$. Correspondingly,
$$I(\mathcal{L}(X(7)))=I_1.\eqno{(1.10)}$$

$(2)$ $p=11$, $V_{I(\mathcal{L}(X(11)))}=\mathbf{10}$ is the discrete
series representation corresponding to $\chi_5$, $\dim V_{I(\mathcal{L}(X(11)))}=10$.
Correspondingly,
$$I(\mathcal{L}(X(11)))=I_{10}.\eqno{(1.11)}$$

$(3)$ $p=13$, $V_{I(\mathcal{L}(X(13)))}=\mathbf{1} \oplus \mathbf{7} \oplus \mathbf{13}$
is the direct sum of the trivial representation, the degenerate principal series
representation and the Steinberg representation, $\dim V_{I(\mathcal{L}(X(13)))}=21$.
Correspondingly,
$$I(\mathcal{L}(X(13)))=I_{1} \cap I_7 \cap I_{13}.\eqno{(1.12)}$$}

\textbf{Theorem 1.3. (Main Theorem 2).} (Surjective realization and modularity).
{\it For a given representation arising from the defining ideals}
$$\rho_p: \text{Gal}(\overline{\mathbb{Q}}/\mathbb{Q}) \rightarrow
  \text{Aut}(\mathcal{L}(X(p))), \quad (p \geq 7),$$
{\it there are surjective realizations by the equations of degree $p+1$ or
$p$ with Galois group isomorphic with} $\text{Aut}(\mathcal{L}(X(p)))$. {\it
These equations are defined over $\mathcal{L}(X(p))$, their coefficients are
invariant under the action of} $\text{Aut}(\mathcal{L}(X(p)))$. {\it Moreover,
it corresponds to the $p$-th order transformation equation of the $j$-function
with Galois group isomorphic to} $\text{PSL}(2, \mathbb{F}_p)$. {\it In particular,}

(1) {\it when $p=7$, there are two such equations of degree eight and seven, respectively.
One comes from the Jacobian multiplier equation of degree eight, the other corresponds
to the decomposition of} $\text{PSL}(2, \mathbb{F}_7)=S_4 \cdot C_7$, {\it where $C_7$
is a cyclic subgroup of order seven.}

(2) {\it When $p=11$, there are two such equations of degree twelve and eleven, respectively.
One comes from the Jacobian multiplier equation of degree twelve, the other corresponds
to the decomposition of} $\text{PSL}(2, \mathbb{F}_{11})=A_5 \cdot C_{11}$, {\it where
$C_{11}$ is cyclic subgroup of order eleven.}

(3) {\it When $p=13$, there are two distinct equations of the same degree fourteen,
respectively. One comes from the Jacobian multiplier equation of degree fourteen,
the other is not the Jacobian multiplier equation.}

{\it Remark.} The significance of the case of $p=7$ is that the
Jacobian multiplier equation of degree eight can be reduced {\it
explicitly} to the the modular equation of degree eight which
corresponds to the seventh-order transformation of $j$-functions.
Hence, they have the same Galois group. Similarly, the equation
of degree seven can also be reduced {\it explicitly} to the the
modular equation of degree seven which also corresponds to the
seventh-order transformation of $j$-functions.

  Our theorems lie at the crossroads of the main areas of number theory
and arithmetic geometry: Weil's Rosetta stone of geometric representation
theory which unifies both motives and anabelian algebraic geometry,
the last mathematical testament of Galois and Klein's elliptic modular
functions, Grothendieck's motives and the global Langlands correspondence
for $\text{GL}(2, \mathbb{Q})$, the Shimura-Taniyama-Weil conjecture,
Artin's cojecture, Serre's modularity conjecture, the Fontaine-Mazur
conjecture, $p$-adic local Langlands correspondence for 
$\text{GL}(2, \mathbb{Q}_p)$, the finite Langlands correspondence, i.e. 
the Macdonald correspondence, Grothendieck's anabelian algebraic geometry, 
``arithmetic theory of $\pi_1$'' expected by Weil, and non-commutative 
Iwasawa theory. More precisely, they can be put into the following nine 
distinct kinds of contexts:

(1) As an anabelian algebraic geometric realization of representations
for $\text{PSL}(2, \mathbb{F}_p)$, this gives rise to a Rosetta stone
for geometric representations of $\text{PSL}(2, \mathbb{F}_p)$
(geometric representation theory in a unified context which includes
both motives and anabelian algebraic geometry). It is well-known
that the genus $g(X(p))$ for the modular curve $X(p)$ is given by
the dimension of its Jacobian variety $J(p)$.  As a consequence of
the following comparison (1.1.13) between Hecke's decomposition
formulas and our decomposition formulas, this gives rise to an
anabelian counterpart of the genus $g(X(p))$ by the dimension of
defining ideals $\dim I(\mathcal{L}(X(p)))$. Hence, for the same
modular curve $X(p)$, there are two distinct kinds of geometric
invariants, one is $g(X(p))$ coming from $H_1$, the other is
$\dim I(\mathcal{L}(X(p)))$ coming from $\pi_1$.
$$\begin{matrix}
  &\textbf{A comparison between}\\
  &\textbf{Hecke's decomposition and our decomposition}
\end{matrix}$$
$$\begin{matrix}
            & \text{motives} & \text{anabelian algebraic geometry}\\
            & \text{Jacobian varieties} & \text{defining ideals}\\
            & \text{(Abel, Riemann)} & \text{(Weierstrass, Klein)}\\
            & \text{abelian part of $X(p)$: $H_1$} &
              \text{anabelian part of $X(p)$: $\pi_1$}\\
            &\text{Hecke's level theory} &\text{our level theory}\\
            & \text{Hecke's decomposition} & \text{our decomposition}\\
            & \dim H^0(X(p), \Omega_{X(p)}^1) & \dim I(\mathcal{L}(X(p)))\\
            & \text{genus of $X(p)$, i.e.,} & \text{dimension of}\\
            & \text{dimension of $J(p)$} & \text{defining ideals}\\
       p=7  & H^0(X(7), \Omega_{X(7)}^1)=V_3 & V_{I(\mathcal{L}(X(7)))}=\mathbf{1}\\
            & \text{degenerate discrete series} & \text{trivial representation}\\
            & J(7) \cong A_3 & I(\mathcal{L}(X(7)))=I_1\\
      p=11  & H^0(X(11), \Omega_{X(11)}^1)=V_5 \oplus V_{10} \oplus V_{11} &
              V_{I(\mathcal{L}(X(11)))}=\mathbf{10}\\
            & \text{degenerate discrete series,} & \text{discrete series}\\
            & \text{discrete series corresponding to $\chi_4$} & \text{corresponding}\\
            & \text{or $\chi_5$ and Steinberg representation} & \text{to $\chi_5$}\\
            & J(11) \cong A_5 \times A_{10} \times A_{11} & I(\mathcal{L}(X(11)))=I_{10}\\
      p=13  & H^0(X(13), \Omega_{X(13)}^1) & V_{I(\mathcal{L}(X(13)))}=\mathbf{1}
              \oplus \mathbf{7} \oplus \mathbf{13}\\
            & =V_{12}^{(1)} \oplus V_{12}^{(2)} \oplus V_{12}^{(3)} \oplus V_{14} & \\
            & \text{three conjugate} & \text{trivial representation,}\\
            & \text{discrete series} & \text{degenerate principal series}\\
            & \text{and principal series} & \text{and Steinberg representation}\\
            & J(13) \cong A_{12}^{(1)} \times A_{12}^{(2)} \times A_{12}^{(3)} \times A_{14}
            & I(\mathcal{L}(X(13)))=I_{1} \cap I_7 \cap I_{13}
\end{matrix}\eqno{(1.1.13)}$$

(2) This leads to a new viewpoint on the last mathematical testament of
Galois by Galois representations arising from the defining ideals of modular
curves, which leads to a connection with Klein's elliptic modular functions.

(3) As an anabelian counterpart of the global Langlands correspondence
for $\text{GL}(2, \mathbb{Q})$ by the \'{e}tale cohomology of modular
curves (Grothendieck's motives and the global Langlands correspondence
for $\text{GL}(2, \mathbb{Q})$) (see Theorem 5.2, i.e. theorem of 
Eichler-Shimura-Deligne-Ihara-Langlands-Piatetskii-Shapiro-Carayol and 
others). Combining these two kinds of correspondences (i.e., Theorem 5.2 
and Theorem 1.1) leads to a complete picture on the correspondence between 
algebraic varieties (both motives and anabelian version) and automorphic 
forms for $\text{GL}(2, \mathbb{Q})$. That is, there are two distinct 
approaches to realize the non-abelian class field theory for 
$\text{GL}(2, \mathbb{Q})$: one is the global Langlands correspondence 
in the context of motives, the other is our correspondence in the context 
of anabelian algebraic geometry.

(4) As an anabelian counterpart of the Shimura-Taniyama-Weil conjecture. Let
$E$ be an elliptic curve over $\mathbb{Q}$. Suppose $\ell$ is a prime number
such that the group of $\ell$-torsion points $E[\ell]$ is irreducible as a
representation of the absolute Galois group
$G_{\mathbb{Q}}=\text{Gal}(\overline{\mathbb{Q}}/\mathbb{Q})$:
$$\overline{\rho}_{E, \ell}: G_{\mathbb{Q}} \rightarrow \text{GL}(2, \mathbb{F}_{\ell}).$$
Then we say that $E[\ell]$ is modular of level $N$ if there exists a normalized
simultaneous eigen-cuspform $f=\sum_{m=1}^{\infty} a_m(f) q^m$ of level dividing
$N$ and a prime ideal $\lambda$ of the integer ring of
$\mathbb{Q}(f):=\mathbb{Q}(a_m(f), m \in \mathbb{N})$ containing $\ell$ such that
$$p+1-\# E(\mathbb{F}_p) \equiv a_p(f) \quad \text{(mod $\lambda$)}\eqno{(1.14)}$$
for all primes $p$ not dividing $N$. By Theorem 1.3, we have the following comparison:
$$\begin{array}{ccc}
  \text{motives} & \longleftrightarrow & \text{anabelian algebraic geometry}\\
  \text{linear realization $E[p]$}  & \longleftrightarrow &
  \text{nonlinear realization $I(\mathcal{L}(X(p)))$}\\
  \text{Aut}(E[p]) \cong \text{GL}(2, \mathbb{F}_p) & \longleftrightarrow  &
  \text{Aut}(\mathcal{L}(X(p))) \cong \text{PSL}(2, \mathbb{F}_p)\\
\end{array}\eqno{(1.15)}$$
In particular, Theorem 1.3 shows that the modularity condition (1.14) (in the
context of motives) should be replaced by the following correspondence (in the
context of anabelian algebraic geometry):
$$\begin{array}{ccc}
  \text{the Jacobian multiplier}  & \longleftrightarrow & \text{$p$-th order transformation}\\
  \text{equation of degree $p+1$} &                     & \text{equation of the $j$-function}\\
\end{array}\eqno{(1.16)}$$
which are connected by the fact that both of the two sides of (1.16) have
Galois groups isomorphic to $\text{PSL}(2, \mathbb{F}_p)$. In fact, (1.16)
corresponds to the exceptional case $\ell=p=N$ in the case of (1.14), which
is not included in (1.14). Hence, the two conditions (1.14) and (1.16) are
complementary to each other. Note that the representation theory of
$\text{GL}(2, \mathbb{F}_p)$ plays a crucial role in the proof of the
Shimura-Taniyama-Weil conjecture (see \cite{W} and \cite{CDT}). In
particular, the cuspidal representation of $\text{GL}(2, \mathbb{F}_3)$
appears in Wiles' proof \cite{W}. In our case, Corollary 1.2 shows that,
besides cuspidal representations (i.e., discrete series), much more
representations such as principal series representations as well as
Steinberg representations of $\text{PSL}(2, \mathbb{F}_p)$ do appear
in our theory.

(5) As an anabelian counterpart of Artin's conjecture (Artin motives),
Serre's modularity conjecture (mod $p$ Langlands correspondence), and
the Fontaine-Mazur conjecture ($p$-adic Galois representations coming
from algebraic geometry, i.e., motives, and $p$-adic local Langlands
correspondence). In fact, together with the two-dimensional Artin conjecture
as well as Serre's modularity conjecture, both of them are realizations
by motives (Artin motives and Grothendieck motives), our correspondence
(anabelian algebraic geometric realization) give a complete picture on
representations of finite Galois groups $\text{PSL}(2, \mathbb{F}_p)$.
Together with the work of Eichler, Shimura, Deligne, Serre and Fontaine-Mazur,
our correspondence leads to a complete picture on Galois representations
associated with modular curves of level $p$ (including both $\ell \neq p$
and $\ell=p$), from the perspective of both motives and anabelain algebraic
geometry. In particular, our correspondence is an integral and anabelian
algebraic geometric counterpart of Serre's modularity conjecture which
gives a mod $p$ modular form realization of the mod $p$ Galois representations.

(6) As an anabelian counterpart of the $p$-adic local Langlands correspondence
for $\text{GL}(2, \mathbb{Q}_p)$. Besides the comparison between Theorem 5.2 
and Theorem 1.1 as above, Theorem 1.1 can also be regarded as an anabelian
counterpart of Theorem 8.4, i.e., theorem of Deligne, Langlands and Carayol.
For simplicity, let us consider the case of $\text{GL}(2)$. Let us denote by
$\rho(f)$ the $p$-adic representation of $\text{Gal}(\overline{\mathbb{Q}}/\mathbb{Q})$
attached to the newform $f$ of weight $k \geq 2$ and level $M$. The philosophy
of the $p$-adic Langlands suggests that one should be able to recover this
Galois representation from purely automorphic data associated to $f$. Similarly,
the local $p$-adic Langlands program should relate the local representation
$\rho_p(f):=\rho(f)|_{\text{Gal}(\overline{\mathbb{Q}}_p/\mathbb{Q}_p)}$
to the local automorphic component $\pi_p(f)$ of the representation of
$\text{GL}(2)$ determined by $f$. From the classical theory, we know that
when the level of $f$ is exactly divisible by $p$, the local automorphic
component $\pi_p(f)$ does not contain enough information to isolate
$\rho_p$. Indeed, in this case, $\pi_p(f)$ is always the Steinberg
representation and it is exactly the $\mathcal{L}$ invariant of $f$
that provide the additional information necessary to identify
the local Galois representation associated to $f$. More precisely,
let $p$ be a prime number and $f$ a cuspidal newform of weight
$k \geq 2$ on $\Gamma_0(Np)$ with $(N, p)=1$. For simplicity, let
$f$ be an eigenvector of the Hecke operators with eigenvalues in
$\mathbb{Z}$. We know that the local component $\rho_p$ is a
non-crystalline semistable representation. Thus, via the Fontaine
functor, it corresponds to a filtered $(\varphi, N)$-module $D_{\text{st}}(\rho_p)$
of dimension $2$, and the invariant $\mathcal{L}$ of $f$, denoted by
$\mathcal{L}(f)$, is defined by the parameter of the Hodge filtration on
$D_{\text{st}}(\rho_p)$ (see \cite{Breuil2004}). The form $f$ is also
associated with an automorphic representation
$\pi=\otimes_{\ell}^{\prime} \pi_{\ell}$ of $\text{GL}(2, \mathbb{A}_{\mathbb{Q}})$,
where $\mathbb{A}_{\mathbb{Q}}$ denotes the adeles of $\mathbb{Q}$. The local
component $\pi_p$ is, up to an unramified torsion, the special representation
or Steinberg representation of $\text{GL}(2, \mathbb{Q}_p)$.

  The significance of Corollary 1.2 comes from that it gives a comparison
between Fontaine's theory on $p$-adic Galois representations and our theory
on $\mathbb{Q}(\zeta_p)$-rational Galois representations:
$$\begin{matrix}
 &\text{Fontaine's theory} &\text{representation theory} &\text{our theory}\\
 &\text{de Rham} &\text{supercuspidal representations} & \rho_{11, 10}\\
 &\text{representations} &\text{(cuspidal or discrete series)} & \text{(see (2.4.31) et seq.)}\\
 &\text{semistable} &\text{Steinberg representations} & \rho_{13, 13}\\
 &\text{representations} & \text{(special series)} & \text{(see (2.4.42) et seq.)}\\
 &\text{crystalline} &\text{principal series} & \rho_{13, 7}\\
 &\text{representations} &\text{(degenerate principal series)} & \text{(see (2.4.42) et seq.)}\\
\end{matrix}\eqno{(1.17)}$$
In particular, it leads to the following comparison:
$$\begin{matrix}
 &\text{motives} &\longleftrightarrow &\text{anabelian algebraic geometry}\\
 &\text{cuspidal eigenforms} &\longleftrightarrow &\text{defining ideals of modular}\\
 &\text{of level $p$} &  & \text{curves of level $p$}\\
 &f \in S_k(\Gamma_0(p)) &\longleftrightarrow & I(\mathcal{L}(X(p)))\\
 &\text{$\pi_p(f)$: only} &\longleftrightarrow & \text{$\pi_p$: trivial representations,}\\
 &\text{Steinberg} &  & \text{principal series, Steinberg}\\
 &\text{representations} & & \text{representations, and discrete series}\\
 &\text{$\rho_p(f)$: only} & \longleftrightarrow & \text{$\rho_p$: ``crystalline'',}\\
 &\text{non-crystalline} & & \text{``semistable'',}\\
 &\text{semistable} &  & \text{and ``de Rham''}\\
 &\text{representations} &  & \text{representations}\\
 &\text{Hodge-Tate weights $(0, k-1)$} & \longleftrightarrow & \text{our decomposition Theorem}\\
 &\text{and $\mathcal{L}$-invariants} & & \text{1.1 and Corollary 1.2}\\
 &\text{$\rho_p(f)$ can {\it not} be} & \longleftrightarrow & \text{$\rho_p$ can be}\\
 &\text{recovered from $\pi_p(f)$} &  & \text{recovered from $\pi_p$}\\
\end{matrix}\eqno{(1.18)}$$
Hence, from the viewpoint of (1.18), for the same level $p$, the defining
ideals of modular curves $X(p)$ can provide much more representations than
the cuspidal eigenforms $S_k(\Gamma_0(p))$.

  As an important step towards the proof of the Fontaine-Mazur conjecture,
Emerton proved the theorem 8.6 by the $p$-adic local Langlands correspondence 
for $\text{GL}(2, \mathbb{Q}_p)$. As an anabelian counterpart of Emerton's 
theorem 8.6, Theorem 1.1 and Theorem 1.3 imply the following:

\textbf{Corollary 1.4.} {\it Let} $\rho_p: G_{\mathbb{Q}} \rightarrow
\text{PSL}(2, \mathbb{F}_p)$ {\it be a Galois representation. Suppose
that}

(1) {\it $\rho_p$ appears in the defining ideals of modular curves $X(p)$.}

(2) {\it $\rho_p|_{G_{\mathbb{Q}}}$ has a decomposition as given by Theorem 1.1.}

{\it Then $\rho_p$ is modular in the sense of Theorem 1.3.}

  In (8.1), we give a more detailed comparison between Emerton's
theorem 8.6 and Corollary 1.4. In particular, Emerton's theorem only
studies two-dimensional $p$-adic Galois representations, while
Corollary 1.4 includes higher dimensional $\mathbb{Q}(\zeta_p)$-rational
Galois representations.

(7) Comparison with the Macdonald correspondence. In 1980, Macdonald (see
\cite{Macdonald}) established an explicit bijection between the isomorphism
classes of the irreducible representations of $\text{GL}(n, k)$, where $k$
is a finite field, and inertia equivalence classes of admissible tamely
ramified $n$-dimensional Weil-Deligne representations of $W_K$, where $K$ is
a non-archimedean local field with residue field $k$ and $W_K$ the absolute
Weil group of $K$. More precisely, let $\Pi(\text{GL}(n, k))$ denote the set 
of equivalence classes of irreducible complex representations of $\text{GL}(n, k)$ 
and $\Phi(\text{GL}(n))$ denote the set of equivalence classes of $n$-dimensional
admissible representations of $W_K$. Let $\Phi^{t}(\text{GL}(n))$ denote the 
set of $\phi \in \Phi(\text{GL}(n))$ which are tamely ramified. Let 
$\Phi_I^{t}(\text{GL}(n))$ denote the set of $I$-equivalent classes (see 
section nine and \cite{Macdonald} for more details). Macdonald established 
a canonical bijection of $\Pi(\text{GL}(n, k))$ onto $\Phi_I^{t}(\text{GL}(n))$, 
by showing that both of theses sets are parameterized by the same set $P_n(\Gamma)$.
As a corollary of Theorem 1.1, we give a comparison between the Macdonald
correspondence in the case of $n=2$, $k=\mathbb{F}_p$ and our correspondence
(1.9) as follows:
$$\begin{array}{ccc}
  \text{the Macdonald correspondence} & \longleftrightarrow & \text{our correspondence}\\
  \pi_{\lambda} \in \Pi(\text{GL}(n, k)) & \longleftrightarrow & \pi_p\\
  \text{complex representations} & \longleftrightarrow &
  \text{$\mathbb{Q}(\zeta_p)$-rational representations}\\
  (\phi_{\lambda}) \in \Phi_I^t(\text{GL}(n)) & \longleftrightarrow & \rho_p\\
  \text{Weil-Deligne representations} & \longleftrightarrow & \text{Galois representations}\\
  \text{the set $P_n(\Gamma)$} & \longleftrightarrow & \text{the defining ideal $I(\mathcal{L}(X(p)))$}\\
  \text{both are parameterized} & \longleftrightarrow & \text{both are parameterized}\\
  \text{by the same $P_n(\Gamma)$} & & \text{by the same $I(\mathcal{L}(X(p)))$}\\
\end{array}\eqno{(1.19)}$$

(8) As an ideal theoretic (i.e. nonlinear) counterpart of Grothendieck's
section conjecture and an ideal theoretic (i.e. nonlinear) reification
of ``arithmetic theory of $\pi_1$'' expected by Weil in \cite{Weil1938}
for modular curves (anabelian algebraic geometry).

(9) As an anabelian counterpart of the theory of Kubert-Lang and Mazur-Wiles
on the cuspidal divisor class groups and the Eisenstein ideals of modular
curves (cuspidal Iwasawa theory and its anabelian counterpart). In fact,
our decomposition theorem (Theorem 1.1 and Corollary 1.2), which can be
regarded as an anabelian counterpart of the eigenspace decomposition
appearing in the context of both cyclotomic Iwasawa theory and cuspidal
Iwasawa theory (the theory of Kubert-Lang and Mazur-Wiles) (see (1.22)
and (1.27)). In particular, Corollary 1.2, i.e., the dimension formula
for defining ideals of modular curve can be regarded as an anabelian
counterpart of the cuspidal divisor class number formula (see (1.23))
as well as a modular curve counterpart of the non-abelian Herbrand-Ribet
theorem (see (1.26)).

  In particular, three different perspectives on non-abelian phenomena
in number theory: the Langlands program and motives, non-commutative
Iwasawa theory and anabelian algebraic geometry fit together within
the single framework of Theorem 1.1, Corollary 1.2 and Theorem 1.3,
the difference lying only in the degree of non-commutativity or
nonlinearity that is preserved in its study.
$$\begin{matrix}
 &\text{Jacobian $J(X)$}      &\text{Selmer variety (cohomology)} &\text{defining ideal}\\
 &\text{(abelian variety/$\mathbb{Q}$)} &\text{(algebraic variety/$\mathbb{Q}_p$)}
                              &\text{(algebraic variety/$\mathbb{Q}$)}\\
 &\text{Tate module of}       &\text{linearized} & \text{ideal theoretic} \\
 &\text{the Jacobian}         &\text{section conjecture} & \text{counterpart of the}\\
 &                            &                          & \text{section conjecture}\\
 &\text{abelianization of $\pi_1$} &\text{unipotent quotient of $\pi_1$}
 & \text{truly $\pi_1$}\\
 &\text{abelian} &\text{nilpotent}           &\text{utmost non-abelian}\\
 &\text{motives} &\text{motivic fundamental} &\text{fundamental groups}\\
 &               &\text{groups}              &                         \\
 &\text{$H_1$ (Grothendieck)} & \text{motivic $\pi_1$ (Deligne)} & \text{$\pi_1$ (ours)}\\
 &\text{linear} &\text{mildly nonlinear} &\text{fully nonlinear}\\
 &\text{abelian} & \text{non-abelian}  &\text{anabelian}\\
 &\text{Iwasawa theory} & \text{Iwasawa theory} &\text{Iwasawa theory}\\
 &\text{Chabauty-Coleman} &\text{Chabauty-Kim} & \text{ours}\\
 &                         &\text{(non-abelian Chabauty)} &
\end{matrix}\eqno{(1.20)}$$
The above is a comparison among homology groups, unipotent motivic
fundamental groups and the ideal theoretic (i.e. nonlinear)
reification of the fundamental groups. In particular, we give a
comparison between motives for $\text{GL}(2, \mathbb{Q})$ and anabelian
algebraic geometry for $\text{GL}(2, \mathbb{Q})$ by the following
anabelianization:
$$\begin{array}{ccc}
  \text{abelian side} & \longrightarrow & \text{anabelian side}\\
  \text{Jacobian varieties of} & \longrightarrow & \text{defining ideals of}\\
  \text{modular curves} &  &  \text{modular curves}\\
  \text{a reification of the} & \longrightarrow & \text{a reification of the}\\
  \text{homology group $H_1$} &  & \text{fundamental group $\pi_1$}\\
  \text{abelian properties of} & \longrightarrow & \text{anabelian properties of}\\
  \text{the modular curves} & & \text{the modular curves}\\
  \text{partial information} & \longrightarrow & \text{complete information}\\
  \text{of modular curves} &  &  \text{of modular curves}
\end{array}\eqno{(1.21)}$$
as an ideal theoretic (i.e. nonlinear) realization of the following picture:
$$\text{Jacobian} \longrightarrow \text{non-abelian Jacobian (Weil, Grothendieck)}.$$
$$\begin{array}{ccc}
  \text{abelian objects} & \longrightarrow & \text{anabelian objects}\\
  \text{motives $H_1$} & \longrightarrow & \text{anabelian geometry $\pi_1$}\\
  \text{Tate conjecture for} & \longrightarrow & \text{Grothendieck's anabelian}\\
  \text{Jacobian varieties of curves} &    & \text{Tate conjecture for curves}\\
  \text{Hecke's decomposition} & \longrightarrow & \text{our decomposition}\\
  \text{of $H^0(X(p), \Omega_{X(p)}^1)$} &    & \text{of $I(\mathcal{L}(X(p)))$}\\
  \text{comparison between de Rham} & \longrightarrow & \text{comparison between de Rham}\\
  \text{cohomology and \'{e}tale cohomology} &  & \text{cohomology and defining ideals}\\
  \text{Hecke's theory v.s.} & \longrightarrow & \text{Hecke's theory v.s.}\\
  \text{Drinfeld-Deligne-Lusztig's theory} &     & \text{our theory}\\
  \text{global Langlands correspondence} & \longrightarrow &
  \text{our correspondence}\\
  \text{for $\text{GL}(2, \mathbb{Q})$ by the \'{e}tale} &  &
  \text{for $\text{GL}(2, \mathbb{Q})$ by the defining}\\
  \text{cohomology of modular curves} &  & \text{ideals of modular curves}\\
  \text{Galois representations arising} & \longrightarrow &
  \text{Galois representations arising}\\
  \text{from the $p$-torsion points of} &  & \text{from the defining ideals of}\\
  \text{elliptic curves $E/\mathbb{Q}$} &  & \text{modular curves $X(p)/\mathbb{Q}$}\\
  \rho_{E, p}: \text{Gal}(\overline{\mathbb{Q}}/\mathbb{Q}) \rightarrow \text{Aut}(E[p])
  & \longrightarrow & \rho_p: \text{Gal}(\overline{\mathbb{Q}}/\mathbb{Q}) \rightarrow
  \text{Aut}(\mathcal{L}(X(p)))\\
  \text{$p$-adic Galois representations} & \longrightarrow &
  \text{$\mathbb{Q}(\zeta_p)$-rational Galois}\\
  \text{coming from algebraic geometry} &  & \text{representations coming from}\\
  \text{(Fontaine-Mazur conjecture)} & & \text{anabelian algebraic geometry}\\
 \text{Hecke operators act on the} & \longrightarrow &
  \text{$\text{PSL}(2, \mathbb{F}_p)$ acts on the}\\
  \text{divisor groups of modular curves} &   &
  \text{the defining ideals $I(\mathcal{L}(X(p)))$}\\
  \text{which induces an} &   & \text{of modular curves $X(p)$}\\
  \text{endomorphism of the Jacobian} &   & \text{as the automorphism group}\\
  \text{Hecke's theory} & \longrightarrow & \text{our theory}\\
  \text{cuspidal divisor class} & \longrightarrow & \text{the locus of}\\
  \text{groups on modular curves} &  & \text{modular curves}\\
  \text{Mazur's Eisenstein ideals} & \longrightarrow & \text{our defining ideals}\\
  \text{for modular curves} &  & \text{for modular curves}\\
  \text{Kubert-Lang theory,} & \longrightarrow & \text{our theory}\\
  \text{Ribet-Mazur-Wiles theory} &    &  \\
  \text{cusps, modular symbols} & \longrightarrow & \text{the locus of}\\
  \text{and Beilinson elements} &  & \text{modular curves}\\
  \text{in $K_2$ of modular curves}  &  &  \\
  \text{$\text{GL}(2, \mathbb{F}_p)$ acts on the} & \longrightarrow &
  \text{$\text{PSL}(2, \mathbb{F}_p)$ acts on the}\\
  \text{Siegel units and $K_2(Y(p))$} &  & \text{defining ideals $I(\mathcal{L}(X(p)))$}\\
  \text{Fukaya-Kato-Sharifi theory} & \longrightarrow & \text{our theory}\\
\end{array}\eqno{(1.22)}$$

  In particular, we should point out that one of the essential distinction
between the abelian objects and the anabelian objects is the following
comparison:
$$\begin{array}{ccc}
  \text{irreducible Galois represen-} & \longrightarrow & \text{Galois representations}\\
  \text{tations $\rho_f$ attached to the} &  & \text{$\rho_p$ arising from the}\\
  \text{cusp forms $f$ (Eichler-Shimura,} &  & \text{defining ideals of}\\
  \text{Deligne, Deligne-Serre, Serre's} &  & \text{modular curves}\\
  \text{modularity conjecture and} &  & X(p)\\
  \text{Fontaine-Mazur conjecture)} &  & \text{(ours)}\\
  \downarrow & & \downarrow\\
  \rho_f: \text{Gal}(\overline{\mathbb{Q}}/\mathbb{Q}) \rightarrow
  \text{GL}(2, \mathbb{F}_p)
  & \longrightarrow & \rho_p: \text{Gal}(\overline{\mathbb{Q}}/\mathbb{Q})
  \rightarrow \text{PSL}(2, \mathbb{F}_p)\\
  \text{occurs in the torsion}  &  & \text{occurs in the}\\
  \text{points of the Jacobian} &  & \text{defining ideal of}\\
  \text{of a modular curve} &   &  \text{a modular curve}\\
  \downarrow & & \downarrow\\
  \text{the geometry of} & \longrightarrow & \text{the geometry of}\\
  \text{a modular curve near $\infty$} &  &\text{a full modular curve}\\
  \text{(local and partial information)} &  & \text{(global and complete information)}\\
\end{array}\eqno{(1.23)}$$
$$\begin{array}{ccc}
  \text{Eisenstein ideals} & \longrightarrow & \text{defining ideals}\\
  \downarrow & & \downarrow\\
  \text{congruences between cusp} & \longrightarrow &\text{transformation equations}\\
  \text{forms and Eisenstein seres} &  & \text{of elliptic functions}\\
  \text{e.g., Ramanujan's congruence} & \longrightarrow &
  \text{e.g., Galois and Klein's}\\
  \text{$\tau(n) \equiv \sum_{0<d | n} d^{11}$ (mod $691$)} &  &
  \text{special transformation equations}\\
  \text{for Ramanujan's $\tau$-function} &  & \text{of elliptic functions}\\
  \downarrow & & \downarrow\\
  \text{residual reducible, indecom-} & \longrightarrow & \text{Galois representations}\\
  \text{posable (upper triangular)} &  & \text{arising from the}\\
  \text{Galois representations} &  & \text{defining ideals of}\\
  \text{attached to cusp forms} &  & \text{modular curves $X(p)$}\\
  \text{(Serre, Ribet, Mazur-Wiles)} &  & \text{(ours)}\\
  \downarrow & & \downarrow\\
  \text{the geometry of} & \longrightarrow & \text{the geometry of}\\
  \text{a modular curve near $\infty$} &  &\text{a full modular curve}\\
  \text{(local and partial information)} &  & \text{(global and complete information)}\\
\end{array}\eqno{(1.24)}$$
The left hand side come from the abelian object, which gives the geometry
of a modular curve near the cusps and thus leads to the Jacobian of the
modular curve. On the other hand, the right hand side come from the
non-abelian object, which gives the geometry of a full modular curve,
i.e., local and partial information versus global and complete information.
It is well-known that the Jacobian variety has always been the
corner-stone in the analysis of algebraic curves and compact Riemann
surfaces. Its power lies in the fact that it abelianizes the curve
and is a reification of the homology group $H_1$. In particular, the
significance of the Eisenstein ideal of a modular curve comes from
the fact that it is the annihilator of the cuspidal divisor class
group, i.e., the rational torsion subgroup in the Jacobian variety of
this modular curve. Correspondingly, the significance of the defining
ideal of a modular curve comes from the fact that it is an anabelian
counterpart of the Eisenstein ideal of this modular curve, i.e., it is
the annihilator of the locus of this (whole) modular curve. It is the
anabelianization of the Jacobian variety of this modular curve. Thus,
it is a reification of the fundamental group $\pi_1$. In fact, by the
uniformization theorem, any compact complex Riemann surface $X$ of genus
$g \geq 2$ can be represented by $X=\mathbb{H}/\pi_1(X)$, where $\pi_1(X)$
is the fundamental group of $M$. Therefore, given a hyperbolic compact
Riemann surface $X$ ($g \geq 2$), it is equivalent to give its fundamental
group $\pi_1(X)$. On the other hand, given a Riemann surface $X$, it is
equivalent to give its defining ideal in some projective space. Hence,
given the defining ideal for a hyperbolic Riemann surface, it is equivalent
to give its fundamental group. Thus, the defining ideal of a hyperbolic
curve is utmost nonabelian (as a fundamental group) and nonlinear (as a
system of highly nonlinear algebraic equations). In order to show this
viewpoint, we give the following comparisons (1.25), (1.26), (1.27),
(1.28), (1.29), (1.30), (1.31), (1.32), (1.33) and (1.34). On the other
hand, the Eisenstein ideals also come from congruences between cusp
forms and Eisenstein series. A typical example is given by Ramanujan's
congruence: if one expands $\Delta(z)=q \prod_{n=1}^{\infty} (1-q^n)^{24}$
as a series $\sum_{n=1}^{\infty} \tau(n) q^n$ with $q=e^{2 \pi i z}$,
$\text{Im}(z)>0$, then $\tau(n) \equiv \sum_{0<d | n} d^{11}$ (mod $691$).
Correspondingly, the defining ideals also come from the transformation
equations of elliptic functions of order $p$. The typical examples are
given by Galois and Klein's special transformation equations of elliptic
functions with order seven and eleven, and our transformation equation of
elliptic functions with order thirteen (see (1.24)).

  We have the following anabelianization from motives to anabelian algebraic
geometry:
$$\begin{array}{ccc}
  \begin{matrix}
  \text{the cusps on}\\
  \text{the modular curve}
  \end{matrix}
 &\left\{\aligned
  &\text{\quad Manin-Drinfeld theorem:}\\
  &\text{(1) cuspidal divisor class groups}\\
  &\text{\qquad in the modular Jacobian}\\
  &\text{(2) homology relative to the}\\
  &\text{\quad cusps of the modular curve}\\
  \endaligned\right\}
 &\begin{matrix}
  \text{the Eisenstein}\\
  \text{ideal}
  \end{matrix}\\
  \downarrow & \downarrow & \downarrow\\
  \begin{matrix}
  \text{the full}\\
  \text{modular curve}
  \end{matrix}
 &\longrightarrow \begin{matrix}
  \text{the locus of}\\
  \text{the modular curve}
  \end{matrix} \longrightarrow
 &\begin{matrix}
  \text{the defining}\\
  \text{ideal}
  \end{matrix}
\end{array}\eqno{(1.25)}$$
$$\begin{array}{ccc}
  \begin{matrix}
  \text{the geometry of modular}\\
  \text{curves near the cusps}
  \end{matrix}
 &\longrightarrow
 &\begin{matrix}
  \text{the arithmetic of}\\
  \text{cyclotomic fields}
  \end{matrix}\\
  \downarrow &  & \downarrow\\
  \begin{matrix}
  \text{the geometry of the full}\\
  \text{modular curves}
  \end{matrix}
 &\longrightarrow
 &\begin{matrix}
  \text{the arithmetic of division fields}\\
  \text{of elliptic curves over $\mathbb{Q}$}
  \end{matrix}
\end{array}\eqno{(1.26)}$$
$$\begin{array}{ccc}
  \text{Eisenstein ideals} & \longleftrightarrow &
  \text{defining ideals $I(\mathcal{L}(X(p)))$}\\
  \text{which annihilate } & \longleftrightarrow & \text{which annihilate}\\
  \text{cuspidal divisor class groups} & \longleftrightarrow &
  \text{the locus $\mathcal{L}(X(p)) \subset \mathbb{P}^{(p-3)/2}$}\\
  \text{in the Jacobian of modular curves} &  & \text{of modular curves $X(p)$}\\
  \text{modular Jacobians $H_1$} & \longleftrightarrow & \pi_1\\
  \text{motives} & \longleftrightarrow & \text{anabelian algebraic geometry}\\
\end{array}\eqno{(1.27)}$$
Moreover, we have the following comparison:
$$\begin{array}{ccc}
  \text{$\text{GL}(1)$: cyclotomic fields} & \longleftrightarrow &
  \text{cuspidal divisor class groups}\\
                                           &  & \text{in the modular Jacobian ($H_1$)}\\
  \updownarrow &  & \updownarrow\\
  \text{$\text{GL(2)}$: division fields}   & \longleftrightarrow &
  \text{the locus of}\\
  \text{of elliptic curves} &  & \text{modular curves ($\pi_1$)}\\
\end{array}\eqno{(1.28)}$$
i.e.,
$$\begin{array}{ccc}
  \text{cyclotomic Iwasawa theory} & \longleftrightarrow &
  \text{cuspidal Iwasawa theory}\\
  \updownarrow &  & \updownarrow\\
  \text{elliptic Iwasawa theory}   & \longleftrightarrow &
  \text{anabelian Iwasawa theory}
\end{array}\eqno{(1.29)}$$
where
$$\begin{array}{ccc}
  \text{cuspidal divisor class groups} & \longleftrightarrow & \text{the locus of}\\
  \text{in the modular Jacobian ($H_1$)} &  & \text{modular curves ($\pi_1$)}\\
  \updownarrow &  & \updownarrow\\
  \text{ideals of quadratic relations} & \longleftrightarrow &
  \text{ideals of quartic relations}
\end{array}$$
In other words,
$$\begin{array}{ccc}
  \text{abelian number fields} & \longleftrightarrow & \text{motives}\\
  \updownarrow &  & \updownarrow\\
  \text{non-abelian number fields} & \longleftrightarrow & \text{anabelian algebraic geometry}\\
  \text{number theory} & \longleftrightarrow & \text{modular curves}\\
  \text{class number formula} & \longleftrightarrow & \text{cuspidal divisor class}\\
  \text{for cyclotomic fields (8.1)} & & \text{number formula (8.3)}\\
  \text{$B_{1, \chi} (L\text{-values})$ in (8.1)} & \longleftrightarrow &
  \text{$B_{2, \chi} (L\text{-values})$ in (8.3)}\\
  \updownarrow &  & \updownarrow\\
  \text{non-abelian class number} & \longleftrightarrow & \text{dimension formula}\\
  \text{formula for division fields} &  & \text{for defining ideals}\\
  \text{of elliptic curves} &  & \text{of modular curves}\\
  \text{(see \S 8.4 for the conjectural form)} &  & \text{(Corollary 1.2)}\\
  \text{division equations} & \longleftrightarrow & \text{transformation equations}\\
  \text{of elliptic functions} &  & \text{of elliptic functions}
\end{array}\eqno{(1.30)}$$
$$\begin{array}{ccc}
  \text{number fields} & \longleftrightarrow & \text{modular curves}\\
  \text{cyclotomic fields} & \longleftrightarrow & \text{modular curves}\\
                           &                     & \text{near the cusps}\\
  \text{Gal}(\mathbb{Q}(\zeta_p)/\mathbb{Q}) \cong \mathbb{F}_p^{*} & \longleftrightarrow
                       & \text{$C(p)/ \pm 1$, where}\\
                       &                         & \text{$C(p)$ is a Cartan group}\\
  \updownarrow &  & \updownarrow\\
  \text{division fields of} & \longleftrightarrow & \text{the full}\\
  \text{elliptic curves} &  & \text{modular curves}\\
  \text{Gal}(\mathbb{Q}(E[p])/\mathbb{Q}) \cong \text{GL}(2, \mathbb{F}_p)
  & \longleftrightarrow & \text{PSL}(2, \mathbb{F}_p)\\
\end{array}\eqno{(1.31)}$$

  The left hand side comes from some specific number fields (number theory),
the right hand side comes from modular curves (arithmetic geometry). The
corresponding ideals appear as follows:
$$\begin{array}{ccc}
  \text{Stickelberger ideals} & \longleftrightarrow & \text{Eisenstein ideals}\\
  \updownarrow &  & \updownarrow\\
  \text{Stickelberger elements} & \longleftrightarrow & \text{defining ideals}
\end{array}\eqno{(1.32)}$$
In particular, we have the following four kinds of theorems/theories:
$$\begin{array}{ccc}
  \text{Herbrand-Ribet} & \longleftrightarrow & \text{theory of Kubert-Lang}\\
  \text{theorem} &  & \text{and Mazur-Wiles}\\
  \updownarrow &  & \updownarrow\\
  \text{non-abelian Herbrand-Ribet} & \longleftrightarrow & \text{Theorem 1.1}\\
  \text{(Prasad-Shekhar theorem)} & & \text{and Corollary 1.2}\\
\end{array}\eqno{(1.33)}$$
More precisely, that is:
$$\begin{array}{ccc}
  \text{eigenspace decomposition} & \longleftrightarrow & \text{eigenspace decomposition}\\
  \text{of the $p$-primary part} &  & \text{of the $p$-primary part of}\\
  \text{of the ideal class groups} &  & \text{the cuspidal divisor class}\\
  \text{for cyclotomic fields $\mathbb{Q}(\zeta_p)$} &  &\text{groups on $X(p)$}\\
  \text{(see (8.1.1) and (8.1.2))} &  & \text{(the theory of Kubert-Lang)}\\
  \updownarrow &  & \updownarrow\\
  \text{eigenspace decomposition} & \longleftrightarrow & \text{eigenspace decomposition}\\
  \text{of the $p$-primary part} &  & \text{of the space $V_{I(\mathcal{L}(X(p)))}$} \\
  \text{of the ideal class groups} & & \text{associated with the defining}\\
  \text{for division fields $\mathbb{Q}(E[p])$} &  & \text{ideals $I(\mathcal{L}(X(p)))$ of}\\
  \text{of elliptic curves $E$ over $\mathbb{Q}$} &  &  \text{modular curves $X(p)$}\\
  \text{(see \S 8.4)} &  & \text{(Theorem 1.1 and Corollary 1.2)}\\
\end{array}\eqno{(1.34)}$$

  Hence, for specific number fields and modular curves, we have both the abelian
side and the non-abelian side in the sense of Weil \cite{Weil1938}:

(1) abelian side: cyclotomic fields versus cuspidal divisor class groups
    on modular curves;

(2) non-abelian side: division fields of elliptic curves versus defining
    ideals of modular curves.

  This paper consists of twelve sections. In section two, we give Weil's
Rosetta stone for geometric representations of $\text{PSL}(2, \mathbb{F}_p)$.
This includes: (1) Hecke's work on the analytic theory and de Rham realization
of representations for $\text{PSL}(2, \mathbb{F}_p)$, i.e., the first aspect
of the Rosetta stone: Riemann surfaces. (2) Drinfeld and Deligne-Lusztig's
work on the algebraic theory and $\ell$-adic realization of representations
for $\text{PSL}(2, \mathbb{F}_p)$, i.e., the second aspect of the Rosetta
stone: curves over $\mathbb{F}_q$. (3) The arithmetic theory and
$\mathbb{Q}(\zeta_p)$-rational realization of representations for
$\text{PSL}(2, \mathbb{F}_p)$ by the defining ideals of modular curves
$X(p)$, i.e., the third aspect of the Rosetta stone: number fields. All of
them are related with the modular curves. This includes the proof of the
Theorem 1.1 and Corollary 1.2. As a consequence, we provide a comparison
between Hecke's decomposition formulas and our decomposition formulas, i.e.,
motives versus anabelian geometry. In section three, we revisit the division
of elliptic functions, the transformation of elliptic functions and their
Galois groups. This includes the division equations of elliptic functions,
their Galois groups and their relation with transformation equations, the
theory of the special transformation equations, the historic background of
transformation theory of elliptic functions, the Galois groups of the special
transformation equations, Galois' theorem, the resolvent of the seventh degree,
and the two resolvents of the eleventh degree. In section four, we study the
Galois representations arising from the defining ideals $I(\mathcal{L}(X(p)))$
of modular curves $X(p)$, their surjective realization and their modularity.
In particular, we give the proof of the Theorem 1.3. In particular, for $p=7$,
$11$ and $13$, we give the explicit construction, respectively. In section
five, we show that our theory is an anabelian counterpart of the global
Langlands correspondence for $\text{GL}(2, \mathbb{Q})$ by the \'{e}tale
cohomology of modular curves. In section six, we show that our theory is
an anabelian counterpart of the Shimura-Taniyama-Weil conjecture. In section 
seven, we show that our theory is an anabelian counterpart of Artin's conjecture, 
Serre's modularity conjecture and the Fontaine-Mazur-Langlands conjecture. In 
section eight, we we show that our theory is an anabelian counterpart of the
$p$-adic local Langlands corresponde3nce for $\text{GL}(2, \mathbb{Q}_p)$.
In section nine, we give a comparison between the Macdonald correspondence
and our correspondence. In section ten, we show that our theory is an ideal 
theoretic (i.e. nonlinear) counterpart of Grothendieck's section conjecture 
and an ideal theoretic (i.e. nonlinear) reification of ``arithmetic theory 
of $\pi_1$'' expected by Weil for modular curves. In section eleven, we show 
that our theory is also an anabelian counterpart of the theory of Kubert-Lang 
and Mazur-Wiles on the cuspidal divisor class groups and the Eisenstein ideals 
of modular curves. In particular, the defining ideal of a modular curve is an 
anabelian counterpart of the Eisenstein ideal of this modular curve. In section 
twelve, we give some relations with other various topics.

\textbf{Acknowledgements}. The author would like to thank P. Deligne and
J.-P. Serre for their very detailed and helpful comments. The author thanks
B. Mazur for his very helpful comments. He author thanks Dick Gross for
pointing out the reference \cite{Gro87}. He thanks I. Fesenko for helpful
comments.

\begin{center}
{\large\bf 2. A Rosetta stone for geometric representations of
              $\text{PSL}(2, \mathbb{F}_p)$}
\end{center}

  The theory of elliptic functions in Gauss, Abel, and Jacobi is
different from the Weierstrass form of this theory. These two different
presentations of the theory of elliptic functions are unified by Klein's
level theory. Fricke gave the analytic, algebraic and arithmetic foundations
of this level theory created by Klein as the ordering principle of the whole
treatment of elliptic functions, especially in order to place the newer form
of the theory, originating from Weierstrass, in the correct relation to the
older form, chiefly created by Jacobi (see \cite{Fricke}). This leads to the
analytic, algebraic and arithmetic theory of modular curves with higher levels
respectively, which plays a central role in the development of number theory
and arithmetic geometry. In particular, it leads to two distinct geometric
realizations of representations for $\text{PSL}(2, \mathbb{F}_p)$: analytic
theory, de Rham cohomology and complex representations coming from modular
curves over $\mathbb{C}$ (Hecke (see \cite{Hec1}, \cite{Hec2}, \cite{Hec3},
\cite{Hec4} and \cite{Hec5}), Frobenius-Schur (see \cite{Fr1} and \cite{S2})),
and algebraic theory, $\ell$-adic cohomology and $\ell$-adic representations
coming from Drinfeld-Deligne-Lusztig curves over $\mathbb{F}_q$ and modular
curves over $\mathbb{Z}_p$ (Drinfeld (see \cite{Dr1}, \cite{Dr2}), Deligne-Lusztig
(see \cite{DeLu})). Both of them can be considered in the context of
Grothendieck motives. That is, they come from the de Rham realization and
the $\ell$-adic realization of Grothendieck motives, respectively.

  By Theorem 1.1, we give the third geometric realization of representations
for $\text{PSL}(2, \mathbb{F}_p)$: arithmetic theory, defining ideals
$I(\mathcal{L}(X(p)))$ and rational representations over the cyclotomic
field $\mathbb{Q}(\zeta_p)$ coming from the rational models $\mathcal{L}(X(p))$
of modular curves over $\mathbb{Q}$. This leads to a new viewpoint on the
last mathematical testament of Galois by Galois representations arising from
the defining ideals of modular curves, which leads to a connection with Klein's
elliptic modular functions, i.e., the $j$-functions. It is a nonlinear and
anabelian counterpart of the global Langlands correspondence (Eichler-Shimura,
Deligne, Langlands and Carayol) between motives and automorphic representations
for $\text{GL}(2, \mathbb{Q})$, namely a correspondence among the $\ell$-adic
\'{e}tale cohomology of modular curves over $\mathbb{Q}$, i.e., Grothendieck
motives ($\ell$-adic system), automorphic representations of
$\text{GL}(2, \mathbb{Q})$ and $\ell$-adic representations. Both of them are
associated with the same object: modular curves over $\mathbb{Q}$, but with
distinct geometrical realizations and Galois representations: one comes from
anabelian algebraic geometry (defining ideals) and $\mathbb{Q}(\zeta_p)$-rational
representations, the other comes from motives and automorphic representations
($\ell$-adic representations). It can also be regarded as a counterpart of
the global Langlands correspondence for function fields (Drinfeld's analogue
for Eichler-Shimura) between motives and automorphic representations for
$\text{GL}(2, F(X))$ for a function field $F(X)$ over a curve $X$ over
$\mathbb{F}_q$, namely a correspondence among the $\ell$-adic \'{e}tale
cohomology of Drinfeld modular curves over $F(X)$, i.e., $\ell$-adic system,
automorphic representations of $\text{GL}(2, F(X))$ and $\ell$-adic
representations.

  In particular, these three kinds of geometric realizations provide a
complete dictionary (1.1) according to the following three parallel
objects. Moreover, (1.1) provides two distinct kinds of comparisons.

(1) The first comparison: $H_{\text{\'{e}t}}^1$ $\leftrightarrow$
    $H_{\text{dR}}^1$ (distinct realizations of motives).

(2) The second comparison: $J(X(p))$ $\leftrightarrow$ $\mathcal{L}(X(p))$
   (motives versus anabelian algebraic geometry as well as linear structure
   versus nonlinear structure).

  The first comparison between $H_{\text{dR}}^1$ and $H_{\text{\'{e}t}}^1$
coming from de Rham cohomology and $\ell$-adic \'{e}tale cohomology is in the
framework of motives, which can be considered as a counterpart of Grothendieck's
mysterious functor. The second comparison between the Jacobian variety
$J(p):=J((X(p))$ of modular curves $X(p)$ and the locus $\mathcal{L}(X(p))$
of modular curves $X(p)$ comes from motives (linear structure) and anabelian
algebraic geometry (defining ideals, nonlinear structure), respectively. In
fact, the second comparison is much more important than the first one because
it provides a framework which unifies both motives and anabelian algebraic
geometry as well as a comparison between them.

  The significance of the dictionary (1.1) comes from that there are three
different languages: the Riemannian theory of algebraic functions, the
Galoisian theory of algebraic functions over a Galois field, and the
arithmetic theory of algebraic numbers. One can bridge the gap between
the right column and the left column, i.e., the arithmetic and the
Riemannian using the middle column, i.e., the Galoisian curves over finite
fields as the intermediary, so as to transfer constructions from one side
to the other. This dictionary (1.1) has a long history. In the 1850s, Riemann
initiated a revolution in algebraic geometry by interpreting algebraic
curves as surfaces covering the sphere with the help of analysis and
topology in non-rigorous fashion. By introducing methods from algebraic
number theory, Dedekind and Weber (see \cite{DeWe}) revealed the deep
analogy between number fields and function fields. They developed the
theory of algebraic functions in analogue with Dedekind's theory of algebraic
numbers, where the concept of ideal plays a central role. By introducing
such concepts into the theory of algebraic curves, Dedekind and Weber
paved the way for modern algebraic geometry. In the 1940s, the focus of
Weil's letter \cite{Weil} is just this analogue between number fields and
the field of algebraic functions of a complex variable. Weil described his
idea about this analogy using a third, intermediate subject, that of function
fields over a finite field, which he thought of as a bridge or Rosetta stone.
Hence, Weil's description of a Rosetta stone is a translation between aspects
of number theory, aspects of Riemann surface theory, and a third intermediate
parallel theory, that of algebraic curves over a finite field.

  For a given algebraic variety, there are two distinct aspects of its geometry:
one comes from its abelian version, i.e., motives, the other comes from its
utmost non-abelian version, i.e., anabelian algebraic geometry (fundamental
groups, defining ideals). A remarkable fact about the geometric realizations
of representations for $\text{PSL}(2, \mathbb{F}_p)$ in the first and second
column of the Rosetta stone (1.1) comes from different kinds of cohomology
theories. That is, for an algebraic variety, they give rise not to one, but
to many cohomology theories, among them the $\ell$-adic theories, one for
each prime $\ell$ different from the characteristic, and in characteristic
zero, the algebraic de Rham cohomology. These theories seem to tell the same
story, over and over again, each in a different language. This leads to the
philosophy of motives that there should exist a universal cohomology theory,
with values in a category of motives to be defined, from which all these
theories could be derived (see \cite{CoSe}). The theory of motives is a
generalization of the various cohomology theories of algebraic varieties.
In fact, the first and second column of the Rosetta stone (1.1) provide
the de Rham realization and the $\ell$-adic realization of such motives,
respectively. Classically, the relation between the Betti cohomology and
de Rham cohomology is expressed in terms of the periods and via Hodge theory.
Here, the first and second column of (1.1) give a comparison between de
Rham cohomology and $\ell$-adic \'{e}tale cohomology in the framework of
geometric representation theory for $\text{PSL}(2, \mathbb{F}_p)$ (or
$\text{SL}(2, \mathbb{F}_p)$). In particular, this gives the comparison
between $H_{\text{dR}}^1$ and $H_{\text{\'{e}t}}^1$. This leads us to
recall Grothendieck's ``mysterious functor'' which was to pass between
the Galois representation and the filtered Dieudonn\'{e} module by using
the theory of $p$-divisible groups as an intermediate (see \cite{Groth1970}
and \cite{Groth1974}). More precisely, given an abelian variety over a local
field, how to go from its Dieudonn\'{e} module to its \'{e}tale cohomology?
When $\Gamma$ is the $p$-divisible group of an abelian variety $A/\mathbb{Q}_p$
with good reduction (i.e., an abelian scheme $\mathbb{A}/\mathbb{Z}_p$), then
the Galois representation is the $p$-adic \'{e}tale cohomology group
$H_{\text{\'{e}t}}^1$ of $A/\overline{\mathbb{Q}}_p$, while the Dieudonn\'{e}
module is the de Rham cohomology group $H_{\text{dR}}^1$ of $\mathbb{A}/\mathbb{Z}_p$.
Therefore, Grothendieck's ``mysterious functor'' was to pass between the $p$-adic
\'{e}tale cohomology and the filtered de Rham cohomology. In particular, it
explains how to construct the Tate module of a $p$-divisible group over the
ring of integers of a $p$-adic field from its filtered Dieudonn\'{e} module
and vise versa. In fact, Tate and Grothendieck discovered that for an abelian
variety over a $p$-adic field, its de Rham cohomology and its $p$-adic
\'{e}tale cohomology both determine and are determined by its $p^{\infty}$-torsion
(more precisely, its $p$-divisible group). Hence, we have the following two
different kinds of comparisons between \'{e}tale cohomology and de Rham
cohomology:

(1) $\ell$-adic \'{e}tale cohomology group $H_{\text{\'{e}t}}^1$ of the
Drinfeld-Deligne-Lusztig curves $\text{DL}_p/\mathbb{F}_p$ $\longleftrightarrow$
de Rham cohomology group $H_{\text{dR}}^1$ of the Jacobian variety
$J(p)/\mathbb{C}$ of modular curves $X(p)$ by using the representation
theory of $\text{SL}(2, \mathbb{F}_p)$ (or $\text{PSL}(2, \mathbb{F}_p)$)
as an intermediate.

(2) $p$-adic \'{e}tale cohomology group $H_{\text{\'{e}t}}^1$ of an abelian
variety $A/\overline{\mathbb{Q}}_p$ $\longleftrightarrow$ de Rham cohomology
group $H_{\text{dR}}^1$ of an abelian scheme $\mathbb{A}/\mathbb{Z}_p$ by using
the theory of $p$-divisible groups as an intermediate.

  In particular, the theory of motives systematically generalizes the idea of
using the Jacobian variety of an algebraic curve $X$ as a replacement for the
cohomology group $H^1(X, \mathbb{Q})$ in the classical theory of correspondences.
In fact, the groups $H^1(X)$ are all deduced from the Jacobian variety
$J(X)=\text{Pic}^{0}(X)$ of $X$. This suggests identification of abelian
varieties with a class of motives, and to look upon $J(X)$ (which is an
abelian variety $A$) as being the motivic cohomology $H_{\text{mot}}^1(X)$
(see \cite{De1976(1)}). For each usual cohomology theory $H$, the dual
$H_1(X)$ of $H^1(X)$ is deduced from $J(X)$ by application of a realization
functor (see \cite{De1994}):
$$H^1(X)=\text{real $H^1_{\text{mot}}(X)$}.$$
Note that $H_1(X)$ is also the Tate twist $H^1(X)(1)$ of $H^1(X)$. For
$\ell$-adic cohomology, the realization of $A$: $\text{real}_{\ell}(A)$ is
the Tate module $V_{\ell}(A)=T_{\ell}(A) \otimes \mathbb{Q}_{\ell}$. For a
field $k$ of characteristic zero and the theory of de Rham, the realization
of $A$: $\text{real}_{\text{DR}}(A)$ is the Lie algebra of the universal
additive extension of $A$. Roughly speaking, $H^1$ are motives that are
abelian varieties. On the other hand, the third column of the Rosetta stone
(1.1) comes from anabelian algebraic geometry.  In particular, our theory of
defining ideals for modular curves $X(p)$ systematically realizes the idea
of using the defining ideal of an algebraic curve $X$ as a replacement for the
fundamental group $\pi_1(X)$ in the theory of modular curves $X(p)$. In particular,
the defining ideal provide a reification of the fundamental group $\pi_1$. Therefore,
the geometric realization of representations for $\text{PSL}(2, \mathbb{F}_p)$
involves both motives and anabelian algebraic geometry (defining ideals).
According to \cite{GrothRS}, the theme of anabelian algebraic geometry is
of comparable depth with the theme of motives. The significance of geometric
realizations of representations for $\text{PSL}(2, \mathbb{F}_p)$ is that it
provides a framework which unifies both motives and anabelian algebraic geometry
as its three kinds of distinct geometric realizations. At the same time, it
provides a comparison between motives and anabelian algebraic geometry as well
as a comparison between $\ell$-adic \'{e}tale cohomology and de Rham cohomology
as two distinct kinds of realizations of the motives. In particular, the comparison
between anabelian algebraic geometry and motives can be regarded as a counterpart
of Grothendieck's anabelian Tate conjecture and Tate conjecture for hyperbolic
curves and their Jacobian varieties, respectively.

  Note that the Langlands correspondence can also be formulated in three
different scenarios in the framework of Weil's Rosetta stone. That is, automorphic
forms over (i) number field; (ii) function field over a finite field; (iii)
function field over a compact Riemann surface. In the first two cases, the
objects are automorphic forms, which are, roughly speaking, functions on the
quotient of the form $G(F) \backslash G(\mathbb{A}_F)/K$. Here $F$ is a number
field or a function field, $G$ is a reductive algebraic group over $F$,
$\mathbb{A}_F$ is the ring of adeles of $F$, and $K$ is a compact subgroup
of $G(\mathbb{A}_F)$. However, in the third case, i.e., the geometric Langlands
program, is quite different from the Langlands program in its original formulation
for number fields and function fields. The most striking difference is that
it is formulated in terms of sheaves rather than functions. More precisely, in
the geometric theory the vector space of automorphic functions on the double quotient
$G(F) \backslash G(\mathbb{A}_F)/K$ is replaced by a derived category of sheaves
on an algebraic stack whose set of $\mathbb{C}$-points is this quotient.
$$\begin{matrix}
 &\text{Riemann surfaces} &\text{Curves over $\mathbb{F}_q$} &\text{Number fields}\\
 &\text{automorphic sheaves} &\text{automorphic forms} &\text{automorphic forms}\\
 &\text{(geometric theory)} &\text{over function fields} & \text{over number fields}
\end{matrix}\eqno{(2.1)}$$
Recently, Langlands proposed an analogous analytic form of the geometric theory
of automorphic forms, i.e., automorphic functions for complex algebraic curves,
in line with his program for global fields (see \cite{Langlands}, \cite{BK},
\cite{Frenkel}, \cite{EFK1}, \cite{EFK2}, \cite{EFK3} and \cite{EFK4}, as well
as \cite{GaWi} (see also \cite{KaWi}) in the gauge theory framework).

  From the viewpoint of representation theory and its geometric realization,
for the group $\text{SL}(2)$ and its variants $\text{GL}(2)$, there are six
theories that deserve to be studied in parallel. These are:

$\bullet$ The representation theory of $\text{SL}(2, \mathbb{F}_q)$
          and its geometric realization;

$\bullet$ The representation theory of $\text{SL}(2, \mathbb{C})$
          and its geometric realization;

$\bullet$ The representation theory of $\text{SL}(2, \mathbb{R})$
          and its geometric realization;

$\bullet$ The representation theory of $\text{SL}(2, \mathbb{Q}_p)$
          and its geometric realization;

$\bullet$ The automorphic representation theory of $\text{GL}(2, \mathbb{Q})$
          and its geometric realization;

$\bullet$ The automorphic representation theory of $\text{GL}(2, F(X))$
          and its geometric realization.

\noindent The second, third and fourth cases can be summarized as $\text{SL}(2, F)$,
where $F$ is a local field, i.e., $F=\mathbb{C}$, $\mathbb{R}$ or $\mathbb{Q}_p$. The
fifth and sixth can be summarized as $\text{SL}(2, F)$, where $F$ is a global field,
i.e., $F=\mathbb{Q}$ or a function field $F(X)$ of a curve $X$ over $\mathbb{F}_q$.
For each of the first four cases, there are the character theory and the geometric
realizations of their representations. In this way, representations of complex,
real, $p$-adic, adelic, and finite groups become a part of a more general scheme.

  Historically, the example of $\text{SL}(2, \mathbb{F}_q)$ played a seminal
role. This subject started in 1896 with the work of Frobenius on the irreducible
characters of $\text{PSL}(2, \mathbb{F}_p)$ (see \cite{Fr1}). In 1927, Weyl and
his student Peter published their proof of the Peter-Weyl theorem for compact
connected Lie groups and a generalization of the Schur orthogonality relations
for finite groups. Weyl was able to compute the irreducible characters of any
such Lie groups, i.e., Weyl's character formula. At the same time, using modular
curves $X(p)$ with $p$ a prime, Hecke constructed a finite dimensional (reducible)
representation $H^0(X(p), \Omega_{X(p)}^{1})$ of the finite group
$\text{PSL}(2, \mathbb{F}_p)$ whose dimension is equal to the genus of $X(p)$, i.e.
the number of linearly independent differentials of the first kind:
$$\dim H^0(X(p), \Omega_{X(p)}^{1})=g(X(p)).$$
This gives the first kind of geometric realization by de Rham cohomology of
complex representations for $\text{PSL}(2, \mathbb{F}_p)$. In particular,
it is also the first constructions of linear representations using cohomology.
Later on, the Borel-Weil theorem gives a geometric realization of each
irreducible representation of a compact connected semisimple Lie group, or
equivalently, of each irreducible holomorphic representation of a complex
connected semisimple Lie group. The realization is in the space of holomorphic
sections of a holomorphic line bundle over the flag variety of the group.
Because of the theorem of the highest weight, it produces a concrete geometric
realization for every irreducible holomorphic representation. Bott's theorem
directly implies a description of the cohomology groups of every holomorphic
line bundle on every complex projective manifold on which the group acts
transitively. In particular, the Borel-Weil-Bott theorem for the group
$\text{SL}(2, \mathbb{C})$ provides a geometric realization by holomorphic
sections of line bundles or cohomology on $\mathbb{P}^1(\mathbb{C})$. The
Atiyah-Bott fixed point formula gives a geometric interpretation of Weyl's
character formula. For a noncompact semisimple Lie group, the discrete series
appear. Harish-Chandra give a classification theorem of discrete series
representations in terms of their global characters. In order to obtain
geometric realizations of discrete series and the proof of Harish-Chandra's
theorem, Langlands observed that a purely formal application of the Atiyah-Bott
fixed point formula produced Harish-Chandra's formula for the discrete series
character. Combining this observation with the vanishing theorem for $L^2$-cohomology
in the very regular case, Langlands arrived at his conjecture on the geometric
realization of discrete series, which was proved by Schmid. In fact, in the
compact case, the Langlands conjecture reduces to the Borel-Weil-Bott Theorem
(see \cite{Schmid}). In particular, for $\text{SL}(2, \mathbb{R})$, whose geometric
realization is given by $L^2$-sections of line bundles or $L^2$-cohomology on
$\mathbb{P}^1(\mathbb{R})$. In 1974, Drinfeld (\cite{Dr1}) constructed a Langlands
correspondence for $\text{GL}(2, K)$, where $K$ is a global field of equal
characteristic. In the course of this work, he remarks that the cuspidal characters
of $G=\text{SL}(2, \mathbb{F}_q)$ may be found in the first $\ell$-adic cohomology
group (where $\ell$ is a prime number different from $p$) of the curve $X$
with equation $x y^q-y x^q=1$, on which $G$ acts naturally by linear changes
of coordinates (we call $X$ the Drinfeld curve). In fact, the polynomial $x y^q-y x^q$
is a determinant used in the construction of the Dickson invariant of the general
linear group and is an invariant of the special linear group. This example inspired
Deligne and Lusztig who then, in their fundamental article \cite{DeLu}, established
the basis of what has come to be known as Deligne-Lusztig theory. This gives the
second kind of geometric realization by $\ell$-adic cohomology of $\ell$-adic
representations for $\text{SL}(2, \mathbb{F}_q)$. For $\text{SL}(2, \mathbb{Q}_p)$,
whose representations involve local Langlands correspondence, its geometric realization
is given by sections of line bundles or $\ell$-adic cohomology on $\mathbb{P}^1(\mathbb{Q}_p)$.
In fact, Drinfeld \cite{Dr1976} introduced $p$-adic symmetric spaces embedded in
projective space. They parameterize certain types of formal groups and thus admit
coverings defined by division points of these. He conjectured that the $\ell$-adic
cohomology groups with compact supports gives a realization of all supercuspidal
representations of $\text{GL}(n, K)$ over the $p$-adic local field $K$ which would
give a construction of these representations analogous to the construction of the
discrete series representations of semisimple Lie groups through $L^2$-cohomology.
He also conjectured that the cohomology of these coverings realizes in some way the
Jacquet-Langlands correspondence between irreducible discrete series representation
of $\text{GL}(n, K)$ and irreducible representations of $D^{\times}$, the multiplicative
group of the division algebra $D$ of invariant $1/n$. This correspondence has been
established for $n>2$ by Deligne and Kazhdan (see \cite{BDKV}) and for $n=2$ by
Carayol (see \cite{Carayol1990}) and Faltings (see \cite{Faltings1994}).

\begin{center}
{\large\bf 2.1. Representations of $\text{PSL}(2, \mathbb{F}_p)$ and modular curves}
\end{center}

  The irreducible representations of $\text{SL}(2, \mathbb{F}_q)$ (where
we assume that $q$ is odd) are well known for a very long time (see
\cite{Fr1}, \cite{J} and \cite{S2}) and are a prototype example in
many introductory courses on the subject (see \cite{Bon}). For
$G=\text{SL}(2, \mathbb{F}_q)$, there are at most
$$1+1+\frac{q-3}{2}+2+\frac{q-1}{2}+2=q+4$$
conjugacy classes which can appear in the following formula on the
dimensions of ordinary representations:
$$|G|=1^2+q^2+\frac{q-3}{2} (q+1)^2+2 \cdot \left(\frac{q+1}{2}\right)^2
     +\frac{q-1}{2} (q-1)^2+2 \cdot \left(\frac{q-1}{2}\right)^2.$$
In particular, Schur (see \cite{S2}) derived the character tables of
the groups $\text{SL}(2, \mathbb{F}_q)$ for all values of $q$.

  In order to calculate the character table of $G$, algebraic methods
(in particular Harish-Chandra induction) give roughly half of the
irreducible characters (non-cuspidal characters with number $\frac{q+5}{2}$).
The other half (the cuspidal characters with number $\frac{q+3}{2}$)
can be obtained by the $\ell$-adic cohomology of the Drinfeld curve,
i.e., the Deligne-Lusztig induction (see \cite{Bon}).

  In general, let $G=G(\mathbb{F}_q)$ be the group of $\mathbb{F}_q$-rational
points of a reductive algebraic group defined over $\mathbb{F}_q$ where $q$
is a power of a prime number $p$. Let $\widehat{G}$ be the set of isomorphism
classes of irreducible representations of $G$. The Gel'fand-Harish-Chandra's
philosophy of cusp forms asserts that each member $\rho$ of $\widehat{G}$ can
be realized (see \cite{GH}), in some effective manner, inside the induction
$\text{Ind}_P^G(\sigma)$ from a cuspidal representation (i.e., one that is not
induced from a smaller parabolic subgroup) of a parabolic subgroup $P$ of $G$
$$\widehat{G}=\left\{\aligned
  &\rho < \text{Ind}_P^G(\sigma)\\
  &G \supset P-\text{parabolic}\\
  &\text{$\sigma$ cuspidal representation of $P$}
  \endaligned\right\}.$$
The above philosophy is a central one and leads to important developments in
representation theory of reductive groups over local and finite fields, in
particular the work of Deligne-Lusztig on the construction of representations
of finite reductive groups and Lusztig's classification of these representations.

  Let us see how this philosophy manifests itself in the cases when $G$ is
$\text{GL}_2(\mathbb{F}_q)$ of $\text{SL}_2(\mathbb{F}_q)$ (see \cite{GH}).
For simplicity, we ignore for the one dimensional representations. In both
cases around half of the irreducible representations are cuspidal, i.e.,
$P=G$, and form the so called discrete series, and half are the so called
principal series, i.e., induced from the Borel subgroup $P=B$ of upper
triangular matrices in $G$. Moreover, it is well known that the discrete
and principal series representations of $\text{SL}_2(\mathbb{F}_q)$ can be
obtained by restriction from the corresponding irreducible representations
of $\text{GL}_2(\mathbb{F}_q)$. These restriction are typically irreducible,
with only two exceptions: one discrete series representation of dimension
$q-1$ and one principal series representation of dimension $q+1$ split
into two pieces. Hence one has two discrete and two principal series
representations of dimensions $\frac{q-1}{2}$ and $\frac{q+1}{2}$,
respectively. These representations are called degenerate discrete series
and degenerate principal series. In general, For $G=\text{SL}_2(\mathbb{F}_q)$,
we have the trivial representation of dimension $1$, the Steinberg representation
of dimension $q$, induced representations $\text{Ind}_B^G \chi$ ($\chi$ a
character of the split torus) of dimension $q+1$ and discrete series (or
super-cuspidal) representations of dimension $q-1$. This discrete series
are associated to the characters of the non-split torus.
$$\begin{matrix}
 &G  & \text{discrete series ($P=G$)} \quad & \text{principal series ($P=B$)}\\
 & \text{GL}(2, \mathbb{F}_q) &\text{dim}(\rho)=q-1 \quad
 & \text{dim($\rho$)$=q+1$ or $q$}\\
 & \text{SL}(2, \mathbb{F}_q) &\text{dim}(\rho)=\left\{\aligned
                               &\text{$q-1$, or}\\
                               &\frac{q-1}{2}
                               \endaligned\right.
                              &\text{dim}(\rho)=\left\{\aligned
                               &\text{$q+1$ or $q$}\\
                               &\frac{q+1}{2}
                               \endaligned\right.
\end{matrix}$$

  Let $N$ be a positive integer and let $X(N)(\mathbb{C})=\Gamma(N) \backslash
\overline{\mathbb{H}}$ be the modular curve of level $N$ with respect to the
principal congruence subgroup $\Gamma(N)$. Here, $\overline{\mathbb{H}}=\mathbb{H}
\cup \{ \infty \} \cup \mathbb{Q}$ is the compactified upper half-plane and
$$\Gamma(N):=\text{ker}(\text{SL}(2, \mathbb{Z}) \rightarrow \text{SL}(2,
  \mathbb{Z}/N \mathbb{Z}))$$
is the principal congruence subgroup of level $N$.

  Recall that the space of holomorphic differentials $\Omega_{X(N)(\mathbb{C})}^1$
can naturally be identified with the space of weight two cusp forms on $\Gamma(N)$,
$$\Omega_{X(N)(\mathbb{C})}^1 \cong S_2(\Gamma(N)).$$
The natural map $X(N)(\mathbb{C}) \rightarrow X(1)(\mathbb{C})$ is a Galois cover
of curves with group $G=G_N=\text{SL}(2, \mathbb{Z}/N \mathbb{Z})/ \{ \pm 1 \}$
and thus $G$ acts on $X(N)(\mathbb{C})$ as a group of automorphisms. Thus, $G$
also acts on $\Omega_{X(N)(\mathbb{C})}^1$. The curve $X(N)(\mathbb{C})$ has a
canonical model $X(N)(\mathbb{Q})$ over $\mathbb{Q}$ with the property that
its space of holomorphic differentials $\Omega_{X(N)(\mathbb{Q})}^1$ can be
identified with the subspace of $S_2(\Gamma(N))$ consisting of forms with
rational Fourier coefficients
$$\Omega_{X(N)(\mathbb{Q})}^1 \cong S_2(\Gamma(N), \mathbb{Q}).$$
The curve  $X(N)(\mathbb{Q})$ is called Shimura's canonical model of $X(N)$.
The morphism $X(N)(\mathbb{C}) \rightarrow X(1)(\mathbb{C})$ descends to a
morphism $X(N)(\mathbb{Q}) \rightarrow X(1)(\mathbb{Q})$ on Shimura's
$\mathbb{Q}$-models but this cover is not Galois. Thus, $G$ does not act
on $X(N)(\mathbb{Q})$, but only on
$$X(N)(\mathbb{Q}(\zeta_N)):=X(N)(\mathbb{Q}) \otimes \mathbb{Q}(\zeta_N).$$

  Put
$$f_a(z)=\frac{g_2(\omega_1, \omega_2) g_3(\omega_1, \omega_2)}{\Delta(\omega_1, \omega_2)}
  \wp\left( a \left[\begin{matrix} \omega_1\\ \omega_2 \end{matrix}\right];
  \omega_1, \omega_2\right),$$
where $z=\omega_1/\omega_2 \in \mathbb{H}$, $a \in \mathbb{Q}^2$,
$a \notin \mathbb{Z}^2$.

\textbf{Proposition 2.1.1.} (see \cite{Sh}, p. 134) {\it For every positive integer
$N$, $\mathbb{C}(j, f_a | a \in N^{-1} \mathbb{Z}^2, \notin \mathbb{Z}^2)$ is
the field of all modular functions of level $N$.}

  Put
$$\mathfrak{F}_N=\mathbb{Q}(j, f_a | a \in N^{-1}
  \mathbb{Z}^2, \notin \mathbb{Z}^2).$$
The proposition 3.1 shows that $\mathbb{C} \cdot \mathfrak{F}_N$ is the field
of all modular functions of level $N$. We call an element of $\mathfrak{F}_N$
a modular function of level $N$ rational over $\mathbb{Q}(\zeta_N)$. The following
theorem will justify this definition.

\textbf{Theorem 2.1.2.} (see \cite{Sh}, pp. 137-138) {\it The field $\mathfrak{F}_N$
has the following properties.

$(1)$ $\mathfrak{F}_N$ is a Galois extension of $\mathbb{Q}(j)$.

$(2)$ For every $\beta \in \text{GL}(2, \mathbb{Z}/N \mathbb{Z})$,
$f_a \mapsto f_{a \beta}$ gives an element of the Galois group
$\text{Gal}(\mathfrak{F}_N/ \mathbb{Q}(j))$, which we write $\tau(\beta)$.
Then $\beta \mapsto \tau(\beta)$ gives an isomorphism of
$\text{GL}(2, \mathbb{Z}/N \mathbb{Z})/ \{ \pm 1 \}$
to $\text{Gal}(\mathfrak{F}_N/\mathbb{Q}(j))$.

$(3)$ If $\zeta_N$ is a primitive $N$-th root of unity, then
$\zeta \in \mathfrak{F}_N$, and $\zeta^{\tau(\beta)}=\zeta^{\det(\beta)}$.

$(4)$ $\mathbb{Q}(\zeta_N)$ is algebraically closed in $\mathfrak{F}_N$.

$(5)$ $\mathfrak{F}_N$ contains the functions $j \circ \alpha$ for all
$\alpha \in M_2(\mathbb{Z})$ such that $\det(\alpha)=N$.}

\textbf{Proposition 2.1.3.} (see \cite{Sh}, p. 140) {\it $(1)$ $\mathfrak{F}_N$
coincides with the field of all the modular functions of level $N$ whose
Fourier expansions with respect to $e^{2 \pi i z/N}$ have coefficients in
$\mathbb{Q}(\zeta_N)$.

$(2)$ The field $\mathbb{Q}(j(z), j(Nz), f_{a_1}(z))$, with $a_1=(N^{-1}, 0)$,
coincides with the field of all the modular functions of level $N$ whose
Fourier expansions with respect to $e^{2 \pi i z/N}$ have rational coefficients.

$(3)$ The field of $(2)$ corresponds to the subgroup
$$\left\{ \left(\begin{matrix} \pm 1 & 0\\ 0 & x \end{matrix}\right): x \in
  (\mathbb{Z}/N \mathbb{Z})^{\times} \right\}/ \{ \pm 1 \}$$
of $\text{GL}(2, \mathbb{Z}/N \mathbb{Z})/ \{ \pm 1 \}$.}

\begin{center}
{\large\bf 2.2. Analytic theory and de Rham realization of representations
                for $\text{PSL}(2, \mathbb{F}_p)$: Riemann surfaces}
\end{center}

  In order to establish an algebraic theory of algebraic functions,
in a series of remarkable papers \cite{Hec1}, \cite{Hec2}, \cite{Hec3},
\cite{Hec4} and \cite{Hec5}, Hecke studied the arithmetic theory of
elliptic modular functions by the methods of function theory
(Riemann-Weierstrass), the methods of arithmetic (Dedekind-Weber),
Galois theory and Frobenius's character theory.

  In the introduction of his work \cite{Hec5}, Hecke stated that
``If an algebraic structure in a variable of genus $p$ has a group $\mathfrak{G}$
of one-to-one analytical mappings (automorphisms), then the $p$ linearly
independent differentials of the first kind experience linear homogeneous
substitutions in these mappings, which together form a representation
$\mathfrak{S}$ of the abstract group $\mathfrak{G}$. A general method of
determining this representation, i.e. determining its character, was first
given by Hurwitz. For the special structure, which is defined by the elliptic
modular functions of a prime number level $q$, I determined the decomposition
of the above representation into its irreducible components a few years ago
in a somewhat different way, whereby a strange connection with the class number
of the quadratic number field $K(\sqrt{-q})$ yielded. Since then, the question
has been taken up again by C. Chevalley and A. Weil, and a very general theorem
about these and related groups has been proved by algebraic methods.''

  In the other work \cite{Hec3}, page 551, Hecke stated that ``We now consider
the algebraic structure, which is defined by the elliptic modular functions of
the prime level $q$ ($q \geq 3$). The structure has a group of mappings into itself,
which are explained by the totality of the modular substitutions applied to the
argument $\tau$; the mappings are thus reduced to the system of modular substitutions,
mod $q$, and isomorphic with the modular group $\mathfrak{M}(q)$ of order
$\frac{q(q^2-1)}{2}$. The group of integrals of the first kind induced in
this way is decomposed into its irreducible components
$$\mathfrak{S}=x \mathfrak{G}_q+y_1 \mathfrak{G}_{\frac{q+\varepsilon}{2}}+
  y_2 \mathfrak{G}_{\frac{q+\varepsilon}{2}}^{\prime}+\sum_{i} u_i
  \mathfrak{G}_{q+1}^{(i)}+\sum_{i} v_i \mathfrak{G}_{q-1}^{(i)}.\eqno{(2.2.1)}$$
Here the $\mathfrak{G}_n$, $\mathfrak{G}_n^{(i)}$ are the various possible
irreducible representations of $\mathfrak{M}(q)$, $n$ being the degree,
$\varepsilon=(-1)^{\frac{q-1}{2}}$; the multiplicities $x$, $u_i$, $v_i$
and also $y_1$, $y_2$ have already been calculated by me.''

  In modern terminology, the group
$\text{PSL}(2, \mathbb{Z}/N \mathbb{Z})=\text{SL}(2, \mathbb{Z}/N \mathbb{Z})/\{\pm I\}$
acts on the holomorphic cusp forms of weight $2$ on $\mathbb{H}$ with
respect to the principal congruence group $\Gamma(N)$ of level $N$:
$$\pi(\gamma) f=f | [\gamma^{-1}]_2$$
where (using Shimura's notation)
$$(f | [g]_k)(z)=f(g(z)) (cz+d)^{-k} \quad \text{if} \quad
  g=\left(\begin{matrix} a & b\\ c & d \end{matrix}\right).$$
The group $\Gamma(N)$ itself acts trivially by definition, and so does $\pm I$,
so the action passes to the quotient $\text{SL}(2, \mathbb{Z})/\Gamma(N) \cdot
\{ \pm I \}$, which is canonically isomorphic to $\text{PSL}(2, \mathbb{Z}/N \mathbb{Z})$.
A natural question apparently first investigated by Hecke is: what is the
irreducible decomposition of this representation? At the time he asked this,
the classification of the irreducible representations of $\text{PSL}(2, \mathbb{Z}/N)$
had been worked out by Frobenius in the case $N=p$, a prime, and this was the
case Hecke looked at.

  It is know that the forms of weight two with respect to $\Gamma(p)$ may be
identified with the space of holomorphic differential forms on $X(p)$. Let $\pi$
be the representation of $\text{PSL}(2, \mathbb{F}_p)$ on this space. It was
known to Hecke that the representation $\pi \oplus \overline{\pi}$ is that of
$\text{PSL}(2, \mathbb{F}_p)$ on the rational cohomology of $X(p)$. Therefore
the sum of $\pi$ and its complex conjugate has to be rational. But the
representation $\pi$ itself is not necessarily even real, and so it makes sense
to ask, what can one say about the difference between $\pi$ and $\overline{\pi}$?
This is the question that most intrigued Hecke (see \cite{Casselman}).

  Using modular curves $X(p)$ with $p$ a prime, Hecke constructed a finite
dimensional (reducible) representation $H^0(X(p), \Omega_{X(p)}^{1})$ of
the finite group $\text{PSL}(2, \mathbb{F}_p)$ whose dimension is equal
to the genus of $X(p)$, i.e. the number of linearly independent differentials
of the first kind:
$$\dim H^0(X(p), \Omega_{X(p)}^{1})=g(X(p)).$$
This construction has the rich structure: the finite dimensional representation
of $\text{PSL}(2, \mathbb{F}_p)$ can be obtained by the methods from geometry,
topology and complex analysis, i.e. Riemann-Weierstrass function theory. In
particular, with the help of the character theory established by Frobenius
(see \cite{Fr1}), this leads Hecke to study the decomposition of the vector
space of regular differentials on $X(p)$: $H^0(X(p), \Omega_{X(p)}^{1})$ under
the action of $\text{PSL}(2, \mathbb{F}_p)$ as the direct sum of spaces of
irreducible representations for $\text{PSL}(2, \mathbb{F}_p)$ defined over
some number fields (see \cite{Hec1}, \cite{Hec2}, \cite{Hec3}, \cite{Hec4}, and
\cite{Hec5}). In particular, he proved that (see \cite{Hec5}, p.761, Satz 14)
that every irreducible representation of $\text{PSL}(2, \mathbb{F}_p)$ is
equivalent to one whose coefficients lie in the field of its character
$\mathbb{Q}(\chi)$.

  According to \cite{Se2}, Chap. 12, let $K$ denote a field of characteristic
zero and $C$ an algebraic closure of $K$. If $V$ is a $K$-vector space, we let
$V_C$ denote the $C$-vector space $C \otimes_K V$ obtained from $V$ by extending
scalars from $K$ to $C$. If $G$ is a finite group, each linear representation
$\rho: G \rightarrow \text{GL}(V)$ over the field $K$ defines a representation
$$\rho_C: G \rightarrow \text{GL}(V) \rightarrow \text{GL}(V_C)$$
over the field $C$. In terms of modules, we have
$$V_C=C[G] \otimes_{K[G]} V.$$
The character $\chi_{\rho}=\text{Tr}(\rho)$ of $\rho$ is the same as for $\rho_C$.
It is a class function on $G$ with values in $K$. We denote by $R_K(G)$ the group
generated by the characters of the representations of $G$ over $K$. It is a sub-ring
of the ring $R(G)=R_C(G)$.

  A linear representation of a finite group $G$ over $C$ is said to be realizable
over $K$ (or rational over $K$) if it is isomorphic to a representation of the form
$\rho_C$, where $\rho$ is a linear representation of $G$ over $K$; this amounts to
saying that it can be realized by matrices having coefficients in $K$.

  Now, let us consider the realizability over cyclotomic fields. Denote by $m$ the
least common multiple of the orders of the elements of $G$. It is a divisor of $g$,
the order of $G$. We have Brauer's theorem (see \cite{Se2}, 12.3, Theorem 24):
If $K$ contains the $m$-th roots of unity, then $R_K(G)=R(G)$. Consequently,
each linear representation of $G$ can be realized over $K$.

  In his arithmetic study of complex representations of finite groups, Schur
(see \cite{S1}) introduced the notion of the index. Roughly speaking, suppose
$G$ is a finite group, $K$ is a splitting field for $G$, and $\chi$ is the
character of an irreducible linear representation $\rho$ of $G$ over $K$.
Suppose $k$ is the subfield of $K$ generated by the character values $\chi(g)$,
$g \in G$. The Schur index of $\chi$ or the Schur index of $\rho$ is defined as
the smallest positive integer $m$ such that there exists a degree $m$ extension
$L$ of $k$ such that $\rho$ can be realized over $L$, i.e., we can change basis
so that all the matrix entries are from $L$. The Schur index of a character
$\chi$ is denoted $m(\chi)$ (see also \cite{Se2} for the definition using the
theory of semi-simple algebras). Many important results about Schur indices appear
to depend on deep facts about division algebras and number theory (see \cite{Fei}
and \cite{I}, Chap. 10). Note that if the representation can be realized over the field
generated by the character values for that representation, the Schur index is one.
In particular, all irreducible complex characters of $\text{PSL}(2, \mathbb{F}_q)$
have Schur index $1$ (see \cite{Jan}, also \cite{S2}, \cite{Dor}, \S 38 and
\cite{Se2}, \S 12.6).

  In fact, let $G$ be a finite group and $\rho: G \rightarrow \text{GL}(n, \mathbb{C})$
be an irreducible complex representation. Then $\rho$ descends to a representation
over $\overline{\mathbb{Q}}$ uniquely up to isomorphism. Hence, we can consider
$\rho$ as a representation with values in $\text{GL}(n, \overline{\mathbb{Q}})$.
Under $\overline{\mathbb{Q}}$, any field over which $\rho$ is definable must
contain the trace field of $\rho$, defined as the extension of $\mathbb{Q}$
generated by the traces of all the matrices $\rho(g)$, $g \in G$, that is,
$\mathbb{Q}(\text{tr}(\rho(g)), \rho \in G)$, which is just the character field
$\mathbb{Q}(\chi)$. The trace field of $\rho$ is characterized as the fixed field
of Galois automorphisms $\sigma \in \text{Gal}(\overline{\mathbb{Q}}/\mathbb{Q})$
where $\rho \cong \rho^{\sigma}$, the Galois twist of $\rho$ by $\sigma$. Thus,
we can denote the trace field of $\rho$ as $\mathbb{Q}(\rho)$. The obstruction to
the representation descending to its trace field can be described by a Galois
$2$-cocycle, which represents an element of the Brauer group of the trace field,
i.e., the Brauer obstruction of $\rho$. More precisely, the Brauer obstruction
$[\psi_{\rho}] \in H^2(\text{Gal}(\overline{\mathbb{Q}}/\mathbb{Q}(\rho)),
\overline{\mathbb{Q}}^{\times})=\text{Br}(\mathbb{Q}(\rho))$ is defined as follows:
for $\sigma \in \text{Gal}(\overline{\mathbb{Q}}/\mathbb{Q}(\rho))$, the isomorphism
$\rho \cong \rho^{\sigma}$ is unique up to $\overline{\mathbb{Q}}^{\times}$ by
Schur lemma. When we choose isomorphisms $\rho \cong \rho^{\sigma}$, we obtain
the Brauer obstruction in $H^2(\text{Gal}(\overline{\mathbb{Q}}/\mathbb{Q}(\rho)),
\overline{\mathbb{Q}}^{\times})$, which is well-defined as the cohomology class is independent
of the choice of isomorphisms.

  Despite of this, we still need the realization in the field $\mathbb{Q}(\chi)$
and the exact knowledge of it. Thanks to Deligne (see \cite{De2} and \cite{De4},
also \cite{De1} and \cite{De3})), who gives an elegant and different proof using
the Kirillov model in representation theory.

\textbf{Theorem 2.2.1. (Deligne's realization theorem).} (see \cite{De4}). {\it Any
representation of $\text{GL}(2, \mathbb{F}_q)$ can be realized over the field of
its character.}

  In fact, Deligne proved the more general result: The split (resp. non-split)
torus $T_0$ of $\text{SL}(2, \mathbb{F}_q)$ are, up to $x \mapsto x^{-1}$,
naturally isomorphic to $\mathbb{F}_q^{\times}$ (diagonal matrices $\left(\begin{matrix}
x & 0\\ 0 & x^{-1} \end{matrix}\right)$), resp. Ker(N: $\mathbb{F}_{q^2}^{\times}
\rightarrow \mathbb{F}_q^{\times})$. Let $\chi_0$ be a nontrivial character
of $T_0$, $t_0$ be its order, and $X_{\chi_0}$ be the corresponding
representation of $\text{SL}(2, \mathbb{F}_q)$, well defined up to
isomorphism. For $\chi_0$ of order $2$, we mean the sum of the two
representations of dimension $\frac{q+1}{2}$ (resp. $\frac{q-1}{2}$).
Assume that $\chi_0(-1)=1$. This is the case for $q$ even. For $q$ odd,
it is equivalent to $t_0| \frac{q-1}{2}$ (resp. $t_0 | \frac{q+1}{2}$).
It is also equivalent to $X_{\chi_0}$ factoring through a representation
of $\text{PSL}(2, \mathbb{F}_q)$.

  Under this assumption, Hecke states that $X_{\chi_0}$ can be realized
over the field (generated by values) of his character, which is the real
subfield $F_0$ of the $t_0$-cyclotomic field.

   Let $H$ be the following quaternion algebra over $F_0$: it is generated
by the quadratic extension $F$ of $F_0$, and by an element $u$, such that
conjugation by $u$ induces on $F$ the generator of $\text{Gal}(F/F_0)$, and
that
$$\left\{\aligned
 &\text{in the split case (principal series)} \quad : u^2=q,\\
 &\text{in the non-split case (discrete series)} : u^2=-q.
\endaligned\right.$$
This algebra is a matrix algebra if and only if $q$ (resp. $-q$) is
in the image of the norm $\text{N}_{F/F_0}$.

  In particular, Deligne proved that the obstruction to realize the
representation of $\text{PSL}(2, \mathbb{F}_q)$ corresponding to $\chi_0$
over the field of its character, that is over $F_0$, is the class of the
quaternion algebra $H$ in the Brauer group of $F_0$.

  Hence, in order to prove (at least for $t_0 \neq 2$) Hecke's statement
that representations of $\text{PSL}(2, \mathbb{F}_q)$ are realizable over
the field of their characters, it remains to prove that the obstruction made
explicit as above vanishes. In particular, for $X$ a discrete series representation
of $\text{GL}(2, \mathbb{F}_q)$, it essentially leads to the Kirillov model of
the representation $X$.

  Recall that the main result of \cite{Hec5} is expressed in the following
simple formula: The multiplicity $r$ of an irreducible representation $\rho$ of
$\text{PSL}(2, \mathbb{F}_p)$ with the character $\chi$ within the group
$H^0(X(p), \Omega_{X(p)}^{1})+\overline{H^0(X(p), \Omega_{X(p)}^{1})}$
(where $\rho$ is not the identical representation) is
$$\kappa(\rho)+\kappa(\overline{\rho})=r=
 f-\frac{1}{p} \sum_{\text{$n$ mod $p$}} \chi(T^n)-\frac{1}{2} \sum_{\text{$n$ mod $2$}}
  \chi(S^n)-\frac{1}{3} \sum_{\text{$n$ mod $3$}} \chi((ST)^n).\eqno{(2.2.2)}$$
Here $T$, $S$ are the two modular substitutions $\tau^{\prime}=\tau+1$,
$\tau^{\prime}=-\frac{1}{\tau}$, $f$ is the degree of $\rho$. Therefore,
Hecke obtained the following beautiful decomposition formula:
$$H^0(X(p), \Omega_{X(p)}^{1})=m V_p \oplus n_1 V_{\frac{p+\varepsilon}{2}}
  \oplus n_2 V_{\frac{p+\varepsilon}{2}}^{\prime} \oplus \bigoplus_i u_i
  V_{p+1}^{(i)} \oplus \bigoplus_j v_j V_{p-1}^{(j)},\eqno{(2.2.3)}$$
where $\varepsilon=(-1)^{\frac{p-1}{2}}$, $i=1, 2, \ldots, \frac{q-\varepsilon}{4}-1$
and $j=1, 2, \ldots, \frac{q+\varepsilon}{4}-\frac{1}{2}$. The multiplicities
$m$, $n_1$, $n_2$, $u_i$ and $v_j$ can be calculated by the above formula on $r$.
Note that for small primes $p$, the irreducible representations appearing in the
decomposition formula  of $H^0(X(p), \Omega_{X(p)}^{1})$ can not exhaust all
of the irreducible representations of $\text{PSL}(2, \mathbb{F}_p)$, especially
for $p=7$, $11$ and $13$.

  Hecke's really interesting observation was that the difference between the
degenerate discrete series representations $\pi_{\pm}$ of dimension $\frac{p-1}{2}$
has global arithmetic significance for every $p>3$, $p \equiv 3$ (mod $4$).
They do not occur with the same multiplicity in the representation of
$\text{SL}(2, \mathbb{F}_p)$ on differential forms on $X(p)$. Instead, the
difference in multiplicities is accounted for by the cusp forms contributed
by Gr\"{o}ssen-characters associated (through binary $\theta$-functions) to
$\mathbb{Q}(\sqrt{-p})$. These had been constructed in \cite{Hec1}. Hecke
proved this by an application of Riemann-Roch theorem, but it can be more cleanly
proved by applying the Atiyah-Bott trace formula (which was in fact motivated by
an earlier result of Eichler about algebraic curves, brought to the attention of
Atiyah and Bott by Shimura at the Woods Hole conference of 1965). Let $\pi$ be
the representation of $G=\text{PSL}(2, \mathbb{F}_p)$. Let $h(-p)$ be the ideal
class number of the quadratic imaginary quadratic extension $\mathbb{Q}(\sqrt{-p})$.
Suppose that as a representation of $G$
$$\pi=\sum m_{\pi_k} \pi_k,$$
where $m_{\pi_k}$ denotes the multiplicity of $\pi_k$.

\textbf{Theorem 2.2.2.} (Hecke). {\it The difference $m_{\pi_{+}}-m_{\pi_{-}}$
is equal to $h(-p)$.}

  Which of the two representations $\pi_{\pm}$ occurs more often than the
other is implicitly related, as Hecke himself observed, to the sign of
Gauss sums. This sort of problem, i.e. relation with class field theory,
is ubiquitous in this business. Hecke did not tie all these facts together
systematically, but this was done by Labesse and Langlands (see \cite{LL}
and \cite{Casselman}) around 1971.  One of the main results of that paper
is that subspaces of automorphic forms on the ad\`{e}lic $\text{SL}(2, \mathbb{Q})
\backslash \text{SL}(2, \mathbb{A})$ arising from quadratic extensions of
$\mathbb{Q}$ via theta functions can be characterized in terms of local
invariants and quadratic reciprocity. This leads to the concept of endoscopy,
which decomposes conjugate-invariant distributions on reductive groups over
local fields as well as on arithmetic quotients into what Langlands
calls stable components. In particular, the Fundamental Lemma is part
of a general program on the endoscopy. It is a statement about
conjugate-invariant functionals on the Hecke algebra associated to
reductive groups defined over $p$-adic fields. It is fundamental because
it will undoubtedly play an important role in any final proof of Langlands
functoriality conjectures. As the paper \cite{LL} already briefly mentions,
certain simple cases of endoscopy were known already to Hecke, in about 1930.

  In fact, for every $p>3$, $p \equiv 3$ (mod $4$), Hecke (see \cite{Hec5},
p. 768, Satz 16) described a $\text{SL}(2, \mathbb{F}_p)$-invariant subspace
of dimension $h \cdot \frac{p-1}{2}$ in $H^0(X(p), \Omega_{X(p)}^{1})$, where $h$
is the class number of the imaginary quadratic field $K=\mathbb{Q}(\sqrt{-p})$,
the group $\text{SL}(2, \mathbb{F}_p)$ acts on this subspace by $h$ copies of $W$,
one of the two irreducible degenerate discrete series representations of
$\text{SL}(2, \mathbb{F}_p)$ of dimension $\frac{p-1}{2}$ (see also \cite{Sh1},
\cite{Sh2}, \cite{Gro87} and \cite{Gro2} for the later development). This subspace
is spanned by weighted binary theta series, and Hecke identified certain periods
of these differentials as the periods of elliptic curves with complex multiplication
by $K$. Namely, Hecke showed that (see \cite{Hec5}, p.768, Satz 16) among the
$h \cdot \frac{p-1}{2}$ integrals of the first kind of this family, the independent
ones can be selected in such a way that their periods at $\Gamma(p)$ are always
integers from $\mathbb{Q}(\sqrt{-p})$. These elliptic curves appear as factors
of the Jacobian of the modular curve $X(p)$.

  In particular, for the smallest such primes $p=7$ and $p=11$ with class number
one, Hecke gave the following remarkable decompositions:
$$H^0(X(7), \Omega_{X(7)}^1)=V_3,\eqno{(2.2.4)}$$
$$H^0(X(11), \Omega_{X(11)}^1)=V_5 \oplus V_{11} \oplus V_{10}.\eqno{(2.2.5)}$$

  Here, when $p=7$, the curve $X(7)$ has genus three and is isomorphic to the
Klein quartic curve. In this case, Hecke showed that the multiplicity $m(W)$
for the irreducible representations $W$ of $\text{PSL}(2, \mathbb{F}_7)$ is
equal to one, for $W$ is one of the degenerate discrete series representations
of dimension $3$, and equal to zero for all other irreducible representations of
$\text{PSL}(2, \mathbb{F}_7)$. Moreover, the integrals of $V_3$ only have periods
that are integers from $\mathbb{Q}(\sqrt{-7})$. The degenerate discrete series
representation of $\text{PSL}(2, \mathbb{F}_7)$ is realizable over
$\mathbb{Q}(\sqrt{-7})$. The Jacobian variety $J(7)(\mathbb{Q}(\sqrt{-7}))$
of $X(7)$ over $\mathbb{Q}(\sqrt{-7})$ is isogenous to the product $E^3$
corresponding to the degenerate discrete series. Here the elliptic curve
$E$ is just the modular curve $X_0(49)$ with conductor $49$, of equation
$$y^2+xy=x^3-x^2-2x-1, \quad \text{with} \quad j(E)=-3^3 \cdot 5^3.$$

  When $p=11$ the curve $X(11)$ has genus $26$. In this case, Hecke showed
that the multiplicity $m(W)$ is equal to one for $W$ is the degenerate discrete
series representation of dimension $5$, the Steinberg representation of
dimension $11$, and the discrete series representation of dimension $10$,
and equal to zero for all other irreducible representations of
$\text{PSL}(2, 11)$. Hecke also pointed out that (see \cite{Hec5}, p.769,
Satz 18) for the modular curve $X(11)$ of level $11$ and genus $g=26$, the
$26$ integrals of the first kind can be selected as elliptic integrals.
Moreover, the integrals of $V_5$ only have periods that are integers from
$\mathbb{Q}(\sqrt{-11})$, and $V_{11}$, $V_{10}$ have a rational character
and the multiplicity $1$. This implies that $V_{10}$ corresponds to the
character $\chi_4$ or $\chi_5$.

  Here, the Steinberg representation and discrete series representation
of $\text{PSL}(2, \mathbb{F}_{11})$ appear in Hecke's decomposition are
realizable over $\mathbb{Q}$. While, the degenerate discrete series
representation of $\text{PSL}(2, \mathbb{F}_{11})$ is realizable over
$K=\mathbb{Q}(\sqrt{-11})$. Let $X_{\text{split}}(11)$, resp.
$X_{\text{nonsplit}}(11)$, resp. $X_{\mathfrak{S}_4}(11)$ the algebraic
curves associated with congruence groups whose image in
$\text{PSL}(2, \mathbb{F}_{11})$ is the normalizer of a split Cartan
group, resp. nonsplit, resp. from a group isomorphic to group
$\mathfrak{A}_4$, an alternate group of four elements (see \cite{Ligozat},
\cite{Se1972} and \cite{Mazur1977}). There is a canonical way to endow
these curves with a curve structure on $\mathbb{Q}$.

  Let us give some background (see \cite{Ligozat}). Let $G_0$ be the
split Cartan group of $\text{GL}(2, \mathbb{F}_{11})$ formed by the
diagonal matrices. Let $G_{\text{split}}$ be the normaliser of $G_0$
in $\text{GL}(2, \mathbb{F}_{11})$. Let $G_{\text{nonsplit}}$ be the
normaliser of the nonsplit Cartan group formed by the matrices of
$\text{GL}(2, \mathbb{F}_{11})$ of the form
$$\left(\begin{matrix}
  a & b\\
 -b & a
 \end{matrix}\right), a, b \in \mathbb{F}_{11}, a^2+b^2 \neq 0.$$
Let $H_{\mathfrak{S}_4}$ be the subgroup of $\text{SL}(2, \mathbb{F}_{11})$
generated by the matrices
$$h^{\prime}=\left(\begin{matrix}
              0 & -3\\
              4 & 0
             \end{matrix}\right), \quad
  h^{\prime \prime}=\left(\begin{matrix}
              5 & 5\\
             -4 & 5
             \end{matrix}\right), \quad
  h^{\prime \prime \prime}=\left(\begin{matrix}
              2 & -4\\
             -2 & -1
             \end{matrix}\right).$$
We have the relations:
$$h^{\prime 2}=-1, \quad h^{\prime \prime 3}=1, \quad
  h^{\prime \prime \prime 3}=-1, \quad
  h^{\prime} h^{\prime \prime} h^{\prime \prime \prime}=-1.$$
It follows that $H_{\mathfrak{S}_4}$ contains $-1$ , and that the
image $H_{\mathfrak{S}_4}/\{ \pm 1 \}$ of $H_{\mathfrak{S}_4}$ in
$\text{PSL}(2, \mathbb{F}_{11})$ is isomorphic to the alternating
group $\mathfrak{A}_4$. The matrix $M=\left(\begin{matrix} 3 & 0\\
0 & -3 \end{matrix}\right)$ normalises $H_{\mathfrak{S}_4}$. Let
$G_{\mathfrak{S}_4}$ be the subgroup of $\text{GL}(2, \mathbb{F}_{11})$
generated by $H_{\mathfrak{S}_4}$ and $M$. Since $-M^2$ generates
the homothety group, $G_{\mathfrak{S}_4}$ contains the homotheties.
On the other hand, $\det(M)=2$ generates $\mathbb{F}_{11}^{\times}$;
this means that the image of $G_{\mathfrak{S}_4}$ in
$\text{PSL}(2, \mathbb{F}_{11})$ contains $\mathfrak{A}_4$ as a normal
subgroup of index $2$, and is therefore isomorphic to $\mathfrak{S}_4$.
Let $H_0$, resp. $H_{\text{split}}$, resp. $H_{\text{nonsplit}}$,
be the intersection of $G_0$, resp. $G_{\text{split}}$, resp.
$G_{\text{nonsplit}}$ with $\text{SL}(2, \mathbb{F}_{11})$.
Let $X_{\text{split}}(11)$, $X_{\text{nonsplit}}(11)$,
$X_{\mathfrak{S}_4}(11)$ resp. be the canonical models associated
with the groups $G_{\text{split}}$, $G_{\text{nonsplit}}$,
$G_{\mathfrak{S}_4}$. These are curves over $\mathbb{Q}$.
Note that each of the three groups $G_{\text{split}}$, $G_{\text{nonsplit}}$,
$G_{\mathfrak{S}_4}$ contains the matrix $M$. The associated curves
therefore appear as quotients over $\mathbb{Q}$ of the $M$-model $X_M(11)$
defined by $M$. Let $J_{\text{split}}(11)$, $J_{\text{nonsplit}}(11)$ and
$J_{\mathfrak{S}_4}(11)$ denote the Jacobian varieties corresponding to
$X_{\text{split}}(11)$, $X_{\text{nonsplit}}(11)$, and $X_{\mathfrak{S}_4}(11)$,
respectively. The Jacobian $J(11)$ of a suitable model over $\mathbb{Q}$
of $X(11)$ is isogenous over $K$ to the product of of $26$ elliptic curves.
The decomposition of $J(11)(K)$ into isogeny classes over $K$ contains $d$
copies of $J_H$, where $H$ is one of the following three cases according
to Hecke's decomposition for $p=11$:

(1) The group $H$ is the Borel group of $\text{SL}(2, \mathbb{F}_{11})$
    formed by the upper triangular matrices, which corresponds to the
    Steinberg representation.

(2) The group $H$ is $H_{\text{nonsplit}}$, which corresponds to the
    degenerate discrete series.

(3) The group $H$ is $H_{\mathfrak{S}_4}$, which corresponds the the
    discrete series associated with the character $\chi_4$.

\noindent Correspondingly, the decomposition of $J(11)(K)$ into isogeny
classes over $K$ contains:

(1) $11$ copies of $X_0(11)$;

(2) $5$ copies of $J_{\text{nonsplit}}(11)$;

(3) $10$ copies of $J_{\mathfrak{S}_4}(11)$.

\noindent Correspondingly,

(1) $J_{\text{split}}(11)$ is isogenous over $K$ to the product $X_0(11) \times E_6$;

(2) $J_{\text{nonsplit}}(11)$ is isogenous over $K$ to $E_6$;

(3) $J_{\mathfrak{S}_4}(11)$ is isogenous over $K$ to $E_4$.

  In fact, it is shown that $J(11)(K)$ is isogenous to the product
$$E_1^{11} \times E_6^5 \times E_4^{10},$$
and that $J_{\text{split}}(11)$, $J_{\text{nonsplit}}(11)$ and
$J_{\mathfrak{S}_4}(11)$ are isogenous over $K$ to $E_1 \times E_6$,
$E_6$ and $E_4$ respectively. Here, $X_0(11)$ is isomorphic to $E_1$
over $\mathbb{Q}$ with conductor $11$, of equation:
$$y^2+y=x^3-x^2-10 x-20, \quad \text{with} \quad
  j(E_1)=-\frac{2^{12} \cdot 31^3}{11^5}.$$
The curve $X_{\text{nonsplit}}(11)$ is isomorphic over $\mathbb{Q}$
to the curve $E_6$, of equation:
$$y^2+y=x^3-x^2-7x+10, \quad \text{with} \quad j(E_6)=-2^{15}.$$
The curve $X_{\mathfrak{S}_4}(11)$ is isomorphic over $\mathbb{Q}$ to the
curve $E_4^{\prime}=E_5/C_5$, of equation:
$$y^2+xy+y=x^3+x^2-305 x+7888, \quad \text{with} \quad j(E_4^{\prime})=-11^2,$$
where $E_4$ is given by the equation:
$$y^2+xy=x^3+x^2-2x-7, \quad \text{with} \quad j(E_4)=-11^2,$$
and $E_5$ is given by the equation:
$$y^2+xy+y=x^3+x^2-30 x-76, \quad \text{with} \quad j(E_5)=-11 \cdot 131^3.$$
Note that $E_4$, $E_5$ and $E_6$ are three elliptic curves of conductor
$121$. In the notation of \cite{Cremona}, they are C1(F), A1(H) and B1(D),
respectively. The elliptic curve $E_4^{\prime}$ with conductor $121$ is
A2(I) in the notation of \cite{Cremona}.

  However, for every $p>5$, $p \equiv 1$ (mod $4$), it becomes completely
different. All of the conjugacy classes in $\text{SL}(2, \mathbb{F}_p)$ are
real and all of the irreducible characters take real values. More precisely,
Hecke (see \cite{Hec5}, p.762, p.763 and p.766) showed that for the Steinberg
representation of dimension $p$, the only such representation is of a rational
character. For the degenerate principal series representation of dimension
$\frac{p+1}{2}$, the two representations of this degree are first known in a
form from $\mathbb{Q}(\zeta_p)$, where $\zeta_p=e^{\frac{2 \pi i}{p}}$, while
the field $\mathbb{Q}(\chi)$ is $\mathbb{Q}(\sqrt{p})$. For the principal series
representation of dimension $p+1$, the characters of such a representation are
certain numbers from the real subfield of the field of the primitive $t$-th roots
of unity $\varrho$, where $t$ passes through all divisors of $\frac{p-1}{2}$
except $1$ and $2$. In fact, there are exactly $\frac{1}{2} \varphi(t)$ different
simple characters of degree $p+1$ for each of these $t$. The field
$\mathbb{Q}(\chi)$ is the real subfield $\mathbb{Q}(\varrho+\varrho^{-1})$.
For the discrete series representation of dimension $p-1$, the characters of
such a representation is again a primitive $t$-th roots of unity $\sigma$, where
now $t$ must be a divisor of $\frac{p+1}{2}$ and $t>2$, such that a total of
$\frac{1}{2} \varphi(t)$ different representations belong to the different
$\sigma$ with the same $t$, where $\sigma$ and $\sigma^{-1}$ give the same
representation. The character $\chi$ generates the real field
$\mathbb{Q}(\sigma+\sigma^{-1})$.

  According to \cite{Hec5}, p. 770), the general case of the representations
$V_f^{(l)}(t)$ is that the field $\mathbb{Q}(\chi)$ of their character is a
totally real (abelian) number field of degree $n=\frac{1}{2} \varphi(t)$.
First let $\kappa=1$, such as $p=13$. Each of the three representations
$V_{12}^{(l)}(7)$ ($l=1, 2, 3$) with conjugate coefficients from the real
cubic subfield of the $7$th roots of unity $\mathbb{Q}(\sigma+\sigma^{-1})$
is realized by an integral vector of the first kind with $12$ components. In
general, Hecke (see \cite{Hec5}, p.771, Satz 19) prove that the period problem
for the $n=\frac{1}{2} \varphi(t)$ algebraic conjugate representations $V_f^{(l)}(t)$
$(l=1, \ldots, n)$ is, if the multiplicity $\kappa$ of them in the integral group
is $1$, equivalent to determining the values of Hilbert modular functions of $n$
variables over the totally real field of the $n$-th degree which is determined
by the characters of the representations $V_f^{(l)}(t)$.

  In particular, for the smallest such prime $p=13$, Hecke showed that the
following deep result:
$$H^0(X(13), \Omega_{X(13)}^1)=V_{14} \oplus V_{12}^{(1)} \oplus V_{12}^{(2)}
                        \oplus V_{12}^{(3)},\eqno{(2.2.6)}$$
where $q=13$, $\frac{q-1}{2}=6$.  Hence, $t=3$ or $6$. This implies that
the field $\mathbb{Q}(\chi)$ is $\mathbb{Q}$ and $V_{14}$ is defined over
$\mathbb{Q}$. On the other hand, $V_{12}^{(i)}$ $(i=1, 2, 3)$ are conjugate
over the totally real field $\mathbb{Q}(\zeta_7+\zeta_7^{-1})$ with $\zeta_7$
the seventh root of unity. In this case, Hecke showed that the multiplicity
$m(W)$ is equal to one for $W$ is the principal series representation of
dimension $14$, and the three conjugate discrete series representations of
dimension $12$, and equal to zero for all other irreducible representations
of $\text{SL}(2, \mathbb{F}_{13})$. In particular, the trivial representation
of dimension one, the degenerate principal series of dimension $7$ and the
Steinberg representation of dimension $13$ do not appear in the above
decomposition corresponding to the modular curve $X(13)$. A natural problem is
{\it why these three irreducible representations do not appear and where do
they appear?}

  What makes Hecke's decomposition so appealing is the wide range of mathematics
it touches: it can be regarded as a result in modular forms and modular curves
(for $X(p)$), number theory (for the field of definition), representation theory
(for $\text{SL}(2, \mathbb{F}_p)$), algebraic geometry (for the genus
$g=\dim H^0(X(p), \Omega_{X(p)}^1)$), and homological algebra (for the direct
sum decomposition for the cohomology groups $H^0(X(p), \Omega_{X(p)}^1)$).

\begin{center}
{\large\bf 2.3. Algebraic theory and $\ell$-adic realization of representations
                for $\text{PSL}(2, \mathbb{F}_p)$: Drinfeld-Deligne-Lusztig curves}
\end{center}

  Let us consider a connected, reductive algebraic group $G$, defined over
a finite field $\mathbb{F}_q$, with Frobenius map $F$. We shall be concerned
with the representation theory of the finite group $G^F$, over fields of
characteristic zero. In \cite{DeLu}, Deligne and Lusztig proved that there
is a well defined correspondence which, to any $F$-stable maximal torus $T$
of $G$ and a character $\theta$ of $T^F$ in general position, associates an
irreducible representation of $G^F$; moreover, if $T$ modulo the centre of
$G$ is anisotropic over $\mathbb{F}_q$, the corresponding representation of
$G^F$ is cuspidal. In particular, for $T$ as above and $\theta$ an arbitrary
character of $T^F$ they constructed virtual representations $R_T^{\theta}$
which have all the required properties. These are defined as the alternating
sum of the cohomology with compact support of the variety of Borel subgroups
of $G$ which are in a fixed relative position with their $F$-transform,
with coefficients in certain $G^F$-equivariant locally constant $\ell$-adic
sheaves of rank one. In particular, for $G=\text{SL}(2)$, this gives a
construction made by Drinfeld (see \cite{Lu1978} or \cite{Bon}), who proved
that the discrete series representations of $\text{SL}(2, \mathbb{F}_q)$ occur
in the cohomology of the affine curve $x y^q-x^q y=1$ (the form
$x y^q-x^q y=\begin{vmatrix} x & x^q\\ y & y^q \end{vmatrix}$ is
$\text{SL}(2, \mathbb{F}_q)$-invariant). For simplicity, let now $q=p$
be a prime, assume that $G=\text{SL}(2, \overline{\mathbb{F}}_p)$ with its
standard $\mathbb{F}_p$-rational structure. Then $G^F=\text{SL}(2, \mathbb{F}_p)$.
There is a unique $G^F$-conjugacy class of $F$-stable maximal tori in $G$
which are not contained in any $F$-stable Borel subgroup. If $T$ is such
a torus there is an isomorphism $T \stackrel{\sim}{\longrightarrow}
\overline{\mathbb{F}}_p^{*}$ under which the action of $F$ on $T$ becomes
the action $\lambda \mapsto \lambda^{-p}$ on $\overline{\mathbb{F}}_p^{*}$.
Thus $T^F$ is isomorphic to $H=\{ \lambda \in \overline{\mathbb{F}}_p^{*}:
\lambda^{p+1}=1 \}$ (the non-split torus or anisotropic torus $\mu_{p+1}$).
Let $U$ be the unipotent radical of a Borel subgroup containing $T$. The
variety
$$L^{-1}(U)=\{ g \in G: g^{-1} F(g) \in U \}$$
with the action of $G^F \times T^F$ can be identified with the following
Drinfeld-Deligne-Lusztig curve:
$$\widetilde{X}_p=\left\{ (x, y) \in \overline{\mathbb{F}}_p^2: x y^p-x^p y=1
  \right\},$$
on which $\text{SL}(2, \mathbb{F}_p)$ acts by linear change of coordinates
and $H$ acts by homothety. $\widetilde{X}_p$ is an irreducible affine curve;
it is known that this implies that $H_c^0(\widetilde{X}_p, \mathbb{Q}_{\ell})=0$
and $H_c^2(\widetilde{X}_p, \mathbb{Q}_{\ell})$ is one dimensional with $H$
acting trivially on it. Here, $H_c^i(*, \mathbb{Q}_{\ell})$ denotes $\ell$-adic
\'{e}tale cohomology with compact support. It follows that for any nontrivial
character $\theta: T^F \rightarrow \overline{\mathbb{Q}}_{\ell}^{*}$ or
equivalently, $\theta: H \rightarrow \overline{\mathbb{Q}}_{\ell}^{*}$, the virtual
$\text{SL}(2, \mathbb{F}_p)$-module $\varepsilon_T \varepsilon_G R_T^G(\theta)
=-R_T^G(\theta)$ is an actual representation: it can be realized in the subspace
$H_{\theta^{-1}}$ of $H_c^1(\widetilde{X}_p, \mathbb{Q}_{\ell}) \otimes
\overline{\mathbb{Q}}_{\ell}$ on which $H$ acts according to
$\theta^{-1}$. In particular, $H_{\theta^{-1}}$ is an irreducible cuspidal
$\text{SL}(2, \mathbb{F}_p)$-module (discrete series representation of dimension
$p-1$) provided that $\theta$ is in general position, i.e. when $\theta^2 \neq 1$.
Assume now that $p$ is odd and that $\theta=\epsilon$ where $\epsilon^2=1$,
$\epsilon \neq 1$. It can be shown that the decomposition of $H_{\epsilon}$ into
two irreducible $\text{SL}(2, \mathbb{F}_p)$-module is achieved by taking the two
eigenspaces of the Frobenius $F$ for the obvious $\mathbb{F}_p$-rational structure
of $\widetilde{X}_p$. In particular, these two $\text{SL}(2, \mathbb{F}_p)$-modules
are cuspidal (degenerate discrete series representations of dimension $\frac{1}{2}(p-1)$).

  In general, for $w$ in the Weyl group $W$ of $G$, The Deligne-Lusztig $X(w) \subset X$
is the locally closed subscheme of $X=G/B$ consisting of all Borel subgroups $B$
of $G$ such that $B$ and $F(B)$ are in relative position $w$. In other words,
$$X(w) \cong \{ gB \in G/B: g^{-1} F(g) \in B w B \} \subset X.$$
Let $U$ be the unipotent radical of $B$. Then $\widetilde{X}=G/U \rightarrow
X=G/B$ is a $T$-torsor: $T$ normalizes $U$ and acts on $\widetilde{X}$ from the
right. We define similarly a locally closed subvariety
$$\widetilde{X}(w)=\{ gU \in G/U: g^{-1} F(g) \in U w U \} \subset \widetilde{X}.$$
Deligne-Lusztig in \cite{DeLu} proved that any irreducible representation $\rho$
of $G^F$ appears in the virtual representation
$$R^{\theta}(w)=\sum_{i} (-1)^i H_c^i(\widetilde{X}(w), \overline{\mathbb{Q}}_{\ell})_{\theta}$$
for some $w$, $\theta$. In particular, for $w=e$, $R^{\theta}(e)$ is simply the
usual parabolic induction of a character $\theta$ of $T^F$ (the induced representation).
In the case $G=\text{SL}(2)$, the Weyl group $W=\{ e, w \}$, the nontrivial element
$w$ acts by $\theta \mapsto \theta^{-1}$. The Deligne-Lusztig varieties are given
as follows:
$$X(e)=\mathbb{P}^1(\mathbb{F}_p), \quad X(w)=\mathbb{P}^1-\mathbb{P}^1(\mathbb{F}_p),$$
In particular, $\widetilde{X}(w)=\widetilde{X}_p$ is just the Drinfeld-Deligne-Lusztig
curve. Correspondingly, $R^{\theta}(e)$ are the principal series representations and
$-R^{\theta}(w)$ are the discrete series representations which are irreducible
whenever $\theta^2 \neq 1$.

  Now, we give the relation between the Drinfeld-Deligne-Lusztig curves $\widetilde{X}_p$
and the bad reduction of the modular curves $X(p)$. Let $p$ be a prime, and $K$ be a
finite extension of the $p$-adic field $\mathbb{Q}_p$, with the ring of integers
$\mathcal{O}$ and the residue field $k$ of cardinality $q$. The proof of the local
Langlands correspondence for $\text{GL}(n, K)$, by Harris-Taylor \cite{HT}, was
achieved by showing that the desired correspondence is realized in the $\ell$-adic
vanishing cycle cohomology groups of the deformation spaces of formal $\mathcal{O}$-modules
of height $n$ with Drinfeld level structures (known as non-abelian Lubin-Tate theory
or the conjecture of Deligne-Carayol \cite{Carayol1990}). As these deformation
spaces occur as complete local rings of certain unitary Shimura varieties at the
supersingular points, they made an essential use of the fact that global Langlands
correspondences are realized in the $\ell$-adic \'{e}tale cohomology groups of
these Shimura varieties over CM-fields. In \cite{Yo}, Yoshida gave a purely
local approach to this non-abelian Lubin-Tate theory, in the special case of
depth $0$ or level $\mathfrak{p}$, by computing the local equation of the deformation
space and constructing its suitable resolution to calculate the vanishing
cycle cohomology directly. He showed that, in this case, the non-abelian Lubin-Tate
theory for supercuspidal representations of $\text{GL}(n, K)$ is essentially
equivalent to the Deligne-Lusztig theory for $\text{GL}(n)$ of the residue
field $k$, which realizes the cuspidal representations of $\text{GL}(n, k)$
in the $\ell$-adic cohomology groups of certain varieties over an algebraic
closure $\overline{k}$ of $k$.

  More precisely, let $K$, $\mathcal{O}$, $k$ as above and fix $n \geq 1$. Let
$K^{\text{ur}}$ be the maximal unramified extension of $K$, and let $W$ be the
completion of the ring of integers $\mathcal{O}^{\text{ur}}$ of $K^{\text{ur}}$.
Let $\eta$, $\overline{\eta}$ be the spectra of $\text{Frac}$ $W$ and its fixed
algebraic closure, respectively. Firstly, let $X$ be the spectrum of the deformation
ring of formal $\mathcal{O}$-module of height $n$ with level $\mathfrak{p}$
structure (\cite{Dr1}), which is a scheme of relative dimension $n-1$ over $W$.
We are interested in the $\ell$-adic \'{e}tale cohomology groups
$H^i(X_{\overline{\eta}}, \overline{\mathbb{Q}}_{\ell})$ ($\ell \neq \text{char}(k)$)
of the geometric generic fiber
$X_{\overline{\eta}}:=X \times_{\text{Spec}(W)} \overline{\eta}$, which are
finite dimensional $\text{GL}(n, k) \times I_K$-modules, where $I_K$ is the
inertia group of $K$. Secondly, let DL be the Deligne-Lusztig variety for
$\text{GL}(n, k)$, associated to the element of the Weyl group of $\text{GL}(n)$
that corresponds to the cyclic permutation $(1, \ldots, n)$ in the symmetric
group of $n$ letters, or equivalently to a non-split torus $T$ with
$T(k) \cong k_n^{\times}$ where $k_n$ is the extension of $k$ of degree $n$
(\cite{DeLu}). This DL is a smooth affine variety over $\overline{k}$ with
actions of $\text{GL}(n, k)$ and $T(k) \cong k_n^{\times}$, hence we can
regard $H_c^i(\text{DL}, \overline{\mathbb{Q}}_{\ell})$ as a
$\text{GL}(n, k) \times I_K$-module by the canonical surjection
$I_K \rightarrow k_n^{\times}$. We denote the alternating sums of these
cohomology groups as follows:
$$H^{*}(X_{\overline{\eta}}):=\sum_{i} (-1)^i H^i(X_{\overline{\eta}},
  \overline{\mathbb{Q}}_{\ell}), \quad
  H_c^{*}(\text{DL}):=\sum_{i} (-1)^i H_c^i(\text{DL},
  \overline{\mathbb{Q}}_{\ell}).$$
Then
$$H^{*}(X_{\overline{\eta}})=H_c^{*}(\text{DL}).$$
In particular, Yoshida (see \cite{Yo}) obtained the information concerning
the geometry of $X$ as the following: Let $\varpi$ be a uniformizer of $\mathcal{O}$.
Then the $W$-scheme $X$ is isomorphic to
$$\text{Spec $W[[X_1, \ldots, X_n]]/(P(X_1, \ldots, X_n)-\varpi)$},$$
where $P \in W[[X_1, \ldots, X_n]]$ is a polynomial in $X_1$, $\ldots$, $X_n$
(see \cite{Yo} for more details). This gives the integral local equations of the
relevant unitary Shimura varieties at supersingular points, and in the special
case $K=\mathbb{Q}_p$ and $n=2$, this gives the integral version of Katz-Mazur's
description of the bad reduction of modular curves $X(p)$ (see \cite{KaMa}). That
is to say that the $\ell$-adic \'{e}tale cohomology group $H^{*}(X_{\overline{\eta}})$
($\ell \neq p$) of the geometric generic fiber $X_{\overline{\eta}}$, where $X$ is
the integral local equation of modular curves $X(p)$ at supersingular points, i.e.,
the integral version of Katz-Mazur's description of the bad reduction of modular
curves $X(p)/\mathbb{Z}_p$ is equal to the $\ell$-adic \'{e}tale cohomology group
$H_c^{*}(\widetilde{X}_p)$ of the Drinfeld-Deligne-Lusztig curve $\widetilde{X}_p$.

  Let us recall the basics of integral models of modular curves. Let
$\Gamma$ be one of the congruence subgroups $\Gamma_0(N)$, $\Gamma_1(N)$,
or $\Gamma(N)$. Then the modular curve $X(\Gamma)$ defined over $\mathbb{Q}$,
admits a smooth model over $\mathbb{Z}[1/N]$. That is, modular curves have
good reduction modulo primes which do not divide the level. If $p \nmid N$,
and if $x \in X(\Gamma)(\overline{\mathbb{F}}_p)$ is a geometric point of
the special fiber, then the complete local ring of $X(\Gamma)$ at $x$ is
$$\widehat{\mathcal{O}}_{X(\Gamma), x} \cong W[[t]],$$
where $W=W(\overline{\mathbb{F}}_p)$. When a prime $p$ does divide the level,
one has to construct an integral model of $X(\Gamma)$ over $\mathbb{Z}_p$, and
singularities begin to appear in the special fiber. The first investigation of
the bad reduction of modular curves was given by Deligne-Rapoport \cite{DeRa},
who constructed a model of $X_0(Np)$ over $W$ whose reduction is the union of
two copies of $X_0(N)_{\overline{\mathbb{F}}_p}$, which meet transversely at
the supersingular points. If $x \in X_0(Np)(\overline{\mathbb{F}}_p)$ is a
supersingular point, then the complete local ring of $X_0(Np)$ at $x$ is
$$\widehat{\mathcal{O}}_{X_0(Np), x} \cong W[[t, u]]/(ut-p).$$
The book of Katz-Mazur \cite{KaMa} constructs integral models of the modular
curves $X(\Gamma)$ by defining moduli problems of elliptic curves with the
notion of Drinfeld level structure \cite{Dr1}. Let us start with modular
model over $\mathbb{Z}[\zeta_p]$, as given by Katz and Mazur (see \cite{KaMa},
Chapter 13, or \cite{EdP1} and \cite{EdP2})), for modular curves with full
level $p$ structure plus some additional level structure $\mathcal{P}$ with
nice properties at $p$, where $p>3$ is a prime number. Consider some representable
finite \'{e}tale moduli problem $\mathcal{P}$ over $(\text{Ell})_{\mathbb{Z}_p}$.
The moduli problem $(\mathcal{P}, [\Gamma(p)])$ classifies triples $(E/S, a, \phi)$
for $S$ a $\mathbb{Z}_p$-scheme, $E/S$ an elliptic curve, $a \in \mathcal{P}(E/S)$
and $\phi \in [\Gamma(p)](E/S)$ in the sense of \cite{KaMa}. To state
Katz-Mazur's central result in the simplest way, we will work with so-called
exotic Igusa structure, where the moduli problem on the elliptic stack
$(\text{Ell})_{\mathbb{F}_p}$ associated with the Igusa structure of level
$p$ is denoted by $[\text{Ig}(p)]$. When restrict to level $p$, we define
the moduli problem $[\text{ExIg}(p, 1)]$ that sends $(E/S/\mathbb{F}_p)$ to:
$$[\text{ExIg}(p, 1)](E/S):=\{\text{$P \in E(S)$, $(0, P)$ is a Drinfeld
  $p$-basis of $E/S$}\}$$
and we can check there is an (exotic) isomorphism
$$[\text{Ig}(p)] \stackrel{\sim}{\longrightarrow} [\text{ExIg}(p, 1)], \quad
  (E/S, P \in [\text{Ig}(p)](E/S)) \mapsto (E^{(p)}/S, (0, P)).$$
Kazt-Mazur's theorem about $\Gamma(p)$-structures assert that
$(\mathcal{P}, [\Gamma(p)])$ is representable by a regular $\mathbb{Z}_p$-scheme
$\mathcal{M}(\mathcal{P}, [\Gamma(p)])$, having a compactification we denote
by $\overline{\mathcal{M}}(\mathcal{P}, [\Gamma(p)])$. The existence of Weil's
pairing $e_p(\cdot, \cdot)$ shows that the morphism
$\overline{\mathcal{M}}(\mathcal{P}, [\Gamma(p)]) \rightarrow \text{Spec}(\mathbb{Z}_p)$
factorizes through $\text{Spec}(\mathbb{Z}_p[\zeta_p])$, with
$\mathbb{Z}_p[\zeta_p]:=\mathbb{Z}_p[x]/(x^{p-1}+\cdots+x+1)$. For all integers
$i$ in $\{ 1, \ldots, p-1\}$, set
$$X_i:=\overline{\mathcal{M}}(\mathcal{P}, [\Gamma(p)^{\zeta_p^i-\text{can}}])$$
for the sub-moduli problem over $(\text{Ell})_{\mathbb{Z}_p[\zeta_p]}$ representing
triples
$$(E/S/\mathbb{Z}_p[\zeta_p], a, \phi) \quad \text{such that} \quad
  e_p(\phi(1, 0), \phi(0, 1))=\zeta_p^i.$$
The obvious morphism
$$\coprod_{i \in \mathbb{F}_p^{*}} X_i \rightarrow
  \overline{\mathcal{M}}(\mathcal{P}, [\Gamma(p)])_{\mathbb{Z}_p[\zeta_p]}$$
induces, by normalization, an isomorphism of schemes over $\mathbb{Z}_p[\zeta_p]$:
$$\coprod_{i \in \mathbb{F}_p^{*}} X_i \stackrel{\sim}{\longrightarrow}
  \overline{\mathcal{M}}(\mathcal{P}, [\Gamma(p)])_{\mathbb{Z}_p[\zeta_p]}^{\sim},$$
where $\overline{\mathcal{M}}(\mathcal{P}, [\Gamma(p)])_{\mathbb{Z}_p[\zeta_p]}^{\sim}
\rightarrow \overline{\mathcal{M}}(\mathcal{P}, [\Gamma(p)])_{\mathbb{Z}_p[\zeta_p]}$
is the normalization. The triviality of $p$-th roots of unity in characteristic
$p$ shows that, after the base change $\mathbb{Z}_p[\zeta_p] \rightarrow \mathbb{F}_p$,
the $X_{i, \mathbb{F}_p}$ are not only isomorphic to each other but actually equal.
Moreover, the modular interpretation of a $\Gamma(p)$-structure
$\phi: (\mathbb{Z}/p \mathbb{Z})^2 \rightarrow E(k)$, in the generic case of an
ordinary elliptic curve $E$ over a field $k$ of characteristic $p$, amounts to
choosing some line $L$ in $(\mathbb{Z}/p \mathbb{Z})^2$ which plays the role of
$\text{Ker}(\phi)$, then some point $P$ in $E(k)$ which defines the induced
isomorphism $(\mathbb{Z}/p \mathbb{Z})^2/L  \stackrel{\sim}{\longrightarrow} E[p](k)$.
More precisely, we have

\textbf{Theorem 2.3.1.} (Katz-Mazur). (see \cite{KaMa}, \cite{EdP1} and \cite{EdP2})
{\it Each curve $X_{i, \mathbb{F}_p}$ obtained from $X_i$ over $\mathbb{Z}_p[\zeta_p]$
via $\mathbb{Z}_p[\zeta_p] \rightarrow \mathbb{F}_p$, is the disjoint union,
with crossings at the supersingular points, of $p+1$ copies of the
$\overline{\mathcal{M}}(\mathcal{P})$-schemes}
$\overline{\mathcal{M}}(\mathcal{P}, [\text{ExIg}(p, 1)])$.

  In particular, the $X_i$ now have two kinds of parts: the vertical (or Igusa)
ones, which already showed up in Katz-Mazur model, and the new horizontal
(Drinfeld) ones, with projective models
$$x^p y-x y^p=z^{p+1},$$
which is just the Drinfeld-Deligne-Lusztig curve $\widetilde{X}_p$. In fact,
the situation at the supersingular points can be described as follows. Let
$\mathbf{x}$ be a point of $X_{i, \mathbb{F}_p}$ whose underlying elliptic
curve is supersingular, and let $k$ be the residue field of $\mathbf{x}$.
Then, we have the following:

\textbf{Theorem 2.3.2.} (Katz-Mazur). (see \cite{KaMa} and \cite{EdP1}) {\it
The complete local ring of arithmetic surface $X_i$ at a supersingular point
$\mathbf{x}$ is isomorphic
$$W(k^{\prime})[\zeta_p][[x, y]]/(x^p y-x y^p+g+(1-\zeta_p) f_1),$$
with $k^{\prime}$ the residue field at $\mathbf{x}$, $W(k^{\prime})$ its
ring of Witt vectors, $g$ belongs to the ideal $(x, y)^{p+2}$ and $f_1$
is a unit of $W(k^{\prime})[\zeta_p][[x, y]]$.}

  The horizontal part and the vertical part correspond to the two kinds
of irreducible components of the stable reduction at $p$ of $X(p)$. In
terms of representation theory, the above two parts correspond to the
cuspidal (discrete series) and the principal series representations of
$\text{SL}(2, \mathbb{F}_p)$, respectively.

\begin{center}
{\large\bf 2.4. Arithmetic theory and $\mathbb{Q}(\zeta_p)$-rational
                realization of representations for
                $\text{PSL}(2, \mathbb{F}_p)$: $\mathcal{L}(X(p))$}
\end{center}

  Let $p \geq 7$ be a prime number. According to \cite{Ad2}, denote
by $L^2(\mathbb{F}_p)$ the $p$-dimensional complex vector space of
all square-integrable complex valued functions on $\mathbb{F}_p$
with respect to counting measure, that is, all functions from
$\mathbb{F}_p$ to the complex numbers. We can decompose
$L^2(\mathbb{F}_p)$ as the direct sum of the space $V^{+}$ of even
functions and the space $V^{-}$ of odd functions. The space $V^{-}$
has dimension $\frac{p-1}{2}$ and its associated projective space
$\mathbb{P}(V^{-})$ has dimension $\frac{p-3}{2}$. If $f$ is a
nonzero element of $V^{-}$, we will denote by $[f]$ the
corresponding element of $\mathbb{P}(V^{-})$, in keeping with the
classical notation for homogeneous coordinates.

  Klein discovered the following general result:

\textbf{Theorem 2.4.1.} {\it The modular curve $X(p)$ is isomorphic to
the locus $\mathcal{L}(X(p))$ of all $[f]$ in $\mathbb{P}(V^{-})$ which
for all $w$, $x$, $y$, $z$ in $\mathbb{F}_p$ satisfy the identities
$$\aligned
  0 &=f(w+x) f(w-x) f(y+z) f(y-z)\\
    &+f(w+y) f(w-y) f(z+x) f(z-x)\\
    &+f(w+z) f(w-z) f(x+y) f(x-y).
\endaligned\eqno{(2.4.1)}$$}

  Let us give some background. In his lectures \cite{Schwarz}, Weierstrass
placed the existence of an addition theorem of the elliptic functions at
the top of the whole development and shares the proposition that the
existence of such a theorem is a characteristic property of elliptic
functions and their degenerations. As an entry into the development of
the addition theorems, we choose a tripartite sigma relation (see
\cite{Fricke}, Vol. II, p. 157, \S 1), also established by Weierstrass,
from which the main formulae of the addition theorems are to be obtained
essentially only by means of analytical transformations. More precisely,
the tripartite Weierstrass sigma relation is given as follows:
$$\aligned
  &\sigma(u+u_1)\sigma(u-u_1)\sigma(u_2+u_3)\sigma(u_2-u_3)\\
 +&\sigma(u+u_2)\sigma(u-u_2)\sigma(u_3+u_1)\sigma(u_3-u_1)\\
 +&\sigma(u+u_3)\sigma(u-u_3)\sigma(u_1+u_2)\sigma(u_1-u_2)=0,
\endaligned\eqno{(2.4.2)}$$
which can be used as the common source of the further addition formulae.
From this formula, Klein (see \cite{K4}) derived the following formula
about theta values:
$$\aligned
  &\vartheta_1(v+w) \vartheta_1(v-w) \vartheta_1(t+u) \vartheta_1(t-u)\\
 +&\vartheta_1(w+u) \vartheta_1(w-u) \vartheta_1(t+v) \vartheta_1(t-v)\\
 +&\vartheta_1(u+v) \vartheta_1(u-v) \vartheta_1(t+w) \vartheta_1(t-w)=0.
\endaligned\eqno{(2.4.3)}$$

  In \cite{K4}, Klein set:
$$z_{\alpha}=(-1)^{\alpha} \left(\frac{\omega_2}{\pi}\right)^{\frac{3n-1}{2}}
  \cdot \frac{q^{\frac{\alpha^2}{n}} \cdot \vartheta_1(\alpha \omega \pi, q^n)}
  {\vartheta_1^{\prime}(0, q)^n}.\eqno{(2.4.4)}$$
For the $z$ normalised in this way, the following behaviour now results
with a linear transformation of $\omega_1$, $\omega_2$:

1. If $\omega_1^{\prime}=\omega_1+\omega_2$, $\omega_2^{\prime}=\omega_2$,
then immediately comes:
$$z_{\alpha}(\omega_1+\omega_2, \omega_2)=\varepsilon^{\frac{\alpha(\alpha-n)}{2}}
  \cdot z_{\alpha}(\omega_1, \omega_2),\eqno{(2.4.5)}$$
where $\varepsilon=e^{\frac{2 \pi i}{n}}$ is set.

2. If $\omega_1^{\prime}=-\omega_2$, $\omega_2^{\prime}=+\omega_1$, then first
$$z_{\alpha}(-\omega_2, \omega_1)=(-1)^{\alpha}
  \left(\frac{\omega_1}{\pi}\right)^{\frac{3n-1}{2}} \cdot
  \frac{e^{-\frac{\alpha^2 \pi i}{n \omega}} \cdot
  \vartheta_1\left(-\frac{\alpha \pi}{\omega}, e^{-\frac{n \pi i}{\omega}}\right)}
  {\vartheta_1^{\prime}\left(0, e^{-\frac{\pi i}{\omega}}\right)^n}$$
results. By means of known formulae, this is converted into the following
equation
$$z_{\alpha}(-\omega_2, \omega_1)=\frac{(-1)^{\alpha+1}}{\sqrt{(-1)^{\frac{n-1}{2}} n}}
  \left(\frac{\omega_2}{\pi}\right)^{\frac{3n-1}{2}} \cdot
  \frac{\vartheta_1\left(\frac{\alpha \pi}{n}, q^{\frac{1}{n}}\right)}
  {\vartheta_1^{\prime}(0, q)^n}.$$
On the other hand, one easily states the existence of the following
$\vartheta$-relation, which, it seems, has not yet been noticed:
$$\aligned
  &(-1)^{\alpha} \vartheta_1\left(\frac{\alpha \pi}{n}, q^{\frac{1}{n}}\right)\\
 =&(-1)^{\frac{n-1}{2}} \sum_{\beta=1}^{\frac{n-1}{2}}
   (\varepsilon^{\alpha \beta}-\varepsilon^{-\alpha \beta})
   \left((-1)^{\beta} \cdot q^{\frac{\beta^2}{n}} \cdot
   \vartheta_1(\beta \omega \pi, q^n)\right).
\endaligned$$
Therefore it follows:
$$\sqrt{(-1)^{\frac{n-1}{2}} n} \cdot z_{\alpha}(-\omega_2, \omega_1)
 =(-1)^{\frac{n+1}{2}} \sum_{\beta=1}^{\frac{n-1}{2}}
  (\varepsilon^{\alpha \beta}-\varepsilon^{-\alpha \beta}) \cdot
  z_{\beta}(\omega_1, \omega_2).\eqno{(2.4.6)}$$
But this is essentially the behaviour Klein used for $n=7$ and $11$.
In fact, the formulas $(2.4.5)$ and $(2.4.6)$ give an explicit realization
of the two generators $T$ and $S$ of the group $\text{PSL}(2, \mathbb{F}_p)$
in the cyclotomic field $\mathbb{Q}(\zeta_p)$.

  By now assuming $n$ to be an odd prime number, the ratios of the
$\frac{n-1}{2}$ quantities $z_{\alpha}$ shall be regarded as homogeneous
coordinates of a space of $\frac{n-3}{2}$ dimensions. If the absolute
invariant $J$, and thus $\omega$, moves in the complex plane, then the
point $z$ runs through a curve. First of all, note that the ratios of
any two $z_{\alpha}$ can be expanded to integral powers of $t=q^{\frac{2}{n}}$,
where at most a finite number of negative exponents occur. Therefore,
considering well-known properties of the $\vartheta$-functions, it follows
that (see \cite{K4}): {\it The curve of $z$ is an algebraic curve}. In fact,
the curve is thus uniquely related to the Riemann surface spread out over
the $J$-plane, which is the image of the Galois resolvent of the modular
equation for transformation of the $n$-th order. Klein showed that (see
\cite{K4}) the $\frac{n(n^2-1)}{2}$ sheets of this surface are $3$ each at
$J=0$, at $J=1$ to $2$ each, at $J=\infty$ to $n$ each, and nowhere else.
Accordingly Klein showed that: {\it The genus of our curve is
$p=\frac{(n+2)(n-3)(n-5)}{24}$. The curve merges into itself through
$\frac{n(n^2-1)}{2}$ collineations. By virtue of this, their points are
grouped together by $\frac{n(n^2-1)}{2}$ each. The points that belong
together are generally different; only for $J=0$ they coincide to $3$ each,
for $J=1$ to $2$ each, for $J=\infty$ to $n$ each. This characterises $J$
as a rational function of the coordinates. The curve of $z$ is of order
$\frac{(n-3)(n^2-1)}{48}$.}  In particular, for special cases $n=7$ and
$11$, Klein found $4$ and $20$ as the relevant order. Now, we will give a
modern interpretation. Let us consider the special case of line bundles.
Here we can say much more. First let us denote by $\text{Pic}(X)^G$ the
group of $G$-invariant line bundles and by $\text{Pic}(G; X)$ the group
of $G$-linearized line bundles. The following theorem shows that the
group of $G$-invariant line bundles on $X(p)$ is generated by a line
bundle $\lambda$ of degree $\frac{p^2-1}{24}$ which is a $(2p-12)$-th
root of the canonical bundle.

\textbf{Theorem 2.4.2.} (see \cite{Dol}, Theorem 2.4) {\it Assume $p \geq 5$
is prime. Let} $G=\text{PSL}(2, \mathbb{F}_p)$. {\it Then}
$$\text{Pic}(X(p))^G=\text{Pic}(\text{SL}(2, \mathbb{F}_p); X(p))=\mathbb{Z} \lambda,$$
{\it where}
$$\lambda^{2p-12}=K_{X(p)}$$
{\it and}
$$\text{deg} \lambda=\frac{p^2-1}{24}.$$
{\it Moreover}, $\text{Pic}(G; X(p))$ {\it is the subgroup of} $\text{Pic}(X(p))^G$
{\it generated by $\lambda^2$.}

  Theorem 2.4.2 implies that the group of $\text{PSL}(2, \mathbb{F}_p)$-invariant
divisors on $X(p)$ modulo principal divisors is generated by a divisor of
degree $(p^2-1)/12$. The tensor powers of the line bundle $\lambda$
generating the group $\text{Pic}(X(p))^G$ allows one to embed $X(p)$
$\text{SL}(2, \mathbb{F}_p)$-equivariantly in the projective space.

\textbf{Theorem 2.4.3.} (see \cite{Dol}, Theorem 2.6) {\it Assume $p \geq 5$ is prime.
Denote by $V^{-}$ the irreducible representations of} $\text{SL}(2, \mathbb{F}_p)$
{\it of dimension $\frac{p-1}{2}$. Then a base-point-free linear subsystem of
$|\lambda^{(p-3)/2}|$ maps $X(p)$ in $\mathbb{P}(V^{-})=\mathbb{P}^{(p-3)/2}$
onto a curve of degree $(p-3)(p^2-1)/48$.}

  In fact, the image of $X(p)$ in $\mathbb{P}(V^{-})$ described in Theorem
2.4.3 is the $z$-curve of Klein. Now, let us consider the following three
cases:

(1) When $p=7$. Then $\deg \lambda=2$, $K_X=\lambda^2$ is of degree $4$
so that $X(7)$ is a curve of genus $3$. The divisor class $\lambda$ is an even
theta-characteristic on $X(7)$. It is the unique theta characteristic invariant
with respect to the group of automorphisms $G=\text{PSL}(2, \mathbb{F}_7)$ of
$X(7)$. The $z$-curve is a canonical model of $X(7)$, a plane quartic, i.e.,
the Klein quartic curve.

(2) When $p=11$. Then $\deg \lambda=5$, $K_X=\lambda^{10}$ is of degree $50$
so that $X(11)$ is a curve of genus $26$. The $z$-curve is a curve of degree
$20$ in $\mathbb{P}^4$.

(3) When $p=13$. Then $\deg \lambda=7$, $K_X=\lambda^{14}$ is of degree $98$
so that $X(13)$ is a curve of genus $50$. The $z$-curve is a curve of degree
$35$ in $\mathbb{P}^5$.

  From the modern point of view these embeddings can be described as follows.
Recall that $X(p)$ is a compactification of the moduli space of isomorphism
classes of pairs $(E, \phi)$, where $E$ is an elliptic curve and $(e_1, e_2)$
is a basis of its group $E_p \cong \mathbb{F}_p \times \mathbb{F}_p$ of
$p$-torsion points such that the Weil pairing is equal to $e^{2 \pi i/p}$.
Let $O$ be the origin of $E$. Deligne pointed out to us (see \cite{De5}) that
there is a unique way, compatible with $x \stackrel{\sigma}{\mapsto} -x: E \rightarrow E$,
and the obvious action of $\sigma$ on $\mathcal{O}(pO)$ to have each $a \in E_p$
act on $\Gamma(E, \mathcal{O}(pO))$, above the action $x \mapsto x+a$ on $E$,
so that for those actions $\sigma_a$
$$\sigma_a \sigma_b=\sigma_{a+b} \cdot [\text{Weil pairing $(a, b)$}]^{\frac{1}{2}}$$
(square root taken in $\mu_p$, the $p$-th unit root). There is a special basis
$(X_0, \ldots, X_{p-1})$ in the space $\Gamma(E, \mathcal{O}_E(pO))$ which
defines a map
$$f: E \rightarrow \mathbb{P}^{p-1}, \quad x \mapsto (X_0(x), \ldots, X_{p-1}(x))$$
satisfying the following properties:

1. $(X_0(x-e_1), \ldots, X_{p-1}(x-e_1))=(X_1(x), X_2(x), \ldots, X_{p-1}(x), X_0(x))$;

2. $(X_0(x+e_2), \ldots, X_{p-1}(x+e_2))=(X_0(x), \zeta X_1(x), \ldots, \zeta^{p-1} X_{p-1}(x))$,
where $\zeta=e^{2 \pi i/p}$;

3. $(X_0(-x), X_1(-x), \ldots, X_{p-1}(-x))=(X_0(x), X_{p-1}(x), \ldots, X_1(x))$.

\noindent Let $\pi: \mathcal{X}(p) \rightarrow X(p)$ be the universal family of
elliptic curves $(E, e_1, e_2)$ (its fibres over cusps are certain degenerate
curves, $p$-gons of lines). The $p$-torsion points of the fibres determine $p^2$
sections of the elliptic surface $\mathcal{X}(p)$. The functions $X_m$ define a
morphism $\mathcal{X}(p) \rightarrow \mathbb{P}^{p-1}$ whose restriction to the
fibre $(E, e_1, e_2)$ is equal to the map $f$. The functions $z_m=X_m-X_{p-m}$,
$m=1, \ldots, \frac{p-1}{2}$, define a projection of the image of $\mathcal{X}(p)$
in $\mathbb{CP}^{(p-3)/2}$. The image of the section of $\pi$ defined by the $0$-point
is the $z$-curve of Klein. In fact, in his paper \cite{K5}, p. 226, formulas (64)
and (65), Klein introduced the following combinations:
$$Y_0=X_0, \quad Y_1=X_1+X_{n-1}, \ldots, Y_{\frac{n-1}{2}}=X_{\frac{n-1}{2}}+X_{\frac{n+1}{2}},$$
$$Z_0=0, \quad Z_1=X_1-X_{n-1}, \ldots, Z_{\frac{n-1}{2}}=X_{\frac{n-1}{2}}-X_{\frac{n+1}{2}},$$
where the $Y_{\alpha}$'s are {\it odd}, the $Z_{\alpha}$'s {\it even} functions
of $u$. It is proven in \cite{Ve}, Theorem 10.6, that the $z$-curve of Klein is
always nonsingular.

  Thus, $X(p)$ is defined by a collection of quartics which
we can write down explicitly. In the case $p=7$, we obtain the
defining equation of the Klein quartic curve $x^3 y+y^3 z+z^3 x=0$.
In the case $p=11$, we recover Klein's theorem about $X(11)$ (see also
\cite{Ad1}). In \cite{AR}, Adler and Ramanan gave a simple geometric
interpretation: the modular curve $X(p)$ is the intersection of a
Grassmannian and a $2$-uply embedded projective space.

  More precisely, consider the Weil representation of
$\text{SL}(2, \mathbb{F}_p)$ on $L^2(\mathbb{F}_p)$. Tensor this
representation with itself and identify $L^2(\mathbb{F}_p) \otimes
L^2(\mathbb{F}_p)$ with $L^2(\mathbb{F}_p^2)$. Define the operator
$T$ from $L^2(\mathbb{F}_p^2)$ to itself by
$$(T \Phi)(x, y)=\Phi\left(\frac{x+y}{2}, \frac{x-y}{2}\right).$$
Then one can show that $T$ normalizes $\text{SL}(2, \mathbb{F}_p)$
as a group of operators on $L^2(\mathbb{F}_p^2)$ and maps
$\bigwedge\!^2(V^{+})$ isomorphically onto $\text{Sym}^2(V^{-})$.
Passing to projective spaces, we can use $T$ to identify
$\mathbb{P}(\bigwedge\!^2(V^{+}))$ with
$\mathbb{P}(\text{Sym}^2(V^{-}))$. Now, in
$\mathbb{P}(\bigwedge\!^2(V^{+}))$ we have the Grassmannian
$\text{Gr}$ of complex $2$-planes in $V^{+}$ and in
$\mathbb{P}(\text{Sym}^2(V^{-}))$ we have the image $\text{Ver}$
of $\mathbb{P}(V^{-})$ under the $2$-uple embedding. Klein's
equations say precisely that $X(p)$ is the intersection of
$\text{Gr}$ and $\text{Ver}$.

  The significance of (2.4.1) comes from the fundamental relation
among the the theta constants (see \cite{AR}, p. 48, (18.7) or
\cite{Mumford1966}, p. 332, Riemann's theta-relation on
null-functions $q$, (A), or \cite{Igusa1972}, p. 146, \S 3. Theta
relations (under the new formalism)). In particular, Mumford (see
\cite{Mumford1966}) pointed out how Riemann's theta-relation can
be deduced from the null-functions $q$. In fact, this relation
will not be exactly the same as in the classical case, because
our $q$'s are not the same as the usual theta-null values. But it
can be shown that they are related by a non-singular linear
transformation with roots of unity as coefficients: so the
formulas are trivial modifications of each other. In the special
case of elliptic curves, i.e., one-dimensional abelian variety, it
takes a simple form (see \cite{AR}, p. 53, (19.6)):
$$q_L(u+v) \cdot q_L(u-v)=\begin{vmatrix}
  q_{L^2}(j_{\delta}(u)) & q_{L^2}(j_{\delta}(v))\\
  q_{L^2}(j_{\delta}(u)+x) & q_{L^2}(j_{\delta}(v)+x)
  \end{vmatrix}.\eqno{(2.4.7)}$$
Let $u$, $v$, $w$, $z$ be any elements of $K(\delta)$. Applying the
Pl\"{u}cker relations to the right hand side of (2.4.7) we obtain the
following quartic relation:
$$\aligned
  0 &=q_L(u+v) \cdot q_L(u-v) \cdot q_L(w+z) \cdot q_L(w-z)\\
    &+q_L(u+w) \cdot q_L(u-w) \cdot q_L(z+v) \cdot q_L(z-v)\\
    &+q_L(u+z) \cdot q_L(u-z) \cdot q_L(v+w) \cdot q_L(v-w).
\endaligned\eqno{(2.4.8)}$$

  Write $V^{+}$ and $V^{-}$ to denote respectively the spaces of even
and odd functions on $K(\delta)$. Let $\Lambda$ denote the subspace of
$V(\delta, \delta)$ consisting of all functions $f$ on $K(\delta, \delta)
  =K(\delta) \times K(\delta)$ such that for all $u$, $v$ in $K(\delta)$
we have $f(u, v)=-f(v, u)$ and $f(-u, v)=f(u, v)$. Let $\Sigma$ denote
the subspace of $V(\delta, \delta)$ consisting of all functions $g$ on
$K(\delta, \delta)$ such that for all $u$, $v$ in $K(\delta)$ we have
$g(u, v)=g(v, u)$ and $g(-u, v)=-g(u, v)$. If $h$ is any function on
$K(\delta, \delta)$, denote by $\mathcal{T}(h)$ the function on
$K(\delta, \delta)$ defined by
$$\mathcal{T}(h)(u, v)=h\left(\frac{u+v}{2}, \frac{u-v}{2}\right).$$
Then $\mathcal{T}$ is an automorphism of the vector space $V(\delta, \delta)$.
Furthermore one can easily show that $\mathcal{T}$ induces an isomorphism
of $\Lambda$ onto $\Sigma$. If $\alpha$, $\beta$ belong to $V(\delta)$
denote by $\alpha \wedge \beta$ and $\alpha \vee \beta$ the functions
on $K(\delta, \delta)$ defined respectively by
$$(\alpha \wedge\beta)(u, v)=\alpha(u) \beta(v)-\alpha(v) \beta(u)
 =\begin{vmatrix}
  \alpha(u) & \alpha(v)\\
  \beta(u) & \beta(v)
  \end{vmatrix},$$
$$(\alpha \vee \beta)(u, v)=\alpha(u) \beta(v)+\alpha(v) \beta(u).$$
It is straightforward to verify that if $\alpha$ and $\beta$ are both
even functions on $K(\delta)$ then $\alpha \wedge \beta$ will belong
to the space $\Lambda$, while if $\alpha$ and $\beta$ are both odd
then $\alpha \vee \beta$ will belong to the space $\Sigma$. The
expressions $\alpha \wedge \beta$ and $\alpha \vee \beta$ are
bilinear in $\alpha$ and $\beta$. Furthermore, one sees very easily
that $\alpha \wedge \beta$ is alternating in $\alpha$ and $\beta$
and $\alpha \vee \beta$ is symmetric. Denote by $\wedge$ and $\vee$
the linear mappings
$$\wedge: \bigwedge\!^2(V^{+}) \rightarrow \Lambda$$
$$\vee: \text{Sym}^2(V^{-}) \rightarrow \Sigma$$
determined by the pairs $\wedge$ and $\vee$ respectively. They are
easily seen to be isomorphisms. it follows that there is one and only
one isomorphism
$$\tau: \bigwedge\!^2(V^{+}) \rightarrow \text{Sym}^2(V^{-})$$
of $\bigwedge\!^2(V^{+})$ onto $\text{Sym}^2(V^{-})$ such that
$$\vee \circ \tau=\mathcal{T} \circ \wedge.$$

  The group $\text{GL}(V^{+})$ acts on $\bigwedge\!^2(V^{+})$ and
therefore on $\Lambda$ via $\wedge$. There is one and only one
closed orbit for $\text{GL}(V^{+})$ in $\mathbb{P}(\Lambda)$,
namely the image under $\wedge$ of the Grassmannian of all
$2$-planes in $V^{+}$. We will denote that closed orbit by
$\text{Gr}$. By means of the isomorphism $\mathcal{T}$ of
$\Lambda$ onto $\Sigma$, we may identify $\text{Gr}$ with a
subvariety of $\mathbb{P}(\Sigma)$.

  The group $\text{GL}(V^{-})$ acts on $\text{Sym}^2(V^{-})$
and therefore on $\Sigma$ via $\vee$. There is one and only one
closed orbit for $\text{GL}(V^{-})$ in $\mathbb{P}(\Sigma)$,
namely the image of $\mathbb{P}(V^{-})$ in $\mathbb{P}(\Sigma)$
under the mapping $[h] \mapsto [h \vee h]$. Denote the image
of this mapping by $\mathcal{V}$. (This notation, in which
$\mathcal{V}$ is the initial of the name Veronese, is motivated
by the case of the Veronese surface in $\mathbb{P}^5$, which is
a special case of this construction. In general, the variety
$\mathcal{V}$ would be referred to as a $2$-uply embedded
projective space and the embedding $[h] \mapsto [h \vee h]$
as the $2$-uple embedding.) Since $\mathcal{V}$ is a subvariety
of $\mathbb{P}(\Sigma)$ and since $\text{Gr}$ has been identified
with a subvariety of $\mathbb{P}(\Sigma)$, we can consider the
intersection of $\mathcal{V}$ and $\text{Gr}$. Denote by
$\mathcal{L}$ the preimage in $\mathbb{P}(V^{-})$ of
$\mathcal{V} \cap \text{Gr}$ under the isomorphism of
$\mathbb{P}(V^{-})$ onto $\mathcal{V}$ induced by $\vee$.
Then $\mathcal{L}$ is the locus of all $[h]$ in $\mathbb{P}(V^{-})$
such that for some $\alpha$, $\beta$ in $V^{+}$ we have
$$h \vee h=\mathcal{T}(\alpha \wedge \beta).$$
It is the same to say that $[h]$ lies in $\mathcal{L}$ if and only if
for some $\alpha$, $\beta$ in $V^{+}$ we have
$$h(u+v) \cdot h(u-v)=\begin{vmatrix}
                      \alpha(u) & \alpha(v)\\
                      \beta(u) & \beta(v)
                      \end{vmatrix}\eqno{(2.4.9)}$$
for all $u$, $v$ in $K(\delta)$. In particular, if the function $q_L$
is not identically $0$ the point $[q_L]$ of $\mathbb{P}(V^{-})$ lies
in $\mathcal{L}$. Applying the Pl\"{u}cker relations to (2.4.9) we
obtain
$$\aligned
  0 &=h(u+v) \cdot h(u-v) \cdot h(w+z) \cdot h(w-z)\\
    &+h(u+w) \cdot h(u-w) \cdot h(z+v) \cdot h(z-v)\\
    &+h(u+z) \cdot h(u-z) \cdot h(v+w) \cdot h(v-w)
\endaligned\eqno{(2.4.10)}$$
for every $[h]$ in $\mathcal{L}$ and all $u$, $v$, $w$, $z$ in $K(\delta)$.
Here, let us recall the Pl\"{u}cker relations which are used to represent
lines in $\mathbb{P}^3$. The coordinates of the line $u$ passing through
two points $p_1(X_1, Y_1, Z_1, W_1)$ and $p_2(X_2, Y_2, Z_2, W_2)$ is
given by the determinants of the six $2 \times 2$ sub-matrices of the
following matrix:
$$[p_1, p_2]=\left[\begin{matrix}
              X_1 & X_2\\
              Y_1 & Y_2\\
              Z_1 & Z_2\\
              W_1 & W_2
             \end{matrix}\right].$$
In other words, $u=(l_{41}; l_{42}, l_{43}; l_{23}; l_{31}; l_{12})$,
where
$$l_{41}=\begin{vmatrix} W_1 & W_2\\ X_1 & X_2 \end{vmatrix}, \quad
  l_{42}=\begin{vmatrix} W_1 & W_2\\ Y_1 & Y_2 \end{vmatrix}, \quad
  l_{43}=\begin{vmatrix} W_1 & W_2\\ Z_1 & Z_2 \end{vmatrix},$$
$$l_{23}=\begin{vmatrix} Y_1 & Y_2\\ Z_1 & Z_2 \end{vmatrix}, \quad
  l_{31}=\begin{vmatrix} Z_1 & Z_2\\ X_1 & X_2 \end{vmatrix}, \quad
  l_{12}=\begin{vmatrix} X_1 & X_2\\ Y_1 & Y_2 \end{vmatrix}.$$
These coordinates $l_{ij}$ are called the Pl\"{u}cker coordinates of
the line. The six Pl\"{u}cker coordinates are sufficient to describe
the line. The coordinates are not independent, they satisfy the
Pl\"{u}cker relations
$$l_{41} l_{23}+l_{42} l_{31}+l_{43} l_{12}=0.$$

  Conversely, suppose that $[h]$ is a point of $\mathbb{P}(V^{-})$
such that (2.4.10) holds for all $u$, $v$, $w$, $z$ in $K(\delta)$. Then
$h$ does not vanish identically so we can choose an element $t$ in
$K(\delta)$ such that $h(t) \neq 0$. Let $\alpha(u)=c \cdot h(u)^2$
where $c$ is a constant to be chosen later and let $\beta(u)=h(t+u)
\cdot h(t-u)$. Then the function $\alpha \wedge \beta$ is given by
$$(\alpha \wedge \beta)(r, s)=c [h(r)^2 \cdot h(t+s) \cdot h(t-s)
                                -h(s)^2 \cdot h(t+r) \cdot h(t-r)].$$
If we take $u=0$, $v=r$, $w=s$ and $z=t$ in (2.4.10) we obtain
$$\aligned
  0 &=h(r) \cdot h(-r) \cdot h(s+t) \cdot h(s-t)\\
    &+h(s) \cdot h(-s) \cdot h(t+r) \cdot h(t-r)\\
    &+h(t) \cdot h(-t) \cdot h(r+s) \cdot h(r-s).
\endaligned$$
Using the fact that $h$ is an odd function, we conclude that
$$(\alpha \wedge \beta)(r, s)=c \cdot h(t)^2 \cdot h(r+s) \cdot h(r-s).$$
If we now take $c=h(t)^{-2}$, we see that
$$h(r+s) \cdot h(r-s)=(\alpha \wedge \beta)(r, s)
 =\begin{vmatrix}
  \alpha(r) & \alpha(s)\\
  \beta(r) & \beta(s)
  \end{vmatrix}$$
for all $r$, $s$ in $K(\delta)$, which proves that $h$ lies in
$\mathcal{L}$. We have therefore proved the following result.

\textbf{Theorem 2.4.4.} (see \cite{AR}, Theorem 19.17). {\it The locus defined
by} (2.4.10) {\it is the same as $\mathcal{L}$. In particular, the fundamental
relation among the theta constants for one-dimensional abelian varieties
says that the point $[q_L]$ of $\mathbb{P}(V^{-})$ lies in the intersection
of a Veronese variety and a Grassmannian, or what is the same, in the
intersection of the minimal orbits for} $\text{GL}(V^{-})$ {\it and}
$\text{GL}(V^{+})$ {\it in $\mathbb{P}(\Sigma)$.}

  Now, we give the geometry of the modular curve via the fundamental relation
(see \cite{AR}). We will take
$$\delta=(p),$$
so that
$$K(\delta)=\mathbb{Z}/p \mathbb{Z}.$$
For every element $t$ of $K(\delta)=\mathbb{Z}/p \mathbb{Z}$,
let $E_t$ denote the linear form on the space $V^{-}$ of odd functions
on $K(\delta)$ given by $E_t(h)=h(t)$. If $w$, $x$, $y$, $z$ are
any elements of $K(\delta)$, denote by $\Phi_{w, x, y, z}$ the
quartic form given by
$$\aligned
  \Phi_{w, x, y, z} &=E_{w+x} \cdot E_{w-x} \cdot E_{y+z} \cdot E_{y-z}+\\
                    &+E_{w+y} \cdot E_{w-y} \cdot E_{z+x} \cdot E_{z-x}+\\
                    &+E_{w+z} \cdot E_{w-z} \cdot E_{x+y} \cdot E_{x-y}.
\endaligned\eqno{(2.4.11)}$$
Then $\mathcal{L}$ is defined by the equations
$$\Phi_{w, x, y, z}=0$$
with $w$, $x$, $y$, $z$ in $\mathbb{Z}/p \mathbb{Z}$. Since $E_{-t}=-E_t$ for
every $t \in K(\delta)$, we have $\Phi_{-w, x, y, z}=\Phi_{w, x, y, z}$.
Furthermore, we have
$$\Phi_{x, w, y, z}=\Phi_{x, y, z, w}=-\Phi_{w, x, y, z}.\eqno{(2.4.12)}$$
Since the odd permutations
$$\left(\begin{matrix}
  1 & 2 & 3 & 4\\
  2 & 1 & 3 & 4
  \end{matrix}\right) \quad \text{and} \quad
  \left(\begin{matrix}
  1 & 2 & 3 & 4\\
  2 & 3 & 4 & 1
  \end{matrix}\right)$$
generate the group of all permutations on four objects, it follows that
if $(a, b, c, d)=\sigma(w, x, y, z)$ is a permutation of $(w, x, y, z)$
then
$$\Phi_{a, b, c, d}=(-1)^{\sigma} \Phi_{w, x, y, z},\eqno{(2.4.13)}$$
where $(-1)^{\sigma}$ denotes the sign of the permutation $\sigma$. It
follows from this and from the sentence preceding (2.4.12) that for all
choices of signs we have
$$\Phi_{\pm w, \pm x, \pm y, \pm z}=\Phi_{w, x, y, z}.\eqno{(2.4.14)}$$
It follows that each quartic $\Phi_{w, x, y, z}$ is equal, up to a sign,
to a quartic $\Phi_{a, b, c, d}$ with
$$0 \leq a < b < c < d \leq \frac{p-1}{2}.$$
The locus $\mathcal{L}$ is therefore defined by
$\left(\begin{matrix} m\\ 4 \end{matrix}\right)$ quartics, where $2m-1=p$.
In general, these quartics are not distinct.

  For every $t \neq 0$ in $K(\delta)$, let $h_t$ denote the element of
$V^{-}$ defined by
$$h_t(x)=\left\{\aligned
         1 \quad &\text{if $x=t$},\\
         -1 \quad &\text{if $x=-t$},\\
         0 \quad &\text{if $x \neq \pm t$},
\endaligned\right.$$
and let $\kappa_t=[h_t]$ be the point of $\mathbb{P}(V^{-})$ determined
by $h_t$. Then we clearly have
$$h_{-t}=-h_t$$
and
$$\kappa_{-t}=\kappa_t$$
for all $t \neq 0$ in $K(\delta)$.

\textbf{Proposition 2.4.5.} {\it Let $\kappa=[h]$ be a point of $\mathcal{L}$.
If
$$h(s)=0$$
for some $s \neq 0$ in $K(\delta)$, then
$$\kappa=\kappa_t$$
for some $t \neq 0$ in $K(\delta)$.}

{\it Proof}. Let
$$Z=\{ x \in K(\delta): h(x)=0\}$$
and let
$$Z^{\prime}=\{ x \in K(\delta): h(x) \neq 0 \}$$
denote the complement of $Z$ in $K(\delta)$. Since $h$ is odd, $Z=-Z$
and by hypothesis $Z \neq \{ 0 \}$. If we put $y=0$ and $z=w+x$ we find
that (2.4.11) becomes
$$0=-h(w+x)^3 \cdot h(w-x)+h(w)^3 \cdot h(w+2x)-h(x)^3 \cdot h(2w+x).$$
If $h(x)=0$ then for all $w$ in $K(\delta)$ we have
$$h(w)^3 \cdot h(w+2x)=h(w+x)^3 \cdot h(w-x).$$
In particular,
$$h(w) \cdot h(w+2x)=0$$
if and only if
$$h(w+x) \cdot h(w-x)=0.$$
Since $h(0)=0$ we certainly have
$$h(0) \cdot h(2x)=0.$$
It follows by induction that
$$h(rx) \cdot h((r+2)x)=0$$
for every positive integer $r$. Since $p$ is prime, $x$ generates
$K(\delta)$. Consequently every element $w$ of $K(\delta)$ is of the
form $rx$ for some positive integer $r$ and we have
$$h(w) \cdot h(w+2x)=0$$
for all $w$ in $K(\delta)$. Therefore, if $w$ lies in $Z^{\prime}$ then
$w+2x$ does not. Since $x$ was an arbitrary nonzero element of $Z$ we
conclude that if $a$, $b$ lies in $Z^{\prime}$ and $a \neq b$ then
$(a-b)/2$ lies in $Z^{\prime}$. Furthermore if $a$, $b$ are elements
of $Z^{\prime}$ such that $a \neq -b$ then since $-b$ lies in
$Z^{\prime}$ we see that $(a+b)/2$ lies in $Z^{\prime}$. If $\kappa$
is not of the form $\kappa_t$ then we can find $a$, $b$ in $Z^{\prime}$
such that $a \neq b$ and $a \neq -b$. Then $(a+b)/2$ and $(a-b)/2$ lie
in $Z^{\prime}$ and $(a+b)/2 \neq (a-b)/2$. Therefore
$$\frac{\frac{a+b}{2}+\frac{a-b}{2}}{2}=\frac{a}{2}$$
lies in $Z^{\prime}$. It follows that if $\kappa$ is not of the form
$\kappa_t$ then
$$\frac{1}{2} Z^{\prime}=Z^{\prime},$$
or what is the same,
$$Z^{\prime}=2 Z^{\prime}.$$
In particular, since $(a+b)/2$ lies in $Z^{\prime}$, so do $a+b$ and
$2a$. Suppose that $r>1$ and that for $1 \leq j \leq r$ the element
$ja$ lies in $Z^{\prime}$. If $(r+1)a=0$ then since $p$ is prime the
elements $ja$ with $1 \leq j \leq r$ comprise all of the nonzero
elements of $K(\delta)$ and so $h(t) \neq 0$ for $t \neq 0$ which
contradicts our hypotheses. Therefore $(r+1)a \neq 0$ and $ra \neq -a$,
so $a+ra=(r+1)a$ belongs to $Z^{\prime}$. It follows that for all
$r>0$ the element $ra$ lies in $Z^{\prime}$ and in particular $pa=0$
lies in $Z^{\prime}$, which is a contradiction. Therefore, $\kappa$
is of the form $\kappa_t$ for some $t \neq 0$.

  Now, we will compute the group of collineations of $\mathbb{P}(V^{-})$
which leave the locus $\mathcal{L}$ invariant. That group will be denoted
$\text{Aut}(\mathcal{L})$. We begin by showing that $\text{PSL}(2, \mathbb{F}_p)$
is isomorphic to a subgroup of $\text{Aut}(\mathcal{L})$. We identify the
underlying sets of $\mathbb{F}_p$ and $K(\delta)$. If $t$ is a nonzero
element of $\mathbb{F}_p$, let $A_t$ and $B_t$ denote the operators on
$V=V(\delta)$ given by
$$\aligned
  (A_t f)(x) &=\zeta^{t x^2} f(x),\\
  (B_t f)(x) &=\left\{\aligned
               \pm \frac{1}{\sqrt{p}} \sum_{y \in \mathbb{F}_p} f(-y) \zeta^{txy},
               \quad \text{$p \equiv 1$ (mod $4$)},\\
               \pm i \frac{1}{\sqrt{p}} \sum_{y \in \mathbb{F}_p} f(-y) \zeta^{txy},
               \quad \text{$p \equiv 3$ (mod $4$)},
               \endaligned\right.
\endaligned$$
for all $f$ in $V$ and all $x$ in $\mathbb{F}_p$, where $\zeta$ is a
primitive $p$-th root of unity and the sign $\pm$ depends on whether
$t$ is multiplied by a quadratic residue or quadratic non-residue.
In fact, the Gaussian sum is given by
$$G_{\sigma}=\sum_{\text{$x$ mod $p$}} e^{2 \pi i \sigma \frac{x^2}{p}}
 =\sqrt{p} \left(\frac{\sigma}{p}\right) \varepsilon_0,$$
where $\left(\frac{\sigma}{p}\right)$ is the Legendre symbol,
$$\varepsilon_0=\left\{\aligned
                1, & \quad \left(\frac{-1}{p}\right)=1,\\
                i, & \quad \left(\frac{-1}{p}\right)=-1.
\endaligned\right.$$
and $\sigma$ is an integer such that $p \nmid \sigma$. Thus,
$$G_{\sigma}=\left\{\aligned
             \pm \sqrt{p}, \quad & \text{$p \equiv 1$ (mod $4$)},\\
             \pm i \sqrt{p}, \quad & \text{$p \equiv 3$ (mod $4$)}.
\endaligned\right.$$
It is known that there is one and only one representation $\rho_t$
of $\text{SL}(2, \mathbb{F}_p)$ on $V$ such that
$$\aligned
  \rho_t\left(\left(\begin{matrix}
               1 & 1\\
               0 & 1
             \end{matrix}\right)\right) &=A_t,\\
  \rho_t\left(\left(\begin{matrix}
               0 & 1\\
              -1 & 0
             \end{matrix}\right)\right) &=B_t.
\endaligned$$
Furthermore, the center of $\text{SL}(2, \mathbb{F}_p)$ acts by a scalar
on $V^{+}$ and on $V^{-}$, although the scalars are different values of
$\pm 1$. The action of $\text{SL}(2, \mathbb{F}_p)$ on $V$ induces an
action on $\bigotimes^2 V$. Call this representation $\rho_t^{\prime}$.
If we identify $\bigotimes^2 V$ with the space of all functions on
$K(\delta, \delta)$, then for all $\gamma$ in $\text{SL}(2, \mathbb{F}_p)$
we have
$$(\rho_t^{\prime}(\gamma)h)(x, y)
 =(\rho_t(\gamma)f)(x) \cdot (\rho_t(\gamma)g)(y)$$
where
$$h(x, y)=f(x) \cdot g(y).$$
It follows at once that the representation $\rho_t^{\prime}$ leaves
$\Lambda$ and $\Sigma$ invariant.

\textbf{Lemma 2.4.6.} {\it The operator $\mathcal{T}$ on $\bigotimes^2 V$
which associates to $f(x, y)$ the function
$$(\mathcal{T} f)(x, y)=f(x+y, x-y)$$
is an intertwining operator from $\rho_1^{\prime}$ to $\rho_2^{\prime}$.}

{\it Proof}. We have to prove that if $\gamma$ belongs to
$\text{SL}(2, \mathbb{F}_p)$ then
$$\mathcal{T} \circ (\rho_1^{\prime}(\gamma))
 =(\rho_2^{\prime}(\gamma)) \circ \mathcal{T}.$$
It is enough to verify this on a set of generators of $\otimes^2 \mathcal{T}$.
If $f$, $g$ belong to $V$, let
$$h(x, y)=f(x) \cdot g(y).$$
Then we must show that
$$\mathcal{T} \circ (\rho_1^{\prime}(\gamma))(h)
 =(\rho_2^{\prime}(\gamma)) \circ \mathcal{T}(h)$$
or, what is the same, that
$$(\rho_1(\gamma)f)(x+y) \cdot (\rho_1(\gamma)g)(x-y)
 =(\rho_2^{\prime}(\gamma) \circ \mathcal{T})(h)(x, y).$$
To prove this, it is enough to verify it for
$$\gamma=\gamma_1=\left(\begin{matrix}
                   1 & 1\\
                   0 & 1
                  \end{matrix}\right)$$
and for
$$\gamma=\gamma_2=\left(\begin{matrix}
                   0 & 1\\
                  -1 & 0
                  \end{matrix}\right),$$
since these elements generate $\text{PSL}(2, \mathbb{F}_p)$. For
$\gamma=\gamma_1$ and any function $M(x, y)$ we have
$$(\rho_t^{\prime}(\gamma)M)(x, y)=M(x, y) \zeta^{t(x^2+y^2)}.$$
In particular, we have
$$\aligned
   (\mathcal{T} \circ \rho_1^{\prime}(\gamma))(h)(x, y)
 &=\rho_1^{\prime}(\gamma)(h)(x+y, x-y)\\
 &=c h(x+y, x-y) \zeta^{(x+y)^2+(x-y)^2}\\
 &=c h(x+y, x-y) \zeta^{2(x^2+y^2)}\\
 &=(\rho_2^{\prime}(\gamma) \circ \mathcal{T})(h)(x, y).
\endaligned$$
If $\gamma=\gamma_2$ then
$$(\rho_t^{\prime}(\gamma)M)(x, y)=c \sum_{u, v \in \mathbb{F}_p}
  M(u, v) \zeta^{t(ux+vy)}$$
for every $M$ in $\otimes^2 V$. In particular, we see that
$$\aligned
   (\mathcal{T} \circ \rho_1^{\prime}(\gamma))(h)(x, y)
 &=(\rho_1^{\prime}(\gamma)(h))(x+y, x-y)\\
 &=c \sum_{u, v \in \mathbb{F}_p} h(u, v) \zeta^{u(x+y)+v(x-y)}\\
 &=c \sum_{u, v \in \mathbb{F}_p} h(u, v) \zeta^{(u+v)x+(u-v)y}\\
 &=c \sum_{a, b \in \mathbb{F}_p} h(a+b, a-b) \zeta^{2(ax+by)}\\
 &=(\rho_2^{\prime}(\gamma) \circ \mathcal{T})(h)(x, y).
\endaligned$$
This completes the proof of Lemma 2.4.6.

  We will refer to the operator $\mathcal{T}$ constructed in the
Lemma 2.4.6 as the fundamental intertwining operator. In fact,
Lemma 2.4.6 can be proved generally by the Weil representations
(see \cite{Ad2}, section 8 and \cite{Weil1964}) which do not use
the computation of the $\text{SL}(2)$ action.

  Since $\text{Gr} \subset \mathbb{P}(\Lambda)$ is invariant under
$\text{GL}(V^{+})$, it is invariant under $\rho_2^{\prime}(\gamma)$
for every $\gamma$ in $\text{SL}(2, \mathbb{F}_p)$. Since $\mathcal{T}$
intertwines $\rho_1^{\prime}$ and $\rho_2^{\prime}$, it follows that
when we identify $\text{Gr}$ with a subvariety of $\mathbb{P}(\Sigma)$
via $\mathcal{T}$, that subvariety of $\mathbb{P}(\Sigma)$ is invariant
under $\rho_2^{\prime}(\gamma)$ for every $\gamma$ in
$\text{SL}(2, \mathbb{F}_p)$. On the other hand, the variety
$\mathcal{V}$ is invariant under $\text{GL}(V^{-})$ and therefore
under $\rho_2^{\prime}(\gamma)$ for every $\gamma$ in
$\text{SL}(2, \mathbb{F}_p)$. Therefore the intersection
$\text{Gr} \cap \mathcal{V}$ is invariant under
$\rho_2^{\prime}(\gamma)$ for every $\gamma$ in
$\text{SL}(2, \mathbb{F}_p)$. For any $h$ in $V$, we have
$$\rho_2^{\prime}(\gamma)(h \vee h)=(\rho_2(\gamma)(h)) \vee
                                    (\rho_2(\gamma)(h))$$
for all $\gamma$ in $\text{SL}(2, \mathbb{F}_p)$. Therefore
$\mathcal{L}$, being the locus of all $h$ in $\mathbb{P}(V^{-})$
such that $h \vee h$ lies in $\text{Gr} \cap \mathcal{V}$, is
invariant under $\rho_2^{\prime}(\gamma)$ for every $\gamma$
in $\text{SL}(2, \mathbb{F}_p)$. Since the center of
$\text{SL}(2, \mathbb{F}_p)$ acts trivially on $\mathbb{P}(V^{-})$,
it follows that the group of collineations of $\mathbb{P}(V^{-})$
of the form $\rho_2(\gamma)$ is isomorphic to
$\text{PSL}(2, \mathbb{F}_p)$. Since $\mathcal{L}$ contains the
generators $\kappa_t$ for $\mathbb{P}(V^{-})$, it follows that the
mapping of $\text{PSL}(2, \mathbb{F}_p)$ into $\text{Aut}(\mathcal{L})$
is injective. Therefore we can identify $\text{PSL}(2, \mathbb{F}_p)$
with a subgroup of $\text{Aut}(\mathcal{L})$. In characteristic $0$,
it is actually true that $\text{PSL}(2, \mathbb{F}_p)$ is equal to
$\text{Aut}(\mathcal{L})$. To see this one can appeal to results of
Brauer \cite{Brauer} characterizing $\text{PSL}(2, \mathbb{F}_p)$ as
a maximal finite subgroup of the collineation group of $\mathbb{P}(V^{-})$.

\textbf{Theorem 2.4.7.} (see \cite{Brauer}, Theorem 4) {\it Let $\mathfrak{G}$
be a group of order $g=p g^{\prime}$ with $(p, g^{\prime})=1$, which has
no normal subgroup of order $p$. If $\mathfrak{G}$ has a
$(1$--$1)$-representation $ \mathfrak{Z}$ of degree $n=(p-1)/2$, then
the factor group of $\mathfrak{G}$ modulo the center $\mathfrak{C}$ of
$\mathfrak{G}$ is isomorphic with} $\text{PSL}(2, \mathbb{F}_p)$. {\it
In other words: $\mathfrak{Z}$, considered as a collineation group,
represents} $\text{PSL}(2, \mathbb{F}_p)$ {\it isomorphically.}

  On the other hand, we can see directly that
$\text{PSL}(2, \mathbb{F}_p)=\text{Aut}(\mathcal{L})$ as soon as
we determine the one-dimensional component of $\mathcal{L}$. That
$\mathcal{L}$ has a one dimensional component follows from Corollary
2.4.9 below. We will show that if $\mathcal{L}$ contains a curve then
the curve is irreducible. To see this, let $\mathcal{C}$ denote the
union of all of the curves lying in $\mathcal{L}$. It will be enough
to show that $\mathcal{C}$ is irreducible.

\textbf{Lemma 2.4.8.} {\it The tangent variety to $\mathcal{L}$ at the
point $\kappa_t$ is the line joining $\kappa_t$ and $\kappa_{3t}$.}

{\it Proof}. Since $\text{Aut}(\mathcal{L})$ contains an automorphism
which maps $\kappa_1$ to $\kappa_t$ and $\kappa_3$ to $\kappa_{3t}$,
we may assume that $t=1$. If a cubic monomial $E_i E_j E_k$ with
$1 \leq i, j, k \leq (p-1)/2$ does not vanish at $\kappa_1$, we
must have $i=j=k=1$. Therefore the first partial derivative of
$\Phi_{w, x, y, z}$ at $\kappa_1$ will all vanish unless one of
the terms of $\Phi_{w, x, y, z}$ is of the form $E_1^3 E_u$.
Suppose that some partial derivative of $\Phi_{w, x, y, z}$ does
not vanish at $\kappa_1$. Then without loss of generality we may
suppose that the first term of $\Phi_{w, x, y, z}$ is $E_1^3 E_u$.
Furthermore, we may suppose that $w+x=w-x=y+z=1$ and $y-z=u$. Then
$w=1$, $x=0$, $u=2y-1=1-2z$ and $\Phi_{w, x, y, z}$ is equal to
$$E_1^3 \cdot E_{2y-1}-E_{y+1} \cdot E_{y-1}^2+E_y^3 \cdot E_{y-2}.$$
$\Phi_{w, x, y, z}$ will be nonzero if and only if $y$ is not equal
to $0$, $1$, $2$ or $(p+1)/2$. Furthermore, if $\Phi_{w, x, y, z}$ is
nonzero then the first term of $\Phi_{w, x, y, z}$ will be the only
term of the form $E_1^3 E_u$. Therefore the only partial derivative
which will not vanish at $\kappa_1$ is the partial with respect to
$E_u$. It follows that if $\kappa$ is a point lying in the tangent
linear variety to $\mathcal{L}$ at $\kappa_1$, we must have
$E_u(\kappa) \neq 0$. If $y$ is not equal to $0$, $1$, $2$ or $(p+1)/2$
then $u$ is not equal to $-1$, $1$, $3$ or $0$. In particular we must
have $E_u(\kappa)=0$ for all $1 \leq u \leq (p+1)/2$ except $u=1$ and
$u=3$. This shows that $\kappa$ lies in the tangent linear variety to
$\mathcal{L}$ at $\kappa_1$ i and only if $\kappa$ lies on the line
joining $\kappa_1$ and $\kappa_3$.

\textbf{Corollary 2.4.9.} {\it Assume that the curve part $\mathcal{C}$ of
the locus $\mathcal{L}$ is non-empty. Then the curve $\mathcal{C}$ is
irreducible. Furthermore, the point $\kappa_t$ lies on $\mathcal{C}$
and is a simple point of $\mathcal{C}$.}

{\it Proof}. By the above Lemma 2.4.8, the tangent linear variety to
$\mathcal{L}$ at $\kappa_t$ has dimension one. Since the curve part
$\mathcal{C}$ of $\mathcal{L}$ is nonempty, the intersection of any
component of $\mathcal{C}$ with the hyperplane $E_x=0$ is nonempty.
On the other hand, we know from the above Lemma 2.4.8 that the hyperplane
$E_x=0$ can only meet $\mathcal{L}$ in points of the form $\kappa_t$.
Since $\mathcal{C}$ is contained in $\mathcal{L}$, it follows that
$E_x=0$ can only meet a component of $\mathcal{C}$ in points of the
form $\kappa_t$. Therefore, every component of $\mathcal{C}$ contains
a point of the form $\kappa_u$. By the above Lemma 2.4.8, we know that
if $\kappa_t$ lies on $\mathcal{C}$ then $\kappa_t$ is a simple point
of $\mathcal{C}$. Consequently, only one irreducible component of
$\mathcal{C}$ can pass through each point $\kappa_u$. Since there
are exactly $(p-1)/2$ points $\kappa_u$, it follows that the number
of irreducible components of $\mathcal{C}$ is at most $(p-1)/2$. On
the other hand, the orbit of $\kappa_t$ under the action of
$\text{PSL}(2, \mathbb{F}_p)$ contains all of the points of
the form $\kappa_u$. Since $\mathcal{C}$ is invariant under
$\text{PSL}(2, \mathbb{F}_p)$ and contains some point of the
form $\kappa_u$, it must contain $\kappa_t$. In particular,
$\kappa_t$ is a simple point of $\mathcal{C}$. If $\mathcal{C}_t$
is the component of $\mathcal{C}$ passing through $\kappa_t$ and
if $\gamma$ is an element of $\text{PSL}(2, \mathbb{F}_p)$ mapping
$\kappa_t$ to $\kappa_u$, then $\gamma$ must map the component
$\mathcal{C}_t$ onto the component $\mathcal{C}_u$ passing through
$\kappa_u$. It follows that $\text{PSL}(2, \mathbb{F}_p)$ acts
transitively on the set of irreducible components of $\mathcal{C}$.
However, it is well-known that $\text{PSL}(2, \mathbb{F}_p)$ has no
nontrivial permutation representation of degree $<p$. Therefore
$\mathcal{C}$ has only one irreducible component, that is to say,
$\mathcal{C}$ is irreducible.

  In his thesis \cite{Ve}, V\'{e}lu investigated the universal
elliptic curves. He was able to use the family of elliptic curves
over $\mathcal{C}$ to prove that $\mathcal{C}$ is nonsingular.
More precisely, in \cite{Ve}, he constructed for any prime
$n \geq 5$ a pair of schemes $(\mathcal{C}, \mathcal{V})$ over
$\mathbb{Z}[1/n]$ equipped with a morphism
$p: \mathcal{V} \rightarrow \mathcal{C}$ which has the following
properties.

(i) $\mathcal{C}$ is a smooth projective curve over $\mathbb{Z}[1/n]$.

(ii) $\mathcal{V}$ is a generalised elliptic curve over $\mathcal{C}$
in the sense of \cite{DeRa}. More precisely,

(a) $p$ is proper and flat, of finite presentation, of relative
    dimension $1$,

(b) the geometric fibres are either proper, smooth and connected curves
of genus $1$, or $n$-sided N\'{e}ron polygons ($n$-gones). We denote
$\widetilde{\mathcal{V}}$ as the smoothness open of $p$.

(c) there exists a morphism
$$+: \widetilde{\mathcal{V}} \times_{\mathcal{C}} \mathcal{V}
     \rightarrow \mathcal{V}$$
whose restriction to $\widetilde{\mathcal{V}}$ makes it a scheme in
commutative groups and defines an action of the scheme in groups
$\widetilde{\mathcal{V}}$ on $\mathcal{V}$. Moreover, translations
operate by rotation on the components of the singular geometric
fibres of $\mathcal{V}$.

(iii) There is an isomorphism $s$ between the kernel $\mathcal{V}_n$
of the multiplication by $n$ in $\widetilde{\mathcal{V}}$ and the
$\mathbb{Z}[1/n]$-scheme $\mathcal{C} \times \mathbb{Z}/n \mathbb{Z}
\times \mu_n$ and this isomorphism preserves the ``Weil form''.

(iv) If $(S, E, \lambda)$ denotes a generalised elliptic curve over
a basis $S$ on $\mathbb{Z}[1/n]$ possessing the properties (a), (b),
(c) above and if $\lambda$ is an isomorphism between $E_n$ and
$S \times \mathbb{Z}/n \mathbb{Z} \times \mu_n$ preserving the ``Weil
form'' there exists a unique pair of morphisms making the diagram
$$\begin{matrix}
  &E &\longrightarrow & \mathcal{V}\\
  &\downarrow &    & \downarrow\\
  &S &\longrightarrow & \mathcal{C}
\end{matrix}$$
Cartesian compatible with the isomorphisms $\lambda$ and $s$. In
other words, the triplet $(\mathcal{C}, \mathcal{V}, s)$ is universal
for the properties listed above.

  The scheme $\mathcal{C}$ is a closed sub-scheme of $\mathbb{P}^{n-1}$
while the scheme $\mathcal{V}$ is a closed sub-scheme of $\mathcal{C}
\times \mathbb{P}^{n-1}$. They are defined by equations (quartic for
$\mathcal{C}$, quadratic for $\mathcal{V}$). Hence, V\'{e}lu gave an
interpretation of Klein's results in \cite{K4} and \cite{K5} and reproved
them.

\textbf{Theorem 2.4.10. (Main Theorem 1).} {\it There is the following
decomposition as}
$\text{Gal}(\overline{\mathbb{Q}}/\mathbb{Q}) \times \text{PSL}(2,
\mathbb{F}_p)$-{\it representations, that is, there is an isomorphism}
$$V_{I(\mathcal{L}(X(p)))}=\bigoplus_{j} m_{p, j} V_{p, j},\eqno{(2.4.15)}$$
{\it where $j$ denotes the $j$-dimensional irreducible representation
of} $\text{PSL}(2, \mathbb{F}_p)$ {\it appearing in the decomposition} (2.4.15),
{\it and $m_{p, j}$ denotes its multiplicities. Correspondingly, the
defining ideal $I(\mathcal{L}(X(p)))$ has the following decomposition as
the intersection of invariant ideals under the action of}
$\text{Gal}(\overline{\mathbb{Q}}/\mathbb{Q}) \times \text{PSL}(2, \mathbb{F}_p)$:
$$I(\mathcal{L}(X(p)))=\bigcap_{j} I_{p, j},\eqno{(2.4.16)}$$
{\it where the ideals $I_{p, j}$ correspond to the representations $V_{p, j}$.
This decomposition} (2.4.15) {\it gives rise to the following two kinds of
representations:}
$$\pi_{p, j}: \text{PSL}(2, \mathbb{F}_p) \longrightarrow \text{Aut}(V_{p, j})
  \eqno{(2.4.17)}$$
{\it is an irreducible $\mathbb{Q}(\zeta_p)$-rational representation
of} $\text{PSL}(2, \mathbb{F}_p)$, {\it and}
$$\rho_{p, j}: \text{Gal}(\overline{\mathbb{Q}}/\mathbb{Q}) \longrightarrow
  \text{Aut}(V_{p, j})\eqno{(2.4.18)}$$
{\it is an irreducible $\mathbb{Q}(\zeta_p)$-rational representation
of} $\text{Gal}(\overline{\mathbb{Q}}/\mathbb{Q})$. {\it On the other
hand},
$$\pi_p: \text{PSL}(2, \mathbb{F}_p) \longrightarrow \text{Aut}(\mathcal{L}(X(p)))
  \eqno{(2.4.19)}$$
{\it is a reducible $\mathbb{Q}(\zeta_p)$-rational representation of}
$\text{PSL}(2, \mathbb{F}_p)$, {\it and}
$$\rho_p: \text{Gal}(\overline{\mathbb{Q}}/\mathbb{Q}) \longrightarrow
  \text{Aut}(\mathcal{L}(X(p)))\eqno{(2.4.20)}$$
{\it is a reducible $\mathbb{Q}(\zeta_p)$-rational representation of}
$\text{Gal}(\overline{\mathbb{Q}}/\mathbb{Q})$. {\it Locally},
$\pi_{p, j}$ {\it mat-ches} $\rho_{p, j}$ {\it by the following
one-to-one correspondence}
$$\pi_{p, j} \longleftrightarrow \rho_{p, j}.\eqno{(2.4.21)}$$
{\it Globally}, $\pi_p$ {\it matches} $\rho_p$ {\it by the following
one-to-one correspondence}
$$\pi_p \longleftrightarrow \rho_p.\eqno{(2.4.22)}$$
{\it Local-to-globally, the set $\{ \pi_{p, j} \}$ correspond to the
set $\{ \rho_{p, j} \}$, where $\{ \pi_{p, j} \}$ are joint irreducible
$\mathbb{Q}(\zeta_p)$-rational representations of}
$\text{PSL}(2, \mathbb{F}_p)$, {\it and} $\{ \rho_{p, j} \}$ {\it are
joint irreducible} $\mathbb{Q}(\zeta_p)$-{\it rational representations
of} $\text{Gal}(\overline{\mathbb{Q}}/\mathbb{Q})$. {\it In particular,
the Galois representations $\rho_{p, j}$ are realized by the algebraic
cycles corresponding to the ideals $I_{p, j}$.}

{\it Proof}. By the above arguments, we give three kinds of proofs (one
comes from finite group theory, the other comes from algebraic geometry,
and the third comes from V\'{e}lu \cite{Ve}) that
$$\text{PSL}(2, \mathbb{F}_p)=\text{Aut}(\mathcal{L}(X(p))).\eqno{(2.4.23)}$$
Hence, $\text{Aut}(\mathcal{L}(X(p)))$ gives rise to a reducible
representation of $\text{PSL}(2, \mathbb{F}_p)$:
$$\pi_p: \text{PSL}(2, \mathbb{F}_p) \longrightarrow \text{Aut}(\mathcal{L}(X(p))).
  \eqno{(2.4.24)}$$
It must be a direct sum of the irreducible representations of
$\text{PSL}(2, \mathbb{F}_p)$. Each irreducible component consists of a
system of quartic equations, which are invariant under the action of
$\text{PSL}(2, \mathbb{F}_p)$. Hence, we obtain a decomposition
$$V_{I(\mathcal{L}(X(p)))}=\bigoplus_{j} m_{p, j} V_{p, j},\eqno{(2.4.25)}$$
where $V_{p, j}$ denotes the $j$-dimensional irreducible representation
of the group $\text{PSL}(2, \mathbb{F}_p)$. Correspondingly, in the context
of algebraic geometry, it can be expressed as a intersection of ideals
which corresponds to the decomposition (2.4.25):
$$I(\mathcal{L}(X(p)))=\bigcap_{j} I_{p, j},\eqno{(2.4.26)}$$
where the ideals $I_{p, j}$ correspond to the representations $V_{p, j}$,
and these ideals are generated by a system of quartic equations, which are
linear combinations of the forms given by (2.4.11). According to the
formulas (2.4.5) and (2.4.6), the smallest irreducible representation of
$\text{PSL}(2, \mathbb{F}_p)$ has dimension $\frac{p-1}{2}$, which is
defined over $\mathbb{Q}(\zeta_p)$. Thus, as a representation,
$\text{Aut}(\mathcal{L}(X(p)))$ can be realized over $\mathbb{Q}(\zeta_p)$
by this smallest irreducible representation on the system of quartic
equations given by (2.4.11). This gives rise to the reducible
$\mathbb{Q}(\zeta_p)$-rational representation
$$\pi_p: \text{PSL}(2, \mathbb{F}_p) \longrightarrow \text{Aut}(\mathcal{L}(X(p)))$$
as well as its irreducible $\mathbb{Q}(\zeta_p)$-rational components
$$\pi_{p, j}: \text{PSL}(2, \mathbb{F}_p) \longrightarrow \text{Aut}(V_{p, j}).$$

  On the other hand, (2.4.23) also gives rise to a Galois representation
$$\rho_p: \text{Gal}(\overline{\mathbb{Q}}/\mathbb{Q}) \rightarrow
  \text{Aut}(\mathcal{L}(X(p))),\eqno{(2.4.27)}$$
since $\mathcal{L}(X(p))$ is defined over $\mathbb{Q}$. At the same time,
we get the correspondence between the $\mathbb{Q}(\zeta_p)$-rational
representation of $\text{PSL}(2, \mathbb{F}_p)$ and the
$\mathbb{Q}(\zeta_p)$-rational Galois representation:
$$\pi_p \longleftrightarrow \rho_p.$$
According to the above decomposition of the $\mathbb{Q}(\zeta_p)$-rational
representation of $\text{PSL}(2, \mathbb{F}_p)$ into irreducible components,
we obtain a corresponding decomposition of the Galois representation (2.4.27)
into the following:
$$\rho_{p, j}: \text{Gal}(\overline{\mathbb{Q}}/\mathbb{Q}) \longrightarrow
  \text{Aut}(V_{p, j}),$$
as well as the correspondence between the irreducible $\mathbb{Q}(\zeta_p)$-rational
representation of $\text{PSL}(2, \mathbb{F}_p)$ and the irreducible
$\mathbb{Q}(\zeta_p)$-rational Galois representation:
$$\pi_{p, j} \longleftrightarrow \rho_{p, j}.$$
This completes the proof of our Main Theorem 1.

\noindent
$\qquad \qquad \qquad \qquad \qquad \qquad \qquad \qquad \qquad
 \qquad \qquad \qquad \qquad \qquad \qquad \qquad \boxed{}$

  Now, we will give the proof of Corollary 1.2 according to $p=7$, $11$
and $13$, respectively.

(1) When $p=7$, we have the following correspondence:
$$\pi_{7} \longleftrightarrow \rho_{7},\eqno{(2.4.28)}$$
which is the same as
$$\pi_{7, 1} \longleftrightarrow \rho_{7, 1}$$
because that the representation
$$\pi_7: \text{PSL}(2, \mathbb{F}_7) \longrightarrow \text{Aut}(\mathcal{L}(X(7)))$$
is irreducible, i.e., the local part is equal to the global part.
Here both the trivial representation
$$\pi_{7, 1}: \text{PSL}(2, \mathbb{F}_7) \longrightarrow \text{Aut}(V_{7, 1}),$$
and the corresponding Galois representation
$$\rho_{7, 1}: \text{Gal}(\overline{\mathbb{Q}}/\mathbb{Q})
  \longrightarrow \text{Aut}(V_{7, 1})$$
have the multiplicity-one properties, that is,
$$m_{7, 1}=1.$$
In his paper \cite{K2}, Klein studied the invariant theory for
$\text{PSL}(2, \mathbb{F}_7)$ by its three-dimensional irreducible complex
representation and gave an explicit realization of this representation.
Recall that the three-dimensional representation of the finite group
$\text{PSL}(2, \mathbb{F}_7)$ of order $168$, which acts on the projective
plane (in Klein's notation)
$$\mathbb{CP}^2=\{ (\lambda, \mu, \nu): \lambda, \mu, \nu \in \mathbb{C} \}.$$
This representation is defined over the cyclotomic field
$\mathbb{Q}(\gamma)$ with $\gamma=e^{\frac{2 \pi i}{7}}$. Put
$$S=\frac{1}{\sqrt{-7}} \left(\begin{matrix}
  \gamma^5-\gamma^2 & \gamma^3-\gamma^4 & \gamma^6-\gamma\\
  \gamma^3-\gamma^4 & \gamma^6-\gamma & \gamma^5-\gamma^2\\
  \gamma^6-\gamma & \gamma^5-\gamma^2 & \gamma^3-\gamma^4
  \end{matrix}\right), \quad
 T=\left(\begin{matrix}
    \gamma &          &    \\
           & \gamma^4 &    \\
           &          & \gamma^2
   \end{matrix}\right).\eqno{(2.4.29)}$$
We have
$$S^2=T^7=(ST)^3=I.$$
Let $G=\langle S, T \rangle$, then
$$G \cong \text{PSL}(2, \mathbb{F}_7).$$
The algebraic cycle corresponding to the ideal $I_{7, 1}$ is given by the
Klein's quartic curve in $\mathbb{P}^2$:
$$\lambda^3 \mu+\mu^3 \nu+\nu^3 \lambda=0,\eqno{(2.4.30)}$$
which is invariant under the action of the group $\text{PSL}(2, \mathbb{F}_7)$.
Hence, the action of $\text{PSL}(2, \mathbb{F}_7)$ on the Klein quartic curve
gives a realization of the trivial representation of the group
$\text{PSL}(2, \mathbb{F}_7)$. The quotient of the Klein quartic curve by this
group is isomorphic to $\mathbb{P}^1(\mathbb{C})$, and the natural projection
defines a clean Belyi map (see \cite{SV}). We do not attempt to draw the
corresponding dessin and refer to Klein's paper (see \cite{K2}, p. 126), where
one finds a very beautiful triangular tesselation of the fundamental domain of
the Klein quartic on the universal covering, this was the figure from which the
uniformization started. After suitable identifications, we get some triangular
dessin $D$ on a curve $X$ of genus three. Thus the Klein quartic itself has
been drawn; indeed, the automorphism group of $D$ and consequently of the
curve $X$ has $168$ elements, and the Klein quartic is the only curve of
genus three with this highest possible number of automorphisms.

(2) When $p=11$, we have the following correspondence:
$$\pi_{11} \longleftrightarrow \rho_{11},\eqno{(2.4.31)}$$
which is the same as
$$\pi_{11, 10} \longleftrightarrow \rho_{11, 10},$$
because that the representation
$$\pi_{11}: \text{PSL}(2, \mathbb{F}_{11}) \longrightarrow
  \text{Aut}(\mathcal{L}(X(11)))$$
is irreducible, i.e., the local part is equal to the global part.
Here both the discrete series representation
$$\pi_{11, 10}: \text{PSL}(2, \mathbb{F}_{11})
  \longrightarrow \text{Aut}(V_{11, 10}),$$
and the corresponding Galois representation
$$\rho_{11, 10}: \text{Gal}(\overline{\mathbb{Q}}/\mathbb{Q})
  \longrightarrow \text{Aut}(V_{11, 10})$$
have the multiplicity-one properties, that is,
$$m_{11, 10}=1.$$
In his paper \cite{K3}, Klein studied the invariant theory for
$\text{PSL}(2, \mathbb{F}_{11})$ by its five-dimensional irreducible
complex representation and gave an explicit realization of this
representation. Recall that the five-dimensional representation of the
finite group $\text{PSL}(2, \mathbb{F}_{11})$ of order $660$, which
acts on the four-dimensional projective space (in Klein's notation)
$$\mathbb{CP}^4=\{ (y_1, y_4, y_5, y_9, y_3): y_i \in \mathbb{C}
  \quad (i=1, 4, 5, 9, 3) \}.$$
This representation is defined over the cyclotomic field
$\mathbb{Q}(\rho)$ with $\rho=e^{\frac{2 \pi i}{11}}$. Put
$$S=\frac{1}{\sqrt{-11}} \begin{pmatrix}
  \rho^9-\rho^2 & \rho^4-\rho^7 & \rho^3-\rho^8 & \rho^5-\rho^6 & \rho-\rho^{10}\\
  \rho^4-\rho^7 & \rho^3-\rho^8 & \rho^5-\rho^6 & \rho-\rho^{10} & \rho^9-\rho^2\\
  \rho^3-\rho^8 & \rho^5-\rho^6 & \rho-\rho^{10} & \rho^9-\rho^2 & \rho^4-\rho^7\\
  \rho^5-\rho^6 & \rho-\rho^{10} & \rho^9-\rho^2 & \rho^4-\rho^7 & \rho^3-\rho^8\\
  \rho-\rho^{10} & \rho^9-\rho^2 & \rho^4-\rho^7 & \rho^3-\rho^8 & \rho^5-\rho^6
\end{pmatrix}\eqno{(2.4.32)}$$
and
$$T=\text{diag}(\rho, \rho^4, \rho^5, \rho^9, \rho^3).\eqno{(2.4.33)}$$
We have
$$S^2=T^{11}=(ST)^3=I.$$
Let $G=\langle S, T \rangle$, then
$$G \cong \text{PSL}(2, \mathbb{F}_{11}).$$
The algebraic cycle corresponding to the ideal $I_{11, 10}$ is given
by the system of quartic equations in $\mathbb{P}^4$:
$$\left\{\aligned
 &\mathbb{B}_1^{(1)}=0, \quad \mathbb{B}_4^{(1)}=0, \quad
  \mathbb{B}_5^{(1)}=0, \quad \mathbb{B}_9^{(1)}=0, \quad
  \mathbb{B}_3^{(1)}=0,\\
 &\mathbb{B}_1^{(2)}=0, \quad \mathbb{B}_4^{(2)}=0, \quad
  \mathbb{B}_5^{(2)}=0, \quad \mathbb{B}_9^{(2)}=0, \quad
  \mathbb{B}_3^{(2)}=0,
\endaligned\right.$$
which is invariant under the action of the group $\text{PSL}(2, \mathbb{F}_{11})$.
In particular, the action of $\text{PSL}(2, \mathbb{F}_{11})$
on the basis
$$(\mathbb{B}_1^{(1)}, \mathbb{B}_4^{(1)}, \mathbb{B}_5^{(1)},
   \mathbb{B}_9^{(1)}, \mathbb{B}_3^{(1)}, \mathbb{B}_1^{(2)},
   \mathbb{B}_4^{(2)}, \mathbb{B}_5^{(2)}, \mathbb{B}_9^{(2)},
   \mathbb{B}_3^{(2)})$$
gives a realization of the discrete series representation of
$\text{PSL}(2, \mathbb{F}_{11})$. Here,
$$\left\{\aligned
  \mathbb{B}_1^{(1)} &=y_4^3 y_9+y_9^3 y_5+y_3^3 y_1,\\
  \mathbb{B}_4^{(1)} &=y_5^3 y_3+y_3^3 y_9+y_1^3 y_4,\\
  \mathbb{B}_5^{(1)} &=y_9^3 y_1+y_1^3 y_3+y_4^3 y_5,\\
  \mathbb{B}_9^{(1)} &=y_3^3 y_4+y_4^3 y_1+y_5^3 y_9,\\
  \mathbb{B}_3^{(1)} &=y_1^3 y_5+y_5^3 y_4+y_9^3 y_3.
\endaligned\right.\eqno{(2.4.34)}$$
$$\left\{\aligned
  \mathbb{B}_1^{(2)} &=y_9^2 y_1 y_4-y_3^2 y_1 y_5-y_5^2 y_9 y_4,\\
  \mathbb{B}_4^{(2)} &=y_3^2 y_4 y_5-y_1^2 y_4 y_9-y_9^2 y_3 y_5,\\
  \mathbb{B}_5^{(2)} &=y_1^2 y_5 y_9-y_4^2 y_5 y_3-y_3^2 y_1 y_9,\\
  \mathbb{B}_9^{(2)} &=y_4^2 y_9 y_3-y_5^2 y_9 y_1-y_1^2 y_4 y_3,\\
  \mathbb{B}_3^{(2)} &=y_5^2 y_3 y_1-y_9^2 y_3 y_4-y_4^2 y_5 y_1.
\endaligned\right.\eqno{(2.4.35)}$$
Then
$$\left\{\aligned
  T(\mathbb{B}_1^{(1)}) &=\rho^{10} \mathbb{B}_1^{(1)},\\
  T(\mathbb{B}_4^{(1)}) &=\rho^7 \mathbb{B}_4^{(1)},\\
  T(\mathbb{B}_5^{(1)}) &=\rho^6 \mathbb{B}_5^{(1)},\\
  T(\mathbb{B}_9^{(1)}) &=\rho^2 \mathbb{B}_9^{(1)},\\
  T(\mathbb{B}_3^{(1)}) &=\rho^8 \mathbb{B}_3^{(1)},
\endaligned\right. \quad
  \left\{\aligned
  T(\mathbb{B}_1^{(2)}) &=\rho \mathbb{B}_1^{(2)},\\
  T(\mathbb{B}_4^{(2)}) &=\rho^4 \mathbb{B}_4^{(2)},\\
  T(\mathbb{B}_5^{(2)}) &=\rho^5 \mathbb{B}_5^{(2)},\\
  T(\mathbb{B}_9^{(2)}) &=\rho^9 \mathbb{B}_9^{(2)},\\
  T(\mathbb{B}_3^{(2)}) &=\rho^3 \mathbb{B}_3^{(2)}.
\endaligned\right.\eqno{(2.4.36)}$$
Without loss of generality, we can only study the action of $S$ on
$\mathbb{B}_1^{(1)}$ and $\mathbb{B}_1^{(2)}$, and then by the cyclic
permutation $(14593)$. We find that
$$\aligned
   11 S(\mathbb{B}_1^{(1)})
 &=3 (\rho^9+\rho^2-\rho^5-\rho^6) \mathbb{B}_1^{(2)}+
   3 (\rho^4+\rho^7-\rho-\rho^{10}) \mathbb{B}_4^{(2)}\\
 &+3 (\rho^3+\rho^8-\rho^9-\rho^2) \mathbb{B}_5^{(2)}+
   3 (\rho^5+\rho^6-\rho^4-\rho^7) \mathbb{B}_9^{(2)}\\
 &+3 (\rho+\rho^{10}-\rho^3-\rho^8) \mathbb{B}_3^{(2)}\\
 &+[2+(\rho+\rho^{10})-(\rho^9+\rho^2)-(\rho^5+\rho^6)] \mathbb{B}_1^{(1)}\\
 &+[2+(\rho^9+\rho^2)-(\rho^4+\rho^7)-(\rho+\rho^{10})] \mathbb{B}_4^{(1)}\\
 &+[2+(\rho^4+\rho^7)-(\rho^3+\rho^8)-(\rho^9+\rho^2)] \mathbb{B}_5^{(1)}\\
 &+[2+(\rho^3+\rho^8)-(\rho^5+\rho^6)-(\rho^4+\rho^7)] \mathbb{B}_9^{(1)}\\
 &+[2+(\rho^5+\rho^6)-(\rho+\rho^{10})-(\rho^3+\rho^8)] \mathbb{B}_3^{(1)}.
\endaligned\eqno{(2.4.37)}$$
$$\aligned
   11 S(\mathbb{B}_1^{(2)})
 &=[2+(\rho+\rho^{10})-(\rho^9+\rho^2)-(\rho^5+\rho^6)] \mathbb{B}_1^{(2)}\\
 &+[2+(\rho^9+\rho^2)-(\rho^4+\rho^7)-(\rho+\rho^{10})] \mathbb{B}_4^{(2)}\\
 &+[2+(\rho^4+\rho^7)-(\rho^3+\rho^8)-(\rho^9+\rho^2)] \mathbb{B}_5^{(2)}\\
 &+[2+(\rho^3+\rho^8)-(\rho^5+\rho^6)-(\rho^4+\rho^7)] \mathbb{B}_9^{(2)}\\
 &+[2+(\rho^5+\rho^6)-(\rho+\rho^{10})-(\rho^3+\rho^8)] \mathbb{B}_3^{(2)}\\
 &+(\rho^9+\rho^2-\rho^5-\rho^6) \mathbb{B}_1^{(1)}+(\rho^4+\rho^7-\rho-\rho^{10})
  \mathbb{B}_4^{(1)}\\
 &+(\rho^3+\rho^8-\rho^9-\rho^2) \mathbb{B}_5^{(1)}+(\rho^5+\rho^6-\rho^4-\rho^7)
  \mathbb{B}_9^{(1)}\\
 &+(\rho+\rho^{10}-\rho^3-\rho^8) \mathbb{B}_3^{(1)}.
\endaligned\eqno{(2.4.38)}$$
Therefore,
$$\text{Tr}(\widetilde{S})=2, \quad \text{Tr}(\widetilde{T})=-1,\eqno{(2.4.39)}$$
where $\widetilde{S}$ and $\widetilde{T}$ are the corresponding matrices with
respect to the actions of $S$ and $T$ on the basis
$(\mathbb{B}_1^{(1)}, \mathbb{B}_4^{(1)}, \mathbb{B}_5^{(1)}, \mathbb{B}_9^{(1)},
  \mathbb{B}_3^{(1)}, \mathbb{B}_1^{(2)}, \mathbb{B}_4^{(2)}, \mathbb{B}_5^{(2)},
  \mathbb{B}_9^{(2)}, \mathbb{B}_3^{(2)})$.
This gives rise to a ten-dimensional representation of $\text{PSL}(2, \mathbb{F}_{11})$,
which is a discrete series representation of $\text{PSL}(2, \mathbb{F}_{11})$.
By (2.4.39), this discrete series representation corresponds to the character
$\chi_5$ (see \cite{CC}).

  The notion of dessin d'enfant is also latent in Klein's work (1879)
(see \cite{K3}), and these objects are called there ``Linienzuges'' (English
translation: line train). The drawing was appeared in Figure 3 of Klein's paper (see
\cite{K3}, Fig. 3, p. 144). It is the pre-image of the interval $[0, 1]$ by a degree-$11$
covering of the Riemann sphere ramified over the three points $0$, $1$ and $\infty$.
Klein labels the preimage of $0$ by $\bullet$ and those of $1$ by $+$. He identified
the monodromy group of this covering as the group $\text{PSL}(2, \mathbb{F}_{11})$.
This work is in the lineage of his Lectures on the icosahedron (1884) (see
\cite{K}). We mention by the way that the study of regular polyhedra is part
of Grothendieck's program which is titled regular polyhedra over finite fields
(see \cite{Groth}).

(3) When $p=13$, we have the following correspondence:
$$\text{globally:} \quad \pi_{13} \longleftrightarrow \rho_{13},\eqno{(2.4.40)}$$
$$\text{local-to-globally:} \quad \{ \pi_{13, 1}, \pi_{13, 7}, \pi_{13, 13} \}
  \longleftrightarrow \{\rho_{13, 1}, \rho_{13, 7}, \rho_{13, 13} \},\eqno{(2.4.41)}$$
where
$$\text{locally:} \quad \pi_{13, 1} \longleftrightarrow \rho_{13, 1}, \quad
  \pi_{13, 7} \longleftrightarrow \rho_{13, 7}, \quad
  \pi_{13, 13} \longleftrightarrow \rho_{13, 13}.\eqno{(2.4.42)}$$
Here,
$$\pi_{13}: \text{PSL}(2, \mathbb{F}_{13}) \longrightarrow
  \text{Aut}(\mathcal{L}(X(13))),$$
$$\left\{\aligned
 &\pi_{13, 1}: \text{PSL}(2, \mathbb{F}_{13})
   \longrightarrow \text{Aut}(V_{13, 1}),\\
 &\pi_{13, 7}: \text{PSL}(2, \mathbb{F}_{13})
  \longrightarrow \text{Aut}(V_{13, 7}),\\
 &\pi_{13, 13}: \text{PSL}(2, \mathbb{F}_{13})
  \longrightarrow \text{Aut}(V_{13, 13}),
\endaligned\right.$$
are the trivial representation, the degenerate principal series
representation and the Steinberg representation of $\text{PSL}(2, \mathbb{F}_{13})$,
respectively, and
$$\rho_{13}: \text{Gal}(\overline{\mathbb{Q}}/\mathbb{Q})
  \longrightarrow \text{Aut}(\mathcal{L}(X(13))),$$
$$\left\{\aligned
 &\rho_{13, 1}: \text{Gal}(\overline{\mathbb{Q}}/\mathbb{Q})
  \longrightarrow \text{Aut}(V_{13, 1}),\\
 &\rho_{13, 7}: \text{Gal}(\overline{\mathbb{Q}}/\mathbb{Q})
  \longrightarrow \text{Aut}(V_{13, 7}),\\
 &\rho_{13, 13}: \text{Gal}(\overline{\mathbb{Q}}/\mathbb{Q})
  \longrightarrow \text{Aut}(V_{13, 13}).
\endaligned\right.$$
are the corresponding Galois representations. They all have
the multiplicity-one properties, that is,
$$m_{13, 1}=m_{13, 7}=m_{13, 13}=1.$$
Let us begin with the invariant theory for $\text{SL}(2, \mathbb{F}_{13})$,
which we developed in \cite{Y1}, \cite{Y2}, \cite{Y3}, \cite{Y4}, \cite{Y5}
and \cite{Y6}. The representation of $\text{SL}(2, \mathbb{F}_{13})$ we will
consider is the unique six-dimensional irreducible complex representation for
which the eigenvalues of $\left(\begin{matrix} 1 & 1\\ 0 & 1 \end{matrix}\right)$
are the $\exp (\underline{a} . 2 \pi i/13)$ for $\underline{a}$ a non-square
mod $13$. We will give an explicit realization of this representation. This
explicit realization will play a major role for giving a complete system of
invariants associated to $\text{SL}(2, \mathbb{F}_{13})$. Recall that the
six-dimensional representation of the finite group
$\text{SL}(2, \mathbb{F}_{13})$ of order $2184$, which acts on the
five-dimensional projective space
$$\mathbb{CP}^5=\{ (z_1, z_2, z_3, z_4, z_5, z_6): z_i \in \mathbb{C}
  \quad (i=1, 2, 3, 4, 5, 6) \}.$$
This representation is defined over the cyclotomic field
$\mathbb{Q}(e^{\frac{2 \pi i}{13}})$. Put
$$S=-\frac{1}{\sqrt{13}} \begin{pmatrix}
  \zeta^{12}-\zeta & \zeta^{10}-\zeta^3 & \zeta^4-\zeta^9
& \zeta^5-\zeta^8 & \zeta^2-\zeta^{11} & \zeta^6-\zeta^7\\
  \zeta^{10}-\zeta^3 & \zeta^4-\zeta^9 & \zeta^{12}-\zeta
& \zeta^2-\zeta^{11} & \zeta^6-\zeta^7 & \zeta^5-\zeta^8\\
  \zeta^4-\zeta^9 & \zeta^{12}-\zeta & \zeta^{10}-\zeta^3
& \zeta^6-\zeta^7 & \zeta^5-\zeta^8 & \zeta^2-\zeta^{11}\\
  \zeta^5-\zeta^8 & \zeta^2-\zeta^{11} & \zeta^6-\zeta^7
& \zeta-\zeta^{12} & \zeta^3-\zeta^{10} & \zeta^9-\zeta^4\\
  \zeta^2-\zeta^{11} & \zeta^6-\zeta^7 & \zeta^5-\zeta^8
& \zeta^3-\zeta^{10} & \zeta^9-\zeta^4 & \zeta-\zeta^{12}\\
  \zeta^6-\zeta^7 & \zeta^5-\zeta^8 & \zeta^2-\zeta^{11}
& \zeta^9-\zeta^4 & \zeta-\zeta^{12} & \zeta^3-\zeta^{10}
\end{pmatrix}\eqno{(2.4.43)}$$
and
$$T=\text{diag}(\zeta^7, \zeta^{11}, \zeta^8, \zeta^6, \zeta^2, \zeta^5)$$
where $\zeta=\exp(2 \pi i/13)$. We have
$$S^2=-I, \quad T^{13}=(ST)^3=I.$$
In \cite{Y1}, we put $P=S T^{-1} S$ and $Q=S T^3$. Then $(Q^3 P^4)^3=-I$ (see
\cite{Y1}, the proof of Theorem 3.1). Let $G=\langle S, T \rangle$, then
$$G \cong \text{SL}(2, \mathbb{F}_{13}).$$

  The algebraic cycle corresponding to the ideal $I(\mathcal{L}(X(13)))$
is given by the system of quartic equations in $\mathbb{P}^5$:
$$I_{13}=\langle \mathbf{B}_0^{(i)}, \mathbf{B}_1^{(j)}, \mathbf{B}_3^{(j)},
    \mathbf{B}_9^{(j)}, \mathbf{B}_{12}^{(j)}, \mathbf{B}_{10}^{(j)},
    \mathbf{B}_4^{(j)}, \mathbf{B}_5, \mathbf{B}_2, \mathbf{B}_6,
    \mathbf{B}_8, \mathbf{B}_{11}, \mathbf{B}_7 \rangle$$
with $i=0. 1, 2$ and $j=1, 2$, where
$$\left\{\aligned
  \mathbf{B}_0^{(0)} &=z_1 z_2 z_4 z_5+z_2 z_3 z_5 z_6+z_3 z_1 z_6 z_4,\\
  \mathbf{B}_0^{(1)} &=z_1 z_5^3+z_2 z_6^3+z_3 z_4^3,\\
  \mathbf{B}_0^{(2)} &=z_1^3 z_6+z_2^3 z_4+z_3^3 z_5,
\endaligned\right.\eqno{(2.4.44)}$$
$$\left\{\aligned
  \mathbf{B}_1^{(1)} &=z_3 z_5^3+z_1^3 z_4-z_1 z_2^3,\\
  \mathbf{B}_1^{(2)} &=z_2 z_4 z_6^2-z_3^2 z_6 z_4-z_1^2 z_2 z_5,\\
  \mathbf{B}_3^{(1)} &=z_2 z_4^3+z_3^3 z_6-z_3 z_1^3,\\
  \mathbf{B}_3^{(2)} &=z_1 z_6 z_5^2-z_2^2 z_5 z_6-z_3^2 z_1 z_4,\\
  \mathbf{B}_9^{(1)} &=z_1 z_6^3+z_2^3 z_5-z_2 z_3^3,\\
  \mathbf{B}_9^{(2)} &=z_3 z_5 z_4^2-z_1^2 z_4 z_5-z_2^2 z_3 z_6,\\
  \mathbf{B}_{12}^{(1)} &=z_1 z_4^3+z_2^3 z_6+z_4 z_5^3,\\
  \mathbf{B}_{12}^{(2)} &=z_2 z_5 z_4^2-z_3^2 z_1 z_5-z_6^2 z_3 z_1,\\
  \mathbf{B}_{10}^{(1)} &=z_3 z_6^3+z_1^3 z_5+z_6 z_4^3,\\
  \mathbf{B}_{10}^{(2)} &=z_1 z_4 z_6^2-z_2^2 z_3 z_4-z_5^2 z_2 z_3,\\
  \mathbf{B}_4^{(1)} &=z_2 z_5^3+z_3^3 z_4+z_5 z_6^3,\\
  \mathbf{B}_4^{(2)} &=z_3 z_6 z_5^2-z_1^2 z_2 z_6-z_4^2 z_1 z_2,
\endaligned\right.\eqno{(2.4.45)}$$
and
$$\left\{\aligned
  \mathbf{B}_5 &=-z_2^2 z_1 z_5+z_4 z_5 z_6^2+z_2 z_3 z_4^2,\\
  \mathbf{B}_2 &=-z_1^2 z_3 z_4+z_6 z_4 z_5^2+z_1 z_2 z_6^2,\\
  \mathbf{B}_6 &=-z_3^2 z_2 z_6+z_5 z_6 z_4^2+z_3 z_1 z_5^2,\\
  \mathbf{B}_8 &=z_2 z_4 z_5^2+z_1 z_2 z_3^2+z_1^2 z_5 z_6,\\
  \mathbf{B}_{11} &=z_1 z_6 z_4^2+z_3 z_1 z_2^2+z_3^2 z_4 z_5,\\
  \mathbf{B}_7 &=z_3 z_5 z_6^2+z_2 z_3 z_1^2+z_2^2 z_6 z_4.
\endaligned\right.\eqno{(2.4.46)}$$

  In order to explain what the above system of quartic equations come from,
we give the geometry of the modular curve via the fundamental relation (see
\cite{AR}). For every element $t$ of $K(\delta)=\mathbb{Z}/p \mathbb{Z}$,
let $E_t$ denote the linear form on the space $V^{-}$ of odd functions
on $K(\delta)$ given by $E_t(h)=h(t)$. If $w$, $x$, $y$, $z$ are
any elements of $K(\delta)$, denote by $\Phi_{w, x, y, z}$ the
quartic form given by
$$\aligned
  \Phi_{w, x, y, z} &=E_{w+x} \cdot E_{w-x} \cdot E_{y+z} \cdot E_{y-z}+\\
                    &+E_{w+y} \cdot E_{w-y} \cdot E_{z+x} \cdot E_{z-x}+\\
                    &+E_{w+z} \cdot E_{w-z} \cdot E_{x+y} \cdot E_{x-y}.
\endaligned\eqno{(2.4.47)}$$
Then $\mathcal{L}$ is defined by the equations
$$\Phi_{w, x, y, z}=0$$
with $w$, $x$, $y$, $z$ in $\mathbb{Z}/p \mathbb{Z}$. Since $E_{-t}=-E_t$ for
every $t \in K(\delta)$, we have $\Phi_{-w, x, y, z}=\Phi_{w, x, y, z}$.
Furthermore, we have
$$\Phi_{x, w, y, z}=\Phi_{x, y, z, w}=-\Phi_{w, x, y, z}.\eqno{(2.4.48)}$$
Since the odd permutations
$$\left(\begin{matrix}
  1 & 2 & 3 & 4\\
  2 & 1 & 3 & 4
  \end{matrix}\right) \quad \text{and} \quad
  \left(\begin{matrix}
  1 & 2 & 3 & 4\\
  2 & 3 & 4 & 1
  \end{matrix}\right)$$
generate the group of all permutations on four objects, it follows that
if $(a, b, c, d)=\sigma(w, x, y, z)$ is a permutation of $(w, x, y, z)$
then
$$\Phi_{a, b, c, d}=(-1)^{\sigma} \Phi_{w, x, y, z},\eqno{(2.4.49)}$$
where $(-1)^{\sigma}$ denotes the sign of the permutation $\sigma$. It
follows from this and from the sentence preceding (2.4.48) that for all
choices of signs we have
$$\Phi_{\pm w, \pm x, \pm y, \pm z}=\Phi_{w, x, y, z}.\eqno{(2.4.50)}$$
It follows that each quartic $\Phi_{w, x, y, z}$ is equal, up to a sign,
to a quartic $\Phi_{a, b, c, d}$ with
$$0 \leq a < b < c < d \leq \frac{p-1}{2}.$$
The locus $\mathcal{L}$ is therefore defined by
$\left(\begin{matrix} m\\ 4 \end{matrix}\right)$ quartics, where $2m-1=p$.
In general, these quartics are not distinct.

  For example, when $p=13$, the cardinality is $35$ but there are only
$21$ distinct quartics, namely:

(1) $\Phi_{0123}=E_1^3 E_5-E_2^3 E_4+E_3^3 E_1$.

(2) $\Phi_{0124}=\Phi_{3456}=E_1^2 E_2 E_6-E_2^2 E_3 E_5+E_4^2 E_1 E_3$.

(3) $\Phi_{0125}=\Phi_{1234}=-E_1^2 E_3 E_6-E_2^2 E_4 E_6+E_5^2 E_1 E_3$.

(4) $\Phi_{0126}=\Phi_{2345}=-E_1^2 E_4 E_5+E_2^2 E_5 E_6+E_6^2 E_1 E_3$.

(5) $\Phi_{0134}=-E_1^3 E_6-E_3^3 E_5+E_4^3 E_2$.

(6) $\Phi_{0135}=\Phi_{2356}=-E_1^2 E_2 E_5-E_3^2 E_4 E_6+E_5^2 E_2 E_4$.

(7) $\Phi_{0136}=\Phi_{1245}=-E_1^3 E_3 E_4+E_3^2 E_5 E_6+E_6^2 E_2 E_4$.

(8) $\Phi_{0145}=-E_1^3 E_4-E_4^3 E_6+E_5^3 E_3$.

(9) $\Phi_{0146}=\Phi_{1256}=-E_1^2 E_2 E_3+E_4^2 E_5 E_6+E_6^2 E_3 E_5$.

(10) $\Phi_{0156}=-E_1^3 E_2+E_5^3 E_6+E_6^3 E_4$.

(11) $\Phi_{0234}=\Phi_{2456}=-E_2^2 E_1 E_6-E_3^2 E_2 E_6+E_4^2 E_1 E_5$.

(12) $\Phi_{0235}=-E_2^3 E_5+E_3^3 E_6+E_5^3 E_1$.

(13) $\Phi_{0236}=\Phi_{1346}=-E_2^2 E_3 E_4+E_3^2 E_4 E_5+E_6^2 E_1 E_5$.

(14) $\Phi_{0245}=\Phi_{1356}=-E_2^2 E_1 E_4+E_4^2 E_3 E_6+E_5^2 E_2 E_6$.

(15) $\Phi_{0246}=-E_2^3 E_3+E_4^3 E_5+E_6^3 E_2$.

(16) $\Phi_{0256}=-E_2^3 E_1+E_5^3 E_4-E_6^3 E_3$.

(17) $\Phi_{0345}=\Phi_{1236}=-E_3^2 E_1 E_4+E_4^2 E_2 E_5-E_5^2 E_1 E_6$.

(18) $\Phi_{0346}=-E_3^3 E_2+E_4^3 E_3-E_6^3 E_1$.

(19) $\Phi_{0356}=\Phi_{1246}=-E_3^2 E_1 E_2+E_5^2 E_3 E_4-E_6^2 E_2 E_5$.

(20) $\Phi_{0456}=\Phi_{1345}=-E_4^2 E_1 E_2+E_5^2 E_2 E_3-E_6^2 E_1 E_4$.

(21) $\Phi_{1235}=\Phi_{1456}=\Phi_{2346}=-E_1 E_2 E_3 E_5+E_2 E_3 E_4 E_6+E_1 E_4 E_5 E_6$.

  This proves the following:

\textbf{Theorem 2.4.11.} {\it The modular curve $X(13)$ is isomorphic to the
above $21$ quartic equations ($\Phi_{abcd}$-terms).}

  Now, we will give a connection between the above quartic $21$ equations
with the invariant theory for the group $G$. Let $ST^{\nu}$ act on the following
quartic term: $z_1 z_2 z_4 z_5+z_2 z_3 z_5 z_6+z_3 z_1 z_6 z_4$. We obtain
$$\aligned
  &13 [ST^{\nu}(z_1) \cdot ST^{\nu}(z_4) \cdot ST^{\nu}(z_2) \cdot ST^{\nu}(z_5)+\\
  &+ST^{\nu}(z_2) \cdot ST^{\nu}(z_5) \cdot ST^{\nu}(z_3) \cdot ST^{\nu}(z_6)+\\
  &+ST^{\nu}(z_3) \cdot ST^{\nu}(z_6) \cdot ST^{\nu}(z_1) \cdot ST^{\nu}(z_4)]\\
 =&\mathbf{B}_0+\zeta^{\nu} \mathbf{B}_1+\zeta^{3 \nu} \mathbf{B}_3
  +\zeta^{9 \nu} \mathbf{B}_9+\zeta^{12 \nu} \mathbf{B}_{12}+\zeta^{10 \nu}
   \mathbf{B}_{10}+\zeta^{4 \nu} \mathbf{B}_4+\\
  &+2 \zeta^{5 \nu} \mathbf{B}_5+2 \zeta^{2 \nu} \mathbf{B}_2+2 \zeta^{6 \nu}
   \mathbf{B}_6+2 \zeta^{8 \nu} \mathbf{B}_8+2 \zeta^{11 \nu} \mathbf{B}_{11}
   +2 \zeta^{7 \nu} \mathbf{B}_7,
\endaligned\eqno{(2.4.51)}$$
where
$$\left\{\aligned
  \mathbf{B}_0 &=5(z_1 z_2 z_4 z_5+z_2 z_3 z_5 z_6+z_3 z_1 z_6 z_4)
                +2(z_1 z_5^3+z_2 z_6^3+z_3 z_4^3)+\\
               &\quad -2(z_1^3 z_6+z_2^3 z_4+z_3^3 z_5),\\
  \mathbf{B}_1 &=z_1^3 z_4-z_1 z_2^3+z_3 z_5^3-z_3^2 z_6 z_4
                +z_2 z_4 z_6^2-z_1^2 z_2 z_5,\\
  \mathbf{B}_3 &=z_3^3 z_6-z_3 z_1^3+z_2 z_4^3-z_2^2 z_5 z_6
                +z_1 z_6 z_5^2-z_3^2 z_1 z_4,\\
  \mathbf{B}_9 &=z_2^3 z_5-z_2 z_3^3+z_1 z_6^3-z_1^2 z_4 z_5
                +z_3 z_5 z_4^2-z_2^2 z_3 z_6,\\
  \mathbf{B}_{12} &=-z_1 z_4^3-z_4 z_5^3-z_2^3 z_6-z_3^2 z_1 z_5
                   -z_3 z_1 z_6^2+z_4^2 z_2 z_5,\\
  \mathbf{B}_{10} &=-z_3 z_6^3-z_6 z_4^3-z_1^3 z_5-z_2^2 z_3 z_4
                   -z_2 z_3 z_5^2+z_6^2 z_1 z_4,\\
  \mathbf{B}_4 &=-z_2 z_5^3-z_5 z_6^3-z_3^3 z_4-z_1^2 z_2 z_6
                -z_1 z_2 z_4^2+z_5^2 z_3 z_6,\\
\endaligned\right.\eqno{(2.4.52)}$$
and $\mathbf{B}_5$, $\mathbf{B}_2$, $\mathbf{B}_6$, $\mathbf{B}_8$,
$\mathbf{B}_{11}$ and $\mathbf{B}_7$ are given by (2.4.46). This leads
us to define the following decompositions:
$$\mathbf{B}_0=5 \mathbf{B}_0^{(0)}+2 \mathbf{B}_0^{(1)}-2 \mathbf{B}_0^{(2)},
  \eqno{(2.4.53)}$$
$$\mathbf{B}_1=\mathbf{B}_1^{(1)}+\mathbf{B}_1^{(2)}, \quad
  \mathbf{B}_3=\mathbf{B}_3^{(1)}+\mathbf{B}_3^{(2)}, \quad
  \mathbf{B}_9=\mathbf{B}_9^{(1)}+\mathbf{B}_9^{(2)}, \eqno{(2.4.54)}$$
$$\mathbf{B}_{12}=-\mathbf{B}_{12}^{(1)}+\mathbf{B}_{12}^{(2)}, \quad
  \mathbf{B}_{10}=-\mathbf{B}_{10}^{(1)}+\mathbf{B}_{10}^{(2)}, \quad
  \mathbf{B}_4=-\mathbf{B}_4^{(1)}+\mathbf{B}_4^{(2)}, \eqno{(2.4.55)}$$
where $\mathbf{B}_0^{(0)}$, $\mathbf{B}_0^{(1)}$ and $\mathbf{B}_0^{(2)}$
are given by (2.4.44), $\mathbf{B}_1^{(1)}$, $\mathbf{B}_1^{(2)}$,
$\mathbf{B}_3^{(1)}$, $\mathbf{B}_3^{(2)}$, $\mathbf{B}_9^{(1)}$,
$\mathbf{B}_9^{(2)}$, $\mathbf{B}_{12}^{(1)}$, $\mathbf{B}_{12}^{(2)}$,
$\mathbf{B}_{10}^{(1)}$, $\mathbf{B}_{10}^{(2)}$, $\mathbf{B}_4^{(1)}$,
and $\mathbf{B}_4^{(2)}$ are given by (2.4.45). The significance of
these quartic polynomials come from the following:

\textbf{Theorem 2.4.12.} {\it There is a one-to-one correspondence
between the above quartic equations $($$\Phi_{abcd}$-terms$)$
and the quartic polynomials $($$\mathbf{B}$-terms$)$.}

{\it Proof}. Under the following correspondence
$$z_1 \longleftrightarrow -e^{\frac{\pi i}{4}} E_1, \quad
  z_2 \longleftrightarrow -e^{\frac{\pi i}{4}} E_3, \quad
  z_3 \longleftrightarrow e^{\frac{\pi i}{4}} E_4,\eqno{(2.4.56)}$$
and
$$z_4 \longleftrightarrow e^{\frac{\pi i}{4}} E_5, \quad
  z_5 \longleftrightarrow e^{\frac{\pi i}{4}} E_2, \quad
  z_6 \longleftrightarrow e^{\frac{\pi i}{4}} E_6,\eqno{(2.4.57)}$$
we have

(1) $\mathbf{B}_0^{(0)} \longleftrightarrow \Phi_{1235}=\Phi_{1456}=\Phi_{2346}$.

(2) $\mathbf{B}_0^{(1)} \longleftrightarrow -\Phi_{0256}$.

(3) $\mathbf{B}_0^{(2)} \longleftrightarrow -\Phi_{0134}$.

(4) $\mathbf{B}_1^{(1)} \longleftrightarrow \Phi_{0123}$.

(5) $\mathbf{B}_1^{(2)} \longleftrightarrow \Phi_{0146}=\Phi_{1256}$.

(6) $\mathbf{B}_3^{(1)} \longleftrightarrow \Phi_{0145}$.

(7) $\mathbf{B}_3^{(2)} \longleftrightarrow -\Phi_{0234}=-\Phi_{2456}$.

(8) $\mathbf{B}_9^{(1)} \longleftrightarrow -\Phi_{0346}$.

(9) $\mathbf{B}_9^{(2)} \longleftrightarrow -\Phi_{0135}=-\Phi_{2356}$.

(10) $\mathbf{B}_{12}^{(1)} \longleftrightarrow \Phi_{0235}$.

(11) $\mathbf{B}_{12}^{(2)} \longleftrightarrow \Phi_{0456}=\Phi_{1345}$.

(12) $\mathbf{B}_{10}^{(1)} \longleftrightarrow -\Phi_{0156}$.

(13) $\mathbf{B}_{10}^{(2)} \longleftrightarrow \Phi_{0236}=\Phi_{1346}$.

(14) $\mathbf{B}_4^{(1)} \longleftrightarrow -\Phi_{0246}$.

(15) $\mathbf{B}_4^{(2)} \longleftrightarrow \Phi_{0125}=\Phi_{1234}$.

(16) $\mathbf{B}_5 \longleftrightarrow \Phi_{0356}=\Phi_{1246}$.

(17) $\mathbf{B}_2 \longleftrightarrow -\Phi_{0126}=-\Phi_{2345}$.

(18) $\mathbf{B}_6 \longleftrightarrow -\Phi_{0245}=-\Phi_{1356}$.

(19) $\mathbf{B}_8 \longleftrightarrow -\Phi_{0124}=-\Phi_{3456}$.

(20) $\mathbf{B}_{11} \longleftrightarrow -\Phi_{0345}=-\Phi_{1236}$.

(21) $\mathbf{B}_7 \longleftrightarrow -\Phi_{0136}=-\Phi_{1245}$.

This completes the proof of Theorem 2.4.12.

\noindent
$\qquad \qquad \qquad \qquad \qquad \qquad \qquad \qquad \qquad
 \qquad \qquad \qquad \qquad \qquad \qquad \qquad \boxed{}$

  Theorem 2.4.12 gives a connection between Klein's quartic system,
i.e., the locus for the modular curve $X(13)$ (algebraic geometry)
and our quartic invariants associated with $\text{SL}(2, \mathbb{F}_{13})$
(invariant theory). This leads us to study Klein's quartic system
from the viewpoint of representation theory.

  Let $I=I(Y)$ be an ideal generated by the above twenty-one quartic
polynomials ($\mathbf{B}$-terms):
$$I=\langle \mathbf{B}_0^{(i)}, \mathbf{B}_1^{(j)}, \mathbf{B}_3^{(j)},
    \mathbf{B}_9^{(j)}, \mathbf{B}_{12}^{(j)}, \mathbf{B}_{10}^{(j)},
    \mathbf{B}_4^{(j)}, \mathbf{B}_5, \mathbf{B}_2, \mathbf{B}_6,
    \mathbf{B}_8, \mathbf{B}_{11}, \mathbf{B}_7 \rangle\eqno{(2.4.58)}$$
with $i=0. 1, 2$ and $j=1, 2$. The corresponding curve associated to
the ideal $I(Y)$ is denoted by $Y$.

\textbf{Theorem 2.4.13.} {\it The modular curve $X=X(13)$ is isomorphic
to the curve $Y$ in $\mathbb{CP}^5$.}

{\it Proof}. It is the consequence of Theorem 2.4.11 and Theorem 2.4.12.

\noindent
$\qquad \qquad \qquad \qquad \qquad \qquad \qquad \qquad \qquad
 \qquad \qquad \qquad \qquad \qquad \qquad \qquad \boxed{}$

\textbf{Theorem 2.4.14.} {\it The curve $Y$ can be constructed from the
invariant quartic Fano four-fold $\Phi_4=0$ in $\mathbb{CP}^5$, i.e.,
$$(z_3 z_4^3+z_1 z_5^3+z_2 z_6^3)-(z_6 z_1^3+z_4 z_2^3+z_5 z_3^3)
  +3(z_1 z_2 z_4 z_5+z_2 z_3 z_5 z_6+z_3 z_1 z_6 z_4)=0.\eqno{(2.4.59)}$$}

{\it Proof}. Note that the invariant quartic polynomial $\Phi_4$ (see (4.3.11)
for a proof of its $G$-invariance) has the following decomposition:
$$\Phi_4=3 \mathbf{B}_0^{(0)}+\mathbf{B}_0^{(1)}-\mathbf{B}_0^{(2)}.\eqno{(2.4.60)}$$
Let $S$ act on the three terms $\mathbf{B}_0^{(0)}$, $\mathbf{B}_0^{(1)}$
and $\mathbf{B}_0^{(2)}$, respectively. This gives that
$$\aligned
   13 S(\mathbf{B}_0^{(0)})
 =&5 \mathbf{B}_0^{(0)}+2 \mathbf{B}_0^{(1)}-2 \mathbf{B}_0^{(2)}
   +2(\mathbf{B}_5+\mathbf{B}_2+\mathbf{B}_6+\mathbf{B}_8+\mathbf{B}_{11}+\mathbf{B}_7)+\\
  &+(\mathbf{B}_1^{(2)}+\mathbf{B}_3^{(2)}+\mathbf{B}_9^{(2)})
   +(\mathbf{B}_{12}^{(2)}+\mathbf{B}_{10}^{(2)}+\mathbf{B}_4^{(2)})+\\
  &+(\mathbf{B}_1^{(1)}+\mathbf{B}_3^{(1)}+\mathbf{B}_9^{(1)})
   -(\mathbf{B}_{12}^{(1)}+\mathbf{B}_{10}^{(1)}+\mathbf{B}_4^{(1)}).
\endaligned\eqno{(2.4.61)}$$
On the other hand,
$$\aligned
   13 S(\mathbf{B}_0^{(1)})
 =&12 \mathbf{B}_0^{(0)}+\frac{7+\sqrt{13}}{2} \mathbf{B}_0^{(1)}+
   \frac{-7+\sqrt{13}}{2} \mathbf{B}_0^{(2)}+\\
  &-3 (\mathbf{B}_5+\mathbf{B}_2+\mathbf{B}_6+\mathbf{B}_8+\mathbf{B}_{11}
   +\mathbf{B}_7)+\\
  &+\frac{-3-3 \sqrt{13}}{2} (\mathbf{B}_1^{(2)}+\mathbf{B}_3^{(2)}
   +\mathbf{B}_9^{(2)})+\\
  &+\frac{-3+3 \sqrt{13}}{2} (\mathbf{B}_{12}^{(2)}+\mathbf{B}_{10}^{(2)}
   +\mathbf{B}_4^{(2)})+\\
  &+\frac{-3+\sqrt{13}}{2} (\mathbf{B}_1^{(1)}+\mathbf{B}_3^{(1)}+\mathbf{B}_9^{(1)})+\\
  &+\frac{3+\sqrt{13}}{2} (\mathbf{B}_{12}^{(1)}+\mathbf{B}_{10}^{(1)}+\mathbf{B}_4^{(1)}).
\endaligned\eqno{(2.4.62)}$$
$$\aligned
   13 S(\mathbf{B}_0^{(2)})
 =&-12 \mathbf{B}_0^{(0)}+\frac{-7+\sqrt{13}}{2} \mathbf{B}_0^{(1)}+
   \frac{7+\sqrt{13}}{2} \mathbf{B}_0^{(2)}+\\
  &+3 (\mathbf{B}_5+\mathbf{B}_2+\mathbf{B}_6+\mathbf{B}_8+\mathbf{B}_{11}
   +\mathbf{B}_7)+\\
  &+\frac{3-3 \sqrt{13}}{2} (\mathbf{B}_1^{(2)}+\mathbf{B}_3^{(2)}
   +\mathbf{B}_9^{(2)})+\\
  &+\frac{3+3 \sqrt{13}}{2} (\mathbf{B}_{12}^{(2)}+\mathbf{B}_{10}^{(2)}
   +\mathbf{B}_4^{(2)})+\\
  &+\frac{3+\sqrt{13}}{2} (\mathbf{B}_1^{(1)}+\mathbf{B}_3^{(1)}+\mathbf{B}_9^{(1)})+\\
  &+\frac{-3+\sqrt{13}}{2} (\mathbf{B}_{12}^{(1)}+\mathbf{B}_{10}^{(1)}+\mathbf{B}_4^{(1)}).
\endaligned\eqno{(2.4.63)}$$
Combining with (2.4.51), this leads to all of the $\mathbf{B}$-terms (2.4.44),
(2.4.45), (2.4.46) and (2.4.52) which we have defined above. Hence, the curve
$Y$ can be constructed from the invariant quartic Fano four-fold $\Phi_4=0$
by the action of the transformation $S$. Moreover,
$$\aligned
  &13 [S(\mathbf{B}_0^{(1)})-S(\mathbf{B}_0^{(2)})]\\
 =&24 \mathbf{B}_0^{(0)}+7 \mathbf{B}_0^{(1)}-7 \mathbf{B}_0^{(2)}
  -6(\mathbf{B}_5+\mathbf{B}_2+\mathbf{B}_6+\mathbf{B}_8+\mathbf{B}_{11}+\mathbf{B}_7)+\\
  &-3(\mathbf{B}_1^{(2)}+\mathbf{B}_3^{(2)}+\mathbf{B}_9^{(2)})
   -3(\mathbf{B}_{12}^{(2)}+\mathbf{B}_{10}^{(2)}+\mathbf{B}_4^{(2)})+\\
  &-3(\mathbf{B}_1^{(1)}+\mathbf{B}_3^{(1)}+\mathbf{B}_9^{(1)})
   +3(\mathbf{B}_{12}^{(1)}+\mathbf{B}_{10}^{(1)}+\mathbf{B}_4^{(1)}),
\endaligned$$
and
$$\aligned
  &13 S(3 \mathbf{B}_0^{(0)}+\mathbf{B}_0^{(1)}-\mathbf{B}_0^{(2)})\\
 =&15 \mathbf{B}_0^{(0)}+6 \mathbf{B}_0^{(1)}-6 \mathbf{B}_0^{(2)}
  +6(\mathbf{B}_5+\mathbf{B}_2+\mathbf{B}_6+\mathbf{B}_8+\mathbf{B}_{11}+\mathbf{B}_7)+\\
  &+24 \mathbf{B}_0^{(0)}+7 \mathbf{B}_0^{(1)}-7 \mathbf{B}_0^{(2)}
  -6(\mathbf{B}_5+\mathbf{B}_2+\mathbf{B}_6+\mathbf{B}_8+\mathbf{B}_{11}+\mathbf{B}_7)\\
 =&13 (3 \mathbf{B}_0^{(0)}+\mathbf{B}_0^{(1)}-\mathbf{B}_0^{(2)}),
\endaligned$$
 i.e.,
 $$S(3 \mathbf{B}_0^{(0)}+\mathbf{B}_0^{(1)}-\mathbf{B}_0^{(2)})
  =3 \mathbf{B}_0^{(0)}+\mathbf{B}_0^{(1)}-\mathbf{B}_0^{(2)}.$$
This proves again that $\Phi_4$ is invariant under the action of $S$.

\noindent
$\qquad \qquad \qquad \qquad \qquad \qquad \qquad \qquad \qquad
 \qquad \qquad \qquad \qquad \qquad \qquad \qquad \boxed{}$

\textbf{Theorem 2.4.15.} {\it The ideal $I(Y)$ is invariant under the action
of $G$, which gives a twenty-one dimensional representation of $G$.}

{\it Proof}. By the proof of Theorem 2.4.14, we have obtained the expression
of $S(\mathbf{B}_0^{(0)})$, $S(\mathbf{B}_0^{(1)})$ and $S(\mathbf{B}_0^{(2)})$.
The other $18$ terms can be divided into six triples:
$$(S(\mathbf{B}_1^{(1)}), S(\mathbf{B}_3^{(1)}), S(\mathbf{B}_9^{(1)})),$$
$$(S(\mathbf{B}_{12}^{(1)}), S(\mathbf{B}_{10}^{(1)}), S(\mathbf{B}_4^{(1)})),$$
$$(S(\mathbf{B}_1^{(2)}), S(\mathbf{B}_3^{(2)}), S(\mathbf{B}_9^{(1)})),$$
$$(S(\mathbf{B}_{12}^{(2)}), S(\mathbf{B}_{10}^{(2)}), S(\mathbf{B}_4^{(2)})),$$
$$(S(\mathbf{B}_5), S(\mathbf{B}_2), S(\mathbf{B}_6)),$$
$$(S(\mathbf{B}_8), S(\mathbf{B}_{11}), S(\mathbf{B}_7)).$$
Without loss of generality, we begin with the computation
of $S(\mathbf{B}_1^{(1)})$, $S(\mathbf{B}_{12}^{(1)})$, $S(\mathbf{B}_1^{(2)})$,
$S(\mathbf{B}_{12}^{(2)})$, $S(\mathbf{B}_5)$, and $S(\mathbf{B}_8)$. We have
$$\aligned
  &13 S(\mathbf{B}_1^{(1)})\\
 =&6 \mathbf{B}_0^{(0)}+\frac{-3+\sqrt{13}}{2} \mathbf{B}_0^{(1)}+
   \frac{3+\sqrt{13}}{2} \mathbf{B}_0^{(2)}+\\
  &+3 [(\zeta^9+\zeta^4+\zeta^6+\zeta^7+\zeta^2+\zeta^{11}) \mathbf{B}_5+\\
  &+(\zeta^3+\zeta^{10}+\zeta^2+\zeta^{11}+\zeta^5+\zeta^8) \mathbf{B}_2+\\
  &+(\zeta+\zeta^{12}+\zeta^5+\zeta^8+\zeta^6+\zeta^7) \mathbf{B}_6+\\
  &+(\zeta^9+\zeta^4+\zeta^6+\zeta^7+\zeta^3+\zeta^{10}) \mathbf{B}_8+\\
  &+(\zeta^3+\zeta^{10}+\zeta^2+\zeta^{11}+\zeta+\zeta^{12}) \mathbf{B}_{11}+\\
  &+(\zeta+\zeta^{12}+\zeta^5+\zeta^8+\zeta^9+\zeta^4) \mathbf{B}_7]+\\
  &+3 [(1+\zeta^9+\zeta^4) \mathbf{B}_1^{(2)}
   +(1+\zeta^3+\zeta^{10}) \mathbf{B}_3^{(2)}
   +(1+\zeta+\zeta^{12}) \mathbf{B}_9^{(2)}+\\
  &+(1+\zeta^6+\zeta^7) \mathbf{B}_{12}^{(2)}
   +(1+\zeta^2+\zeta^{11}) \mathbf{B}_{10}^{(2)}
   +(1+\zeta^5+\zeta^8) \mathbf{B}_4^{(2)}]+\\
  &+[1+(\zeta^9+\zeta^4)+2 (\zeta+\zeta^{12})+(\zeta^2+\zeta^{11})] \mathbf{B}_1^{(1)}+\\
  &+[1+(\zeta^3+\zeta^{10})+2 (\zeta^9+\zeta^4)+(\zeta^5+\zeta^8)] \mathbf{B}_3^{(1)}+\\
  &+[1+(\zeta+\zeta^{12})+2 (\zeta^3+\zeta^{10})+(\zeta^6+\zeta^7)] \mathbf{B}_9^{(1)}+\\
  &-[1+(\zeta^6+\zeta^7)+2 (\zeta^5+\zeta^8)+(\zeta^3+\zeta^{10})] \mathbf{B}_{12}^{(1)}+\\
  &-[1+(\zeta^2+\zeta^{11})+2 (\zeta^6+\zeta^7)+(\zeta+\zeta^{12})] \mathbf{B}_{10}^{(1)}+\\
  &-[1+(\zeta^5+\zeta^8)+2 (\zeta^2+\zeta^{11})+(\zeta^9+\zeta^4)] \mathbf{B}_4^{(1)}.
\endaligned\eqno{(2.4.64)}$$
$$\aligned
  &13 S(\mathbf{B}_{12}^{(1)})\\
 =&-6 \mathbf{B}_0^{(0)}+\frac{3+\sqrt{13}}{2} \mathbf{B}_0^{(1)}+
   \frac{-3+\sqrt{13}}{2} \mathbf{B}_0^{(2)}+\\
  &-3 [(\zeta^9+\zeta^4+\zeta^6+\zeta^7+\zeta^3+\zeta^{10}) \mathbf{B}_5+\\
  &+(\zeta^3+\zeta^{10}+\zeta^2+\zeta^{11}+\zeta+\zeta^{12}) \mathbf{B}_2+\\
  &+(\zeta+\zeta^{12}+\zeta^5+\zeta^8+\zeta^9+\zeta^4) \mathbf{B}_6+\\
  &+(\zeta^9+\zeta^4+\zeta^6+\zeta^7+\zeta^2+\zeta^{11}) \mathbf{B}_8+\\
  &+(\zeta^3+\zeta^{10}+\zeta^2+\zeta^{11}+\zeta^5+\zeta^8) \mathbf{B}_{11}+\\
  &+(\zeta+\zeta^{12}+\zeta^5+\zeta^8+\zeta^6+\zeta^7) \mathbf{B}_7]+\\
  &-3 [(1+\zeta^6+\zeta^7) \mathbf{B}_1^{(2)}
   +(1+\zeta^2+\zeta^{11}) \mathbf{B}_3^{(2)}
   +(1+\zeta^5+\zeta^8) \mathbf{B}_9^{(2)}+\\
  &+(1+\zeta^9+\zeta^4) \mathbf{B}_{12}^{(2)}
   +(1+\zeta^3+\zeta^{10}) \mathbf{B}_{10}^{(2)}
   +(1+\zeta+\zeta^{12}) \mathbf{B}_4^{(2)}]+\\
  &-[1+(\zeta^6+\zeta^7)+2 (\zeta^5+\zeta^8)+(\zeta^3+\zeta^{10})] \mathbf{B}_1^{(1)}+\\
  &-[1+(\zeta^2+\zeta^{11})+2 (\zeta^6+\zeta^7)+(\zeta+\zeta^{12})] \mathbf{B}_3^{(1)}+\\
  &-[1+(\zeta^5+\zeta^8)+2 (\zeta^2+\zeta^{11})+(\zeta^9+\zeta^4)] \mathbf{B}_9^{(1)}+\\
  &+[1+(\zeta^9+\zeta^4)+2 (\zeta+\zeta^{12})+(\zeta^2+\zeta^{11})] \mathbf{B}_{12}^{(1)}+\\
  &+[1+(\zeta^3+\zeta^{10})+2 (\zeta^9+\zeta^4)+(\zeta^5+\zeta^8)] \mathbf{B}_{10}^{(1)}+\\
  &+[1+(\zeta+\zeta^{12})+2 (\zeta^3+\zeta^{10})+(\zeta^6+\zeta^7)] \mathbf{B}_4^{(1)}.
\endaligned\eqno{(2.4.65)}$$
$$\aligned
  &13 S(\mathbf{B}_1^{(2)})\\
 =&2 \mathbf{B}_0^{(0)}+\frac{-1-\sqrt{13}}{2} \mathbf{B}_0^{(1)}+
   \frac{1-\sqrt{13}}{2} \mathbf{B}_0^{(2)}+\\
  &+(\zeta^9+\zeta^4+\zeta^6+\zeta^7+\zeta^2+\zeta^{11}) \mathbf{B}_5+\\
  &+(\zeta^3+\zeta^{10}+\zeta^2+\zeta^{11}+\zeta^5+\zeta^8) \mathbf{B}_2+\\
  &+(\zeta+\zeta^{12}+\zeta^5+\zeta^8+\zeta^6+\zeta^7) \mathbf{B}_6+\\
  &+(\zeta^9+\zeta^4+\zeta^6+\zeta^7+\zeta^3+\zeta^{10}) \mathbf{B}_8+\\
  &+(\zeta^3+\zeta^{10}+\zeta^2+\zeta^{11}+\zeta+\zeta^{12}) \mathbf{B}_{11}+\\
  &+(\zeta+\zeta^{12}+\zeta^5+\zeta^8+\zeta^9+\zeta^4) \mathbf{B}_7+\\
  &+[-1-(\zeta^9+\zeta^4)+2 (\zeta+\zeta^{12})+(\zeta^2+\zeta^{11})] \mathbf{B}_1^{(2)}+\\
  &+[-1-(\zeta^3+\zeta^{10})+2 (\zeta^9+\zeta^4)+(\zeta^5+\zeta^8)] \mathbf{B}_3^{(2)}+\\
  &+[-1-(\zeta+\zeta^{12})+2 (\zeta^3+\zeta^{10})+(\zeta^6+\zeta^7)] \mathbf{B}_9^{(2)}+\\
  &+[-1-(\zeta^6+\zeta^7)+2 (\zeta^5+\zeta^8)+(\zeta^3+\zeta^{10})] \mathbf{B}_{12}^{(2)}+\\
  &+[-1-(\zeta^2+\zeta^{11})+2 (\zeta^6+\zeta^7)+(\zeta+\zeta^{12})] \mathbf{B}_{10}^{(2)}+\\
  &+[-1-(\zeta^5+\zeta^8)+2 (\zeta^2+\zeta^{11})+(\zeta^9+\zeta^4)] \mathbf{B}_4^{(2)}+\\
  &+(1+\zeta^9+\zeta^4) \mathbf{B}_1^{(1)}+(1+\zeta^3+\zeta^{10}) \mathbf{B}_3^{(1)}
   +(1+\zeta+\zeta^{12}) \mathbf{B}_9^{(1)}+\\
  &-(1+\zeta^6+\zeta^7) \mathbf{B}_{12}^{(1)}-(1+\zeta^2+\zeta^{11}) \mathbf{B}_{10}^{(1)}
   -(1+\zeta^5+\zeta^8) \mathbf{B}_4^{(1)}.
\endaligned\eqno{(2.4.66)}$$
$$\aligned
  &13 S(\mathbf{B}_{12}^{(2)})\\
 =&2 \mathbf{B}_0^{(0)}+\frac{-1+\sqrt{13}}{2} \mathbf{B}_0^{(1)}+
   \frac{1+\sqrt{13}}{2} \mathbf{B}_0^{(2)}+\\
  &+(\zeta^9+\zeta^4+\zeta^6+\zeta^7+\zeta^3+\zeta^{10}) \mathbf{B}_5+\\
  &+(\zeta^3+\zeta^{10}+\zeta^2+\zeta^{11}+\zeta+\zeta^{12}) \mathbf{B}_2+\\
  &+(\zeta+\zeta^{12}+\zeta^5+\zeta^8+\zeta^9+\zeta^4) \mathbf{B}_6+\\
  &+(\zeta^9+\zeta^4+\zeta^6+\zeta^7+\zeta^2+\zeta^{11}) \mathbf{B}_8+\\
  &+(\zeta^3+\zeta^{10}+\zeta^2+\zeta^{11}+\zeta^5+\zeta^8) \mathbf{B}_{11}+\\
  &+(\zeta+\zeta^{12}+\zeta^5+\zeta^8+\zeta^6+\zeta^7) \mathbf{B}_7+\\
  &+[-1-(\zeta^6+\zeta^7)+2 (\zeta^5+\zeta^8)+(\zeta^3+\zeta^{10})] \mathbf{B}_1^{(2)}+\\
  &+[-1-(\zeta^2+\zeta^{11})+2 (\zeta^6+\zeta^7)+(\zeta+\zeta^{12})] \mathbf{B}_3^{(2)}+\\
  &+[-1-(\zeta^5+\zeta^8)+2 (\zeta^2+\zeta^{11})+(\zeta^9+\zeta^4)] \mathbf{B}_9^{(2)}+\\
  &+[-1-(\zeta^9+\zeta^4)+2 (\zeta+\zeta^{12})+(\zeta^2+\zeta^{11})] \mathbf{B}_{12}^{(2)}+\\
  &+[-1-(\zeta^3+\zeta^{10})+2 (\zeta^9+\zeta^4)+(\zeta^5+\zeta^8)] \mathbf{B}_{10}^{(2)}+\\
  &+[-1-(\zeta+\zeta^{12})+2 (\zeta^3+\zeta^{10})+(\zeta^6+\zeta^7)] \mathbf{B}_4^{(2)}+\\
  &+(1+\zeta^6+\zeta^7) \mathbf{B}_1^{(1)}+(1+\zeta^2+\zeta^{11}) \mathbf{B}_3^{(1)}
   +(1+\zeta^5+\zeta^8) \mathbf{B}_9^{(1)}+\\
  &-(1+\zeta^9+\zeta^4) \mathbf{B}_{12}^{(1)}-(1+\zeta^3+\zeta^{10}) \mathbf{B}_{10}^{(1)}
   -(1+\zeta+\zeta^{12}) \mathbf{B}_4^{(1)}.
\endaligned\eqno{(2.4.67)}$$
$$\aligned
  &13 S(\mathbf{B}_5)\\
 =&4 \mathbf{B}_0^{(0)}-\mathbf{B}_0^{(1)}+\mathbf{B}_0^{(2)}+\\
  &-[(\zeta^3+\zeta^{10})+2 (\zeta^2+\zeta^{11})+2 (\zeta^9+\zeta^4)
   +2 (\zeta+\zeta^{12})] \mathbf{B}_5+\\
  &-[(\zeta+\zeta^{12})+2 (\zeta^5+\zeta^8)+2 (\zeta^3+\zeta^{10})
   +2 (\zeta^9+\zeta^4)] \mathbf{B}_2+\\
  &-[(\zeta^9+\zeta^4)+2 (\zeta^6+\zeta^7)+2 (\zeta+\zeta^{12})
   +2 (\zeta^3+\zeta^{10})] \mathbf{B}_6+\\
  &-[(\zeta^2+\zeta^{11})+2 (\zeta^3+\zeta^{10})+2 (\zeta^6+\zeta^7)
   +2 (\zeta^5+\zeta^8)] \mathbf{B}_8+\\
  &-[(\zeta^5+\zeta^8)+2 (\zeta+\zeta^{12})+2 (\zeta^2+\zeta^{11})
   +2 (\zeta^6+\zeta^7)] \mathbf{B}_{11}+\\
  &-[(\zeta^6+\zeta^7)+2 (\zeta^9+\zeta^4)+2 (\zeta^5+\zeta^8)
   +2 (\zeta^2+\zeta^{11})] \mathbf{B}_7+\\
  &+(\zeta^9+\zeta^4+\zeta^6+\zeta^7+\zeta^2+\zeta^{11}) \mathbf{B}_1^{(2)}+\\
  &+(\zeta^3+\zeta^{10}+\zeta^2+\zeta^{11}+\zeta^5+\zeta^8) \mathbf{B}_3^{(2)}+\\
  &+(\zeta+\zeta^{12}+\zeta^5+\zeta^8+\zeta^6+\zeta^7) \mathbf{B}_9^{(2)}+\\
  &+(\zeta^9+\zeta^4+\zeta^6+\zeta^7+\zeta^3+\zeta^{10}) \mathbf{B}_{12}^{(2)}+\\
  &+(\zeta^3+\zeta^{10}+\zeta^2+\zeta^{11}+\zeta+\zeta^{12}) \mathbf{B}_{10}^{(2)}+\\
  &+(\zeta+\zeta^{12}+\zeta^5+\zeta^8+\zeta^9+\zeta^4) \mathbf{B}_4^{(2)}+\\
  &+(\zeta^9+\zeta^4+\zeta^6+\zeta^7+\zeta^2+\zeta^{11}) \mathbf{B}_1^{(1)}+\\
  &+(\zeta^3+\zeta^{10}+\zeta^2+\zeta^{11}+\zeta^5+\zeta^8) \mathbf{B}_3^{(1)}+\\
  &+(\zeta+\zeta^{12}+\zeta^5+\zeta^8+\zeta^6+\zeta^7) \mathbf{B}_9^{(1)}+\\
  &-(\zeta^9+\zeta^4+\zeta^6+\zeta^7+\zeta^3+\zeta^{10}) \mathbf{B}_{12}^{(1)}+\\
  &-(\zeta^3+\zeta^{10}+\zeta^2+\zeta^{11}+\zeta+\zeta^{12}) \mathbf{B}_{10}^{(1)}+\\
  &-(\zeta+\zeta^{12}+\zeta^5+\zeta^8+\zeta^9+\zeta^4) \mathbf{B}_4^{(1)}.
\endaligned\eqno{(2.4.68)}$$
$$\aligned
  &13 S(\mathbf{B}_8)\\
 =&4 \mathbf{B}_0^{(0)}-\mathbf{B}_0^{(1)}+\mathbf{B}_0^{(2)}+\\
  &-[(\zeta^2+\zeta^{11})+2 (\zeta^3+\zeta^{10})+2 (\zeta^6+\zeta^7)
   +2 (\zeta^5+\zeta^8)] \mathbf{B}_5+\\
  &-[(\zeta^5+\zeta^8)+2 (\zeta+\zeta^{12})+2 (\zeta^2+\zeta^{11})
   +2 (\zeta^6+\zeta^7)] \mathbf{B}_2+\\
  &-[(\zeta^6+\zeta^7)+2 (\zeta^9+\zeta^4)+2 (\zeta^5+\zeta^8)
   +2 (\zeta^2+\zeta^{11})] \mathbf{B}_6+\\
  &-[(\zeta^3+\zeta^{10})+2 (\zeta^2+\zeta^{11})+2 (\zeta^9+\zeta^4)
   +2 (\zeta+\zeta^{12})] \mathbf{B}_8+\\
  &-[(\zeta+\zeta^{12})+2 (\zeta^5+\zeta^8)+2 (\zeta^3+\zeta^{10})
   +2 (\zeta^9+\zeta^4)] \mathbf{B}_{11}+\\
  &-[(\zeta^9+\zeta^4)+2 (\zeta^6+\zeta^7)+2 (\zeta+\zeta^{12})
   +2 (\zeta^3+\zeta^{10})] \mathbf{B}_7+\\
  &+(\zeta^9+\zeta^4+\zeta^6+\zeta^7+\zeta^3+\zeta^{10}) \mathbf{B}_1^{(2)}+\\
  &+(\zeta^3+\zeta^{10}+\zeta^2+\zeta^{11}+\zeta+\zeta^{12}) \mathbf{B}_3^{(2)}+\\
  &+(\zeta+\zeta^{12}+\zeta^5+\zeta^8+\zeta^9+\zeta^4) \mathbf{B}_9^{(2)}+\\
  &+(\zeta^9+\zeta^4+\zeta^6+\zeta^7+\zeta^2+\zeta^{11}) \mathbf{B}_{12}^{(2)}+\\
  &+(\zeta^3+\zeta^{10}+\zeta^2+\zeta^{11}+\zeta^5+\zeta^8) \mathbf{B}_{10}^{(2)}+\\
  &+(\zeta+\zeta^{12}+\zeta^5+\zeta^8+\zeta^6+\zeta^7) \mathbf{B}_4^{(2)}+\\
  &+(\zeta^9+\zeta^4+\zeta^6+\zeta^7+\zeta^3+\zeta^{10}) \mathbf{B}_1^{(1)}+\\
  &+(\zeta^3+\zeta^{10}+\zeta^2+\zeta^{11}+\zeta+\zeta^{12}) \mathbf{B}_3^{(1)}+\\
  &+(\zeta+\zeta^{12}+\zeta^5+\zeta^8+\zeta^9+\zeta^4) \mathbf{B}_9^{(1)}+\\
  &-(\zeta^9+\zeta^4+\zeta^6+\zeta^7+\zeta^2+\zeta^{11}) \mathbf{B}_{12}^{(1)}+\\
  &-(\zeta^3+\zeta^{10}+\zeta^2+\zeta^{11}+\zeta^5+\zeta^8) \mathbf{B}_{10}^{(1)}+\\
  &-(\zeta+\zeta^{12}+\zeta^5+\zeta^8+\zeta^6+\zeta^7) \mathbf{B}_4^{(1)}.
\endaligned\eqno{(2.4.69)}$$

  The other terms: $S(\mathbf{B}_3^{(1)})$, $S(\mathbf{B}_{10}^{(1)})$,
$S(\mathbf{B}_3^{(2)})$, $S(\mathbf{B}_{10}^{(2)})$, $S(\mathbf{B}_2)$,
$S(\mathbf{B}_{11})$, $S(\mathbf{B}_9^{(1)})$, $S(\mathbf{B}_4^{(1)})$,
$S(\mathbf{B}_9^{(2)})$, $S(\mathbf{B}_4^{(2)})$, $S(\mathbf{B}_6)$, and
$S(\mathbf{B}_7)$ can be obtained by the permutation on the above six terms.
For each of the above triples, we can obtain the other two terms by the
following permutation
$$\zeta^9 \mapsto \zeta^3 \mapsto \zeta, \quad
  \zeta^4 \mapsto \zeta^{10} \mapsto \zeta^{12}, \quad
  \zeta^6 \mapsto \zeta^2 \mapsto \zeta^5, \quad
  \zeta^7 \mapsto \zeta^{11} \mapsto \zeta^8\eqno{(2.4.70)}$$
on the first term. For example, for the first triple
$(S(\mathbb{B}_1^{(1)}), S(\mathbb{B}_3^{(1)}), S(\mathbb{B}_9^{(1)}))$,
we can get the expression of $S(\mathbb{B}_3^{(1)})$ and $S(\mathbb{B}_9^{(1)})$
by the above permutation (12.49) on $S(\mathbb{B}_1^{(1)})$ as follows:
$$\aligned
  &13 S(\mathbf{B}_3^{(1)})\\
 =&6 \mathbf{B}_0^{(0)}+\frac{-3+\sqrt{13}}{2} \mathbf{B}_0^{(1)}+
   \frac{3+\sqrt{13}}{2} \mathbf{B}_0^{(2)}+\\
  &+3 [(\zeta^3+\zeta^{10}+\zeta^2+\zeta^{11}+\zeta^5+\zeta^8) \mathbf{B}_5+\\
  &+(\zeta+\zeta^{12}+\zeta^5+\zeta^8+\zeta^6+\zeta^7) \mathbf{B}_2+\\
  &+(\zeta^9+\zeta^4+\zeta^6+\zeta^7+\zeta^2+\zeta^{11}) \mathbf{B}_6+\\
  &+(\zeta^3+\zeta^{10}+\zeta^2+\zeta^{11}+\zeta+\zeta^{12}) \mathbf{B}_8+\\
  &+(\zeta+\zeta^{12}+\zeta^5+\zeta^8+\zeta^9+\zeta^4) \mathbf{B}_{11}+\\
  &+(\zeta^9+\zeta^4+\zeta^6+\zeta^7+\zeta^3+\zeta^{10}) \mathbf{B}_7]+\\
  &+3 [(1+\zeta^3+\zeta^{10}) \mathbf{B}_1^{(2)}
   +(1+\zeta+\zeta^{12}) \mathbf{B}_3^{(2)}
   +(1+\zeta^9+\zeta^4) \mathbf{B}_9^{(2)}+\\
  &+(1+\zeta^2+\zeta^{11}) \mathbf{B}_{12}^{(2)}
   +(1+\zeta^5+\zeta^8) \mathbf{B}_{10}^{(2)}
   +(1+\zeta^6+\zeta^7) \mathbf{B}_4^{(2)}]+\\
  &+[1+(\zeta^3+\zeta^{10})+2 (\zeta^9+\zeta^4)+(\zeta^5+\zeta^8)] \mathbf{B}_1^{(1)}+\\
  &+[1+(\zeta+\zeta^{12})+2 (\zeta^3+\zeta^{10})+(\zeta^6+\zeta^7)] \mathbf{B}_3^{(1)}+\\
  &+[1+(\zeta^9+\zeta^4)+2 (\zeta+\zeta^{12})+(\zeta^2+\zeta^{11})] \mathbf{B}_9^{(1)}+\\
  &-[1+(\zeta^2+\zeta^{11})+2 (\zeta^6+\zeta^7)+(\zeta+\zeta^{12})] \mathbf{B}_{12}^{(1)}+\\
  &-[1+(\zeta^5+\zeta^8)+2 (\zeta^2+\zeta^{11})+(\zeta^9+\zeta^4)] \mathbf{B}_{10}^{(1)}+\\
  &-[1+(\zeta^6+\zeta^7)+2 (\zeta^5+\zeta^8)+(\zeta^3+\zeta^{10})] \mathbf{B}_4^{(1)}.
\endaligned\eqno{(2.4.71)}$$
$$\aligned
  &13 S(\mathbf{B}_9^{(1)})\\
 =&6 \mathbf{B}_0^{(0)}+\frac{-3+\sqrt{13}}{2} \mathbf{B}_0^{(1)}+
   \frac{3+\sqrt{13}}{2} \mathbf{B}_0^{(2)}+\\
  &+3 [(\zeta+\zeta^{12}+\zeta^5+\zeta^8+\zeta^6+\zeta^7) \mathbf{B}_5+\\
  &+(\zeta^9+\zeta^4+\zeta^6+\zeta^7+\zeta^2+\zeta^{11}) \mathbf{B}_2+\\
  &+(\zeta^3+\zeta^{10}+\zeta^2+\zeta^{11}+\zeta^5+\zeta^8) \mathbf{B}_6+\\
  &+(\zeta+\zeta^{12}+\zeta^5+\zeta^8+\zeta^9+\zeta^4) \mathbf{B}_8+\\
  &+(\zeta^9+\zeta^4+\zeta^6+\zeta^7+\zeta^3+\zeta^{10}) \mathbf{B}_{11}+\\
  &+(\zeta^3+\zeta^{10}+\zeta^2+\zeta^{11}+\zeta+\zeta^{12}) \mathbf{B}_7]+\\
  &+3 [(1+\zeta+\zeta^{12}) \mathbf{B}_1^{(2)}
   +(1+\zeta^9+\zeta^4) \mathbf{B}_3^{(2)}
   +(1+\zeta^3+\zeta^{10}) \mathbf{B}_9^{(2)}+\\
  &+(1+\zeta^5+\zeta^8) \mathbf{B}_{12}^{(2)}
   +(1+\zeta^6+\zeta^7) \mathbf{B}_{10}^{(2)}
   +(1+\zeta^2+\zeta^{11}) \mathbf{B}_4^{(2)}]+\\
  &+[1+(\zeta+\zeta^{12})+2 (\zeta^3+\zeta^{10})+(\zeta^6+\zeta^7)] \mathbf{B}_1^{(1)}+\\
  &+[1+(\zeta^9+\zeta^4)+2 (\zeta+\zeta^{12})+(\zeta^2+\zeta^{11})] \mathbf{B}_3^{(1)}+\\
  &+[1+(\zeta^3+\zeta^{10})+2 (\zeta^9+\zeta^4)+(\zeta^5+\zeta^8)] \mathbf{B}_9^{(1)}+\\
  &-[1+(\zeta^5+\zeta^8)+2 (\zeta^2+\zeta^{11})+(\zeta^9+\zeta^4)] \mathbf{B}_{12}^{(1)}+\\
  &-[1+(\zeta^6+\zeta^7)+2 (\zeta^5+\zeta^8)+(\zeta^3+\zeta^{10})] \mathbf{B}_{10}^{(1)}+\\
  &-[1+(\zeta^2+\zeta^{11})+2 (\zeta^6+\zeta^7)+(\zeta+\zeta^{12})] \mathbf{B}_4^{(1)}.
\endaligned\eqno{(2.4.72)}$$
The other five triples can be obtained in the similar way.

  On the other hand, we have
$$T(\mathbf{B}_0^{(0)})=\mathbf{B}_0^{(0)}, \quad
  T(\mathbf{B}_0^{(1)})=\mathbf{B}_0^{(1)}, \quad
  T(\mathbf{B}_0^{(2)})=\mathbf{B}_0^{(2)}.$$
$$T(\mathbf{B}_1^{(1)})=\zeta \mathbf{B}_1^{(1)}, \quad
  T(\mathbf{B}_1^{(2)})=\zeta \mathbf{B}_1^{(2)}.$$
$$T(\mathbf{B}_3^{(1)})=\zeta^3 \mathbf{B}_3^{(1)}, \quad
  T(\mathbf{B}_3^{(2)})=\zeta^3 \mathbf{B}_3^{(2)}.$$
$$T(\mathbf{B}_9^{(1)})=\zeta^9 \mathbf{B}_9^{(1)}, \quad
  T(\mathbf{B}_9^{(2)})=\zeta^9 \mathbf{B}_9^{(2)}.$$
$$T(\mathbf{B}_{12}^{(1)})=\zeta^{12} \mathbf{B}_{12}^{(1)}, \quad
  T(\mathbf{B}_{12}^{(2)})=\zeta^{12} \mathbf{B}_{12}^{(2)}.$$
$$T(\mathbf{B}_{10}^{(1)})=\zeta^{10} \mathbf{B}_{10}^{(1)}, \quad
  T(\mathbf{B}_{10}^{(2)})=\zeta^{10} \mathbf{B}_{10}^{(2)}.$$
$$T(\mathbf{B}_4^{(1)})=\zeta^4 \mathbf{B}_4^{(1)}, \quad
  T(\mathbf{B}_4^{(2)})=\zeta^4 \mathbf{B}_4^{(2)}.$$
$$T(\mathbf{B}_5)=\zeta^5 \mathbf{B}_5, \quad
  T(\mathbf{B}_8)=\zeta^8 \mathbf{B}_8.$$
$$T(\mathbf{B}_2)=\zeta^2 \mathbf{B}_2, \quad
  T(\mathbf{B}_{11})=\zeta^{11} \mathbf{B}_{11}.$$
$$T(\mathbf{B}_6)=\zeta^6 \mathbf{B}_6, \quad
  T(\mathbf{B}_7)=\zeta^7 \mathbf{B}_7.$$

  Let $\widetilde{S}$ and $\widetilde{T}$ denote the corresponding matrices
of order $21$ corresponding to $S$ and $T$. Then they have the following trace:
$$\text{Tr}(\widetilde{S})=1, \quad \text{Tr}(\widetilde{T})=\frac{3+\sqrt{13}}{2}.
  \eqno{(2.4.73)}$$
This completes the proof of Theorem 2.4.15.

\noindent
$\qquad \qquad \qquad \qquad \qquad \qquad \qquad \qquad \qquad
 \qquad \qquad \qquad \qquad \qquad \qquad \qquad \boxed{}$

  Now, we have constructed a $21$-dimensional reducible representation
of $\text{SL}(2, \mathbb{F}_{13})$. We will give its decomposition:
$\mathbf{21}=\mathbf{1} \oplus \mathbf{7} \oplus \mathbf{13}$. For the
group $\text{SL}(2, \mathbb{F}_q)$, when $q=13$, $|G|=2184$. There are
at most $q+4=17$ conjugacy classes and
$$1^2+13^2+5 \times 14^2+2 \times 7^2+6 \times 12^2+2 \times 6^2=2184$$
corresponding to the decomposition of conjugacy classes:
$$1+1+5+2+6+2=17.$$
For $G=\text{SL}(2, \mathbb{F}_{13})$. There are the following
representations $\rho$:

(1) discrete series
$$\text{dim}(\rho)=\left\{\aligned
  &12 \quad \text{(super-cuspidal representations) or}\\
  &6 \quad \text{(two degenerate discrete series representations)}
\endaligned\right.$$

(2) principal series
$$\text{dim}(\rho)=\left\{\aligned
  &14 \quad \text{(principal series representations) or}\\
  &13 \quad \text{(Steinberg representation) or}\\
  &7 \quad \text{(two degenerate principal series representations)}
\endaligned\right.$$

\textbf{Theorem 2.4.16}. {\it The $21$-dimensional representation is reducible,
which can be decomposed as the direct sum of a $1$-dimensional representation
$($the trivial representation$)$, a $7$-dimensional representation $($the
degenerate principal series representation) and a $13$-dimensional representation
$($the Steinberg representation):
$$\mathbf{21}=\mathbf{1} \oplus \mathbf{7} \oplus \mathbf{13}.\eqno{(2.4.74)}$$
More precisely, let $V$ be a complex vector space generated by the $21$ quartic
polynomials $\mathbf{B}_0^{(0)}$, $\mathbf{B}_0^{(1)}$, $\ldots$, $\mathbf{B}_7$.
Then $V$ has the following decomposition:
$$V=V_1 \oplus V_7 \oplus V_{13},\eqno{(2.4.75)}$$
where the $1$-dimensional subspace
$$V_1=\langle 3 \mathbf{B}_0^{(0)}+\mathbf{B}_0^{(1)}-\mathbf{B}_0^{(2)} \rangle,
      \eqno{(2.4.76)}$$
the $7$-dimensional subspace
$$\aligned
  V_7=\langle &\mathbf{B}_0^{(1)}+\mathbf{B}_0^{(2)},
               \mathbf{B}_1^{(1)}-3 \mathbf{B}_1^{(2)},
               \mathbf{B}_3^{(1)}-3 \mathbf{B}_3^{(2)},
               \mathbf{B}_9^{(1)}-3 \mathbf{B}_9^{(2)},\\
              &\mathbf{B}_{12}^{(1)}+3 \mathbf{B}_{12}^{(2)},
               \mathbf{B}_{10}^{(1)}+3 \mathbf{B}_{10}^{(2)},
               \mathbf{B}_4^{(1)}+3 \mathbf{B}_4^{(2)}
\rangle\endaligned\eqno{(2.4.77)}$$
and the $13$-dimensional subspace
$$\aligned
  V_{13}=\langle &4 \mathbf{B}_0^{(0)}-\mathbf{B}_0^{(1)}+\mathbf{B}_0^{(2)},
               \mathbf{B}_5, \mathbf{B}_2, \mathbf{B}_6, \mathbf{B}_8,
               \mathbf{B}_{11}, \mathbf{B}_7,\\
              &\mathbf{B}_1^{(1)}+\mathbf{B}_1^{(2)},
               \mathbf{B}_3^{(1)}+\mathbf{B}_3^{(2)},
               \mathbf{B}_9^{(1)}+\mathbf{B}_9^{(2)},\\
              &-\mathbf{B}_{12}^{(1)}+\mathbf{B}_{12}^{(2)},
               -\mathbf{B}_{10}^{(1)}+\mathbf{B}_{10}^{(2)},
               -\mathbf{B}_4^{(1)}+\mathbf{B}_4^{(2)} \rangle
\endaligned\eqno{(2.4.78)}$$
are stable under the action of} $G \cong \text{SL}(2, \mathbb{F}_{13})$.

{\it Proof}. Our $21$-dimensional representation does not belong to the
discrete series or the principal series. By the proof of Theorem 2.4.12 and
Theorem 2.4.13, we have
$$S(3 \mathbf{B}_0^{(0)}+\mathbf{B}_0^{(1)}-\mathbf{B}_0^{(2)})
 =3 \mathbf{B}_0^{(0)}+\mathbf{B}_0^{(1)}-\mathbf{B}_0^{(2)}.\eqno{(2.4.79)}$$
$$T(3 \mathbf{B}_0^{(0)}+\mathbf{B}_0^{(1)}-\mathbf{B}_0^{(2)})
 =3 \mathbf{B}_0^{(0)}+\mathbf{B}_0^{(1)}-\mathbf{B}_0^{(2)}.\eqno{(2.4.80)}$$
Let $V$ be a complex vector space generated by the $21$ quartic polynomials
$\mathbf{B}_0^{(0)}$, $\mathbf{B}_0^{(1)}$, $\ldots$, $\mathbf{B}_7$. Put
$$V_1=\langle 3 \mathbf{B}_0^{(0)}+\mathbf{B}_0^{(1)}-\mathbf{B}_0^{(2)} \rangle.
  \eqno{(2.4.81)}$$
Then $V_1$ is a one-dimensional subspace which is stable under the action of
$G \cong \text{SL}(2, \mathbb{F}_{13})$. Thus,
$$V=V_1 \oplus V_1^0,\eqno{(2.4.82)}$$
where $V_1^0$ is the complement of $V_1$ in $V$. This shows that this $21$-dimensional
representation is reducible. It can be decomposed as a direct sum of a
$1$-dimensional representation and a $20$-dimensional representation
$$\mathbf{21}=\mathbf{1} \oplus \mathbf{20},\eqno{(2.4.83)}$$
the latter one can be decomposed into the sum of irreducible representations
listed as above. More precisely, we have
$$\aligned
  &\sqrt{13} S(\mathbf{B}_0^{(1)}+\mathbf{B}_0^{(2)})\\
 =&(\mathbf{B}_0^{(1)}+\mathbf{B}_0^{(2)})
  +(\mathbf{B}_1^{(1)}-3 \mathbf{B}_1^{(2)})
  +(\mathbf{B}_3^{(1)}-3 \mathbf{B}_3^{(2)})
  +(\mathbf{B}_9^{(1)}-3 \mathbf{B}_9^{(2)})\\
 +&(\mathbf{B}_{12}^{(1)}+3 \mathbf{B}_{12}^{(2)})+
   (\mathbf{B}_{10}^{(1)}+3 \mathbf{B}_{10}^{(2)})+
   (\mathbf{B}_4^{(1)}+3 \mathbf{B}_4^{(2)}).
\endaligned\eqno{(2.4.84)}$$
$$T(\mathbf{B}_0^{(1)}+\mathbf{B}_0^{(2)})
 =\mathbf{B}_0^{(1)}+\mathbf{B}_0^{(2)}.\eqno{(2.4.85)}$$
Without loss of generality, we can only consider the action of $S$
and $T$ on $\mathbf{B}_1^{(1)}-3 \mathbf{B}_1^{(2)}$ and
$\mathbf{B}_{12}^{(1)}+3 \mathbf{B}_{12}^{(2)}$. The action of $S$
and $T$ on
$\mathbf{B}_3^{(1)}-3 \mathbf{B}_3^{(2)}$,
$\mathbf{B}_9^{(1)}-3 \mathbf{B}_9^{(2)}$,
$\mathbf{B}_{10}^{(1)}+3 \mathbf{B}_{10}^{(2)}$,
and $\mathbf{B}_4^{(1)}+3 \mathbf{B}_4^{(2)}$
can be done similarly. We have
$$\aligned
  &13 S(\mathbf{B}_1^{(1)}-3 \mathbf{B}_1^{(2)})\\
 =&2 \sqrt{13} (\mathbf{B}_0^{(1)}+\mathbf{B}_0^{(2)})\\
  &+[-2 (1+\zeta^9+\zeta^4)+2 (\zeta+\zeta^{12})+(\zeta^2+\zeta^{11})]
  (\mathbf{B}_1^{(1)}-3 \mathbf{B}_1^{(2)})\\
  &+[-2 (1+\zeta^3+\zeta^{10})+2 (\zeta^9+\zeta^4)+(\zeta^5+\zeta^8)]
  (\mathbf{B}_3^{(1)}-3 \mathbf{B}_3^{(2)})\\
  &+[-2 (1+\zeta+\zeta^{12})+2 (\zeta^3+\zeta^{10})+(\zeta^6+\zeta^7)]
  (\mathbf{B}_9^{(1)}-3 \mathbf{B}_9^{(2)})\\
  &+[2 (1+\zeta^6+\zeta^7)-2 (\zeta^5+\zeta^8)-(\zeta^3+\zeta^{10})]
  (\mathbf{B}_{12}^{(1)}+3 \mathbf{B}_{12}^{(2)})\\
  &+[2 (1+\zeta^2+\zeta^{11})-2 (\zeta^6+\zeta^7)-(\zeta+\zeta^{12})]
  (\mathbf{B}_{10}^{(1)}+3 \mathbf{B}_{10}^{(2)})\\
  &+[2 (1+\zeta^5+\zeta^8)-2 (\zeta^2+\zeta^{11})-(\zeta^9+\zeta^4)]
  (\mathbf{B}_4^{(1)}+3 \mathbf{B}_4^{(2)}).
\endaligned\eqno{(2.4.86)}$$
$$\aligned
  &13 S(\mathbf{B}_{12}^{(1)}+3 \mathbf{B}_{12}^{(2)})\\
 =&2 \sqrt{13} (\mathbf{B}_0^{(1)}+\mathbf{B}_0^{(2)})\\
  &+[2 (1+\zeta^6+\zeta^7)-2 (\zeta^5+\zeta^8)-(\zeta^3+\zeta^{10})]
  (\mathbf{B}_1^{(1)}-3 \mathbf{B}_1^{(2)})\\
  &+[2 (1+\zeta^2+\zeta^{11})-2 (\zeta^6+\zeta^7)-(\zeta+\zeta^{12})]
  (\mathbf{B}_3^{(1)}-3 \mathbf{B}_3^{(2)})\\
  &+[2 (1+\zeta^5+\zeta^8)-2 (\zeta^2+\zeta^{11})-(\zeta^9+\zeta^4)]
  (\mathbf{B}_9^{(1)}-3 \mathbf{B}_9^{(2)})\\
  &+[-2 (1+\zeta^9+\zeta^4)+2 (\zeta+\zeta^{12})+(\zeta^2+\zeta^{11})]
  (\mathbf{B}_{12}^{(1)}+3 \mathbf{B}_{12}^{(2)})\\
  &+[-2 (1+\zeta^3+\zeta^{10})+2 (\zeta^9+\zeta^4)+(\zeta^5+\zeta^8)]
  (\mathbf{B}_{10}^{(1)}+3 \mathbf{B}_{10}^{(2)})\\
  &+[-2 (1+\zeta+\zeta^{12})+2 (\zeta^3+\zeta^{10})+(\zeta^6+\zeta^7)]
  (\mathbf{B}_4^{(1)}+3 \mathbf{B}_4^{(2)}).
\endaligned\eqno{(2.4.87)}$$
$$T(\mathbf{B}_1^{(1)}-3 \mathbf{B}_1^{(2)})
 =\zeta (\mathbf{B}_1^{(1)}-3 \mathbf{B}_1^{(2)}).\eqno{(2.4.88)}$$
$$T(\mathbf{B}_{12}^{(1)}+3 \mathbf{B}_{12}^{(2)})
 =\zeta^{12} (\mathbf{B}_{12}^{(1)}+3 \mathbf{B}_{12}^{(2)}).\eqno{(2.4.89)}$$
Put
$$\aligned
  V_7=\langle &\mathbf{B}_0^{(1)}+\mathbf{B}_0^{(2)},
               \mathbf{B}_1^{(1)}-3 \mathbf{B}_1^{(2)},
               \mathbf{B}_3^{(1)}-3 \mathbf{B}_3^{(2)},
               \mathbf{B}_9^{(1)}-3 \mathbf{B}_9^{(2)},\\
              &\mathbf{B}_{12}^{(1)}+3 \mathbf{B}_{12}^{(2)},
               \mathbf{B}_{10}^{(1)}+3 \mathbf{B}_{10}^{(2)},
               \mathbf{B}_4^{(1)}+3 \mathbf{B}_4^{(2)}
\rangle.\endaligned\eqno{(2.4.90)}$$
Then $V_7$ is a seven-dimensional subspace which is stable under the
action of $G \cong \text{SL}(2, \mathbb{F}_{13})$. Moreover,
$$\aligned
  &13 S(4 \mathbf{B}_0^{(0)}-\mathbf{B}_0^{(1)}+\mathbf{B}_0^{(2)})\\
 =&-(4 \mathbf{B}_0^{(0)}-\mathbf{B}_0^{(1)}+\mathbf{B}_0^{(2)})+14 \mathbf{B}_5
   +14 \mathbf{B}_2+14 \mathbf{B}_6+14 \mathbf{B}_8+14 \mathbf{B}_{11}+14 \mathbf{B}_7\\
  &+7 (\mathbf{B}_1^{(1)}+\mathbf{B}_1^{(2)})
   +7 (\mathbf{B}_3^{(1)}+\mathbf{B}_3^{(2)})
   +7 (\mathbf{B}_9^{(1)}+\mathbf{B}_9^{(2)})\\
  &+7 (-\mathbf{B}_{12}^{(1)}+\mathbf{B}_{12}^{(2)})
   +7 (-\mathbf{B}_{10}^{(1)}+\mathbf{B}_{10}^{(2)})
   +7 (-\mathbf{B}_4^{(1)}+\mathbf{B}_4^{(2)}).
\endaligned\eqno{(2.4.91)}$$
$$T(4 \mathbf{B}_0^{(0)}-\mathbf{B}_0^{(1)}+\mathbf{B}_0^{(2)})
 =4 \mathbf{B}_0^{(0)}-\mathbf{B}_0^{(1)}+\mathbf{B}_0^{(2)}.\eqno{(2.4.92)}$$
Without loss of generality, we can only consider the action of $S$
and $T$ on $\mathbf{B}_5$, $\mathbf{B}_8$, $\mathbf{B}_1^{(1)}+\mathbf{B}_1^{(2)}$,
$-\mathbf{B}_{12}^{(1)}+\mathbf{B}_{12}^{(2)}$. The action of $S$
and $T$ on $\mathbf{B}_2$, $\mathbf{B}_6$, $\mathbf{B}_{11}$, $\mathbf{B}_7$,
$\mathbf{B}_3^{(1)}+\mathbf{B}_3^{(2)}$, $\mathbf{B}_9^{(1)}+\mathbf{B}_9^{(2)}$,
$-\mathbf{B}_{10}^{(1)}+\mathbf{B}_{10}^{(2)}$ and $-\mathbf{B}_4^{(1)}+\mathbf{B}_4^{(2)}$
can be done similarly. We have
$$\aligned
  &13 S(\mathbf{B}_5)\\
 =&4 \mathbf{B}_0^{(0)}-\mathbf{B}_0^{(1)}+\mathbf{B}_0^{(2)}+\\
  &-[(\zeta^3+\zeta^{10})+2 (\zeta^2+\zeta^{11})+2 (\zeta^9+\zeta^4)
   +2 (\zeta+\zeta^{12})] \mathbf{B}_5+\\
  &-[(\zeta+\zeta^{12})+2 (\zeta^5+\zeta^8)+2 (\zeta^3+\zeta^{10})
   +2 (\zeta^9+\zeta^4)] \mathbf{B}_2+\\
  &-[(\zeta^9+\zeta^4)+2 (\zeta^6+\zeta^7)+2 (\zeta+\zeta^{12})
   +2 (\zeta^3+\zeta^{10})] \mathbf{B}_6+\\
  &-[(\zeta^2+\zeta^{11})+2 (\zeta^3+\zeta^{10})+2 (\zeta^6+\zeta^7)
   +2 (\zeta^5+\zeta^8)] \mathbf{B}_8+\\
  &-[(\zeta^5+\zeta^8)+2 (\zeta+\zeta^{12})+2 (\zeta^2+\zeta^{11})
   +2 (\zeta^6+\zeta^7)] \mathbf{B}_{11}+\\
  &-[(\zeta^6+\zeta^7)+2 (\zeta^9+\zeta^4)+2 (\zeta^5+\zeta^8)
   +2 (\zeta^2+\zeta^{11})] \mathbf{B}_7+\\
  &+(\zeta^9+\zeta^4+\zeta^6+\zeta^7+\zeta^2+\zeta^{11})
  (\mathbf{B}_1^{(1)}+\mathbf{B}_1^{(2)})+\\
  &+(\zeta^3+\zeta^{10}+\zeta^2+\zeta^{11}+\zeta^5+\zeta^8)
  (\mathbf{B}_3^{(1)}+\mathbf{B}_3^{(2)})+\\
  &+(\zeta+\zeta^{12}+\zeta^5+\zeta^8+\zeta^6+\zeta^7)
  (\mathbf{B}_9^{(1)}+\mathbf{B}_9^{(2)})+\\
  &+(\zeta^9+\zeta^4+\zeta^6+\zeta^7+\zeta^3+\zeta^{10})
  (-\mathbf{B}_{12}^{(1)}+\mathbf{B}_{12}^{(2)})+\\
  &+(\zeta^3+\zeta^{10}+\zeta^2+\zeta^{11}+\zeta+\zeta^{12})
  (-\mathbf{B}_{10}^{(1)}+\mathbf{B}_{10}^{(2)})+\\
  &+(\zeta+\zeta^{12}+\zeta^5+\zeta^8+\zeta^9+\zeta^4)
  (-\mathbf{B}_4^{(1)}+\mathbf{B}_4^{(2)}).
\endaligned\eqno{(2.4.93)}$$
$$\aligned
  &13 S(\mathbf{B}_8)\\
 =&4 \mathbf{B}_0^{(0)}-\mathbf{B}_0^{(1)}+\mathbf{B}_0^{(2)}+\\
  &-[(\zeta^2+\zeta^{11})+2 (\zeta^3+\zeta^{10})+2 (\zeta^6+\zeta^7)
   +2 (\zeta^5+\zeta^8)] \mathbf{B}_5+\\
  &-[(\zeta^5+\zeta^8)+2 (\zeta+\zeta^{12})+2 (\zeta^2+\zeta^{11})
   +2 (\zeta^6+\zeta^7)] \mathbf{B}_2+\\
  &-[(\zeta^6+\zeta^7)+2 (\zeta^9+\zeta^4)+2 (\zeta^5+\zeta^8)
   +2 (\zeta^2+\zeta^{11})] \mathbf{B}_6+\\
  &-[(\zeta^3+\zeta^{10})+2 (\zeta^2+\zeta^{11})+2 (\zeta^9+\zeta^4)
   +2 (\zeta+\zeta^{12})] \mathbf{B}_8+\\
  &-[(\zeta+\zeta^{12})+2 (\zeta^5+\zeta^8)+2 (\zeta^3+\zeta^{10})
   +2 (\zeta^9+\zeta^4)] \mathbf{B}_{11}+\\
  &-[(\zeta^9+\zeta^4)+2 (\zeta^6+\zeta^7)+2 (\zeta+\zeta^{12})
   +2 (\zeta^3+\zeta^{10})] \mathbf{B}_7+\\
  &+(\zeta^9+\zeta^4+\zeta^6+\zeta^7+\zeta^3+\zeta^{10})
   (\mathbf{B}_1^{(1)}+\mathbf{B}_1^{(2)})+\\
  &+(\zeta^3+\zeta^{10}+\zeta^2+\zeta^{11}+\zeta+\zeta^{12})
   (\mathbf{B}_3^{(1)}+\mathbf{B}_3^{(2)})+\\
  &+(\zeta+\zeta^{12}+\zeta^5+\zeta^8+\zeta^9+\zeta^4)
   (\mathbf{B}_9^{(1)}+\mathbf{B}_9^{(2)})+\\
  &+(\zeta^9+\zeta^4+\zeta^6+\zeta^7+\zeta^2+\zeta^{11})
   (-\mathbf{B}_{12}^{(1)}+\mathbf{B}_{12}^{(2)})+\\
  &+(\zeta^3+\zeta^{10}+\zeta^2+\zeta^{11}+\zeta^5+\zeta^8)
   (-\mathbf{B}_{10}^{(1)}+\mathbf{B}_{10}^{(2)})+\\
  &+(\zeta+\zeta^{12}+\zeta^5+\zeta^8+\zeta^6+\zeta^7)
   (-\mathbf{B}_4^{(1)}+\mathbf{B}_4^{(2)}).
\endaligned\eqno{(2.4.94)}$$
$$\aligned
  &13 S(\mathbf{B}_1^{(1)}+\mathbf{B}_1^{(2)})\\
 =&2 (4 \mathbf{B}_0^{(0)}-\mathbf{B}_0^{(1)}+\mathbf{B}_0^{(2)})+\\
  &+4 (\zeta^9+\zeta^4+\zeta^6+\zeta^7+\zeta^2+\zeta^{11}) \mathbf{B}_5+\\
  &+4 (\zeta^3+\zeta^{10}+\zeta^2+\zeta^{11}+\zeta^5+\zeta^8) \mathbf{B}_2+\\
  &+4 (\zeta+\zeta^{12}+\zeta^5+\zeta^8+\zeta^6+\zeta^7) \mathbf{B}_6+\\
  &+4 (\zeta^9+\zeta^4+\zeta^6+\zeta^7+\zeta^3+\zeta^{10}) \mathbf{B}_8+\\
  &+4 (\zeta^3+\zeta^{10}+\zeta^2+\zeta^{11}+\zeta+\zeta^{12}) \mathbf{B}_{11}+\\
  &+4 (\zeta+\zeta^{12}+\zeta^5+\zeta^8+\zeta^9+\zeta^4) \mathbf{B}_7]+\\
  &+[2 (1+\zeta^9+\zeta^4)+2 (\zeta+\zeta^{12})+(\zeta^2+\zeta^{11})
   (\mathbf{B}_1^{(1)}+\mathbf{B}_1^{(2)})+\\
  &+[2 (1+\zeta^3+\zeta^{10})+2 (\zeta^9+\zeta^4)+(\zeta^5+\zeta^8)]
   (\mathbf{B}_3^{(1)}+\mathbf{B}_3^{(2)})+\\
  &+[2 (1+\zeta+\zeta^{12})+2 (\zeta^3+\zeta^{10})+(\zeta^6+\zeta^7)]
   (\mathbf{B}_9^{(1)}+\mathbf{B}_9^{(2)})+\\
  &+[2 (1+\zeta^6+\zeta^7)+2 (\zeta^5+\zeta^8)+(\zeta^3+\zeta^{10})]
   (-\mathbf{B}_{12}^{(1)}+\mathbf{B}_{12}^{(2)})+\\
  &+[2 (1+\zeta^2+\zeta^{11})+2 (\zeta^6+\zeta^7)+(\zeta+\zeta^{12})]
   (-\mathbf{B}_{10}^{(1)}+\mathbf{B}_{10}^{(2)})+\\
  &+[2 (1+\zeta^5+\zeta^8)+2 (\zeta^2+\zeta^{11})+(\zeta^9+\zeta^4)]
   (-\mathbf{B}_4^{(1)}+\mathbf{B}_4^{(2)}).
\endaligned\eqno{(2.4.95)}$$
$$\aligned
  &13 S(-\mathbf{B}_{12}^{(1)}+\mathbf{B}_{12}^{(2)})\\
 =&2 (4 \mathbf{B}_0^{(0)}-\mathbf{B}_0^{(1)}+\mathbf{B}_0^{(2)})+\\
  &+4 (\zeta^9+\zeta^4+\zeta^6+\zeta^7+\zeta^3+\zeta^{10}) \mathbf{B}_5+\\
  &+4 (\zeta^3+\zeta^{10}+\zeta^2+\zeta^{11}+\zeta+\zeta^{12}) \mathbf{B}_2+\\
  &+4 (\zeta+\zeta^{12}+\zeta^5+\zeta^8+\zeta^9+\zeta^4) \mathbf{B}_6+\\
  &+4 (\zeta^9+\zeta^4+\zeta^6+\zeta^7+\zeta^2+\zeta^{11}) \mathbf{B}_8+\\
  &+4 (\zeta^3+\zeta^{10}+\zeta^2+\zeta^{11}+\zeta^5+\zeta^8) \mathbf{B}_{11}+\\
  &+4 (\zeta+\zeta^{12}+\zeta^5+\zeta^8+\zeta^6+\zeta^7) \mathbf{B}_7+\\
  &+[2 (1+\zeta^6+\zeta^7)+2 (\zeta^5+\zeta^8)+(\zeta^3+\zeta^{10})]
   (\mathbf{B}_1^{(1)}+\mathbf{B}_1^{(2)})+\\
  &+[2 (1+\zeta^2+\zeta^{11})+2 (\zeta^6+\zeta^7)+(\zeta+\zeta^{12})]
   (\mathbf{B}_3^{(1)}+\mathbf{B}_3^{(2)})+\\
  &+[2 (1+\zeta^5+\zeta^8)+2 (\zeta^2+\zeta^{11})+(\zeta^9+\zeta^4)]
   (\mathbf{B}_9^{(1)}+\mathbf{B}_9^{(2)})+\\
  &+[2 (1+\zeta^9+\zeta^4)+2 (\zeta+\zeta^{12})+(\zeta^2+\zeta^{11})]
   (-\mathbf{B}_{12}^{(1)}+\mathbf{B}_{12}^{(2)})+\\
  &+[2 (1+\zeta^3+\zeta^{10})+2 (\zeta^9+\zeta^4)+(\zeta^5+\zeta^8)]
   (-\mathbf{B}_{10}^{(1)}+\mathbf{B}_{10}^{(2)})+\\
  &+[2 (1+\zeta+\zeta^{12})+2 (\zeta^3+\zeta^{10})+(\zeta^6+\zeta^7)]
   (-\mathbf{B}_4^{(1)}+\mathbf{B}_4^{(2)}).
\endaligned\eqno{(2.4.96)}$$
$$T(\mathbf{B}_5)=\zeta^5 \mathbf{B}_5, \quad
  T(\mathbf{B}_8)=\zeta^8 \mathbf{B}_8.\eqno{(2.4.97)}$$
$$T(\mathbf{B}_1^{(1)}+\mathbf{B}_1^{(2)})
 =\zeta (\mathbf{B}_1^{(1)}+\mathbf{B}_1^{(2)}).\eqno{(2.4.98)}$$
$$T(-\mathbf{B}_{12}^{(1)}+\mathbf{B}_{12}^{(2)})
 =\zeta^{12} (-\mathbf{B}_{12}^{(1)}+\mathbf{B}_{12}^{(2)}).\eqno{(2.4.99)}$$
Let
$$\aligned
  V_{13}=\langle &4 \mathbf{B}_0^{(0)}-\mathbf{B}_0^{(1)}+\mathbf{B}_0^{(2)},
               \mathbf{B}_5, \mathbf{B}_2, \mathbf{B}_6, \mathbf{B}_8,
               \mathbf{B}_{11}, \mathbf{B}_7,\\
              &\mathbf{B}_1^{(1)}+\mathbf{B}_1^{(2)},
               \mathbf{B}_3^{(1)}+\mathbf{B}_3^{(2)},
               \mathbf{B}_9^{(1)}+\mathbf{B}_9^{(2)},\\
              &-\mathbf{B}_{12}^{(1)}+\mathbf{B}_{12}^{(2)},
               -\mathbf{B}_{10}^{(1)}+\mathbf{B}_{10}^{(2)},
               -\mathbf{B}_4^{(1)}+\mathbf{B}_4^{(2)} \rangle.
\endaligned\eqno{(2.4.100)}$$
Then $V_{13}$ is a thirteen-dimensional subspace which is stable under
the action of $G \cong \text{SL}(2, \mathbb{F}_{13})$. It is easy to see
that each $v \in V$ can be represented by a linear combination of the basis
of $V_1 \oplus V_7 \oplus V_{13}$. This gives the following decomposition
of $V$:
$$V=V_1 \oplus V_7 \oplus V_{13},\eqno{(2.4.101)}$$
i.e.,
$$\mathbf{21}=\mathbf{1} \oplus \mathbf{7} \oplus \mathbf{13}.\eqno{(2.4.102)}$$
This completes the proof of Theorem 2.4.16.

\noindent
$\qquad \qquad \qquad \qquad \qquad \qquad \qquad \qquad \qquad
 \qquad \qquad \qquad \qquad \qquad \qquad \qquad \boxed{}$

   Let $\mathfrak{a}_i$ be the ideals generated by the basis of $V_i$
$(i=1, 7, 13)$ and $Y_i$ be the varieties corresponding to the ideals
$\mathfrak{a}_i$. Theorem 2.4.14 implies that $\mathfrak{a}_i$ $(i=1, 7, 13)$
are $\text{SL}(2, \mathbb{F}_{13})$-invariant ideals and $Y_i$ $(i=1, 7, 13)$
are $\text{SL}(2, \mathbb{F}_{13})$-invariant varieties.

\textbf{Corollary 2.4.17.} {\it The ideal $I(Y)$ can be decomposed as the sum of
ideals $\mathfrak{a}_i$ $(i=1, 7, 13)$. The curve $Y$ can be expressed as the
intersection of $Y_1$, $Y_7$ and $Y_{13}$:
$$I(Y)=\mathfrak{a}_1+\mathfrak{a}_7+\mathfrak{a}_{13}.\eqno{(2.4.103)}$$
$$Y=Y_1 \cap Y_7 \cap Y_{13}.\eqno{(2.4.104)}$$}

  This gives the proof of Corollary 1.2.

\noindent
$\qquad \qquad \qquad \qquad \qquad \qquad \qquad \qquad \qquad
 \qquad \qquad \qquad \qquad \qquad \qquad \qquad \boxed{}$

\begin{center}
{\large\bf 2.5. A comparison between Hecke's decomposition formulas
                and our decomposition formulas, i.e., motives versus
                anabelian geometry}
\end{center}

  According to Grothendieck (see \cite{Groth}), the starting point of
anabelian algebraic geometry is exactly a study of the action of absolute
Galois groups (particularly the groups $\text{Gal}(\overline{K}/K)$, where
$K$ is an extension of finite type of the prime field) on geometric
fundamental groups of algebraic varieties defined over $K$, and more
particularly fundamental groups which are very far from abelian groups
and which for this reason he calls anabelian. The central theme of
anabelian algebraic geometry is to reconstitute certain so-called
anabelian varieties $X$ over an absolute field $K$ from their mixed
fundamental group, the extension of $\text{Gal}(\overline{K}/K)$ by
$\pi_1(X_{\overline{K}})$. This is when Grothendieck discovered the
fundamental conjecture of anabelian algebraic geometry, close to the
conjectures of Mordell and Tate which have been proved by Faltings
(\cite{Faltings})

\textbf{Conjecture 2.5.1. (Anabelian Tate conjecture).} {\it If $C$, $C^{\prime}$
are two hyperbolic curves over a field $k$ finitely generated over the rationals,
then the natural mapping}
$$\Phi_{(C, C^{\prime})}: \text{Hom}_k(C, C^{\prime}) \longrightarrow
  \text{Hom}_{\text{Gal}(\overline{k}/k)}(\pi_1(C), \pi_1(C^{\prime}))/
  \text{Int} \pi_1(\overline{C^{\prime}})$$
{\it gives a one-to-one correspondence between the dominant $k$-morphisms
and the classes of} $\text{Gal}(\overline{k}/k)$-{\it compatible open
homomorphisms. Here}, ``$/\text{Int} \pi_1$'' {\it means the quotient by
the natural actions by the inner automorphisms of $\pi_1$}.

  This Anabelian Tate conjecture is closely related to the Tate conjecture
which can also be regarded as an anabelian result for abelian varieties-where
the fundamental groups are replaced by their abelian counterparts, i.e.,
$\ell$-adic Tate modules $T_{\ell}$:

\textbf{Theorem 2.5.2. (Tate conjecture (Faltings' theorem).} (see \cite{Tate},
\cite{Tate1966} and \cite{Faltings}). {\it For $A$ and $B$ abelian varieties
over a number field $k$, the natural map
$$\text{Hom}_k(A, B) \otimes_{\mathbb{Z}} \mathbb{Z}_{\ell} \longrightarrow
  \text{Hom}_{\text{Gal}(\overline{k}/k)}(T_{\ell}(A), T_{\ell}(B))$$
is an isomorphism.}

  As Grothendieck himself observed, this conjecture bears some resemblance
to the above Tate conjecture which can be reformulated by one-dimensional
\'{e}tale homology groups of abelian varieties:
$$\text{Hom}_k(A, B) \otimes_{\mathbb{Z}} \widehat{\mathbb{Z}} \cong
  \text{Hom}_{\text{Gal}(\overline{k}/k)}(H_1(\overline{A}, \widehat{\mathbb{Z}}),
  H_1(\overline{B}, \widehat{\mathbb{Z}})).$$
Here, $A$ and $B$ are abelian varieties defined over a global field $k$,
and $\widehat{\mathbb{Z}}$ is the profinite completion of $\mathbb{Z}$.
Note that $H_1$ is just the abelianization of $\pi_1$. In the case of
abelian varieties, the above exact sequence becomes that
$$1 \longrightarrow H_1(\overline{A}) \longrightarrow \pi_1(A)
  \longrightarrow \text{Gal}(\overline{k}/k) \longrightarrow 1,$$
where $\pi_1(\overline{A})$ can be replaced by the $\ell$-adic
component, namely the Tate module $T_{\ell}(\overline{A})$.

  If $J$ is the Jacobian variety of the smooth compactification
of the curve $C$, then the one-dimensional homology group
$H_1(\overline{J})$ appears as a natural abelian quotient of
$\pi_1(\overline{C})$. The Tate conjecture asserts that for
Jacobian varieties $J$, $J^{\prime}$,
$$\text{Hom}_k(J, J^{\prime}) \otimes \mathbb{Z}_{\ell} \cong
  \text{Hom}_{\text{Gal}(\overline{k}/k)}(H_1(\overline{J}),
  H_1(\overline{J^{\prime}})).$$
So, the Anabelian Tate conjecture can be regarded as a non-abelian
lifting of the Tate conjecture.

  Research on Grothendieck's conjecture was begun at the end of the
1980s by Nakamura, given significant impetus (including the case of
positive characteristic) by Tamagawa, and brought to a final solution
by means of a new ($p$-adic) interpretation of the problem due to
Mochizuki (see \cite{NTM}).

  We will give the following results which provide a comparison
between modular curves with higher level, i.e., hyperbolic curves
(anabelian algebraic geometry) and their Jacobian varieties (motives).
Hence, it can be regarded as a counterpart of Grothendieck's anabelian
Tate conjecture and Tate conjecture for hyperbolic curves and their
Jacobian varieties, respectively. It can also be considered as a
counterpart of the comparison theorem between $p$-adic \'{e}tale
cohomology and de Rham cohomology (crystalline cohomology). Recall
that the $p$-adic comparison theorems (or the $p$-adic periods
isomorphisms) are isomorphisms, analog to the complex periods isomorphism
$$H_{\text{dR}}^i(X/\mathbb{C}) \cong H^i(X(\mathbb{C}), \mathbb{Q}) \otimes \mathbb{C}$$
for a smooth and projective variety over $\mathbb{C}$, between the $p$-adic
cohomology and the de Rham cohomology (plus some additional structure) of
smooth and projective varieties over a finite extension of $\mathbb{Q}_p$.
The theory started with the work of Grothendieck's ``mysterious period functor''
relating directly $p$-adic \'{e}tale cohomology and de Rham cohomology as
well as Tate and Fontaine on abelian varieties and $p$-divisible groups,
and continues with the work of Fontaine, Bloch, Kato, Messing, Faltings,
Hyodo, Tsuji, Breuil and others by various methods: syntomic topology,
de Rham-Witt complexes and the method of almost \'{e}tale extensions. In
particular, this comparison leads to the $p$-adic Hodge theory.

  Let $K$ be a local field of characteristic $0$ and residue characteristic
$p$, i.e., $K$ is a $p$-adic field. Let $\overline{K}$ be an algebraic
closure of $K$. The first significant result towards Grothendieck's
question on the mysterious functor is due to Tate and appears in his
paper on $p$-divisible groups (see \cite{Tate1967}), i.e., the Hodge-Tate
decomposition, which is the starting point of $p$-adic Hodge theory.
It states that, when $A$ is a smooth abelian scheme over
$\text{Spec}(\mathbb{Z}_p)$, we have a
$\text{Gal}(\overline{\mathbb{Q}}_p/\mathbb{Q}_p)$-equivariant isomorphism:
$$H_{\text{\'{e}t}}^1(A_{\overline{\mathbb{Q}}_p}, \mathbb{Q}_p) \otimes_{\mathbb{Q}_p}
  \mathbb{C}_p \simeq (H^1(A, \mathcal{O}_A) \otimes_K \mathbb{C}_p) \oplus
  (H^0(A, \Omega_{A/K}^1) \otimes_K \mathbb{C}_p(-1)),$$
where $\mathbb{C}_p$ denotes the completion of $\overline{\mathbb{Q}}_p$ and
$\mathbb{C}_p(-1)$ is its twist by the inverse of the cyclotomic character.
This also works for a proper smooth curve over $K$ by considering its Jacobian.
We have the following comparison:

\noindent (1) Hecke's decomposition for $J(p)$ in the following (2.5.2):
              $\text{PSL}(2, \mathbb{F}_p)$-equivariant.

\noindent (2) Hodge-Tate decomposition for Jacobian:
              $\text{Gal}(\overline{\mathbb{Q}}_p/\mathbb{Q}_p)$-equivariant.

\noindent (3) Our decomposition for $I(\mathcal{L}(X(p)))$ in the following (2.5.2):
              $\text{PSL}(2, \mathbb{F}_p) \times
              \text{Gal}(\overline{\mathbb{Q}}/\mathbb{Q})$-equivariant.

\noindent Here, (1) and (2) come from motives, while (3) comes from
anabelian algebraic geometry.

  It is known that the Tate conjecture is a central problem in arithmetical
algebraic geometry, which predicts a description of all subvarieties of a given
algebraic variety over a field $k$ in terms of the representations of the Galois
group of $k$ on \'{e}tale cohomology in the context of motives. In contrast, we
have the following comparison:

$$\begin{array}{ccc}
  & \text{Galois actions on} & \\
  & \text{algebraic cycles} & \\
  \swarrow &  & \searrow\\
  \text{in motives} & \longleftrightarrow & \text{in anabelian}\\
                    &  & \text{algebraic geometry}\\
  \text{Tate conjecture for} & \longleftrightarrow & \text{Grothendieck's anabelian}\\
  \text{Jacobians of curves} &  & \text{Tate conjecture}\\
  \text{Tate conjecture} & \longleftrightarrow & \text{Our Theorem 1.1}\\
  \text{$H^{2d}(X)_{\ell}$} & \longleftrightarrow & \text{$I(\mathcal{L}(X(p)))$}\\
  \text{linear combination of} & \longleftrightarrow & \text{intersection of}\\
  \text{cohomology classes} &  & \text{invariant ideals}\\
  \text{Jacobians of curves} & \longleftrightarrow & \text{curves}\\
  \text{semisimplicity for} & \longleftrightarrow & \text{semisimplicity for}\\
  \text{$T_{\ell}(A)$} &  & \text{$V_{I(\mathcal{L}(X(p)))}$}\\
\end{array}\eqno{(2.5.1)}$$
In particular, the decomposition of a defining ideal as a representation
into irreducible representations provides a relation between the defining
ideals and the irreducible algebraic cycles in the context of anabelian
algebraic geometry, which can be regarded as a counterpart of the
semi-simplicity conjecture in the context of motives.

  In the context of motives, let us consider algebraic cycles of
codimension $d$ in a smooth projective variety $X$ over a field $k$ which
is of finite type over the prime field. Let $\ell$ be a prime number not
dividing the characteristic of $k$, and let $H^{2d}(X)_{\ell}$ be the
\'{e}tale cohomology group of $X/\overline{k}$ with coefficients in
$\mathbb{Q}_{\ell}(d)$. Any codimension $d$ cycle defines an element
$c(z)$ of $H^{2d}(X)_{\ell}$ which is invariant under the Galois group
$\text{Gal}(\overline{k}/k)$. By considering number-theoretic analogue
of the classical conjectures of Hodge concerning algebraic cycles, the
Tate conjecture (see \cite{Tate}) states that, conversely, any element
of $H^{2d}(X)_{\ell}$ invariant under $\text{Gal}(\overline{k}/k)$ is a
$\mathbb{Q}_{\ell}$-linear combination of classes of type $c(z)$. In
particular, in the case of abelian varieties as well as Jacobian varieties
of curves, Tate formulated his conjecture which links the problem of
classifying abelian varieties up to isogeny to that of classifying the
Galois representations to which they give rise (see \cite{Mazur1987}):
Let $\ell$ be a prime number. Let $K$ be a number field, and $\overline{K}$
an algebraic closure of $K$. An abelian variety $A$ over $K$ is determined
up to isogeny over $K$ by the natural representation of $\text{Gal}(\overline{K}/K)$
on the Tate module
$$T_{\ell}(A)=\varprojlim_{n} A[\ell^n](\overline{K}).$$
An important adjunct to the conjecture of Tate concerning the representations
of the Galois group $\text{Gal}(\overline{K}/K)$ on the Tate module of an
abelian variety $A$ defined over a number field $K$ is the assertion of
semisimplicity of the representation of $\text{Gal}(\overline{K}/K)$ on
$T_{\ell}(A)$.

\textbf{Theorem 2.5.3. (Semisimplicity Theorem).} (see \cite{Tate} and
\cite{Faltings}) {\it Let $A$ be an abelian variety over a number field $K$.
Then $T_{\ell}(A) \otimes_{\mathbb{Z}_{\ell}} \mathbb{Q}_{\ell}$ is semisimple
as a representation of} $\text{Gal}(\overline{K}/K)$ {\it $($that is, it is
a direct sum of irreducible representations$)$.}

  In the context of anabelian algebraic geometry, our Theorem 1.1 states
that the defining ideal $I(\mathcal{L}(X(p)))$ which is invariant under
the action of both $\text{Gal}(\overline{\mathbb{Q}}/\mathbb{Q})$ and
$\text{PSL}(2, \mathbb{F}_p)$ can de represented as an intersection of
both $\text{Gal}(\overline{\mathbb{Q}}/\mathbb{Q})$ and
$\text{PSL}(2, \mathbb{F}_p)$ invariant ideals, which are parameterized
by some irreducible representations of $\text{PSL}(2, \mathbb{F}_p)$
defined over $\mathbb{Q}(\zeta_p)$.

  We find geometric and arithmetic properties of modular curves $X(p)$ that
are not determined by their Jacobian $J(p):=\text{Jac}(X(p))$, which is
the motive of the curve $X(p)$. In fact, by Serre duality theorem, we have
$$H^1(X(p), \mathcal{O}_{X(p)}) \cong H^0(X(p), \Omega_{X(p)}^1),$$
where $H^0(X(p), \Omega_{X(p)}^1)$ denotes the complex vector space of
global sections of the sheaf of holomorphic differentials on the Riemann
surface $X(p)$. In fact, this complex vector space is the same as the
complex vector space of global sections of the sheaf of holomorphic
differentials on its Jacobian variety $J(p)$. The following theorem
is fundamental.

\textbf{Theorem 2.5.4.} (see \cite{Saito2014}, Theorem D.3) {\it Let $S$ be a scheme,
and let $f:X \rightarrow S$ be a proper smooth curve over $S$ of genus $g$
such that each geometric fiber is connected.}

(1) $\text{Pic}_{X/S}^{0}$ {\it is represented by an abelian scheme $j: J \rightarrow S$
over $S$ of relative dimension $g$.}

(2) {\it There is a natural isomorphism}
$$j_{*} \Omega_{J/S} \rightarrow f_{*} \Omega_{X/S}.\eqno{(2.5.2)}$$

(3) {\it Let $n \geq 1$ be an integer. The Weil pairing defines a bilinear
form $J[n] \times J[n] \rightarrow \mu_n$ and defines an isomorphism to the
Cartier dual}
$$J[n] \rightarrow J[n]^{\vee}.$$

  The moduli space $J$ of the functor $\text{Pic}_{X/S}^{0}$ is called the
Jacobian of $X$.

\textbf{Corollary 2.5.5.} (see \cite{Saito2014}, Corollary D.5) {\it Let $k$ be
a field, and let $X$ be a proper smooth curve over $k$ of genus $g$ such that
the geometric fiber $X_{\overline{k}}$ is connected.}

(1) {\it The Jacobian} $J=\text{Pic}_{X/k}^{0}$ {\it is an abelian variety over
$k$ of dimension $g$.}

(2) {\it There is a natural isomorphism}
$$\Gamma(J, \Omega_{J/k}) \rightarrow \Gamma(X, \Omega_{X/k}).\eqno{(2.5.3)}$$

(3) {\it Let $n \geq 1$ be an integer. The Weil pairing defines a bilinear
form $J[n] \times J[n] \rightarrow \mu_n$ and the isomorphism $J[n] \rightarrow
J[n]^{\vee}$ to its Cartier dual.}

  By the above natural isomorphism (2.5.3), we have
$$\Gamma(J(p), \Omega_{J(p)}^1) \cong \Gamma(X(p), \Omega_{X(p)}^1).\eqno{(2.5.4)}$$

  It is well-known that the concept of genus was introduced by Abel and Riemann
in two versions. This fundamental concept arises in problems of integration,
through algebraic curves and their associated Riemann surfaces. Its importance
in analysis, algebraic geometry, number theory and topology is emphasized through
many theorems. In fact, the genus is closely related to the Jacobian variety of
algebraic curves. As a consequence of the following comparison between Hecke's
decomposition and our decomposition (i.e., the first comparison between motives
and anabelian algebraic geometry), it gives an anabelian counterpart of the
genus for $X(p)$ which is the dimension of its Jacobian variety $J(p)$, by the
dimension of defining ideals $\dim I(\mathcal{L}(X(p)))$:
$$\begin{matrix}
  &\textbf{A comparison between}\\
  &\textbf{Hecke's decomposition and our decomposition}
\end{matrix}$$
$$\begin{matrix}
            & \text{motives} & \text{anabelian algebraic geometry}\\
            & \text{Jacobian varieties} & \text{defining ideals}\\
            & \text{(Abel, Riemann)} & \text{(Weierstrass, Klein)}\\
            & \text{abelian part of $X(p)$: $H_1$} &
              \text{anabelian part of $X(p)$: $\pi_1$}\\
            &\text{Hecke's level theory} &\text{our level theory}\\
            & \text{Hecke's decomposition} & \text{our decomposition}\\
            & \dim H^0(X(p), \Omega_{X(p)}^1) & \dim I(\mathcal{L}(X(p)))\\
            & \text{genus of $X(p)$, i.e.,} & \text{dimension of}\\
            & \text{dimension of $J(p)$} & \text{defining ideals}\\
       p=7  & H^0(X(7), \Omega_{X(7)}^1)=V_3 & V_{I(\mathcal{L}(X(7)))}=\mathbf{1}\\
            & \text{degenerate discrete series} & \text{trivial representation}\\
            & J(7) \cong A_3 & I(\mathcal{L}(X(7)))=I_1\\
      p=11  & H^0(X(11), \Omega_{X(11)}^1)=V_5 \oplus V_{10} \oplus V_{11} &
              V_{I(\mathcal{L}(X(11)))}=\mathbf{10}\\
            & \text{degenerate discrete series,} & \text{discrete series}\\
            & \text{discrete series corresponding to $\chi_4$} & \text{corresponding}\\
            & \text{or $\chi_5$ and Steinberg representation} & \text{to $\chi_5$}\\
            & J(11) \cong A_5 \times A_{10} \times A_{11} & I(\mathcal{L}(X(11)))=I_{10}\\
      p=13  & H^0(X(13), \Omega_{X(13)}^1) & V_{I(\mathcal{L}(X(13)))}=\mathbf{1}
              \oplus \mathbf{7} \oplus \mathbf{13}\\
            & =V_{12}^{(1)} \oplus V_{12}^{(2)} \oplus V_{12}^{(3)} \oplus V_{14} & \\
            & \text{three conjugate} & \text{trivial representation,}\\
            & \text{discrete series} & \text{degenerate principal series}\\
            & \text{and principal series} & \text{and Steinberg representation}\\
            & J(13) \cong A_{12}^{(1)} \times A_{12}^{(2)} \times A_{12}^{(3)} \times A_{14}
            & I(\mathcal{L}(X(13)))=I_{1} \cap I_7 \cap I_{13}
\end{matrix}\eqno{(2.5.5)}$$

  In particular, the cases $p=7$ or $11$ are especially interesting, because
they are closely related to the last mathematical testament of Galois \cite{Ga},
and these two cases are the most important examples in Klein-Fricke's lectures
on the elliptic modular functions which go beyond the icosahedron (see \cite{KF1},
Section three, Chapter seven and \cite{KF2}, Section five, Chapter five). Now,
let us give some background. The transformation equations of elliptic functions
were started by Euler, Landen, Lagrange, Gauss and Legendre, and then by Abel,
Jacobi and Galois, i.e., his final letter on May 29, 1832. Galois' idea stimulated
the work of Hermite, Kronecker, and in particular, Klein, who combined ideas of
Riemann and Galois, which leads to the theory of modular curves $X(p)$. We will
give a modern framework and interpretation by Galois representations arising from
anabelian algebraic geometry which correspond to the Galois groups of modular
equations which appear in the transformation theory of elliptic functions.

  In fact, Abel approached the theory of elliptic functions mainly from
an algebraic point of view, centered on the various algebraic equations
that the transformation theory provided in such abundance. The rich source
of irreducible algebraic equations, in particular, the division equations
that the transformation theory of elliptic functions provided, was instrumental
for both Abel and Galois as suggestive examples that led to the Galois theory.
In particular, in his analysis of the division equations, Abel saw how far
Gauss' method for solving cyclotomic equations could be generalized. For Abel
the transformation theory of elliptic functions can be stated as follows: Find
the conditions under which there exists an algebraic function $y=y(x)$, i.e.,
$R(x, y)=0$ for some rational function $R$ in two variables over $\mathbb{C}$,
which is a solution of the separable differential equation
$$\frac{dy}{\sqrt{(1-y^2)(1-\lambda^2 y^2)}}
 =a \frac{dx}{\sqrt{(1-x^2)(1-\kappa^2 x^2)}},$$
where $a$ is a constant (the multiplier), and $\lambda$ and $\kappa$ are called
moduli. In \cite{Abel2}, a much more general problem is posed, it is proved by
algebraic methods that one can reduce the situation to $y$ being a rational
function of $x$. The solution of the above differential equation is transferred
by Abel to finding all solutions of the equation
$$(1-y^2)(1-\lambda^2 y^2)=t^2 (1-x^2)(1-\kappa^2 x^2),$$
where $y=\frac{U(x)}{V(x)}$ and $t=t(x)$ are rational functions of $x$. Now
the general case can be traced back to the special case where the degrees
of the polynomials $U(x)$ and $V(x)$ are $p$ and $p-1$, respectively, where
$p$ is a prime. We say that the transformation is of degree $p$. In this case
there exists an algebraic equation of degree $p+1$ between $\lambda$ and $\kappa$
with integer coefficients, the so-called modular equation. The multiplier $a$
is determined by this differential equation. Now, we briefly describe a modern
way to look at the transformation theory of elliptic functions, a viewpoint that
was also initiated by Abel (see \cite{Abel2}). Instead of considering the above
differential equation, we consider two elliptic functions $\varphi(u)=\varphi(u |
\omega_1, \omega_2)$ and $\psi(u^{\prime})=\psi(u^{\prime} | \omega_1^{\prime},
\omega_2^{\prime})$ with primitives periods $\omega_1$, $\omega_2$ and
$\omega_1^{\prime}$, $\omega_2^{\prime}$, respectively. The transformation
theory is to investigate under which conditions there exists an algebraic
relation between $\varphi(u)$ and $\psi(u^{\prime})$ if $u^{\prime}=m u$ for
some constant $m$. Assuming that both $\varphi$ and $\psi$ are homogeneous
functions in three variables, one can assume that $m=1$ and hence $u^{\prime}=u$.
Moreover, one can reduce the investigation into studying the relation between
$\varphi(u | \omega_1, \omega_2)$ and $\varphi(u | \omega_1^{\prime}, \omega_2^{\prime})$,
where $\omega_1^{\prime}=a \omega_1+b \omega_2$, $\omega_2^{\prime}=c \omega_1+
d \omega_2$, and $a, b, c, d \in \mathbb{Z}$, $ad-bc=n>0$. We say that
$\varphi(u | \omega_1^{\prime}, \omega_2^{\prime})$ arises from
$\varphi(u | \omega_1, \omega_2)$ by a transformation of degree $n$. The
significance comes from the close relation between transformation theory
and division equations. That is, as for the modular equation associated to
transformations of prime degree $p$ of elliptic functions, it leads to
essentially the same problem as for the division equation. Note that the
special division equations are closely related to the torsion points of an
elliptic curve, where the elliptic curve is associated to an elliptic
function by the Weierstrass parametrization. This is intimately related
to the multiplication maps of an elliptic curve to itself, which are
important for the study of the arithmetic of elliptic curves. Similarly,
the transformation theory of elliptic functions plays much more important
role in the Galois theory. In particular, there are two different kinds of
viewpoints, one is the algebraic version, i.e., algebraic (rational)
solutions of the above differential equation, the other is the transcendental
version, i.e, the theory of elliptic modular functions. This leads to the
last mathematical testament of Galois (see \cite{Ga}), who put together
these two distinct kinds of theories, i.e., the Galois theory of algebraic
equations as well as the Galois theory of elliptic modular functions
(i.e., transcendental functions).

  Let us recall some ideas of Galois in \cite{Ga}. The main applications of
Galois' ambiguity theory concern algebraic numbers and algebraic functions.
Half a century after Galois, Klein's work blended Galois with Riemann, which
started a tradition through Klein and Poincar\'{e}'s uniformization theorem
of Riemann surfaces as well as Poincar\'{e}'s view about fundamental groups,
eventually led to Grothendieck's fusion between Galois theory and the theory
of coverings, a cornerstone of arithmetic geometry. But Galois also foresaw
applications of ambiguity theory beyond the domains of algebraic numbers and
algebraic functions. In fact, in the end of his last mathematical testament
\cite{Ga}, Galois pointed out that ``My principal meditations, for some time
now, were directed on the application of the theory of ambiguity to transcendental
analysis. It was to see, a priori, in a relation between transcendental
quantities or functions, what exchanges can be done, what quantities we could
substitute to the given quantities, without changing the relation.'' The most
interesting as well as important example of such a relation between transcendental
quantities is the transformation equation between elliptic functions and its
Galois group. This leads to the modular equation of degree $p$ and its Galois
group. Galois showed that (see \cite{Ga}) ``The last application of the theory
of equations is related to the modular equation of elliptic functions. We show
that the group of the equation which has for roots the sine of the amplitude of
$p^2-1$ divisions of a period is:
$$x_{k, l}, \quad x_{ak+bl, ck+dl};$$
consequently, the corresponding modular equation has for its group
$$x_{\frac{k}{l}}, \quad x_{\frac{ak+bl}{ck+dl}},$$
in which $\frac{k}{l}$ can take the $p+1$ values
$$\infty, 0, 1, 2, \cdots, p-1.$$
Thus, by agreeing that $k$ can be infinity, we can simply write
$$x_k, \quad x_{\frac{ak+b}{ck+d}}.$$
By giving to $a$, $b$, $c$, $d$ all their values, we obtain $(p+1)p(p-1)$
permutations.

\noindent Now this group is decomposable properly into two groups, for
which the substitutions are
$$x_k, \quad x_{\frac{ak+b}{ck+d}},$$
$ad-bc$ being a quadratic residue of $p$.

\noindent The group thus simplified is of
$$(p+1)p \frac{p-1}{2}$$
permutations.

\noindent For $p=7$ we find a group of $(p+1)(p-1)/2$ permutations, where
$$\infty \quad 1 \quad 2 \quad 4$$
are respectively related to
$$0 \quad 3 \quad 6 \quad 5.$$
For $p=11$, the same substitutions take place with the same notations,
$$\infty \quad 1 \quad 3 \quad 4 \quad 5 \quad 9$$
are respectively related to
$$0 \quad 2 \quad 6 \quad 8 \quad 10 \quad 7.$$
Thus, for the case of $p=5$, $7$, $11$, the modular equation is reduced
to degree $p$.''

  From modern viewpoint of group theory, Galois' result is that the group
$\text{PSL}(2, \mathbb{F}_p)$ which acts transitively on the $p+1$ points
of $\mathbb{P}^1(\mathbb{F}_p)$ only acts transitively on $p$ points if
$p=2$, $3$, $5$, $7$, $11$. Geometrically the most interesting cases are
$p=7$ and $p=11$. In particular it can not operate non-trivially on a set
with fewer than $p+1$-elements when $p>11$. This may be expressed as follows
(see \cite{Kostant}). The cyclic group $C_p$ of order $p$ embeds uniquely
in $\text{PSL}(2, \mathbb{F}_p)$, up to conjugacy, as the Sylow $p$-subgroup.
It is the unipotent radical of a Borel subgroup and pulled up to
$\text{SL}(2, \mathbb{F}_p)$ may be taken to be
$$\mathbb{Z}/p \mathbb{Z}=\left\{ \left(\begin{matrix}
  1 & x\\ 0 & 1 \end{matrix}\right): x \in \mathbb{F}_p \right\}.$$
The result of Galois is that if $p>11$, there exists no subgroup of
$\text{PSL}(2, \mathbb{F}_p)$ which is complementary to $C_p$. By this
we mean that there exists no subgroup $F$ such that set theoretically
$\text{PSL}(2, \mathbb{F}_p)=F \cdot C_p$. The notation means that every
element $g$ can be uniquely written $g=fz$ where $f \in F$ and $z \in C_p$.
Implicit in Galois' statement is certainly the knowledge that his statement
is not true for the three cases of a simple $\text{PSL}(2, \mathbb{F}_p)$
where $p \leq 11$, namely $p=5$, $7$, $11$. It is surely a marvelous fact,
that for the three exceptional cases, the groups $F$ which run counter to
Galois' statement are precisely the symmetry groups of the Platonic solids.
Namely one has
$$\aligned
  \text{PSL}(2, \mathbb{F}_5) &=A_4 \cdot C_5,\\
  \text{PSL}(2, \mathbb{F}_7) &=S_4 \cdot C_7,\\
  \text{PSL}(2, \mathbb{F}_{11}) &=A_5 \cdot C_{11}.
\endaligned$$
These exceptional cases merit elaboration. It was taken up by Conway (see
\cite{Conway}) in connection with exceptional properties of the known
simple groups and the discovery of the sporadic simple groups. Conway
examined more closely the cases $p=2$, $3$, $5$, $7$, $11$ in Galois' theorem,
leading to exceptional properties of the groups $\text{PSL}(2, \mathbb{F}_p)$
in these cases. For example, these include isomorphisms
$$\text{PSL}(2, \mathbb{F}_2) \cong S_3, \quad
  \text{PSL}(2, \mathbb{F}_3) \cong A_4, \quad
  \text{PSL}(2, \mathbb{F}_5) \cong A_5,$$
$$\text{PSL}(2, \mathbb{F}_7) \cong \text{PSL}(3, \mathbb{F}_2),$$
and a presentation of the Mathieu simple group $M_{12}$ related to a
presentation of the simple group $\text{PSL}(2, \mathbb{F}_{11})$. Conway
continued with an illustration of the way the exceptional permutation
representations of the small simple groups $\text{PSL}(2, \mathbb{F}_p)$
can arise in other situations by considering the permutation representation
of Janko's sporadic simple group $J_1$ on $266$ letters. The stabilizer of
a point is an $\text{PSL}(2, \mathbb{F}_{11})$, with orbits of sizes $1$,
$11$, $110$, $132$, $12$. Knowing this much, it is easy to construct $J_1$.

  Now, we will give a geometric viewpoint. As a Riemann surface, the Klein
quartic curve is the first Hurwitz surface, which is uniformised by a
surface subgroup of the $(2, 3, 7)$-triangle group, and has automorphism
group $\text{PSL}(2, \mathbb{F}_7)$ of order $168$. In his paper \cite{K2},
Klein drew the famous main figure (see \cite{K2}, p. 126), showing this
surface as a $14$-sided hyperbolic polygon, tessellated by $336$ triangles
with internal angles $\frac{\pi}{2}$, $\frac{\pi}{3}$ and $\frac{\pi}{7}$.
The following year Klein also wrote about the surface with
$\text{PSL}(2, \mathbb{F}_{11})$ as automorphism group as an image of the
$(2, 3, 11)$-triangle group (see \cite{K3}). This paper introduced some
early ideas about dessins: for example, Klein drew the ten plane trees
(see \cite{K3}, p. 143) which correspond to dessins of type $(3, 2, 11)$
and genus $0$, two of them having monodromy group $\text{PSL}(2, \mathbb{F}_{11})$
in its representation of degree $11$, and the other eight having $A_{11}$
as monodromy group. In particular, these ten plane trees anticipated
Grothendieck's concept of dessins d'enfants (see \cite{Groth}) by over a
century, classifying the ten dessins of degree $11$ and type $(3, 2, 11)$.
More precisely, Klein drew a diagram showing ten plane trees, each with
eleven edges and twelve vertices. In each tree there is a bipartite partition
of the vertices, with five vertices (coloured white) of valencies dividing
$3$, and seven vertices (indicated by short cross-bars, perpendicular to
their incident edges) of valencies dividing $2$; in each case the unique
face is an $11$-gon. Among these plane trees there are four chiral
(mirror-image) pairs, numbered I to IV, and two others, number V and VI,
which exhibit bilateral symmetry. The chiral pair I plays a major role
in Klein's paper \cite{K3}. Klein also considered several quotient curves
and maps, very much in the spirit of some of today's work on dessins:
these include the modular curve $X_0(11)$ of genus $1$, arising from
the action of $\text{PSL}(2, \mathbb{F}_{11})$ on the projective line
$\mathbb{P}^1(\mathbb{F}_{11})$, and a chiral pair of cubic plane trees,
arising from its action on the $11$ cosets of a subgroup isomorphic to
$A_5$ (see \cite{JZ}).

  Now, we will give a unified geometric framework for both $p=7$ and
$p=11$ (see \cite{BJS}). Let us have a brief look at the geometry of
Klein quartic curve as a Riemann surface $S$. The group
$\text{PSL}(2, \mathbb{F}_7)$ contains $28$ $C_3$'s and each $C_3$
fixes two points of $S$. Each $S_4$ in $\text{PSL}(2, \mathbb{F}_7)$
contains $4$ $C_3$'s and together these fix $8$ points which Klein
called $A$, $A^{\prime}$, $B$, $B^{\prime}$, $C$, $C^{\prime}$, $D$,
$D^{\prime}$ which he regarded as being vertices of an inscribed cube.
Each vertex of this cube is the centroid of a triangle whose vertices
have valency $7$. This gives $8 \times 3=24$ vertices of valency $7$.
In this way we see that these $24$ points are vertices of a truncated
cube which are fixed by the $S_4$. Note that as $S_4$ has index $7$ so
that $\text{PSL}(2, \mathbb{F}_7)$ has a transitive permutation
representation on $7$ points. For $p=11$, instead of $(2, 3, 11)$
triangle group used as above, it is more interesting to study
$\text{PSL}(2, \mathbb{F}_{11})$ as an image of the $(2, 5, 11)$ triangle
group because we now get calculations and geometry that are in the spirit
of the work on the Klein quartic curve described above. In the Klein quartic
curve we have automorphism group $\text{PSL}(2, \mathbb{F}_7)$ which contains
an $S_4$ which contains $4$ $C_3$'s each with $2$ fixed points. We then have
$8$ points which form a cube and then we see that the $24$ fixed points of
the elements of order $7$ give a truncated cube stabilized by the $S_4$. Very
similar considerations applied to $\text{PSL}(2, \mathbb{F}_{11})$ as an image
of $\Delta(2, 5, 11)$. Let $L$ be the kernel of the homomorphism from
$\Delta(2, 5, 11)$ to $\text{PSL}(2, \mathbb{F}_{11})$ and let
$\mathcal{B}=\mathbb{H}/L$. Now $\text{PSL}(2, \mathbb{F}_{11})$ contains two
conjugacy classes of subgroups isomorphic to $A_5$. Each $A_5$ contains $6$
elements of order $5$ (each fixing two points of $\mathcal{B}$ and so
together we have $12$ fixed points stabilized by an $A_5$). These $12$ points
then form the vertices of an icosahedron embedded into $\mathcal{B}$. Also just
as there are $24$ fixed points of elements of order $7$ in the Klein surface,
there are $60$ points fixed by elements of order $11$. By considering the
$(2, 5, 11)$ triangles we see that each of these fixed points of a $C_5$ is
the center of a pentagon all of whose vertices is one of the $60$ fixed points
of a $C_{11}$, and thus these $60$ points are the vertices of a truncated
icosahedron or buckyball (see \cite{MS}). This is why we call $\mathcal{B}$
the Buckyball surface.

  An interesting fact about the Klein quartic curve is its relationship to the
Fano plane. This is the projective plane of order two and has $7$ points and
$7$ lines. Its collineation group is $\text{PSL}(2, \mathbb{F}_7)$. It is shown
that there are regular embeddings of the Fano plane into the Klein surface $S$.
This can be described purely algebraically. As Klein observed, there are $14$
subgroups isomorphic to $S_4$ in $\text{PSL}(2, \mathbb{F}_7)$ in two conjugacy
classes of $7$ each. Each $S_4$ corresponds to a truncated cube. Thus we have
two sets $P_1$, $P_2$, $\ldots$, $P_7$ and $L_1$, $L_2$, $\ldots$, $L_7$ of
truncated cubes in the Klein surface. If we regard the $P_i$ as points and
the $L_j$ as lines in a finite geometry and study their intersections then
we get back the Fano plane. Just as the Fano plane, the simplest finite
projective plane has an embedding into the Klein quartic curve, the simplest
biplane has an embedding into the Buckyball surface. In fact, note that there
are two conjugacy classes of $A_5$'s in $\text{PSL}(2, \mathbb{F}_{11})$. We
have two sets $P_1$, $P_2$, $\ldots$, $P_{11}$ and $L_1$, $L_2$, $\ldots$,
$L_{11}$ of buckyballs in $\mathcal{B}$. If we study their intersections we
get back the biplane with $11$ points. Also, we see that this gives
$\text{PSL}(2, \mathbb{F}_{11})$ as a transitive permutation representation
on $11$ points (see \cite{MS} and \cite{BJS}).

  Klein came to his work on elliptic functions, especially modular
functions, which was to assume such a large scope, originally from
the investigations on the icosahedron, with the intention of
understanding the methods of solving the equation of the fifth
degree by elliptic functions, as given in 1858 by Hermite, Kronecker
and Brioschi, from the icosahedron. Then he was thus forced of his
own accord to take a general interest in the transformation theory
of elliptic functions and thus in particular to take up Hermite's
task anew: to really set up in the simplest form the resolvents of
the $5$th, $7$th and $11$th degree, which, according to the famous
theorem by Galois as above, exist for transformations of the $5$th,
$7$th and $11$th order. By combining the group-theoretical and
invariant-theoretical approaches familiar to him with the geometric
and descriptive methods of Riemann's function theory, this was
achieved by Klein in a surprisingly simple way. At the same time, a
complete overview of the various types of modular equations and
multiplier equations that existed in the literature was obtained,
as well as a clear insight into the nature of those algebraic
equations that can be solved by elliptic functions.

  When $p=7$, in the context of motives, the degenerate discrete series
representation of $\text{PSL}(2, \mathbb{F}_7)$ is realizable over
$\mathbb{Q}(\sqrt{-7})$. The Jacobian variety $J(7)(\mathbb{Q}(\sqrt{-7}))$
of $X(7)$ over $\mathbb{Q}(\sqrt{-7})$ is isogenous to the product $E^3$
corresponding to the degenerate discrete series. Here the elliptic curve $E$
is just the modular curve $X_0(49)$ with conductor $49$, of equation
$$y^2+xy=x^3-x^2-2x-1.$$
In the context of anabelian algebraic geometry, the trivial representation
leads to the Galois representation $\rho_7: \text{Gal}(\overline{\mathbb{Q}}/\mathbb{Q})
\rightarrow \text{Aut}(\mathcal{L}(X(7)))$ and the last mathematical
testament of Galois \cite{Ga}. Namely, there is a decomposition
$\text{PSL}(2, \mathbb{F}_7)=S_4 \cdot C_7$, where $C_7$ is a cyclic group
of order $7$. Correspondingly, the Galois representation $\rho_7$ can be
realized by the transformation equation of elliptic functions of degree seven,
and there is a Galois resolvent of degree seven, whose coefficients are defined
over $\mathbb{Q}(\sqrt{-7})$:
$$c^7+\frac{7}{2} (-1 \mp \sqrt{-7}) \nabla c^4-7 \left(\frac{5 \mp
  \sqrt{-7}}{2}\right) \nabla^2 c-C=0.\eqno{(2.5.6)}$$
(See (4.1.24), (4.1.22) et seq. for more details).

  When $p=11$, in the context of motives, the Steinberg representation
and discrete series representation of $\text{PSL}(2, \mathbb{F}_{11})$
appear in Hecke's decomposition are realizable over $\mathbb{Q}$. While,
the degenerate discrete series representation of $\text{PSL}(2, \mathbb{F}_{11})$
is realizable over $\mathbb{Q}(\sqrt{-11})$. In fact, it is shown that
(see \cite{Ligozat}) the Jacobian variety $J(11)(\mathbb{Q}(\sqrt{-11}))$
of $X(11)$ over $\mathbb{Q}(\sqrt{-11})$ is isogenous to the product
$$E_1^{11} \times E_6^5 \times E_4^{10},$$
corresponding to the Steinberg representation, degenerate discrete
series, and discrete series, respectively. Here the elliptic curve
$E_1$ is just the modular curve $X_0(11)$ with conductor $11$, of
equation:
$$y^2+y=x^3-x^2-10 x-20.$$
The elliptic curve $E_6$ with conductor $121$ is given by the equation:
$$y^2+y=x^3-x^2-7x+10,$$
and the elliptic curve $E_4$ with conductor $121$ is given by the equation:
$$y^2+xy=x^3+x^2-2x-7.$$
From the viewpoint of modular curves with level $11$, the Jacobian variety
$J_{\text{split}}(11)$ is isogenous over $\mathbb{Q}(\sqrt{-11})$ to the
product $X_0(11) \times E_6$, $J_{\text{nonsplit}}(11)$ is isogenous over
$\mathbb{Q}(\sqrt{-11})$ to $E_6$, and $J_{\mathfrak{S}_4}(11)$ is isogenous
over $\mathbb{Q}(\sqrt{-11})$ to $E_4$. Here, the modular curve $X_0(11)$ is
isomorphic to $E_1$ over $\mathbb{Q}$, $X_{\text{nonsplit}}(11)$ is isomorphic
over $\mathbb{Q}$ to the curve $E_6$, and $X_{\mathfrak{S}_4}(11)$ is
isomorphic over $\mathbb{Q}$ to the curve $E_4^{\prime}=E_5/C_5$, of equation:
$$y^2+xy+y=x^3+x^2-305 x+7888,$$
where $E_5$ is given by the equation:
$$y^2+xy+y=x^3+x^2-30 x-76.$$
In the context of anabelian algebraic geometry, the discrete series
representation leads to the Galois representation $\rho_{11}:
\text{Gal}(\overline{\mathbb{Q}}/\mathbb{Q}) \rightarrow \text{Aut}(\mathcal{L}(X(11)))$
and the last mathematical testament of Galois \cite{Ga} again. Namely, there
is a decomposition $\text{PSL}(2, \mathbb{F}_{11})=A_5 \cdot C_{11}$, where
$C_{11}$ is a cyclic group of order $11$. Correspondingly, the Galois
representation $\rho_{11}$ can be realized by the transformation equation
of elliptic functions of degree eleven, and there is a Galois resolvent of
degree eleven, whose coefficients are defined over $\mathbb{Q}(\sqrt{-11})$:
$$\varphi^{11}+\alpha \nabla^2 \cdot \varphi^8+\beta \nabla^4 \cdot \varphi^5
 +\gamma \nabla C \cdot \varphi^4+\delta \nabla^6 \cdot \varphi^2+\varepsilon
  \nabla^3 C \cdot \varphi+\zeta \cdot C^2=0.\eqno{(2.5.7)}$$
(See (4.2.16), (4.2.11) et seq. for more details).

  The significance of the comparison (2.5.5) is that there are two different kinds
of viewpoints to look at it. Now, we will give the second viewpoint, in which the
cases of $p=7$ and $13$ are equally important as the above cases of $p=7$ and $11$.
More precisely, two Galois representations $\rho_7: \text{Gal}(\overline{\mathbb{Q}}/
\mathbb{Q}) \rightarrow \text{Aut}(\mathcal{L}(X(7)))$ and $\rho_{13}: \text{Gal}
(\overline{\mathbb{Q}}/\mathbb{Q}) \rightarrow \text{Aut}(\mathcal{L}(X(13)))$
can be (surjectively) realized by Jacobian multiplier equations of degree eight and
fourteen, respectively. Their modular realizations correspond to the transformations
between $J$ and $J^{\prime}$ of order seven and thirteen, respectively. This leads
to Klein's construction of transformation equations of order seven and thirteen,
respectively. In particular, these two equations can be solved by elliptic modular
functions which come from the multipliers of the transformation equations between
periods of elliptic integrals. Now, we need the different viewpoint: the periods
of elliptic integrals and the transformation equations between periods of elliptic
integrals. Recall that the elliptic differential
$$\frac{dx}{\sqrt{f(x)}}=\frac{dx}{\sqrt{a_0 x^4+4 a_1 x^3+6 a_2 x^2+4 a_3 x+a_4}},$$
which takes the following form when written homogeneously:
$$\frac{x_2 dx_1-x_1 dx_2}{\sqrt{a_0 x_1^4+4 a_1 x_1^3 x_2+6 a_2 x_1^2 x_2^2
       +4 a_3 x_1 x_2^3+a_4 x_2^4}},$$
can be regarded as a covariant of the biquadratic binary form in the
denominator; for if new variables are introduced instead of $x_1$, $x_2$
by a linear substitution, the substitution determinant appears as a factor.
Accordingly, it is essentially dependent on the rational invariants of
this binary form. According to Weierstrass, we denote these invariants
by $g_2$, $g_3$ and write accordingly:
$$g_2=a_0 a_4-4 a_1 a_3+3 a_2^2,$$
$$g_3=\begin{vmatrix}
      a_0 & a_1 & a_2\\
      a_1 & a_2 & a_3\\
      a_2 & a_3 & a_4
\end{vmatrix}.$$
The discriminant $\Delta$ of the biquadratic form is composed of
$g_2$ and $g_3$ in the known way
$$\Delta=g_2^3-27 g_3^2.$$
They can be used to form the absolute invariant $J$:
$$J=\frac{g_2^3}{\Delta}, \quad J-1=\frac{27 g_3^2}{\Delta}.$$

  We will express the periods of the elliptic integral
$$\int \frac{dx}{\sqrt{f(x)}},$$
arising from the above differential, which are to be called $\omega_1$
and $\omega_2$, by $g_2$, $g_3$, or else, which is more appropriate
for the following investigations, to represent the {\it ratio} of
the periods $\frac{\omega_1}{\omega_2}=\omega$ by the absolute invariant
$J$. One could call $\omega_1$ and $\omega_2$ the transcendental invariants.
For with $\omega_1$, $\omega_2$, as is known, every other pair of values
$$\left\{\aligned
  \omega_1^{\prime} &=\alpha \omega_1+\beta \omega_2,\\
  \omega_2^{\prime} &=\gamma \omega_1+\delta \omega_2
\endaligned\right.$$
is conjugated, provided that $\alpha$, $\beta$, $\gamma$, $\delta$ denote
integers whose determinant $\alpha \delta-\beta \gamma$ is equal to one.
The periods are invariants because they only change by the substitution
determinant as a factor if new variables are introduced into the given
integral instead of $x_1$, $x_2$ by linear substitution. It is often
useful to normalise the elliptic integral in such a way that its periods
are absolute invariants without further ado. The integral is therefore
written in the following form, which is also occasionally referred to
as the normal form (first level):
$$\int \frac{\root 12 \of{\Delta} \cdot dx}{\sqrt{f(x)}}.$$
Now, for the sake of a more definite expression, let $n$ be a prime number,
and we require the integral thus normalised to be converted into a likewise
normalised integral
$\int \frac{\root 12 \of{\Delta_1} \cdot dx_1}{\sqrt{f_1(x_1)}}$
by transformation of the $n$-th order. This results in a multiplier that
is defined by the equation:
$$M \int \frac{\root 12 \of{\Delta} \cdot dx}{\sqrt{f(x)}}=
  \int \frac{\root 12 \of{\Delta_1} \cdot dx_1}{\sqrt{f_1(x_1)}}.$$
If we now form a period on both sides, such as $\omega_2$, the following
follows, depending on the type of transformation (or the selected period):
$$M=\frac{\root 12 \of{\Delta_1} \cdot \omega_2^{\prime}}{\root 12 \of{\Delta}
    \cdot \omega_2}, \quad \text{or} \quad
   =\frac{1}{n} \frac{\root 12 \of{\Delta_1} \cdot \omega_2^{\prime}}{\root 12 \of{\Delta}
    \cdot \omega_2}.$$
In fact (by Jacobi's formula),
$$\root 12 \of{\Delta} \cdot \omega_2=2 \pi q^{\frac{1}{6}} \prod (1-q^{2 \nu})^2,
  \quad \text{for $q=e^{\pi i \omega}$}.$$
The same is true for
$$\root 12 \of{\Delta_1} \cdot \omega_2^{\prime}=2 \pi q_1^{\frac{1}{6}}
  \prod (1-q_1^{2 \nu})^2, \quad \text{for $q_1=e^{\pi i \omega^{\prime}}$},$$
and here, as we know, we have to set the following values for $q_1$:
$$q^{\frac{1}{n}}, \alpha q^{\frac{1}{n}}, \alpha^2 q^{\frac{1}{n}}, \ldots,
  \alpha^{n-1} q^{\frac{1}{n}}, q^n,$$
where $\alpha$ is a primitive $n$-th root of unity. The following $(n+1)$
expressions are thus obtained for $M$, which are initially defined only up
to the sixth roots of unity and, incidentally, are to be distinguished by
indices in a frequently used way:
$$\left\{\aligned
  M_{\varrho} &=\frac{1}{n} \cdot \frac{\alpha^{\frac{\varrho}{6}} q^{\frac{1}{6n}}
                \cdot \prod (1-\alpha^{2 \varrho \nu} q^{\frac{2 \nu}{n}})}
                {q^{\frac{1}{6}} \cdot \prod (1-q^{2 \nu})^2} \quad \text{for
                $\varrho=0$, $1$, $2$, $\ldots$, $(n-1)$},\\
  M_{\infty} &=\frac{q^{\frac{n}{6}} \cdot \prod (1-q^{2 \nu n})^2}
               {q^{\frac{1}{6}} \cdot \prod (1-q^{2 \nu})^2}.
\endaligned\right.$$

  In his paper \cite{K1}, using Riemann's idea about algebraic functions,
Klein studied the following question: How must $s$ be branched as a function
of $J$ if the equation $$\varphi(s, J)=0$$
is to be solved by elliptic modular functions? We mean that the equation
should be solvable in such a way that the period ratio $\omega$ is calculated
from $J$ and we now have unique functions of $\omega$ defined in the whole
positive $\omega$-half-plane, which represent the roots of $\varphi=0$. For
this, it is necessary and sufficient that the single root $s$, understood
as a function of $\omega$, does not branch within the positive $\omega$-half-plane.
Klein proved the following theorem: {\it Branch points may only lie in the
Riemann surface, which represents $s$ as a function of $J$, at $J=0$, $1$,
$\infty$. With $J=0$, three sheets can be connected any number of times,
with $J=1$, two sheets can be connected any number of times. For $J=\infty$,
the branching can be any.} All transcendental functions of $J$, which can
be uniquely represented by modular functions, are determined in exactly the
same way. At the same time, the problem is solved: {\it to set up all subgroups
which are contained in the entirety of the [integer] substitutions}
$$\omega^{\prime}=\frac{\alpha \omega+\beta}{\gamma \omega+\delta}, \quad
  (\alpha \delta-\beta \gamma=1).$$

  In particular, if we are looking for equations of genus zero and put
them in the form:
$$J=\frac{\varphi(\tau)}{\psi(\tau)},$$
then $\varphi$ may only contain triple factors in addition to single factors,
$\varphi-\psi$ may only contain double factors in addition to single factors,
and $\psi$ may only contain factors of any multiplicity. But no equation
$\lambda \varphi+\mu \psi=0$ that is different from $\varphi=0$, $\psi=0$,
$\varphi-\psi=0$ may have multiple roots. These equations are remarkable
because they have nothing to do with the transformation problem of elliptic
functions, and thus provide a first calculated example of the more general
equations that can be solved by elliptic modular functions.

  Let us restrict to the transformation problem of elliptic functions.
Let $J$ and $J^{\prime}$ be the absolute invariants of two elliptic
integrals which are obtained by transformation of the $n$-th order,
where $n$ may be a prime number. What we need to do is to set up the
equation of the $(n+1)$-th degree, which links $J^{\prime}$ with $J$.
The equations between $J^{\prime}$ and $J$ were first treated by Felix
M\"{u}ller in his dissertation published in 1867 following Weierstrass'
lectures. He started from the study of double-periodic functions and
arrived at ready equations for $n=2$, $3$, $4$, $5$, $7$. Then in 1874
Brioschi tackled the question from the algebraic side by directly
considering the transformation of the elliptic integral. However,
Klein's derivation of the transformation equations for $n=2$, $3$, $4$,
$5$, $7$, $13$ differs essentially in that it only makes use of the
variables $J$ and $\omega$, but not of the integration variable of the
elliptic integral or of this integral itself. And even the variable
$\omega$ only enters the consideration by means of which the branch
of $J^{\prime}$ with respect to $J$, but not in the calculation.
For the sake of simplicity, we will assume that the prime number $n$
is greater than $3$. We also think of $J^{\prime}$ as being calculated
by picking out any value of $\omega$ belonging to $J$ and now equating
$\omega^{\prime}$ in sequence:
$$\frac{\omega}{n}, \frac{\omega+1}{n}, \cdots, \frac{\omega+(n-1)}{n},
  -\frac{1}{n \omega}.$$
According to the former, a branch of $J^{\prime}$ with respect to $J$
can only take place at $J=0$, $1$, $\infty$. Klein studied the branch
of $J^{\prime}$ with respect to $J$, and showed that there is a unique
transformation of the Riemann surface into itself, which corresponds to
the interchange of $J$ and $J^{\prime}$. Moreover, Klein calculated the
genus of the equation between $J$ and $J^{\prime}$. In particular, the
genus is zero for $n=5$, $7$, $13$. Klein explained in more detail for
those cases which result in $g=0$. The Riemann surface can then be
extended into a plane and thus its division into areas can be visualised
directly by drawing. Klein set up the equations between $J$ and $J^{\prime}$
in the cases $n=2$, $3$, $4$, $5$, $7$, $13$ in such a way that he represented
both rationally by the variable $\tau$, which assumes each value only once
in the Riemann surface. In doing so, Klein used only the position and
multiplicity of the branch points of $J^{\prime}$ in relation to $J$. We need
only two cases: $n=7$ and $n=13$:

(1) Transformation of the seventh order:
$$\left\{\aligned
  J:J-1:1 =&(\tau^2+13 \tau+49)(\tau^2+5 \tau+1)^3\\
          :&(\tau^4+14 \tau^3+63 \tau^2+70 \tau-7)^2\\
          :&1728 \tau,\\
  \text{$J^{\prime}$ also in $\tau^{\prime}$},\\
  \tau \tau^{\prime} =&49.
\endaligned\right.\eqno{(2.5.8)}$$

(2) Transformation of the thirteenth order.
$$\left\{\aligned
  J&:J-1:1 \\
  =&(\tau^2+5 \tau+13)(\tau^4+7 \tau^3+20 \tau^2+19 \tau+1)^3\\
  :&(\tau^2+6 \tau+13)(\tau^6+10 \tau^5+46 \tau^4+108 \tau^3+122 \tau^2+38 \tau-1)^2\\
  :&1728 \tau,\\
  &\text{$J^{\prime}$ also in $\tau^{\prime}$},\\
  &\tau \tau^{\prime}=13.
\endaligned\right.\eqno{(2.5.9)}$$

  It is these formulae for $M_{\varrho}$ which readily give the resolution
of the above equations of order seven and thirteen by making $\tau$
identical with a power of $M$, apart from a numerical factor.  In order
to see this directly, Klein (see \cite{K1}) considered a value of $M$, say
$$M_0=\frac{1}{n} \cdot \frac{q^{\frac{1}{6n}} \cdot \prod
      (1-q^{\frac{2 \nu}{n}})^2}{q^{\frac{1}{6}} \cdot \prod (1-q^{2 \nu})^2},$$
as a function of the location in our Riemann surface $J$, $J^{\prime}$,
or, what amounts to the same thing, in our fundamental polygon located
in the $\omega$-plane. Obviously, $M_0$ becomes zero only once in the
fundamental polygon, where $\omega=0$, $q=1$. Similarly, it only becomes
infinite once, where $\omega=i \infty$, $q=0$. If we denote by $\lambda$
the smallest multiple of $\frac{n-1}{12}$, which is equal to an integer,
then $M_0^{\lambda}$ is unique in our Riemann surface [since this has the
genus zero in the cases we are considering] and is therefore a rational
function of $J$ and $J^{\prime}$, which becomes zero and infinite at only
one point. For $n=2$, $3$, $5$, $7$, $13$, $\lambda$ becomes equal to $12$,
$6$, $3$, $2$, $1$. At the same time, $\lambda \cdot \frac{n-1}{12}$ is not
only an integer, but also equal to one. Therefore, in these cases $M_0^{\lambda}$
is not only zero or infinite once, but also just zero or infinite. Apart
from one numerical factor, $M_0^{\lambda}$ therefore corresponds to the
previous $\tau$, which becomes zero and infinite at the same points and
in the same way. The numerical factor is determined as follows. Just as
Klein previously introduced a quantity $\tau^{\prime}$ next to $\tau$, he
considered the other multiplier $M^{\prime}$ next to $M$, which occurs
during the transition from the transformed integral back to the original.
Then we have in the familiar way:
$$M M^{\prime}=\frac{1}{n}.$$
If we compare this relation with the equations
$$\tau \tau^{\prime}=C$$
which appeared in the seven and thirteen order transformation, we obtain
the numerical factor we are looking for. We thus obtain the following formulae,
which contain the resolution of the equations of order $7$ and $13$:
$$\left\{\aligned
  n &=7, \quad \tau= 49 M^2,\\
  n &=13, \quad \tau=13M.
\endaligned\right.$$

  Now, we will give a unified way to look at the above three cases $p=7$, $11$
and $13$ by the Galois resolvents of the equations between $J$ and $J^{\prime}$.
It is known that the Galois group of the transformation equation comprises
six substitutions for $p=2$, twenty-four for $p=4$ and $\frac{p(p^2-1)}{2}$
for each prime number $p$ which is $>2$. The Riemann surface, which represents
the Galois resolvent of the transformation equation, will therefore have
just as many sheets. Klein (see \cite{K1}) showed that these sheets are related
to $3$ each for $J=0$, to $2$ each for $J=1$ and to $n$ each for $J=\infty$.
This therefore gives the genus $g=0$ for $p=2$, $g=0$ again for $p=4$, for the
other $p$, $g=\frac{(p-3)(p-5)(p+2)}{24}$. For $n=7$, $11$ and $13$, the genus
is equal to $3$, $26$ and $50$.

  Let $X$ be an algebraic curve of genus $g$ over the field of complex
numbers. The Jacobian of $X$ is said to be completely decomposable if
it is isogenous to a product of $g$ elliptic curves (see \cite{ES}).
Then, according to the above results, the Jacobian varieties $J(7)$ and
$J(11)$ are completely decomposable. Hence, their geometric and arithmetic
properties are completely determined by elliptic curves. On the other
hand, we will show that the anabelian part, $I(\mathcal{L}(X(7)))$ and
$I(\mathcal{L}(X(11)))$ contain much more information which do not come
from elliptic curves, i.e., the abelian part. In particular, the anabelain
part and the abelian part come from distinct kinds of representations,
which are defined over different algebraic number fields.

  In the general case, note that the modular curves $X(p)$ have two
kinds of geometry and topology: the abelian part $H_1$ and the
anabelian part $\pi_1$. The Jacobian varieties $J(p)$ of modular
curves $X(p)$ realize $H_1$, while the defining ideals
$I(\mathcal{L}(X(p))$ realize $\pi_1$. The above result shows that
$H_1$ and $\pi_1$ give different kinds of representations of
$\text{PSL}(2, \mathbb{F}_p)$, which are defined over distinct
algebraic number fields, hence with distinct arithmetic properties.
Moreover, $H_1$ corresponds to abelian varieties, while $\pi_1$
corresponds to ideals.

\begin{center}
{\large\bf 3. The division of elliptic functions, the transformation
              of elliptic functions and their Galois groups}
\end{center}

\begin{center}
{\bf \S 3.1. The division equations of elliptic functions, their Galois groups
             and their relation with transformation equations}
\end{center}

  For the development of division theorems of elliptic functions, Abel's
creations have been fundamental (see \cite{Abel1}). In particular, Abel
recognised the algebraic character of the division equations, which led
him to the general investigations on the resolution of algebraic equations.
More precisely, the division theorems of elliptic functions were discovered
by Abel and presented for the first time in the ``Recherches sur les fonctions
elliptiques'' (see \cite{Abel1}). Abel's realisation that the general division
equations are solvable by drawing roots must have been of considerable influence
on his general investigations of algebraically solvable equations. Given the
interest that the division equations present as examples for Galois' equation
theory, the division theorems come into play more extensively in those
presentations in which algebraic questions are in the foreground. Thus,
C. Jordan devoted a special chapter in his well-known work ``Trait\'{e} des
substitutions et des \'{e}quations algebriques'' (Paris, 1870) to the groups
of division equations, both the ``general'' division equations as well as the
``special'' division equations. H. Weber treated the algebraic side of the
division theorems in particular detail in his book ``Elliptic Functions and
Algebraic Numbers'' (Braunschweig, 1891). A second edition was published in
1908 as the third volume of the ``Lehrbuch der Algebra''. The theorem, that
the general division equation is an Abelian equation, presupposes the
adjunction and thus the knowledge of the division values $\wp_{\lambda \mu}$,
$\wp_{\lambda \mu}^{\prime}$. These in turn satisfy equations which are called
``special division equations'' and whose theory belongs to the most interesting
objects of Fricke's presentation \cite{Fricke}, Vol. II. The Galois group of
the special division equation belonging to degree $n$ was considered by
C. Jordan in his well-known group work mentioned as above without arriving
at definitive results. However, such results are achieved in the investigations
by Sylow and Kronecker. From a function-theoretical point of view, Kiepert
in particular had examined the division values of the functions of the
first level and traced their relationship to the transformation theory.

  The number of incongruent pairs of numbers mod $n$ that are coprime to
$n$ is denoted by the symbol $\chi(n)$. It is shown that (see \cite{Fricke},
Vol. II, p. 219) if the prime factorisation of the number $n$ is given by
$n=p_1^{\nu_1} \cdot p_2^{\nu_2} \cdot p_3^{\nu_3} \ldots$, then the number
$\chi(n)$ of the pairs $l$, $m$ coprime to $n$ and thus the number of
substitutions of period $n$ contained in $G_{n^2}$ is:
$$\chi(n)=n^2 \left(1-\frac{1}{p_1^2}\right)\left(1-\frac{1}{p_2^2}\right)
          \left(1-\frac{1}{p_3^2}\right) \cdots. \eqno{(3.1.1)}$$
In \cite{Fricke}, Vol. II, pp. 260-261, it is shown that the Galois group of the
special division equation is a $G_{\frac{1}{2} n \chi(n) \varphi(n)}$ of order
$\frac{1}{2} n \chi(n) \varphi(n)$, the group $G_{\frac{1}{2} n \chi(n)}$ of
order $\frac{1}{2} n \chi(n)$, the monodromy group of our equation, as a normal
subgroup. The associated quotient group:
$$G_{\frac{1}{2} n \chi(n) \varphi(n)}/G_{\frac{1}{2} n \chi(n)}=G_{\varphi(n)},$$
however, is the group of the cyclotomic equation of the $n$-th degree of division.
Using equation (2) in \cite{Fricke}, Vol. II, p. 255, the monodromy group would
refer to the orbits of $J$ in its plane. According to the explanations in
\cite{Fricke}, Vol. I, 299 et seq., we achieve these orbits by exercising the
substitutions of the modular group on $\omega$ or $\omega_1$, $\omega_2$.
{\it As the permutation group of the roots $\wp_{\lambda \mu}$ of the special
division equation, the monodromy group is thus given by the $\frac{1}{2} n \chi(n)$
incongruent substitutions} (7) in \cite{Fricke}, Vol. II, p. 246, {\it where
two substitutions, which merge by simultaneous sign change of the four
coefficients, are not different.} Of course, this is to be understood in such
a way that one thinks of the $\frac{1}{2} \chi(n)$ roots of our equation
written in a row and then successively exercises the $\frac{1}{2} n \chi(n)$
substitutions (7) in \cite{Fricke}, Vol. II, p.246 on the indices $\lambda$,
$\mu$, thus producing the permutations of $\wp_{\lambda \mu}$. Incidentally,
after \cite{Fricke}, Vol. II, p. 246 we can state the theorem {\it that the
monodromy group of the special division equation is isomorphic with the
$G_{\frac{1}{2} n \chi(n)}$ to which the non-homogeneous modular group mod
$n$ reduces.}

  To arrive at the Galois group, we have to replace $\varepsilon$ in turn
by all $\varphi(n)$ primitive unit roots of $n$-th degree. This
substitution effects on the indices $\lambda$, $\mu$ of the division
values $\wp_{\lambda \mu}$ the substitution
$$\lambda^{\prime} \equiv \lambda, \quad
  \mu^{\prime} \equiv \kappa \mu \quad \text{(mod $n$)},\eqno{(3.1.2)}$$
where $\kappa$ has to pass through all $\varphi(n)$ incongruent residues
mod $n$ which are coprime to of $n$. If we combine the substitutions (3.1.2)
with the $\frac{1}{2} n \chi(n)$ substitutions (7) in \cite{Fricke}, Vol.
II, p. 246, we arrive at Galois' $G_{\frac{1}{2} n \chi(n) \varphi(n)}$
(see also \cite{Weber}, \S 63, especially, p. 212):

\textbf{Theorem 3.1.1.} {\it The Galois group $G_{\frac{1}{2} n \chi(n) \varphi(n)}$
of the special division equation is obtained as the permutation group of
$\wp_{\lambda \mu}$ by the $\frac{1}{2} n \chi(n) \varphi(n)$ incongruent
substitutions:}
$$\lambda^{\prime} \equiv \alpha \lambda+\gamma \mu, \quad
  \mu^{\prime} \equiv \beta \lambda+\delta \mu \quad
  \text{(mod $n$)}, \eqno{(3.1.3)}$$
{\it whose determinants $(\alpha \delta-\beta \gamma)$ are all incongruent
residues $\kappa$ which are coprime to $n$, and for which two substitutions
which can be converted into one another by simultaneous sign change of the
coefficients are considered to be not different.}

  We therefore have to focus our attention solely on the prime degrees $n=p$
and, after adjunction of the $p$-th root of unity $\varepsilon$, we now consider
the monodromy group $G_{\frac{1}{2}p(p^2-1)}$ of the special division equation,
which we present after \cite{Fricke}, Vol. II, p. 246  in the form of the mod $n$
reduced non-homogeneous modular group $\Gamma$. The basic theorem now holds:
{\it For all prime numbers $p>3$, the group $G_{\frac{1}{2}p(p^2-1)}$ is simple,
i.e. it has no normal subgroup $($except for $G_1$ and $G_{\frac{1}{2}p(p^2-1)}$$)$}.
Accordingly, for $p>3$ the index series (cf. \cite{Fricke}, Vol. II, p. 13)
of our group $G_{\frac{1}{2}p(p^2-1)}$ consists only of the single member
$\frac{1}{2}p(p^2-1)$; and since this number is not a prime number, {\it
the special division equation of the $p$-th degree of division for $p>3$
cannot be solved algebraically according to the theorem of} \cite{Fricke},
Vol. II, p. 75. For $p=2$ and $p=3$, the division equations belong to degrees
$3$ and $4$: {\it In the two lowest cases, $p=2$ and $p=3$, the division
values can be calculated by root extraction alone.}

  In this situation, it was a particularly important aim of Klein's
theory of elliptic modular functions to solve the peculiar algebraic
problems which, corresponding to the monodromy groups for $p=5$, $7$,
$11$, $\ldots$ are called ``Galois problems'' of degree $60$, $168$,
$660$, $\ldots$, generally of degree $\frac{1}{2} p(p^2-1)$. For $p=5$
one arrives at the ``icosahedron theory'' (see \cite{K} and \cite{K1}).
In the case of $p=7$, Klein developed his particularly beautiful theory
of the underlying group $G_{168}$ by direct algebraic methods without
the aid of series expansions of the elliptic functions (see \cite{K2}),
and also in the case of $p=11$ he succeeded in carrying out a corresponding
theory (see \cite{K3}). For the higher cases Klein, however, used the
analytic tools of the elliptic functions.

  The treatment of these Galois problems is given in great detail in
the work \cite{KF1} and \cite{KF2}. In contrast, the present development
takes the turn, that it more closely follows the algebraic points of view
originating in the older theory. It is not a question of Galois resolvents,
but the ``resolvents of the lowest degree'' of the special division equations.
These are, at least in general, the modular and multiplier equations or, as
we shall say, the ``special transformation equations'', which, in the
transformation of higher degree of the elliptic functions. In this area, too,
by the way, as will be explained in more detail below, Klein's methods have
had many groundbreaking effects. Before we go into the transformation theory
of elliptic functions in general, a few remarks on the division values of
the Jacobian functions sn, cn, dn should be added.

  The investigations by Sylow and Kronecker mentioned as above refer to
the division values of the sn function and, incidentally, only concern
{\it odd} degrees $n$ of division. {\it The task of calculating the
squares} sn$_{\lambda \mu}^2$ {\it differs from the problem of calculating
the $\wp_{\lambda \mu}$ only in that the adjunction of the modular form
of the second level $(e_2-e_1)$ still has to be carried out.}  The main
task is now to determine the group of this special division equation.

\textbf{Theorem 3.1.2.} (see \cite{Fricke}, Vol. II, pp. 266-268) {\it The
individual permutation of the Galois group transform} $\text{sn}_{\lambda \mu}$
{\it into} $\text{sn}_{\lambda^{\prime} \mu^{\prime}}$, {\it where}
$$\lambda^{\prime} \equiv \alpha \lambda+\gamma \mu, \quad
  \mu^{\prime} \equiv \beta \lambda+\delta \mu \quad \text{(mod $n$)}\eqno{(3.1.4)}$$
{\it holds and $\alpha$, $\beta$, $\gamma$, $\delta$ are four integers, of which
at least $\alpha$, $\beta$ and $\gamma$, $\delta$ form pairs coprime to $n$.}
{\it $G_{n \chi(n) \varphi(n)}$ certainly contains the Galois group of the
special division equation of the} $\text{sn}$-{\it function}. {\it The Galois
group of the division equation is therefore really $G_{n \chi(n) \varphi(n)}$}.
{\it After adjunction of $\varepsilon$, the group of the equation indeed comes
back to the $G_{n \chi(n)}$.}

  The investigation of the resolvents of the lowest degree of the specific
division equation was described as the aim of further development on
\cite{Fricke}, Vol. II, p. 265. Fricke extended this task in such a way
that one can also look for remarkable resolvents of lower degree for the
general division equation. In fact, such resolvents are given to us by the
equations of transformation theory, the ``general'' and the ``special
transformation equations''. This is an extraordinarily versatile developed
theory, whose beginnings go back to the early 19th century, whose most
important principles were established by Abel and Jacobi, but which only
gained its present form through the development of the theory of elliptic
modular functions.

  The basic task of transformation theory already appears with Euler and
Lagrange, in the treatment of certain mechanical tasks. In the language
of elliptic functions, it is a question of whether an elliptic differential,
which has $z$ as a variable, can be transformed into another elliptic
differential with $z^{\prime}$ as a variable by prescribing an {\it
algebraic} relation between $z$ and $z^{\prime}$. This is the algebraic
version of the transformation problem, while the transcendental form is
connected to the concept of the double-periodic function. This transcendent
figure foreshadows in the simplest way the general solution to the problem.
In the fundamental investigations of Abel and Jacobi, both sides of the
problem naturally come into play, the algebraic side mainly in the ``Pr\'{e}cis
d'une th\'{e}orie des fonctions elliptiques'' \cite{Abel2} and in the first
part of the ``Fundamenta nova'' \cite{Ja}, the transcendental in the
``Recherches sur les fonctions elliptiques'' \cite{Abel1} and the later
developments of the ``Fundamenta nova'' \cite{Ja}. In more recent times,
the transformation problem is usually presented in a transcendental form
and the algebraic version is only mentioned in passing (see \cite{Fricke},
Vol. II, p. 270 et seq.)

  Let $\varphi(u| \omega_1, \omega_2)$ and $\psi(u^{\prime} | \omega_1^{\prime},
\omega_2^{\prime})$ be two elliptic functions with fixed periods. The basic
problem of the general transformation is then in transcendental form: {\it
Under what circumstances does the assumption of a linear relation
$u^{\prime}=mu+\mu$ with constant $m$, $\mu$ entail an ``algebraic''
relation $F(\varphi, \psi)=0$ between the elliptic functions $\varphi$
and $\psi$?}

  Since $\varphi(u| \omega_1, \omega_2)$ and $\wp(u| \omega_1, \omega_2)$ are
algebraically related, and likewise $\psi(u^{\prime}| \omega_1^{\prime},
\omega_2^{\prime})$ and $\wp(u^{\prime}| \omega_1^{\prime}, \omega_2^{\prime})$,
we can refer our question to $\wp(u| \omega_1, \omega_2)$ and $\wp(u^{\prime}|
\omega_1^{\prime}, \omega_2^{\prime})$ instead of the functions $\varphi(u|
\omega_1,\omega_2)$ and $\psi(u^{\prime}| \omega_1^{\prime}, \omega_2^{\prime})$
without any restriction of generality. Further, $\wp(mu+\mu| \omega_1^{\prime},
\omega_2^{\prime})$ is algebraically related to $\wp(mu| \omega_1^{\prime},
\omega_2^{\prime})$ (on the basis of the addition theorem), and, since the
$\wp$-function is homogeneous of dimension $-2$ in its three arguments,
$$\wp(mu| \omega_1^{\prime}, \omega_2^{\prime})=m^{-2} \cdot
  \wp\left(u \mid \frac{\omega_1^{\prime}}{m}, \frac{\omega_2^{\prime}}{m}\right)$$
holds. If one therefore immediately sets $\omega_1^{\prime}$, $\omega_2^{\prime}$
again for $\frac{\omega_1^{\prime}}{m}$, $\frac{\omega_2^{\prime}}{m}$, the
problem posed has taken the following form without restriction of generality:
{\it Under which circumstances, i.e. for which pairs of periods $\omega_1$,
$\omega_2$ and $\omega_1^{\prime}$, $\omega_2^{\prime}$, does an algebraic
relation
$$F\left(\wp(u | \omega_1, \omega_2), \wp(u | \omega_1^{\prime},
  \omega_2^{\prime})\right)=0\eqno{(3.1.5)}$$
exist between $\wp(u| \omega_1, \omega_2)$ and $\wp(u| \omega_1^{\prime},
\omega_2^{\prime})$ with coefficients which are naturally supposed to be
independent of $u$?}

  Let us write $z$ and $z^{\prime}$ for the two $\wp$-functions connected
in relation (3.1.5), and to introduce the algebraic notation of the differential
$du$ we denote $g_2(\omega_1^{\prime}, \omega_2^{\prime})$ and $g_3(\omega_1^{\prime},
\omega_2^{\prime})$ briefly by $g_2^{\prime}$, $g_3^{\prime}$, then we can
dress up our problem in the following algebraic form: {\it Under what
circumstances can we in general integrate the differential equation:
$$\frac{dz}{\sqrt{4 z^3-g_2 z-g_3}}=\frac{dz^{\prime}}{\sqrt{4
  z^{\prime 3}-g_2^{\prime} z^{\prime}-g_3^{\prime}}}\eqno{(3.1.6)}$$
by an algebraic relation $F(z, z^{\prime})=0$, or how can the differential
on the left in $(8.1.6)$ be transformed into a differential of the same
construction by means of an algebraic function $z^{\prime}$ of $z$?}

  The solution of the problem is suggested by the following
consideration: The group $\Gamma^{(u)}$ introduced in \cite{Fricke},
Vol. I, p. 229 for the periods $\omega_1$, $\omega_2$ will be called
$\Gamma$ for short; the corresponding group for the periods $\omega_1^{\prime}$,
$\omega_2^{\prime}$ of the second $\wp$-function is called $\Gamma^{\prime}$.
In \cite{Fricke}, Vol. II, p. 271, it is proved that: {\it In any case, an
algebraic relation $(3.1.5)$ can only exist between $\wp(u | \omega_1, \omega_2)$
and $\wp(u | \omega_1^{\prime}, \omega_2^{\prime})$ if the groups $\Gamma$
and $\Gamma^{\prime}$ are commensurable. Conversely the commensurability
of the two groups $\Gamma$ and $\Gamma^{\prime}$ is also sufficient for the
existence of an algebraic relation} (3.1.5). Furthermore, the following theorem
holds:

\textbf{Theorem 3.1.3.} {\it For the two functions $\wp(u | \omega_1, \omega_2)$
and $p(u | \omega_1^{\prime}, \omega_2^{\prime})$ with commensurable groups
$\Gamma$, $\Gamma^{\prime}$, representations
$$\left\{\aligned
  \wp(u | \omega_1, \omega_2) &=R(\wp(u | A \omega_1+B \omega_2, D \omega_2)),\\
  \wp(u | \omega_1^{\prime}, \omega_2^{\prime}) &=R^{\prime}(\wp(u | A \omega_1
  +B \omega_2, D \omega_2))
\endaligned\right.\eqno{(3.1.7)}$$
hold as rational functions of the same argument $\wp(u | A \omega_1+B \omega_2,
D \omega_2)$, by whose elimination from the equations $(3.1.7)$ an algebraic
relation $(3.1.5)$ between $\wp(u | \omega_1, \omega_2)$ and $\wp(u | \omega_1^{\prime},
\omega_2^{\prime})$ is obtained. Here, if one chooses all divisors of $n$ for
$D$ in turn, $A=\frac{n}{D}$ and enter $0$, $1$, $2$, $\ldots$, $D-1$ for $B$
in succession.}

  The commensurability of the groups is therefore, as already claimed above,
is also sufficient for the existence of a relation (3.1.5).

  A subgroup $I_n$ of index $n$ in the group $\Gamma$ has as its discontinuity
region a parallelogram whose vertices $0$, $\omega_2^{\prime}$, $\omega_1^{\prime}
+\omega_2^{\prime}$, $\omega_1^{\prime}$ are lattice points of the original network,
and whose content is $n$ times as large as the content of the parallelogram of
the vertices $0$, $\omega_2$, $\omega_1+\omega_2$, $\omega_1$. We assume the
labels $\omega_1^{\prime}$, $\omega_2^{\prime}$ distributed among the vertices
in such a way that the quotient $\omega^{\prime}=\omega_1^{\prime}: \omega_2^{\prime}$
has positive imaginary component. {\it Conversely, each such parallelogram
provides the discontinuity region of a subgroup $\Gamma_n$ of index $n$ in
$\Gamma$, which is producible from the two substitutions:}
$$u^{\prime}=S_i^{\prime}(u)=u+\omega_i^{\prime}, \quad (i=1, 2).\eqno{(3.1.8)}$$
The $\omega_1^{\prime}$, $\omega_2^{\prime}$ represent themselves as lattice
points of the original network in the form:
$$\omega_1^{\prime}=a \omega_1+b \omega_2, \quad
  \omega_2^{\prime}=c \omega_1+d \omega_2 \eqno{(3.1.9)}$$
by means of integers $a$, $b$, $c$, $d$. Since the quotient $\omega^{\prime}
=\omega_1^{\prime}: \omega_2^{\prime}$ has a positive imaginary component,
the determinant $(ad-bc)$ is positive, and since the content of the parallelogram
of the vertices $0$, $\omega_2^{\prime}$, $\omega_1^{\prime}+\omega_2^{\prime}$,
$\omega_1^{\prime}$ is equal to the $n$-fold content of the original parallelogram,
in particular
$$ad-bc=n \eqno{(3.1.10)}$$
holds. Accordingly, one can say: {\it If one calculates the $\omega_1^{\prime}$,
$\omega_2^{\prime}$ from the $\omega_1$, $\omega_2$ according to $(3.1.9)$ by
means of any integers $a$, $b$, $c$, $d$ of the determinant $n$, then the
substitutions $S_1^{\prime}$, $S_2^{\prime}$ given in $(3.1.8)$ produce a
subgroup of the index $n$ in $\Gamma$; at the same time one arrives in this
way at all subgroups $\Gamma_n$ of $\Gamma$.}  We now say that in (3.1.9)
there is a ``{\it transformation of the $n$-th degree}'' of the periods
$\omega_1$, $\omega_2$, if $a$, $b$, $c$, $d$ are four integers of determinant
$n$. According to \cite{Fricke}, Vol. II, p. 273, each class must contain one
and only one transformation of the form:
$$T=\left(\begin{matrix} A & B\\ 0 & D \end{matrix}\right), \quad
  A \cdot D=n, \quad 0 \leq B < D,\eqno{(3.1.11)}$$
where $A$, $B$, $D$ are positive integers satisfying the conditions given
in (3.1.11). In \cite{Fricke}, Vol. II, p. 220, Fricke introduced a third
number $\psi(n)$ dependent on $n$ by:
$$\psi(n)=n\left(1+\frac{1}{p_1}\right)\left(1+\frac{1}{p_2}\right)
          \left(1+\frac{1}{p_3}\right) \cdots, \eqno{(3.1.12)}$$
then $\chi(n)=\varphi(n) \cdot \psi(n)$. Among the $\psi(n)$ actually
belonging to the degree $n$ essentially different transformations we call
the one given by $A=n$, $B=0$, $D=1$ the ``{\it first principal transformation}'',
while the numbers $A=1$, $B=0$, $D=n$ may provide the ``{\it second principal
transformation}'' (see \cite{Fricke}, Vol. II, p. 278-279).

\begin{center}
{\bf \S 3.2. The theory of the special transformation equations}
\end{center}

  Indeed, in \cite{Fricke}, Vol. II, p. 296, Fricke pointed out that the main
interest of transformation theory has essentially deviated from the ``general
transformation equations'', and turned to the ``special transformation equations'',
which are far more interesting from an algebraic and arithmetic point of view,
which are satisfied by the transformed modular functions and modular forms. In
particular, the Jacobian modular equations and multiplier equations also belong
to these equations. If one applies the transformation of the $n$-th degree to
the modular functions and modular forms which depend on the periods $\omega_1$,
$\omega_2$ alone, the result is quantities which are again algebraically related
to the original functions or forms. We call these algebraic relations ``special
transformation equations''. They correspond to the special division equations
and are, as has already been occasionally remarked, in general, i.e., apart from
some special degrees of transformation, the resolvents of the lowest degree of
the special division equations. Because of the richly related and extended theory
of these special transformation equations, it seems expedient to focus the
consideration first on the transformation of the special division equations, to
restrict the consideration at first to the transformation of the functions of the
first level $J(\omega)$, $g_2(\omega_1, \omega_2)$, $g_3(\omega_1, \omega_2)$
and $\Delta(\omega_1, \omega_2)$; but at the same time also the roots of
the discriminant $\Delta$, in so far as they yield unique modular forms,
are to be admitted at the same time.

  In the algebraic individual explanations of these special transformation
equations we still need knowledge of certain systems of elliptic modular
forms of the $n$-th level, which are advantageously used in all deeper
investigations of the transformation of the functions of the first level.
Those systems of modular forms are provided by certain systems of entire
elliptic functions of the $n$-th level of the third kind, whose theory
is remarkable in itself, especially from the point of view of group
theory (see \cite{Fricke}, Vol. II, p. 297 et seq.). According to this,
the algebraic equations which are satisfied by the quantities arising from
$g_2$, $g_3$, $\Delta$ and $J$ on transformation of the $n$-th degree should
be called ``{\it special transformation equations}''. In particular,
we call them special transformation equations of the ``first level''
and will also calculate those equations for them for the roots of
the discriminant of the $8$-th, $12$-th and $24$-th degree of the
discriminant $\Delta$. The special transformation equations correspond,
as noted above, the special division equations, to which they belong
as resolvents. First of all the general theory of these first-level
transformation equations has been developed in \cite{Fricke}, Vol. II,
p. 335 et seq. Concluding individual investigations on lower transformation
degrees follow in \cite{Fricke}, Vol. II, p. 370 et seq. In particular,
in \cite{Fricke}, Vol. II, p. 338, an essential property of a special
transformation equation of the first level is proved as follows:

\textbf{Theorem 3.2.1.} {\it The special transformation equations of
the $\psi(n)$-th degree of the modular forms $g_2$, $g_3$, $\Delta$
and the modular function $J$ are rational resolvents of the special
division equation of the $\wp$-function; the coefficients of those
transformation equations therefore belong to the field
$(\mathfrak{R}, g_2, g_3)$.}

  Note that in \cite{Weber}, \S 65, p. 223, a similar relation between
division equations and transformation equations is given as follows:
{\it If $n$ is a prime number or a power of a prime number, the square
roots of the division equation are rationally expressible by the roots
of the transformation equation.}

  The most important transformation equations of the first level are
those for the modular function of the first level $J(\omega)$. Incidentally,
we do not want to write the equations for $J(\omega)$ here, but according
to H. Weber's procedure for the function
$$12^3 J(\omega)=j(\omega),$$
because in this way, from an arithmetical point of view, significant advantages
are obtained. These advantages are based on the power series for $j(\omega)$,
which we first establish. From (3) in \cite{Fricke}, Vol. I, p. 274 and (9) in
\cite{Fricke}, Vol. I, p. 433 follows
$$\aligned
  12 g_2 &=\left(\frac{2 \pi}{\omega_2}\right)^4 (1+240 q^2+2160 q^4+6720 q^6+\cdots),\\
  \frac{1}{\root 3 \of{\Delta}} &=\left(\frac{\omega_2}{2 \pi}\right)^4
  (1+8 q^2+44 q^4+192 q^6+\cdots),
\endaligned$$
where the power series in the brackets have integer coefficients throughout.
By multiplying the last two equations it follows:
$$\aligned
  \frac{12 g_2}{\root 3 \of{\Delta}} &=12 \root 3 \of{J}=\root 3 \of{j(\omega)}\\
  &=q^{-\frac{2}{3}} (1+248 q^2+4124 q^4+34752 q^6+\cdots).
\endaligned$$
By raising to the third power, the result for $j(\omega)$ itself is:
$$j(\omega)=q^{-2}+744+196884 q^2+21493760 q^4+\cdots$$
or by splitting the coefficients into prime factors
$$j(\omega)=q^{-2}+2^3 \cdot 3 \cdot 31+2^2 \cdot 3^3 \cdot 1823 q^2
  +2^{11} \cdot 5 \cdot 2099 q^4+\cdots.$$
The important thing here is {\it that the coefficient of the first member
is equal to $1$ and the coefficients of all the links are integers.}

  The $\psi(n)$ to be obtained from $j(\omega)$ in transformation of the
$n$-th degree are now
$$j\left(\frac{A \omega+B}{D}\right)=e^{-\frac{2 \pi i B}{D}} q^{-\frac{2A}{D}}
  +744+196884 e^{\frac{2 \pi i B}{D}} q^{\frac{2 A}{D}}+\cdots,\eqno{(3.2.1)}$$
where $A$, $B$, $D$ are used in the familiar sense. If $\omega$ lies in the
interior of the $\omega$-half plane, the same applies to all transformed
values $\frac{A \omega+B}{D}$. The symmetric basic functions of the $\psi(n)$
transformed functions (3.2.1), in which we first recognised rational functions
of $j(\omega)$, are therefore ``entire'' rational functions of $j(\omega)$,
since they cannot become infinite for any finite value of $j$. For the
transformation equation we thus have the approach
$$j^{\prime \psi}+G_1(j) j^{\prime \psi-1}+G_2(j) j^{\prime \psi-2}+\cdots+
  +G_{\psi}(j)=0, \eqno{(3.2.2)}$$
where the coefficients are rational entire functions of $j(\omega)$ with
rational number coefficients. It can now be proved {\it that the number
coefficients of the functions $G(j)$ are all integers.} Another important
theorem is, {\it that the transformation equation $(3.2.2)$ also increases
in $j$ to the degree $\psi(n)$.} Concerning the algebraic nature of the
transformation equations, we first note only the following theorem to be
used immediately: {\it The entire transformation equations of the $\psi$-th
degree are irreducible in the field $(\mathfrak{R}, g_2, g_3)$ and remain
irreducible even if one extends the field by adjunction of some ``numerical''
irrationality} (see \cite{Fricke}, Vol. II, p. 344-347).

  The whole integer function of $j^{\prime}$ and $j$ on the left side of
the transformation equation (3.2.2) may be called $F(j^{\prime}, j)$ for
short. If one sets in particular $j^{\prime}=j(n \omega)$, corresponding
to the first principal transformation, then equation
$$F(j(n \omega), j(\omega))=0$$
exists identically in $\omega$ and thus also remains valid if one enters
$\frac{\omega}{n}$ instead of $\omega$:
$$F\left(j(\omega), j\left(\frac{\omega}{n}\right)\right)=0.$$
Summarising all the individual results, Fricke wrote down the theorem (see
\cite{Fricke}, Vol. II, p. 348):

\textbf{Theorem 3.2.2.} {\it The irreducible transformation equation of
the $\psi(n)$-th degree for $j(\omega)=12^3 J(\omega)$ occurring with
transformation of the $n$-th degree has the form:
$$j^{\prime \psi}+j^{\psi}+\sum_{k, l} c_{kl} j^{\prime k} j^l=0,\eqno{(3.2.3)}$$
where the sum refers to the combinations of integers $k$, $l$ of the
series $0$, $1$, $2$, $\ldots$, $\psi-1$; the $c$ are integers satisfying
the condition $c_{kl}=c_{lk}$, and especially in the case of prime number
$n$, where $\psi(n)=n+1$, $c_{nn}=-1$, while in the case of a composite
degree of transformation $c_{\psi-1, n}=c_{n, \psi-1}=-1$ holds.}

  The theory of special transformation equations underwent a further
development of fundamental importance through the theory of modular
functions developed by Klein. Compared to the various types of transformation
equations, this theory first of all provided the possibility of sorting and
appropriate arrangement by providing a judgement on which of the transformed
quantities is to be regarded as the simplest in the individual case and what
the relations of the others are to these simplest quantities. In addition to
the power series, Klein used the tools of the Riemann theory of algebraic
functions, which had grown up on a group-theoretical and geometrical basis,
to set up the transformation equations. The following presentation will show
the great successes of Klein's methods. The transformation equations for
$J(\omega)$ now come to the fore. The prime degrees up to $n=13$ were treated
by Klein himself; following this, J. Gierster carried out
investigations on the transformation equations of $J(\omega)$ for some lower
composite degrees. Later, these investigations were taken up by Kiepert
and were also substantially promoted with the aid of series expansions. A. Hurwitz
gave a treatment of the transformation equations for the roots of the discriminant
from the point of view of the theory of modular functions.

\begin{center}
{\bf \S 3.3. The transformation theory of elliptic functions and modular equations}
\end{center}

  Since the theory of transformation for elliptic functions plays a central
role in our theory, we will give its connection with modern concepts such as
Kronecker congruence relation, Eichler-Shimura congruence relation,
Langlands-Kottwitz's method, Mazur's Eisenstein ideals, Weil, Igusa,
Deligne-Rapoport, Drinfeld and Katz-Mazur's work as follows. In fact, the theory
of transformations of elliptic integrals and elliptic functions was the historic
predecessor, in the nineteenth century, of today's theory of isogenies between
elliptic curves, and of modular curves. It has been customary to call the theory
of isogeny of elliptic curves with variable moduli as transformation theory of
elliptic functions.

  Kronecker's congruence relation is of great importance, and which is the
principal aim of his efforts in the paper \cite{Kronecker3}. As he himself
confirmed, this congruence relation has a fundamental meaning for the
entire transformation theory of elliptic functions as well as its
arithmetic applications. The congruence relation (64) just says that the
reduction modulo $p$ of some $p$-isogeny coinciding with the Frobenius map.
This leads to the Eichler-Shimura congruence relation (see \cite{Sh}, p. 181
or \cite{Vladut}, p. 76) which comes from the study of reductions mod $p$ of
$p$-isogenies of an elliptic curve $E$ such that $j(E)=j$, $j$ being an independent
variable. This relation says that exactly one $p$-isogeny of $E$ is mapped
to the Frobenius map by the reduction modulo $p$, and all the others are
mapped to its transpose, the latter being equivalent to the separability
of their reductions. The statement on separability of reductions of all
the other $p$-isogenies can also be deduced from Kronecker's complements
to his main result. In the context of modular curves $X_0(N)$ with
$p \nmid N$, the Eichler-Shimura congruence relation expresses the local
$L$-function of a modular curve at a prime $p$ in terms of the eigenvalues
of Hecke operators. That is, the Hecke correspondence $T_p$ on the pair
of modular curves $X_0(N)$ is congruent mod $p$ to the sum of the Frobenius
map $F$ and its transpose $V=p/F$. In other words,
$$T_p=F+V$$
as endomorphisms of the Jacobian $J_0(N)(\mathbb{F}_p)$ of the modular
curve $X_0(N)$ over the finite field $\mathbb{F}_p$. It plays the crucial
role in the reduction of the Ramanujan conjecture to the Weil conjecture.
In fact, Eichler-Shimura made use of this relation to prove that for
all but finitely many $p$ the Ramanujan conjecture held for any given
cusp form of weight two and level $N$ which is a simultaneous
eigenfunction of all $T_p$ with $p$ not dividing $N$. In 1968, Deligne
(see \cite{De1968}) completed the general reduction of Ramanujan
conjecture for forms of arbitrary weight by generalizing the congruence
formula to any weight, and in particular for $\Delta$ with weight $12$,
to the Weil conjecture. The Eichler-Shimura congruence relation and its
generalizations to Shimura varieties, i.e., the Langlands-Kottwitz
method, can be used to show that the Hasse-Weil $L$-function of a
modular curve or more general a Shimura variety is a product of
automorphic $L$-functions. On the other hand, in the above
Eichler-Shimura relation, $F$ is the identity on $J_0(N)(\mathbb{F}_p)$,
this implies that $T_p=1+p$ on $J_0(N)(\mathbb{F}_p)$. This leads to
Mazur's Eisenstein ideal \cite{Mazur1977}. For a prime $N$, Mazur
proved Ogg's conjecture via a careful study of subgroups of the
Jacobian variety $J_0(N)$ annihilated by the Eisenstein ideal.
As a consequence, Mazur obtained multiplicity-one results for mod
$p$ Galois representations (see \cite{Mazur1977} for more details).
The Eisenstein ideals are closely related to the congruence module
that arises measures congruences between cusp forms and an Eisenstein
series (see \cite{MW}). Indeed, the ideas of \cite{MW} provide the
means to estimate the order of this module.

  In a 1950 congress address \cite{Weil1950}, Andr\'{e} Weil gave a geometric
formulation of the transformation theory and suggested that ``it is very probable
that a reconsideration of this splendid work from a modern point of view would
not merely enrich our knowledge of elliptic function-fields, but would also
reveal principles of great importance for any further development of algebraic
geometry over integers.'' In a series of papers \cite{Igusa1}, \cite{Igusa2},
\cite{Igusa3} and \cite{Igusa4}, Igusa gave a geometric construction of the
field of elliptic modular functions of an arbitrary level $N$ for any
characteristic which does not divide $N$. He put together the fields of
modular functions of level $N$ in different characteristics so that one
gets a fibre system parameterized by the ring of rational integers. In
particular, he formulated the transformation theory of elliptic functions
as follows (\cite{Igusa2}): Suppose that $A$ is an elliptic curve whose
absolute invariant $j(A)$ is a variable over a prime field $F$. Consider a
homomorphism of $A$ onto another elliptic curve $B$. Since every homomorphism
can be factored into homomorphisms of prime degrees, we can assume that the
degree of the homomorphism is a prime number $p$. Also we can assume that $p$
is different from the characteristic of $F$, for otherwise the situation will
become trivial. Then the transformation theory of elliptic functions will
answer the following problems:

(1) How is $j(B)$ related to $j(A)$?

(2) What is the arithmetic nature of the homomorphism?

  The first part is the transformation of the parameter space while the
second part is the transformation of the entire fibre system of elliptic
curves. At any rate the homomorphism is determined up to an isomorphism
by its kernel, which is a cyclic subgroup of $A$ of order $p$. We know
that the points of $A$ of order $p$ form a finite Abelian group of type
$(p, p)$. Therefore we have $p+1$ distinct homomorphisms from $A$ to
another elliptic curves of which $B$ is a member. This already suggests
the importance of points of finite orders. In fact the Galois theory
of these points is a necessary preliminary for the transformation
theory. Note that the method of Eichler-Shimura was intrinsically
incapable of specifying the exceptional $p$, which consist only of
primes dividing $N$. In order to settle definitely this question
of exceptional $p$ for weight two forms, Igusa gave a complete and
definitive account of the level $N$ moduli scheme over $\mathbb{Z}[1/N]$.
This leads to an interpretation of a modular curve as a (coarse)
moduli space for elliptic curves endowed with a so-called level
structure and its relation with modular forms. In \cite{DeRa},
Deligne and Rapoport gave the modular interpretation of the cusps
by compactification of modular curves and their reduction modulo
$p$. Given level $H$ structures for congruence subgroups $\Gamma_H$
inverse images of subgroups $H$ of $\text{GL}(2, \mathbb{Z}/N \mathbb{Z})$,
such as $\Gamma_0(N)=\Gamma_H$ for $H=\left(\begin{matrix} * & *\\ 0 & *
\end{matrix}\right)$, the reduction modulo a prime number $p$ of the
corresponding modular stacks $\mathcal{M}_H[1/N]$ is examined. For $p$
dividing $N$, a model $\mathcal{M}_H$ of $\mathcal{M}_H[1/N]$ over
$\text{Spec} \mathbb{Z}$ is needed: Deligne and Rapoport defined
$\mathcal{M}_H$ as the normalization of $\mathcal{M}_1$ in
$\mathcal{M}_H^{0}[1/N]$. The stack $\mathcal{M}_{\Gamma_0(p)}$ has
a modular interpretation: it classifies pairs $(C/S, A)$ of a generalized
elliptic curve over $S$ and a rank $p$ locally free subgroup $A$ meeting
each irreducible component of any geometric fiber of $C$. In particular,
they proved a refinement of the Eichler-Shimura congruence formula (see
also \cite{Illusie}):

\textbf{Theorem 3.3.1.} {\it The stack $\mathcal{M}_{\Gamma_0(p)}$ is regular,
proper and flat over $\text{Spec} \mathbb{Z}$, of relative dimension $1$,
smooth outside the supersingular points of characteristic $p$, and with
semistable reduction at these points; $\mathcal{M}_{\Gamma_0(p)} \otimes
\mathbb{F}_p$ is the union of two irreducible components crossing
transversally at the supersingular points. Moreover, the $($open$)$
coarse moduli space $M_{\Gamma_0(p)}^{0}$ is the spectrum of the
normalization of $\mathbb{Z}[j, j^{\prime}]/(\Phi_p(j, j^{\prime}))$,
where $\Phi_p(j, j^{\prime})$ is the modular equation, a polynomial
congruent to $(j-j^{\prime p})(j^{\prime}-j^p)$ modulo $p$.}

  For more general groups $H$ the definition of $\mathcal{M}_H$ as a
normalization made it difficult to study its reduction mod $p$. Drinfeld's
notion of full level $N$ structures, providing a simple modular interpretation
of $\mathcal{M}_H$, solved the problem. In their treatise \cite{KaMa},
Katz and Mazur gave a systematic exposition of this theory.

  In the present paper, we will give a different geometric construction
by the modular equations and their Galois theory, which arise from the
solution to the above problem (1) of Igusa. In fact, the above theorem
of Deligne-Rapoport has already shown the importance of modular equations
even in the Eichler-Shimura theory.

  Let us recall some basic facts about modular equations (see \cite{dSG}).
Given an integer $N \geq 2$, we seek an equation linking $j(\tau)$ and
$j^{\prime}(\tau)=j(N \tau)$ for $\tau \in \mathbb{H}$. It is easy to
check that $j^{\prime}$ is left invariant by the group
$$\Gamma_0(N)=\left\{ \left(\begin{matrix}
              a & b\\
              c & d
             \end{matrix}\right) \in \text{SL}(2, \mathbb{Z}):
             c \equiv 0 (\text{mod $N$}) \right\},$$
which is, in fact, precisely the stabilizer of $j^{\prime}$.

  On the other hand, $j^{\prime}$ is meromorphic at the cusps of $\Gamma_0(N)$.
Indeed, by means of the action of $\Gamma(1)=\text{SL}(2, \mathbb{Z})$, one
reduces the situation to the cusp $\infty$ and to a function of the form
$j \circ \left(\begin{matrix} a & b\\ c & d \end{matrix}\right)$ with $a$,
$b$ and $d$ integers; for sufficiently large $k$ the product of the latter
with $q^{k/m}$ is bounded in a neighborhood of $q^{1/m}=0$. The extension
$K(\Gamma_0(N))/\mathbb{C}(j)$ being finite, this implies the existence of
an algebraic relation between $j$ and $j^{\prime}$. In order to exhibit such
a relation, one considers the transforms of $j^{\prime}$ by the elements of
$\Gamma(1)$, that is, the $j \circ \alpha$ with $\alpha$ ranging over the
orbit $O_N$ of the point
$$p_N=\Gamma(1) \left(\begin{matrix} N & 0\\ 0 & 1 \end{matrix}\right)
      \in \Gamma(1) \backslash \Delta_N,$$
under the action of $\Gamma(1)$ on the right; here $\Delta_N$ denotes the
set of integer matrices of determinant $N$. One can check that the stabilizer
of the point $p_N$ in $\Gamma(1)$ is $\Gamma_0(N)$, so that the orbit $O_N$
may be identified with the quotient $\Gamma_0(N) \backslash \Gamma(1)$.
Denote $d_N$ for the index of $\Gamma_0(N)$ in $\Gamma(1)$ and
$\alpha_k \in \Delta_N$ $(k=1, \ldots, d_N)$ for a system of representatives
of the orbit $O_N$. Then the coefficients of the polynomial
$\prod_{k=1}^{d_N} (X-j \circ \alpha_k)$ are invariant under $\Gamma(1)$,
holomorphic on $\mathbb{H}$ and meromorphic at the cusp $\infty$. We have
thus found a polynomial $\Phi_N \in \mathbb{C}[X, Y]$ of degree $d_N$ in $X$
such that
$$\Phi_N(j^{\prime}, j)=0.\eqno{(3.3.1)}$$
This is the modular equation associated with transformations of order $N$.
The stabilizer of $j \circ \alpha_k$ is conjugate to $\Gamma_0(N)$ (the
stabilizer of $j^{\prime}$), whence the subgroup fixing all the $j \circ
\alpha_k$ coincides with $\Gamma(N)=\cap_{\gamma \in \Gamma(1)} \gamma
\Gamma_0(N) \gamma^{-1}$. It follows that the splitting field of
$\Phi_N \in \mathbb{C}[j][X]$ is $K(\Gamma(N))$. Moreover $\Gamma(1)$
acts as a set of automorphisms of $K(\Gamma(N))$ in permuting transitively
the roots of this polynomial, which is therefore irreducible, whence, in
particular, $K(\Gamma_0(N))=\mathbb{C}(j, j^{\prime})$. When $N=p$ is a
prime, one sees that the matrices $\left(\begin{matrix} 1 & k\\ 0 & p
\end{matrix}\right)$ $(0 \leq k < p)$ and $\left(\begin{matrix} p & 0\\
0 & 1 \end{matrix}\right)$ form a system of representatives of
$O_p=\Gamma(1) \backslash \Delta_p$; the index of $\Gamma_0(p)$ is thus
$d_p=p+1$.

  An elementary calculation shows that $\Gamma_0(N)$ is normalized by the
matrix
$$\left(\begin{matrix}
  0 & N^{-1/2}\\
  -N^{1/2} & 0
  \end{matrix}\right)$$
which induces an involutary automorphism of the modular curve $X_0(N)$ and
its function field: this is the Fricke involution interchanging $j$ and
$j^{\prime}$. One infers from its existence that $\Phi_N \in \mathbb{C}[X, Y]$
is symmetric. Klein relies on this symmetry in his investigation of the modular
equation for $N=2$, $3$, $4$, $5$, $7$ and $13$. For these values of $N$ the
modular curve $X_0(N)$ is of genus zero and there exists $\tau \in K(\Gamma_0(N))$
such that $K(\Gamma_0(N))=\mathbb{C}(j, j^{\prime})=\mathbb{C}(\tau)$; one then
has $j=F(\tau)$ and $j^{\prime}=F(\tau^{\prime})$ with $F \in \mathbb{C}(Z)$, the
function $\tau^{\prime}$ being linked to $\tau$ by the Fricke involution. In each
of these cases Klein describes a fundamental region for the action of $\Gamma_0(N)$
on the half-plane, then deduces from ramification data an expression for $F$ and
gives the relation between $\tau$ and $\tau^{\prime}$. Note that for $N \in \{2,
3, 4, 5 \}$, the modular curve $X(N)$ is also of genus zero, with respective
automorphism groups (leaving the set of cusps globally fixed) the dihedral group,
the tetrahedral group $A_4$, the group $S_4$ of the cube and the octahedron, and
the group $A_5$ of the dodecahedron and the icosahedron. Over the two year period
1878-1879, Klein published a series of papers on modular equations, devoted
respectively to transformations of order $p=5$, $7$ and $11$. In each case he
constructs by geometric means a Galois resolvent, gives its roots explicitly-using
modular forms-and shows how to find the modular equation itself of degree $p+1$
as well as a resolvent of degree $p$ for each of these particular values of $p$.

  In fact, the complex function field of $X=X_0(p)$ with $p$ a prime
consists of the modular functions $f(z)$ for $\Gamma_0(p)$ which are
meromorphic on the extended upper half-plane. A function $f$ lies in
the rational function field $\mathbb{Q}(X)$ if and only if the Fourier
coefficients in its expansion at $\infty$: $f(z)=\sum a_n q^n$ are all
rational numbers. The field $\mathbb{Q}(X)$ is known to be generated over
$\mathbb{Q}$ by the classical $j$-functions
$$\left\{\aligned
  j &=j(z)=q^{-1}+744+196884 q+\cdots,\\
  j_p &=j \left(\frac{-1}{pz}\right)=j(pz)=q^{-p}+744+\cdots.
\endaligned\right.$$
A further element in the function field $\mathbb{Q}(X)=\mathbb{Q}(j,
j_p)$ is the modular unit $u=\Delta(z)/\Delta(pz)$ with divisor
$(p-1) \{ (0)-(\infty) \}$ where $\Delta(z)$ is the discriminant.
If $m=\text{gcd}(p-1, 12)$, then an $m$th root of $u$ lies in
$\mathbb{Q}(X)$. This function has the Fourier expansion
$$t=\root m \of{u}=q^{(1-p)/m} \prod_{n \geq 1} \left(\frac{1-q^n}
    {1-q^{np}}\right)^{24/m}=\left(\frac{\eta(z)}{\eta(pz)}\right)^{24/m}.$$
When $p-1$ divides $12$, so $m=p-1$, the function $t$ is a Hauptmodul
for the curve $X$ which has genus zero. It is well-known that the genus
of the modular curve $X$ for prime $p$ is zero if and only if
$p=2, 3, 5, 7, 13$. In his paper, Klein studied the modular equations
of orders $2$, $3$, $5$, $7$, $13$ with degrees $3$, $4$, $6$, $8$, $14$,
respectively. They are uniformized, i.e., parameterized, by so-called
Hauptmoduln (principal moduli). A Hauptmodul is a function $J_{\Gamma}$
that is a modular function for some subgroup of $\Gamma(1)=\text{SL}(2, \mathbb{Z})$,
with any other modular function expressible as a rational function of it. In this
case, $\Gamma=\Gamma_0(p)$.

  The curve $X_0(p)$ can be given as a plane curve by the modular polynomial
$\Phi_p(X, Y)$. These can quickly get very complicated. For $p=2$ we have
$$\aligned
  \Phi_2(X, Y)=&X^3+Y^3-X^2 Y^2+1488 (X^2 Y+X Y^2)-162000 (X^2+Y^2)\\
               &+40773375 XY+8748000000 (X+Y)-157464000000000,
\endaligned$$
where $X$, $Y$ are the $j$-invariants of the two elliptic curves involved.
It is not so easy to guess that this is a genus $0$ curve. Hence, it is much
better, for conceptual understanding, to parameterize the curve by a different
modular function, and then write $X$ and $Y$ in terms of the parameter. This
leads to the Hauptmodul
$$\tau=\left(\frac{\eta(z)}{\eta(2z)}\right)^{24}.$$
Here,
$$X=j(z)=\frac{(\tau+256)^3}{\tau^2}=\frac{(\tau^{\prime}+16)^3}{\tau^{\prime}},$$
$$Y=j(2z)=\frac{(\tau+16)^3}{\tau}.$$
The Fricke involution involution $\tau^{\prime}=\tau(-\frac{1}{2z})$ satisfies
that
$$\tau \tau^{\prime}=2^{12}.$$

  For $p=3$, the modular polynomial is given by
$$\aligned
  \Phi_3(X, Y)=&X^4+Y^4-X^3 Y^3+2232 (X^3 Y^2+X^2 Y^3)\\
               &-1069956 (X^3 Y+X Y^3)+36864000 (X^3+Y^3)\\
               &+2587918086 X^2 Y^2+8900222976000 (X^2 Y+X Y^2)\\
               &+452984832000000 (X^2+Y^2)-770845966336000000 XY\\
               &+1855425871872000000000(X+Y),
\endaligned$$
We have
$$X=j(z)=\frac{(\tau+27) (\tau+243)^3}{\tau^3}
        =\frac{({\tau^{\prime}}+27) (\tau^{\prime}+3)^3}{\tau^{\prime}},$$
$$Y=j(3z)=\frac{(\tau+27) (\tau+3)^3}{\tau},$$
where the Hauptmodul is given by
$$\tau=\left(\frac{\eta(z)}{\eta(3 z)}\right)^{12},$$
and the Fricke involution $\tau^{\prime}=\tau(-\frac{1}{3z})$ satisfies
that
$$\tau \tau^{\prime}=3^6.$$

  For $p=5$, the geometric model of $X(5)$ he uses is the regular icosahedron,
the resolvent of degree five being linked, as had been shown by Hermite,
to the general quintic equation. In fact, Klein shows that the morphism
$X(5) \rightarrow X(1)$ is isomorphic to that taking the quotient of the
unit sphere in $\mathbb{R}^3$ by the action of the symmetry group of the
regular icosahedron. More precisely, in the modular equation $\Phi_5(X, Y)=0$
of level five, we have
$$X=j(z)=-\frac{(\tau^2-250 \tau+3125)^3}{\tau^5}
      =-\frac{({\tau^{\prime}}^2-10 \tau^{\prime}+5)^3}{\tau^{\prime}},\eqno{(3.3.2)}$$
$$Y=j(5z)=-\frac{(\tau^2-10 \tau+5)^3}{\tau},\eqno{(3.3.3)}$$
where the Hauptmodul is given by
$$\tau=\left(\frac{\eta(z)}{\eta(5 z)}\right)^6,\eqno{(3.3.4)}$$
and the Fricke involution $\tau^{\prime}=\tau(-\frac{1}{5z})$ satisfies
that
$$\tau \tau^{\prime}=125.$$

  For $p=7$, Klein shows that the modular curve $X(7)$ is isomorphic to
the smooth plane quartic curve $C_4$ with equation $x^3 y+y^3 z+z^3 x=0$
in $\mathbb{CP}^2$, invariant under the action of a group $G$ isomorphic
to $\text{PSL}(2, 7)$, which is the automorphism group of $X(7)$. In this
projective model the natural morphism from $X(7)$ onto
$X(1) \simeq \mathbb{CP}^1$ is made concrete as the projection of $C_4$ on
$C_4/G$ (identified with $\mathbb{CP}^1$). This is a Galois covering whose
generic fibre is considered by Klein as the Galois resolvent of the modular
equation $\Phi_7(\cdot, j)=0$ of level $7$, which means that the function
field of $C_4$ is the splitting field of this modular equation over $\mathbb{C}(j)$.
Moreover, just as the sphere has a regular tiling induced from the faces of
an inscribed icosahedron, so also does the modular curve $X(7)$ admits a
regular tiling by triangles. This tiling is inherited combinatorially from
a tiling of $\mathbb{H}$ of type $(2, 3, \infty)$, and its triangles are
of type $(2, 3, 7)$. Here, a tiling of $\mathbb{H}$ by triangles is said
to be of type $(a, b, c)$ if it is realized by hyperbolic triangles
$(a, b, c)$, that is, with angles $(\frac{2 \pi}{a}, \frac{2 \pi}{b},
\frac{2 \pi}{c})$. More precisely, in the modular equation $\Phi_7(X, Y)=0$
of level seven, we have
$$\aligned
  X=j(z) &=\frac{(\tau^2+13 \tau+49) (\tau^2+245 \tau+2401)^3}{\tau^7}\\
         &=\frac{({\tau^{\prime}}^2+13 \tau^{\prime}+49)({\tau^{\prime}}^2
         +5 \tau^{\prime}+1)^3}{\tau^{\prime}},
\endaligned\eqno{(3.3.5)}$$
$$Y=j(7z)=\frac{(\tau^2+13 \tau+49)(\tau^2+5 \tau+1)^3}{\tau},\eqno{(3.3.6)}$$
where the Hauptmodul is given by
$$\tau=\left(\frac{\eta(z)}{\eta(7 z)}\right)^4,\eqno{(3.3.7)}$$
and the Fricke involution $\tau^{\prime}=\tau(-\frac{1}{7z})$ satisfies
that
$$\tau \tau^{\prime}=49.$$

  Now, we come to the group $\text{PSL}(2, 13)$ and the modular curves
$X(13)$. The function field $K(\Gamma(13))$ of the modular curve $X(13)$
is the splitting field of the polynomial $\Phi_{13} \in \mathbb{C}(j)[X]$
associated with transformations of order $13$. Let
$\mathbb{F}_{13}=\mathbb{Z}/13 \mathbb{Z}$ be the field of thirteen elements.
Since $\text{SL}(2, 13)$ is generated by
$\left(\begin{matrix} 1 & 1\\ 0 & 1 \end{matrix}\right)$ and
$\left(\begin{matrix} 1 & 0\\ 1 & 1 \end{matrix}\right)$, the reduction
morphism modulo $13$ from $\text{SL}(2, \mathbb{Z})$ to $\text{SL}(2, 13)$
is surjective, whence the exact sequence
$$1 \rightarrow \overline{\Gamma}(13) \rightarrow \text{PSL}(2, \mathbb{Z})
  \rightarrow \text{PSL}(2, 13) \rightarrow 1,$$
where $\overline{\Gamma}:=\Gamma/ \{ \pm I \}$. In particular, the quotient
$\overline{G}=\text{PSL}(2, \mathbb{Z})/\overline{\Gamma}(13)$ is isomorphic
to $\text{PSL}(2, 13)$, a simple group of order $1092$. The group $\overline{G}$
acts on $X(13)$ via automorphisms and $X(13)/\overline{G}$ can be identified
with $X(1)$. Thus, the fibres of the projection $X(13) \rightarrow X(1)$ are
the orbits of the action of $\overline{G}$ on $X(13)$. There are therefore
three singular fibres corresponding to the values $J=\infty$, $0$ and $1$
(recall that $J=j/1728$), that is, $j=\infty$, $0$ and $1728$, whose elements
are called $A$-points, $B$-points and $C$-points in Klein's terminology, with
stabilizers of orders $13$, $3$ and $2$, respectively. These fibres have
cardinality $84$, $364$ and $546$; all others have $1092$ elements. By the
Riemann-Hurwitz formula, the genus of $X(13)$ satisfies the relation
$$2-2g=2 \cdot 1092-12 \cdot 84-2 \cdot 364-546,$$
whence $g=50$. Hence, the modular curve $X(13)$ admits a regular tiling by
triangles. This tiling is inherited combinatorially from a tiling of
$\mathbb{H}$ of type $(2, 3, \infty)$, and its triangles are of type
$(2, 3, 13)$.

  The following facts about the modular subgroups of level $13$ should
be noted. One has $\overline{\Gamma}(13)<\overline{\Gamma}_0(13)<\overline{\Gamma}(1)$,
the respective subgroup indices being $78=6 \cdot 13$ and $14=13+1$. The
quotient $\overline{\Gamma}(1)/\overline{\Gamma}(13)$ is of order
$1092=78 \cdot 14$, and $\overline{\Gamma}_0(13)/\overline{\Gamma}(13)$
is a subgroup of order $78$, which is isomorphic to a semidirect product
of $\mathbb{Z}_{13}$ by $\mathbb{Z}_6$. The respective quotients of
$\mathbb{H}$ by these three groups (compactified) are the modular curves
$X(13)$, $X_0(13)$, $X(1)$, and the coverings $X(13) \rightarrow X_0(13)
\rightarrow X(1)$ are respectively $78$-sheeted and $14$-sheeted. The curve
$X_0(13)$ like $X(1)$ being of genus zero.

  In the modular equation $\Phi_{13}(X, Y)=0$ of level thirteen, we have
$$\aligned
  X=j(z) &=\frac{(\tau^2+5\tau+13)(\tau^4+247 \tau^3+3380 \tau^2
         +15379 \tau+28561)^3}{\tau^{13}}\\
         &=\frac{({\tau^{\prime}}^2+5 \tau^{\prime}+13)({\tau^{\prime}}^4
         +7 {\tau^{\prime}}^3+20 {\tau^{\prime}}^2+19 \tau^{\prime}+1)^3}
         {\tau^{\prime}},
\endaligned\eqno{(3.3.8)}$$
and
$$Y=j(13z)=\frac{(\tau^2+5 \tau+13)(\tau^4+7 \tau^3+20 \tau^2+19 \tau+1)^3}{\tau},
  \eqno{(3.3.9)}$$
where the Hauptmodul is given by
$$\tau=\left(\frac{\eta(z)}{\eta(13z)}\right)^2,\eqno{(3.3.10)}$$
and the Fricke involution $\tau^{\prime}=\tau(-\frac{1}{13 z})$ satisfies
that
$$\tau \tau^{\prime}=13.$$

  In section 12, we will show that the modular curve $X(13)$ is isomorphic to
the algebraic space curve $Y$ with the $21$ equations in $\mathbb{CP}^5$,
invariant under the action of a group $G$ isomorphic to $\text{SL}(2, 13)$,
which is the automorphism group of $X(13)$. In this projective model the
natural morphism from $X(13)$ onto $X(1) \simeq \mathbb{CP}^1$ is made
concrete as the projection of $Y$ onto $Y/G$ (identified with $\mathbb{CP}^1$).
This is a Galois covering whose generic fibre is considered as the Galois
resolvent of the modular equation $\Phi_{13}( \cdot, j)=0$ of level $13$,
which means that the function field $K(Y)$ of $Y$ is the splitting field of
this modular equation over $\mathbb{C}(j)$. This proves the following:

\textbf{Theorem 3.3.2.} (Galois covering, Galois resolvent and their geometric
realizations). {\it In the projective model, $X(13)$ is isomorphic to $Y$ in
$\mathbb{CP}^5$. The natural morphism from $X(13)$ onto $X(1) \simeq \mathbb{CP}^1$
is realized as the projection of $Y$ onto $Y/G$ $($identified with $\mathbb{CP}^1$$)$.
This is a Galois covering whose generic fibre is interpreted as the Galois
resolvent of the modular equation $\Phi_{13}( \cdot, j)=0$ of level $13$, i.e.,
the function field of $Y$ is the splitting field of this modular equation
over $\mathbb{C}(j)$.}

\begin{center}
{\bf \S 3.4. The Galois groups of the special transformation equations}
\end{center}

  In his lectures on elliptic functions (see \cite{Fricke}, Vol. II, pp.
459-491), Fricke studied the Galois groups of the special transformation
equations and the three resolvents of degrees $5$, $7$ and $11$. Fricke
pointed out that this topic comes from the solution of the special division
equations of the $\wp$-function and he intended to give the continuation
of these developments for the special transformation equations. It is known
that there is the theorem that if the division values of the prime division
degrees $n$ are known, the division values of all further degrees can be
calculated by rational calculations and root extractions alone. However,
in all prime number cases $n>3$, the special division equations of the
$\wp$-function could no longer be solved by root extractions alone, so
that it is precisely these cases $n$ that are of further interest.

  In fact, the special transformation equations can be regarded as
resolvents of the special division equations. It will now be a question
of developing more precisely the algebraic theory of these equations,
whereby we can restrict ourselves to the prime cases $n>3$ according to what
has been said. The first question will be that of the Galois group of the
special transformation equation. By then going into the structure of this
group in more detail, we will address the question whether the special
transformation equations, which have degree $(n+1)$ for the single $n$,
are the resolvents of the lowest degree of the special division equation
or not. This is answered by a famous theorem discovered by Galois, according
to which for $n>11$ the transformation equations are the lowest resolvents,
but that in the first three cases $n=5$, $7$ and $11$ resolvents of the
$n$-th degree exist. The theorem is given by Galois in his letter to A. Chevalier
of 29 May 1832; cf. the collection of Galois' papers in the Journ. de Math.,
Vol. 11 (1846) (see \cite{Ga}).

  The transition from the division equations to the transformation
equations was accomplished by the considerations from \cite{Fricke},
Vol. II, p. 335 et seq. We arrange the $\frac{1}{2}(n^2-1)$ division
values $\wp_{\lambda, \mu}$, of prime degree $n$ in $(n+1)$ systems:
$$\wp_{\lambda, \mu}, \wp_{2 \lambda, 2 \mu}, \wp_{3 \lambda, 3 \mu},
  \cdots, \wp_{\frac{n-1}{2} \lambda, \frac{n-1}{2} \mu} \eqno{(3.4.1)}$$
to $\frac{1}{2}(n-1)$ each, which are inverse to the substitutions of
the Galois group $G_{\frac{1}{2} n (n-1) (n^2-1)}$ of the division
equation. We were able to build up this group from all incongruent
substitutions
$$\lambda^{\prime} \equiv \alpha \lambda+\gamma \mu, \quad
  \mu^{\prime} \equiv \beta \lambda+\delta \mu, \quad
  \text{(mod $n$)} \eqno{(3.4.2)}$$
to be exercised on the indices $\lambda$, $\mu$, whose determinants
$(\alpha \delta-\beta \gamma)$ were not divisors of $n$, and for which
two substitutions merging into one another by simultaneous sign change
of $\alpha$, $\beta$, $\gamma$, $\delta$ were considered as not different
(see \cite{Fricke}, Vol. II, p. 261). By the individual of these substitutions,
the systems (3.4.1), apart from rearrangements of the $\wp$-division values
in the individual systems, are in fact only permuted among themselves. We
can characterise the individual system (3.4.1) with $\lambda \neq 0$ by
that number $\kappa$ of the series $0$, $1$, $2$, $\ldots$, $n-1$ which
satisfies the congruence $\kappa \lambda \equiv \mu$ (mod $n$), a number
which we may also denote by $\frac{\mu}{\nu}$ or $\mu \cdot \lambda^{-1}$
in a known manner. For the notation of the system (3.4.1) with $\lambda=0$
we use the symbol $\kappa=\infty$ accordingly.

  As the most important example of the transformation equations, we now
consider that of $j(\omega)=12^3 J(\omega)$. According to \cite{Fricke}, Vol. II,
p. 338, the transformed function $j(n \omega)$ occurring in the first principal
transformation is a symmetric function of the $\frac{1}{2}(n-1)$ division
values of the system (3.4.1) belonging to $\kappa=\infty$ with coefficients
of the field $(\mathfrak{R}, g_2, g_3)$. We write
$$j_{\infty}=j(n \omega)=R\left(\wp_{01}, \wp_{02}, \cdots,
  \wp_{0, \frac{n-1}{2}}\right).\eqno{(3.4.3)}$$
with the inclusion of $\kappa$ as index. This is followed by the other
functions
$$j_{\kappa}=j\left(\frac{\omega+\kappa}{n}\right)=R\left(\wp_{1, \kappa},
  \wp_{2, 2 \kappa}, \wp_{3, 3 \kappa}, \cdots, \wp_{\frac{n-1}{2},
  \frac{n-1}{2} \kappa}\right) \eqno{(3.4.4)}$$
with $\kappa=0$, $1$, $2$, $\ldots$, $n-1$, which are conjugate among
themselves and with $j_{\infty}$ and correspond to the other systems (3.4.1),
where $R$ has the same meaning as in (3.4.3). In particular, $j_0$ corresponds
to the second principal transformation of the $n$-th degree.

  In $j_{\infty}$, $j_0$, $j_1$, $\ldots$, $j_{n-1}$ we now have before us
the roots of the transformation equation, to which the general principles
of Galois theory concerning the groups of resolvents are now to be applied
(cf. \cite{Fricke}, Vol. II, p. 52 et seq.). For this purpose we first have
to determine the subgroup of those substitutions (3.4.2) which cause the identical
permutation of $j_{\infty}$, $j_0$, $j_1$, $\ldots$, $j_{n-1}$. Now $j_{\infty}$
is transformed into itself by the substitutions (3.4.2) with $\gamma \equiv 0$
(mod $n$), but $j_0$ by the substitutions with $\beta \equiv 0$ (mod $n$).
Similarly, we find $\alpha \equiv \delta$ (mod $n$) as the condition for
$j_1$ to transform into itself. The substitutions (3.4.2), which satisfy these
three conditions, transform all roots (3.4.3) and (3.4.4) into themselves and
thus provide the average of all $(n+1)$ subgroups belonging to the individual
roots: {\it The substitutions of the Galois group $G_{\frac{1}{2}n(n-1)(n^2-1)}$
of the division equation, which provide the identical permutation of the
$j_{\infty}$, $j_0$, $j_1$, $\ldots$, $j_{n-1}$ form the normal subgroup
$G_{\frac{1}{2}(n-1)}$ of the $\frac{1}{2}(n-1)$ substitutions:}
$$\lambda^{\prime} \equiv \alpha \lambda, \quad
  \mu^{\prime} \equiv \alpha \mu \quad (\text{mod $n$}). \eqno{(3.4.5)}$$

  Proceeding to the transformation equation, the reduction of the Galois
group to the quotient group $G_{\frac{1}{2}n(n-1)(n^2-1)}/G_{\frac{1}{2}(n-1)}$
now occurs according to \cite{Fricke}, Vol. II, p. 53, {\it so that the
Galois group of the special transformation equation is a $G_{n(n^2-1)}$ of
order $n(n^2-1)$.} We can produce them from the group (3.4.2) simply by considering
all substitutions (3.4.2) with proportional number quadruples $\alpha$, $\beta$,
$\gamma$, $\delta$ as not different from each other. {\it It is thus sufficient
allow all $\frac{1}{2}n(n^2-1)$ substitutions of $(3.4.2)$ with determinant
$1$ and also all substitutions $(3.4.2)$, whose determinant is equal to a
quadratic non-residue of $n$ to be chosen arbitrarily, i.e. in the
case of $n=4h+3$ approximately equal to the non-residue $-1$.} By
introducing the indices $\kappa \equiv \mu \lambda^{-1}$, we can also
express the obtained theorem in this way:

\textbf{Theorem 3.4.1.} {\it The Galois group of the special transformation
equation consists of all $n(n^2-1)$ permutations of the $j_{\infty}$, $j_0$,
$j_1$, $\ldots$, $j_{n-1}$, which are obtained by the substitutions:}
$$\kappa^{\prime} \equiv \frac{\delta \kappa+\beta}{\gamma \kappa+\alpha}
  \quad \text{(mod $n$)} \eqno{(3.4.6)}$$
{\it to be exercised on the index $\kappa$, where $\alpha \delta-\beta \gamma
\equiv 1$ and congruent to a quadratic non-residue of $n$ to be chosen
arbitrarily, and of two quadruples merging by quadruples $\alpha$, $\beta$,
$\gamma$, $\delta$, which merge into one another by changing signs, only
one is to be allowed.}

  In the Galois group $G_{n(n^2-1)}$ contains${}$ (this is the group of
all substitutions of (3.4.2) with $\alpha \delta -\beta \gamma \equiv 1$ (mod $n$))
the {\it monodromy group} $G_{\frac{1}{2}n(n^2-1)}$ as a normal subgroup
of the index $2$, which in turn is simple in all cases $n>3$ according to
\cite{Fricke}, Vol. II, p. 262 that are in question here. The
$G_{\frac{1}{2}n(n^2-1)}$ is isomorphic with the mod $n$ reduced
non-homogeneous modular group $\Gamma^{(\omega)}$ and can therefore be
generated from two permutations, which correspond to the two substitutions
$S=\left(\begin{matrix} 1 & 1\\ 0 & 1 \end{matrix} \right)$ and
$T=\left(\begin{matrix} 0 & 1\\ -1 & 0 \end{matrix}\right)$.
These two permutations are, as can easily be seen,
$$j_{\infty}^{\prime}=j_{\infty}, \quad
  j_{\kappa}^{\prime}=j_{\kappa+1}, \eqno{(S)}$$
$$j_{\infty}^{\prime}=j_0, \quad j_0^{\prime}=j_{\infty}, \quad
  j_{\kappa}^{\prime}=j_{-\kappa^{-1}}, \eqno{(T)}$$
where for the permutation $S$, of course, $j_n=j_0$ is to be taken.
Then for the Galois group $G_{n(n^2-1)}$ there is added as a third
generating permutation:
$$j_{\infty}^{\prime}=j_{\infty}, \quad j_0^{\prime}=j_0, \quad
  j_{\kappa}^{\prime}=j_{\nu \kappa},$$
where $\nu$ is some quadratic non-residue of $n$.

  For the reduction of the Galois group to the monodromy group, the
adjunction of the unit root $\varepsilon=e^{\frac{2 \pi i}{n}}$ is
sufficient, but not necessary. Since the index of the monodromy group
in the $G_{n(n^2-1)}$ is equal to $2$, the adjunction of a single
numerical irrationality of the second degree is sufficient. Since it
belongs to the cyclotomic field $(\mathfrak{R}, \varepsilon)$, it is
the root of the quadratic resolvent of the cyclotomic equation for the
$n$-th degree of division. This quadratic resolvent is the equation
for the $\frac{1}{2}(n-1)$-membered sum:
$$\varepsilon+\varepsilon^4+\varepsilon^9+\cdots+
  \varepsilon^{\left(\frac{n-1}{2}\right)^2}=\frac{1}{2}\left(-1+
  i^{\frac{n-1}{2}} \sqrt{n}\right), \eqno{(3.4.7)}$$
where on the left in the exponents are the $\frac{1}{2}(n-1)$ quadratic
residues of $n$.

  According to \cite{Fricke}, Vol. II, p. 260, $\varepsilon$ is a natural
irrationality of the $\wp$-division equation, that is, it can be represented
as a rational function of $\wp_{\lambda \mu}$ with coefficients of the field
$(\mathfrak{R}, g_2, g_3)$. This function remains unchanged in the
substitutions of the monodromy group and changes to $\varepsilon^d$ in
the single substitution (3.4.2) of the determinant $\alpha \delta-\beta
\gamma \equiv d$. The sum (3.4.7) is therefore a quantity that is invariant
with respect to the $G_{\frac{1}{2}(n-1)}$ of the substitutions (3.4.5),
i.e. it belongs to the transformation equation as a natural
irrationality according to \cite{Fricke}, Vol. II, p. 49 et seq.
Now, Fricke proved the following (see \cite{Fricke}, Vol. II, p. 462,
also see \cite{Weber}, \S 78, especially, p. 288 and p. 289):

\textbf{Theorem 3.4.2.} {\it The Galois group $G_{n(n^2-1)}$ of the
special transformation equation reduces after adjunction of the
natural irrationality $i^{\frac{n-1}{2}} \sqrt{n}$ to it to the
monodromy group $G_{\frac{1}{2}n(n^2-1)}$ which is isomorphic
with the} mod $n$ {\it reduced non-homogeneous modular group.}

\begin{center}
{\bf \S 3.5. Galois' theorem}
\end{center}

  The Galois group of the transformation equation was, in the prime case
$n$, a more definite $G_{n(n^2-1)}$ in \S 8.4, which, after adjunction of
the square root $i^{\frac{n-1}{2}} \sqrt{n}$, reduced to the $G_{\frac{1}{2}
n (n^2-1)}$ decomposed in \cite{Fricke}, Vol. II, p. 465, \S 3 et seq. The
transformation equation was the resolvent, which belonged to the $(n+1)$
meta-cyclic $G_{\frac{1}{2}n(n-1)}$. The index $(n+1)$ of these groups in
the total group $G_{\frac{1}{2} n(n^2-1)}$ provided the degree $(n+1)$ of
the transformation equation.

  We first hold to the adjunction of $i^{\frac{n-1}{2}} \sqrt{n}$, so that
the group of each resolvent is again $G_{\frac{1}{2}n(n^2-1)}$, since this
group is simple. An equation with a degree $<n$ cannot have a group, whose
order is divisible by the prime number $n$. Accordingly, shall we discuss
the possibility, whether resolvents of degree $<n+1$ occur, then it can at
most be a question of resolvents of the $n$-th degree itself. The answer to
the question raised is contained in Galois' theorem, {\it according to which,
although in the three lowest cases $n=5$, $7$ and $11$, but for no other
prime $n$ do resolvents of the $n$-th degree occur.}

  In \cite{Fricke}, Vol. II, pp. 465-475, Fricke gave a classification of the
subgroups occurring in $G_{\frac{1}{2}n(n^2-1)}$, which consists of cyclic
subgroups, meta-cyclic subgroups, dihedral subgroups and exceptional subgroups.
{\it In $G_{\frac{1}{2}n(n^2-1)}$, $(n+1)$ conjugate cyclic subgroups $G_n$ of
order $n$ occur and, correspondingly, just as many conjugate meta-cyclic subgroups
$G_{\frac{1}{2}n(n-1)}$ of order $\frac{1}{2} n(n-1)$.} The individual $G_n$ as
a group of prime order can be generated from each of its substitutions, apart
from the identical substitution $1$. From this it follows that two different
$G_n$, apart from the substitution $1$, cannot contain any common substitutions.
{\it In total, we therefore find $(n+1)(n-1)=n^2-1$ substitutions of period $n$.}
{\it Apart from the subgroups contained in the $G_{\frac{1}{2}n(n-1)}$ and the
total group $G_{\frac{1}{2}n(n^2-1)}$, the meta-cyclic $G_{\frac{1}{2}n(n-1)}$
are the only groups in which the cyclic $G_n$ are involved.} {\it The cyclic
$G_{\frac{1}{2}(n-1)}$ is normally contained in a most comprehensive $G_{n-1}$
of the dihedral type}, {\it so that we obtain in total $\frac{1}{2}n(n+1)$
cyclic subgroups $G_{\frac{1}{2}(n-1)}$ of order $\frac{1}{2}(n-1)$, one of
which consists of the substitutions}
$$\left(\begin{matrix} g & 0\\ 0 & g^{-1} \end{matrix}\right),
  \left(\begin{matrix} g^2 & 0\\ 0 & g^{-2} \end{matrix}\right),
  \left(\begin{matrix} g^3 & 0\\ 0 & g^{-3} \end{matrix}\right), \cdots,
  \left(\begin{matrix} g^{\frac{1}{2}(n-1)} & 0\\ 0 & g^{-\frac{1}{2}(n-1)}
  \end{matrix}\right) \equiv 1.\eqno{(3.5.1)}$$
{\it No two of the cyclic groups of order $\frac{1}{2}(n-1)$ can have a
substitution in common, apart from substitution $1$.} {\it The
$\frac{1}{2}n(n+1)$ cyclic subgroups $G_{\frac{1}{2}(n-1)}$ contain
different substitutions except $1$ in the whole $\frac{1}{4}n(n+1)(n-3)$.}
{\it For $n>5$ we find $\frac{1}{2}n(n+1)$ conjugate dihedral groups
$G_{n-1}$, one of which consists of the substitutions $(3.5.1)$ and}
$$\left(\begin{matrix} 0 & -g\\ g^{-1} & 0 \end{matrix}\right),
  \left(\begin{matrix} 0 & -g^2\\ g^{-2} & 0 \end{matrix}\right),
  \left(\begin{matrix} 0 & g^{-\frac{1}{2}(n-1)}\\ g^{-\frac{1}{2}(n-1)} & 0
  \end{matrix}\right) \equiv \left(\begin{matrix} 0 & 1\\ -1 & 0 \end{matrix}
  \right).\eqno{(3.5.2)}$$
For $n=5$ there is an exceptional case. {\it We find $\frac{1}{2}n(n-1)$
conjugate cyclic groups $G_{\frac{1}{2}(n+1)}$ of order $\frac{1}{2}(n+1)$
and, corresponding to them, just as many conjugate dihedral groups $G_{n+1}$,
each of which is the most comprehensive group in which the corresponding
$G_{\frac{1}{2}(n+1)}$ is normally contained.} {\it Two cyclic
$G_{\frac{1}{2}(n+1)}$ cannot have any substitution in common apart from
the identical substitution $1$.} We then count off {\it that in the
$G_{\frac{1}{2}(n+1)}$, apart from the identical substitution, contain
$\frac{1}{4}n(n-1)^2$ different substitutions in total.}

  Including the identical substitution, there are now a total of
$$1+(n^2-1)+\frac{1}{4}n(n+1)(n-3)+\frac{1}{4}n(n-1)^2=\frac{1}{2}n(n^2-1)$$
different substitutions, i.e. all the substitutions of our $G_{\frac{1}{2}n(n^2-1)}$
are categorised into cyclic groups. If we add the cyclic subgroups contained
in $G_{\frac{1}{2}(n-1)}$ and $G_{\frac{1}{2}(n+1)}$, we have obtained all
the cyclic subgroups contained in $G_{\frac{1}{2}n(n^2-1)}$ at all. The cyclic
subgroups contained in $G_{\frac{1}{2}(n-1)}$ and $G_{\frac{1}{2}(n+1)}$ also
give rise to the formation of dihedral groups by combining them with substitutions
of period $2$ from the corresponding dihedral $G_{n-1}$ or $G_{n+1}$. In this
way, all subgroups of the dihedral type contained in the $G_{\frac{1}{2}n(n^2-1)}$
are obtained.

  The subgroups of $G_{\frac{1}{2}n(n^2-1)}$ that are still missing can be
set up on the basis of a consideration developed by C. Jordan. The
subgroups still contained in the meta-cyclic $G_{\frac{1}{2}n(n-1)}$, which
can easily be obtained from the cyclic subgroups of $G_{\frac{1}{2}(n-1)}$,
may be left aside, since they will still only be ``most comprehensive''
subgroups. We therefore only have only have to search for those subgroups
in which the cyclic $G_n$ are not involved. {\it Apart from the groups
mentioned as above, only three types of subgroups $G_{12}$, $G_{24}$ and
$G_{60}$ of orders $12$, $24$ and $60$ can occur in $G_{\frac{1}{2}n(n^2-1)}$,
which correspond to the following solutions of equation}
$$\frac{m-1}{m}=\frac{\mu_1-1}{e_1 \mu_1}+\frac{\mu_2-1}{e_2 \mu_2}+\cdots+
                \frac{\mu_k-1}{e_k \mu_k}:\eqno{(3.5.3)}$$
$$G_{12}, k=2, e_1=2, e_2=1, \mu_1=2, \mu_2=3, m=12,\eqno{(3.5.4)}$$
$$G_{24}, k=3, e_1=e_2=e_3=2, \mu_1=2, \mu_2=3, \mu_3=4, m=24,\eqno{(3.5.5)}$$
$$G_{60}, k=3, e_1=e_2=e_3=2, \mu_1=2, \mu_2=3, \mu_3=5, m=60.\eqno{(3.5.6)}$$
The three permutation groups obtained are now known to be isomorphic
with the tetrahedral group, the octahedral group and the icosahedral
group. Incidentally, since the order $\frac{1}{2}n(n^2-1)$ of the
total group $G_{\frac{1}{2}n(n^2-1)}$ must of course be divisible by
$24$ or $60$ when subgroups $G_{24}$ and $G_{60}$ occur, the final
theorem is as follows: {\it Apart from the groups already mentioned as
above, only tetrahedral groups $G_{12}$, octahedral groups $G_{24}$ and
icosahedral groups $G_{60}$ can occur in $G_{\frac{1}{2}n(n^2-1)}$,
whereby the occurrence of octahedral groups requires prime numbers $n$
of the form $n=8h \pm 1$, while that of icosahedral groups requires
prime numbers of the form $n=10 h \pm 1$}. An exception to the last
statement concerning $G_{60}$ only exists in the lowest case $n=5$,
{\it insofar as the total group $G_{\frac{1}{2}n(n^2-1)}$ belonging
to $n=5$ is itself an icosahedral group.} As above, this is easily
deduced from the occurrence of cyclic groups and groups of four in
$G_{60}$.

  The negative part of the theorem is immediately apparent from the
results mentioned as above. The cyclic subgroups and the dihedral groups
of $G_{\frac{1}{2}n(n^2-1)}$ never reach the order of the meta-cyclic
$G_{\frac{1}{2}n(n-1)}$. If, therefore, even subgroups whose order is
$>\frac{1}{2} n(n-1)$ are to occur, they would have to belong to the
$G_{12}$, $G_{24}$ or $G_{60}$ established as above. If $n \geq 13$,
then $\frac{1}{2}n(n-1) \geq 78$ holds, so that the $G_{\frac{1}{2}n(n-1)}$
are certainly the subgroups of highest order. In contrast, for $n=5$
the order $12$, for $n=7$ the order $24$ and for $n=11$ the order $60$
is greater than the respective order $10$, $21$ and $55$ of the meta-cyclic
groups.

   In the lowest case $n=5$, the actual occurrence of groups $G_{12}$ is
easy to prove. {\it There are actually five conjugate $G_{12}$ in the
$G_{60}$, which are tetrahedral groups. These $G_{12}$ correspond to five
conjugate congruence groups of the fifth level $\Gamma_5$ of index $5$.}
{\it $G_{120}$ contains five conjugate subgroups $G_{24}$ of index $5$.}
The $G_{168}$ belonging to $n=7$ contains a total of $21$ cyclic $G_4$,
all of which are conjugate, and the $21$ conjugate $G_2$ within them.
Since the individual $G_4$ can now be combined with $24$ substitutions
forming an octahedron-$G_{24}$, {\it the $G_{168}$ in fact also contains
two systems of seven conjugate octahedral groups $G_{24}$.} These groups
correspond to {\it two systems of seven conjugate congruence groups of
the $7$-th level $\Gamma_7$ of index $7$}. The search for any icosahedron-$G_{60}$
in the $G_{660}$ belonging to $n=11$ should be preceded by the following
consideration concerning the generation of the $G_{60}$, in which we use
the $G_{60}$ in the form of the $60$ mod $5$ incongruent substitutions
$\left(\begin{matrix} \alpha & \beta\\ \gamma & \delta \end{matrix}\right)$
of determinant $1$. {\it The individual $G_{60}$ is always one of eleven
conjugate subgroups.} {\it If icosahedra $G_{60}$ occur at all in the
$G_{660}$, then there is either a system of eleven conjugate $G_{60}$ or
two such systems.} {\it This means that in $G_{660}$ there are in fact two
systems of eleven conjugate icosahedral groups $G_{60}$ each.} {\it The
$G_{60}$ correspond to two systems of eleven conjugate congruence groups
$\Gamma_{11}$ of the eleventh degree of index $11$.} We arrive at the theorem,
which we immediately express again for Galois $G_{n(n^2-1)}$: {\it Within the
Galois group $G_{336}$ there is a system of $14$ conjugate octahedral groups
$G_{24}$ of index $14$, and likewise within Galois $G_{1320}$ there is a
system of $22$ conjugate icosahedral groups $G_{60}$ of index $22$.} In
accordance with the general theorems of Galois equation theory, we now
immediately draw the algebraic conclusions from this: {\it First, in
the case of $n=5$ we arrive at a numerically rational resolvent of the
transformation equation of the sixth degree.} If we disregard the
adjunction of the irrationality $i \sqrt{n}$ for the time being in the
case of $n=7$ and $11$, the following continues:

\textbf{Theorem 3.5.1.} {\it For the transformation degrees $n=7$ and $11$,
a resolvent of degree $14$ and $22$, respectively, with rational number
coefficients occurs.}

  However, when proceeding to the monodromy groups and thus after adjunction
of $i \sqrt{7}$ or $i \sqrt{11}$, these resolvents become reducible, and one
arrives at the theorem:

\textbf{Theorem 3.5.2.} {\it For the degrees of transformation $n=7$ and
$n=11$ there are two resolvents each of the degrees $7$ and $11$, in whose
coefficients the irrationalities $i \sqrt{7}$ and $i \sqrt{11}$ occur.}

  Now, we will give a modern viewpoint according to Serre (\cite{Se1972})
and Mazur (\cite{Mazur1977}). Let $p$ be a prime number, if
$H \subset \text{GL}(2, \mathbb{F}_p)$ is any subgroup such that
$\det H=\mathbb{F}_p^{*}$ there is a projective curve $X_H$ over
$\mathbb{Q}$ parameterizing elliptic curves with level $H$-structures.
The determination of the rational points of $X_H$ amounts to a
classification of elliptic curves over $\mathbb{Q}$ satisfying the
property that the associated representation of
$\text{Gal}(\overline{\mathbb{Q}}/\mathbb{Q})$ on $p$-division points
factors through a conjugate of $H$ (see \cite{Mazur1977}). It is well-known
that if $p \geq 5$ is a prime number, any proper subgroup
$H \subset \text{GL}(2, \mathbb{F}_p)$ is conjugate to a subgroup of one
of the entries of the following table, in which $\mathfrak{S}_n$ denotes
the symmetric group on $n$ letters, $\mathfrak{A}_n$ the alternating group
and $\mathbb{Q}(\chi_p)$ denotes the quadratic subfield of
$\mathbb{Q}(e^{\frac{2 \pi i}{p}})$.
$$\begin{matrix}
  H & \text{Notation for $X_p$} & \text{Field of definition}\\
  \text{(1) The Borel subgroup} & X_0(p) & \mathbb{Q}\\
  \left(\begin{matrix} * & *\\ 0 & * \end{matrix}\right) &  & \\
  \text{(2) The normalizer of} & X_{\text{split}}(p) & \mathbb{Q}\\
  \text{a split Cartan subgroup} &  & \\
  \left(\begin{matrix} * & 0\\ 0 & * \end{matrix}\right) \sqcup
  \left(\begin{matrix} 0 & * \\ * & 0 \end{matrix}\right)  &  & \\
  \text{(3) The normalizer of} & X_{\text{non-split}}(p) & \mathbb{Q}\\
  \text{a non-split Cartan subgroup} &  & \\
  \mathbb{F}_{p^2}^{*} \subset \text{GL}(2, \mathbb{F}_p) &  &  \\
  \text{(4) The inverse image in} & X_{\mathfrak{S}_4}(p) & \text{$\mathbb{Q}$ if
  $p \equiv \pm 3$ (mod $8$)}\\
  \text{$\text{GL}(2, \mathbb{F}_p)$ of $\mathfrak{S}_4 \subset \text{PGL}(2, \mathbb{F}_p)$}
  &  & \text{$\mathbb{Q}(\chi_p)$ if otherwise}\\
  \text{(5) The inverse image in} & X_{\mathfrak{A}_5}(p) & \mathbb{Q}(\chi_p)\\
  \text{$\text{GL}(2, \mathbb{F}_p)$ of $\mathfrak{A}_5 \subset
  \text{PGL}(2, \mathbb{F}_p)$} & & \\
  \text{(possibly only if $p \equiv \pm 1$ (mod $5$))}  &  & \\
  \text{(6) The inverse image in} & X_{\mathfrak{A}_4}(p) & \mathbb{Q}(\chi_p)\\
  \text{$\text{GL}(2, \mathbb{F}_p)$ of $\mathfrak{A}_4 \subset
  \text{PGL}(2, \mathbb{F}_p)$} & &
\end{matrix}$$
In fact, any proper subgroup $H \subset \text{GL}(2, \mathbb{F}_p)$ is
contained in one of the following four types of subgroups (see
\cite{Mazur1977}):

(i) $H=\text{a Borel subgroup}$. Then $X_H=X_0(p)$.

(ii) $H=\text{the normalizer of a split Cartan subgroup}$. In this case,
denote $X_H=X_{\text{split}}(p)$. There is a natural isomorphism between
$X_{\text{split}}(p)$ and $X_0(p^2)/w_{p^2}$ as projective curves over
$\mathbb{Q}$, where $w_{p^2}$ is the canonical involution induced from
$z \mapsto -\frac{1}{p^2 z}$ on the upper half plane.

(iii) $H=\text{the normalizer of a nonsplit cartan subgroup}$. In this
case write $X_H=X_{\text{nonsplit}}(p)$.

(iv) $H=\text{an exceptional subgroup}$, $H$ is the inverse image in
$\text{GL}(2, \mathbb{F}_p)$ of an exceptional subgroup of
$\text{PGL}(2, \mathbb{F}_p)$. It is known that an exceptional subgroup
of $\text{PGL}(2, \mathbb{F}_p)$ is a subgroup isomorphic to the symmetric
group $\mathfrak{S}_4$, or alternating groups $\mathfrak{A}_4$ or
$\mathfrak{A}_5$.

  The further requirement $\det H=\mathbb{F}_p^{*}$ insures that the image
of $H$ in $\text{PGL}(2, \mathbb{F}_p)$ be isomorphic to $\mathfrak{S}_4$.
Moreover, if such an $H$ (with surjective determinant) exists when $K=\mathbb{Q}$,
then $p \equiv \pm 3$ (mod $8$). For such $p$ write $X_H=X_{\mathfrak{S}_4}(p)$.

  Let $G$ be a subgroup of $\text{GL}(V)$ of order prime to $p$, and let
$H$ be its image in $\text{PGL}(V)$. According to \cite{Se1972}, Prop. 16,
applies to the field $k=\mathbb{F}_p$, shows that we have the following
possibilities:

(i) $H$ is cyclic, hence contained in a Cartan subgroup of $\text{PGL}(V)$,
unique if $H \neq \{ 1 \}$; we conclude that $G$ is contained in a Cartan
subgroup of $\text{GL}(V)$.

(ii) $H$ is dihedral, so contains a non-trivial cyclic subgroup $C^{\prime}$
of index $2$. The group $C^{\prime}$ is contained in a unique Cartan subgroup
$C$ of $\text{PGL}(V)$; as $H$ normalises $C^{\prime}$, it also normalises $C$.
Returning to $\text{GL}(V)$, we conclude that $G$ is contained in the normaliser
of a Cartan subgroup of $\text{GL}(V)$.

(iii) $H$ is isomorphic to $\mathfrak{A}_4$, $\mathfrak{S}_4$ or $\mathfrak{A}_5$,
the latter only being possible if $p \equiv \pm 1$ (mod $5$). The elements of $H$
are then of order $1$, $2$, $3$, $4$ or $5$.

\begin{center}
{\bf \S 3.6. The resolvent of the seventh degree}
\end{center}

  To obtain the fifth-degree resolvent, which exists in the case of
the fifth degree transformation, one would have, following the general
approaches of Galois theory, first form the Galois field of the
transformation equation and choose for a single one of the groups
$G_{12}$ a suitable function of that field, which would then be the
root of the sought equation of the fifth degree. These are extended
developments, which are exhaustively presented by Klein in \cite{K}.
The corresponding investigation for the case of the seventh degree
transformation is given by Klein in the treatise \cite{K2}. Now,
however, at least in the two mentioned lowest cases $n=5$ and $n=7$,
one can arrive at the real knowledge of the resolvents with circumvention
of Galois fields directly from the discontinuity regions of the groups
$\Gamma_5$ and $\Gamma_7$, as was shown by Klein in a special treatise
\cite{Klein1879}. Here one has to make use of the purely
algebraic method, which Klein also already used for the transformation
equations themselves. With regard to the relationship of the resolvents
of the fifth and seventh degree to the Galois fields and the transformation
equations, we refer to the representations mentioned above, in the
following we will only give the actual list of resolvents according to
Klein's algebraic method (see \cite{Fricke}, Vol. II, pp. 483-486.)

  Klein's method is sufficient in the case $n=7$ to obtain the two
resolvents of the $7$-th degree. The discontinuity region of one of the
groups $\Gamma_7$ shown in \cite{Fricke}, Vol. II, Fig. 31 in p. 477,
is mapped by $J(\omega)$ onto a seven-sheeted Riemann surface above the
$J$-plane, which again is branched only at $J=0$, $1$ and $\infty$; and
there are now two three-sheeted branch points at $J=0$, at $J=1$ two
two-sheeted and at $J=\infty$ one seven-sheeted one. {\it This surface
also belongs to the genus $0$, so that we can again use a single-valued
function $\zeta(\omega)$.} Let
$$\zeta\left(\frac{-3+i \sqrt{3}}{2}\right)=0, \quad
  \zeta(i \infty)=\infty$$
hold so that $\zeta$ vanishes at $J=0$ in the sheet running isolated there
and becomes infinite at $J=\infty$ in the seven-sheeted branch point. By
these two conditions $\zeta$ is then only determined up to a constant
factor; a further determination remains reserved.

  In $\zeta$, $J$ is an entire function of the $7$-th degree, and in fact
from the distribution of values of $\zeta$ at $J=0$ and $J=1$ one derives
the approaches
$$J=a \zeta (\zeta^2+b \zeta+7c)^3, \quad
  J-1=a (\zeta^3+d \zeta^2+e \zeta+f)(\zeta^2+g \zeta+h)^2,\eqno{(3.6.1)}$$
where, to shorten the calculation, the factor $7$ is included in the
absolute member of the first bracket. By equating the two expressions
to be obtained for $\frac{dJ}{d \zeta}$, the omission of superfluous
factors results in the identical equation
$$\aligned
  &\left(\zeta^2+\frac{4}{7}b \zeta+c\right)(\zeta^2+b \zeta+7c)^2
   =(\zeta^2+g\zeta+h)\\
  &\times \left(\zeta^4+\frac{6d+5g}{7} \zeta^3+\frac{4dg+5e+3h}{7} \zeta^2
  +\frac{2dh+3eg+4f}{7} \zeta+\frac{eh+2fg}{7}\right).
\endaligned$$
Now it follows from (3.6.1) that the functions $(\zeta^2+b \zeta+7c)$ and
$(\zeta^2+g \zeta+h)$ have no linear factor in common. Thus the functions
$(\zeta^2+\frac{4}{7}b \zeta+c)$ and $(\zeta^2+g \zeta+h)$ are identical,
and by comparing coefficients one finds in total
$$7g=4b, \quad h=c, \quad 6d+5g=14b, \quad 4dg+5e+3h=7b^2+98c,$$
$$2dh+3eg+4f=98bc, \quad eh+2fg=343c^2.$$

  If now $b=0$, then $g=0$, $d=0$ and because of the penultimate equation
$f=0$ would follow. But it is $f \neq 0$, because for $\zeta=0$ also $J=0$
and not $J=1$. According to this, $b \neq 0$. But then we can fix the factor
still available for the unique determination of $\zeta$ in such a way that
$b=7$. The equations to be solved now shorten to
$$b=7, \quad g=4, \quad h=c, \quad d=13, \quad e=27+19c, \quad 3e+f=165c,$$
$$ce+8f=343c^2.$$
By eliminating $e$ and $f$ from the last three equations, the quadratic equation
$$4 c^2-11 c+8=0$$
follows for $c$, in the solution of which, as it must be, the irrationality
$i \sqrt{7}$ is found. We obtain, corresponding to our two resolvents, the
two systems of coefficients
$$c=h=\frac{11 \pm i \sqrt{7}}{8}, \quad
  e=\frac{425 \pm 19 i \sqrt{7}}{8}, \quad
  f=\frac{135 \pm 27 i \sqrt{7}}{2},$$
while the coefficients $b$, $g$, $d$ have the rational values already
indicated. Note that the Riemann surface leading to the second resolvent
of degree $7$ above the $J$-plane, which is symmetric to the first surface,
has branch points of the same number of sheets at $J=0$, $1$ and $\infty$.
The reasoning carried out had to lead to the second resolvent at the same
time. Finally, the value $a$ is obtained by entering $\zeta=0$ and $J=0$
into the second equation of (3.6.1); one finds:
$$-1=a \cdot f \cdot h^2=a \frac{135 \pm 27 i \sqrt{7}}{2} \left(\frac{11
  \pm i \sqrt{7}}{8}\right)^2, \quad
  -a^{-1}=\frac{351 \pm 189 i \sqrt{7}}{4}.$$

  If one inserts the calculated values of the coefficients into the approaches
(8.6.1), then these can be combined into the following proportion:
$$\aligned
  &J: (J-1): 1=\zeta \left(\zeta^2+7 \zeta+\frac{77 \pm 7 i \sqrt{7}}{8}\right)^3\\
  &:\left(\zeta^3+13 \zeta^2+\frac{425 \pm 19 i \sqrt{7}}{8} \zeta+
   \frac{135 \pm 27 i \sqrt{7}}{2}\right)\left(\zeta^2+4\zeta+\frac{11
   \pm i \sqrt{7}}{8}\right)^2\\
 &:-\frac{351 \pm 189 i \sqrt{7}}{4}.
\endaligned\eqno{(3.6.2)}$$
{\it As equations of the $7$-th degree for $\zeta$, we have here the two
resolvents of the $7$-th degree in a first form.}

  From (8.6.2) one can also easily derive equations for the simplest modular forms
of $\Gamma_7$. First calculate from the first and third member of the proportion
of (8.6.2):
$$(12 g_2)^3=\left(\frac{1 \pm i \sqrt{7}}{2}\right)^7 \Delta \zeta
  \left(\zeta^2+7 \zeta+\frac{77 \pm 7 i \sqrt{7}}{8}\right)^3.$$
From this it can be seen, {\it that the cube roots
$$f(\omega_1, \omega_2)=\root 3 \of{\frac{1 \pm i \sqrt{7}}{2} \cdot \Delta \zeta}
  \eqno{(3.6.3)}$$
are also unique and indeed entire modular forms of the $(-4)$-th dimension
of $\Gamma_7$.} If one introduces $f$ instead of $\zeta$ as an unknown into
the penultimate equation, then, after subtracting the cube root, we obtain
the {\it new form of our resolvent of the $7$-th degree:}
$$f^7+7 \frac{1 \pm i \sqrt{7}}{2} \Delta f^4-7 \frac{5 \mp i \sqrt{7}}{2} \Delta^2 f
  -12 g_2 \Delta^2=0. \eqno{(3.6.4)}$$

\begin{center}
{\bf \S 3.7. The two resolvents of the eleventh degree}
\end{center}

  For obtaining the fifth-degree resolvent it was sufficient to state, that
the five-sheeted Riemann surface above the $J$-plane has a three-sheeted
branch point at $J=0$, at $J=1$ two two-sheeted branch points and at
$J=\infty$ one five-sheeted branch point, and that further branch points
did not occur. Which sheets are connected with each other at $J=0$ and
$J=1$, however, did not need to be determined in more detail. Accordingly,
there is only one single five-sheeted surface, which has the branch in
question. Correspondingly, in the case of $n=7$ there are only two
(mutually symmetrical) surfaces, which have the branch described in more
detail above for this case. In contrast, Klein found when extending his
investigations to the case $n=11$ (see \cite{K3}), that there are not only
two, but ten eleven-sheeted Riemann surfaces above the $J$-plane, which
the region of \cite{Fricke}, Vol. II, Fig. 32 in p. 480, have the numbers
of branch points at $J=0$, $1$ and $\infty$. The purely algebraic method
of the previous paragraph can be applied with a single-valued function
$\zeta$ in exactly the same way as above, since the genus $0$ is again
present, but would have to lead to ten different systems of coefficients
here, and it would therefore still have to be investigated in more detail,
which of the ten equations are the two sought-after resolvents we are
looking for. It therefore seems more appropriate here, to link up with the
Galois field in accordance with the consideration given at the beginning of
the previous paragraph, and from there to produce a suitable puncture for
the individual icosahedral group $G_{60}$ in $G_{660}$.

  For the construction of Galois field Klein uses the system of five
modular forms $x_{\lambda}$ belonging to $n=11$, which are given by
\cite{K3}, \S 11. Leaving out the common factor $i$ and using the five
quadratic residues of $11$ as indices, one has the representations
$$\left\{\aligned
  x_1 &=\sqrt{\frac{2 \pi}{\omega_2}} \root 8 \of{\Delta^5} q^{\frac{81}{44}}
        (1-q^2-q^{20}+\cdots),\\
  x_3 &=\sqrt{\frac{2 \pi}{\omega_2}} \root 8 \of{\Delta^5} q^{\frac{25}{44}}
        (1-q^6-q^{16}+\cdots),\\
  x_9 &=\sqrt{\frac{2 \pi}{\omega_2}} \root 8 \of{\Delta^5} q^{\frac{49}{44}}
        (1-q^4-q^{18}+\cdots),\\
  x_5 &=\sqrt{\frac{2 \pi}{\omega_2}} \root 8 \of{\Delta^5} q^{\frac{1}{44}}
        (1-q^{10}-q^{12}+\cdots),\\
  x_4 &=-\sqrt{\frac{2 \pi}{\omega_2}} \root 8 \of{\Delta^5} q^{\frac{9}{44}}
        (1-q^8-q^{14}+\cdots),
\endaligned\right.\eqno{(3.7.1)}$$
for these quantities belonging to the dimension $(-8)$. In contrast to the
substitutions $S=\left(\begin{matrix} 1 & 1\\ 0 & 1 \end{matrix}\right)$ and
$T=\left(\begin{matrix} 0 & 1\\ -1 & 0 \end{matrix}\right)$, the
$x_{\lambda}$ according to \cite{Fricke}, Vol. II, p. 315 undergo
the linear transformations
$$\left\{\aligned
  (S) \quad \qquad &x_{\lambda}^{\prime}=\varepsilon^{-\frac{\lambda(11-\lambda)}{2}}
   x_{\lambda},\\
  (T) \quad i \sqrt{11} &x_{\lambda}^{\prime}=\sum_{\kappa} (\varepsilon^{\kappa
   \lambda}-\varepsilon^{-\kappa \lambda}) x_{\kappa},
\endaligned\right.\eqno{(3.7.2)}$$
where the sum related to $\kappa$ here and furthermore relates to the five
quadratic residues of $11$ and is $\varepsilon=e^{\frac{2 \pi i}{11}}$. The
$G_{660}$ of the transformations of the Galois field in itself presents
itself here as a group of linear $x_{\lambda}$-substitutions, which can be
generated from the two substitutions (3.7.2).

  We now consider the icosahedron group $G_{60}$ belonging to the discontinuity
region of \cite{Fricke}, Vol. II, Fig. 32 in p. 480, and have to form associated
functions of $x_{\lambda}$. The $G_{60}$ contains the substitution of the period
five $S^{\prime} \equiv \left(\begin{matrix} 3 & 0\\ 0 & 4 \end{matrix}\right)$
(mod $11$), to be denoted by $S^{\prime}$, in the exercise of which the
$x_{\lambda}$ permute themselves according to (18) and (19) in \cite{Fricke},
Vol II, p. 313:
$$(S^{\prime}) \quad \quad x_{\lambda}^{\prime}=x_{3 \lambda}. \eqno{(3.7.3)}$$
The sum of the $x_{\lambda}$, to be abbreviated as $\sigma_{\infty}$:
$$\sigma_{\infty}=x_1+x_3+x_9+x_5+x_4 \eqno{(3.7.4)}$$
thus represents a modular form which is invariant with respect to $S^{\prime}$.
However, as can be shown with the help of relation
$$\varepsilon+\varepsilon^3+\varepsilon^9+\varepsilon^5+\varepsilon^4
 -\varepsilon^{10}-\varepsilon^8-\varepsilon^2-\varepsilon^6-\varepsilon^7=i \sqrt{11},$$
the sum is also invariant with respect to $T$, i.e. belongs to the metacyclic
subgroup $G_{10}$ of $G_{60}$ to be generated from $S^{\prime}$ and $T$. Through
the entire substitutions of $G_{60}$, the modular form (3.7.4) is thus transformed
into six conjugate forms, which may be denoted by $\sigma_{\infty}$, $\sigma_0$,
$\sigma_1$, $\ldots$, $\sigma_4$. In particular, $\sigma_0$ emerges from
$\sigma_{\infty}$ through the substitution of the period two
$S^{-1} \cdot T \cdot S=\left(\begin{matrix} 1 & 2\\ -1 & -1 \end{matrix}\right)$
belonging to $G_{60}$. The effect of this substitution on the $x_{\lambda}$ is:
$$i \sqrt{11} x_{\lambda}^{\prime}=\sum_{\kappa} (\varepsilon^{\kappa \lambda}
  -\varepsilon^{-\kappa \lambda}) \varepsilon^{6 (\kappa^2-\lambda^2)} x_{\kappa}.$$
From this one calculates for $\sigma_0$ the equation:
$$i \sqrt{11} \sigma_0=\sum_{\kappa} (\varepsilon^{\kappa^2}+2 \varepsilon^{2 \kappa^2}
  -2 \varepsilon^{5 \kappa^2}-\varepsilon^{6 \kappa^2}) x_{\kappa}.$$
The remaining four forms $\sigma_1$, $\sigma_2$, $\ldots$ follow from $\sigma_0$
by repeated exercise of the substitution $S^{\prime}$; one finds:
$$i \sqrt{11} \sigma_{\nu}=\sum_{\kappa} (\varepsilon^{\kappa^2}+2 \varepsilon^{2 \kappa^2}
  -2 \varepsilon^{5 \kappa^2}-\varepsilon^{6 \kappa^2}) x_{\kappa \cdot 3^{\nu}}.
  \quad \nu=0, 1, \cdots, 4.\eqno{(3.7.5)}$$

  As functions of $G_{60}$ we now use symmetric expressions of the six $\sigma$.
Since their sum vanishes identically, the two simplest expressions are the
second and the third power sum, for which one finds the representations in
the $x_{\lambda}$:
$$\left\{\aligned
  \frac{-1+i \sqrt{11}}{12} \sum \sigma^2 &=\sum_{\kappa} \left(x_{\kappa}^2
  -x_{\kappa} x_{4 \kappa}-\frac{1-i \sqrt{11}}{2} x_{\kappa} x_{5 \kappa}\right),\\
  -\frac{i \sqrt{11}}{6} \sum \sigma^3 &=\sum_{\kappa} \left(x_{\kappa}^3-3 i \sqrt{11}
  x_{\kappa}^2 x_{3 \kappa}+3 x_{\kappa} x_{5 \kappa} (x_{\kappa}-x_{4 \kappa})\right.\\
  &\quad \left.+\frac{1+i \sqrt{11}}{2} x_{\kappa} (x_{\kappa} x_{4 \kappa}-x_{4 \kappa}
  x_{9 \kappa}-2 x_{9 \kappa} x_{\kappa})\right).
\endaligned\right.\eqno{(3.7.6)}$$

  Since the $x_{\lambda}$ belong to the dimension $(-8)$, we have in the
expressions on the right entire modular forms of the dimensions $-16$ and
$-24$ of the $\Gamma_{11}$ belonging to the $G_{60}$. According to (3.7.1)
the quotient of the first of these forms and the discriminant $\Delta$ is
also an entire modular form of $\Gamma_{11}$, namely of the dimension $-4$.
In this form
$$f(\omega_1, \omega_2)=\Delta^{-1} \cdot \sum_{\kappa} \left(x_{\kappa}^2-
  x_{\kappa} x_{4 \kappa}-\frac{1-i \sqrt{11}}{2} x_{\kappa} x_{5 \kappa}\right)
  \eqno{(3.7.7)}$$
to be denoted by $f(\omega_1, \omega_2)$, whose series representation is
calculated from (3.7.1) to:
$$\aligned
 &f(\omega_1, \omega_2)\\
=&\left(\frac{2 \pi}{\omega_2}\right)^4 q^{\frac{6}{11}}
  \left(1+q^{\frac{2}{11}}+q^{\frac{4}{11}}-\frac{1-i \sqrt{11}}{2} q^{\frac{6}{11}}
  +\frac{1-i \sqrt{11}}{2} q^{\frac{8}{11}}+\cdots\right),
\endaligned\eqno{(3.7.8)}$$
we have obtained a particularly simple modular form of $\Gamma_{11}$.

  The form on the right-hand side of the second equation of (3.7.6) gives, divided
by $\Delta^2$, a modular function of the $\Gamma_{11}$, which becomes infinite
only in the cusp $\omega=i \infty$ of the discontinuity region. Since from (3.7.1)
we calculate $q^{-\frac{2}{11}}$ as the initial element of the power series for
this quotient, a first-order pole lies in that cusp, so that we have obtained
a {\it single-valued function of the} $\Gamma_{11}$ in the quotient. To simplify
the following formulae, we introduce the quotient reduced by $3i \sqrt{11}$ as
a single-valued function $\zeta(\omega)$ of the $\Gamma_{11}$:
$$\zeta(\omega)=-3i \sqrt{11}+\Delta^{-2} \sum_{\kappa} (x_{\kappa}^3-3 i
  \sqrt{11} x_{\kappa}^2 x_{3 \kappa}+3 x_{\kappa} x_{5 \kappa}
  (x_{\kappa}-x_{4 \kappa})+\cdots).\eqno{(3.7.9)}$$
The power series of this function is then
$$\zeta(\omega)=q^{-\frac{2}{11}} \left(1+*+\frac{1+i \sqrt{11}}{2} q^{\frac{4}{11}}
  +2q^{\frac{6}{11}}+\frac{1+i \sqrt{11}}{2} q^{\frac{8}{11}}+\cdots\right),
  \eqno{(3.7.10)}$$
where the asterisk emphasises that the member with $q^{\frac{2}{11}}$ drops out.

  The form $f$ as of dimension $(-4)$ has zero points of total order $\frac{11}{3}$
in the discontinuity region. Since according to (3.7.8) there is a zero point of order
$3$ in the cusp $i \infty$, then the order $\frac{2}{3}$ still remains. Thus, there
will either be a zero point of order $\frac{1}{3}$ in each of the two range corners
at $\omega=\frac{-7+i \sqrt{3}}{2}$ and $\frac{9+i \sqrt{3}}{2}$ or a zero point of
order $\frac{2}{3}$ in one of these places. In any case, one recognises in the quotient
of $f^3$ and $\Delta$ an entire function of the second degree of $\zeta$, for which
the series expansions easily produce the representation
$$\frac{f^3}{\Delta}=\zeta^2+3 \zeta+(5-i \sqrt{11}) \eqno{(3.7.11)}$$
without difficulty. Since there is no square on the right here, zero points of order
$\frac{1}{3}$ of $f$ are present in each of the two corners mentioned.

  We further form the quotient of $g_2$ and $f$, which also represents a function
of $\Gamma_{11}$. The two zero points just mentioned have continued in this quotient.
There remain, coming from the numerator $g_2$, three first-order zero points in the
three vertex cycles of the region composed of the vertices $\frac{\pm 1+i \sqrt{3}}{2}$,
$\frac{\pm 3 +i \sqrt{3}}{2}$, $\ldots$ and a third-order pole in the cusp
$\omega=i \infty$, coming from the denominator $f$. The quotient is therefore an entire
function of the third degree of $\zeta$, for which the power series give the
representation:
$$\frac{g_2}{f}=\frac{1}{12} \left(\zeta^3-\zeta^2-3 \frac{1+i \sqrt{11}}{2} \zeta
  -\frac{7-i \sqrt{11}}{2}\right).\eqno{(3.7.12)}$$

  The elimination of $\zeta$ from (3.7.11) and (3.7.12) leads to the result:
{\it One of the two sought resolvents of the $11$-th degree has the form
$$\aligned
  f^{11}-22 \Delta f^8&+11 (9+2i \sqrt{11}) \Delta^2 f^5-132 g_2 \Delta^2 f^4
  +88 i \sqrt{11} \Delta^3 f^2\\
  &+66 (3-i \sqrt{11}) g_2 \Delta^3 f-144 g_2^2 \Delta^3=0
\endaligned\eqno{(3.7.13)}$$
as equation for the modular form $f$ of $\Gamma_{11}$; the second resolvent
results from this by sign change from $\sqrt{11}$.}

  We obtain a second form of the resolvent by multiplying equation (3.7.11)
by the third power of equation (3.7.12):
$$12^3 J=(\zeta^2+3 \zeta+(5-i \sqrt{11})) \left(\zeta^3-\zeta^2
  -3 \frac{1+i \sqrt{11}}{2} \zeta-\frac{7-i \sqrt{11}}{2}\right)^3.\eqno{(3.7.14)}$$
This result can be subjected to a test. If we map the discontinuity region of
$\Gamma_{11}$ onto the $J$-plane, we obtain a Riemann surface, which has, besides
other branch points (at $J=0$ and $J=\infty$), four two-sheeted branch points at
$J=1$. Accordingly, $(J-1)$ must be equal to the product of an entire function of
the third degree and the square of an entire function of the fourth degree of $\zeta$,
so that there exists an equation identical in $\zeta$ of the following form:
$$\aligned
  &(\zeta^2+3 \zeta+(5-i \sqrt{11})) \left(\zeta^3-\zeta^2-3 \frac{1+i \sqrt{11}}{2}
   \zeta-\frac{7-i \sqrt{11}}{2}\right)^3-12^3\\
 =&(\zeta^3+a \zeta^2+b \zeta+c)(\zeta^4+d \zeta^3+e \zeta^2+f \zeta+g)^2.
\endaligned$$
The calculation confirms this indeed. By differentiating the above equation by
$\zeta$ and proceeding as in the previous paragraph, one easily finds the
coefficients $a$, $b$, $c$, $\ldots$. {\it As the second form of the two
resolvents of the eleventh degree, we find the following two equations of the
eleventh degree for $\zeta$:
$$\aligned
 &J:(J-1):1\\
=&(\zeta^2+3 \zeta+(5 \mp i \sqrt{11}))\left(\zeta^3-\zeta^2-3
  \frac{1 \pm i \sqrt{11}}{2} \zeta-\frac{7 \mp i \sqrt{11}}{2}\right)^3\\
 &:\left(\zeta^3-4 \zeta^2+\frac{7 \mp 5 i \sqrt{11}}{2} \zeta-(4 \mp 6 i \sqrt{11})\right)
  \left(\zeta^4+2 \zeta^3+3 \frac{1 \mp i \sqrt{11}}{2} \zeta^2\right.\\
 &\left.-(5 \pm i \sqrt{11}) \zeta-3 \frac{5 \pm i \sqrt{11}}{2}\right)^2: 12^3.
\endaligned\eqno{(3.7.15)}$$}

\begin{center}
{\large\bf 4. Galois representations arising from $\mathcal{L}(X(p))$ and
              their modularity}
\end{center}

  In this section, we will prove the following:

\textbf{Theorem 4.1. (Main Theorem 2).} (Surjective realization and modularity).
{\it For a given representation arising from the defining ideals}
$$\rho_p: \text{Gal}(\overline{\mathbb{Q}}/\mathbb{Q}) \rightarrow
  \text{Aut}(\mathcal{L}(X(p))), \quad (p \geq 7),$$
{\it there are surjective realizations by the equations of degree $p+1$ or
$p$ with Galois group isomorphic with} $\text{Aut}(\mathcal{L}(X(p)))$. {\it
These equations are defined over $\mathcal{L}(X(p))$, their coefficients are
invariant under the action of} $\text{Aut}(\mathcal{L}(X(p)))$. {\it Moreover,
it corresponds to the $p$-th order transformation equation of the $j$-function
with Galois group isomorphic to} $\text{PSL}(2, \mathbb{F}_p)$. {\it In particular,}

(1) {\it when $p=7$, there are two such equations of degree eight and seven, respectively.
One comes from the Jacobian multiplier equation of degree eight, the other corresponds
to the decomposition of} $\text{PSL}(2, \mathbb{F}_7)=S_4 \cdot C_7$, {\it where $C_7$
is a cyclic subgroup of order seven.}

(2) {\it When $p=11$, there are two such equations of degree twelve and eleven, respectively.
One comes from the Jacobian multiplier equation of degree twelve, the other corresponds
to the decomposition of} $\text{PSL}(2, \mathbb{F}_{11})=A_5 \cdot C_{11}$, {\it where
$C_{11}$ is cyclic subgroup of order eleven.}

(3) {\it When $p=13$, there are two distinct equations of the same degree fourteen,
respectively. One comes from the Jacobian multiplier equation of degree fourteen,
the other is not the Jacobian multiplier equation.}

{\it Proof}. Given Galois representation arising from the defining ideals
of modular curves $X(p)$:
$$\rho_p: \text{Gal}(\overline{\mathbb{Q}}/\mathbb{Q}) \rightarrow
  \text{Aut}(\mathcal{L}(X(p))), \quad (p \geq 7)$$
in order to give its surjective realization, we will give the equations of
degree $p+1$ or $p$ satisfying:

(1) they are defined over $\mathcal{L}(X(p))$;

(2) their coefficients are invariant under the action of $\text{Aut}(\mathcal{L}(X(p)))$;

(3) their Galois group is isomorphic to $\text{Aut}(\mathcal{L}(X(p)))$.

  In order to prove the modularity of $\rho_p$, we will show that it corresponds
to the $p$-th order transformation equation of the $j$-function with Galois
group isomorphic to $\text{PSL}(2, \mathbb{F}_p)$.

  In fact, by Theorem 2.4.1, the modular curve $X(p)$ is isomorphic to the
locus $\mathcal{L}(X(p))$. This proves (1).  For (2) and (3), we need Jacobian
multiplier equations of degree $p+1$.

  Let us recall a result which Jacobi had established as early as $1828$ in
his ``Notices sur les fonctions elliptiques'' (see \cite{Ja1828} and \cite{K}).
Jacobi there considered, instead of the modular equation, the so-called
multiplier-equation, together with other equations equivalent to it, and
found that their $(n+1)$ roots are composed in a simple manner of $\frac{n+1}{2}$
elements, with the help of merely numerical irrationalities. Namely, if we
denote these elements by $\mathbf{A}_0$, $\mathbf{A}_1$, $\ldots$,
$\mathbf{A}_{\frac{n-1}{2}}$, and further, for the roots $z$ of the equation
under consideration, apply the indices employed by Galois, we have, with
appropriate determination of the square root occurring on the left-hand side:
$$\left\{\aligned
  \sqrt{z_{\infty}} &=\sqrt{(-1)^{\frac{n-1}{2}} \cdot n} \cdot \mathbf{A}_0,\\
  \sqrt{z_{\nu}} &=\mathbf{A}_0+\epsilon^{\nu} \mathbf{A}_1+\epsilon^{4 \nu} \mathbf{A}_2
  +\cdots+\epsilon^{(\frac{n-1}{2})^2 \nu} \mathbf{A}_{\frac{n-1}{2}}
\endaligned\right.\eqno{(4.1)}$$
for $\nu=0, 1, \cdots, n-1$ and $\epsilon=e^{\frac{2 \pi i}{n}}$.
Jacobi had himself emphasized the special significance of his result by
adding to his short communication: ``C'est un th\'{e}or\`{e}me des plus
importants dans la th\'{e}orie alg\'{e}brique de la transformation et de
la division des fonctions elliptiques.''

  In fact, there is a maximal subgroup $H$ of order $\frac{1}{2}p(p-1)$
of index $p+1$ in the group $\text{PSL}(2, \mathbb{F}_p)$ which acts on
the projective space
$$\mathbb{CP}^{\frac{p-3}{2}}=\left\{ (z_1, z_2, \ldots, z_{\frac{p-1}{2}}):
  z_i \in \mathbb{C} \quad \left(1 \leq i \leq \frac{p-1}{2}\right) \right\}.$$
The two generators of $\text{PSL}(2, \mathbb{F}_p)$ are given by $S$ and $T$
with
$$S^2=T^p=(ST)^3=1.\eqno{(4.2)}$$
For this maximal subgroup $H$, we can construct a polynomial $\phi_{\infty}$
in $\frac{p-1}{2}$ variables $z_1$, $\ldots$, $z_{\frac{p-1}{2}}$, which is
invariant under the action of $H$. By the coset-decomposition of
$\text{PSL}(2, \mathbb{F}_p)$ corresponding to $H$, there are $p+1$ such cosets:
$$\text{PSL}(2, \mathbb{F}_p)=H \cup H \cdot S \cup H \cdot ST \cup \cdots \cup
                              H \cdot ST^{p-1}.\eqno{(4.3)}$$
Correspondingly, we can construct $p+1$ polynomials $\phi_{\infty}$, $\phi_0$,
$\ldots$, $\phi_{p-1}$ by
$$\phi_{\nu}=\phi(ST^{\nu}(z_1, z_2, \ldots, z_{\frac{p-1}{2}})),\eqno{(4.4)}$$
where $\nu=0, \ldots, p-1$. The Jacobian multiplier equation od degree
$p+1$ is just the equation whose $p+1$ roots are $\phi_{\infty}$, $\phi_0$,
$\ldots$, $\phi_{p-1}$. The Galois group is the permutation group of
these $p+1$ roots. By the construction as above, $\text{PSL}(2, \mathbb{F}_p)$
acts on $H$, $H \cdot S$, $\ldots$, $H \cdot ST^{p-1}$ by the coset
decomposition (4.3) as a permutation of these $p+1$ cosets, which corresponds
to the permutation of these $p+1$ roots $\phi_{\infty}$, $\phi_0$, $\ldots$,
$\phi_{p-1}$. This implies that the Galois group is isomorphic with
$\text{PSL}(2, \mathbb{F}_p)$. Moreover, the coefficients of this Jacobian
multiplier equation od degree $p+1$ are the symmetric polynomials of $\phi_{\infty}$,
$\phi_0$, $\ldots$, $\phi_{p-1}$. By the construction, these symmetric polynomials
are $\text{PSL}(2, \mathbb{F}_p)$-invariant polynomials.

  In fact, in \cite{K4}, a close relationship of the aforementioned
quantities to the roots of the new [first-level] multiplier equations
was found. If one calls these roots $z_{\infty}$, $z_0$, $\ldots$,
$z_{n-1}$ in a known manner and sets, since we are dealing with a
Jacobian equation
$$\sqrt{z_{\infty}}=\sqrt{(-1)^{\frac{n-1}{2}} n} \cdot A_0,$$
$$\sqrt{z_{\nu}}=A_0+\varepsilon^{36 \nu} A_1+\varepsilon^{4 \cdot 36 \nu} A_2
  +\cdots+\varepsilon^{\left(\frac{n-1}{2}\right)^2 36 \nu} A_{\frac{n-1}{2}},
  \quad \left(\varepsilon=e^{\frac{2 \pi i}{n}}\right),$$
then one finds, apart from a factor $\varrho$ which is not further considered
here:
$$\left\{\aligned
  \varrho A_0 &=\sum_{\lambda=-\infty}^{+\infty} (-1)^{\lambda} \cdot
                q^{\frac{(6 \lambda+1)^2 n}{12}},\\
  \varrho A_{\alpha} &=(-1)^{\alpha} \cdot \sum_{\lambda=-\infty}^{+\infty}
                (-1)^{\lambda} \left\{q^{\frac{((6 \lambda+1) n+6 \alpha)^2}{12 n}}
                +q^{\frac{((6 \lambda+1)n-6 \alpha)^2}{12 n}}\right\}.
\endaligned\right.\eqno{(4.5)}$$
Keeping to the designation method introduced here, Klein's result for $n=7$ is
as follows:
$$\frac{\mu}{\lambda}=-q^{-\frac{4}{7}}+\cdots=\frac{A_1}{A_0}, \quad
  \frac{\nu}{\mu}=+q^{-\frac{2}{7}}+\cdots=\frac{A_2}{A_0}, \quad
  \frac{\lambda}{\nu}=-q^{+\frac{6}{7}}+\cdots=\frac{A_3}{A_0},$$
and for $n=11$:
$$\frac{y_5}{y_4}=+q^{-\frac{8}{11}}+\cdots=-\frac{A_2}{A_0}, \quad
  \frac{y_9}{y_5}=-q^{-\frac{10}{11}}+\cdots=-\frac{A_4}{A_0},$$
$$\frac{y_3}{y_9}=-q^{+\frac{4}{11}}+\cdots=-\frac{A_3}{A_0}, \quad
  \frac{y_1}{y_3}=+q^{_\frac{6}{11}}+\cdots=-\frac{A_5}{A_0},$$
$$\frac{y_4}{y_1}=+q^{+\frac{20}{11}}+\cdots=-\frac{A_1}{A_0}.$$

  As is well known, one has:
$$\frac{\vartheta(x, q^n) \cdot \vartheta_2(x, q^n) \cdot \vartheta_3(x, q^n)}
  {2 q^{\frac{n}{6}} \prod_{\nu=1}^{\infty} (1-q^{2 n \nu})^2}
 =\sum_{\lambda=-\infty}^{+\infty} (-1)^{\lambda}
  q^{\frac{(6 \lambda+1)^2 n}{12}} \cdot \cos(6 \lambda+1) x.$$
Therefore, for $q=e^{\pi i \omega}$ it follows:
$$\frac{A_{\alpha}}{A_0}=(-1)^{\alpha} \cdot 2 \cdot q^{\frac{3 \alpha^2}{n}}
  \cdot \frac{\vartheta(\alpha \omega \pi, q^n) \cdot
  \vartheta_2(\alpha \omega \pi, q^n) \cdot \vartheta_3(\alpha \omega \pi, q^n)}
  {\vartheta(0, q^n) \cdot \vartheta_2(0, q^n) \cdot \vartheta_3(0, q^n)}.$$
But now in general:
$$2 \cdot \frac{\vartheta(x, q) \cdot \vartheta_2(x, q) \cdot \vartheta_3(x, q)}
  {\vartheta(0, q) \cdot \vartheta_2(0, q) \cdot \vartheta_3(0, q)}
 =\frac{\vartheta_1(2x, q)}{\vartheta_1(x, q)}.$$
Thus we get:
$$\frac{A_{\alpha}}{A_0}=(-1)^{\alpha} \cdot q^{\frac{3 \alpha^2}{n}} \cdot
  \frac{\vartheta_1(2 \alpha \omega \pi, q^n)}{\vartheta_1(\alpha \omega \pi, q^n)}.$$

  For the sake of brevity Klein wrote
$$\sigma z_{\alpha}=(-1)^{\alpha} \cdot q^{\frac{\alpha^2}{n}} \cdot
  \vartheta_1(\alpha \omega \pi, q^n),$$
where $\sigma$ may mean a factor to be determined later and $\alpha$ is to run
from $1$ to $(n-1)$, where, by the way, only $\frac{n-1}{2}$ essentially
different quantities $z_{\alpha}$ arise by becoming (for odd $n$)
$$z_{n-\alpha}=-z_{\alpha}.$$
{\it Then one has:}
$$\frac{A_{\alpha}}{A_0}=\frac{z_{2 \alpha}}{z_{\alpha}}.$$
And from this it follows directly {\it that the ratios Klein used for $n=7$,
$11$ essentially correspond to the $z_{\alpha}$.} In fact, a comparison of
the formulae just given gives:

1. for $n=7$:
$$\lambda: \mu: \nu=z_1: z_2: z_4,$$

2. for  $n=11$:
$$y_1: y_4: y_5: y_9: y_3=z_1: z_9: z_4: z_3: z_5.$$

  By the results of section 3.4, especially Theorem 3.4.1 and Theorem 3.4.2,
for these Jacobian multiplier equation of degree $p+1$, it corresponds
to the $p$-th order transformation equation of the $j$-function with Galois
group isomorphic to $\text{PSL}(2, \mathbb{F}_p)$.

  The proof for the special cases of $p=7$, $11$ and $13$ will be
given in the following three sections.

\begin{center}
{\large\bf 4.1. The case of $p=7$ and Klein quartic curve}
\end{center}

  This process has a long history which dates back to F. Klein. In the
simplest case $p=7$, it originated from Klein's work on the Galois group
$\text{PSL}(2, \mathbb{F}_7)$.  As is well known that the modular equation
corresponding to the seventh-order transformation of elliptic functions
has a Galois group of $168$ substitutions, which is isomorphic to
$\text{PSL}(2, \mathbb{F}_7)$, and it is an old problem to transform equations
of the seventh and eighth degree, which have just this group, to the modular
equation in question by algebraic processes, and thus to solve them by elliptic
functions. Klein (see \cite{Klein1877/1878} and \cite{Klein1879(II)}) succeeded
to give the modular equation such a form, that this reduction is indeed possible.
In doing so, one needs an auxiliary equation of the fourth degree, which one
can show cannot be avoided. This is just the Klein quartic curve.

  Now for $p=7$, consider the Galois representation:
$$\rho_7: \text{Gal}(\overline{\mathbb{Q}}/\mathbb{Q}) \rightarrow
  \text{Aut}(\mathcal{L}(X(7))),$$
there are two such equations with degree $8$ and $7$, respectively. Correspondingly,
there are two distinct kinds of surjective realizations as well as two distinct
kinds of modularities. For the degree eight equation with Galois group isomorphic
to $\text{PSL}(2, \mathbb{F}_7)$, Klein constructed (see \cite{K2} and \cite{KF1},
Section three, Chapter seven) the following equation (4.1.16) of degree eight
corresponding to a maximal subgroup $G_{21}$ of $G_{168}$. In particular, the
equation (4.1.16) is defined over the Klein quartic curve $f=0$, which is just
the locus $\mathcal{L}(X(7))$ of $X(7)$. Its coefficients are invariant under
the action of the automorphism group of the Klein quartic curve. Klein proved
(see \cite{K2}) that the Galois group of (4.1.16) is isomorphic with the
automorphism group of the Klein quartic curve. In fact, (4.1.16) is just the
Jacobian multiplier equation of degree eight. Klein showed that (4.1.16) can be
reduced to the the modular equation of degree eight which corresponding to the
seventh-order transformation of $j$-functions given by (4.1.12). For its modern
form, see (3.3.5), (3.3.6) and (3.3.7).

  For the degree seven equation with Galois group isomorphic to
$\text{PSL}(2, \mathbb{F}_7)$, Klein constructed the following
equation (4.1.24) of degree seven corresponding to a maximal
subgroup $G_{24}$ of $G_{168}$. In particular, the equation (4.1.24)
is defined over the Klein quartic curve $f=0$. Its coefficients
are invariant under the action of the automorphism group of the
Klein quartic curve. Klein proved (see \cite{K2}) that the Galois
group of (4.1.24) is isomorphic with the automorphism group of the
Klein quartic curve. Klein showed that (4.1.24) can be reduced to
the seventh-order transformation of $j$-functions given by (4.1.22).

  In fact, in the above two cases, degree eight and seven, the reduction process
is always given by the theta constants of order seven, and the Klein quartic
curve $\mathcal{L}(X(7))$ is parameterized by the theta constants of order seven.
This is the reason that both (4.1.16) and (4.1.24) are defined over the locus
$\mathcal{L}(X(7))$. Therefore, we have proved the surjective realization as
well as the modularity of the Galois representation $\rho_7$.

  Now, we will give the detailed construction. In \cite{K2} or \cite{KF1},
Section three, Chapter seven, starting with a quartic invariant associated
with $\text{PSL}(2, \mathbb{F}_7)$:
$$f=\lambda^3 \mu+\mu^3 \nu+\nu^3 \lambda,\eqno{(4.1.1)}$$
Klein constructed a system of invariants coming from geometry. The
first covariant of $f$ is the Hessian $\nabla$ of order six:
$$\nabla=\frac{1}{54} \begin{vmatrix}
         \frac{\partial^2 f}{\partial \lambda^2} &
         \frac{\partial^2 f}{\partial \lambda \partial \mu} &
         \frac{\partial^2 f}{\partial \lambda \partial \nu}\\
         \frac{\partial^2 f}{\partial \mu \partial \lambda} &
         \frac{\partial^2 f}{\partial \mu^2} &
         \frac{\partial^2 f}{\partial \mu \partial \nu}\\
         \frac{\partial^2 f}{\partial \nu \partial \lambda} &
         \frac{\partial^2 f}{\partial \nu \partial \mu} &
         \frac{\partial^2 f}{\partial \nu^2}
         \end{vmatrix}
        =5 \lambda^2 \mu^2 \nu^2-(\lambda^5 \nu+\nu^5 \mu+\mu^5 \lambda).
\eqno{(4.1.2)}$$
As is well known, Hesse constructed for a general curve of order four a
curve of order $14$ that goes through the contact points of the bitangents.
In Klein's case, this property holds for any covariant of order $14$ that
is not a multiple of $f^2 \nabla$. Klein chose
$$C=\frac{1}{9} \begin{vmatrix}
    \frac{\partial^2 f}{\partial \lambda^2} &
    \frac{\partial^2 f}{\partial \lambda \partial \mu} &
    \frac{\partial^2 f}{\partial \lambda \partial \nu} &
    \frac{\partial \nabla}{\partial \lambda}\\
    \frac{\partial^2 f}{\partial \mu \partial \lambda} &
    \frac{\partial^2 f}{\partial \mu^2} &
    \frac{\partial^2 f}{\partial \mu \partial \nu} &
    \frac{\partial \nabla}{\partial \mu}\\
    \frac{\partial^2 f}{\partial \nu \partial \lambda} &
    \frac{\partial^2 f}{\partial \nu \partial \mu} &
    \frac{\partial^2 f}{\partial \nu^2} &
    \frac{\partial \nabla}{\partial \nu}\\
    \frac{\partial \nabla}{\partial \lambda} &
    \frac{\partial \nabla}{\partial \mu} &
    \frac{\partial \nabla}{\partial \nu} & 0
    \end{vmatrix}=(\lambda^{14}+\mu^{14}+\nu^{14})+\cdots.\eqno{(4.1.3)}$$
He also formed a function of degree $21$, the functional determinant of
$f$, $\nabla$, and $C$:
$$K=\frac{1}{14} \begin{vmatrix}
         \frac{\partial f}{\partial \lambda} &
         \frac{\partial \nabla}{\partial \lambda} &
         \frac{\partial C}{\partial \lambda}\\
         \frac{\partial f}{\partial \mu} &
         \frac{\partial \nabla}{\partial \mu} &
         \frac{\partial C}{\partial \mu}\\
         \frac{\partial f}{\partial \nu} &
         \frac{\partial \nabla}{\partial \nu} &
         \frac{\partial C}{\partial \nu}
         \end{vmatrix}=-(\lambda^{21}+\mu^{21}+\nu^{21})+\cdots.
\eqno{(4.1.4)}$$
The equation $K=0$ represents the union of the $21$ axes. Klein showed
that under the condition $f=0$ we have, for appropriate values of $k$
and $l$, a relation of the form
$$\nabla^7=k C^3+l K^2,$$
and then it follows further that $f$, $\nabla$, $C$ and $K$, which are
connected by this one equation, generate the whole system of forms under
consideration, and a fortiori the whole system of covariants of $f$. He
determined that the relation among $\nabla$, $C$ and $K$ is
$$(-\nabla)^7=\left(\frac{C}{12}\right)^3-27 \left(\frac{K}{216}\right)^2.
  \eqno{(4.1.5)}$$
Note that $J$ must be equal to $0$ at the contact points of the bitangents,
equal to $1$ at the sextatic points, and equal to $\infty$ at the inflection
points. Thus we have the equation
$$J:J-1:1=\left(\frac{C}{12}\right)^3: 27 \left(\frac{K}{216}\right)^2: -\nabla^7
  \eqno{(4.1.6)}$$
over the Klein quartic curve $f=0$. If we use, instead of $J$, the invariants
$g_2$, $g_3$ and $\Delta$, we can write
$$g_2=\frac{C}{12}, \quad g_3=\frac{K}{216}, \quad \root 7 \of{\Delta}=-\nabla.
  \eqno{(4.1.7)}$$

  The group $\text{PSL}(2, \mathbb{F}_7)$ of $168$ collineations contains
maximal subgroups $G_{21}^{\prime}$ and $G_{24}^{\prime}$ of order $21$
and $24$. Accordingly, there are resolvents of degree eight and of degree
seven. Klein used the expressions in $\lambda$, $\mu$, $\nu$ that, when set
to zero, represent the eight inflection triangles and the two times seven
conics, respectively. Denote by $\delta_{\infty}$ the inflection triangle
to be used as a coordinate triangle, and write
$$\delta_{\infty}=-7 \lambda \mu \nu.\eqno{(4.1.8)}$$
introducing on the right a factor that will later prove convenient. The
following formulas then arise for the remaining inflection triangles,
where $x=0, 1, \ldots, 6$:
$$\aligned
  \delta_{x}
 =&\delta_{\infty}(ST^{x}(\lambda, \mu, \nu))\\
 =&\lambda \mu \nu-(\gamma^{3x} \lambda^3+\gamma^{5x} \mu^3+\gamma^{6x} \nu^3)
  +(\gamma^{6x} \lambda^2 \mu+\gamma^{3x} \mu^2 \nu+\gamma^{5x} \nu^2 \lambda)\\
  &+2(\gamma^{4x} \lambda^2 \nu+\gamma^{x} \nu^2 \mu+\gamma^{2x} \mu^2 \lambda),
\endaligned\eqno{(4.1.9)}$$
where $S$ and $T$ are given by (2.4.29). Next Klein obtained equations
for two of the $14$ conics in terms of $\lambda$, $\mu$, $\nu$:
$$(\lambda^2+\mu^2+\nu^2)+\frac{-1 \pm \sqrt{-7}}{2}
  (\mu \nu+\nu \lambda+\lambda \mu)=0.\eqno{(4.1.10)}$$
Correspondingly, if we denote the left-hand side of the conics by
$c_x$, for $x=0, 1, 2, \ldots, 6$, we get
$$c_x=(\gamma^{2x} \lambda^2+\gamma^{x} \mu^2+\gamma^{4x} \nu^2)
     +\frac{-1 \pm \sqrt{-7}}{2} (\gamma^{6x} \mu \nu+\gamma^{3x}
     \nu \lambda+\gamma^{5x} \lambda \mu).\eqno{(4.1.11)}$$
It is these two expressions that will later lead to the simplest
resolvents of the eighth and seventh degree. For the resolvent of
degree eight, consider the $G_{21}^{\prime}$ generated by the two
substitutions $H$ and $T$, where
$$H=\left(\begin{matrix}
    0 & 1 & 0\\
    0 & 0 & 1\\
    1 & 0 & 0
\end{matrix}\right).$$
It leaves invariant the inflection triangle $\delta_{\infty}$ given
by (4.1.8), and of course $\nabla$, so also the rational function
$\sigma=\delta_{\infty}^2/\nabla$. The latter has the property that
it takes a prescribed value at only $21$ points of the curve $f=0$,
because the pencil of order-six curves $\delta_{\infty}^2-\sigma \nabla$
has three fixed points (the vertices of the coordinate triangle) in
common with $f=0$, each with multiplicity one. Thus, if we use $\sigma$
as a variable, $J$ becomes a rational function of $\sigma$, of degree
eight:
$$J=\frac{\varphi(\sigma)}{\psi(\sigma)}.$$
Klein determined the multiplicity of the individual factors in
$\varphi$, $\psi$, and $\varphi-\psi$ as follows: $J$ becomes
infinite with multiplicity seven at the $24$ inflection points.
At three of these points-the vertices of the coordinate triangle-
$\sigma$ vanishes with multiplicity seven, since $\delta_{\infty}$
has a fourfold zero and $\nabla$ a simple zero. At the remaining
$21$ inflection points $\sigma$ becomes infinite with multiplicity
one because of the denominator $\nabla$. Therefore $\psi(\sigma)$
consists of a simple factor and a sevenfold one, the first vanishing
at $\sigma=0$ and the second at $\sigma=\infty$. Thus, apart from a
constant factor, $\psi(\sigma)$ equals $\sigma$. $J$ vanishes with
multiplicity three at the $56$ contact points of the bitangents. With
respect to the group $G_{21}^{\prime}$ the bitangents fall into two
classes, one with $7$ and one with $21$ elements, so the contact points
are divided into two orbits with $7$ points each and $2$ with $21$ each.
At points of the first kind $\sigma$ takes a certain value with multiplicity
three, and at points of the second with multiplicity one. This means that
$\varphi$ contains two simple and two threefold linear factors. Finally,
$J$ takes the value $1$ with multiplicity two at the $84$ sextatic points.
With respect to $G_{21}^{\prime}$ these points fall into four orbits of
$21$, at each point $\sigma$ takes its value with multiplicity one.
Therefore $\varphi-\psi$ is the square of an expression of degree four
of nonzero discriminant. Klein noted that these are the same condition
on $\varphi$, $\psi$, $\varphi-\psi$ that led him in (\cite{K1}, Section
II, \S 14) to the construction of the modular equation of degree eight:
$$J: J-1:1=(\tau^2+13 \tau+49)(\tau^2+5 \tau+1)^3: (\tau^4+14 \tau^3+63 \tau^2+70 \tau-7)^2:
           1728 \tau,\eqno{(4.1.12)}$$
where for the quartic factor, one finds the decomposition
$$[\tau^2+(7+2 \sqrt{7}) \tau+(21+8 \sqrt{7})] \cdot
  [\tau^2+(7-2 \sqrt{7}) \tau+(21-8 \sqrt{7})].$$
We arrive at this same equation in the present case, if we denote an
appropriate multiple of $\sigma$ by $\tau$. One root $\tau$ of the modular
equation (4.1.12) has the form
$$\tau_{\infty}=-\frac{\delta_{\infty}^2}{\nabla}.\eqno{(4.1.13)}$$
The remaining roots have the values
$$\tau_{x}=-\frac{\delta_x^2}{\nabla}.\eqno{(4.1.14)}$$
Therefore, Klein showed that one can express the roots of the modular equation
of degree eight as a rational function of one points on the Klein quartic curve.

  Now return to the equation:
$$\lambda^3 \mu+\mu^3 \nu+\nu^3 \lambda=0.$$
According to \cite{Radford}, this may be written in the form
$$\frac{\lambda^2}{\nu}+\frac{\mu^2}{\lambda}+\frac{\nu^2}{\mu}=0,$$
and therefore to satisfy it we may take
$$\frac{\lambda^2}{\nu}=z, \quad \frac{\mu^2}{\lambda}=-1, \quad
  \frac{\nu^2}{\mu}=-(z-1),$$
and from these we can express $\lambda$, $\mu$, $\nu$ in terms of $z$.
We have, from the last two,
$$\frac{\nu^4}{\lambda}=-(z-1)^2,$$
and therefore by the first
$$\lambda^8=-\lambda (z-1)^2 z^4,$$
therefore,
$$\lambda=-z^{\frac{4}{7}} (z-1)^{\frac{2}{7}},$$
and the others are
$$\mu=-z^{\frac{2}{7}} (z-1)^{\frac{1}{7}}, \quad
  \nu=z^{\frac{1}{7}} (z-1)^{\frac{4}{7}}.$$

  In fact, the modular equation (4.1.12) can be transformed as follows:
write $z^2$ instead of $\tau$, $27 g_3^2/\Delta$ instead of $J-1$, and
take the square root of both sides, to obtain
$$z^8+14 z^6+63 z^4+70 z^2-\frac{216 g_3}{\sqrt{\Delta}} z-7=0.\eqno{(4.1.15)}$$
We can further replace $216 g_3/\sqrt{\Delta}$ with $K/\sqrt{-\nabla^7}$, by
(4.1.7), and replacing also $z$ by its value $\delta/\sqrt{-\nabla}$, given by
(4.1.13), (4.1.14), the result is
$$\delta^8-14 \delta^6 \nabla+63 \delta^4 \nabla^2-70 \delta^2 \nabla^3-K \delta
  -7 \nabla^4=0.\eqno{(4.1.16)}$$
Klein showed that we would have arrived at the same equation (4.1.16) if we
have taken the polynomial approach. For the simplest polynomial function
of $\lambda$, $\mu$, $\nu$ that is left invariant by $G_{21}^{\prime}$ is
$\delta_{\infty}$ given by (4.1.8). Under the $168$ collineations $\delta$
takes eight distinct values $\delta_{\infty}$ and $\delta_x$ ($x=0, 1,
\ldots, 6$) given by (4.1.9), whose symmetric function must be a polynomial
function of $\nabla$, $C$, $K$ (since $f$ is taken to equal $0$). Therefore
$\delta$ satisfies an equation of the eighth degree, which in view of the
degrees of $\nabla$, $C$, $K$, must have the form
$$\delta^8+a \nabla \delta^6+b \nabla^2 \delta^4+c \nabla^3 \delta^2
  +d K \delta+e \nabla^4=0,$$
and if the coefficients $a$, $b$, $c$, $d$, $e$ are determined by
substituting for $\delta$, $\nabla$, $K$ their values in terms of
$\lambda$, $\mu$, $\nu$ and taking into account that $f=0$, we recover
(4.1.16). The eight roots of (4.1.16) can be expressed as follows, by
virtue of (4.1.8) and (4.1.9):
$$\aligned
  \delta_{\infty} &=-7 \lambda \mu \nu,\\
  \delta_x &=\lambda \mu \nu-\gamma^{-x} (\nu^3-\lambda^2 \mu)
             -\gamma^{-4 x} (\lambda^3-\mu^2 \nu)-\gamma^{-2 x} (\mu^3-\nu^2 \lambda)\\
           &+2 \gamma^{x} \nu^2 \mu+2 \gamma^{4x} \lambda^2 \nu+2 \gamma^{2x} \mu^2 \lambda.
\endaligned$$
Klein pointed out in (\cite{K1}, section II, \S 18) that (4.1.15), and therefore
also (4.1.16), is a Jacobian equation of degree eight, that is, the square roots
of its roots can be written in terms of four quantities $A_0$, $A_1$, $A_2$, $A_3$
as follows:
$$\left\{\aligned
  \sqrt{\delta_{\infty}} &=\sqrt{-7} A_0,\\
  \sqrt{\delta_x} &=A_0+\gamma^{\rho x} A_1+\gamma^{4 \rho x} A_2+\gamma^{2 \rho x} A_3,
\endaligned\right.\eqno{(4.1.17)}$$
where $\rho$ is any integer not divisible by $7$.

  For the resolvent of degree seven, set
$$x_1=\frac{\lambda+\mu+\nu}{\alpha-\alpha^2}, \quad
  x_2=\frac{\lambda+\alpha \mu+\alpha^2 \nu}{\alpha-\alpha^2}, \quad
  x_3=\frac{\lambda+\alpha^2 \mu+\alpha \nu}{\alpha-\alpha^2},$$
where $\alpha=e^{\frac{2 \pi i}{3}}$. The equation of the Klein quartic
curve becomes
$$0=f=\frac{1}{3} \left\{x_1^4+3 x_1^2 x_2 x_3-3 x_2^2 x_3^2+
      x_1 [(1+3 \alpha^2) x_2^3+(1+3 \alpha) x_3^3]\right\}.$$
To get rid of the cubic roots of unity, Klein further set
$$x_1=\frac{y_1}{\root 3 \of{7}}, \quad x_2=y_2 \root 3 \of{3\alpha+1}, \quad
  x_3=y_3 \root 3 \of{3 \alpha^2+1}$$
and got
$$0=f=\frac{1}{21 \root 3 \of{7}} [y_1^4+21 y_1^2 y_2 y_3-147 y_2^2 y_3^2
     +49 y_1 (y_2^3+y_3^3)].$$
Klein wrote the equation of his quartic curve in three different ways in the
form $pqrs-w^2=0$, where $p$, $q$, $r$, $s$ are bitangents and $w$ is the
conic that goes through the contact points. The first such expression is
$$0=\frac{1}{21 \root 3 \of{7}}[49 y_1 (y_1+y_2+y_3)(y_1+\alpha y_2+\alpha^2 y_3)
    (y_1+\alpha^2 y_2+\alpha y_3)-3 (4 y_1^2-7 y_2 y_3)^2],$$
and the other two are
$$\aligned
 0=\frac{1}{21 \root 3 \of{7}} \left\{\frac{y_1}{7 \cdot 8^3}\right.
   &[(-7 \pm 3 \sqrt{-7}) y_1+56 y_2+56 y_3]\\
   \times &[(-7 \pm 3 \sqrt{-7}) y_1+56 \alpha y_2+56 \alpha^2 y_3]\\
   \times &[(-7 \pm 3 \sqrt{-7}) y_1+56 \alpha^2 y_2+56 \alpha y_3]\\
   &\left.-3 \left(\frac{1 \pm 3 \sqrt{-7}}{16} y_1^2-7 y_2 y_3\right)^2\right\},
\endaligned\eqno{(4.1.18)}$$
which yields the substitutions in a $G_{24}^{\prime}$. Set
$$\left\{\aligned
  \mathfrak{z}_1 &=(21 \mp 9 \sqrt{-7}) y_1,\\
  \mathfrak{z}_2 &=(-7 \pm 3 \sqrt{-7}) y_1+56 y_2+56 y_3,\\
  \mathfrak{z}_3 &=(-7 \pm 3 \sqrt{-7}) y_1+56 \alpha y_2+56 \alpha^2 y_3,\\
  \mathfrak{z}_4 &=(-7 \pm 3 \sqrt{-7}) y_1+56 \alpha^2 y_2+56 \alpha y_3,
\endaligned\right.$$
so that $\sum \mathfrak{z}_i=0$. Then (4.1.18) becomes, apart from a scalar
factor,
$$\left(\sum_{i=1}^{4} \mathfrak{z}_i^2\right)^2-(14 \pm 6 \sqrt{-7})
  \mathfrak{z}_1 \mathfrak{z}_2 \mathfrak{z}_3 \mathfrak{z}_4=0.\eqno{(4.1.19)}$$
and this equation is invariant under the $24$ collineations determined
by the permutations of the $\mathfrak{z}_i$. Therefore, these are the
collineations of the $G_{24}^{\prime \prime}$ in question. We see that
the collineation of a $G_{24}^{\prime \prime}$ always leave invariant
a certain conic
$$\sum_{i=1}^{4} \mathfrak{z}_i^2=0,\eqno{(4.1.20)}$$
which goes through the contact points of the corresponding bitangents.
Since there are $2 \cdot 7$ groups $G_{24}^{\prime \prime}$ and all
bitangents have equal right, there are $2 \cdot 7$ such conics, and
by taking any seven together and intersecting with the Klein quartic
curve we get all the contact points of bitangents. Next Klein obtained
equations for two of the $14$ conics, by taking the equation (4.1.20)
and expressing it first in terms of the $y_i$ and from there in terms
of $\lambda$, $\mu$, $\nu$:
$$(\lambda^2+\mu^2+\nu^2)+\frac{-1 \sqrt{-7}}{2} (\mu \nu+\nu \lambda+\lambda \mu)=0,$$
which is just (4.1.10). The substitutions of a $G_{24}^{\prime \prime}$
always leave invariant a conic $c_x$ given by (4.1.11) that goes through
the contact points of four bitangents. Now form the rational function
$$\xi=\frac{c_x^3}{\nabla}.\eqno{(4.1.21)}$$
Since the numerator and denominator are invariant under the substitutions
in the $G_{24}^{\prime \prime}$, and since the pencil of sixth-order curves
$\nabla-\xi c_x^3=0$ has no fixed intersection with the Klein quartic
curve, we conclude that $\xi$ takes a given value at exactly the points
of an orbit of the $G_{24}^{\prime \prime}$. Therefore $J$ is a rational
function of degree seven of $\xi$:
$$J=\frac{\varphi(\xi)}{\psi(\xi)}.$$
Klein considered again the values $J=\infty$, $0$, $1$. The $24$ inflection
points, where $J$ becomes infinite with multiplicity seven, form a single
orbit of the $G_{24}^{\prime \prime}$, each point appearing once. Thus
$\psi(\xi)$ is the seventh power of a linear factor. But $\xi$ itself is
infinite at the inflection points, because of (4.1.21). Therefore $\psi(\xi)$
is a constant. Of the $56$ contact points of the $28$ bitangents eight
lie on $c_x=0$, so $\xi$ vanishes with order three at those points. The
other $48$ split into $2$ orbits of as $24$ (each corresponding to $12$
tangents). Thus $\varphi$ contains the simple factor $\xi$ and the cube
of a quadratic factor of nonzero discriminant. Finally, the $84$ sextatic
points fall into three orbits of $12$ points each and two of $24$ points
each. Thus $\varphi-\psi$ contains a simple cubic factor and the square
of a quadratic factor. Again, these are the requirements on $\varphi$ and
$\psi$ that led in (\cite{Klein1879}, \S 7) to the construction of the
simplest equation of degree seven, which has the following form:
$$\aligned
  &J:J-1:1=\mathfrak{z} [\mathfrak{z}^2-2^2 \cdot 7^2 (7 \mp \sqrt{-7})
           \mathfrak{z}+2^5 \cdot 7^4 (5 \mp \sqrt{-7})]^3\\
         :&[\mathfrak{z}^3-2^2 \cdot 7 \cdot 13 (7 \mp \sqrt{-7})
           \mathfrak{z}^2+2^6 \cdot 7^3 (88 \mp 23 \sqrt{-7}) \mathfrak{z}\\
          &-2^8 \cdot 3^3 \cdot 7^4 (35 \mp 9 \sqrt{-7})][\mathfrak{z}^2-2^4
           \cdot 7 (7 \mp \sqrt{-7}) \mathfrak{z}+2^5 \cdot 7^3 (5 \mp \sqrt{-7})]^2\\
         :&\mp 2^{27} \cdot 3^3 \cdot 7^{10} \cdot \sqrt{-7}.
\endaligned\eqno{(4.1.22)}$$
Klein first transformed (4.1.22) by setting $\mathfrak{z}=z^3$ and $J=g_2^3/\Delta$,
and then he took the cubic root of both sides:
$$z^7-2^2 \cdot 7^2 (7 \mp \sqrt{-7}) z^4+2^5 \cdot 7^4 (5 \mp \sqrt{-7}) z
  \mp 2^9 \cdot 3 \cdot 7^3 \sqrt{-7} \frac{g_2}{\root 3 \of{\Delta}}=0.\eqno{(4.1.23)}$$
Klein substituted $C/(12 \root 3 \of{-\nabla^7})$ for $g_2/\root 3 \of{\Delta}$
and $\pm 2 \sqrt{-7} c/\root 3 \of{\nabla}$ for $z$. The result is
$$c^7+\frac{7}{2} (-1 \mp \sqrt{-7}) \nabla c^4-7
  \left(\frac{5 \mp \sqrt{-7}}{2}\right) \nabla^2 c-C=0.\eqno{(4.1.24)}$$
The roots $z$ of (4.1.23) and $\mathfrak{z}$ of (4.1.22) have the following
values in terms of $\lambda$, $\mu$, $\nu$:
$$\aligned
 &z=\mathfrak{z}^{\frac{1}{3}}\\
=&\frac{\pm 2 \sqrt{-7} [(\gamma^{2x}
  \lambda^2+\gamma^{x} \mu^2+\gamma^{4x} \nu^2)+\frac{-1 \mp \sqrt{-7}}{2}
  (\gamma^{6 x} \mu \nu+\gamma^{3 x} \nu \lambda+\gamma^{5 x} \lambda \mu)]}
  {\root 3 \of{\nabla}},
\endaligned\eqno{(4.1.25)}$$
and Klein had explicitly written the $\mathfrak{z}$'s as rational functions
of one point on his quartic curve. Naturally, the polynomial approach leads
to the same equation. Indeed, the lowest polynomial function of $\lambda$,
$\mu$, $\nu$ that remains invariant under a $G_{24}^{\prime \prime}$ is
exactly the corresponding $c_x$, and this $c_x$ must satisfy an equation
of degree seven, whose coefficients are polynomials in $\nabla$, $C$, $K$,
and which therefore has the form
$$c^7+\alpha \nabla c^4+\beta \nabla^2 c+\gamma C=0,$$
where $\alpha$, $\beta$, $\gamma$ are to be determined by the substitution
of values for $\lambda$, $\mu$, $\nu$.

\begin{center}
{\large\bf 4.2. The case of $p=11$ and Klein cubic threefold}
\end{center}

  For $p=11$, consider the Galois representation:
$$\rho_{11}: \text{Gal}(\overline{\mathbb{Q}}/\mathbb{Q}) \rightarrow
  \text{Aut}(\mathcal{L}(X(11))),$$
there are two such equations with degree $12$ and $11$, respectively. Correspondingly,
there are two distinct kinds of surjective realizations as well as two distinct
kinds of modularities. For the degree twelve equation with Galois group isomorphic
to $\text{PSL}(2, \mathbb{F}_{11})$, Klein constructed (see \cite{K3} and \cite{KF2},
Section five, Chapter five) the following equation (4.2.17) of degree twelve
corresponding to a maximal subgroup $G_{55}$ of $G_{660}$, which is the Jacobian
multiplier equation of degree twelve. He showed that it can be reduced to the the
modular equation of degree twelve which corresponding to the eleven-order
transformation of $j$-functions.

  For the eleventh-degree equation with a Galois group isomorphic with the
$\text{PSL}(2, \mathbb{F}_{11})$, Klein constructed the following equation
of degree eleven corresponding to a maximal subgroup $G_{60}$ of $G_{660}$.
At first, he obtained the eleventh-order transformation (4.2.11) of $j$-functions.
Then, he showed that the invariant-polynomial approach leads to the following
equation (4.2.16), which are the same as in (4.2.11). In particular, Klein
constructed a system of equations, which he referred to as $H_{ik}=0$ given
by (4.2.8). It turns out that this system is just the locus $\mathcal{L}(X(11))$
given by (2.4.34) and (2.4.35). In particular, the equation (4.2.16) is defined
over the locus $\mathcal{L}(X(11))$. Its coefficients are invariant under the
action of the automorphism group $\text{Aut}(\mathcal{L}(X(11)))$. Klein proved
(see \cite{K3}) that the Galois group of (4.2.16) is isomorphic with the
automorphism group of the locus $\mathcal{L}(X(11))$.

  In fact, in the above two cases, degree twelve and eleven, the reduction
process is always given by the theta constants of order eleven, and the locus
$\mathcal{L}(X(11))$ is parameterized by the theta constants of order eleven.
This is the reason that both the Jacobian multiplier equation of degree twelve
and (4.2.16) are defined over the locus $\mathcal{L}(X(11))$. Therefore, we
have proved the surjective realization as well as the modularity of the
Galois representation $\rho_{11}$.

  Now, we will give the detailed construction. In his paper \cite{K3}
or \cite{KF2}, Section five, Chapter five, Klein's intention is not to
communicate all invariant entire functions of $y$; in any case this would
be a far reaching and possibly difficult task. Rather Klein defined three
of them which he would use later. The first of which is the function of
the third degree
$$\nabla=y_1^2 y_9+y_4^2 y_3+y_5^2 y_1+y_9^4 y_4+y_3^2 y_5,\eqno{(4.2.1)}$$
clearly the lowest invariant function. The second is the Hessian
determinant of the fifth degree:
$$\begin{vmatrix}
  y_9 & 0 & y_5 & y_1 & 0\\
  0 & y_3 & 0 & y_9 & y_4\\
  y_5 & 0 & y_1 & 0 & y_3\\
  y_1 & y_9 & 0 & y_4 & 0\\
  0 & y_4 & y_3 & 0 & y_5
\end{vmatrix},\eqno{(4.2.2)}$$
whose first sub-determinant will be of interest later. The third is
the function of the eleventh degree:
$$C=(y_1^{11}+y_4^{11}+y_5^{11}+y_9^{11}+y_3^{11})+\cdots.\eqno{(4.2.3)}$$
Klein constructed a very simple function, that is invariant under the
operations $C$ and $S$, namely the sum of $y$:
$$p_{\infty}=y_1+y_4+y_5+y_9+y_3,\eqno{(4.2.4)}$$
where
$$C=\left(\begin{matrix}
     & 1 &   &   &  \\
     &   & 1 &   &  \\
     &   &   & 1 &  \\
     &   &   &   & 1\\
   1 &   &   &   &
\end{matrix}\right)$$
and $S$ is given by (2.4.32). If one applies the collineation $TST^{-1}$
to this, a short calculation results in
$$\aligned
  &p_0=p_{\infty}(TST^{-1}(y_1, y_4, y_5, y_9, y_3))\\
 =&\frac{1}{\sqrt{-11}}\left\{[2(\rho^7-\rho)+(\rho^9-\rho^{10})] y_1+
   [2(\rho^6-\rho^4)+(\rho^3-\rho^7)] y_4\right.\\
  &+[2(\rho^2-\rho^5)+(\rho-\rho^6)]y_5+[2(\rho^8-\rho^9)+(\rho^4-\rho^2)] y_9\\
  &\left.+[2(\rho^{10}-\rho^3)+(\rho^5-\rho^8)] y_3 \right\},
\endaligned\eqno{(4.2.5)}$$
where $T$ is given by (2.4.33), and if one cyclically permutes the five $y$,
one obtains a further four distinct expressions, which we will refer to as
$$p_1, p_2, p_3, p_4.$$
Since $p_{\infty}$ remains unchanged by $10$ collineations of the subgroups,
it thus is six-valued under action by collection of collineations; that
is the six expressions $p$ are permuted amongst each other by the $60$
collineations of the subgroup. The symmetric functions of the six $p$
are invariant under all the collineations of the subgroup and are thus
eleven-valued under the $660$ collineations generated by $S$ and $T$,
they must therefore belong to the functions that are left unchanged
entirely. Accordingly, in order to have eleven-valued functions that
are as low as possible, one calculates the smallest non-vanishing
symmetric functions of $p$, i.e. the sum of squares and the sum of
cubes. One thus finds the following functions:

1) the function of the second degree:
$$\aligned
  \varphi_0= &(y_1^2+y_4^2+y_5^2+y_9^2+y_3^2)\\
           - &(y_1 y_9+y_4 y_3+y_5 y_1+y_9 y_4+y_3 y_5)\\
  +\frac{-1+\sqrt{-11}}{2} &(y_1 y_4+y_4 y_5+y_5 y_9+y_9 y_3+y_3 y_1).
\endaligned\eqno{(4.2.6)}$$

2) The function of the third degree
$$\aligned
  f_0= &(y_1^3+y_4^3+y_5^3+y_9^3+y_3^3)\\
     +3 &(y_1^2 y_3+y_4^2 y_1+y_5^2 y_4+y_9^2 y_5+y_3^2 y_9)\\
     -3 &(y_1 y_4 y_9+y_4 y_5 y_3+y_5 y_9 y_1+y_9 y_3 y_4+y_3 y_1 y_5)\\
 +\frac{1+\sqrt{-11}}{2} &(y_1^2 y_5+y_4^2 y_9+y_5^2 y_3+y_9^2 y_1+y_3^2 y_4)\\
 -\frac{1+\sqrt{-11}}{2} &(y_1 y_4 y_5+y_4 y_5 y_9+y_5 y_9 y_3+y_9 y_3 y_1+y_3 y_1 y_4)\\
 -(1+\sqrt{-11}) &(y_1^2 y_4+y_4^2 y_5+y_5^2 y_9+y_9^2 y_3+y_3^2 y_1).
\endaligned\eqno{(4.2.7)}$$
The function $\varphi_0$ coincides with $\frac{-1+\sqrt{-11}}{12} \sum p^2$;
the function $f_0$ differs from $\frac{-\sqrt{-11} \sum p^3}{6}$ on only
one term, which is a numerical multiple of $\nabla$ given by (4.2.1). The
eleven values which $\varphi_0$ and $f_0$ take on under action of the
660 collineations, which Klein called $\varphi_{\nu}$ and $f_{\nu}$
respectively, arise from $\varphi_0$ and $f_0$, when one replaces
$y_{\kappa^2}$ with $\rho^{\kappa^2 \nu} \cdot y_{\kappa^2}$ in the
collineation $T^{\nu}$. Klein obtained the following results for the
transformation of the eleventh degree of elliptic functions:

(1) The Galois resolvent of the $660$th degree can be written in the
following way: One takes the five related quantities
$$y_1: y_4: y_5: y_9: y_3$$
subject to the $15$ relations which come from constructing all
sub-determinants of $H$ in (4.2.2) and setting them equal to zero.
One thus get a system of equations, which Klein referred to as
$$H_{ik}=0 \eqno{(4.2.8)}$$
and whose equations arise from the following three
$$\left\{\aligned
  0 &=y_4 y_5 y_9 y_3-y_1^2 y_5 y_3+y_1^2 y_4^2+y_3^3 y_1,\\
  0 &=y_1^2 y_5 y_9-y_4^2 y_5 y_3-y_3^2 y_1 y_9,\\
  0 &=y_4^3 y_9+y_9^3 y_5+y_3^3 y_1,
\endaligned\right.\eqno{(4.2.9)}$$
up to cyclic permutation of the $y$. Set
$$\frac{-C^3}{1728 \nabla^{11}}=J,\eqno{(4.2.10)}$$
where $\nabla$ denotes the function of the third degree in (4.2.1), and
$C$ denotes the function of the eleventh degree in (4.2.3).

(2) There exist two simplest forms of the resolvent of the eleventh degree.
The first one is given by
$$\aligned
  &J:J-1:1=[z^2-3z+(5-\sqrt{-11})] \cdot\\
    \cdot &\left(z^3+z^2-3 \cdot \frac{1+\sqrt{-11}}{2} z+\frac{7-\sqrt{-11}}{2}\right)^3\\
        : &\left[z^3+4 z^2+\frac{7-5 \sqrt{-11}}{2} z+(4-6 \sqrt{-11})\right] \cdot\\
    \cdot &\left[z^4-2 z^3+3 \cdot \frac{1-\sqrt{-11}}{2} z^2+(5+\sqrt{-11}) z
          -3 \cdot \frac{5+\sqrt{-11}}{2}\right]^2\\
        : &-1728,
\endaligned\eqno{(4.2.11)}$$
its eleven roots are given by the formula:
$$z_{\nu}=\frac{f_{\nu}}{\nabla}\eqno{(4.2.12)}$$
where $f_{\nu}$ is defined by the equation (4.2.7). The second form is given by
$$\aligned
  \xi^{11} &-22 \xi^8+11(9-2 \sqrt{-11}) \xi^5-11 \cdot \frac{12 g_2}
  {\root 3 \of{\Delta}} \xi^4+88 \sqrt{-11} \xi^2\\
  &-11 \cdot \frac{-3+\sqrt{-11}}{2} \frac{12 g_2}{\root 3 \of{\Delta}} \xi
   -\frac{144 g_2^2}{\root 3 \of{\Delta^2}}=0,
\endaligned\eqno{(4.2.13)}$$
and its roots are
$$\xi_{\nu}=\frac{\varphi_{\nu}}{\nabla^{\frac{2}{3}}},\eqno{(4.2.14)}$$
where $\varphi_{\nu}$ are understood as the functions (4.2.6).

  It remains only to show that how the equation (4.2.13) is connected with
equation (4.2.11). Klein showed that the following relation holds (naturally
still given the $H_{ik}=0$):
$$\varphi_{\nu}^3=f_{\nu}^2-3 f_{\nu} \nabla+(5-\sqrt{-11}) \nabla^2,
  \eqno{(4.2.15)}$$
and that one can get equation (4.2.11) from equation (4.2.13) by setting
$$\xi^3=z^2-3z+(5-\sqrt{-11}).$$
Naturally, the polynomial approach leads to the same equation:
$$\varphi^{11}+\alpha \nabla^2 \cdot \varphi^8+\beta \nabla^4 \cdot \varphi^5
 +\gamma \nabla C \cdot \varphi^4+\delta \nabla^6 \cdot \varphi^2+\varepsilon
  \nabla^3 C \cdot \varphi+\zeta \cdot C^2=0.\eqno{(4.2.16)}$$

  In fact, in the equation $J=F(z)$, we require only one parameter $J$ to
be present. We are thus faced with properly finding a curve in the projective
space $\mathbb{P}^4=\{ (y_1: y_4: y_5: y_9: y_3) \}$, which is mapped to
itself by the $660$ collineations, and to carry out the problem of $y$
that applies to specifically to it. Regarding this curve, we know that
it must be the image of the Galois resolvent of the transformation equation.
As Klein showed that, we now have that the Galois resolvent is presented
by a Riemann surface, which is $660$ sheeted across the plane of $J$ and
whose sheets are in triples at $J=0$, in pairs at $J=1$, and eleven-fold
at $J=\infty$, and otherwise are not joined, that is their genus is equal
to $26$. Therefore, a rational function $J$ must exist on our curve, such
that it takes on every value at and only at those $660$ points, which arise
from the $660$ collineations. In each of these groups of $660$ points each
there can only be three that consist of a smaller number of multiply counted
points: a group of $220$ points counted three times, a group of $330$ points
counted two times, and a group of $60$ counted eleven times. The genus of
the curve is naturally thus equal to $26$:
$$2 \cdot 26-2=660 \left(1-\frac{1}{2}-\frac{1}{3}-\frac{1}{11}\right).$$
Next, we consider the equation $J=F(z)$. Given everything else it implies
that there exists a rational function of $z$ on our curve, taking on every
value at $60$ and only $60$ points, which stem from a subgroup of collineations
of index $11$. Further properties of the function $F$: That $F(z)$ is a rational
entire function of the eleventh degree, and must contain a cubic factor three
times and $F(z)-1$ must contain a quartic factor twice. That $F(z)$ is an entire
function of the eleventh degree follows from the fact that for the $660$ points
corresponding to a value of $J$ split with respect to the $60$ collineations
of subgroup into $11 \cdot 60$, however to $60$ points $J=\infty$ all arise
from the collineations of the subgroup. The other properties follow from the
behavior of the $220$ points of $J=0$ and the $330$ points of $J=1$. Among the
$660$ collineations, we know that there exist $2 \cdot 55$ of period $3$,
$55$ of period $2$. Thus, each collineation of period $3$ fixes $4$ points of
$J=0$, and each collineation of period $2$ fixes $6$ points of $J=1$. However,
the subgroup of index $11$ contains $2 \cdot 10$ collineations of period $3$,
$15$ of period $2$. The $220$ points of $J=0$ thus separate with respect to
this in $2 \cdot 20+3 \cdot 60$ and the $330$ points of $J=1$ separate into
$3 \cdot 30+4 \cdot 60$. And exactly this is meant by the stated properties
of $F$. Now, we have that the lowest rational function, depending only on
the relations of the $y$ and unchanged by the $60$ collineations of a
subgroup, $z=\frac{f_{\nu}}{\nabla}$ is a function of the third degree. As
we presuppose this function should be the function taking on each value on
$60$ points on our curve; we will thus introduce the hypothesis that our
curve be of the $20$th degree, that it lies on neither $f_{\nu}=0$ nor on
$\nabla=0$, and does not cross the common intersection of $f_{\nu}=0$ and
$\nabla=0$. We will further assume, that it does not lie on $C=0$ and also
does not cross the intersection of $C=0$, $\nabla=0$. Then $\frac{C^3}{\nabla^{11}}$
is a function which takes on each value at $660$ associated points; for $C=0$
one only gets $220$ separated points, and for $\nabla=0$ only $60$. One sees:
One must set $J=k \frac{C^3}{\nabla^{11}}$, where $k$ is a numerical constant.
Now Klein stated that if the genus of our curve is equal to $26$, there directly
exists only one group of multiply counted points besides the group of $220$
triply counted points and the group of $60$ counted eleven times, namely
the $330$ counted twice. It is the matter of finding a curve in the space
of $y$ with degree $20$ and genus $26$, which is mapped to itself by the
$660$ collineations, lies on neither $f_{\nu}=0$, $\nabla=0$, nor finitely
on $C=0$, and does not cross the intersection of $f_{\nu}=0$ or $\nabla=0$,
nor the intersection of $C=0$ or $\nabla=0$. This leads to the double curve
of $H=0$: if our curve of the $20$th degree exists, it must lie on the
threefold $H=0$. This is true since otherwise there would be $100$
points of intersection with $H=0$, and these $100$ points would have
to be permuted amongst each other by the $660$ collineations, which
is impossible. Now one recalls from ordinary space geometry that the
Hessian of a threefold of the third degree has $10$ nodes. For this
purpose one simultaneously sets all first sub-determinants of the
Hessian determinant equal to zero. One can proceed in the same manner
in the case of $5$ variables and obtains the following general theorem
via use of known methods: in the case of $5$ variables the Hessian of
a threefold of the third degree possess a double curve of the $20$th
degree and genus $26$. Klein showed that the curve of the $20$th degree
we seek is the double curve of the Hessina threefold $H=0$. Denote the
five points, where four of the five $y$ vanish, as the points I, IV, V,
IX, III. Then Klein proved the following:

(1) The five points I, IV, V, IX, III belong to our curve.

(2) Our curve lies on neither $C=0$ nor $f_{\nu}=0$.

(3) Every plane $y_{k^2}=0$ intersects our curve at only four points,
namely at those four of the five points I, IV, V, IX, III, which are
not named after the index of $y_{k^2}$.

(4) The five points are simple points of our curve.

(5) The curve is of degree $20$.

(6) $\nabla=0$ has precisely $60$ points of intersection with our curve,
and neither $C$ nor $f_{\nu}$ vanish at these points.

(7) The genus of our curve is equal to $26$.

  Besides the equation (4.2.11) we just obtained we can obtain a second
with a similar amount of ease, given one starts with the function of
the second degree $\varphi_{\nu}$ (4.2.6) instead of the eleven-valued
function of the third degree $f_{\nu}$ (4.2.7). Klein proved the following:
Subjected to the relation $H_{ik}=0$, we can reduce the collection of entire
functions of $y$ invariant under the $660$ collineations to entire functions
of $\nabla$ and $C$. Therefore, by virtue of the relations $H_{ik}=0$, the
eleven $\varphi_{\nu}$ satisfy an equation of the eleventh degree whose
coefficients are entire functions of $\nabla$ and $C$. With respect to the
degree of the functions to be considered, we obtain the equation (4.2.16).

  Connection with equations of the twelfth degree

  Klein also showed how the quantities $y$ related to the multiplier
equation of the twelfth degree, which is a Jacobian equation:
$$\aligned
  &z^{12}-90 \cdot 11 \cdot \root 2 \of{\Delta} \cdot z^6
   +40 \cdot 11 \cdot 12 g_2 \cdot \root 3 \of{\Delta} \cdot z^4
   -15 \cdot 11 \cdot 216 g_3 \cdot \root 4 \of{\Delta} \cdot z^3\\
 &+2 \cdot 11 \cdot (12 g_2)^2 \cdot \root 6 \of{\Delta} \cdot z^2
  -12 g_2 \cdot 216 g_3 \cdot \root 12 \of{\Delta} \cdot z-11 \cdot \Delta=0.
\endaligned\eqno{(4.2.17)}$$
Here, Klein set the following:
$$\left\{\aligned
  \sqrt{z_{\infty}} &=\sqrt{-11} \cdot A_0,\\
  \sqrt{z_{\nu}} &=A_0+\rho^{\nu} A_1+\rho^{4 \nu} A_4+\rho^{5 \nu} A_5
                  +\rho^{9 \nu} A_9+\rho^{3 \nu} A_3, \quad (\nu=0, 1, \ldots, 10)
\endaligned\right.\eqno{(4.2.18)}$$
In fact, Klein only deviated from the Jacobian description by choosing
the indices for $A$ to be quadratic residues modulo $11$.

  In \cite{K4}, Klein pointed that among the $15$ relations in (4.2.9),
one initially obtains only those $10$ which consist of three members.
From these, however, one can derive the five four-membered ones as an
algebraic sequence. Because it is:
$$\aligned
  &y_4 y_5 y_9 y_3-y_1^2 y_5 y_3+y_1^2 y_4^2+y_3^3 y_1\\
 =&\frac{1}{y_9} \left\{-y_3 (y_1^2 y_5 y_9-y_4^2 y_5 y_3-y_3^2 y_1 y_9)
   -y_4 (y_3^2 y_4 y_5-y_1^2 y_4 y_9-y_9^2 y_3 y_5)\right\}.
\endaligned$$
This leads us to put (2.4.34) and (2.4.35).

\begin{center}
{\large\bf 4.3. The case of $p=13$ and our quartic fourfold}
\end{center}

  For $p=13$, consider the Galois representation:
$$\rho_{13}: \text{Gal}(\overline{\mathbb{Q}}/\mathbb{Q}) \rightarrow
  \text{Aut}(\mathcal{L}(X(13))),$$
there are two distinct equations with the same degree fourteen and the same
Galois group isomorphic with $\text{PSL}(2, \mathbb{F}_{13})$. One is the Jacobian
multiplier equation of degree fourteen given by the following equation (4.3.8),
the other has degree fourteen, but it is not the Jacobian multiplier equation (see
(4.3.18)). Both of them are defined over the locus $\mathcal{L}(X(13))$, which is parameterized
by theta constants of order $13$. We have verified this in \cite{Y6}, Theorem 7.6.
Moreover, their coefficients are invariant under the action of
$\text{Aut}(\mathcal{L}(X(13)))$, and their Galois groups are isomorphic to
$\text{Aut}(\mathcal{L}(X(13)))$. Correspondingly, there are two distinct kinds
of surjective realizations as well as two distinct kinds of modularities. On
the other hand, for the degree fourteen equation with Galois group isomorphic
to $\text{PSL}(2, \mathbb{F}_{13})$, Klein constructed the following equation
(4.3.19) of degree fourteen corresponding to a maximal subgroup $G_{78}$ of
$G_{1092}$. For its modern form, see (3.3.8), (3.3.9) and (3.3.10).

  Now, we will give the detailed construction. Put $\theta_1=\zeta+\zeta^3+\zeta^9$,
$\theta_2=\zeta^2+\zeta^6+\zeta^5$, $\theta_3=\zeta^4+\zeta^{12}+\zeta^{10}$,
and $\theta_4=\zeta^8+\zeta^{11}+\zeta^7$. We find that
$$\left\{\aligned
  &\theta_1+\theta_2+\theta_3+\theta_4=-1,\\
  &\theta_1 \theta_2+\theta_1 \theta_3+\theta_1 \theta_4+
   \theta_2 \theta_3+\theta_2 \theta_4+\theta_3 \theta_4=2,\\
  &\theta_1 \theta_2 \theta_3+\theta_1 \theta_2 \theta_4+
   \theta_1 \theta_3 \theta_4+\theta_2 \theta_3 \theta_4=4,\\
  &\theta_1 \theta_2 \theta_3 \theta_4=3.
\endaligned\right.$$
Hence, $\theta_1$, $\theta_2$, $\theta_3$ and $\theta_4$ satisfy
the quartic equation $z^4+z^3+2 z^2-4z+3=0$,
which can be decomposed as two quadratic equations
$$\left(z^2+\frac{1+\sqrt{13}}{2} z+\frac{5+\sqrt{13}}{2}\right)
  \left(z^2+\frac{1-\sqrt{13}}{2} z+\frac{5-\sqrt{13}}{2}\right)=0$$
over the real quadratic field $\mathbb{Q}(\sqrt{13})$. Therefore, the
four roots are given as follows:
$$\left\{\aligned
  \theta_1=\frac{1}{4} \left(-1+\sqrt{13}+\sqrt{-26+6 \sqrt{13}}\right),\\
  \theta_2=\frac{1}{4} \left(-1-\sqrt{13}+\sqrt{-26-6 \sqrt{13}}\right),\\
  \theta_3=\frac{1}{4} \left(-1+\sqrt{13}-\sqrt{-26+6 \sqrt{13}}\right),\\
  \theta_4=\frac{1}{4} \left(-1-\sqrt{13}-\sqrt{-26-6 \sqrt{13}}\right).
\endaligned\right.$$
Moreover, we find that
$$\left\{\aligned
  \theta_1+\theta_3+\theta_2+\theta_4 &=-1,\\
  \theta_1+\theta_3-\theta_2-\theta_4 &=\sqrt{13},\\
  \theta_1-\theta_3-\theta_2+\theta_4 &=-\sqrt{-13+2 \sqrt{13}},\\
  \theta_1-\theta_3+\theta_2-\theta_4 &=\sqrt{-13-2 \sqrt{13}}.
\endaligned\right.$$

  Let us study the action of $S T^{\nu}$ on $\mathbb{P}^5$, where
$\nu=0, 1, \ldots, 12$. Put
$$\alpha=\zeta+\zeta^{12}-\zeta^5-\zeta^8, \quad
   \beta=\zeta^3+\zeta^{10}-\zeta^2-\zeta^{11}, \quad
   \gamma=\zeta^9+\zeta^4-\zeta^6-\zeta^7.$$
We find that
$$\aligned
  &13 ST^{\nu}(z_1) \cdot ST^{\nu}(z_4)\\
=&\beta z_1 z_4+\gamma z_2 z_5+\alpha z_3 z_6+\\
 &+\gamma \zeta^{\nu} z_1^2+\alpha \zeta^{9 \nu} z_2^2+\beta \zeta^{3 \nu} z_3^2
  -\gamma \zeta^{12 \nu} z_4^2-\alpha \zeta^{4 \nu} z_5^2-\beta \zeta^{10 \nu}
  z_6^2+\\
 &+(\alpha-\beta) \zeta^{5 \nu} z_1 z_2+(\beta-\gamma) \zeta^{6 \nu} z_2 z_3
  +(\gamma-\alpha) \zeta^{2 \nu} z_1 z_3+\\
 &+(\beta-\alpha) \zeta^{8 \nu} z_4 z_5+(\gamma-\beta) \zeta^{7 \nu} z_5 z_6
  +(\alpha-\gamma) \zeta^{11 \nu} z_4 z_6+\\
 &-(\alpha+\beta) \zeta^{\nu} z_3 z_4-(\beta+\gamma) \zeta^{9 \nu} z_1 z_5
  -(\gamma+\alpha) \zeta^{3 \nu} z_2 z_6+\\
 &-(\alpha+\beta) \zeta^{12 \nu} z_1 z_6-(\beta+\gamma) \zeta^{4 \nu} z_2 z_4
  -(\gamma+\alpha) \zeta^{10 \nu} z_3 z_5.
\endaligned$$
$$\aligned
  &13 ST^{\nu}(z_2) \cdot ST^{\nu}(z_5)\\
=&\gamma z_1 z_4+\alpha z_2 z_5+\beta z_3 z_6+\\
 &+\alpha \zeta^{\nu} z_1^2+\beta \zeta^{9 \nu} z_2^2+\gamma \zeta^{3 \nu} z_3^2
  -\alpha \zeta^{12 \nu} z_4^2-\beta \zeta^{4 \nu} z_5^2-\gamma \zeta^{10 \nu}
  z_6^2+\\
 &+(\beta-\gamma) \zeta^{5 \nu} z_1 z_2+(\gamma-\alpha) \zeta^{6 \nu} z_2 z_3
  +(\alpha-\beta) \zeta^{2 \nu} z_1 z_3+\\
 &+(\gamma-\beta) \zeta^{8 \nu} z_4 z_5+(\alpha-\gamma) \zeta^{7 \nu} z_5 z_6
  +(\beta-\alpha) \zeta^{11 \nu} z_4 z_6+\\
 &-(\beta+\gamma) \zeta^{\nu} z_3 z_4-(\gamma+\alpha) \zeta^{9 \nu} z_1 z_5
  -(\alpha+\beta) \zeta^{3 \nu} z_2 z_6+\\
 &-(\beta+\gamma) \zeta^{12 \nu} z_1 z_6-(\gamma+\alpha) \zeta^{4 \nu} z_2 z_4
  -(\alpha+\beta) \zeta^{10 \nu} z_3 z_5.
\endaligned$$
$$\aligned
  &13 ST^{\nu}(z_3) \cdot ST^{\nu}(z_6)\\
=&\alpha z_1 z_4+\beta z_2 z_5+\gamma z_3 z_6+\\
 &+\beta \zeta^{\nu} z_1^2+\gamma \zeta^{9 \nu} z_2^2+\alpha \zeta^{3 \nu} z_3^2
  -\beta \zeta^{12 \nu} z_4^2-\gamma \zeta^{4 \nu} z_5^2-\alpha \zeta^{10 \nu}
  z_6^2+\\
 &+(\gamma-\alpha) \zeta^{5 \nu} z_1 z_2+(\alpha-\beta) \zeta^{6 \nu} z_2 z_3
  +(\beta-\gamma) \zeta^{2 \nu} z_1 z_3+\\
 &+(\alpha-\gamma) \zeta^{8 \nu} z_4 z_5+(\beta-\alpha) \zeta^{7 \nu} z_5 z_6
  +(\gamma-\beta) \zeta^{11 \nu} z_4 z_6+\\
 &-(\gamma+\alpha) \zeta^{\nu} z_3 z_4-(\alpha+\beta) \zeta^{9 \nu} z_1 z_5
  -(\beta+\gamma) \zeta^{3 \nu} z_2 z_6+\\
 &-(\gamma+\alpha) \zeta^{12 \nu} z_1 z_6-(\alpha+\beta) \zeta^{4 \nu} z_2 z_4
  -(\beta+\gamma) \zeta^{10 \nu} z_3 z_5.
\endaligned$$
Note that $\alpha+\beta+\gamma=\sqrt{13}$, we find that
$$\aligned
  &\sqrt{13} \left[ST^{\nu}(z_1) \cdot ST^{\nu}(z_4)+
   ST^{\nu}(z_2) \cdot ST^{\nu}(z_5)+ST^{\nu}(z_3) \cdot ST^{\nu}(z_6)\right]\\
 =&(z_1 z_4+z_2 z_5+z_3 z_6)+(\zeta^{\nu} z_1^2+\zeta^{9 \nu} z_2^2+\zeta^{3 \nu}
   z_3^2)-(\zeta^{12 \nu} z_4^2+\zeta^{4 \nu} z_5^2+\zeta^{10 \nu} z_6^2)+\\
  &-2(\zeta^{\nu} z_3 z_4+\zeta^{9 \nu} z_1 z_5+\zeta^{3 \nu} z_2 z_6)
   -2(\zeta^{12 \nu} z_1 z_6+\zeta^{4 \nu} z_2 z_4+\zeta^{10 \nu} z_3 z_5).
\endaligned$$
Let
$$\varphi_{\infty}(z_1, z_2, z_3, z_4, z_5, z_6)=\sqrt{13} (z_1 z_4+z_2 z_5+z_3 z_6)
  \eqno{(4.3.1)}$$
and
$$\varphi_{\nu}(z_1, z_2, z_3, z_4, z_5, z_6)=\varphi_{\infty}(ST^{\nu}(z_1, z_2,
                                              z_3, z_4, z_5, z_6))\eqno{(4.3.2)}$$
for $\nu=0, 1, \ldots, 12$. Then
$$\aligned
  \varphi_{\nu}
=&(z_1 z_4+z_2 z_5+z_3 z_6)+\zeta^{\nu} (z_1^2-2 z_3 z_4)+\zeta^{4 \nu}
  (-z_5^2-2 z_2 z_4)+\\
 &+\zeta^{9 \nu} (z_2^2-2 z_1 z_5)+\zeta^{3 \nu} (z_3^2-2 z_2 z_6)+
   \zeta^{12 \nu} (-z_4^2-2 z_1 z_6)+\\
 &+\zeta^{10 \nu} (-z_6^2-2 z_3 z_5).
\endaligned\eqno{(4.3.3)}$$
This leads us to define the following senary quadratic forms
(quadratic forms in six variables):
$$\left\{\aligned
  \mathbf{A}_0 &=z_1 z_4+z_2 z_5+z_3 z_6,\\
  \mathbf{A}_1 &=z_1^2-2 z_3 z_4,\\
  \mathbf{A}_2 &=-z_5^2-2 z_2 z_4,\\
  \mathbf{A}_3 &=z_2^2-2 z_1 z_5,\\
  \mathbf{A}_4 &=z_3^2-2 z_2 z_6,\\
  \mathbf{A}_5 &=-z_4^2-2 z_1 z_6,\\
  \mathbf{A}_6 &=-z_6^2-2 z_3 z_5.
\endaligned\right.\eqno{(4.3.4)}$$
Hence,
$$\sqrt{13} ST^{\nu}(\mathbf{A}_0)=\mathbf{A}_0+\zeta^{\nu} \mathbf{A}_1
  +\zeta^{4 \nu} \mathbf{A}_2+\zeta^{9 \nu} \mathbf{A}_3+\zeta^{3 \nu}
  \mathbf{A}_4+\zeta^{12 \nu} \mathbf{A}_5+\zeta^{10 \nu} \mathbf{A}_6.
  \eqno{(4.3.5)}$$
Let $H:=Q^5 P^2 \cdot P^2 Q^6 P^8 \cdot Q^5 P^2 \cdot P^3 Q$ where
$P=S T^{-1} S$ and $Q=S T^3$. Then (see \cite{Y2}, p.27)
$$H=\begin{pmatrix}
  0 &  0 &  0 & 0 & 0 & 1\\
  0 &  0 &  0 & 1 & 0 & 0\\
  0 &  0 &  0 & 0 & 1 & 0\\
  0 &  0 & -1 & 0 & 0 & 0\\
 -1 &  0 &  0 & 0 & 0 & 0\\
  0 & -1 &  0 & 0 & 0 & 0
\end{pmatrix}.\eqno{(4.3.6)}$$
Note that $H^6=1$ and $H^{-1} T H=-T^4$. Thus,
$\langle H, T \rangle \cong \mathbb{Z}_{13} \rtimes \mathbb{Z}_6$.
Hence, it is a maximal subgroup of order $78$ of $\text{PSL}(2, \mathbb{F}_{13})$
with index $14$. We find that $\varphi_{\infty}^2$ is invariant under
the action of the maximal subgroup $\langle H, T \rangle$. Note that
$$\varphi_{\infty}=\sqrt{13} \mathbf{A}_0, \quad
  \varphi_{\nu}=\mathbf{A}_0+\zeta^{\nu} \mathbf{A}_1+
  \zeta^{4 \nu} \mathbf{A}_2+\zeta^{9 \nu} \mathbf{A}_3+
  \zeta^{3 \nu} \mathbf{A}_4+\zeta^{12 \nu} \mathbf{A}_5+
  \zeta^{10 \nu} \mathbf{A}_6\eqno{(4.3.7)}$$
for $\nu=0, 1, \ldots, 12$. Let $w=\varphi^2$,
$w_{\infty}=\varphi_{\infty}^2$ and $w_{\nu}=\varphi_{\nu}^2$.
Then $w_{\infty}$, $w_{\nu}$ for $\nu=0, \ldots, 12$ form an
algebraic equation of degree fourteen, which is just the Jacobian
multiplier equation of degree fourteen (see \cite{K}, pp.161-162),
whose roots are these $w_{\nu}$ and $w_{\infty}$:
$$w^{14}+a_1 w^{13}+\cdots+a_{13} w+a_{14}=0.\eqno{(4.3.8)}$$
In particular, the coefficients
$$-a_1=w_{\infty}+\sum_{\nu=0}^{12} w_{\nu}=
  26 ({\bf A}_0^2+{\bf A}_1 {\bf A}_5+{\bf A}_2 {\bf A}_3+{\bf A}_4 {\bf A}_6).
  \eqno{(4.3.9)}$$
This leads to an invariant quadric
$$\Psi_2:={\bf A}_0^2+{\bf A}_1 {\bf A}_5+{\bf A}_2 {\bf A}_3+{\bf A}_4 {\bf A}_6
    =2 \Phi_4(z_1, z_2, z_3, z_4, z_5, z_6),\eqno{(4.3.10)}$$
where
$$\begin{array}{rl}
  \Phi_4:=&(z_3 z_4^3+z_1 z_5^3+z_2 z_6^3)-(z_6 z_1^3+z_4 z_2^3+z_5 z_3^3)+\\
          &+3(z_1 z_2 z_4 z_5+z_2 z_3 z_5 z_6+z_3 z_1 z_6 z_4),
\end{array}\eqno{(4.3.11)}$$
Hence, the variety $\Phi_4=0$ is a quartic four-fold, which is invariant
under the action of the group $G$.

  Let $x_i(z)=\eta(z) a_i(z)$ $(1 \leq i \leq 6)$, where
$$\left\{\aligned
  a_1(z) &:=e^{-\frac{11 \pi i}{26}} \theta \begin{bmatrix}
            \frac{11}{13}\\ 1 \end{bmatrix}(0, 13z),\\
  a_2(z) &:=e^{-\frac{7 \pi i}{26}} \theta \begin{bmatrix}
            \frac{7}{13}\\ 1 \end{bmatrix}(0, 13z),\\
  a_3(z) &:=e^{-\frac{5 \pi i}{26}} \theta \begin{bmatrix}
            \frac{5}{13}\\ 1 \end{bmatrix}(0, 13z),\\
  a_4(z) &:=-e^{-\frac{3 \pi i}{26}} \theta \begin{bmatrix}
            \frac{3}{13}\\ 1 \end{bmatrix}(0, 13z),\\
  a_5(z) &:=e^{-\frac{9 \pi i}{26}} \theta \begin{bmatrix}
            \frac{9}{13}\\ 1 \end{bmatrix}(0, 13z),\\
  a_6(z) &:=e^{-\frac{\pi i}{26}} \theta \begin{bmatrix}
            \frac{1}{13}\\ 1 \end{bmatrix}(0, 13z)
\endaligned\right.\eqno{(4.3.12)}$$
are theta constants of order $13$ and
$\eta(z):=q^{\frac{1}{24}} \prod_{n=1}^{\infty} (1-q^n)$ with
$q=e^{2 \pi i z}$ is the Dedekind eta function which are all
defined in the upper-half plane
$\mathbb{H}=\{ z \in \mathbb{C}: \text{Im}(z)>0 \}$. When parameterized by
theta constants of order $13$, it is shown that $\Phi_4=0$ (see \cite{Y6},
Theorem 9.1), i.e., $a_1=0$. This shows that the Jacobian multiplier equation
of degree fourteen is different from the transformation equation of degree
$13$ for the $j$-function (see (4.3.19) and (4.3.21)).

  Besides the Jacobian multiplier-equation of degree fourteen we have
established, we will need another kind of modular equation. We have
$$\aligned
  &-13 \sqrt{13} ST^{\nu}(z_1) \cdot ST^{\nu}(z_2) \cdot ST^{\nu}(z_3)\\
 =&-r_4 (\zeta^{8 \nu} z_1^3+\zeta^{7 \nu} z_2^3+\zeta^{11 \nu} z_3^3)
   -r_2 (\zeta^{5 \nu} z_4^3+\zeta^{6 \nu} z_5^3+\zeta^{2 \nu} z_6^3)\\
  &-r_3 (\zeta^{12 \nu} z_1^2 z_2+\zeta^{4 \nu} z_2^2 z_3+\zeta^{10 \nu} z_3^2 z_1)
   -r_1 (\zeta^{\nu} z_4^2 z_5+\zeta^{9 \nu} z_5^2 z_6+\zeta^{3 \nu} z_6^2 z_4)\\
  &+2 r_1 (\zeta^{3 \nu} z_1 z_2^2+\zeta^{\nu} z_2 z_3^2+\zeta^{9 \nu} z_3 z_1^2)
   -2 r_3 (\zeta^{10 \nu} z_4 z_5^2+\zeta^{12 \nu} z_5 z_6^2+\zeta^{4 \nu} z_6 z_4^2)\\
  &+2 r_4 (\zeta^{7 \nu} z_1^2 z_4+\zeta^{11 \nu} z_2^2 z_5+\zeta^{8 \nu} z_3^2 z_6)
   -2 r_2 (\zeta^{6 \nu} z_1 z_4^2+\zeta^{2 \nu} z_2 z_5^2+\zeta^{5 \nu} z_3 z_6^2)+\\
  &+r_1 (\zeta^{3 \nu} z_1^2 z_5+\zeta^{\nu} z_2^2 z_6+\zeta^{9 \nu} z_3^2 z_4)
   +r_3 (\zeta^{10 \nu} z_2 z_4^2+\zeta^{12 \nu} z_3 z_5^2+\zeta^{4 \nu} z_1 z_6^2)+\\
  &+r_2 (\zeta^{6 \nu} z_1^2 z_6+\zeta^{2 \nu} z_2^2 z_4+\zeta^{5 \nu} z_3^2 z_5)
   +r_4 (\zeta^{7 \nu} z_3 z_4^2+\zeta^{11 \nu} z_1 z_5^2+\zeta^{8 \nu} z_2 z_6^2)+\\
  &+r_0 z_1 z_2 z_3+r_{\infty} z_4 z_5 z_6+\\
  &-r_4 (\zeta^{11 \nu} z_1 z_2 z_4+\zeta^{8 \nu} z_2 z_3 z_5+\zeta^{7 \nu} z_1 z_3 z_6)+\\
  &+r_2 (\zeta^{2 \nu} z_1 z_4 z_5+\zeta^{5 \nu} z_2 z_5 z_6+\zeta^{6 \nu} z_3 z_4 z_6)+\\
  &-3 r_4 (\zeta^{7 \nu} z_1 z_2 z_5+\zeta^{11 \nu} z_2 z_3 z_6+\zeta^{8 \nu} z_1 z_3 z_4)+\\
  &+3 r_2 (\zeta^{6 \nu} z_2 z_4 z_5+\zeta^{2 \nu} z_3 z_5 z_6+\zeta^{5 \nu} z_1 z_4 z_6)+\\
  &-r_3 (\zeta^{10 \nu} z_1 z_2 z_6+\zeta^{4 \nu} z_1 z_3 z_5+\zeta^{12 \nu} z_2 z_3 z_4)+\\
  &+r_1 (\zeta^{3 \nu} z_3 z_4 z_5+\zeta^{9 \nu} z_2 z_4 z_6+\zeta^{\nu} z_1 z_5 z_6),
\endaligned$$
where
$$r_0=2(\theta_1-\theta_3)-3(\theta_2-\theta_4), \quad
  r_{\infty}=2(\theta_4-\theta_2)-3(\theta_1-\theta_3),$$
$$r_1=\sqrt{-13-2 \sqrt{13}}, \quad r_2=\sqrt{\frac{-13+3 \sqrt{13}}{2}},$$
$$r_3=\sqrt{-13+2 \sqrt{13}}, \quad r_4=\sqrt{\frac{-13-3 \sqrt{13}}{2}}.$$
This leads us to define the following senary cubic forms (cubic forms in
six variables):
$$\left\{\aligned
  \mathbf{D}_0 &=z_1 z_2 z_3,\\
  \mathbf{D}_1 &=2 z_2 z_3^2+z_2^2 z_6-z_4^2 z_5+z_1 z_5 z_6,\\
  \mathbf{D}_2 &=-z_6^3+z_2^2 z_4-2 z_2 z_5^2+z_1 z_4 z_5+3 z_3 z_5 z_6,\\
  \mathbf{D}_3 &=2 z_1 z_2^2+z_1^2 z_5-z_4 z_6^2+z_3 z_4 z_5,\\
  \mathbf{D}_4 &=-z_2^2 z_3+z_1 z_6^2-2 z_4^2 z_6-z_1 z_3 z_5,\\
  \mathbf{D}_5 &=-z_4^3+z_3^2 z_5-2 z_3 z_6^2+z_2 z_5 z_6+3 z_1 z_4 z_6,\\
  \mathbf{D}_6 &=-z_5^3+z_1^2 z_6-2 z_1 z_4^2+z_3 z_4 z_6+3 z_2 z_4 z_5,\\
  \mathbf{D}_7 &=-z_2^3+z_3 z_4^2-z_1 z_3 z_6-3 z_1 z_2 z_5+2 z_1^2 z_4,\\
  \mathbf{D}_8 &=-z_1^3+z_2 z_6^2-z_2 z_3 z_5-3 z_1 z_3 z_4+2 z_3^2 z_6,\\
  \mathbf{D}_9 &=2 z_1^2 z_3+z_3^2 z_4-z_5^2 z_6+z_2 z_4 z_6,\\
  \mathbf{D}_{10} &=-z_1 z_3^2+z_2 z_4^2-2 z_4 z_5^2-z_1 z_2 z_6,\\
  \mathbf{D}_{11} &=-z_3^3+z_1 z_5^2-z_1 z_2 z_4-3 z_2 z_3 z_6+2 z_2^2 z_5,\\
  \mathbf{D}_{12} &=-z_1^2 z_2+z_3 z_5^2-2 z_5 z_6^2-z_2 z_3 z_4,\\
  \mathbf{D}_{\infty}&=z_4 z_5 z_6.
\endaligned\right.\eqno{(4.3.13)}$$
Then
$$\aligned
  &-13 \sqrt{13} ST^{\nu}(\mathbf{D}_0)\\
 =&r_0 \mathbf{D}_0+r_1 \zeta^{\nu} \mathbf{D}_1+
   r_2 \zeta^{2 \nu} \mathbf{D}_2+
   r_1 \zeta^{3 \nu} \mathbf{D}_3+r_3 \zeta^{4 \nu} \mathbf{D}_4+\\
  &+r_2 \zeta^{5 \nu} \mathbf{D}_5+r_2 \zeta^{6 \nu} \mathbf{D}_6+
   r_4 \zeta^{7 \nu} \mathbf{D}_7+r_4 \zeta^{8 \nu} \mathbf{D}_8+\\
  &+r_1 \zeta^{9 \nu} \mathbf{D}_9+r_3 \zeta^{10 \nu} \mathbf{D}_{10}
   +r_4 \zeta^{11 \nu} \mathbf{D}_{11}+r_3 \zeta^{12 \nu} \mathbf{D}_{12}
   +r_{\infty} \mathbf{D}_{\infty}.
\endaligned$$
$$\aligned
  &-13 \sqrt{13} ST^{\nu}(\mathbf{D}_{\infty})\\
 =&r_{\infty} \mathbf{D}_0-r_3 \zeta^{\nu} \mathbf{D}_1-
   r_4 \zeta^{2 \nu} \mathbf{D}_2-r_3 \zeta^{3 \nu} \mathbf{D}_3+
   r_1 \zeta^{4 \nu} \mathbf{D}_4+\\
  &-r_4 \zeta^{5 \nu} \mathbf{D}_5-r_4 \zeta^{6 \nu} \mathbf{D}_6+
   r_2 \zeta^{7 \nu} \mathbf{D}_7+r_2 \zeta^{8 \nu} \mathbf{D}_8+\\
  &-r_3 \zeta^{9 \nu} \mathbf{D}_9+r_1 \zeta^{10 \nu} \mathbf{D}_{10}+
   r_2 \zeta^{11 \nu} \mathbf{D}_{11}+r_1 \zeta^{12 \nu} \mathbf{D}_{12}-
   r_0 \mathbf{D}_{\infty}.
\endaligned$$

  Let
$$\delta_{\infty}(z_1, z_2, z_3, z_4, z_5, z_6)
 =13^2 (z_1^2 z_2^2 z_3^2+z_4^2 z_5^2 z_6^2)\eqno{(4.3.14)}$$
and
$$\delta_{\nu}(z_1, z_2, z_3, z_4, z_5, z_6)
 =\delta_{\infty}(ST^{\nu}(z_1, z_2, z_3, z_4, z_5, z_6))\eqno{(4.3.15)}$$
for $\nu=0, 1, \ldots, 12$. Then
$$\delta_{\nu}=13^2 ST^{\nu}(\mathbf{G}_0)
 =-13 \mathbf{G}_0+\zeta^{\nu} \mathbf{G}_1+\zeta^{2 \nu} \mathbf{G}_2+
  \cdots+\zeta^{12 \nu} \mathbf{G}_{12},\eqno{(4.3.16)}$$
where the senary sextic forms (i.e., sextic forms in six variables) are
given as follows:
$$\left\{\aligned
  \mathbf{G}_0 =&\mathbf{D}_0^2+\mathbf{D}_{\infty}^2,\\
  \mathbf{G}_1 =&-\mathbf{D}_7^2+2 \mathbf{D}_0 \mathbf{D}_1+10 \mathbf{D}_{\infty}
                 \mathbf{D}_1+2 \mathbf{D}_2 \mathbf{D}_{12}+\\
                &-2 \mathbf{D}_3 \mathbf{D}_{11}-4 \mathbf{D}_4 \mathbf{D}_{10}-2
                 \mathbf{D}_9 \mathbf{D}_5,\\
  \mathbf{G}_2 =&-2 \mathbf{D}_1^2-4 \mathbf{D}_0 \mathbf{D}_2+6 \mathbf{D}_{\infty}
                 \mathbf{D}_2-2 \mathbf{D}_4 \mathbf{D}_{11}+\\
                &+2 \mathbf{D}_5 \mathbf{D}_{10}-2 \mathbf{D}_6 \mathbf{D}_9-2
                 \mathbf{D}_7 \mathbf{D}_8,\\
  \mathbf{G}_3 =&-\mathbf{D}_8^2+2 \mathbf{D}_0 \mathbf{D}_3+10 \mathbf{D}_{\infty}
                 \mathbf{D}_3+2 \mathbf{D}_6 \mathbf{D}_{10}+\\
                &-2 \mathbf{D}_9 \mathbf{D}_7-4 \mathbf{D}_{12} \mathbf{D}_4-2
                 \mathbf{D}_1 \mathbf{D}_2,\\
  \mathbf{G}_4 =&-\mathbf{D}_2^2+10 \mathbf{D}_0 \mathbf{D}_4-2 \mathbf{D}_{\infty}
                 \mathbf{D}_4+2 \mathbf{D}_5 \mathbf{D}_{12}+\\
                &-2 \mathbf{D}_9 \mathbf{D}_8-4 \mathbf{D}_1 \mathbf{D}_3-2
                 \mathbf{D}_{10} \mathbf{D}_7,\\
  \mathbf{G}_5 =&-2 \mathbf{D}_9^2-4 \mathbf{D}_0 \mathbf{D}_5+6 \mathbf{D}_{\infty}
                 \mathbf{D}_5-2 \mathbf{D}_{10} \mathbf{D}_8+\\
                &+2 \mathbf{D}_6 \mathbf{D}_{12}-2 \mathbf{D}_2 \mathbf{D}_3-2
                 \mathbf{D}_{11} \mathbf{D}_7,\\
  \mathbf{G}_6 =&-2 \mathbf{D}_3^2-4 \mathbf{D}_0 \mathbf{D}_6+6 \mathbf{D}_{\infty}
                 \mathbf{D}_6-2 \mathbf{D}_{12} \mathbf{}_7+\\
                &+2 \mathbf{D}_2 \mathbf{D}_4-2 \mathbf{D}_5 \mathbf{D}_1-2
                 \mathbf{D}_8 \mathbf{D}_{11},\\
  \mathbf{G}_7 =&-2 \mathbf{D}_{10}^2+6 \mathbf{D}_0 \mathbf{D}_7+4 \mathbf{D}_{\infty}
                 \mathbf{D}_7-2 \mathbf{D}_1 \mathbf{D}_6+\\
                &-2 \mathbf{D}_2 \mathbf{D}_5-2 \mathbf{D}_8 \mathbf{D}_{12}-2
                 \mathbf{D}_9 \mathbf{D}_{11},\\
  \mathbf{G}_8 =&-2 \mathbf{D}_4^2+6 \mathbf{D}_0 \mathbf{D}_8+4 \mathbf{D}_{\infty}
                 \mathbf{D}_8-2 \mathbf{D}_3 \mathbf{D}_5+\\
                &-2 \mathbf{D}_6 \mathbf{D}_2-2 \mathbf{D}_{11} \mathbf{D}_{10}-2
                 \mathbf{D}_1 \mathbf{D}_7,\\
  \mathbf{G}_9 =&-\mathbf{D}_{11}^2+2 \mathbf{D}_0 \mathbf{D}_9+10 \mathbf{D}_{\infty}
                 \mathbf{D}_9+2 \mathbf{D}_5 \mathbf{D}_4+\\
                &-2 \mathbf{D}_1 \mathbf{D}_8-4 \mathbf{D}_{10} \mathbf{D}_{12}-2
                 \mathbf{D}_3 \mathbf{D}_6,\\
  \mathbf{G}_{10} =&-\mathbf{D}_5^2+10 \mathbf{D}_0 \mathbf{D}_{10}-2 \mathbf{D}_{\infty}
                    \mathbf{D}_{10}+2 \mathbf{D}_6 \mathbf{D}_4+\\
                   &-2 \mathbf{D}_3 \mathbf{D}_7-4 \mathbf{D}_9 \mathbf{D}_1-2
                    \mathbf{D}_{12} \mathbf{D}_{11},\\
  \mathbf{G}_{11} =&-2 \mathbf{D}_{12}^2+6 \mathbf{D}_0 \mathbf{D}_{11}+4 \mathbf{D}_{\infty}
                    \mathbf{D}_{11}-2 \mathbf{D}_9 \mathbf{D}_2+\\
                   &-2 \mathbf{D}_5 \mathbf{D}_6-2 \mathbf{D}_7 \mathbf{D}_4-2
                    \mathbf{D}_3 \mathbf{D}_8,\\
  \mathbf{G}_{12} =&-\mathbf{D}_6^2+10 \mathbf{D}_0 \mathbf{D}_{12}-2 \mathbf{D}_{\infty}
                    \mathbf{D}_{12}+2 \mathbf{D}_2 \mathbf{D}_{10}+\\
                   &-2 \mathbf{D}_1 \mathbf{D}_{11}-4 \mathbf{D}_3 \mathbf{D}_9-2
                    \mathbf{D}_4 \mathbf{D}_8.
\endaligned\right.\eqno{(4.3.17)}$$
We have that $\mathbf{G}_0$ is invariant under the action of
$\langle H, T \rangle$, a maximal subgroup of order $78$ of
$\text{PSL}(2, \mathbb{F}_{13})$ with index $14$. Note that $\delta_{\infty}$,
$\delta_{\nu}$ for $\nu=0, \ldots, 12$ form an algebraic equation
of degree fourteen. However, we have
$$\Phi_6=\Phi_{0, 1}=\delta_{\infty}+\sum_{\nu=0}^{12} \delta_{\nu}=0.\eqno{(4.3.18)}$$
Hence, it is not the Jacobian multiplier equation of degree fourteen.

  Now, we have two distinct kinds of modular equations of degree fourteen:
the first one (4.3.8) is formed by $w_{\nu}$ $(\nu=0, 1, \ldots, 12)$ and
$w_{\infty}$, while the second one is formed by $\delta_{\nu}$ $(\nu=0, 1,
\ldots, 12)$ and $\delta_{\infty}$. By their construction, their Galois groups
are isomorphic to $\text{PSL}(2, \mathbb{F}_{13})$. Their coefficients are
invariant under the action of $\text{Aut}(\mathcal{L}(X(13)))$. Both of them
are defined over the locus $\mathcal{L}(X(13))$. In fact, in \cite{Y6}, Theorem 7.6,
we have verified explicitly that the locus $\mathcal{L}(X(13))$ can be parameterized
by the theta constants of order $13$.

  In his paper \cite{K1}, Klein gave a construction of the modular equation
of degree fourteen:
$$\aligned
  J: J-1:1 =&(\tau^2+5 \tau+13)(\tau^4+7 \tau^3+20 \tau^2+19 \tau+1)^3\\
           :&(\tau^2+6 \tau+13)(\tau^6+10 \tau^5+46 \tau^4+108 \tau^3
              +122 \tau^2+38 \tau-1)^2\\
           :&1728 \tau,
\endaligned\eqno{(4.3.19)}$$
whose Galois group is isomorphic with $\text{PSL}(2, \mathbb{F}_{13})$,
where for the quartic factor, one finds the decomposition
$$\left(\tau^2+\frac{7+\sqrt{13}}{2} \tau+\frac{11+3 \sqrt{13}}{2}\right)
  \left(\tau^2+\frac{7-\sqrt{13}}{2} \tau+\frac{11-3 \sqrt{13}}{2}\right),$$
and for the sextic factor, one finds the decomposition
$$\aligned
 &\left(\tau^3+5 \tau^2+\frac{21-\sqrt{13}}{2} \tau+\frac{3+\sqrt{13}}{2}\right)\\
 \times &\left(\tau^3+5 \tau^2+\frac{21+\sqrt{13}}{2} \tau+\frac{3-\sqrt{13}}{2}\right).
\endaligned$$
By the second equation in (4.3.19), we have the following expansion of $J-1$:
$$\aligned
  (J-1) \cdot 1728 \tau=&\tau^{14}+13 (2 \tau^{13}+25 \tau^{12}+196 \tau^{11}
                         +1064 \tau^{10}+4180 \tau^9\\
                        &+12086 \tau^8+25660 \tau^7+39182 \tau^6+41140 \tau^5
                         +27272 \tau^4\\
                        &+9604 \tau^3+1165 \tau^2)-982 \tau+13.
\endaligned\eqno{(4.3.20)}$$
Set
$$\tau=\frac{z}{\root 12 \of{\Delta}}$$
and $J=\frac{g_2^3}{\Delta}$, we have the following modular equation:
$$\aligned
  z^{14} &+13 (2 \Delta^{\frac{1}{12}} \cdot z^{13}
              +25 \Delta^{\frac{2}{12}} \cdot z^{12}
              +196 \Delta^{\frac{3}{12}} \cdot z^{11}
              +1064 \Delta^{\frac{4}{12}} \cdot z^{10}\\
             &+4180 \Delta^{\frac{5}{12}} \cdot z^9
              +12086 \Delta^{\frac{6}{12}} \cdot z^8
              +25660 \Delta^{\frac{7}{12}} \cdot z^7
              +39182 \Delta^{\frac{8}{12}} \cdot z^6\\
             &+41140 \Delta^{\frac{9}{12}} \cdot z^5
              +27272 \Delta^{\frac{10}{12}} \cdot z^4
              +9604 \Delta^{\frac{11}{12}} \cdot z^3
              +1165 \Delta \cdot z^2)\\
             &+[746 \Delta-(12 g_2)^3] \Delta^{\frac{1}{12}} \cdot z
              +13 \Delta^{\frac{14}{12}}=0.
\endaligned\eqno{(4.3.21)}$$
This equation appeared in \cite{Hurwitz1881}, p. 65.

  Therefore, on the side of algebraic equations, we have two distinct kinds
of equations with Galois group isomorphic to $\text{PSL}(2, \mathbb{F}_{13})$.
On the side of elliptic modular functions, we have the modular equations (4.3.19)
or (4.3.21) with Galois group isomorphic to $\text{PSL}(2, \mathbb{F}_{13})$.

\begin{center}
{\large\bf 5. An anabelian counterpart of the global Langlands correspondence
              for $\text{GL}(2, \mathbb{Q})$ by the cohomology of modular curves}
\end{center}

  A central theme in modern algebraic number theory is the connection
between modular forms and Galois representations. This link has its origins
in class field theory, and evolved from the work of Artin and Weil before
its scope was vastly expanded by Langlands in the late 1960s and 1970s,
linking the emerging theory of automorphic representations to Grothendieck's
conjectural theory of motives, where Galois representations appear as
realizations. In particular, Langlands recognized the connection between
non-abelian Galois representations and automorphic representations on reductive
algebraic groups, and told us essentially what to prove. More precisely,
according to \cite{De1976(2)}, the main goal of the global theory in
the non-abelian class field theory is to understand the relation between
motives defined over a number field, modular forms and their $L$-functions.
A motive $M$ over $\mathbb{Q}$ defines a compatible system of $\ell$-adic
representations
$\text{Gal}(\overline{\mathbb{Q}}/\mathbb{Q}) \rightarrow \text{GL}(M_{\ell})$:
the typical case is that in which $M$ is the $i$-th cohomology group
of a nonsingular projective variety $X$ defined over $\mathbb{Q}$, and
for each $\ell$, $\text{Gal}(\overline{\mathbb{Q}}/\mathbb{Q})$ acts
on $M_{\ell}=H^i(X, \mathbb{Q}_{\ell})$; $i$ is the weight of $M$. In
fact, the Galois representations associated with the \'{e}tale
cohomology of a variety defined over a number field is a central
subject of study in number theory. This gives rise to a function
$$L(M, s)=\prod_{p} (\det(1-F_p \cdot p^{-s}, M_{\ell}))^{-1}$$
(product over all primes, $F_p$ being the Frobenius endomorphism at
$p$; the product is absolutely convergent for $\text{Re}(s)$ large
enough; it is well defined except for a finite number of factors).
The only known results are proved by first linking $L(M, s)$ to a
modular form. Rather than modular forms, one should speak of the
irreducible components of the automorphic representation of
$G(\mathbb{A})$ ($\mathbb{A}$ being the adele ring of $\mathbb{Q}$)
in $L_{\chi}(G(\mathbb{A})/G(\mathbb{Q}))$ ($G$ reductive algebraic
group over $\mathbb{Q}$, $\chi$ a quasi-character of its center,
$L_{\chi}$ a sensible function space). Such a component $\pi$ defines
functions $L(\pi, \rho, s)$ ($\rho$ representation of the dual group
$\widehat{G}$). It is conjectured that every function $L(M, s)$ is
of type $L(\pi, \rho, s)$. It has two steps: attaching Galois
representations or motives to certain automorphic representations,
and conversely, showing that all Galois representations of motives
arise in this way. The converse direction incorporates the Fontaine-Mazur
conjecture as well, which posits a purely Galois-theoretic criterion
for when a Galois representation should arise from a motive.

  By Theorem 1.1 and Theorem 1.3, we give a distinct connection between
modular forms and Galois representations, linking $\mathbb{Q}(\zeta_p)$-rational
representations to defining ideals, where Galois representations
also appear as realizations. This implies that there is another
distinct approach to prove, which is different from the Langlands's
approach by automorphic representations. For simplicity, let $G$ be
a reductive group defined over $\mathbb{Q}$, $K$ a finite extension
of $\mathbb{Q}$, and $q$ a power of a prime number $p$. We have the
following two distinct connections between modular forms and
$\text{Gal}(\overline{\mathbb{Q}}/\mathbb{Q})$ according to the
comparison between motives and anabelian algebraic geometry.
$$\begin{array}{ccc}
  & \text{algebraic varieties} & \\
  \swarrow &  & \searrow\\
  \text{motives} & \longleftrightarrow & \text{defining ideals}\\
  \text{abelian} & \longleftrightarrow & \text{anabelian}\\
  \text{algebraic geometry} &          & \text{algebraic geometry}
\end{array}\eqno{(5.1)}$$

(1) The first kind of connection in the context of motives:

$$\begin{matrix}
  \text{modular}\\
  \text{forms}
  \end{matrix}
  \left\{\begin{matrix}
  \text{automorphic}  & & & &\text{$\ell$-adic}\\
  \text{representations} & \leftrightarrow &\text{motives} &
  \leftrightarrow &\text{representations}\\
  \text{of $G(\mathbb{A}_{\mathbb{Q}})$} & & & &\text{of}
  \end{matrix}\right\} \text{Gal}(\overline{\mathbb{Q}}/\mathbb{Q}).\eqno{(5.2)}$$
In particular, the reciprocity between motives and (arithmetic)
automorphic representations of $G(\mathbb{A}_{\mathbb{Q}})$ is the central
problem.

(2) The second kind of connection in the context of anabelian algebraic
geometry:

$$\begin{matrix}
  \text{modular}\\
  \text{forms}
  \end{matrix}
  \left\{\begin{matrix}
  \text{$K$-rational} & &\text{defining} & &\text{$K$-rational}\\
  \text{representations} & \leftrightarrow &\text{ideals} &
  \leftrightarrow &\text{representations}\\
  \text{of $G(\mathbb{F}_q)$} & & & &\text{of}
  \end{matrix}\right\} \text{Gal}(\overline{\mathbb{Q}}/\mathbb{Q}).\eqno{(5.3)}$$
In particular, the reciprocity between defining ideals (objects in anabelian
algebraic geometry) and $K$-rational representations of $G(\mathbb{F}_q)$ is
also a central problem.

  In fact, in the context of the usual Langlands reciprocity, $L$-functions
can classify the motives. However, in classical nonlinear Diophantine geometry,
we are interested in schemes, not motives, in particular, actual maps between
schemes. Hence, there is a need for nonlinear reciprocity of some sort. Our
correspondence (5.3) provides such a fully nonlinear theory.

  For the simplest group $\text{GL}(2, \mathbb{Q})$, the corresponding
simplest Shimura varieties are the elliptic modular curves, we give the
explicit realization of these two kinds of connections. Note that in
this case, the first kind of connection is classical, which is due
to Eichler-Shimura-Deligne-Serre-Langlands. On the other hand, in
this case, the second kind of connection is new. It is well-known that
the main tool to establish a Langlands correspondence is the construction
of geometric objects (e.g. Shimura varieties) that have so many symmetries
such that their cohomology carries a Galois action and an automorphic
action. And then one hopes to prove that this cohomology realizes the
Langlands correspondence. It turns out that the automorphic $L$-functions
are Mellin transforms of modular forms which represent the de Rham cohomology
of modular curves, while Galois representations embody their \'{e}tale
cohomology. It gives rise to a third comparison between \'{e}tale cohomology
and de Rham cohomology. However, linear algebraic and geometric methods used
in Langlands correspondence can not fully cover or substitute various profound
arithmetic issues, in particular those revealed by class field theory and its
other generalizations. Langlands correspondence is non-abelian but not utmost
non-abelian in comparison to anabelian algebraic geometry. This leads to our
correspondence. In our case, the geometric objects are still Shimura varieties
(modular curves over $\mathbb{Q}$). However, instead of their cohomology, we
use their defining ideals which carry both a Galois action and an action of
finite groups $\text{PSL}(2, \mathbb{F}_p)$. In particular, this defining
ideal realizes our correspondence. That is, there are two distinct approaches
to realize the non-abelian class field theory for $\text{GL}(2, \mathbb{Q})$:
one is the global Langlands correspondence, the other is our correspondence.
Hence, for $\text{GL}(2, \mathbb{Q})$, we have unified two generalizations of
class field theory: Langlands correspondence and anabelian algebraic geometry
by modular curves. In fact, for the same Shimura varieties, Langlands made the
cohomology of Shimura varieties a central part of his program, both as a source
of $\ell$-adic representations of Galois groups and as tests for his conjecture
that all motivic $L$-functions are automorphic. On the other hand, in our case,
we make the defining ideals of Shimura varieties a central part, both as a
source of $\mathbb{Q}(\zeta_p)$-rational representations of Galois groups and
the prediction that these Galois representations come from modular forms. It is
expected that the $K$-rational Galois representations
$$\rho: \text{Gal}(\overline{\mathbb{Q}}/\mathbb{Q}) \rightarrow G(\mathbb{F}_q)$$
can be realized by the defining ideals $I(S_K(G, X))$ of the corresponding
Shimura varieties $S_K(G, X)$ (which is dependent on $q$):
$$\rho: \text{Gal}(\overline{\mathbb{Q}}/\mathbb{Q}) \rightarrow
  \text{Aut}(I(S_K(G, X)))$$
as well as their surjective and modular realizations. As a consequence,
this gives an answer to the following question (see \cite{Fesenko} and
\cite{Fesenko2021}):

\textbf{Question 5.1.} Can the conjectures in arithmetic Langlands correspondence
be fully established remaining solely inside the use of representation theory
for adelic objects and Galois groups and class field theory? Can non-linear
methods help new fundamental developments in Langlands correspondence? Or
should one use more information about the absolute Galois group of global
and local fields, which can not be reached via representation theory, e.g.
such as in anabelian geometry?

  In fact, whereas the Langlands correspondence provides some conjectural
answers using representation theory, while anabelian geometry provides
very different insights into the full structure of the absolute Galois
group. In particular, our correspondence is in the context of anabelian
algebraic geometry. It does not use representation theory for adelic
objects, i.e., it can not be reached via automorphic representation theory
for $\text{GL}(2, \mathbb{A}_{\mathbb{Q}})$. Note that since Langlands
correspondence is a linear theory, it inevitably misses various important
features of the absolute Galois group that are not of linear representation
type, whereas class field theory is not, the asked enhancements of Langlands
correspondence have to be nonlinear (see \cite{Fesenko} and \cite{Fesenko2021}).
Nonlinear methods will help new fundamental developments in class field theory
for $\text{GL}(2, \mathbb{Q})$. Anabelian algebraic geometry is a sort of utmost
nonabelian and nonlinear theory, working with the absolute Galois groups and
fundamental groups of hyperbolic curves (smooth projective geometrically connected
curve whose Euler characteristic is negative).

  Therefore, for the group $\text{GL}(2, \mathbb{Q})$ and the associated
elliptic modular curves, there are two different approaches. The first
method is to embed the modular curve, for example, $X_0(p)$ with $p$ a
prime, into its Jacobian variety $J_0(p)$, in which the consideration of
explicit equations can largely be avoided. Recall that Abel's extension
of the addition theorem for elliptic integrals to the general case of
Abelian integrals led Jacobi to formulate the inversion problem for
hyperelliptic integrals. Through the works of Riemann and Clebsch, it
became the base of a new method to study the geometry of algebraic
curves. Abel's method to prove this theorem contains in germ the notions
of divisors and of linear families of divisors on an algebraic curves
and Riemann's interpretation of Abel's result leads to the notion of
Jacobian variety of an algebraic curve. It was Weil, in his proof of
the Riemann hypothesis for curves over finite fields (see \cite{Weil1948}),
who first made essential use of the passage from curves to abelian varieties
to derive important consequences for the arithmetic of curves. Also, Weil
(see \cite{Weil1948}) used his theory of the Jacobian to express the reciprocity
law for unramified class field theory in the form of an isomorphism between
the Galois group of the maximal abelian geometric unramified covering of a
smooth curve $X$ over a finite field and the group of rational points of the
Jacobian of $X$. The geometric insight, that in pursuing questions about curves
it sometimes pays to appeal to their Jacobians leads to a close relationship
between curves and their Jacobians. In particular, modular curves are endowed
with a plentiful supply of algebraic correspondences over $\mathbb{Q}$, which
emerge naturally from their moduli description and are geometric
incarnations of Hecke operators. The resulting endomorphisms break up
$J_0(p)$ into arithmetically simpler pieces with a large endomorphism
algebra, whose Tate modules give rise to compatible systems of
two-dimensional $\ell$-adic representations of
$\text{Gal}(\overline{\mathbb{Q}}/\mathbb{Q})$. These abelian variety
quotients of $\text{GL}(2)$ type provide a testing ground for the
general program of understanding $\ell$-adic or mod $\ell$ representations
of the Galois groups of number fields, which is a cornerstone of the
Langlands program. The two-dimensional representations of
$\text{Gal}(\overline{\mathbb{Q}}/\mathbb{Q})$ represent a prototypical
first step in this program, which go beyond the abelian setting of
global class field theory. The second method is to use the locus of
modular curves, for example, $\mathcal{L}(X(p))$, and its defining ideals
$I(\mathcal{L}(X(p)))$ which consist of a system of explicit algebraic
equations defined over $\mathbb{Q}$. This defining ideal $I(\mathcal{L}(X(p)))$
has the action of the automorphism group which is isomorphic to the
special projective linear group $\text{PSL}(2, \mathbb{F}_p)$. This
automorphism group breaks up $I(\mathcal{L}(X(p)))$ into arithmetically
simpler pieces with the same automorphism group, i.e., as the
intersection of some $\text{PSL}(2, \mathbb{F}_p)$-invariant ideals,
the corresponding algebraic cycles give rise to a family of
$\mathbb{Q}(\zeta_p)$-rational representations of
$\text{Gal}(\overline{\mathbb{Q}}/\mathbb{Q})$. These algebraic
cycles together with their automorphism groups provide a testing
ground of understanding $\mathbb{Q}(\zeta_p)$-rational
representations of the Galois groups of number fields. Note that
these representations are usually not two-dimensional, but instead
higher-dimensional. These higher-dimensional representations of
$\text{Gal}(\overline{\mathbb{Q}}/\mathbb{Q})$ represent a different
realization of global class field theory for $\text{GL}(2, \mathbb{Q})$.

  In particular, we have the following comparison between motives and
anabelian algebraic geometry:

$$\begin{array}{ccc}
  \text{motives} & \longleftrightarrow & \text{anabelian algebraic geometry}\\
  \text{Various cohomology theory} &  & \text{Defining ideals}\\
  \text{on modular curves with} & \longleftrightarrow & \text{of modular curves with}\\
  \text{$\text{Gal}(\overline{\mathbb{Q}}/\mathbb{Q}) \times
  \text{GL}(2, \mathbb{A}_{\mathbb{Q}})$-action} & &
  \text{$\text{Gal}(\overline{\mathbb{Q}}/\mathbb{Q}) \times
  \text{PSL}(2, \mathbb{F}_p)$-action}\\
  \text{Galois representations} &  & \text{Galois representations}\\
  \text{coming from division points} & \longleftrightarrow & \text{coming from}\\
  \text{on Jacobian varieties of} & & \text{defining ideals of}\\
  \text{modular curves} & & \text{modular curves}\\
  \text{Cuspidal representations} & \longleftrightarrow &
  \text{$\mathbb{Q}(\zeta_p)$-rational representations}\\
  \text{of $\text{GL}(2, \mathbb{A}_{\mathbb{Q}})$} &  &
  \text{of $\text{PSL}(2, \mathbb{F}_p)$}
\end{array}\eqno{(5.4)}$$

  Now, let us study the relation between our result with Langlands
correspondence. For $\text{SL}(2)$ and its variants $\text{GL}(2)$,
there are: global Langlands correspondence for
$\text{GL}(2, \mathbb{A}_F)$, where $F=\mathbb{Q}$ or $F(X)$, local
Langlands correspondence for $\text{GL}(2, F)$, where $F=\mathbb{R}$,
$\mathbb{C}$ or $\mathbb{Q}_p$, and mod $p$ Langlands correspondence
for $\text{GL}(2, \mathbb{F}_p)$ (Serre's modularity conjecture),
where $F=\mathbb{F}_p$, respectively. According to the above comparison
(5.4) between motives and anabelian algebraic geometry, the above three
kinds of Langlands correspondences come from motives. In contrast with
them, our result give a distinct correspondence which comes from
anabelian algebraic geometry. Because for the irreducible representations
of $\text{PSL}(2, \mathbb{F}_p)$, whose realizations can be formulated
in these three different scenarios (1.1). In particular, for its
$\mathbb{Q}(\zeta_p)$-rational realization, we have
$$\begin{matrix}
  \text{our correspondence: a distinct realization of}\\
  \text{the non-abelian class field theory for $\text{GL}(2, \mathbb{Q})$}\\
  \text{(in the context of anabelian algebraic geometry)}
\end{matrix}$$
$$\begin{array}{ccc}
  \text{defining ideals} & \longleftrightarrow & \text{$\mathbb{Q}(\zeta_p)$-rational}\\
  \text{$I(\mathcal{L}(X(p)))$} &  & \text{representations}\\
                                &  & \text{of $\text{PSL}(2, \mathbb{F}_p)$}\\
  \searrow &              & \swarrow\\
           & \text{$\mathbb{Q}(\zeta_p)$-rational Galois} &\\
           & \text{representations} &\\
           & \rho_p: \text{Gal}(\overline{\mathbb{Q}}/\mathbb{Q})
             \rightarrow \text{Aut}(\mathcal{L}(X(p)))\\
           & \updownarrow &\\
           & \text{their modular and} &\\
           & \text{surjective realizations} &
\end{array}\eqno{(5.5)}$$
Its formulation can be considered as a counterpart of the above global Langlands
correspondence over number fields and global Langlands correspondence over
function fields of curves over $\mathbb{F}_q$ as well as Serre's modularity
conjecture (\cite{Se1987}) (mod $p$ Langlands correspondence). In particular,
our correspondence (5.5) can be realized by the defining ideals of modular
curves. More precisely, the correspondence (5.5) are given by the Theorem 1.1
and Theorem 1.3.

  In general, the Langlands conjectures predict that automorphic forms
for a global field $F$ admit a spectral decomposition. Each piece of this
decomposition should be a motive, i.e., a part of the cohomology of an
algebraic variety over $F$. This is far from being proved in general. In
the case of number fields, cohomology of Shimura varieties give the answer,
but only under strong additional conditions. In the case of function fields,
varieties of Drinfeld shtukas give the answer in general, at the level of
cohomology, but not yet in a motivic way.

  In his paper \cite{Lafforgue1} (see also \cite{Lafforgue2}), Lafforgue
constructed a Langlands decomposition of the space of cuspidal automorphic
functions for function fields. More precisely, he constructed a canonical
decomposition of $\overline{\mathbb{Q}}_{\ell}$-vector spaces
$$C_c^{\text{cusp}}(\text{Bun}_G(\mathbb{F}_q)/\Xi, \overline{\mathbb{Q}}_{\ell})
 =\bigoplus_{\sigma} \mathfrak{H}_{\sigma},$$
where the direct sum is taken over global Langlands parameters
$\sigma: \pi_1(X, \overline{\eta})$ $\rightarrow$
$\widehat{G}(\overline{\mathbb{Q}}_{\ell})$. This decomposition is respected
by and compatible with the action of Hecke operators. Recently, it is shown
that (see \cite{Raskin}) a version of the geometric Langlands conjectures
yields a description of the eigenspaces of Lafforgue's decomposition in the
everywhere unramified case. In other words, there are no unramified cusp
forms with trivial Langlands parameters. In particular, it is shown that
the spectral decomposition actually follows from an $\ell$-adic version
of the geometric Langlands conjecture. It is known that Grothendieck's
motives over a given field form a $\overline{\mathbb{Q}}$-linear category
and unify the $\ell$-adic cohomologies of varieties over this field for
different $\ell$: a motive is a factor in a universal cohomology of a
variety. We consider here motives over $F$. In \cite{Lafforgue1} and
\cite{Lafforgue2} it is conjectured that the decomposition constructed above
is defined over $\overline{\mathbb{Q}}$ instead of $\overline{\mathbb{Q}}_{\ell}$,
indexed by motivic Langlands parameters $\sigma$, and independent on $\ell$.
The notion of motivic Langlands parameter is clear if we admit the standard
conjectures. A motivic Langlands parameter defined over $\overline{\mathbb{Q}}$
would give rise to a compatible family of morphisms $\sigma_{\ell, \iota}:
\text{Gal}(\overline{F}/F) \longrightarrow \widehat{G}(\overline{\mathbb{Q}}_{\ell})$
for any $\ell$ not dividing $q$ and any embedding $\iota: \overline{\mathbb{Q}}
\hookrightarrow \overline{\mathbb{Q}}_{\ell}$.

  In contrast with the above Lafforgue's decomposition for function fields
which is in the context of Grothendieck's motives, our decomposition (1.2)
in Theorem 1.1 can be defined over $\overline{\mathbb{Q}}$, which is
independent of $p$ in the context of anabelian algebraic geometry. In
particular, a Galois representation defined over $\mathbb{Q}(\zeta_p)$:
$$\rho_p: \text{Gal}(\overline{\mathbb{Q}}/\mathbb{Q}) \longrightarrow
  \text{Aut}(\mathcal{L}(X(p)))$$
will give rise to a Galois representation defined over $\overline{\mathbb{Q}}$
which is independent of $p$ by the embedding $\iota: \mathbb{Q}(\zeta_p)
\hookrightarrow \overline{\mathbb{Q}}$. This leads to the third comparison
between the second column and the third column of our Rosetta stone (1.1).
In the second column, the $\ell$-adic cohomology of the Drinfeld curves
gives the $\ell$-adic representations of $\text{SL}(2, \mathbb{F}_p)$ for
any $\ell \neq p$, which is dependent on $\ell$. In contrast with it, a
representation of $\text{PSL}(2, \mathbb{F}_p)$ defined over
$\mathbb{Q}(\zeta_p)$:
$$\pi_p: \text{PSL}(2, \mathbb{F}_p) \longrightarrow \text{Aut}(\mathcal{L}(X(p)))$$
will give a representation defined over $\overline{\mathbb{Q}}$ which is
independent of $p$ by the embedding $\iota: \mathbb{Q}(\zeta_p) \hookrightarrow
\overline{\mathbb{Q}}$.

  Note that Lafforgue's decomposition is respected by and compatible
with the action of Hecke operators. Our decomposition is respected
by and compatible with the action of the finite simple group
$\text{PSL}(2, \mathbb{F}_p)$. In our decomposition (1.2) in Theorem
1.1, the global Langlands parameters are replaced by the
$\mathbb{Q}(\zeta_p)$-rational Galois representations. In particular,
besides discrete series (cuspidal) representations, other types of
representations, such as degenerate principal series representations,
Steinberg representations, and even trivial representations do appear
in our decomposition (1.2).

  The central theme of the Langlands correspondence relates the spectra
of reductive algebraic groups over local and global fields and
parameterizes their representations by Galois theoretic data. It gives
rise to the global Langlands correspondence, local Langlands correspondences,
both Archimedean and non-Archimedean, which includes $\ell$-adic
($\ell \neq p$) and $p$-adic Langlands correspondences, mod $\ell \neq p$
and mod $p$ Langlands correspondences. Indeed, the global Langlands
conjectures over number fields predict that there exists a correspondence
among the cohomology of any algebraic variety over a number field,
automorphic representations and Galois representations, with the
simplest case being the cohomology of Shimura varieties.
$$\text{Global Langlands correspondence over number fields}$$
$$\text{(motives)}$$
$$\begin{array}{ccc}
  \text{cohomology of} & \longleftrightarrow & \text{automorphic}\\
  \text{algebraic varieties} &     &  \text{representations}\\
  \searrow &              & \swarrow\\
           & \text{Galois representations} &
\end{array}\eqno{(5.6)}$$
Here, the cohomology of algebraic varieties gives geometric realizations
of automorphic representations as well as Galois representations. In
particular, Shimura varieties are algebraic varieties defined over number
fields and equipped with many symmetries, which often provide a geometric
realization of the Langlands correspondence over number fields. One is a
Hecke symmetry coming from varying the level and considering various
transition morphisms between Shimura varieties at different levels. The
other is a Galois symmetry coming from the natural action of the Galois
groups on the \'{e}tale cohomology of Shimura varieties. It is a folklore
belief that the cohomology of Shimura varieties should decompose in terms
of certain automorphic representations and their Langlands parameters.

  In the modern theory of automorphic forms, one starts from a semisimple
algebraic group $G$ over $\mathbb{Q}$, and from an irreducible automorphic
representation $\pi$ of the adelic group $G(\mathbb{A})$. It is well known
that $\pi$ decomposes as an infinite tensor product of $p$-adic representations
$\pi_p$ (one for each prime $p$) and one archimedean component $\pi_{\infty}$,
which is a representation of $G(\mathbb{R})$. If the representation $\pi_{\infty}$
is algebraic (see \cite{Clozel} and \cite{BG}) , then an $\ell$-adic Galois
representation should correspond to $\pi$. This should be the case in particular
when $\pi_{\infty}$ is a discrete series representation or a limit of discrete
series. Let us consider the simplest case $G=\text{SL}(2)$. In the 1950s Eichler
(\cite{Eichler}) and Shimura (\cite{Shimura58}) showed that, to every classical
modular forms of weight two which is cuspidal and an eigenform for the action
of the Hecke operators, there corresponds a two-dimensional Galois representation,
in the sense that the Hecke eigenvalue at a prime is equal to the trace of a
Frobenius element at that prime. The work of Deligne and Serre in the 1960s
and 1970s (see \cite{De1968} and \cite{DeSe}) generalized the theorem of
Eichler and Shimura to cuspidal modular eigenforms of arbitrary weight $k \in
\mathbb{Z}$ and $k \geq 1$. A consequence of this construction of Deligne
(\cite{De1968}) is the proof of the Ramanujan conjecture:
$$|\tau(p)| \leq 2 p^{\frac{11}{2}}, \quad \text{where} \quad
  q \prod_{n=1}^{\infty} (1-q^n)^{24}=\sum_{n=1}^{\infty} \tau(n) q^n$$
via the \'{e}tale cohomology of modular curves. Deligne obtained the desired
bound as a consequence of his proof \cite{De1974} of the Weil conjectures
for smooth projective varieties over finite fields. In fact, two different
situations arise corresponding to the two possibilities $k \geq 2$ and $k=1$.
For weight $k \geq 2$, the construction of the Galois representation is
straightforward, based on a decomposition of the $\ell$-adic \'{e}tale cohomology
of modular curves under Hecke operators (see \cite{De1968}). In fact, the
realization of automorphic forms most directly related to Galois representations
is in \'{e}tale cohomology, because the latter carries an action of the absolute
Galois group as part of its definition. However, the \'{e}tale realization directly
applies only to the automorphic representation $\pi$ such that the archimedean
component $\pi_{\infty}$ is discrete series. In the $\text{SL}(2)$ case, only
modular forms of weight at least two can be realized in the \'{e}tale cohomology
of modular curves, which is why Deligne's construction only works in that case.
In particular, in the case of weight two, one has a triangle:
$$\begin{array}{ccc}
  \text{newforms of weight two} & \longrightarrow & \text{elliptic curves}\\
  \text{on $\Gamma_0(N)$ with rational} &     &\text{over $\mathbb{Q}$}\\
  \text{Fourier coefficients} &   &\text{up to isogeny}\\
  \searrow &              & \swarrow\\
           & \text{compatible systems} &   \\
           & \text{of $\ell$-adic representations} &
\end{array}\eqno{(5.7)}$$
In particular, it includes the Eichler-Shimura theory and the
Shimura-Tani-yama-Weil conjecture. This leads to the concept of regular algebraic
cuspidal automorphic representations. The condition that $\pi$ is regular
algebraic implies that the Hecke eigenvalues of a twist of $\pi$ appear in the
cohomology of the arithmetic locally symmetric spaces attached to the group
$\text{GL}(n, L)$, where $L$ is a number field. In particular, when $n=2$ and
$L=\mathbb{Q}$, this is just Deligne's construction. When $n>2$ or when $L$
is not totally real, the arithmetic locally symmetric spaces attached to the
group $\text{GL}(n, L)$ do not arise from Shimura varieties. However, the
self-duality condition implies that $\pi$ or one of its twists can be shown
to descend to another reductive (unitary or unitary similitude) group $G$
which does admit a Shimura variety. In the non-self-dual case, the automorphic
Galois representation will never occur in the Betti or \'{e}tale cohomology
of a Shimura variety. It can be constructed as a $\ell$-adic limit of
representations which do occur in such cohomology groups (see \cite{HLTT}
and \cite{Scholze}). By contrast, for weight $k=1$, one uses congruences between
forms of different weights to obtain a two-dimensional Artin representation, i.e.,
one with finite image (see \cite{DeSe}):
$$\begin{array}{ccc}
  \text{newforms of weight one} & \longrightarrow & \text{two-dimensional complex}\\
  \text{on $\Gamma_0(N)$ with character $\chi$} &     &\text{(Artin) representations}
\end{array}\eqno{(5.8)}$$
Although coherent cohomology lacks a Galois action, it may be used to construct
congruences between modular forms of weight one and modular forms of higher weight,
which allowed to reduce the construction of Galois representations in weight one
to that in higher weight. Generalizing this, Galois representations have been
constructed corresponding to Hecke eigenvalues appearing in the coherent cohomology
of Shimura varieties. This relies on finding congruences to Hecke eigenvalues for
regular algebraic automorphic representations. Note that the limitations of coherent
cohomology are twofold: first, coherent cohomology lacks a Galois action; second,
many $\pi$ that are conjectured to correspond to a Galois representation do not
appear in the coherent cohomology of any Shimura variety; all $\pi$ with $\pi_{\infty}$
a degenerate limit of discrete series are of this type, but there are still even
wider $\pi$, for example, the $\pi$ associated to classical Maass forms with
eigenvalue $\frac{1}{4}$. This difference between weight $\geq 2$ and weight
one is mirrored by the archimedean component of the corresponding automorphic
representation of $\text{SL}(2)$, which belongs to the discrete series in the
first case and is a limit of discrete series in the second case. In general,
for $G$ associated to a Hermitian symmetric domain and $\pi_{\infty}$ in the
discrete series, a construction analogous to that explained as above, but with
Shimura varieties replacing modular curves, should in principle work, again
according to the Langlands philosophy. According to \cite{Harris}, all of the
Galois representations associated to automorphic representations have been
constructed, either directly or by $\ell$-adic interpolation, using the
cohomology of Shimura varieties. This source of Galois representations has
been or soon will be exhausted, and new methods will need to be invented in
order to find the Galois representations attached to automorphic representations
that can not be related in any way to cohomology of Shimura varieties. In fact,
only a restricted class of Galois representations can be obtained using the cohomology
of Shimura varieties, and only those that can be realized directly in the cohomology
are associated to motives that admit an automorphic interpretation. In particular,
for those Galois representations which never be obtained in the cohomology of
Shimura varieties, although they are expected to be geometric no one has the
slightest idea where they might arise in the cohomology of algebraic varieties
(see \cite{Harris}, also \cite{Calegari}). For the mod $\ell$ representations,
it is impossible to prove the Galois representations are geometric using Shimura
varieties. It is not even clear what sense to give the question (see \cite{Harris2013}).
For example, the Galois representations associated with the cohomology of the
locally symmetric spaces attached to congruence subgroups of $\text{GL}(n)$
(including the Galois representations associated to torsion classes) do not
occur as a summand of the \'{e}tale cohomology of some smooth proper algebraic
variety over a number field (see \cite{AC}). Nevertheless, we can still define
cohomology groups associated to an arbitrary weight and level and also define
the Hecke algebras that act on these cohomology groups. However, for representations
$\Pi$ which do not act on the cohomology, in the end of his lecture \cite{Harris2013},
Harris pointed out:``for $\Pi$ not cohomological, a serious barrier; practically
no ideas''. In fact, from Matsushima's formula and computations of
$(\mathfrak{g}, K_{\infty})$-cohomology, we know that tempered automorphic
representations contributions in the cohomology of Shimura varieties with\
complex coefficients, is concentrated in middle degree. If one considers
cohomology with coefficients in a very regular local system, then only
tempered representations can contribute so that all of the cohomology
is concentrated in middle degree. On the other hand, for
$\overline{\mathbb{Z}}_{\ell}$-coefficients and Shimura varieties of
Kottwitz-Harris-Taylor types, whatever the weight of the coefficients
is, when the level is large enough at $\ell$, there are always nontrivial
torsion cohomology classes, so that the $\overline{\mathbb{F}}_{\ell}$-cohomology
can not be concentrated in middle degree.

  Recall that a Shimura datum is a pair $(G, X)$ where $G/\mathbb{Q}$
is a connected reductive group and $X \simeq G(\mathbb{R})/K_{\infty}$
is a Hermitian symmetric domain. Given such a pair $(G, X)$ and
$K \subset G(\mathbb{A}_f)$ a sufficiently small open compact subgroup,
the locally symmetric manifold
$$G(\mathbb{Q}) \backslash (X \times G(\mathbb{A}_f))/K \cong
  \text{Sh}(G, X)_K(\mathbb{C})$$
for a certain smooth quasi-projective algebraic variety
$\text{Sh}(G, X)_K$ defined over a number field $E=E(G, X)$. This
gives a tower $\{ \text{Sh}(G, X)_K \}_K$ with $G(\mathbb{A}_f)$-action.
In particular, take
$(G, X)=(\text{GL}(2), \mathbb{H}^{\pm}=\mathbb{C}-\mathbb{R})$.
Then $E=\mathbb{Q}$ and $Y_K=\text{Sh}(G, X)_K$ is the usual tower
of modular curves. When
$$K=K(N)=\{ g \in \text{GL}(2, \widehat{\mathbb{Z}}):
            g \equiv I (\text{mod $N$}) \},$$
$Y_{K(N)}$ is the moduli space of elliptic curves $E$ with a
trivialization $(\mathbb{Z}/N \mathbb{Z})^2 \simeq E[N]$. The
cohomology
$$\varinjlim_{K} H_{\text{\'{e}t}}^1(Y_{K, \overline{\mathbb{Q}}},
  \overline{\mathbb{Q}}_{\ell})$$
has a natural action of $\text{Gal}(\overline{\mathbb{Q}}/\mathbb{Q})
\times \text{GL}(2, \mathbb{A}_f)$. This action encodes the following
deep information, which gives the classical Langlands correspondence
for $\text{GL}(2, \mathbb{Q})$ in both the local and global contexts.
In particular, it provides the realization of the global Langlands
correspondence in the cohomology of modular curves.

\textbf{Theorem 5.2.} (Eichler-Shimura-Deligne-Ihara-Langlands-
Piatetskii-Shapiro-Carayol and others) {\it As}
$\text{Gal}(\overline{\mathbb{Q}}/\mathbb{Q}) \times \text{GL}(2,
\mathbb{A}_f)$-{\it representations, there is an isomorphism}
$$\varinjlim_{K} H_{\text{\'{e}t}}^1(Y_{K, \overline{\mathbb{Q}}},
  \overline{\mathbb{Q}}_{\ell})=\bigoplus_{f} \rho_f \otimes
  \bigotimes_p \pi_{f, p}+\cdots,\eqno{(5.9)}$$
{\it where the sum runs over all cuspidal automorphic representations
(satisfying a condition pertaining to $k$ at the infinite place), or
equivalently over all cuspidal newforms $f$ of weight $k$. For each
such $f$, $\rho_f$ is a two-dimensional representations of}
$\text{Gal}(\overline{\mathbb{Q}}/\mathbb{Q})$. {\it Moreover,}
$\pi_{f, p} \in \text{Irr}_{\overline{\mathbb{Q}}_{\ell}}
(\text{GL}(2, \mathbb{Q}_p))$ {\it matches} $\rho_f|_{W_{\mathbb{Q}_p}}$
{\it via the local Langlands correspondence. Here, ``\ldots'' are all
irreducible representations of} $\text{Gal}(\overline{\mathbb{Q}}/\mathbb{Q})$
{\it which occur there are one-dimensional. In particular, $\rho_f$ is
realized in the $\ell$-adic Tate module of elliptic curves.}

  In fact, Theorem 1.1 can be regarded as a counterpart of the above Theorem
5.2. In particular, (1.2) and (1.3) in Theorem 1.1 can be regarded as a
counterpart of (5.9). The \'{e}tale cohomology $H^1$ of modular curves
(motives) appearing in (5.9) is replaced by the defining ideals of modular
curves (anabelian algebraic geometry) in (1.2) and (1.3). In particular,
in Theorem 5.2, the Langlands correspondence for $\text{GL}(2, \mathbb{Q})$
gives the following:
$$\text{automorphic representations of $\text{GL}(2, \mathbb{A}_{\mathbb{Q}})$}
  \longleftrightarrow \text{Galois representations}.\eqno{(5.10)}$$
In contrast, our correspondence in Theorem 1.1 gives the following:
$$\text{$\mathbb{Q}(\zeta_p)$-rational representations of
  $\text{PSL}(2, \mathbb{F}_p)$} \longleftrightarrow
  \text{Galois representations}.\eqno{(5.11)}$$
In particular, our decomposition (1.2) in Theorem 1.1 can be regarded
as an anabelian counterpart of the decomposition of motives. That is,
the motives can be decomposed into irreducible parts which enable
us to analyze them piece by piece. The study of irreducible motives
should admit a simplicity analogous to the study of irreducible
representations. In our case, the defining ideals of $X(p)$ can
de decomposed into irreducible parts, which enable us to analyze
them piece by piece, and the study of irreducible parts admit a
simplicity leads to the study of irreducible representations of
$\text{PSL}(2, \mathbb{F}_p)$.

  Let us give some background (see \cite{Fargues-Mantovan}). For the
group $\text{GL}(2, \mathbb{Q})$, let
$f=\sum_{n=1}^{\infty} a_n q^n \in S_k(\Gamma_0(N))$ be a cuspidal
holomorphic modular form of weight $k \geq 2$, which are eigenfunctions
for Hecke operators outside $N$ and primitive. The form $f(z) dz^{k/2}$
defines a class in the de Rham cohomology with coefficients of the modular
curves $X_0(N)$ (Eichler-Shimura isomorphism). Applying the comparison
theorems between cohomological theories, we deduce a cohomology class
in the $\ell$-adic \'{e}tale cohomology with coefficients of $X_0(N)$.
This is the eigenfunction for Hecke operators outside $N$ and defines
a Galois representation
$$\rho_f: \text{Gal}(\overline{\mathbb{Q}}/\mathbb{Q}) \longrightarrow
  \text{GL}(2, \overline{\mathbb{Q}}_{\ell})$$
such that for every $p \nmid N \ell$, $\rho_f$ is unramified at $p$ and
$$\text{tr} \rho_f(\text{Frob}_p)=a_p,$$
where $\text{Frob}_p$ denotes a Frobenius substitution at $p$. A proof
of this equality consists in applying a Lefschetz trace formula to the
modulo $p$ reduction of smooth integral models of $X_0(N)$ over
$\mathbb{Z}_p$ and then comparing it with the Selberg trace formula.

  Associated with $f$ is an automorphic form
$$\phi \in L_{\text{cusp}}^2(\text{GL}(2, \mathbb{Q}) \backslash
 \text{GL}(2, \mathbb{A}_{\mathbb{Q}})/K_0(N))$$
whose iterates by the Hecke operators at all places give rise to an
automorphic representation
$$\pi=\pi_{\infty} \otimes \bigotimes_p \pi_p$$
where the equality of traces above means that for every $p \nmid N \ell$,
the representation $\pi_p$ of $\text{GL}(2, \mathbb{Q}_p)$ is associated
with $\rho_{f|D_p}$ via the unramified Langlands correspondence (the
Hecke eigenvalues correspond to the Frobenius eigenvalues) where $D_p
\subset \text{Gal}(\overline{\mathbb{Q}}/\mathbb{Q})$ denotes the
decomposition group at $p$. In fact, for every $p \nmid \ell$,
$(\rho_{f|D_p})^{\text{Frob}_p\text{-ss}}$ and $\pi_p$ correspond via
the local Langlands correspondence for $\text{GL}(2, \mathbb{Q}_p)$.

  The question is then to understand geometrically why $\rho_{f|D_p}$
depends only on $\pi_p$. We would like to find a geometric realization
of the correspondence $\pi_p \mapsto \rho_{f|D_p}$ in local cohomology
spaces, analogous at $p$ to $X_0(N)$, which is compatible with the
global correspondence $\pi \mapsto \rho$.

  This is the case when for example $\pi_p$ is supercuspidal (i.e.,
$(\rho_{f|D_p})^{\text{ss}}$ is irreducible) since then $\pi_p \otimes
(\rho_{f|D_p})$ is realized in the \'{e}tale cohomology of the $p$-adic
rigid analytic tower obtained by considering the points specializing at
a supersingular point (the Lubin-Tate tower for $\text{GL}(2)$). Indeed,
on the formal completion along the ordinary locus of regular integral
models over $\mathbb{Z}_p$ of curves $X(N)$ defined by the Drinfeld level
structures, there is an exact sequence
$$0 \longrightarrow E[p^{\infty}]^{0} \longrightarrow E[p^{\infty}]
  \longrightarrow E[p^{\infty}]^{\text{\'{e}t}} \longrightarrow 0$$
where $E$ denotes the universal ordinary elliptic curve and $E[p^{\infty}]^0$
and $E[p^{\infty}]^{\text{\'{e}t}}$ are of height $1$. When $p^n | N$
the level structure of Drinfeld
$$\eta: (p^{-n} \mathbb{Z}/\mathbb{Z})^2 \longrightarrow E[p^n]$$
then defines a splitting of this formal completion indexed by subgroups
$M=\eta^{-1}(E[p^n]^0) \subset (p^{-n} \mathbb{Z}/\mathbb{Z})^2$ direct
factor of rank $1$. This demonstrates, modulo the problems associated
with the points, that the difference between the cohomology of $X(N)$
and that of its super-singular part is induced at $p$ from the Borel
subgroup $\left(\begin{matrix} * & *\\ 0 & * \end{matrix}\right)$ (the
stabilizer of the component indexed by $M=(p^{-n} \mathbb{Z}/\mathbb{Z})
\oplus (0)$) and therefore does not contribute to the super-cuspidal part.

  According to a theorem of Serre-Tate (see \cite{LST}), the Lubin-Tate
tower does not depend on the global objects $X_0(N)$ since it depends
only on the formal group law associated with a super-singular elliptic
curve over $\overline{\mathbb{F}}_p$.

  Let us recall some facts about Lubin-Tate tower. Let $K$ be a
non-archimedean local field with uniformizer $\pi$ and residue field
$k \cong \mathbb{F}_q$, and let $n \geq 1$. The Lubin-Tate tower is a
projective system of formal schemes $\mathcal{M}_m$ which parameterize
deformation with level $\pi^m$ structure (i.e., Drinfeld level-$m$-structure)
of a one-dimensional formal $\mathcal{O}_K$-module of height $n$ over
$\overline{\mathbb{F}}_q$. After extending scalars to a separable closure
of $K$, the Lubin-Tate tower admits an action of the triple product group
$\text{GL}(n, K) \times D^{\times} \times W_K$, where $D/K$ is the central
division algebra of invariant $1/n$, and $W_K$ is the Weil group of $K$.
In particular, the $\ell$-adic \'{e}tale cohomology of the Lubin-Tate
tower realizes both the Jacquet-Langlands correspondence between
$\text{GL}(n, K)$ and $D^{\times}$ and the local Langlands correspondence
between $\text{GL}(n, K)$ and $W_K$. When $n=1$, this statement reduces to
classical Lubin-Tate theory. For $n=2$, the result was proved by Deligne
(see \cite{De1973}) and Carayol (see \cite{Carayol}). Moreover, Carayol
gave a conjectural description of the general phenomenon under the name
non-abelian Lubin-Tate theory (see \cite{Carayol1990}). Non-abelian
Lubin-Tate theory was established for all $n$ by Boyer for $K$ of positive
characteristic (see \cite{Boyer}) and by Harris and Taylor for $p$-adic $K$
(see \cite{HT}). In both cases, the result is established by embedding $K$
into a global field and appealing to results from the cohomology theory of
Shimura varieties or Drinfeld modular varieties.

  Recall that over $\mathbb{Z}$, Deligne and Rapoport initiated the study
of semistable models for the modular curves $X_0(Np)$ (see \cite{DeRa}).
The affine curve $xy=p$ over $W(\overline{\mathbb{F}}_p)$ appeared in
their work as a local model for $X_0(p)$ near a mod $p$ super-singular
point. Recently, Weinstein produced semistable models for Lubin-Tate
curves at height $n=2$ by passing to the infinite $p^{\infty}$-level,
when they each have the structure of a perfectoid space (see
\cite{Weinstein}). Such a Lubin-Tate curve is the rigid space attached
to the $p$-adic completion of a modular curve at one of its mod $p$
super-singular points. Since the super-singular locus is the interesting
part of the special fiber of a modular curve, Weinstein's work essentially
gives semistable models for $X(Np^m)$. In particular, his affine models
include curves with equations $xy^q-x^q y=1$ and $y^q+y=x^{q+1}$ over
$\overline{\mathbb{F}}_q$, where $q$ is a power of $p \neq 2$. In fact,
these two equations determine isomorphic projective curves, that is, the
Drinfeld curves.

  This is closely related to the Kronecker congruence
$$(j-j^{\prime p})(j^{\prime}-j^p) \equiv \text{$0$ (mod $p$)}.$$
It gives an equation for the modular curve $X_0(p)$ which represents
the moduli problem $[\Gamma_0(p)]$ for elliptic curves over a perfect
field of characteristic $p$. This moduli problem associates to such an
elliptic curve its finite flat subgroup schemes of rank $p$. A choice
of such a subgroup scheme is equivalent to a choice of an isogeny from
the elliptic curve with a prescribed kernel. The $j$-invariants of the
source and target curves along this isogeny are parameterized by $j$
and $j^{\prime}$. In fact, the Kronecker congruence gives a local
description for $[\Gamma_0(p)]$ at a super-singular point. For large
primes $p$, the mod $p$ super-singular locus may consist of more than
one closed point. In this case, the modular curve $X_0(p)$ does not
have an equation in the simple form given as above. Only its completion
at a single super-singular point does. What we need is only a suitable
local integral equation for $X_0(p)$ completed at a single mod $p$
super-singular point. This completion of $X_0(p)$ is a Lubin-Tate curve,
which is a moduli space for formal groups. In fact, there is a relation
between the moduli of formal groups and the moduli of elliptic curves.
This is the Serre-Tate theorem mentioned as above (see \cite{LST}), which
states that $p$-adically, the deformation theory of an elliptic curve is
equivalent to the deformation theory of its $p$-divisible group. In
particular, the $p$-divisible group of a super-singular elliptic curve
is formal. Therefore, the local information given by the Kronecker
congruence and its integral lifts are important for understanding
deformations of formal groups of height $n=2$. On the other hand,
the Eichler-Shimura relation $T_p \equiv F+V$ (mod $p$) gives
another interpretation for the moduli problem $[\Gamma_0(p)]$ in
characteristic $p$. However, there are modular equations which
describe $X_0(p)$ as a curve over $\text{Spec}(\mathbb{Z})$. In
particular, the classical modular equations lift and globalize the
Kronecker congruence. This leads to the following comparison between
the local theory and the global theory associated with the modular
equations for $\Gamma_0(p)$:

$$\text{Modular equations for $\Gamma_0(p)$}$$
$$\begin{matrix}
  \text{Local theory:} & & \text{Global theory:}\\
  \text{Kronecker congruence,} &  & \text{Galois groups of modular equations,}\\
  \text{Eichler-Shimura relation,} & & \text{modularity of Galois representations}\\
  \text{Lubin-Serre-Tate theorem.} & & \text{$\rho_p: \text{Gal}(\overline{
  \mathbb{Q}}/\mathbb{Q}) \rightarrow \text{Aut}(\mathcal{L}(X(p)))$.}
\end{matrix}\eqno{(5.12)}$$

  In the general case (see \cite{Carayol1990}), let
$f=\sum_{n=1}^{\infty} a_n q^n \in S_k(N, \varepsilon)$ be a
primitive normalized cuspidal eigenform of weight $k \geq 2$,
level $N$ and character $\varepsilon$. Then a classical construction,
due to Eichler and Shimura in the weight two case and to Deligne (see
\cite{De1968}) in the general situation, associates to $f$ a system
$(\rho_{\lambda})$ of two-dimensional $\lambda$-adic Galois
representations:
$$\rho_{\lambda}: \text{Gal}(\overline{\mathbb{Q}}/\mathbb{Q})
  \longrightarrow \text{GL}(2, E_{\lambda}),$$
where $\lambda$ ranges over the set of primes of the number field $E$,
which is generated by the coefficients $a_n$ and the values of
$\varepsilon$. In classical terms, the relationship between $f$ and
$\rho_{\lambda}$ is as follows: if we denote by $\ell$ the residual
characteristic of $\lambda$, then $\rho_{\lambda}$ is un-ramified
outside $N \ell$; and for $p$ any prime number not dividing $N \ell$,
the trace and determinant of $\rho_{\lambda}$ on the Frobenius
$\text{Frob}_p$ are given by
$$\left\{\aligned
  \text{tr} \rho_{\lambda}(\text{Frob}_p) &=a_p,\\
  \det \rho_{\lambda}(\text{Frob}_p) &=\varepsilon(p) p^{k-1}.
\endaligned\right.$$

  The classical theory of modular forms was not suitable to formulate
a precise theorem describing the behaviour of the representations
$\rho_{\lambda}$ at bad primes: that means, for $p \neq \ell$
a divisor of $N$, giving a recipe to compute the restriction
$\rho_{\lambda, p}$ of $\rho_{\lambda}$ to the local Galois group
$\text{Gal}(\overline{\mathbb{Q}}_p/\mathbb{Q}_p)$. It was only
after Jacquet-Langlands' work that such a precise recipe was
elaborated: a modular form $f$ as above gives rise to an
automorphic representation $\pi=\otimes \pi_{v}$ of the group
$\text{GL}(2, \mathbb{A})$, and the theorem was that the local
restriction $\rho_{\lambda, p}$ should correspond to the local
factor $\pi_p$ via the (suitably normalized) local Langlands
correspondence. More precisely, according to Langlands (see
\cite{Langlands73}), Deligne (see \cite{De1973}) and Carayol
(see \cite{Carayol1990}), let $\pi_f$ be a cuspidal automorphic
representation of $\text{GL}(2, \mathbb{A})$ (satisfying some condition
at the infinite place). Then there is a number field $E$ and a strictly
compatible system $\{ \sigma^{\lambda}(\pi_f)\}_{\lambda}$ of continuous
$E_{\lambda}$-adic Galois representations of
$\text{Gal}(\overline{\mathbb{Q}}/\mathbb{Q})$ such that for all primes
$p>2$, $\lambda \nmid p$,
$$\rho_f|_{W_{\mathbb{Q}_p}}=\sigma^{\lambda}(\pi_f)|_{W_{\mathbb{Q}_p}}
  \cong \text{LLC}(\pi_{f, p}).$$
Here $\text{LLC}$ denotes the (suitably normalized) local Langlands
correspondence between irreducible representations of $\text{GL}(2, \mathbb{Q}_p)$
and representations of $W_{\mathbb{Q}_p}$ with coefficients in
$\overline{\mathbb{Q}}_{\ell}$, where $\ell$ is the residue
characteristic of $\lambda$.

  Recall that the newform $f$ is associated to an automorphic
representation $\pi_{\mathbb{A}}$ of $\text{GL}(2, \mathbb{A})$,
where $\mathbb{A}$ is the ad\`{e}le ring of $\mathbb{Q}$. Note that
we can write an admissible irreducible representation $\pi$ of
$\text{GL}(2, \mathbb{A})$ as a restricted tensor product $\pi=
\bigotimes_{\nu}^{\prime} \pi_{\nu}$ of local components. In
particular, let $\pi_p$ be the component of $\pi_{\mathbb{A}}$
at $p$, so that $\pi_p$ is an admissible representation of
$\text{GL}(2, \mathbb{Q}_p)$. This representation may be classified
as a principal series representation, a special representation, or
a cuspidal representation of $\text{GL}(2, \mathbb{Q}_p)$. The local
Langlands correspondence attaches to $\pi$ a $\lambda$-adic
representation $\rho_{\lambda, \pi}$ of the Weil group of $\mathbb{Q}_p$.
More precisely, the local Langlands correspondence gives the following
explicit dictionary between automorphic representations for
$\text{GL}(2, \mathbb{Q})$ and Galois representations (see
\cite{Langlands73}, \cite{De1973} or \cite{Carayol1990}):
$$\begin{matrix}
 &\text{principal series} &\longleftrightarrow
 &\quad \text{decomposed Galois}\\
 &\text{representations} & &\quad \text{representations}\\
 &\text{special representations} &\longleftrightarrow
 &\quad \text{special Galois representations}\\
 &\text{cuspidal (Weil)} &\longleftrightarrow
 &\quad \text{irreducible (induced)}\\
 &\text{representations} & &\quad \text{Galois representations}\\
\end{matrix}\eqno{(5.13)}$$

  In fact, as soon as $\pi_p$ is not a spherical principal series
representation, the Galois representation $\rho_{\lambda}$ occurs
in the cohomology of a modular curve which has bad reduction at
$p$, and the whole question amounts to computing that cohomology
group: it sits in an exact sequence, with on the left side the
cohomology group of the special fiber and on the right the
cohomology of vanishing cycles. One can explicitly describe the
set of points of the special fiber, together with the Hecke and
Galois actions. Using such a description, Langlands (see
\cite{Langlands73}) was able to compute the cohomology of the
special fiber and (comparing the Selberg and Lefschetz trace
formulas) to prove the above theorem in the case of principal
or special representations; more precisely, it turned out that
principal series representations occurred only in the special
fiber cohomology and cuspidal ones only in the vanishing cycles
cohomology, while special representations contributed to both.

  The case of cuspidal series was solved by Deligne (see \cite{De1973})
at least for $p \neq 2$. Note that in this case there was no way
to compute explicitly the vanishing cycles cohomology. His method
consisted in constructing a local representation of the product group
$\text{GL}(2, \mathbb{Q}_p) \times B_p^{*} \times W_{\mathbb{Q}_p}$,
where $B_p$ denotes the quaternion division algebra over $\mathbb{Q}_p$,
and $W_{\mathbb{Q}_p}$ the Weil group. Using this local representation
and its interplay with the global representation on the vanishing
cycles group of the modular curve, he was able to prove that the
local restriction $\rho_{\lambda, p}$ was expressible in terms of
the local component $\pi_p$ alone. As a consequence of this, when
the cuspidal factor $\pi_p$ is a Weil representation which is always
the case if $p \neq 2$, one can reduce oneself to the situation where
the automorphic representation $\pi$ itself is obtained from the global
Weil construction, and then the theorem is easy to prove.

  For the remaining case of so-called extraordinary cuspidal
representation of $\text{GL}(2, \mathbb{Q}_2)$, in order to prove
the existence of the local Langlands correspondence in that case,
Carayol (see \cite{Carayol}) studied the same question for Hilbert
modular forms. He found that a theory of bad reduction for Shimura
curves existed, similar to the one for modular curves and using
this theory, he was able to generalize the results of Deligne and
Langlands to the Hilbert case. In this more general context, it
became possible to use the base change arguments and finally to
prove the above theorem even in the case of extraordinary cuspidal
representations. Kutzko had proven in the interval the existence
in all cases of the local Langlands correspondence for $\text{GL}(2)$
(see \cite{Kutzko}).

  In particular, when we take as level structure $K$ a principal
congruence subgroup $K=K_{N p^n}$, $p \nmid N$, then $\rho_f=\sigma(\pi_f)$
is the $\text{Gal}(\overline{\mathbb{Q}}/\mathbb{Q})$-representation attached
to the cuspidal automorphic representation (with respect to the cusp form $f$)
$\pi_f$ in the sense of Deligne (see \cite{De1968}), i.e., for all
$p \nmid N \ell$, the representation is un-ramified, and we have the
characteristic polynomial of $\rho_f(\text{Frob}_p)$ is equal to
$$X^2-a_p X+\chi(p) \cdot p^{k-1},$$
where $a_p$ is the Fourier coefficient of the cusp form $f$ corresponding
to $\pi_f$, and $\chi$ is the nebentypus character attached to this cusp
form $f$.

  Now, for $\text{GL}(2, \mathbb{Q})$, the global Langlands correspondence
among the $\ell$-adic \'{e}tale cohomology of modular curves over
$\mathbb{Q}$, automorphic representations of $\text{GL}(2, \mathbb{Q})$
and $\ell$-adic representations is a typical example.
$$\begin{array}{ccc}
  \text{$\ell$-adic cohomology} & \longleftrightarrow & \text{automorphic}\\
  \text{of modular curves} &     &\text{representations}\\
 (\text{$\ell$-adic system}) &   &\text{for $\text{GL}(2, \mathbb{Q})$}\\
  \searrow &              & \swarrow\\
           & \text{$\ell$-adic representations} &
\end{array}\eqno{(5.14)}$$
Here, the $\ell$-adic cohomology of modular curves provides geometric
realizations of automorphic representations as well as $\ell$-adic
representations. More precisely, let $p$ be a prime number. Langlands,
Deligne and Carayol showed in \cite{Langlands73}, \cite{De1973}, \cite{Carayol}
and \cite{Carayol1990} that one can realize a part of the global Langlands
correspondence for $\text{GL}(2)$ in the \'{e}tale cohomology:
$$\varinjlim_{K} H_{\text{\'{e}t}}^1(Y(K) \times_{\mathbb{Q}}
  \overline{\mathbb{Q}}, \mathbb{Q}_p),$$
where the limit is taken on the compact open subgroups $K$ of $\text{GL}(2)$
of finite ad\`{e}les of $\mathbb{Q}$ and where $Y(K)$ is the modular (open)
curve of level $K$. More precisely, if $\mathcal{O}_E$ is the ring of integers
of a finite extension $E$ of $\mathbb{Q}_p$ containing the Hecke eigenvalues of
a parabolic modular eigenform of weight $2$, the space
$$\varinjlim_{K} H_{\text{\'{e}t}}^1(Y(K) \times_{\mathbb{Q}} \overline{\mathbb{Q}},
  \mathcal{O}_E) \otimes_{\mathcal{O}_E} E$$
contains $\rho_f \otimes _E \otimes_{\ell}^{\prime} \pi_{\ell}(\rho_f|_{\text{Gal}(\overline{\mathbb{Q}_{\ell}}/\mathbb{Q}_{\ell})})$
where $\rho_f$ is the $p$-adic representation of
$\text{Gal}(\overline{\mathbb{Q}}/\mathbb{Q})$ associated to $f$ and $\pi_{\ell}(\rho_f|_{\text{Gal}(\overline{\mathbb{Q}_{\ell}}/\mathbb{Q}_{\ell})})$
the smooth representation of $\text{GL}(2, \mathbb{Q}_{\ell})$ corresponding to the
local Galois representation
$\rho_f|_{\text{Gal}(\overline{\mathbb{Q}_{\ell}}/\mathbb{Q}_{\ell})}$
by the local Langlands correspondence. Moreover, $\pi_{\ell}(\rho_f|_{\text{Gal}(\overline{\mathbb{Q}_{\ell}}/\mathbb{Q}_{\ell})})$
determines (essentially)
$\rho_f|_{\text{Gal}(\overline{\mathbb{Q}_{\ell}}/\mathbb{Q}_{\ell})}$
if $\ell \neq p$. Thus, the local Langlands correspondence for each
$\text{GL}(2, \mathbb{Q}_{\ell})$ fits naturally into a global correspondence
which is realized on a cohomology space. This is called local-global compatibility.
When $\ell=p$, it is no longer true in general that the smooth representation
$\pi_p(\rho_f|_{\text{Gal}(\overline{\mathbb{Q}_p}/\mathbb{Q}_p)})$
determines $\rho_f|_{\text{Gal}(\overline{\mathbb{Q}_p}/\mathbb{Q}_p)}$.
The local $p$-adic correspondence for the group $\text{GL}(2, \mathbb{Q}_p)$
associates to the representation
$\rho_f|_{\text{Gal}(\overline{\mathbb{Q}_p}/\mathbb{Q}_p)}$ a $p$-adic Banach
space $B(\rho_f|_{\text{Gal}(\overline{\mathbb{Q}_p}/\mathbb{Q}_p)})$ with
continuous action of $\text{GL}(2, \mathbb{Q}_p)$ which determines it completely.
This local $p$-adic correspondence also fits into a global correspondence
realising on the space of complete \'{e}tale cohomology:
$$\varinjlim_{K^p} \left(\varinjlim_{K_p} H_{\text{\'{e}t}}^1(Y(K^p K_p)
  \times_{\mathbb{Q}} \overline{\mathbb{Q}}, \mathcal{O}_E)\right)^{\wedge}
  \otimes_{\mathcal{O}_E} E,$$
where the hat $\wedge$ denotes the $p$-adic completion and where the inductive
limits are taken respectively on the compact open subgroups $K^p$ (resp. $K_p$)
of $\text{GL}(2)$ of the finite ad\`{e}les outside $p$ (resp. of
$\text{GL}(2, \mathbb{Q}_p)$) (see \cite{Emer2011}).

  For $\text{GL}(2)$, the Langlands correspondence over function fields of curves
over $\mathbb{F}_q$ provides a correspondence among the $\ell$-adic \'{e}tale
cohomology of Drinfeld modular curves over $F(X)$, automorphic representations
of $\text{GL}(2, F(X))$ and $\ell$-adic representations.
$$\text{Global Langlands correspondence over function fields of curves over $\mathbb{F}_q$}$$
$$\text{(motives)}$$
$$\begin{array}{ccc}
  \text{$\ell$-adic cohomology of} & \longleftrightarrow & \text{automorphic}\\
  \text{Drinfeld modular curves} &     &\text{representations for}\\
 (\text{$\ell$-adic system}) &   &\text{$\text{GL}(2, F(X))$}\\
  \searrow &              & \swarrow\\
           & \text{$\ell$-adic representations} &
\end{array}\eqno{(5.15)}$$

  In fact, in his two papers \cite{Dr1} and \cite{Dr2}, Drinfeld introduced
analogues of Shimura varieties for $\text{GL}(n, F(X))$ over a function field
$F(X)$ of characteristic $p>0$. Decomposing their $\ell$-adic cohomology under
the action of the Hecke operators he constructed Galois representations of
$F(X)$. In particular, for $n=2$ he showed that the correspondence which to
an automorphic representation associates the Galois representation on its
eigenspace is, up to a Tate twist, a Langlands correspondence. This is
completely analogous to the classical case of modular curves over $\mathbb{Q}$,
i.e., the Shimura varieties associated to $\text{GL}(2, \mathbb{Q})$. The
corresponding varieties are Drinfeld modular curves defined over the function
field $F(X)$.

  In conclusion, given a variety $X$, there are the following two correspondences:
$$\aligned
  X &\mapsto H^{*}(X) \quad \text{(various cohomology, motives) (linearization of $X$)},\\
  X &\mapsto I(X) \qquad \text{(defining ideals) (nonlinear)}.
\endaligned\eqno{(5.16)}$$
Correspondingly, we have the following two kinds of realizations:
$$\begin{matrix}
 &\text{\{algebraic geometry \}} &\longrightarrow &\text{\{theory of motives\}}\\
 &X & \longmapsto & H^{*}(X)\\
 &\text{\{symmetries of $X$\}} &\longmapsto & \text{\{symmetries of $H^{*}(X)$\}}\\
 &\text{\{algebraic geometry\}} &\longrightarrow &\text{\{defining ideals\}}\\
 &X & \longmapsto & I(X)\\
 &\text{\{symmetries of $X$\}} &\longmapsto & \text{\{symmetries of $I(X)$\}}
\end{matrix}\eqno{(5.17)}$$
In particular, it is well-known that there is a correspondence between ideals
and varieties $V(I(X))=X$. According to Grothendieck-Langlands, there should
be a correspondence:
$$\text{automorphic representations} \longleftrightarrow \text{motives}.$$
Theorem 1.1 and Theorem 1.3 show that there is the following correspondence
$$\text{$\mathbb{Q}(\zeta_p)$-rational representations} \longleftrightarrow
  \text{defining ideals}.$$

  Reciprocity between automorphic forms and algebraic varieties (including
motives (various cohomology) and defining ideals (anabelian algebraic geometry))
can be realized in two distinct kinds of ways:
$$\begin{array}{ccc}
  & \text{reciprocity} & \\
  \swarrow & \quad \quad  & \searrow\\
  \text{automorphic forms} & \longleftrightarrow & \text{algebraic varieties}
\end{array}$$
$$\text{GL}(2, \mathbb{Q}) \quad \text{(modular curves)}$$
$$\text{curves $X$} \longleftrightarrow \text{uniformization theorem}
  \longleftrightarrow X=\pi_1(X) \backslash \mathbb{H}$$
$$\begin{matrix}
 &\text{$H_1(X)$ (Jacobian varieties $J(X)$)} &\longleftrightarrow &
  \text{fundamental groups $\pi_1(X)$}\\
 &\text{of modular curves} &    & \text{of modular curves}\\
 &\text{(abelianization of $\pi_1(X)$} &  & \text{(anabelianization of $H_1(X)$)}\\
\end{matrix}$$
$$\text{(topology)}\eqno{(5.18)}$$
$$\begin{matrix}
 &\text{$H^1(X)$ (various cohomology)} & \longleftrightarrow
 &\text{defining ideals $I(X)$}\\
 &\text{motives of modular curves} &  & \text{of modular curves}\\
 &\text{(linearization of $X$)} &  & \text{(nonlinear counterpart of $H^1(X)$)}
\end{matrix}$$
$$\text{(algebraic geometry)}\eqno{(5.19)}$$
$$\begin{matrix}
 &\text{automorphic representations} &\longleftrightarrow
 &\text{reducible $\mathbb{Q}(\zeta_p)$-rational}\\
 &\text{of $\text{GL}(2, \mathbb{A}_{\mathbb{Q}})$} &
 &\text{representations of $\text{PSL}(2, \mathbb{F}_p)$}\\
 &\text{(adelic version of $\text{GL}(2, \mathbb{Q})$)} &
 &\text{(mod $p$ version of $\text{GL}(2, \mathbb{Q})$)}
\end{matrix}$$
$$\text{(representation theory)}\eqno{(5.20)}$$
$$\begin{matrix}
 &\text{$\mathbb{A}_{\mathbb{Q}}$ (adelic version of $\mathbb{Q}$)}
 &\longleftrightarrow &\text{$\mathbb{Q}(\zeta_p)$ (algebraic number fields)}
\end{matrix}$$
$$\text{(arithmetic)}\eqno{(5.21)}$$

\begin{center}
{\large\bf 6. An anabelian counterpart of the Shimura-Taniyama-Weil conjecture}
\end{center}

\textbf{Definition 6.1.} (see \cite{Saito2013}, \S 2.9) An elliptic curve $E$ over
$\mathbb{Q}$ is called modular exactly of level $N$ if there exists a normalized
simultaneous eigen-newform $f=\sum_{m=1}^{\infty} a_m(f) q^m$ of level $N$
such that
$$p+1-\# E(\mathbb{F}_p)=a_p(f)\eqno{(6.1)}$$
for all prime numbers $p$ not dividing $N$.

\textbf{Definition 6.2.} (see \cite{Saito2013}, Definition 0.14) Let $E$ be
an elliptic curve over $\mathbb{Q}$. Suppose $\ell$ is a prime number such
that the group of $\ell$-torsion points $E[\ell]$ is irreducible as a
representation of the absolute Galois group
$G_{\mathbb{Q}}=\text{Gal}(\overline{\mathbb{Q}}/\mathbb{Q})$:
$$\overline{\rho}_{E, \ell}: G_{\mathbb{Q}} \rightarrow \text{GL}(2, \mathbb{F}_{\ell}).
  \eqno{(6.2)}$$
Then we say that $E[\ell]$ is modular of level $N$ if there exists a normalized
simultaneous eigen-cuspform $f=\sum_{m=1}^{\infty} a_m(f) q^m$ of level dividing
$N$ and a prime ideal $\lambda$ of the integer ring of
$\mathbb{Q}(f):=\mathbb{Q}(a_m(f), m \in \mathbb{N})$ containing $\ell$ such that
$$p+1-\# E(\mathbb{F}_p) \equiv a_p(f) \quad \text{(mod $\lambda$)}\eqno{(6.3)}$$
for all primes $p$ not dividing $N$.

\textbf{Lemma 6.3.} (see \cite{Saito2013}, Lemma 4.1) {\it Let $E$ be a semistable
elliptic curve, and let $\ell$ be an odd prime. Let $N$ be the conductor of $E$.
Suppose the} mod $\ell$ {\it representation} $\overline{\rho}_{E, \ell}: G_{\mathbb{Q}}
\rightarrow \text{GL}(2, \mathbb{F}_{\ell})$ {\it associated with the $\ell$-torsion
points of $E$ is irreducible and modular. Then, $E$ is modular exactly of level $N$.}

  By Lemma 6.3, it suffices to show that the mod $\ell$ representation
$\overline{\rho}_{E, \ell}$ is modular for at least one of $\ell=3$ or $5$.
In particular, the modularity of $\overline{\rho}_{E, 3}$ plays a crucial
role in Wiles' proof (see \cite{W}) of the Shimura-Taniyama-Weil conjecture.
In order to prove its modularity, we lift $\overline{\rho}_{E, 3}$ to a
complex representation to which the theorem of Langlands-Tunnell is applicable
(see \cite{Gelbart}). More precisely, we extend
$\overline{\rho}_{E, 3}: G_{\mathbb{Q}} \rightarrow \text{GL}(2, \mathbb{F}_3)$
to a complex representation $\sigma: G_{\mathbb{Q}} \rightarrow \text{GL}(2, \mathbb{C})$
by composing $\overline{\rho}_{E, 3}$ with a specific (injective) homomorphism
$$\Psi: \text{GL}(2, \mathbb{F}_3) \hookrightarrow \text{GL}(2, \mathbb{Z}[\sqrt{-2}])
  \subset \text{GL}(2, \mathbb{C})\eqno{(6.4)}$$
through the formulas
$$\Psi\left(\begin{matrix} -1 & 1\\ -1 & 0 \end{matrix}\right)
 =\left(\begin{matrix} -1 & 1\\ -1 & 0 \end{matrix}\right), \quad \text{and} \quad
  \Psi\left(\begin{matrix} 1 & -1\\ 1 & 1 \end{matrix}\right)
 =\left(\begin{matrix} 1 & -1\\ -\sqrt{-2} & -1+\sqrt{-2} \end{matrix}\right).$$
Here $\alpha=\left(\begin{matrix} -1 & 1\\ -1 & 0 \end{matrix}\right)$ and
$\beta=\left(\begin{matrix} 1 & -1\\ 1 & 1 \end{matrix}\right)$ are two convenient
generators of $\text{GL}(2, \mathbb{F}_3)$. Once it is checked that the above formulas
indeed preserve the required relations, it is immediately seen that the resulting
homomorphism
$$\Psi: \text{GL}(2, \mathbb{F}_3) \rightarrow \text{GL}(2, \mathbb{Z}[\sqrt{-2}])
  \subset \text{GL}(2, \mathbb{C})$$
is the identity upon reduction mod $(1+\sqrt{-2})$. In particular,
$$\text{Trace}(\Psi(g)) \equiv \text{Trace}(g) \quad \text{(mod $(1+\sqrt{-2})$)}$$
and
$$\det(\Psi(g)) \equiv \det(g) \quad \text{(mod $3=(1+\sqrt{-2})(1-\sqrt{-2})$)}.$$
Moreover, $\sigma=\Psi \circ \overline{\rho}_{E, 3}: G_{\mathbb{Q}} \rightarrow
\text{GL}(2, \mathbb{C})$ is odd, irreducible and solvable. It is important to
point out that the representation
$$\Psi: \text{GL}(2, \mathbb{F}_3) \rightarrow \text{GL}(2, \mathbb{C})\eqno{(6.5)}$$
is really just one of the three cuspidal representations of the group
$\text{GL}(2, \mathbb{F}_3)$.

\textbf{Proposition 6.4.} (see \cite{Saito2013}, Proposition 4.2) {\it Let $E$
be a semistable elliptic curve over $\mathbb{Q}$.}

(1) {\it If the} mod $3$ {\it representation $\overline{\rho}_{E, 3}$ determined
by the $3$-torsion points of $E$ is reducible, then the} mod $5$ {\it representation}
$\overline{\rho}_{E, 5}$ {\it is irreducible.}

(2) {\it If the} mod $5$ {\it representation} $\overline{\rho}_{E, 5}$ {\it is
irreducible, then there exists a semistable elliptic curve} $E^{\prime}$ {\it
such that the} mod $3$ {\it representation} $\overline{\rho}_{E^{\prime}, 3}$
{\it is irreducible and} $\overline{\rho}_{E, 5} \simeq \overline{\rho}_{E^{\prime}, 5}$.

\textbf{Theorem 6.5.} (Shimura-Taniyama-Weil conjecture) (see \cite{BCDT})
{\it If $E$ is an elliptic curve over $\mathbb{Q}$, then $E$ is modular.}

\textbf{Theorem 6.6.} (see \cite{BCDT}) {\it If} $\overline{\rho}:
\text{Gal}(\overline{\mathbb{Q}}/\mathbb{Q}) \rightarrow \text{GL}(2, \mathbb{F}_5)$
{\it is an irreducible continuous representation with cyclotomic determinant,
then $\overline{\rho}$ is modular.}

  The proof of the Shimura-Taniyama-Weil conjecture also needs deformations
of Galois representations (see \cite{W}), which uses (6.2) again but in a
slightly different viewpoint (see \cite{GouMa}). Let
$$\overline{\rho}: G_{\mathbb{Q}} \rightarrow \text{GL}(2, \mathbb{F}_p)\eqno{(6.6)}$$
be a continuous absolutely irreducible representation which is unramified
outside $p$ and infinity and which is modular of non-critical slope in the
sense that it is the residual representation attached to some modular eigenform
of weight $k$ and level $p$ with coefficients in $\mathbb{Z}_p$ whose slope
is not equal to $0$ or $k-1$. This restricts $p$ to be at least $11$. The basic
deformation problem which may be posed concerning $\overline{\rho}$ is the
following: classify all equivalence classes of liftings of $\overline{\rho}$
to representations
$$\rho: G_{\mathbb{Q}} \rightarrow \text{GL}(2, A),\eqno{(6.7)}$$
where $A$ is any commutative complete Noetherian local ring with residue
field $\mathbb{F}_p$ and where the representation $\rho$ is unramified
outside $p$ and infinity.

  As we have shown above, the representation theory of $\text{GL}(2, \mathbb{F}_p)$
plays a crucial role in the proof of the Shimura-Taniyama-Weil conjecture
(see \cite{W} and \cite{CDT}). In particular, the cuspidal representation
(6.5) of $\text{GL}(2, \mathbb{F}_3)$ appears in Wiles' proof \cite{W}. In
our case, much more representations do appear. Corollary 1.2 shows that,
besides cuspidal representations (i.e., discrete series), principal series
representations as well as Steinberg representations of
$\text{PSL}(2, \mathbb{F}_p)$ appear in our theory.
$$\begin{matrix}
 &\text{Wiles' proof} &\longleftrightarrow &\text{our theory}\\
 &\text{only} & \longleftrightarrow &\text{cuspidal representations,}\\
 &\text{cuspidal} & & \text{principal series, and}\\
 &\text{representation} & &\text{Steinberg representations}\\
 &\text{of $\text{GL}(2, \mathbb{F}_3)$} &
 &\text{of $\text{PSL}(2, \mathbb{F}_p)$ ($p \geq 7$)}
\end{matrix}\eqno{(6.8)}$$

  In contrast with the above proof of the Shimura-Taniyama-Weil conjecture
which deals with the cases $p=3$ and $p=5$, our Galois representation (2.4.20)
begins with $p \geq 7$. For the almost same Galois representation
$$\rho: G_{\mathbb{Q}} \rightarrow \text{GL}(2, \mathbb{F}_p) \quad \text{or} \quad
  \text{PSL}(2, \mathbb{F}_p),\eqno{(6.9)}$$
this gives two distinct geometric realizations: one comes from the $p$-torsion
points of elliptic curves, the other comes from the defining ideals of modular
curves $X(p)$. More precisely, we have the following comparison:
$$\begin{array}{ccc}
  \text{motives} & \longleftrightarrow & \text{anabelian algebraic geometry}\\
  \text{linear realization $E[p]$}  & \longleftrightarrow &
  \text{nonlinear realization $I(\mathcal{L}(X(p)))$}\\
  \text{Aut}(E[p]) \cong \text{GL}(2, \mathbb{F}_p) & \longleftrightarrow  &
  \text{Aut}(\mathcal{L}(X(p))) \cong \text{PSL}(2, \mathbb{F}_p)\\
\end{array}\eqno{(6.10)}$$
In particular, Theorem 1.3 shows that the modularity condition
(6.3) should be replaced by the following correspondence:
$$\begin{array}{ccc}
  \text{the Jacobian multiplier}  & \longleftrightarrow & \text{$p$-th order transformation}\\
  \text{equation of degree $p+1$} &                     & \text{equation of the $j$-function}\\
\end{array}\eqno{(6.11)}$$
which are connected by the fact that both of the two sides of (6.11) have
Galois groups isomorphic to $\text{PSL}(2, \mathbb{F}_p)$. In fact, (6.11)
corresponds to the exceptional case $\ell=p=N$ in the case of (6.3), which
is not included in (6.3). Hence, the two conditions (6.3) and (6.11) are
complementary to each other. This leads to the following two kinds of
correspondences. The first one comes from arithmetic geometry:
$$p+1-\# E(\mathbb{F}_p) \longleftrightarrow \text{the Jacobian multiplier
  equation of degree $p+1$}.\eqno{(6.12)}$$
The second one comes from modular forms:
$$a_p(f) \longleftrightarrow \text{$p$-th order transformation
  equation of the $j$-function}.\eqno{(6.13)}$$
The left-hand-side of (6.12) and (6.13) comes from the Shimura-Taniyama-Weil
conjecture (motives), while, the right-hand-side of (6.12) and (6.13) can be
regarded as an anabelian counterpart of the the Shimura-Taniyama-Weil
conjecture (anabelian algebraic geometry).

\begin{center}
{\large\bf 7. An anabelian counterpart of Artin's conjecture, Serre's modularity
              conjecture and the Fontaine-Mazur conjecture}
\end{center}

  In the theory of elliptic functions, given a prime number $p$, there
are two distinct theories: one is the division equation of order $p$,
the other is the transformation equation of order $p$. The first one
leads to the Galois theory of division equations, the second one leads
to the Galois theory of transformation equations. In particular, from
the viewpoint of Galois representations, the first one leads to the
Galois representations on the torsion points (division points) of
elliptic curves or abelian varieties, the second one leads to the
Galois representations arising from the defining ideals
$I(\mathcal{L}(X(p)))$ of modular curves $X(p)$ according to Theorem
1.1 and Theorem 1.3. This leads to the following comparison:
$$\begin{matrix}
 \text{division} & \rightarrow & \text{$p^{\infty}$-torsion $A[p^{\infty}]$}
 & \rightarrow & \text{$p$-adic $H_{\text{\'{e}t}}^1$, $H_{\text{dR}}^1$}
 & \text{(motives)}\\
 \text{equations} &  & \text{($p$-divisible groups)} &  & \text{(Grothendieck's} & \\
   & & \text{(Tate, Grothendieck)} &  & \text{mysterious} & \\
   & &  & & \text{functor)} & \\
 \text{division} & \rightarrow & \text{representations of}
 & \rightarrow & \text{$\ell$-adic $H_{\text{\'{e}t}}^1$, $H_{\text{dR}}^1$}
 & \text{(motives)}\\
 \text{equations} &  & \text{$\text{SL}(2, \mathbb{F}_p)$ or $\text{PSL}(2, \mathbb{F}_p)$}
 &  &  & \\
 & & \text{(Hecke, Drinfeld,} & & & \\
 & & \text{Deligne-Lusztig)} & & & \\
 \text{transfor-} & \rightarrow & \text{the locus $\mathcal{L}(X(p))$} &
 \rightarrow & \pi_1 & \text{(anabelian}\\
 \text{mation} & & \text{(ours)} &  & \text{defining ideals} & \text{algebraic} \\
 \text{equations} &  &  &  & I(\mathcal{L}(X(p))) & \text{geometry)}
\end{matrix}\eqno{(7.1)}$$

  Now, we give the first theory in the context of motives. Let $E$ be an
elliptic curve defined over a field $k$ and let $\ell$ be a prime different
from $\text{char}(k)$. The Tate modules
$T_{\ell}(E)=\varprojlim_{n} E[\ell^n](\overline{k})$ are special cases of
the $\ell$-adic homology groups associated to algebraic varieties. The
Galois group $\text{Gal}(\overline{k}/k)$ acts on $T_{\ell}(E)$ and on the
$\mathbb{Q}_{\ell}$-vector space $V_{\ell}(E)=\mathbb{Q}_{\ell} \otimes T_{\ell}(E)$.
We may consider the associated $\ell$-adic Galois representation
$\rho_{\ell}: \text{Gal}(\overline{k}/k) \rightarrow \text{Aut}(T_{\ell})
 \simeq \text{GL}(2, \mathbb{Z}_{\ell})$. The image $G_{\ell}$ of $\rho_{\ell}$
is an $\ell$-adic Lie subgroup of $\text{Aut}(T_{\ell})$. The Galois extension
associated to $G_{\ell}$ is obtained by adding the coordinates of the points
of $E(\overline{k})$ of order a power of $\ell$ to the field $k$. The main
theorem of \cite{Se1972} states that if $E$ is an elliptic curve defined over
an algebraic number field $k$, which does not have complex multiplication,
then $G_{\ell}=\text{Aut}(T_{\ell}(E))$, for almost all $\ell$. In particular,
we have $\text{Gal}(k(E_{\ell})/k) \simeq \text{GL}(2, \mathbb{F}_{\ell})$,
for almost all $\ell$. In general, let $A$ be an abelian variety of dimension
$d$ defined over a field $k$. Given a prime $\ell \neq \text{char}(k)$, the
Tate module $T_{\ell}(A)=\varprojlim_{n} A[\ell^n](\overline{k})$ is a free
$\mathbb{Z}_{\ell}$-module of rank $2d$. Let $V_{\ell}(A)=T_{\ell}(A) \otimes
\mathbb{Q}_{\ell}$. The action of the absolute Galois group of $k$ on the
Tate module of $A$ gives an $\ell$-adic Galois representation $\rho_{\ell}:
\text{Gal}(\overline{k}/k) \rightarrow \text{GL}(T_{\ell}(A))=\text{GL}(2d,
\mathbb{Z}_{\ell})$. Its image $G_{\ell}$ is a compact subgroup of
$\text{GL}(2d, \mathbb{Z}_{\ell})$, hence it is a Lie subgroup of the
$\ell$-adic Lie group $\text{GL}(T_{\ell})$.

  On the other hand, the second theory is in the context of anabelian
algebraic geometry. We give a distinct realization of the above Galois
representation
$$\rho_p: \text{Gal}(\overline{\mathbb{Q}}/\mathbb{Q}) \rightarrow
  \text{GL}(2, \mathbb{F}_p) \quad \text{or} \quad \text{PSL}(2,
  \mathbb{F}_p)\eqno{(7.2)}$$
by the following
$$\rho_p: \text{Gal}(\overline{\mathbb{Q}}/\mathbb{Q}) \rightarrow
  \text{Aut}(\mathcal{L}(X(p)))\eqno{(7.3)}$$
which comes from the transformation equations of elliptic functions with
order $p$. In particular, Corollary 1.2 and Theorem 1.3 show that the
cases of $p=7$, $11$ and $13$ are especially important from the viewpoint
of algebraic geometry, representation theory and arithmetic. At first, we
give their geometrical as well as representation theoretic aspect as we
have done in section 2.5. The arithmetical significance mainly come from
that Theorem 1.3 gives a connection between Galois representations and
Felix Klein's elliptic modular functions (see \cite{KF1} and \cite{KF2}),
i.e., $j$-functions by the defining ideals of modular curves $X(p)$. In
particular, we have two distinct kinds of viewpoints: one comes from Galois,
the other comes from Riemann. In the first viewpoint, the cases of $p=7$
and $11$ are especially interesting. In the second viewpoint, the cases
of $p=7$ and $13$ are especially interesting. At first, let us consider
the first viewpoint. Because these are the only two cases ($p=7$ and $11$)
that there are equations of degree $p+1$ and $p$ whose Galois groups are
isomorphic with $\text{Aut}(\mathcal{L}(X(p)))$, which are closely related
to the last mathematical testament of Galois \cite{Ga}. Moreover, these two
Galois groups are not only isomorphic, one comes from algebraic equations,
the other comes from elliptic modular functions, but also there is a direct
and explicit reduction process from the side of algebraic equations to the
side of elliptic modular functions. From the second viewpoint, these are the
only two cases ($p=7$ and $13$) that the genus of the corresponding curves
$X_0(p)$ has genus zero. Moreover, the case of $p=13$ is interesting because
that there exist an algebraic equation of degree fourteen which is not the
Jacobian multiplier equation of degree fourteen.

  Now, we give a comparison with the Fontaine-Mazur-Langlands conjecture. In
general, the Langlands reciprocity conjecture articulates a correspondence
between certain automorphic representations (those which are algebraic) of
$G(\mathbb{A}_F)$ for a connected reductive algebraic group $G$ over a number
field $F$, and motives over $F$ whose motivic Galois group is closely related
to the $C$-group ${}^{C}G$ of $G$ (see \cite{Langlands79} and \cite{BG}), where
the $C$-group is a refinement of the $L$-group introduced by Langlands (see
\cite{BG}). Conjectures on Galois representations (especially the Fontaine-Mazur
conjecture (see \cite{FM})) suggest that such motives in turn may be identified
with certain compatible systems of $\ell$-adic Galois representations
$\text{Gal}(\overline{F}/F) \rightarrow {}^{C}G(\overline{\mathbb{Q}}_{\ell})$.
This leads to the construction of $\ell$-adic Galois representations. It is
well-known that the $\ell$-adic \'{e}tale cohomology of algebraic varieties
is much richer than their classical cohomology in the sense that it admits
the action of Galois groups. This opens up a new field of inquiry, even in
the classical case. More precisely, given an algebraic variety $X$ over
$\mathbb{Q}$, Grothendieck's $\ell$-adic \'{e}tale cohomology
$H_{\text{\'{e}t}}^i(X, \mathbb{Q}_{\ell})$ attaches to such an $X$ a
collection of finite dimensional $\mathbb{Q}_{\ell}$-vector spaces with a
continuous action of $\text{Gal}(\overline{\mathbb{Q}}/\mathbb{Q})$
($\ell$-adic representations)
$$\rho_{\ell, i}: \text{Gal}(\overline{\mathbb{Q}}/\mathbb{Q}) \rightarrow
  \text{Aut}(H_{\text{\'{e}t}}^i(X, \mathbb{Q}_{\ell})).\eqno{(7.4)}$$
On the other hand, for modular curves $X(p)$, whose defining ideals of
the locus $\mathcal{L}(X(p))$ provide a collection of finite dimensional
$\mathbb{Q}$-vector spaces with a continuous action of
$\text{Gal}(\overline{\mathbb{Q}}/\mathbb{Q})$ ($\mathbb{Q}(\zeta_p)$-rational
representations)
$$\rho_p: \text{Gal}(\overline{\mathbb{Q}}/\mathbb{Q}) \rightarrow
  \text{Aut}(\mathcal{L}(X(p))) \quad (p \geq 7).\eqno{(7.5)}$$
In particular, both $H_{\text{\'{e}t}}^i(X, \mathbb{Q}_{\ell})$ and
$\mathcal{L}(X(p))$ come from geometry: one comes from motives, the
other comes from anabelian algebraic geometry. This leads to the
following comparison:
$$\begin{array}{ccc}
  & \text{prime $p$} & \\
  \swarrow &  & \searrow\\
  \text{division equation} & \longleftrightarrow & \text{transformation equation}\\
  \text{of order $p$} &  & \text{of order $p$}\\
  \downarrow &  & \downarrow\\
  \text{motives} & \longleftrightarrow & \text{anabelian geometry}\\
  \downarrow &  & \downarrow\\
  \text{Galois groups of} & \longleftrightarrow & \text{Galois groups of}\\
  \text{division equations} &  & \text{transformation equations}\\
  \searrow &  & \swarrow\\
  & \text{modular curves} & \\
  & \text{of level $p$} & \\
  \swarrow &  & \searrow\\
  X_0(p) & \longleftrightarrow & X(p)\\
  \downarrow &  & \downarrow\\
  \text{$H_1$ (abelian)} & \longleftrightarrow & \text{$\pi_1$ (anabelian)}\\
  \downarrow &  & \downarrow\\
  \text{Jacobian variety} & \longleftrightarrow & \text{defining ideals}\\
  J_0(p)=J(X_0(p)) &  & \text{$I(\mathcal{L}(X(p)))$}\\
  \downarrow &  & \downarrow\\
  \text{Galois action on} & \longleftrightarrow & \text{Galois action}\\
  \text{division points of $J_0(p)$} &  & \text{on $\mathcal{L}(X(p))$}\\
  \downarrow &  & \downarrow\\
  \text{$\text{End}(J_0(p)) \otimes \mathbb{Q}=\mathbf{T}$} & \longleftrightarrow &
  \text{$\text{Aut}(\mathcal{L}(X(p)))$}\\
  \text{(Hecke algebras)} &  & \text{(as representations)}\\
  \downarrow &  & \downarrow\\
  \mathbf{T}=K_1 \times \cdots \times K_t & \longleftrightarrow
  & \text{Aut}(\mathcal{L}(X(p)))=V_1 \oplus \cdots \oplus V_t\\
  \text{defined over} &  & \text{defined over}\\
  \text{totally real fields} & \longleftrightarrow & \text{cyclotomic fields}\\
  \text{or CM-fields} &  & \mathbb{Q}(\zeta_p)\\
  \downarrow &  & \downarrow\\
  J_0(p)=A_1 \times \cdots \times A_t & \longleftrightarrow
  & I(\mathcal{L}(X(p)))=I_1 \cap \cdots \cap I_t\\
  \downarrow &  & \downarrow\\
  \text{multiplicity-one} & \longleftrightarrow & \text{multiplicity-one}\\
  \text{phenomenon} & & \text{phenomenon}
\end{array}\eqno{(7.6)}$$

  At first, we give some interpretation on the left-hand side of (7.6).
To understand representations arising from modular forms, it is helpful
to realize these representations inside of geometric objects such as the
Jacobian variety $J$ of modular curves. These representations are constructed
geometrically with the help of the Hecke algebra $\mathbf{T}$. In particular,
the Jacobian $J_0(N)$ of $X_0(N)$ is an abelian variety defined over
$\mathbb{Q}$. It has a number of remarkable arithmetical properties. Many
of these are related to the structure of the endomorphism ring of $J_0(N)$
and its simple factors. Suppose that $N$ is prime. Then the endomorphism
algebra of $J_0(N)$ is the algebra $\mathbf{T}$ generated by the Hecke
operators $T_n$ with $n$ prime to $N$ (see \cite{Ribet1975}). More precisely,
the algebra $(\text{End} J_0(N)) \otimes \mathbb{Q}$ of all endomorphisms of
$J_0(N)$ is equal to the algebra generated by the Hecke operators, regarded
as $\mathbb{Q}$-endomorphisms of $J_0(N)$. In general, let $N$ be a positive
integer. Let $R$ be the subring of $\text{End}(J_0(N))$ generated by the
Hecke operators $T_n$ with $n$ prime to $N$. The theory of newforms shows
that $E=R \otimes \mathbb{Q}$ is a product of totally real algebraic number
fields $E_{\alpha}$:
$$E=E_1 \times \cdots \times E_t$$
and that the degree $[E: \mathbb{Q}]$ is the number of (normalized) newforms
of weight two, trivial character, and level dividing $N$. Note that the ring
$R$ operates faithfully on the abelian variety $A:=\prod_{M | N} J_0(M)_{\text{new}}$,
where $J_0(M)_{\text{new}}$ is the new subvariety of $J_0(M)$. The dimension
of $A$ is the degree $[E: \mathbb{Q}]$, and the decomposition of $E$ into the
product $\prod E_{\alpha}$:
$$E=E_1 \times \cdots \times E_t$$
decomposes $A$, up to isogeny, as a product
$$A=A_1 \times \cdots \times A_t$$
of abelian varieties $A_{\alpha}$ with real multiplication by the factors
$E_{\alpha}$. On the other hand, for the Nebentypus character, Ribet proved
the following result (see \cite{Ribet1975}). Let $N$ be a prime congruent
to $1$ mod $4$. Let $S$ be the vector space of cusp forms of weight two
and Nebentypus $\left(\frac{}{N}\right)$ on $\Gamma_0(N)$, where
$\left(\frac{}{N}\right)$ is the quadratic character mod $N$. Let $J$
be the abelian variety over $\mathbb{Q}$ associated with $S$. According
to \cite{Shimura1972}, this is an abelian variety of dimension
$d=\dim_{\mathbb{C}} S=g(X_0(N))$ (the genus of $X_0(N)$) whose complex
tangent space is $S$ and whose $\mathbb{Q}$-endomorphism algebra contains
the algebra $\mathbf{T}$ generated over $\mathbb{Q}$ by the Hecke operators
$T_n$ ($n$ prime to $N$) acting on $S$. The rank $[\mathbf{T}: \mathbb{Q}]$
is again equal to $d$. The algebra $\mathbf{T}$ breaks up into a product
$$K_1 \times \cdots \times K_t$$
of fields which are no longer totally real but instead are CM-fields.
This decomposition of $\mathbf{T}$ forces $J$ to split over $\mathbb{Q}$
(up to isogeny) as a product
$$A_1 \times \cdots \times A_t$$
with
$$\dim A_i=[K_i: \mathbb{Q}].$$
The significance of the above decompositions comes from the fact that the
structure of $J_0(N)$ to within isogeny is determined by the structure
of the space of cusp forms of weight two on $\Gamma_0(N)$ as of a module
over the Hecke algebra $\mathbf{T}$. More precisely, we have (see \cite{Sh},
p. 183, Theorem 7.14 or \cite{Vladut}, p. 272, Theorem 2.1.1):

\textbf{Theorem 7.1.} {\it Let $f=\sum_{n=1}^{\infty} a_n q^n$, $a_1=1$,
be a normalized cusp form of weight two on $\Gamma_0(N)$, which is
an eigenform for Hecke operators $\{ T_n \}$, $(n, N)=1$, and let
$K=\mathbb{Q}(\{ a_n \})$ be the subfield of $\overline{\mathbb{Q}}$
generated by its coefficients. Then there exists an abelian subvariety
$A_f \subset J_0(N)$ defined over $\mathbb{Q}$ and an isomorphism
$\theta: K \simeq \text{End}(A_f) \otimes \mathbb{Q}$ with the following
properties:

(1) $\dim A_f=[K: \mathbb{Q}]$;

(2)  for $(n, N)=1$ the element $\theta(a_n)$ coincides with the
restriction of $T_n \in \text{End}(J_0(N)) \otimes \mathbb{Q}$ to $A_f$.
these requirements uniquely determine the pair $(A_f, \theta)$.}

  The space of cusp forms of weight two is generated by eigenforms for
the Hecke operators, therefore the above theorem 6.1 gives a complete
description of the isogeny type of $J_0(N)$. Let
$\rho: \text{Gal}(\overline{\mathbb{Q}}/\mathbb{Q}) \rightarrow
\text{GL}(2, \overline{\mathbb{F}}_{\ell})$ be an irreducible Galois
representation that arises from a weight-two newform $f$. The next step,
after having attached a maximal ideal $\mathfrak{m}$ to $f$, is to find
a $\mathbf{T}/\mathfrak{m}$-vector space affording $\rho$ inside of the
group of $\ell$-torsion points of $J$. Following \cite{Mazur1977}, we
consider the $\mathbf{T}/\mathfrak{m}$-vector space
$$J[\mathfrak{m}]:=\{ \text{$P \in J(\overline{\mathbb{Q}})$: $tP=0$
  for all $t \in \mathfrak{m}$}\} \subset J(\overline{\mathbb{Q}})[\ell]
  \cong (\mathbb{Z}/\ell \mathbb{Z})^{2g}.$$
Since the endomorphisms in $\mathbf{T}$ are $\mathbb{Q}$-rational,
$J[\mathfrak{m}]$ comes equipped with a linear action of
$\text{Gal}(\overline{\mathbb{Q}}/\mathbb{Q})$. The fact that
$\text{tr}(\rho(\text{Frob}_p))$ and $\det(\rho(\text{Frob}_p))$ both
lie in the subfield $\mathbf{T}/\mathfrak{m}$ of $\overline{\mathbb{F}}_{\ell}$
suggests that $\rho$ has a model over $\mathbf{T}/\mathfrak{m}$, in the
sense that $\rho$ is equivalent to a representation taking values in
$\text{GL}(2, \mathbf{T}/\mathfrak{m}) \subset \text{GL}(2, \overline{\mathbb{F}}_{\ell})$
(see \cite{RS}).

\textbf{Proposition 7.2.} {\it The representation $\rho$ has a model
$\rho_{\mathfrak{m}}$ over the finite field $\mathbf{T}/\mathfrak{m}$.}

  In fact, we have the following stronger result: Let $N$ be a positive
integer. Let $\mathbf{T}=\mathbf{T}_N$ be the ring generated by the Hecke
operators $T_n$ $(n \geq 1)$ on the space of weight-two cusp forms on
$\Gamma_0(N)$. Let $\mathfrak{m}$ be a maximal ideal of $\mathbf{T}$. The
following result is a variant of \cite{DeSe}, Th\'{e}or\`{e}me 6.7, which
attaches a mod $\ell$ representation of
$\text{Gal}(\overline{\mathbb{Q}}/\mathbb{Q})$ to mod $\ell$ modular forms
which are eigenvectors for the Hecke operators $T_n$.

\textbf{Proposition 7.3.} (\cite{Ribet1990}, Proposition 5.1). {\it
There is a unique semisimple representation}
$$\rho_{\mathfrak{m}}: \text{Gal}(\overline{\mathbb{Q}}/\mathbb{Q})
  \rightarrow \text{GL}(2, \mathbf{T}/\mathfrak{m})$$
{\it satisfying}
$$\text{$\text{tr}(\rho_{\mathfrak{m}}(\text{Frob}_r))=T_r$ (mod $\mathfrak{m})$,
  \quad $\det(\rho_{\mathfrak{m}}(\text{Frob}_r))=r$ (mod $\mathfrak{m}$)}$$
{\it for almost all primes $r$. The representation $\rho_{\mathfrak{m}}$
is unramified at all primes $r$ prime to $\mathfrak{m} N$, and the
displayed relations hold for all such primes.} ({\it We say that
$\rho_{\mathfrak{m}}$ is the representation of}
$\text{Gal}(\overline{\mathbb{Q}}/\mathbb{Q})$ {\it attached to
$\mathfrak{m}$.})

  Let $f$ be an eigenform and fix a nonzero prime $\lambda$ of the ring
generated by the Fourier coefficients of $f$ such that $\rho_{f, \lambda}$
is absolutely irreducible. View the Hecke algebra $\mathbf{T}$ as a subring
of $\text{End}(J_0(p))$, and let $\mathfrak{m}$ be the maximal ideal
associated to $f$ and $\lambda$. Let $V_{\mathfrak{m}}$ again be a
two-dimensional $\mathbf{T}/\mathfrak{m}$-vector space that affords
$\rho_{\mathfrak{m}}: \text{Gal}(\overline{\mathbb{Q}}/\mathbb{Q})
\rightarrow \text{GL}(2, \mathbf{T}/\mathfrak{m})$. Mazur proved (see
\cite{Mazur1977}, Proposition 14.2) the following multiplicity one
theorem:
$$J[\mathfrak{m}] \cong V_{\mathfrak{m}}.$$

  Now, we give some interpretation on the right-hand side of (7.6). The
transformation theory of elliptic functions leads to the following two
different theories which are given by (7.7): one is the Eichler-Shimura
congruence relation (or more general, the Langlands-Kottwitz method)
and the Eisenstein ideal, which comes from Kronecker's congruence
formula, the other is the modular equation and its Galois theory, which
gives rise to the $\mathbb{Q}(\zeta_p)$-rational Galois representation
$\rho_p: \text{Gal}(\overline{\mathbb{Q}}/\mathbb{Q}) \rightarrow
 \text{Aut}(\mathcal{L}(X(p)))$ and its surjective realization
and modularity. In particular, the defining ideal of the modular
curve is an anabelian counterpart of the Eisenstein ideal of this
modular curve. It is the anabelianization of the Jacobian variety
of this modular curve.

$$\begin{array}{ccc}
  & \text{Jacobi's} & \\
  & \text{transformation} & \\
  & \text{equations of} & \\
  & \text{order $p$} & \\
    \swarrow &  & \searrow\\
  \text{Kronecker}  & \longleftrightarrow & \text{Galois groups of}\\
  \text{congruence} &  & \text{transformation}\\
  \text{relation}   &  & \text{equations}\\
  \text{(Kronecker)} &  & \text{(Galois, Klein)}\\
  \text{(local theory)} &  & \text{(global theory)}\\
  \downarrow &  & \downarrow\\
  \text{Eichler-Shimura relation,} & \longleftrightarrow &
  \text{$\mathbb{Q}(\zeta_p)$-rational Galois}\\
  \text{Serre-Tate theorem} & & \text{representations}\\
  \text{for modular curves} & & \text{for modular curves}\\
  \text{(modular curves mod $p$)} & & \text{and their modularity}\\
  \downarrow &  & \downarrow\\
  \text{Mazur's Eisenstein ideal} & \longleftrightarrow & \text{our defining ideal}\\
  \text{multiplicity-one} & \longleftrightarrow & \text{multiplicity-one}\\
  \text{for mod $p$ Galois} & & \text{for $\mathbb{Q}(\zeta_p)$-rational}\\
  \text{representations} & & \text{Galois representations}\\
  \downarrow &  & \downarrow\\
  \text{Langlands-Kottwitz} & \longleftrightarrow & \text{$K$-rational Galois}\\
  \text{methods for general} & & \text{representations}\\
  \text{Shimura varieties or} & & \rho: \text{Gal}(\overline{\mathbb{Q}}/\mathbb{Q})
  \rightarrow \text{Aut}(I(\text{Sh}_K))\\
  \text{local Shimura varieties:} & & \text{for general Shimura}\\
  \text{Lubin-Tate spaces,} & & \text{varieties and}\\
  \text{Rapoport-Zink spaces} & & \text{their modularity}
\end{array}\eqno{(7.7)}$$

  Now, we point out to the relation with multiplicity-one mod $p$
Galois representations \cite{MR} appearing in the left-hand side
of the above comparison (7.7). In the theory of automorphic
representations of a reductive algebraic group $G$ over a number
field $K$, it is usually true that irreducible representations
occurring in $L^2(G(\mathbb{A})/G(K))$ occur with multiplicity one.
In a classical special case ($G=\text{GL}(2)$, $K=\mathbb{Q}$, and
where we restrict attention to automorphic representations which
are holomorphic, cuspidal, and of weight two), the Galois-theoretic
counterpart of the above multiplicity one phenomenon is the assertion
that given a newform of the above type, of level $N$, the two-dimensional
$p$-adic $\text{Gal}(\overline{\mathbb{Q}}/\mathbb{Q})$-representation
associated to it occurs with multiplicity one in the $p$-adic Tate
module of $J_1(N)$.

  In the case of weight two, Mazur (see \cite{Mazur1977}) discussed
the Eisenstein ideal, which detect the congruences between Eisenstein
series and cusp forms. An ideal of a Hecke ring is called Eisenstein
if it contains $T_{p}-p-1$ for all primes $p$ not dividing the level.
In his paper \cite{Mazur1977}, he proved that $\mathbf{T}(N)/I \simeq
\mathbb{Z}/n \mathbb{Z}$, where $\mathbf{T}(N)$ is the Hecke ring of
level $N$, $I$ is the Eisenstein ideal of $\mathbf{T}(N)$, and $n$ is
the numerator of $\frac{N-1}{12}$. Moreover, for each prime $\ell | n$,
he further proved (see \cite{Mazur1977}, Proposition 14.2) the following
multiplicity one property:
$$\dim J_0(N)[\mathfrak{m}]=2,$$
where $\mathfrak{m}$ is an Eisenstein maximal ideal generated by $I$
and $\ell$, and
$$J_0(N)[\mathfrak{m}]:=\{ x \in J_0(N)(\overline{\mathbb{Q}}):
  \text{$tx=0$ for all $t \in \mathfrak{m}$} \}.$$
Using this result, by studying the action of Galois groups on torsion
points of $J_0(N)$, he proved a classification theorem for the rational
torsion subgroups of elliptic curves over $\mathbb{Q}$. On the other hand,
assume that $\mathfrak{m}$ is a non-Eisenstein maximal ideal of the Hecke
ring $\mathbf{T}$ of level $N$. Then it is shown that (see \cite{BLR})
$J_0(N)[\mathfrak{m}] \simeq V^{\oplus \mu_{\mathfrak{m}}}$, where $V$ is
the underlying irreducible module for the Galois representation
$$\rho_{\mathfrak{m}}: \text{Gal}(\overline{\mathbb{Q}}/\mathbb{Q})
  \rightarrow \text{GL}(2, \mathbf{T}/\mathfrak{m})$$
associated to $\mathfrak{m}$. For example, when the characteristic
of $\mathbf{T}/\mathfrak{m}$ is prime to $2N$, the multiplicity
$\mu_{\mathfrak{m}}$ is one (see \cite{Ribet1990},  Theorem 5.2 (b)).
The motivation of the interest in the mod $p$ multiplicity one question
for Galois representations comes from the Serre's modularity conjecture
(see \cite{Se1987}). Assume that $p$ is an odd prime, and suppose that
$\rho$ is an irreducible mod $p$ representation of
$\text{Gal}(\overline{\mathbb{Q}}/\mathbb{Q})$ which arises from the
space of weight-two modular forms on $\Gamma_0(N)$. We then say that
$\rho$ is modular of level $N$. Assume that $\ell \neq p$ is a prime
factor of $N$ for which $\rho$ is unramified at $\ell$. Then Serre's
conjectures predict that $\rho$ is modular of level $N_0$, where $N_0$
is the prime-to-$\ell$ part of $N$. Ribet (see \cite{Ribet1990} and
\cite{MR}) showed that $\rho$ is modular of level $N_0$ whenever $\ell$
exactly divides $N$ (i.e., $N_0=N/\ell$) and $\rho$ occurs with
multiplicity one in the Jacobian $J_0(N)$. The main results of
\cite{Ribet1990} rely on a multiplicity-one principle (\cite{Ribet1990},
Theorem 5.2 (b)). In his paper \cite{Edixhoven}, Edixhoven extended some
results of Mazur \cite{Mazur1977}, \cite{MR}, Ribet \cite{Ribet1990}
and Gross \cite{Gro90} that state that in certain cases the Galois
representation $\rho_f$ occurs with multiplicity one in a certain
subspace of the $p$-torsion of the Jacobian of the modular curve that
is used to construct $\rho_f$ (see \cite{Edixhoven}, p. 590, Theorem 9.2):

\textbf{Theorem 7.4.} {\it Let $f$ be a cuspidal eigenform of type
$(N, k, \varepsilon)$, defined over $\overline{\mathbb{F}}_p$, with
$p \nmid N$ and $2 \leq k \leq p+1$. Let $J_{\mathbb{Q}}$ be the
Jacobian of the curve $X_1(pN)_{\mathbb{Q}}$ if $k>2$ and let
$J_{\mathbb{Q}}$ be the Jacobian of $X_1(N)_{\mathbb{Q}}$ if $k=2$.
Let $H$ be the subring of} $\text{End}(J_{\mathbb{Q}})$ {\it generated
by all $T_{\ell}$ and $\langle a \rangle_N$, and $\langle b \rangle_p$
if $k>2$. Let $m$ be the maximal ideal of $H$ corresponding to $f$,
and let $\mathbb{F}=H/m \subset \overline{\mathbb{F}}_p$. Suppose
that the representation} $\rho_f: \text{Gal}(\overline{\mathbb{Q}}/\mathbb{Q})
\rightarrow \text{GL}(2, \overline{\mathbb{F}}_p)$ {\it is irreducible.
Then $J_{\mathbb{Q}}(\overline{\mathbb{Q}})[m]$ is an $\mathbb{F}$-vector
space of dimension two in each of the following cases:

$1$. $2 \leq k<p$,

$2$. $\text{$k=p$ and $a_p^2 \neq \varepsilon(p)$, where $T_p^{*}f=a_p f$}$,

$3$. $\text{$k=p$ and $\rho_f$ is ramified at $p$}$,

$4$. $\text{$k=p+1$ and there is no form $g$ of type $(N, 2, \varepsilon)$
   with $\rho_g \cong \rho_f$}$}

  Such multiplicity-one results as well as the Gorenstein property of
Hecke rings played a crucial role in the work of Wiles on Fermat's Last
Theorem (see \cite{W}, p. 482, Theorem 2.1). In general, we have the
following (see \cite{RS}, Theorem 3.5):

\textbf{Theorem 7.5.} {\it An irreducible modular Galois representation}
$$\rho: \text{Gal}(\overline{\mathbb{Q}}/\mathbb{Q}) \longrightarrow
  \text{GL}(2, \overline{\mathbb{F}}_{\ell})$$
{\it is a multiplicity one representation, except perhaps when all of the
following hypothesis on $\rho$ are simultaneously satisfied:

$1$. $k(\rho)=\ell$;

$2$. $\text{$\rho$ is unramified at $\ell$}$;

$3$. $\text{$\rho$ is ordinary at $\ell$}$;

$4$. $\text{$\rho \mid_{D_{\ell}} \sim \left(\begin{matrix}
     \alpha & *\\ 0 & \beta \end{matrix}\right)$ with $\alpha=\beta$}$.}

  In fact, the multiplicity-one phenomenon also appears in the context
of open elliptic modular curves, for example, $Y_1(N)$, by crystalline
cohomology. Fix a level $N$ and a weight $k \geq 2$. Denote by
$S_k(\Gamma_1(N))$ the space of cusp forms of weight $k$ for
$\Gamma_1(N)$ and by $\mathbf{T}$ the Hecke algebra acting on this
space. In \cite{FJ}, Faltings and Jordan showed that for a maximal
ideal $\mathfrak{m} \subset \mathbf{T}$ of residue characteristic
$p>k$ with $(p, N)=1$ the modular Galois representation
$\rho_{\mathfrak{m}}: \text{Gal}(\overline{\mathbb{Q}}/\mathbb{Q})
\rightarrow \text{GL}(2, \overline{\mathbb{F}}_p)$ is crystalline
with the weights $0$ and $(k-1)$ each occurring with multiplicity
one. In this case, the comparison theorem between $p$-adic \'{e}tale
and crystalline cohomology is Hecke-equivariant. The application is
the result of multiplicity one for maximal ideals $\mathfrak{m}
\subset \mathbf{T}$ of residue characteristic $p>k$, $p \nmid N$,
such that the associated modular Galois representation
$\rho_{\mathfrak{m}}$ is irreducible.

\textbf{Theorem 7.6.} (see \cite{FJ}, Theorem 2.1). {\it Let $\mathfrak{m}
\subset \mathbf{T}$ be a maximal ideal with $\mathbf{k}=\mathbf{T}/\mathfrak{m}$
of characteristic $p$. Suppose $\rho_{\mathfrak{m}}$ is irreducible.
If $p>k$, then}:

(1) $H_{\text{par}}^1(Y_1(N)_{\overline{\mathbb{Q}}}, \overline{\vartheta}_p)
[\mathfrak{m}]^{\vee}$ {\it is isomorphic to the}
$\mathbf{k}[\text{Gal}(\overline{\mathbb{Q}}/\mathbb{Q})]$-{\it module
corresponding to $\rho_{\mathfrak{m}}$. Here, for a prime $p$ denote by
$\vartheta_p$ the \'{e}tale sheaf} $\text{Symm}^{k-2}(\mathbb{V})$
{\it on $Y_1(N)/\mathbb{Q}$ with}
$\mathbb{V}=R^1 \pi_{*, \text{\'{e}t}}(\mathbb{Z}_p)$. {\it Define}
$\overline{\vartheta}_p=\vartheta_p/p \vartheta_p$. {\it In particular}
$\dim_{\mathbf{k}} H_{\text{par}}^1(Y_1(N)_{\overline{\mathbb{Q}}},
 \overline{\vartheta}_p)[\mathfrak{m}]=2$.

(2) {\it The local ring $\mathbf{T}_{\mathfrak{m}}$ is Gorenstein.}

  The multiplicity $\mu_{\mathfrak{m}}$ is typically equal to one. To
cite the simplest possible example, take $N=11$. Then $J_0(11)$ is an
elliptic curve, $\mathbf{T}=\mathbb{Z}$, and the ideals $\mathfrak{m}$
are the prime ideals $(p)$ of $\mathbb{Z}$. The kernel $J_0(11)[p]$ is
then an $\mathbb{F}_p$-vector space of dimension two. To make a concrete
example of a maximal ideal $\mathfrak{m}$, take $p=11$. The ring
$\mathbf{T}_{11}$ is isomorphic to $\mathbb{Z}$, and there is a unique
ideal $\mathfrak{m}=(11)$ of residue characteristic $p$. The associated
representation $\rho_{\mathfrak{m}}$ is the two-dimensional representation
$J_0(11)[11]$ of $\text{Gal}(\overline{\mathbb{Q}}/\mathbb{Q})$ over
$\mathbb{F}_{11}$. This representation is known to be absolutely irreducible;
indeed, the associated map $\text{Gal}(\overline{\mathbb{Q}}/\mathbb{Q})
\rightarrow \text{Aut}(J_0(11)[11])$ is surjective (see \cite{MR}).

  Now, we give some interpretation on the right-hand side of (7.7) about the
multiplicity-one property. By Theorem 1.1, Corollary 1.2 and Theorem 1.3, we
have constructed the decompositions of the $\mathbb{Q}(\zeta_p)$-rational
reducible representations $\pi_p$ of $\text{PSL}(2, \mathbb{F}_p)$ as
well as the Galois representations $\rho_p$ given by (7.2) for the
simplest three cases $p=7$, $11$ and $13$. For these cases, both
$\pi_p$ and $\rho_p$ have the multiplicity-one property. In
particular, for $p=11$, we can compare the above Galois representation
$\text{Gal}(\overline{\mathbb{Q}}/\mathbb{Q}) \rightarrow \text{Aut}(J_0(11)[11])$
with our $\rho_{11}: \text{Gal}(\overline{\mathbb{Q}}/\mathbb{Q}) \rightarrow
\text{Aut}(\mathcal{L}(X(11)))$. The first one comes from the $11$-division points
of the Jacobian variety of the modular curve $X_0(11)$ in the context of motives,
while the second one come from the defining ideals of the modular curve $X(11)$
in the context of anabelian algebraic geometry. In fact, the first one comes
from the division equations of elliptic functions, while the second one comes from
the transformation equations of elliptic functions.

  Thus, for modular curves $X(p)$, we have the following two distinct theories:

(1) Eichler-Shimura-Deligne-Serre-Langlands:
$$\begin{array}{ccc}
  \text{eigenvalues of Hecke} & \longleftrightarrow & \text{trace of Frobenius}\\
  \text{operators $a_{\ell}$} &  & \text{elements $\text{Frob}_{\ell}$}\\
  \searrow &  & \swarrow\\
  & \text{Tr$(\rho_f(\text{Frob}_{\ell}))=a_{\ell}$} &\\
  & \text{(\text{prime} $\ell \neq p)$} &
\end{array}\eqno{(7.8)}$$

(2) Ours:
$$\begin{array}{ccc}
  \text{$\mathbb{Q}(\zeta_p)$-rational} & \longleftrightarrow &
  \text{$\mathbb{Q}(\zeta_p)$-rational}\\
  \text{representations} & & \text{representations}\\
  \text{of $\text{PSL}(2, \mathbb{F}_p)$} &  &\text{of $\text{Gal}(\overline{\mathbb{Q}}/\mathbb{Q})$}\\
  \searrow &  & \swarrow\\
  & \text{surjective realization} & \\
  & \text{and modularity} & \\
  & \text{(\text{prime} $\ell=p)$} &
\end{array}\eqno{(7.9)}$$

  Let us give some background. Let $f$ be a newform of weight $k$ and
level $N$ with integer coefficients. Deligne-Serre theorem says that there
exist an associated representation
$$\rho_f^{(\ell)}: \text{Gal}(\overline{\mathbb{Q}}/\mathbb{Q}) \rightarrow
  \text{GL}(2, \mathbb{Q}_{\ell}),$$
for any prime $\ell$ not dividing $N$. In particular, this induces a
representation
$$\rho_{f, \ell}: \text{Gal}(\overline{\mathbb{Q}}/\mathbb{Q}) \rightarrow
  \text{GL}(2, \mathbb{F}_{\ell}).$$
In particular, due to Serre's open image theorem, both adelic and mod
$\ell$ representations are surjective for all but finitely many $\ell$,
when $f$ has weight two and comes from an elliptic curve. More precisely,
we have the following:

\textbf{Theorem 7.7.} (Eichler-Shimura, Deligne, Deligne-Serre). {\it
Let $f$ be a normalized Hecke eigenform of weight $k$, level $N$, and
nebentype character $\chi$. There exists an irreducible Galois
representation}
$$\rho_f: \text{Gal}(\overline{\mathbb{Q}}/\mathbb{Q}) \rightarrow
  \text{GL}(2, \overline{\mathbb{Q}}_p)$$
{\it such that

$(1)$ for each $\ell \nmid p N$, $\rho_f$ is unramified at $\ell$, and
for the Frobenius} $\text{Frob}_{\ell}$ {\it at $\ell$,}
$$\text{Tr}(\rho_f(\text{Frob}_{\ell}))=a_{\ell}, \quad \text{and} \quad
  \det(\rho_f(\text{Frob}_{\ell}))=\ell^{k-1} \chi(\ell);$$

{\it $(2)$ if $p \nmid N$, $\rho_f$ is crystalline at $p$ with Hodge-Tate
weight $(0, k-1)$.}

  Combining these two theories leads to a complete picture on Galois
representations associated with modular curves of level $p$ (including
the cases that the level $p$ is coprime to $\ell$ or not).

  Now, we continue the comparison between division equations of order $p$
and transformation equations of order $p$, which will give rise to a
comparison between Serre's modularity conjecture and our Theorem 1.3.
In his paper \cite{Abel1}, Abel applied complex multiplication to the
division problem of elliptic functions, which is an analogue of Gauss's
theorem on the division of the circle. Moreover, Abel connected the
solvability of the division equations with complex multiplication,
assuming that they are unsolvable whenever the elliptic function does
not enjoy complex multiplication. This leads to the works of Serre (see
\cite{Seabel}, \cite{Se1972} and \cite{Se1987}). Let $E/K$ be an elliptic
curve defined over a number field $K$ (we need only the cases
$K=\mathbb{Q}$ and $K=k$, where $k$ is a complex quadratic field),
let $G=\text{Gal}(\overline{\mathbb{Q}}/K)$ be the Galois group of $K$,
and let $\rho_{\ell}: G \rightarrow \text{Aut}(E_{\ell}) \simeq
\text{GL}(2, \mathbb{F}_{\ell})$ be the representation of $G$ by
automorphisms of $E_{\ell}$, where $E_{\ell}$ is the group of points
of an $\ell$-division on $E$ and $\ell$ is a fixed rational prime. The
equation defining the value $\phi\left(\frac{\omega}{\ell}\right)$,
where $\phi$ is an elliptic function corresponding to $E$, is solvable
by radicals if and only if the group $\rho_{\ell}(G)=G_{\ell}$ is
solvable. If $E$ enjoys complex multiplication, then $G_{\ell}$ is a
subgroup of the abelian group $(\mathcal{O}/ \ell \mathcal{O})^{\times}$,
where $\mathcal{O}=\text{End}_{\mathbb{C}}(E)$ is the ring of complex
multiplications of $E$. When $E$ does not enjoy complex multiplication,
the situation is described by the main theorem of Serre's paper
(\cite{Se1972}):

\textbf{Theorem 7.8.} {\it If $E$ does not enjoy complex multiplication, then
for all but a finite number of primes $\ell$ one has}
$\rho_{\ell}(G)=\text{Aut}(E_{\ell})$.

  This fact, together with the simplicity of $\text{PSL}(2, \mathbb{F}_{\ell})$
for $\ell \geq 5$, shows that if $\ell$ is not contained in some finite exceptional
set $S_E$ of primes $p$ with $\rho_p(G) \neq \text{Aut}(E_p)$, then
$\phi\left(\frac{\omega}{\ell}\right)$ can not be expressed by radicals (see
\cite{Vladut}). The problem of explicit determination of the exceptional set
$S_E$ is very delicate. it has been solved in quite a few cases; for example,
the main results of Mazur's paper \cite{Mazur1978} provides its solution in the
simplest case of a semi-stable $E$. In this case one has
$\rho_{\ell}(G)=\text{Aut}(E_{\ell})$ whenever $\ell \geq 11$. It follows
that for a semi-stable $\phi$ without complex multiplication the values
$\phi\left(\frac{\omega}{\ell}\right)$, $\ell \geq 11$ are not expressible
by radicals. This result can be considered as a far-reaching development of
Abel's ideas. In particular, this leads to the theory of mod $\ell$ and
$\ell$-adic Galois representations.

  In fact, the level $N$ modular curve $X(N)$ is the moduli space for
elliptic curves with a basis for the $N$-torsion. This leads to two
distinct kinds of Galois representations. The first kind of Galois
representations are the following (see \cite{Seabel}, \cite{Se1972}):
let $E/\mathbb{Q}$ be an elliptic curve and $\ell$ a prime. The Galois
group $\text{Gal}(\overline{\mathbb{Q}}/\mathbb{Q})$ acts on
$E[\ell] \cong \mathbb{Z}/ \ell \mathbb{Z} \times \mathbb{Z}/\ell \mathbb{Z}$
and gives rise to the mod $\ell$ Galois representation of $E$:
$$\overline{\rho}_{E, \ell}: \text{Gal}(\overline{\mathbb{Q}}/\mathbb{Q})
  \rightarrow \text{GL}(2, \mathbb{F}_{\ell}).$$
Write $T_{\ell}(E)=\varprojlim_{n} E[\ell^n]$ for the $\ell$-adic Tate
module of $E$. Then
$T_{\ell}(E) \cong \mathbb{Z}_{\ell} \times \mathbb{Z}_{\ell}$,
and by choosing a $\mathbb{Z}_{\ell}$-basis for $T_{\ell}(E)$ we obtain
the $\ell$-adic Galois representation of $E$, which describes the action
of the absolute Galois group $\text{Gal}(\overline{\mathbb{Q}}/\mathbb{Q})$
on $T_{\ell}(E)$:
$$\rho_{E, \ell}: \text{Gal}(\overline{\mathbb{Q}}/\mathbb{Q}) \rightarrow
  \text{GL}(2, \mathbb{Q}_{\ell}).$$
In particular, Serre (\cite{Seabel} and \cite{Se1972}) and Mazur
(\cite{Mazur1977} and \cite{Mazur1978}) obtain semi-stable elliptic
curve realization
$$\overline{\rho}_{E, \ell}: \text{Gal}(\overline{\mathbb{Q}}/\mathbb{Q})
  \rightarrow \text{Aut}(E[\ell])\eqno{(7.10)}$$
of the mod $\ell$ Galois representation
$$\overline{\rho}_{\ell}: \text{Gal}(\overline{\mathbb{Q}}/\mathbb{Q})
  \rightarrow \text{GL}(2, \mathbb{F}_{\ell})\eqno{(7.11)}$$
by the $\ell$-torsion points of $E/\mathbb{Q}$ such that it is surjective
for $\ell \geq 11$. Moreover, for such a semi-stable elliptic curve $E$
defined over $\mathbb{Q}$, let $\ell$ be a prime number $\geq 5$. We will
focus on the representation
$$\rho_{E, \ell}: \text{Gal}(\overline{\mathbb{Q}}/\mathbb{Q}) \rightarrow
  \text{GL}(2, \mathbb{F}_{\ell})\eqno{(7.12)}$$
given by the $\ell$-torsion points of $E$. We have the following
irreducibility result:

\textbf{Proposition 7.9.} (Mazur). {\it The representation $\rho_{E, \ell}$
is irreducible.}

  This proposition is proved as Proposition 6 in \cite{Se1987}, \S 4.1.
The proof given by Serre in \cite{Se1987} is based on results of Mazur
(\cite{Mazur1977}, Theorem 8) and relies on the fact that the $2$-division
points of $E$ are rational over $\mathbb{Q}$. In the same paper \cite{Se1987},
Serre formulated his modularity conjecture which gives a mod $\ell$ modular
form realization of this mod $\ell$ Galois representation. The significance
of Serre's conjecture is as follows: In one direction, one can think about
representations coming from $\ell$-division points on elliptic curves,
or more generally from $\ell$-division points on abelian varieties of
$\text{GL}(2)$-type. Conversely, Serre's conjectures imply that
all odd irreducible two-dimensional mod $\ell$ representations of
$\text{Gal}(\overline{\mathbb{Q}}/\mathbb{Q})$ may be realized in
spaces of $\ell$-division points on such abelian varieties.

  On the other hand, the second kind of Galois representations come from
$\mathcal{L}(X(p))$. The Galois group $\text{Gal}(\overline{\mathbb{Q}}/\mathbb{Q})$
acts on $\mathcal{L}(X(p))$ by automorphisms of $\mathcal{L}(X(p))$ and
gives rise to the $\mathbb{Q}(\zeta_p)$-rational Galois representation
of $\mathcal{L}(X(p))$:
$$\rho_{\mathcal{L}, p}: \text{Gal}(\overline{\mathbb{Q}}/\mathbb{Q})
  \rightarrow \text{Aut}(\mathcal{L}(X(p))) \cong \text{PSL}(2, \mathbb{F}_p).
  \eqno{(7.13)}$$
Theorem 1.3 shows that it has a modular and surjective realization for
$p \geq 7$. This leads us to give the following comparison:
$$\begin{matrix}
 &\text{motives} &\longleftrightarrow & \text{anabelian algebraic geometry}\\
 &\text{$\ell$-adic representations associated} &\longleftrightarrow
 &\text{Galois representations $\rho_p$}\\
 &\text{to abelian varieties (\cite{Seabel}, \cite{Se1972})} &
 &\text{associated to $\mathcal{L}(X(p))$}\\
 &\text{semisimplicity of} &\longleftrightarrow & \text{semisimplicity of}\\
 &\text{$\ell$-adic representations (\cite{Faltings})} & & \text{Galois representations $\rho_p$}\\
 &\text{Hodge cycles behave under the} &\longleftrightarrow & \text{algebraic cycles}\\
 &\text{action of the Galois group as if} & &\text{under the action}\\
 &\text{they were algebraic cycles (\cite{DeHodge})} & & \text{of the Galois group}
\end{matrix}\eqno{(7.14)}$$

  Now, we give the relation with the Artin conjecture. One important aspect
of the present view is that the theory of ordinary complex representations
of finite Galois groups, i.e., the theory of Artin motives, can not be
separated from the theory of Grothendieck motives, i.e., from the theory
of system of $\ell$-adic representations of Galois groups coming from
\'{e}tale cohomology of algebraic varieties. Artin motives are motives of
zero-dimensional varieties: the products of $H^0(L/K)$ for finite extensions
$L/K$. Grothendieck's vision of motives is as a universal cohomology theory
but also as higher dimensional version of Galois theory. On the other hand,
our $\mathbb{Q}(\zeta_p)$-rational representations of finite Galois groups
$\text{PSL}(2, \mathbb{F}_p)$ arise from defining ideals of modular curves
over $\mathbb{Q}$, which is independent of $\ell$. Now, there are three kinds
of geometric constructions associated with Galois representations: Artin
motives, Grothendieck motives and defining ideals:
$$\aligned
 &\text{(1) complex Galois representations, Artin motives}\\
 &\quad \quad \text{of algebraic number fields (independent of $\ell$);}\\
 &\text{(2) $\ell$-adic Galois representations, Grothendieck motives}\\
 &\quad \quad \text{of algebraic varieties ($\ell$-adic system);}\\
 &\text{(3) $\mathbb{Q}(\zeta_p)$-rational Galois representations, defining ideals}\\
 &\quad \quad \text{of algebraic varieties (independent of $\ell$).}
\endaligned\eqno{(7.15)}$$
Correspondingly, we have the following correspondences among Galois
representations, algebraic varieties and automorphic forms:
$$\begin{matrix}
 &\text{$\ell$-adic representations} &\leftrightarrow
 &\text{$\ell$-adic \'{e}tale cohomology} &\leftrightarrow
 &\text{cusp forms of}\\
 &\text{(Deligne)} &  &\text{of modular curves} &  &\text{weight $k \geq 2$}\\
 &\text{Arin representations} &\leftrightarrow &\text{coherent cohomology}
 &\leftrightarrow &\text{cusp forms of}\\
 &\text{(Deligne-Serre)} &  &\text{of modular curves} &  &\text{weight $k=1$}\\
 &\text{$\mathbb{Q}(\zeta_p)$-rational} &\leftrightarrow
 &\text{defining ideals} &\leftrightarrow &\text{cusp forms of}\\
 &\text{representations} &  &\text{of modular curves} &  &\text{weight $k=\frac{1}{2}$}\\
\end{matrix}\eqno{(7.16)}$$
Here, the cusp forms of weight $k=\frac{1}{2}$ refer to the theta constants
of order $p$. Therefore, for the same modular curves over $\mathbb{Q}$, there
exist two distinct kinds of correspondences, one comes from automorphic
representations of $\text{GL}(2, \mathbb{A}_{\mathbb{Q}})$, i.e. Grothendieck
motives (the systems of $\ell$-adic representations), the other comes from
$\mathbb{Q}(\zeta_p)$-rational representations of $\text{PSL}(2, \mathbb{F}_p)$,
i.e., defining ideals $I(\mathcal{L}(X(p)))$ (independent of $\ell$).

  It should be pointed out that there is a dictionary between the classical
and modern theories for $\text{GL}(2)$ (see \cite{Gelbart}).

\textbf{Proposition 7.10.} {\it There is a one-to-one correspondence between
normalized new forms
$$f(z)=\sum_{n=1}^{\infty} a_n e^{2 \pi i n z}$$
in $S_k(\Gamma_0(N), \psi)$, and irreducible automorphic cuspidal representations
$$\pi=\otimes_{p} \pi_p$$
of $G(\mathbb{A})$ such that:}

(a) {\it the central character of $\pi$ is $\chi_{\psi}$,}

(b) {\it $\pi_p$ is unramified for all $p \nmid N$, and}

(c) {\it $\pi_{\infty}$ has ``lowest weight $k$.''}

{\it Moreover, if $N=\prod p_i^{\alpha_i}$, and $p | N$, then $\pi_p$ has
conductor $\alpha_i$; but if $p \nmid N$, then $\pi_p=\pi(\mu_1, \mu_2)$,
with $\mu_1$, $\mu_2$ unramified characters of $\mathbb{Q}_p^{\times}$ such
that}
$$p^{\frac{k-1}{2}}(\mu_1(p)+\nu_2(p))=a_p.$$

  In case $k>1$, then $\pi_{\infty}$ is a discrete series representation
of ``lowest weight $k$,'' sitting inside some $\pi_{\infty}(\mu_1, \mu_2)$;
this discrete series representation has central character trivial on $\mathbb{R}^{+}$,
and is denoted $D_k$. However, if $k=1$, then $\pi_{\infty}$ is not a
sub-representation of any $\pi_{\infty}(\mu_1, \mu_2)$, but rather equals
the full induced (principal series) representation $\pi_{\infty}(1, \text{sgn})$.
In this case, the corresponding representation $\pi=\otimes \pi_p$ (with
$\pi_{\infty}=\pi_{\infty}(1, \text{sgn})$) is called an automorphic cuspidal
representation of weight one. Hence, in (7.16), the case of $k \geq 2$ corresponds
to the discrete series representation, while the case of $k=1$ corresponds to
the principal series representation. In contrast with it, in our case of
$k=\frac{1}{2}$, all of the representations of $\text{PSL}(2, \mathbb{F}_p)$
do appear, namely the discrete series representation, the principal series
representation, the Steinberg representation and the trivial representation.

  It fact, the only known constructions of Galois representations involve
the \'{e}tale cohomology of varieties in some way. Even more, an Artin
representation comes from the degree zero \'{e}tale cohomology of a
zero-dimensional variety. In contrast, our construction of Galois
representations arise from the defining ideals of modular curves (in
the context of anabelian algebraic geometry) which does not involve the
\'{e}tale cohomology of varieties.

  From the analytic viewpoint, we have the following comparison between
transformation equations and division equations of elliptic functions.
Here, the first principal transformation is given by
$\left(\begin{matrix} p & 0\\ 0 & 1 \end{matrix}\right)$,
the second principal transformation is given by
$\left(\begin{matrix} 1 & 0\\ 0 & p \end{matrix}\right)$.
The $p-1$ transformations
$\left(\begin{matrix} 1 & 1\\ 0 & p \end{matrix}\right)$,
$\left(\begin{matrix} 1 & 2\\ 0 & p \end{matrix}\right)$, $\ldots$,
$\left(\begin{matrix} 1 & p-1\\ 0 & p \end{matrix}\right)$
are associated with the second principal transformation. The Hecke operator
$T_p$ is an average of these two kinds of transformations.
$$\begin{matrix}
  \text{transformation} & \rho_p: \text{Gal}(\overline{\mathbb{Q}}/\mathbb{Q})
  \rightarrow \text{Aut}(\mathcal{L}(X(p))) & \leftrightarrow & \text{Galois groups}\\
  \text{equations} & \mathcal{L}(X(p)): \text{nonlinear} &  & \text{of $p$-th degree}\\
                  & \pi_1(X(p)): \text{anabelian}  &  & \text{transformation of}\\
                       &                           &  & \text{$j$-functions (1st and}\\
                       &                           &  & \text{2nd principal}\\
                       &                           &  & \text{transformations)}\\
  \text{by resolvents} & \text{by linearization}   &  & \text{by average}\\
                       & \text{and abelianization} &  &                  \\
           \downarrow  & \downarrow                &  & \downarrow\\
       \text{division} & H_1: \text{the first \'{e}tale} & \leftrightarrow &
       \text{Hecke operators}\\
      \text{equations} & \text{homology group of} &  &
  T_p f=\frac{1}{p} \sum_{i=0}^{p-1} f(\frac{z+i}{p})\\
                       & \text{the elliptic curve $E/\mathbb{Q}$} &  & +pf(pz)\\
                       & \text{$p^n$-torsion points} &  & \text{eigenforms $f$}\\
                       & \rho_E: G_{\mathbb{Q}} \rightarrow \text{Aut}(E[p^n]) & \leftrightarrow &
                       \rho_f: G_{\mathbb{Q}} \rightarrow \text{GL}(2, \overline{\mathbb{Q}}_p)\\
                       & \text{tr}(\rho_E(\text{Frob}_{\ell}))=a_{\ell} \quad (\ell \neq p)
                       & \leftrightarrow & T_{\ell} f=a_{\ell} f \quad (\ell \neq p)
\end{matrix}\eqno{(7.17)}$$

  In conclusion, given a Galois representation (7.2), there are two distinct kinds
of realizations: (1) Realization of motives: $\rho_p$ can be realized by the division
points of elliptic curves, which come from division equations of elliptic functions.
This leads to Serre's modularity conjecture, i.e., mod $p$ Langlands correspondence
for $\text{GL}(2, \mathbb{Q})$. (2) Realization of anabelian algebraic geometry:
$\rho_p$ can be realized by the defining ideals $I(\mathcal{L}(X(p)))$, which come
from transformation equations of elliptic functions. This leads to Theorem 1.1 and
Theorem 1.3, i.e., anabelian algebraic geometric version for $\text{GL}(2, \mathbb{Q})$.

\begin{center}
{\large\bf 7.1. Relation with the Artin conjecture and Langlands' principle
                of functoriality}
\end{center}

  Emil Artin associated to a complex Galois representation
$\rho: \text{Gal}(\overline{\mathbb{Q}}/\mathbb{Q}) \rightarrow
 \text{GL}(n, \mathbb{C})$ an $L$-function $L(\rho, s)$ and conjectured
that it has an analytic continuation to the whole complex plane. Via
work of Deligne and Serre, the Langlands program relates such
representations, when $n=2$, to certain cusp forms of weight one on a
group slightly different from $\Gamma_0(N)$. In fact, using techniques
of base change for the group $\text{GL}(2)$, Langlands (\cite{Langlands80})
showed how to establish the holomorphy of Artin's $L$-series attached
to two-dimensional representations of the solvable Galois group $A_4$
using the following Deligne-Serre theorem, and thus prove that such
representations come from automorphic forms. For general Artin $L$-functions,
Langlands's principle of functoriality can be interpreted as an identity
(see \cite{Arthur2024}):
$$L^{S}(s, r)=L^{S}(s, \pi, St(n))$$
between a general Artin $L$-function and a standard automorphic $L$-function
for $\text{GL}(n)$. This represents a general formulation of nonabelian
class field theory. It identifies purely arithmetic objects, Artin
$L$-functions, with objects associated with harmonic analysis, automorphic
$L$-functions, thereby implying that the arithmetic $L$-functions have
meromorphic continuation and functional equation, and that they are
essentially entire. Functoriality would thus include a proof of the Artin
conjecture.

  For $\text{GL}(2, \mathbb{Q})$, in contrast with Artin's as well as
Langlands' strategies to nonabelian class field theory, Theorem 1.3
gives a connection between Galois representations and Klein's elliptic
modular functions by the defining ideals of modular curves. It states
that for the Galois representations
$\rho_p: \text{Gal}(\overline{\mathbb{Q}}/\mathbb{Q}) \rightarrow
 \text{Aut}(\mathcal{L}(X(p)))$ with $p \geq 7$, there are surjective
realizations by the equations of degree $p+1$ or $p$ with Galois group
isomorphic with $\text{Aut}(\mathcal{L}(X(p)))$. These equations are
defined over $\mathcal{L}(X(p))$, their coefficients are invariant under
the action of $\text{Aut}(\mathcal{L}(X(p)))$. Moreover, it corresponds to
the $p$-th order transformation equation of the $j$-function with Galois
group isomorphic to $\text{PSL}(2, \mathbb{F}_p)$ $(p \geq 7)$. In particular,
this goes beyond the scope of two-dimensional Artin conjecture, which only
includes the Galois groups $A_4$, $S_4$ and $A_5$.

\begin{center}
{\large\bf 7.2. Relation with Deligne's theorem ($k \geq 2$) and the
                Deligne-Serre theorem ($k=1$)}
\end{center}

  In \cite{De1968}, Deligne proved Serre's conjecture on the existence
of $\ell$-adic representations associated to modular forms (\cite{Se1969}).
This result has been essential for the study of modular forms modulo $p$ and
for that of $p$-adic modular forms. A converse theorem over $\mathbb{Q}$ was
established by Deligne-Serre (\cite{DeSe}) for holomorphic modular forms of
weight one: every such form gives rise to a two-dimensional Galois
representation with odd determinant. More precisely, Deligne and Serre proved
that every cusp form of weight one, which is an eigenfunction of the Hecke
operators, corresponds by Mellin's transform to the Artin $L$-function of
an irreducible complex representation $\rho:
\text{Gal}(\overline{\mathbb{Q}}/\mathbb{Q}) \rightarrow \text{GL}(2, \mathbb{C})$.
In order to prove this theorem, they constructed mod $\ell$ Galois
representations, for each prime $\ell$, of sufficiently small image, which
allows them to lift these representations to characteristic zero and to
obtain from them the desired complex representation. In particular, they
used the existence of $\ell$-adic representations associated to modular
forms of weight $k \geq 2$ to deduce an existence theorem for complex
representations associated to weight $k=1$ modular forms.

  A $p$-adic Galois representation $\rho: \text{Gal}(\overline{F}/F) \rightarrow
\text{GL}(d, E)$, $F$ a number field, $E$ a finite extension of $\mathbb{Q}_p$,
which is unramified outside a finite set of primes and potentially semi-stable
at primes above $p$ is called geometric by J.-M. Fontaine and B. Mazur (see
\cite{FM}). A subquotient of the Galois representation given by the $p$-adic
\'{e}tale cohomology of a projective and smooth algebraic variety $X$ over $F$
is geometric. It is unramified outside $p$ and the primes of bad reduction of
$X$, and it is potentially semistable by the $p$-adic comparison theorem.
Conversely, J.-M. Fontaine and B. Mazur conjecture that a $p$-adic geometric
irreducible representation of $\text{Gal}(\overline{F}/F)$ comes from such a
subquotient.

  The following two kinds of Galois representations are geometric (see
\cite{KW}) in the context of motives:

(1) Deligne (\cite{De1968}): $\ell$-adic Galois representations attached
to an $f$ of weight $k \geq 2$ is geometric, as it appears in the cohomology
of a fiber product of the universal generalized elliptic curve over a modular
curve (Grothendieck motives).

(2) Deligne-Serre (\cite{DeSe}): complex Galois representations attached to an
$f$ of weight $k=1$ has finite image, hence is also geometric (Artin motives).

  By Theorem 1.1 and Theorem 1.3, the $\mathbb{Q}(\zeta_p)$-rational Galois
representations
$$\rho_p: \text{Gal}(\overline{\mathbb{Q}}/\mathbb{Q})
  \longrightarrow \text{PSL}(2, \mathbb{F}_p)\eqno{(7.2.1)}$$
can be realized by $\text{Aut}(\mathcal{L}(X(p)))$, which are connected with
the theta constants of order $p$ and weight $k=\frac{1}{2}$. In particular,
$\rho_p$ arises from geometry in the context of anabelian algebraic geometry,
as it both has finite image and comes from the defining ideals of modular curves.

  Therefore, all of the three distinct kinds of Galois representations appearing
in (7.15) and (7.16) come from algebraic geometry (Grothendieck motives, Artin
motives, and anabelian algebraic geometry).

\begin{center}
{\large\bf 7.3. Relation with Serre's modularity conjecture and Serre weights}
\end{center}

  In his paper (\cite{Se1987}), Serre assumed at the start that the
level is coprime to $p$. He pointed out that this is not necessary:
all the stated results remain true in the general case. However, the
gained generality does not supply mod $p$ modular forms that are
genuinely new; indeed we know that any cusp form with coefficients
in $\overline{\mathbb{F}}_p$ of level $p^m N$ is also of level $N$,
at the expense of increasing the weight. A typical example is that
of forms of weight $2$ and level $p$, which are also of weight $p+1$
and level $1$.

  Serre's modularity conjecture gives a relation between mod $p$
Galois representations and mod $p$ modular forms. Let
$\rho: \text{Gal}(\overline{\mathbb{Q}}/\mathbb{Q}) \rightarrow
\text{GL}(2, \overline{\mathbb{F}}_p)$ be a continuous, odd, irreducible
Galois representation. Then there exists a normalized cuspidal mod $p$
eigenform $f$, such that $\rho$ is isomorphic to $\rho_f$, the Galois
representation associated to $f$. Moreover, $f \in S_{k(\rho)}(N(\rho),
\varepsilon(\rho); \overline{\mathbb{F}}_p)$, where $N(\rho)$ and
$k(\rho)$ are the minimal weight and level. Serre's modularity conjecture
can be viewed as a first step in the direction of a mod $p$ analogue of
the Langlands program.

  Note that the conjecture of Artin about odd representations can be
stated that if we have a continuous representation
$\rho: \text{Gal}(\overline{\mathbb{Q}}/\mathbb{Q}) \rightarrow
\text{GL}(2, \mathbb{C})$ which is odd and irreducible, then it arises
from a newform of weight one. In this way, Khare showed that Serre's
conjecture implies Artin's conjecture with odd representations.
Moreover, assume Serre's conjecture. A representation $\rho$ of Serre
type arises from a classical newform of weight one via the construction
in Deligne-Serre \cite{DeSe}, if and only if $\rho$ lifts to a complex
representation (see \cite{Khare1997} as well as \cite{KW1}, \cite{KW2}
and \cite{Kisin2009}).

  Theorem 1.3 gives a relation between $\mathbb{Q}(\zeta_p)$-rational
Galois representations $\rho_p$ given by (7.2.1) and the elliptic modular
functions, i.e., the $j$-functions. Now, we will give a comparison between
Serre's modularity conjecture with Theorem 1.3 in the special cases $p=7$,
$11$ and $13$.

$\bullet$ For $p=7$, at the end of his paper \cite{Se1987}, Serre gave an
example using the simple group $\text{PSL}(2, \mathbb{F}_7)$ of order $168$.
The degree seven extension of $\mathbb{Q}$ defined by the equation
$$X^7-7X+3=0\eqno{(7.3.1)}$$
has Galois group $\text{PSL}(2, \mathbb{F}_7)$. Serre used it to
construct a representation of $\text{Gal}(\overline{\mathbb{Q}}/\mathbb{Q})$
in characteristic $7$ as follows: Let $G$ be the subgroup of
$\text{GL}(2, \mathbb{F}_{49})$ defined by:
$$G=\{ \pm 1, \pm i \} \cdot \text{SL}(2, \mathbb{F}_7)
   =\text{SL}(2, \mathbb{F}_7) \cup i \cdot \text{SL}(2, \mathbb{F}_7),$$
where $i$ is an element of order $4$ of $\mathbb{F}_{49}^{\times}$. We have
$\det G=\{ \pm 1 \}$, and the image of $G$ in $\text{PGL}(2, \mathbb{F}_{49})$
is $\text{PGL}(2, \mathbb{F}_7)$. We get the exact sequence:
$$\{ 1 \} \rightarrow \{ \pm 1 \} \rightarrow G \rightarrow
  \text{PSL}(2, \mathbb{F}_7) \times \{ \pm 1 \} \rightarrow \{ 1 \}.$$
Let $K$ be the field of degree $7$ defined by the above equation, and let
$\alpha^K: \text{Gal}(\overline{\mathbb{Q}}/\mathbb{Q}) \rightarrow
\text{PSL}(2, \mathbb{F}_7)$ be the corresponding homomorphism. On the
other hand, let
$$\varepsilon: \text{Gal}(\overline{\mathbb{Q}}/\mathbb{Q})
  \rightarrow \{ \pm 1 \}$$
be the quadratic character associated with the field $\mathbb{Q}(\sqrt{-3})$.
The pair $(\alpha^K, \varepsilon)$ defines a homomorphism
$$\alpha: \text{Gal}(\overline{\mathbb{Q}}/\mathbb{Q}) \rightarrow
  \text{PSL}(2, \mathbb{F}_7) \times \{ \pm 1 \}.\eqno{(7.3.2)}$$
Let $\text{obs}(\alpha) \in \text{Br}_2(\mathbb{Q})$ be the obstruction
to lifting $\alpha$ to a homomorphism
$$\rho: \text{Gal}(\overline{\mathbb{Q}}/\mathbb{Q}) \rightarrow
  G \subset \text{GL}(2, \mathbb{F}_{49}).$$
It can be shown that
$$\text{obs}(\alpha)=w+(-1)(-3),$$
where $w$ is the Witt invariant of the quadratic form
$\text{Tr}_{K/\mathbb{Q}}(x^2)$. In fact, we have $w=(-1)(-3)$, hence
$\text{obs}(\alpha)=0$. This proves the existence of the representation
$$\rho: \text{Gal}(\overline{\mathbb{Q}}/\mathbb{Q}) \rightarrow
  \text{GL}(2, \mathbb{F}_{49})$$
we are looking for. By construction, we have $\det \rho=\varepsilon$.

  Once again, we choose $\rho$ so that its conductor is as small as
possible. The discriminant of the polynomial $X^7-7X+3$ is $3^8 7^8$
and that of the field $K$ is $3^6 7^8$. It follows that the conductor
of $\rho$ can be chosen to be $3^n$, and a ramification calculation
shows that $n=3$. On the other hand, the study of the ramification at
$7$ shows that the action of the inertia at $7$ is:
$$\text{either} \quad \left(\begin{matrix} \chi & *\\ 0 & \chi^{-1}
  \end{matrix}\right), \quad \text{either} \quad \left(\begin{matrix}
  \chi^4 & *\\ 0 & \chi^{-4} \end{matrix}\right),$$
where $\chi$ is the cyclotomic character.

  After tensoring $\rho$ by $\chi$, or by $\chi^4$, we get a new
representation $\rho^{\prime}$ where the action of inertia at $7$
is given by:
$$\left(\begin{matrix}
  \chi^2 & *\\
  0 & 1
  \end{matrix}\right),$$
which leads to a weight $k$ equal to $3$. We have
$$\det \rho^{\prime}=\varepsilon \cdot \chi^2.$$

  Serre's modularity conjecture states that $\rho^{\prime}$ is of
the form $\rho_f$, where $f=\sum a_n q^n$ is a cusp form of type
$(3^3, 3, \varepsilon)$, with coefficients in $\mathbb{F}_{49}$,
and which is a normalized eigenfunction for the Hecke operators.
The link between the eigenvalues $a_{\ell}$ $(\ell \neq 3, 7)$ and
the decomposition of $\ell$ in $K$ is the following: if we write
$\text{ord}(\ell)$ for the order of the Frobenius automorphism
attached to $\ell$ in $\text{Gal}(K^{\text{gal}}/\mathbb{Q})
\cong \text{PSL}(2, \mathbb{F}_7)$, we must have:
$$\aligned
  &\text{ord$(\ell)=1$ or $7$} & \Leftrightarrow \quad
  &a_{\ell}^2=4 \ell^2 \varepsilon(\ell) \quad &\text{in $\mathbb{F}_7$},\\
  &\text{ord$(\ell)=2$} & \Leftrightarrow \quad
  &a_{\ell}=0 \quad &\text{in $\mathbb{F}_7$},\\
  &\text{ord$(\ell)=3$} & \Leftrightarrow \quad
  &a_{\ell}^2=\ell^2 \varepsilon(\ell) \quad &\text{in $\mathbb{F}_7$},\\
  &\text{ord$(\ell)=4$} & \Leftrightarrow \quad
  &a_{\ell}^2=2 \ell^2 \varepsilon(\ell) \quad &\text{in $\mathbb{F}_7$},
\endaligned$$
with $\varepsilon(\ell)=\left(\frac{\ell}{3}\right)$.

  Indeed, we can find a form $f$ with these properties. It is the
reduction (mod $7$) of a newform $F$ in characteristic $0$:
$$F=q+\sum_{n \geq 2} A_n q^n=q+3i q^2-5 q^4-3i q^5+5 q^7-3i q^8+\cdots.$$
This form has coefficients in $\mathbb{Z}[i]$.

  Calculation of $F$. Let $\theta_1$ be the theta function associated
with the field $\mathbb{Q}(\sqrt{-3})$:
$$\theta_1=\sum_{x, y \in \mathbb{Z}} q^{x^2+xy+y^2}
          =1+6(q+q^3+q^4+2 q^7+q^9+\cdots).$$
It is an Eisenstein series of weight $1$, level $3$ and character
$\varepsilon$. If we set
$$\aligned
  \theta_2 &=\theta_1(3z)=1+6 (q^3+q^9+q^{12}+\cdots),\\
  \theta_3 &=\theta_1(9z)=1+6 (q^9+q^{27}+q^{36}+\cdots),
\endaligned$$
we obtain forms of levels $3^2$ and $3^3$.

  On the other hand, the series
$$g=q \prod_{n \geq 1} \left(1-q^{3n}\right)^2 \left(1-q^{9n}\right)^2
   =q-2 q^4-q^7+5 q^{13}+\cdots$$
is the unique normalized cusp form of weight $2$, level $3^3$ and trivial
character (it corresponds to the elliptic curve $y^2+y=x^3-3$, of
conductor $3^3$).

  The products $g \theta_1$, $g \theta_2$ and $g \theta_3$ are forms
of weight $3$, level $3^3$ and character $\varepsilon$. They form a
basis for the space of cusp forms of type $(3^3, 3, \varepsilon)$. The
series
$$F=\frac{1}{2} i g \theta_1-\frac{1}{2} (1+i) g \theta_2+\frac{3}{2} g \theta_3
   =q+3i q^2-5 q^4+\cdots\eqno{(7.3.3)}$$
is the cusp form we are looking for.

  It is interesting to compare the above result with the result given
in the section 4.1, especially the comparison between the formulas
(7.3.1) and (4.1.24), the Galois representations (7.3.2) and (7.2.1)
for $p=7$, as well as the modular forms given by (7.3.3) and the
seven-order transformation of elliptic modular functions given by
(4.1.22).

  For primes $p \geq 11$, in his proof of Serre's conjecture on companion
forms \cite{Gro90}, Gross used $p$-adic techniques, and specifically the
different $p$-adic cohomology theories (de Rham, crystalline, Washnitzer-Monsky)
of modular curves and their Jacobians. In particular, Gross gave some
examples of cuspidal eigenforms (mod $p$)  of weight $k \leq p+1$ on
$\Gamma_1(N)$, discussed their liftings to forms of weight $2$ on
$\Gamma_1(Np)$ and their Galois representations. The cusp form
$$\Delta(z)=q \prod_{n=1}^{\infty} (1-q^n)^{24}
 =\sum_{n=1}^{\infty} \tau(n) q^n$$
has weight $12$ for $\Gamma_1(1)=\text{SL}(2, \mathbb{Z})$ and is defined
over $\mathbb{Z}$. It gives an eigenform (mod $p$) for all primes $p$, but
we will assume that $p \geq 11$ so as to have the inequality $k \leq p+1$.
When $p \geq 11$, $p \neq 23$ or $691$, the Galois representation
$\rho_{\Delta}: \text{Gal}(\overline{\mathbb{Q}}/\mathbb{Q}) \rightarrow
\text{GL}(2, \mathbb{F}_p)$ has image the subgroup of invertible matrices
$A$ with $\det A \in (\mathbb{Z}/p \mathbb{Z})^{* 11}$. In particular,
$\rho_{\Delta}$ determines a Galois extension of $\mathbb{Q}(\mu_p)$
which is unramified outside $p$ and has Galois group
$\simeq \text{SL}(2, \mathbb{F}_p)$.

$\bullet$ For $p=11$ we have $k=p+1$. There is a lifting of $f=\Delta$ to a form
$F$ of weight $2$ and trivial character on $\Gamma_1(11)$. The unique such
form has $q$-expansion
$$F=q \prod_{n=1}^{\infty} (1-q^n)^2 (1-q^{11 n})^2.$$
The representation $\rho_{\Delta}$ occurs on the $11$-division points of
the elliptic curve $J_0(11)$ or equivalently, on the $11$-division points
of the $5$-isogenous curve $J_1(11)$. Its restriction to the inertia group
at $11$ is tr\`{e}s ramifi\'{e} (highly ramified) in the sense of Serre
\cite{Se1987}.

  It is interesting to compare the above result with the result given
in the section 2.4 and the section 4.2. Note that the main difference
is that Gross used $11$-division points of the Jacobian variety of
modular curves $X_0(11)$ or $X_1(11)$, which comes from the division
equations of elliptic functions. On the other hand, we use the defining
ideal $I(\mathcal{L}(X(11)))$, especially (2.4.34), (2.4.35) and (4.2.9),
which comes from the transformation equations of elliptic functions of
order eleven given by (4.2.11) and (4.2.16).

$\bullet$ When $p=13$ we have $k=p-1$. There is a lifting of $f=\Delta$ to
a form $F$ of weight $2$ and character $\omega^{10}$ on $\Gamma_1(13)$.
The unique such form has $q$-expansion in the subring $\mathbb{Z}[\mu_6]$
of $\mathbb{Z}_{13}$. If $\alpha$ is the unique $6$-th root of $1$ in
$\mathbb{Z}_{13}^{*}$ which satisfies $\alpha \equiv 4$ (mod $13$) then
the Fourier expansion of $F$ begins
$$F=q+(-2+\alpha) q^2+(-2 \alpha) q^3+(1-\alpha) q^4+(2 \alpha-1) q^5
   +(2+2 \alpha) q^6+\cdots.$$
The form $F$ is congruent to the eigenform $f=\Delta_{12}$ (mod $13$).
The Galois representation $\rho_{\Delta}=\rho_f$ occurs in the subspace
of $13$-division points of $J_1(13)$ where the Galois group
$(\mathbb{Z}/13 \mathbb{Z})^{*}/ \langle \pm 1 \rangle$ of
$X_1(13)$ over $X_0(13)$ acts by $\omega^2$. As there is no companion
form, $\rho_f$ gives an $\text{SL}(2, \mathbb{F}_{13})$-extension of
$\mathbb{Q}(\mu_{13})$, which is wildly ramified at $13$.

   It is interesting to compare the above result with the result given
in the section 2.4 and the section 4.3. Note that the main difference
is that Gross used $13$-division points of the Jacobian variety of
modular curves $X_0(13)$ or $X_1(13)$, which comes from the division
equations of elliptic functions. On the other hand, we use the defining
ideal $I(\mathcal{L}(X(13)))$, especially (2.4.44), (2.4.45) and (2.4.46),
which comes from the transformation equations of elliptic functions of
order thirteen given by (4.3.19). More precisely, the Jacobian multiplier
equation (4.3.8), the modular equation formed by $\delta_{\nu}$ $(\nu=0, 1,
\ldots, 12)$ and $\delta_{\infty}$ and (4.3.19) have the same Galois group
which is isomorphic to $\text{SL}(2, \mathbb{F}_{13})$.

  Now we explain a representation-theoretic reformulation of the weight
$k$ in Serre's conjecture (see \cite{GHS}). The Eichler-Shimura isomorphism
allows one to reinterpret Serre's conjecture in terms of the cohomology of
arithmetic groups. If $V$ is an $\overline{\mathbb{F}}_p$-representation
of $\text{GL}(2, \mathbb{F}_p)$ and $N$ is prime to $p$, then we have a
natural action of the Hecke algebra of $\Gamma_1(N)$ on $H^1(\Gamma_1(N), V)$,
and so it makes sense to speak of a continuous representation
$\overline{\rho}: G_{\mathbb{Q}} \rightarrow \text{GL}(2, \overline{\mathbb{F}}_p)$
being associated to an eigenclass in that cohomology group. If $\overline{\rho}$
is odd and irreducible, then the Eichler-Shimura isomorphism implies that
$\overline{\rho}$ is modular of weight $k$ and prime-to-$p$ level $N$ if and
only if $\overline{\rho}$ is associated to an eigenclass in
$H^1(\Gamma_1(N), \text{Sym}^{k-2} \overline{\mathbb{F}}_p^2)$, where
$\text{Sym}^{k-2} \overline{\mathbb{F}}_p^2$ is the $(k-2)$-th symmetric
power of the standard representation of $\text{GL}(2, \mathbb{F}_p)$ on
$\overline{\mathbb{F}}_p^2$. By d\'{e}vissage one deduces that $\overline{\rho}$
is modular of weight $k$ and prime-to-$p$ level $N$ if and only if $\overline{\rho}$
is associated to an eigenclass in $H^1(\Gamma_1(N), V)$ for some
Jordan-H\"{o}lder factor $V$ of $\text{Sym}^{k-2} \overline{\mathbb{F}}_p^2$.
Recall that the representation $\text{Sym}^{k-2} \overline{\mathbb{F}}_p^2$
is reducible as soon as $k>p+1$. It is then natural to associate to $\overline{\rho}$
the set $W(\overline{\rho})$ of irreducible $\overline{\mathbb{F}}_p$-representations
$V$ of $\text{GL}(2, \mathbb{F}_p)$ such that $\overline{\rho}$ is associated to an
eigenclass in $H^1(\Gamma_1(N), V)$ for some prime-to-$p$ level $N$. Such
representations of $\text{GL}(2, \mathbb{F}_p)$ are often referred to as
Serre weights. To give a concrete example, let $\chi$ denote the
cyclotomic character. Suppose that
$$\overline{\rho}|_{I_{\mathbb{Q}_p}} \cong \left(\begin{matrix}
  \chi^{k-1} & *\\
  0 & 1
  \end{matrix}\right),$$
where $I_{\mathbb{Q}_p}$ is the inertia group at $p$, and $2<k<p-1$.
Then the minimal weight predicted by Serre's recipe is $k$. If the
extension class $*$ is non-split, then we have $W(\overline{\rho})=\{
\text{Sym}^{k-2} \overline{\mathbb{F}}_p^2 \}$. In fact, let $k \geq 2$
and let $f$ be a normalised cuspidal eigenform of level $N$ prime to
$p$. Then $\rho_f|_{G_{\mathbb{\mathbb{Q}}_p}}$ is crystalline. Moreover,
any crystalline representation with Hodge-Tate weights $\{ 0, k-1 \}$
which is reducible mod $p$, must be of this form. However, in the
companion forms case where the extension class $*$ is split, we
have
$$W(\overline{\rho})=\{ \text{Sym}^{k-2} \overline{\mathbb{F}}_p^2,
  \text{det}^{k-1} \otimes \text{Sym}^{p-1-k} \overline{\mathbb{F}}_p^2 \}.$$
The decomposition formulas in Theorem 1.1 and Corollary 1.2 can be regarded
as an anabelian counterpart of Serre weights:
$$\begin{matrix}
 &\text{Serre weights} &\longleftrightarrow
 &\text{Our decomposition formulas for}\\
 &\text{$W(\overline{\rho})$} & &\text{$\pi_p \leftrightarrow \rho_p$
                                       in Theorem 1.1 and Corollary 1.2}
\end{matrix}$$

  One of the significance of these representations of
$\text{GL}(2, \mathbb{F}_p)$ is that they are the starting point for
constructing $\mathbb{F}_p$-representations of $\text{GL}(2, \mathbb{Z}_p)$
and $\text{GL}(2, \mathbb{Q}_p)$. These representations play a very
important role in the $p$-adic local Langlands correspondence. Let
$E$ denote a finite extension of $\mathbb{Q}_p$ with residue field $k_E$.
We identify the center of $\text{GL}(2, \mathbb{Q}_p)$ with
$\mathbb{Q}_p^{\times}$. If $r \geq 0$, then $\text{Sym}^r k_E^2$ is
a representation of $\text{GL}(2, \mathbb{F}_p)$ that yields a representation
of $\text{GL}(2, \mathbb{Z}_p)$ by inflation, and we extend it to
$\text{GL}(2, \mathbb{Z}_p) \mathbb{Q}_p^{\times}$ by letting $p$ act
via the identity. The representation
$$\text{Ind}_{\text{GL}(2, \mathbb{Z}_p) \mathbb{Q}_p^{\times}}^{\text{GL}(2, \mathbb{Q}_p)}
  \text{Sym}^r k_E^2$$
is the set of functions $f: \text{GL}(2, \mathbb{Q}_p) \rightarrow \text{Sym}^r k_E^2$
that are locally constant, have a compact support modulo $\text{GL}(2, \mathbb{Z}_p)
\mathbb{Q}_p^{\times}$, and satisfy $f(kg)=\text{Sym}^r(k) f(g)$ if
$k \in \text{GL}(2, \mathbb{Z}_p) \mathbb{Q}_p^{\times}$ and $g \in
\text{GL}(2, \mathbb{Q}_p)$. This space is equipped with the action
of $\text{GL}(2, \mathbb{Q}_p)$ given by $(gf)(h)=f(hg)$. This leads
to the $p$-adic local Langlands correspondence for $\text{GL}(2, \mathbb{Q}_p)$
which is a bijection between some two-dimensional representations of
$\text{Gal}(\overline{\mathbb{Q}}_p/\mathbb{Q}_p)$ and some representations
of $\text{GL}(2, \mathbb{Q}_p)$. This bijection can in fact be constructed
using the theory of $(\varphi, \Gamma)$-modules and some results of $p$-adic
analysis (see \cite{Berger} for more details).

\begin{center}
{\large\bf 7.4. Relation with the Fontaine-Mazur-Langlands conjecture}
\end{center}

  Let us recall the Fontaine-Mazur-Langlands conjecture (see \cite{FM}):
{\it Let}
$$\rho: \text{Gal}(\overline{\mathbb{Q}}/\mathbb{Q}) \rightarrow
  \text{GL}(2, \overline{\mathbb{Q}}_p)$$
{\it be an irreducible representation which is unramified except at a finite
number of primes and which is not the Tate twist of an even representation
which factors through a finite quotient group of}
$\text{Gal}(\overline{\mathbb{Q}}/\mathbb{Q})$.

  {\it Then $\rho$ is associated to an cuspidal newform $f$ if and only if
$\rho$ is potentially semi-stable at $p$.}

  Here, ``associated'' is a technical term which means that one can find an
isomorphism between $\overline{\mathbb{Q}}_p$ and $\mathbb{C}$ such that for
all but a finite number of prime numbers $\ell$, the above isomorphism brings
the trace of $\rho(\text{Frob}_{\ell})$ to the eigenvalues of the Hecke
operator $T_{\ell}$ on the newform $f$. Here, $\text{Frob}_{\ell}$ denotes
the Frobenius element at $\ell$, and of course, one need only consider prime
numbers $\ell$ which do not divide the level of $f$, and at which the
representation $\rho$ is unramified.

  Fontaine and Mazur also gave some refinement of their conjecture. Let us
recall the Conjecture 3c in \cite{FM}: {\it Let $V$ be an irreducible geometric
$\overline{\mathbb{Q}}_p$-representation of}
$\text{Gal}(\overline{\mathbb{Q}}/\mathbb{Q})$ {\it of degree two which is
not a Tate twist of a finite representation. Then there is an integer $i \in
\mathbb{Z}$ such that $V(-i)$ is isomorphic to the representation associated
to a ``new'' modular form.}

  Note that if $X$ is an elliptic curve over $\mathbb{Q}$, if $T_p(X)$ is its
Tate module and if $V=\mathbb{Q}_p \otimes_{\mathbb{Z}_p} T_p(X)$, $V$ is
geometric, conjecture 3c implies the Taniyama-Shimura-Weil conjecture.

  In the context of motives, the Fontaine-Mazur-Langlands conjecture
leads to a conjectural necessary and sufficient condition, stated only
in terms of a local condition at $p$ (i.e., geometric), for an irreducible
representation
$\rho: \text{Gal}(\overline{\mathbb{Q}}/\mathbb{Q}) \rightarrow \text{GL}(2,
 \overline{\mathbb{Q}}_p)$
to come from algebraic geometry, and hence to be the representation associated
to a cuspidal newform. However, the Fontaine-Mazur-Langlands conjecture
excludes the case of finite Galois representations. In contrast with it,
in the context of anabelian algebraic geometry, we prove that the finite
$\mathbb{Q}(\zeta_p)$-rational Galois representations $\rho_p$ given by
(7.2.1) can come from the defining ideals of modular curves:
$\rho_p: \text{Gal}(\overline{\mathbb{Q}}/\mathbb{Q}) \rightarrow
 \text{Aut}(\mathcal{L}(X(p)))$,
i.e., come from anabelian algebraic geometry. Moreover, Theorem 1.3 shows
that it is modular, i.e., $\rho_p$ corresponds to the elliptic modular
functions.

\begin{center}
{\large\bf 8. An anabelian counterpart of the $p$-adic local Langlands
              correspondence for $\text{GL}(2, \mathbb{Q}_p)$}
\end{center}

  We will have a look back at the origins of the $p$-adic Langlands program
(see \cite{Breuil2008}). Let us fix a prime number $p$. We will describe
the background that leads to the discovery of the first instances of the
$p$-adic local correspondence for $\text{GL}(2, \mathbb{Q}_p)$. Let us fix
embeddings $\overline{\mathbb{Q}} \hookrightarrow \mathbb{C}$ and
$\overline{\mathbb{Q}} \hookrightarrow \overline{\mathbb{Q}}_p$. Let $F$
be a number field, $\mathbb{A}_F$ the locally compact group of adeles of $F$,
and $n \geq 1$ an integer. We normalise local reciprocity maps such that the
geometric Frobenius map to the inverse of the uniformisers. The following
classic conjecture is attributed to Langlands, Fontaine, Mazur and Tate:

\textbf{Conjecture 8.1.} {\it There is a unique bijection between the two sets:}
$$\{ \text{\it isomorphism classes of algebraic cuspidal automorphic}$$
$$\text{{\it representations $\pi=\otimes^{\prime} \pi_{\mathfrak{l}}$ of}
  \text{$\text{GL}(n, \mathbb{A}_F)$}}\}$$
$$\updownarrow$$
$$\{\text{\it isomorphism classes of continuous irreducible representations}$$
$$\text{$\rho: \text{Gal}(\overline{\mathbb{Q}}/F) \rightarrow \text{GL}(n, \overline{\mathbb{Q}}_p)$
  {\it which are unramified at almost all places}}$$
$$\text{\it and potentially semistable at places dividing $p$}\}$$
{\it such that, for any finite place $\mathfrak{l} | \ell$, $\ell \neq p$,
$\pi_{\mathfrak{l}}$ and} $\rho|_{\text{Gal}(\overline{\mathbb{Q}}_{\ell}/F_{\mathfrak{l}})}$
{\it determine one another via the $($appropriately normalised$)$ local Langlands
correspondence.}

  It is recalled that an automorphic representation of $\text{GL}(n, \mathbb{A}_F)$
is algebraic if, for every infinite place, the restriction to $\mathbb{C}^{\times}$
of the corresponding Langlands parameter is a direct sum of characters
$z \mapsto z^{a-\frac{n-1}{2}} \overline{z}^{b-\frac{n-1}{2}}$ where
$a$, $b \in \mathbb{Z}$. In fact, the correspondence
$\pi_{\mathfrak{l}} \longleftrightarrow \rho|_{\text{Gal}(\overline{\mathbb{Q}}_{\ell}/F_{\mathfrak{l}})}$
above involves the Weil-Deligne representation, conjectured to be $F$-semisimple,
associated with the $\ell$-adic representation
$\rho|_{\text{Gal}(\overline{\mathbb{Q}}_{\ell}/F_{\mathfrak{l}})}$.
For $\ell \neq p$, this Weil-Deligne representation determines
$\rho|_{\text{Gal}(\overline{\mathbb{Q}}_{\ell}/F_{\mathfrak{l}})}$.
However, for places $\mathfrak{p}$ dividing $p$, it has long been known
that $\pi_{\mathfrak{p}}$ generally does not determine
$\rho|_{\text{Gal}(\overline{\mathbb{Q}}_{p}/F_{\mathfrak{p}})}$, but
only a small part: the Weil-Deligne representation associated with
$\rho|_{\text{Gal}(\overline{\mathbb{Q}}_{p}/F_{\mathfrak{p}})}$.
However, $\pi$ does indeed determine the global representation $\rho$,
and therefore, a fortiori,
$\rho|_{\text{Gal}(\overline{\mathbb{Q}}_{p}/F_{\mathfrak{p}})}$, which
means that whatever is missing from $\pi_{\mathfrak{p}}$ to determine
$\rho|_{\text{Gal}(\overline{\mathbb{Q}}_{p}/F_{\mathfrak{p}})}$ is still
somewhere in $\pi$. The $p$-adic Langlands programme stems from a desire
to understand what needs to be added to $\pi_{\mathfrak{p}}$ in order to
reconstruct $\rho|_{\text{Gal}(\overline{\mathbb{Q}}_{p}/F_{\mathfrak{p}})}$.
In other words, we wish to elucidate the emergence of $p$-adic Hodge theory
on the Galois side in terms of representation theory on the automorphic side.
Here is a slightly less vague way of formulating this objective:

\textbf{Question 8.2.} Can we associate with $\pi$ a ``natural $p$-adic''
representation $\widehat{\pi}=\otimes^{\prime} \widehat{\pi}_{\mathfrak{l}}$
of $\text{GL}(n, \mathbb{A}_F)$ such that, for all finite places $\mathfrak{l}$:
$\widehat{\pi}_{\mathfrak{l}}$ and
$\rho|_{\text{Gal}(\overline{\mathbb{Q}}_{\ell}/F_{\mathfrak{l}})}$
determine one another?

  ``Natural'' means, for example, that $\widehat{\pi}$ or its factors
$\widehat{\pi}_{\mathfrak{l}}$ realise in suitable cohomology spaces.
For $\mathfrak{l} \nmid p$, we can already take $\widehat{\pi}_{\mathfrak{l}}
=\pi_{\mathfrak{l}}$, so the main task is to define $\widehat{\pi}_{\mathfrak{p}}$
for $\mathfrak{p} | p$. More precisely, we therefore not only want
$\widehat{\pi}_{\mathfrak{p}}$ to determine $\pi_{\mathfrak{p}}$, but
also to contain the additional information needed to fully reconstruct
the $p$-adic representation
$\rho_{\mathfrak{p}}=\rho|_{\text{Gal}(\overline{\mathbb{Q}}_p/F_{\mathfrak{p}})}$.
From a local perspective, we are therefore looking for new correspondence
$\widehat{\pi}_{\mathfrak{p}} \longleftrightarrow \rho_{\mathfrak{p}}$.
If $\widehat{\pi}$ exists, the hope is that the correspondence
$\pi \longleftrightarrow \rho$ in Conjecture 8.1 factorialises naturally
via $\pi \longrightarrow \widehat{\pi} \longleftrightarrow \rho$,
where, in the latter case, the local representations determine one
another place by place, and where the analysis at infinite places
on the automorphic side is replaced by $p$-adic analysis at places
dividing $p$.

  Now, let us restrict to $\text{GL}(2)$. From the classical theory,
we know that when the level of $f$ is exactly divisible by $p$, the
local automorphic component $\pi_p(f)$ does not contain enough information
to isolate $\rho_p$. Indeed, in this case, $\pi_p(f)$ is always the
Steinberg representation and it is exactly the $\mathcal{L}$ invariant
of $f$ that provide the additional information necessary to identify
the local Galois representation associated to $f$. More precisely,
let $p$ be a prime number and $f$ a cuspidal newform of weight
$k \geq 2$ on $\Gamma_0(Np)$ with $(N, p)=1$. For simplicity, let
$f$ be an eigenvector of the Hecke operators with eigenvalues in
$\mathbb{Z}$. Let $\rho$ be the $p$-adic representation of
$\text{Gal}(\overline{\mathbb{Q}}/\mathbb{Q})$ associated with $f$;
we know that the local component
$\rho_p=\rho|_{\text{Gal}(\overline{\mathbb{Q}}_p/\mathbb{Q}_p)}$
is a non-crystalline semistable representation. Thus, via the Fontaine
functor, it corresponds to a filtered $(\varphi, N)$-module $D_{\text{st}}(\rho_p)$
of dimension $2$, and the invariant $\mathcal{L}$ of $f$, denoted by
$\mathcal{L}(f)$, is defined by the parameter of the Hodge filtration on
$D_{\text{st}}(\rho_p)$ (see \cite{Breuil2004}). The form $f$ is also
associated with an automorphic representation
$\pi=\otimes_{\ell}^{\prime} \pi_{\ell}$ of $\text{GL}(2, \mathbb{A}_{\mathbb{Q}})$,
where $\mathbb{A}_{\mathbb{Q}}$ denotes the adeles of $\mathbb{Q}$. The local
component $\pi_p$ is, up to an unramified torsion, the special representation--
or Steinberg representation--of $\text{GL}(2, \mathbb{Q}_p)$, denoted by $\text{St}$.
This is an irreducible smooth representation of $\text{GL}(2, \mathbb{Q}_p)$
in which there is no $\mathcal{L}$-invariant whatsoever. For this reason,
whilst $\mathcal{L}(f)$ is a local invariant on the Galois side of $p$, it
has the privilege of being regarded as a global invariant on the $\text{GL}(2)$
side. This motivation is linked to the following open question which is closely
related to Question 8.2: is there a Langlands correspondence (global or local)
between continuous $p$-adic representations of Galois groups and certain
continuous $p$-adic representations of linear groups? Indeed, if such a
correspondence exists, then in our particular case the hypothetical $p$-adic
representation $\Pi_p$ of $\text{GL}(2, \mathbb{Q}_p)$ associated with $\rho_p$
must determine $\rho_p$ (up to isomorphism) and must therefore contain a
structure equivalent to the Hodge filtration on $D_{\text{st}}(\rho_p)$.
In other words, the Hodge-Tate weights $(0, k-1)$ and the $\mathcal{L}$-invariant
must lie somewhere within the local representation $\Pi_p$.

  In order to understand the above theory, we need Fontaine's theory on $p$-adic
Galois representations (see \cite{Mezard}). A $p$-adic Galois representation is
a continuous homomorphism
$$\rho: \text{Gal}(\overline{\mathbb{Q}}/\mathbb{Q}) \longrightarrow \text{GL}(2, E)$$
where $E$ is a finite extension of $\mathbb{Q}_{\ell}$. For any prime number
$p$, by restricting $\rho$ to $\text{Gal}(\overline{\mathbb{Q}}_p/\mathbb{Q}_p)$,
we obtain
$$\rho_p: \text{Gal}(\overline{\mathbb{Q}}_p/\mathbb{Q}_p) \longrightarrow \text{GL}(2, E).$$
For $\ell \neq p$, we can describe all the representations $\rho_p$. But for
$\ell=p$, there are too many $p$-adic representations. Fontaine's idea is to
focus first on those that arise naturally in geometry: semistable representations.
(Fontaine defined rings of period $B_{\text{dR}}$ and $B_{\text{st}}$ equipped
with certain structures, including an action of
$G_{\mathbb{Q}_p}=\text{Gal}(\overline{\mathbb{Q}}_p/\mathbb{Q}_p)$.
The $p$-adic representation $V$ is said to be de Rham if
$\dim_{\mathbb{Q}_p} (B_{\text{dR}} \otimes_{\mathbb{Q}_p} V)^{G_{\mathbb{Q}_p}}
 =\dim_{\mathbb{Q}_p} V$. The $p$-adic representation $V$ is said to be
semistable if $\dim_{\mathbb{Q}_p} (B_{\text{st}} \otimes_{\mathbb{Q}_p} V)^{G_{\mathbb{Q}_p}}
 =\dim_{\mathbb{Q}_p} V$.) We obtain such representations in particular

(1) by taking the \'{e}tale cohomology of a proper, smooth variety over
   $\mathbb{Q}_p$ with semistable reduction;

(2) by constructing the $p$-adic representation associated with a modular
form for $\Gamma_0(pN)$, $(p, N)=1$.

  Thanks to Fontaine's theory, we can characterize semistable representations.
They are, up to isomorphism, parameterized by one or two invariants: the
Hodge-Tate weights $(0, k-1)$ and the invariant $\mathcal{L}$, which is in
$\mathbb{Q}_p$. Specifically,

\textbf{Theorem 8.3.} (Colmez-Fontaine) {\it There is a category equivalence}:
$$\left\{\begin{matrix}
  \text{{\it semistable representations}}\\
  \text{Gal}(\overline{\mathbb{Q}}_p/\mathbb{Q}_p) \longrightarrow \text{GL}(2, \mathbb{Q}_p)
\end{matrix}\right\} \longleftrightarrow
  \left\{\begin{matrix}
  \text{{\it weakly admissible filtered}}\\
  \text{$(\varphi, N)$-{\it modules of dimension two}}
\end{matrix}\right\}.$$

  The parameters $k$ and $\mathcal{L}$ determine the filtration of weakly
admissible $(\varphi, N)$-modules. There are two types of semistable
representations: crystalline representations and non-crystalline semistable
representations. Each of them is associated with the weakly admissible
$(\varphi, N)$-module $D$.

  The local $\ell$-adic Langlands correspondence was proven for dimensions
$n \geq 2$ by Harris and Taylor (see \cite{HT}) in 1998. Henniart (see
\cite{Henniart}) provided a simplified proof of this in 2000. It establishes
the correspondence between
$$\aligned
  \left\{\begin{matrix}
  \text{admissible irreducible representations}\\
  \text{of $\text{GL}(n, \mathbb{Q}_p)$ over $\overline{\mathbb{Q}}_{\ell}$}
\end{matrix}\right\} &\longleftrightarrow
  \left\{\begin{matrix}
  \rho_{\ell}: W_p \longrightarrow \text{GL}(n, \overline{\mathbb{Q}}_{\ell})\\
  \text{irreducible, $\phi$-semisimple}
\end{matrix}\right\}\\
  \pi_{\ell} &\longleftrightarrow \rho_{\ell}
\endaligned$$
where $W_p$ is the Weil subgroup of $\text{Gal}(\overline{\mathbb{Q}}_p/\mathbb{Q}_p)$.
The proof provides a geometric interpretation of the correspondence when
the representations arise from modular forms. The method is global in the
sense that we treat all prime numbers simultaneously by considering

(1) the representations of $\text{GL}(2, \mathbb{A})$ where $\mathbb{A}$ is the
restricted product subring of $\mathbb{R} \times \prod \mathbb{Q}_{\ell}$
(these global representations are called automorphic),

(2) the Galois representations of the global group $\text{Gal}(\overline{\mathbb{Q}}/\mathbb{Q})$.

  Let us clarify the geometric construction strategy. This leads to the
appearance of $\text{SL}(2, \mathbb{F}_p)$ and its relation with
$\text{SL}(2, \mathbb{Z}_p)$ and $\text{SL}(2, \mathbb{Q}_p)$. To an
eigenform $f$, i.e., a cuspidal modular form of weight $k \geq 2$ on
$\Gamma_0(N)$ that is an eigenvector of the Hecke operators, we can
associate an $\ell$-adic representation
$$\rho(f): \text{Gal}(\overline{\mathbb{Q}}/\mathbb{Q}) \longrightarrow
  \text{GL}(2, \overline{\mathbb{Q}}_{\ell})$$
and, for any prime $p$, an irreducible smooth representation $\pi_p(f)$
of $\text{GL}(2, \mathbb{Q}_p)$ (over $\overline{\mathbb{Q}}_{\ell}$)
such that if $p \neq \ell$, $\pi_p(f)$ can be reconstructed from
$\rho(f)|_{\text{Gal}(\overline{\mathbb{Q}}_p/\mathbb{Q}_p)}$ and
vice-versa:
$$\pi_p(f) \longleftrightarrow \rho(f)|_{\text{Gal}(\overline{\mathbb{Q}}_p/\mathbb{Q}_p)}.$$
This correspondence is realized geometrically on a cohomology group (see
\cite{Mezard}): for simplicity, we assume here that $N=p^n$ is a power of $p$.
Let
$$Y(p^n)(\mathbb{C}):=\Gamma(p^n) \backslash \mathbb{H} \quad \text{where}
  \quad \Gamma(p^n):=\text{ker}(\text{SL}(2, \mathbb{Z}) \rightarrow
  \text{SL}(2, \mathbb{Z}/p^n \mathbb{Z})).$$
We have a surjection
$$Y(p^{n+1})(\mathbb{C}) \rightarrow Y(p^n)(\mathbb{C}),$$
and $\text{SL}(2, \mathbb{Z}/p^n \mathbb{Z})$ acts on $Y(p^n)(\mathbb{C})$
via $\overline{z} \mapsto \overline{gz}$, where $z$ is a conjugate of
$\overline{z}$ in $\mathbb{H}$ and $g \in \text{SL}(2, \mathbb{Z})$ (this
action does not depend on the choices since $\Gamma(p^n)$ is distinguished).
Consider
$$\varinjlim_{n} H_{\text{sing}}^1(Y(p^n)(\mathbb{C}), \overline{\mathbb{Q}}_{\ell}),$$
where $H_{\text{sing}}^1$ is the singular cohomology of the topological
space $Y(p^n)(\mathbb{C})$ with coefficients in $\overline{\mathbb{Q}}_{\ell}$.
We can show that this space is naturally equipped with an action of
$\text{SL}(2, \mathbb{Q}_p)$. For example, the action of $\text{SL}(2, \mathbb{Z}_p)$
arises from the action of $\text{SL}(2, \mathbb{Z}/p^n \mathbb{Z})$ at
each level. In particular, when $n=1$, i.e., the level $p$, it gives
the action of $\text{SL}(2, \mathbb{F}_p)$. It is also equipped with
a continuous action of $\text{Gal}(\overline{\mathbb{Q}}/\mathbb{Q})$,
but this is much more subtle: this stems from the fact that $Y(p^n)(\mathbb{C})$
corresponds to the complex points of an algebraic variety defined over
$\mathbb{Q}$, and that in this case we can replace
$H_{\text{sing}}^1(Y(p^n)(\mathbb{C}), \overline{\mathbb{Q}}_{\ell})$
with the $\ell$-adic \'{e}tale cohomology of $Y(p^n)(\overline{\mathbb{Q}})$,
on which there is a Galois group action. In practice, we prefer an action
of $\text{GL}(2, \mathbb{Q}_p)$ to an action of $\text{SL}(2, \mathbb{Q}_p)$,
and at each level we replace $Y(p^n)(\mathbb{C})$ with
$|(\mathbb{Z}/p^n \mathbb{Z})^{*}|$ copies of $Y(p^n)(\mathbb{C})$ on which
we let $(\mathbb{Z}/p^n \mathbb{Z})^{*}$ act by permutation. Let us denote
this new space by $\mathbb{Y}(p^n)(\mathbb{C})$; this time, we have an action
of $\text{GL}(2, \mathbb{Q}_p) \times \text{Gal}(\overline{\mathbb{Q}}/\mathbb{Q})$
on
$$\varinjlim_{n} H_{\text{sing}}^1(\mathbb{Y}(p^n)(\mathbb{C}), \overline{\mathbb{Q}}_{\ell}).$$
We can then state the geometric realization theorem of the Langlands
correspondence for $k=2$:

\textbf{Theorem 8.4.} (Deligne, Langlands, Carayol) {\it There exists a}
$\text{GL}(2, \mathbb{Q}_p) \times \text{Gal}(\overline{\mathbb{Q}}/\mathbb{Q})$-{\it equivariant 
isomorphism}
$$\varinjlim_{n} H_{\text{sing}}^1(\mathbb{Y}(p^n)(\mathbb{C}), \overline{\mathbb{Q}}_{\ell})
  \simeq \oplus^{\prime} \pi_p(f) \otimes \rho(f) \oplus (\text{remainder})$$
{\it where the sum $\oplus^{\prime}$ is over cuspidal eigen-newforms
$f \in S_2(p^{*})$ and the $($re-mainder$)$ term relates to non-cuspidal
forms.}

  For weights $k>2$, the constant coefficients $\overline{\mathbb{Q}}_{\ell}$
must be replaced by the local system arising from the fibre bundle
$\Gamma(p^n) \backslash (\text{Sym}^{k-2} \overline{\mathbb{Q}}_{\ell}^2
 \times \mathbb{H}) \longrightarrow \Gamma(p^n) \backslash \mathbb{H}$.
Now, we consider the $p$-adic case. For $n \in \{ 0, 1 \}$, a Galois
representation associated with an eigen-newform of $S_k(\Gamma_0(p^n))$
is semistable. However, for $n \geq 2$, such a representation is only
de Rham. Suppose that the representation $\rho_p(f)$ arises from a
cuspidal eigen-newform $f$ of $S_k(\Gamma_0(p^n))$. Let us consider
$\widehat{H_{\text{sing}}^{1}}$, the $p$-adic completion of
$$\varinjlim_{n} H_{\text{sing}}^1(\mathbb{Y}(p^n)(\mathbb{C}), \mathbb{Z}_p),$$
which is equipped with a continuous action of $\text{GL}(2, \mathbb{Q}_p)
\times \text{Gal}(\overline{\mathbb{Q}}/\mathbb{Q})$. It is unitary,
i.e., there exists an invariant form for $\text{GL}(2, \mathbb{Q}_p)$.
Let $\Pi_p(f)$ be the unitary $p$-adic Banach space equipped with an
action of $\text{GL}(2, \mathbb{Q}_p)$
$$\Pi_p(f)=\text{Hom}_{\text{Gal}(\overline{\mathbb{Q}}/\mathbb{Q})}
  (\rho(f), \widehat{H_{\text{sing}}^{1}}).$$
On the other hand, there always exists a $\text{GL}(2, \mathbb{Q}_p)$-equivariant
injection
$$\underline{\text{Sym}}^{k-2} E^2 \otimes \pi_p(f) \longrightarrow \Pi_p(f).$$
Note that we have the algebraic representation $\text{sym}^{(k-2)} \overline{\mathbb{Q}}_p^2$
of the highest weight $(0, k-2)$. For technical reasons, $k$ is assumed to be
even, and we will instead work with $\underline{\text{Sym}}^{k-2} \overline{\mathbb{Q}}_p^2
=\text{Sym}^{k-2} \overline{\mathbb{Q}}_p^2 \otimes |\det|_p^{\frac{k-2}{2}}$,
whose central character is the $(k-2)$-th power of the cyclotomic character.
In fact, $\underline{\text{Sym}}^{k-2} E^2 \otimes \pi_p(f)$ is identified
with the locally algebraic vectors of $\Pi_p(f)$. Let $\widehat{\pi_p(f)}$
be the completion of $\underline{\text{Sym}}^{k-2} E^2 \otimes \pi_p(f)$ in
$\Pi_p(f)$. Thus $\widehat{\pi_p(f)}$ is the completion of
$\underline{\text{Sym}}^{k-2} E^2 \otimes \pi_p(f)$ induced by the
invariant norm on $\widehat{H_{\text{sing}}^1}$.

\textbf{Conjecture 8.5.} (Breuil) {\it If
$\rho(f)|_{\text{Gal}(\overline{\mathbb{Q}}_p/\mathbb{Q}_p)}=\rho_p(f)$
is absolutely irreducible, then $\widehat{\pi_p(f)}$ is non-zero, topologically
irreducible, admissible, determines $\rho_p(f)$ and depends only on it.}

  Note that when $\rho_p(f)$ is reducible, $\widehat{\pi_p(f)}$ is not
always sufficient to determine $\rho_p(f)$. When $\rho_p(f)$ is absolutely
irreducible, Breuil also conjectures that $\Pi_p(f)=\widehat{\pi_p(f)}$.
Conjecture 8.5 is proven in the following cases: The case of an eigen
newform $f$ of $S_k(\Gamma_0(p))$ (in which case $\rho_p(f)$ is semistable
non-crystalline and $\pi_p(f)$ is said to be a special series) is established:
Breuil shows that there exists a unique non-zero morphism from $B(k, \mathcal{L})$
to $\widehat{\pi_p(f)}$ if and only if $\mathcal{L}=\mathcal{L}(f)$. The
case of an eigen newform $f$ of $S_k(\Gamma_0(1))$ (in which case $\rho_p(f)$
is crystalline and $\pi_p(f)$ is said to be a principal series) has been
studied by Berger and Breuil. In this case, there is only one local unitary
completion $B(k)$ of $\underline{\text{Sym}}^{k-2} E^2 \otimes \pi_p$,
which is therefore necessarily induced by the global cohomology. The case
of an eigen newform $f$ of $S_k(\Gamma_0(p^n))$ with $n \geq 2$ (in which
case $\pi_p(f)$ is said to be supercuspidal) has not yet been considered.

  This leads to a comparison between Fontaine's theory on $p$-adic Galois
representations and our theory on $\mathbb{Q}(\zeta_p)$-rational Galois
representations by the intermediate of representation theory:
$$\begin{matrix}
 &\text{Fontaine's theory} &\text{representation theory} &\text{our theory}\\
 &\text{de Rham} &\text{supercuspidal representations} & \rho_{11, 10}\\
 &\text{representations} &\text{(cuspidal or discrete series)} & \text{(see (2.4.31) et seq.)}\\
 &\text{semistable} &\text{Steinberg representations} & \rho_{13, 13}\\
 &\text{representations} & \text{(special series)} & \text{(see (2.4.42) et seq.)}\\
 &\text{crystalline} &\text{principal series} & \rho_{13, 7}\\
 &\text{representations} &\text{(degenerate principal series)} & \text{(see (2.4.42) et seq.)}\\
\end{matrix}\eqno{(8.1)}$$
As a consequence, for the simplest case of level $p$, our theory provides
a unified geometric realization for all of the three kinds of representations:
``crystalline'', ``semistable'' and ``de Rham'' representations.

  Fix a prime number $p$. Let $K^p=\prod_{\ell \neq p} K_{\ell}$ be an
open compact subgroup of $\text{GL}(2, \mathbb{A}_f^p)$. Set
$$\widetilde{H}^i(K^p, \mathbb{Z}/p^n):=\varinjlim_{K_p \subseteq
  \text{GL}(2, \mathbb{Q}_p)} H^i(Y_{K^p K_p}(\mathbb{C}), \mathbb{Z}/p^n),$$
where $K_p$ runs through all open compact subgroups of $\text{GL}(2, \mathbb{Q}_p)$,
and $Y_{K^p K_P}$ denotes the modular curve of level $K^p K_p$ over $\mathbb{Q}$
whose $\mathbb{C}$-points are given by the usual quotient $Y_{K^p K_p}(\mathbb{C})
=\text{GL}(2, \mathbb{Q}) \backslash (\mathbb{H}^{\pm} \times \text{GL}(2, \mathbb{A}_f)/K)$,
and $\mathbb{H}^{\pm}=\mathbb{C}-\mathbb{R}$ denotes the union of the upper and lower
half-planes. The completed cohomology of tame level $K^p$, introduced by Emerton, is
defined as
$$\widetilde{H}^i(K^p, \mathbb{Z}_p):=\varprojlim_{n} \widetilde{H}^i(K^p, \mathbb{Z}/p^n).$$
It is $p$-adically complete and equipped with a natural continuous action of
$\text{GL}(2, \mathbb{Q}_p) \times G_{\mathbb{Q}}$. Emerton was interested
in the finite dimensional Galois representation of $G_{\mathbb{Q}}$ inside
of $\widetilde{H}^1(K^p, \mathbb{Z}_p)$ (see \cite{Emer2006} or \cite{Pan}).

\textbf{Theorem 8.6.} {\it Let $E$ be a finite extension of $\mathbb{Q}_p$ and}
$$\rho: G_{\mathbb{Q}} \rightarrow \text{GL}(2, E)$$
{\it be a two-dimensional continuous absolutely irreducible representation
of $G_{\mathbb{Q}}$. Suppose that}

(1) {\it $\rho$ appears in $\widetilde{H}^1(K^p, E):=\widetilde{H}^1(K^p,
\mathbb{Z}_p) \otimes_{\mathbb{Z}_p} E$, i.e.,}
$$\Pi_{\rho}:=\text{Hom}_{E[G_{\mathbb{Q}}]}(\rho, \widetilde{H}^1(K^p, E)) \neq 0.$$

(2) {\it $\rho|_{G_{\mathbb{Q}_p}}$ is de Rham of Hodge-Tate weights $0$, $k$ for
some integer $k>0$, where $G_{\mathbb{Q}_p} \subseteq G_{\mathbb{Q}}$ denotes a
decomposition group at $p$.}

{\it Then $\rho$ arises from a cuspidal eigenform of weight $k+1$.}

  This result was first obtained by Emerton in \cite{Emer2006} (see
also \cite{Pan} for a different proof) and was used by him to study
the Fontaine-Mazur conjecture. To do this, Emerton proved a local-global
compatibility result of the completed cohomology which implies that
$\Pi_{\rho}$ corresponds to $\rho|_{G_{\mathbb{Q}_p}}$ under the $p$-adic
local Langlands correspondence for $\text{GL}(2, \mathbb{Q}_p)$. Since $\rho$
appears in $\widetilde{H}^1(K^p, E)$, one can show that $\rho|_{G_{\mathbb{Q}_{\ell}}}$
is unramified for all but finitely many $\ell$ and $\rho|_{G_{\mathbb{R}}}$
is odd. Hence theorem 8.6 is a special case of the Fontaine-Mazur conjecture.
In fact, in Emerton's approach to the Fontaine-Mazur conjecture \cite{Emer2006},
he first showed that any two-dimensional $p$-adic geometric absolutely
irreducible odd Galois representation of $G_{\mathbb{Q}}$ appears in
$\widetilde{H}^1(K^p, E)$, then he invoked this result to deduce that
it actually comes from a cuspidal eigenform.

  Let us recall some facts about the relation between the completed cohomolgy
and the Hodge-Tate decomposition appearing in (2) of Theorem 8.6 (see
\cite{Pan}). Fix $k \geq 1$. The starting point is a characterization of de
Rham Galois representations in Fontaine's classification of almost
$B_{\text{dR}}$-representations. More precisely, let $V$ be a two-dimensional
continuous representation of $G_{\mathbb{Q}_p}$ over $\mathbb{Q}_p$. Suppose
that $V$ is Hodge-Tate of weights $0$, $k$, i.e., there is a natural
decomposition
$$W:=V \otimes_{\mathbb{Q}_p} C=W_0 \oplus W_k$$
with both $W_0$ and $W_k$ being non-zero and $W_i(i)=W_i(i)^{G_{\mathbb{Q}_p}}
\otimes_{\mathbb{Q}_p} C$ for $i=0,$, $k$, where $C$ is the $p$-adic completion
of $\overline{\mathbb{Q}}_p$. Fontaine proved that there is a natural $C$-linear,
$G_{\mathbb{Q}_p}$-equivariant operator $N: W_0 \rightarrow W_k(k)$ such that
$V$ is de Rham exactly when $N=0$. Let $\Pi_{\rho}^{\text{la}} \subseteq
\Pi_{\rho}$ be the subspace of $\text{GL}(2, \mathbb{Q}_p)$-locally analytic
vectors. Then $\Pi_{\rho}^{\text{la}}$ has an infinitesimal character
$\widetilde{\chi}_k: Z(U(\mathfrak{gl}_2(\mathbb{Q}_p))) \rightarrow
\mathbb{Q}_p$ which is equal to the infinitesimal character of the
$(k-1)$-th symmetric power of the dual of the standard representation.
In fact, the $\widetilde{\chi}_k$-isotypic part
$\widetilde{H}^1(K^p, \mathbb{Q}_p)^{\text{la}, \widetilde{\chi}_k}$ is
Hodge-Tate of weights $0$, $k$, i.e., there is a natural Hodge-Tate
decomposition
$$\widetilde{H}^1(K^p, \mathbb{Q}_p)^{\text{la}, \widetilde{\chi}_k}
  \widehat{\otimes}_{\mathbb{Q}_p} C=W_0 \oplus W_k$$
such that $W_i(i)=W_i(i)^{G_{\mathbb{Q}_p}} \widehat{\otimes}_{\mathbb{Q}_p} C$,
$i=0$, $k$. Fontaine's construction can also be generalized to this
setting: there exists a natural continuous map $N: W_0 \rightarrow
W_k(k)$ such that for any two-dimensional irreducible subrepresentation
$V \subseteq \widetilde{H}^1(K^p, \mathbb{Q}_p)^{\text{la}, \widetilde{\chi}_k}$
of $G_{\mathbb{Q}}$, the restriction $N|_{(V \otimes_{\mathbb{Q}_p} C)_0}$
agrees with the Fontaine operator $(V \otimes_{\mathbb{Q}_p} C)_0 \rightarrow
(V \otimes_{\mathbb{Q}_p} C)_k(k)$ for $V$.

  As an anabelian counterpart of Theorem 8.6, Theorem 1.1 and
Theorem 1.3 imply the following

\textbf{Theorem 8.7.} {\it Let} $\rho_p: G_{\mathbb{Q}} \rightarrow
\text{PSL}(2, \mathbb{F}_p)$ {\it be a Galois representation. Suppose
that}

(1) {\it $\rho_p$ appears in the defining ideals of modular curves $X(p)$.}

(2) {\it $\rho_p|_{G_{\mathbb{Q}}}$ has a decomposition as given by Theorem 1.1.}

{\it Then $\rho_p$ is modular in the sense of Theorem 1.3.}

  In particular, we have the following comparison:
$$\begin{array}{ccc}
  \text{motives}  & \longleftrightarrow & \text{anabelian geometry}\\
  \text{completed cohomology} & \longleftrightarrow & \text{defining ideals}\\
  \text{of modular curves} &  & \text{of modular curves}\\
  \text{$p$-adic local Langlands} & \longleftrightarrow &\text{our correspondence}\\
  \text{correspondence for $\text{GL}(2, \mathbb{Q}_p)$} &
  & \text{for $\text{PSL}(2, \mathbb{F}_p)$}\\
  \text{two-dimensional} & \longleftrightarrow & \text{higher-dimensional}\\
  \text{$p$-adic} & & \text{$\mathbb{Q}(\zeta_p)$-rational}\\
  \text{Galois representation} &  & \text{Galois representation}\\
  \rho: G_{\mathbb{Q}} \rightarrow \text{GL}(2, \mathbb{Q}_p) & &
  \rho_p: G_{\mathbb{Q}} \rightarrow \text{PSL}(2, \mathbb{F}_p)\\
  \text{$\rho$ appears in completed} & \longleftrightarrow & \text{$\rho_p$ appears in defining}\\
  \text{cohomology of modular curves} & & \text{ideals of modular curves}\\
  \text{$\rho|_{G_{\mathbb{Q}_p}}$ is de Rham of} & \longleftrightarrow &
  \text{$\rho_p|_{G_{\mathbb{Q}}}$ has a decomposition}\\
  \text{Hodge-Tate weights $0$, $k$} & & \text{in Theorem 1.1}\\
  \text{Hodge-Tate decomposition} & \longleftrightarrow & \text{semisimple decomposition}\\
  \text{with nilpotent operator $N$} &  & \text{with $\text{PSL}(2, \mathbb{F}_p)$-action}\\
  \text{$\rho$ arises from a cuspidal} & \longleftrightarrow & \text{$\rho_p$ is modular in the sense}\\
  \text{eigenform of weight $k+1$} & & \text{of Theorem 1.3}\\
  \text{GL}(2, \mathbb{Q}_p) \times G_{\mathbb{Q}}\text{-action} & \longleftrightarrow &
  \text{PSL}(2, \mathbb{F}_p) \times G_{\mathbb{Q}}\text{-action}
\end{array}\eqno{(8.2)}$$

\begin{center}
{\large\bf 9. Comparison with the finite Langlands correspondence (the 
              Macdonald correspondence)}
\end{center}

  In 1980, Macdonald (see \cite{Macdonald}) proved an analogous result of
the local Langlands correspondence of $\text{GL}(n)$ for finite fields.
Namely, he established an explicit bijection between the isomorphism
classes of the irreducible representations of $\text{GL}(n, k)$, where
$k$ is a finite field, and inertia equivalence classes of admissible tamely
ramified $n$-dimensional Weil-Deligne representations of $W_K$, where
$K$ is a non-archimedean local field with residue field $k$ and $W_K$ the
absolute Weil group of $K$. More precisely, let $\Pi(\text{GL}(n, k))$ denote
the set of equivalence classes of irreducible complex representations of
$\text{GL}(n, k)$, and let $\Phi(\text{GL}(n))$ denote the set of equivalence
classes of $n$-dimensional admissible representations of $W_K$. On the other
hand, let $I$ and $P$ be respectively the inertia and ramification subgroups
of the Weil group $W_K$ of $K$, and let $\Phi^{t}(\text{GL}(n))$ denote the
set of $\phi \in \Phi(\text{GL}(n))$ which are tamely ramified, i.e. trivial
on $P$. Two elements of $\Phi^{t}(\text{GL}(n))$ will be called $I$-equivalent
if their restrictions to $I$ are equivalent. Let $\Phi_I^{t}(\text{GL}(n))$
denote the set of $I$-equivalent classes. Macdonald established a canonical
bijection of $\Pi(\text{GL}(n, k))$ onto $\Phi_I^{t}(\text{GL}(n))$, by showing
that both of theses sets are parameterized by the same set $P_n(\Gamma)$.

  More precisely, let $k$ be a finite field with $q$ elements,
$\overline{k}$ an algebraic closure of $k$, and $F: x \rightarrow x^q$
the Frobenius automorphism of $\overline{k}$ over $k$. For each $n \geq 1$
let $k_n=\{ x \in \overline{k}: F^n x=x \}$ be the unique extension of $k$
of degree $n$ in $\overline{k}$, and let $k_n^{*}$ be its multiplicative
group. If $m$ divides $n$, the norm homomorphism
$N_{n, m}: k_n^{*} \rightarrow k_m^{*}$ (namely $x \rightarrow x^{(q^n-1)/(q^m-1)}$)
is surjective. Let $\Pi(k_n^{*})$ be the character group of $k_n^{*}$. By
transposition, $N_{n, m}$ defines an injective homomorphism $\Pi(k_m^{*})
\rightarrow \Pi(k_n^{*})$ whenever $m$ divides $n$. Let $\Gamma=\varinjlim
\Pi(k_n^{*})$ be the direct limit; $\Gamma$ is a discrete group in which every
element has finite order, and is (non-canonically) isomorphic to the multiplicative
group of $\overline{k}$. The Frobenius $F$ acts on $\Gamma$: namely $F \gamma=\gamma^q$.
If $f$ is an $F$-orbit in $\Gamma$, we denote by $d(f)$ the number of elements in
$f$, and we call $d(f)$ the degree of $f$. Let $f_1$, $\ldots$, $f_m$ be distinct
$F$-orbits in $\Gamma$, and let $\lambda_1$, $\ldots$, $\lambda_m$ be partitions
such that $\sum |\lambda_i| d(f_i)=n$, where $|\lambda_i$ is the sum of the parts
of the partition $\lambda_i$. Then $(f_1^{\lambda_1}) \circ \cdots \circ (f_m^{\lambda_m})$
is an irreducible representation of $\text{GL}(n, k)$, and all the irreducible
representations of $\text{GL}(n, k)$ are obtained in this way, without repetition.
Macdonald restated this result as follows. Let $P_n(\Gamma)$ denote the set of all
partition-valued functions $\lambda$ on $\Gamma$ such that (i) $\lambda \circ F=\lambda$,
i.e., $\lambda$ is constant on the $F$-orbits; (ii) $\sum_{\lambda \in \Gamma}
|\lambda(\gamma)|=n$. Then there exists a canonical bijection $\lambda \mapsto
\pi_{\lambda}$ of $P_n(\Gamma)$ onto $\Pi(\text{GL}(n, k))$. On the other hand,
Macdonald established a canonical bijection $\lambda \mapsto (\phi_{\lambda})$
of the set $P_n(\Gamma)$ onto $\Phi_I^{t}(\text{GL}(n))$, and hence a canonical
bijection:
$$\aligned
  \Pi(\text{GL}(n, k)) &\rightarrow \Phi_I^t(\text{GL}(n))\\
  \pi_{\lambda} &\mapsto (\phi_{\lambda}).
\endaligned\eqno{(9.1)}$$

  In general, suppose $G$ is a reductive group defined over $\mathbb{F}_q$.
Deligne-Lusztig (\cite{DeLu}) described irreducible representations of
$G(\mathbb{F}_q)$. More precisely, let $G$ be a connected, reductive
algebraic group over an algebraic closure $\mathbb{F}$ of $\mathbb{F}_q$,
with an $\mathbb{F}_q$-rational structure and Frobenius morphism
$F: G \rightarrow G$. They constructed virtual representations over
$\overline{\mathbb{Q}}_{\ell}$ (here $\ell$ is a prime not dividing $q$)
of the finite group $G^F$ of $F$-fixed points of $G$, via the $\ell$-adic
cohomology of algebraic varieties over $\mathbb{F}$ which are now known as
Deligne-Lusztig varieties. In \cite{Lu1979} and \cite{Lu2004}, Lusztig gave
a cohomological construction for certain representations of $G(R)$ over a
finite ring $R$ arising from the ring of integers in a non-archimedean local
field by reduction modulo a power of the maximal ideal, which is an extension
of the construction of the virtual representations in \cite{DeLu} for groups
over a finite field to a local field. Therefore, a natural question arises
(see \cite{Vogan2019} and \cite{Vogan2020}): Can their results be formulated
in spirit of local Langlands correspondence as above? More precisely, what is
the Langlands correspondence say about $\widehat{G(\mathbb{F}_q)}$? Vogan
(\cite{Vogan2019} and \cite{Vogan2020}) gave some conjectural approach in the
context of local Langlands correspondence.

  Recently, Imai and Vogan (see \cite{ImaiVogan} and \cite{Imai}) formulated
a parametrization of representations of a finite Chevalley group $G(\mathbb{F}_q)$
over $\overline{\mathbb{Q}}_{\ell}$ in terms of the Langlands dual group over
$\overline{\mathbb{Q}}_{\ell}$ of $G$. One motivation is to relate this
parametrization for $G(\mathbb{F}_q)$ to the (still conjectural) Langlands
parametrization of irreducible representations of groups over local fields.
The reason for using the dual group over $\overline{\mathbb{Q}}_{\ell}$ is
that Langlands' philosophy suggests that representations of a reductive group
$G$ on $k=\overline{k}$ vector spaces ought to be related to (maps of a Weil
group to) a $k$-dual group. In particular, such a parametrization was given
for representations of $\text{GL}(n, \mathbb{F}_q)$ by Macdonald as above.

  We give a comparison between the Macdonald correspondence (9.1) in the case
of $n=2$, $k=\mathbb{F}_p$ and Theorem 1.1 as follows:
$$\begin{array}{ccc}
  \text{the Macdonald correspondence} & \longleftrightarrow & \text{our correspondence}\\
  \pi_{\lambda} \in \Pi(\text{GL}(n, k)) & \longleftrightarrow & \pi_p\\
  \text{complex representations} & \longleftrightarrow &
  \text{$\mathbb{Q}(\zeta_p)$-rational representations}\\
  (\phi_{\lambda}) \in \Phi_I^t(\text{GL}(n)) & \longleftrightarrow & \rho_p\\
  \text{Weil-Deligne representations} & \longleftrightarrow & \text{Galois representations}\\
  \text{the set $P_n(\Gamma)$} & \longleftrightarrow & \text{the defining ideal $I(\mathcal{L}(X(p)))$}\\
  \text{both are parameterized} & \longleftrightarrow & \text{both are parameterized}\\
  \text{by the same $P_n(\Gamma)$} & & \text{by the same $I(\mathcal{L}(X(p)))$}\\
\end{array}\eqno{(9.2)}$$

\begin{center}
{\large\bf 10. An ideal theoretic counterpart of Grothendieck's section conjecture
               and an ideal theoretic reification of ``arithmetic theory of $\pi_1$''}
\end{center}

  In his paper \cite{Weil1938}, Weil predicted the non-abelian mathematics
in the future. He pointed out that ``Some of the most brilliant advances in
modern mathematics have been made, as we know, in arithmetic, the theory of
algebraic varieties and topology. Now, not only are these fields closely
related, and by analogies that we have not finished exploring and exploiting,
but most of the progress made is marked by a common imprint, namely an
essentially abelian character. Whether we are talking about the class field,
multiple integrals on algebraic varieties, or the homology properties of a
space, commutative groups are used everywhere, and it is most often this
restriction that leads to success. Perhaps the reason for this lies in
duality, so that Pontrjagin's theory should be placed at the centre of all
this research. Riemann's work on algebraic functions is already deeply
marked by this character, and it is no accident that abelian functions
bear the same name as abelian groups, since the division of these functions
leads (as Abel, Galois and Jacobi saw) to abelian equations, and their main
algebraic use is to generate and study the relatively abelian extensions of
fields of algebraic functions of one variable. Without wishing to give other
examples, I would like to recall the main contribution to modern topology is
indeed properties of homology, i.e. abelian properties. This abelian mathematics
is not yet complete, it lacks, for example, the general notion of the multiple
integral, which, when it is found, will perhaps be its main tool. But it already
constitutes an imposing field of doctrine, to which considerable treatises are
devoted. On the other hand, as soon as we want to go further in the fields
listed above, we find ourselves enveloped in darkness; the truth is that we
have already begun to study some special problems, which I shall not attempt
to enumerate; in topology especially, questions relating to homotopy have
been tackled, not without success. But there is still room for a great deal
of research before we can hope to reach an overall view; and this is why I
have undertaken to examine, from the non-abelian point of view, the theory
of algebraic functions of one variable, in the hope also that a generalisation
of theta functions might provide arithmetic with a powerful instrument.''

  Note that Weil emphasized the importance of the division of abelian
functions and the non-abelian viewpoint of the theory of algebraic functions
of one variable. As we have shown in the section 6, in particular (6.6), the
division equations of the elliptic functions lead to the Galois representations
associated with the division points of elliptic curves or abelian varieties,
i.e., the Jacobian varieties, which is of abelian nature. On the other hand,
the transformation equations of elliptic functions lead to the Galois representations
arising from defining ideals of modular curves $X(p)$, which are of anabelian
nature. In particular, the defining ideals of modular curves provide a
non-abelian viewpoint of the theory of algebraic functions of one variable.

  In fact, the success of Weil's thesis arose from the use of abelian varieties,
and in particular of Jacobians; for the Mordell conjecture, Weil felt that it
would be necessary to go beyond the abelian setting (see \cite{Weil1938},
\cite{Se1999} or \cite{Groth1958}). In his paper \cite{Weil1938}, Weil pointed
out the importance of developing non-abelian mathematics, whose ingredients
should involve the moduli of vector bundles and the fundamental groups.
As pointed out by Serre \cite{Se1999}, this paper was ``a text presented
as analysis, whose significance is essentially algebraic, but whose
motivation is arithmetic''. Weil discussed the need to move beyond abelian
objects in the study of arithmetic. He began the theory of vector bundles
over algebraic curves and initiated the study of vector bundles in this
context to generalize the Jacobian. The homological nature of the Jacobian
was emphasized and the footnote contained allusion to the importance of
non-abelian $\pi_1$. From modern perspective, moduli of semi-stable bundles
corresponds to the theory of reductive completions of $\pi_1$. Weil thought
that such theories should have application to arithmetic. It is plausible
that he had the Mordell conjecture in his mind. However, at the time of
Weil, there was no serious arithmetic theory of $\pi_1$ (see also
\cite{Kim2006}). Weil tried to use this parametrization coming from Jacobian
(i.e., homological construction) to prove the Mordell conjecture. However,
it does not quite work because the Jacobian has too many points which is
also related to the linearity. The eventual proof of Mordell conjecture
involves a non-linear parametrization, that is, a family of curves (the
Kodaira-Parshin construction) (see \cite{Faltings} and \cite{Kim2005April}).
According to \cite{Lang1983}, as a curve of genus $\geq 1$ is characterized
by the fact that it is embeddable in an abelian variety, one is led to study
subvarieties of abelian varieties, and especially Mordell's conjecture, in
this light. A curve of genus $\geq 2$ is characterized by the fact that it
is unequal to any non-zero translation of itself in its Jacobian. This is a
very difficult hypotheses to use. In fact, for a general variety $X$, the
category of motives, being of homological nature (i.e., abelian nature),
destroys the information about $X(\mathbb{Q})$ and does not touch on some
very basic Diophantine phenomena. These abelian invariants do not yield in
general information about $X(\mathbb{Q})$ and leave thereby untouched the
most basic questions of Diophantine geometry. This is an artifact of the
fact that the theory of motives as presently developed is implicitly modeled
on the theory of abelian varieties and $H_1$. We have an important contrast:
$$\text{linear} \longleftrightarrow \text{nonlinear}$$
as well as their realization in topology:
$$\text{homology} \longleftrightarrow \text{homotopy},$$
where
$$\text{homology} \longleftrightarrow \text{the study of $L$-functions},$$
$$\text{homotopy} \longleftrightarrow \text{the study of Diophantine sets}.$$
Hence, one expects a fully nonlinear theory, when fully understood, to be
more powerful than a linearized version (see \cite{Kim2005April}).

  In order to remedy this, in a letter to Faltings \cite{GtoF}, Grothendieck
proposed to study rational points on a curve $X$ through the geometric
\'{e}tale fundamental group $\pi_1^{\text{\'{e}t}}(\overline{X}, b)$ of
$X$ with base point $b$. More precisely, the section conjecture of
Grothendieck states that the non-abelian Albanese map
$$\aligned
  \kappa^{\text{na}}: X(\mathbb{Q}) & \rightarrow
  H^1(\text{Gal}(\overline{\mathbb{Q}}/\mathbb{Q}),
  \pi_1^{\text{\'{e}t}}(\overline{X}, b))\\
  x & \mapsto [\pi_1^{\text{\'{e}t}}(\overline{X}; b, x)]
\endaligned\eqno{(10.1)}$$
which attaches to every rational point $x \in X(\mathbb{Q})$ the class of
the Galois representation defined by the corresponding \'{e}tale path torsor
$\pi_1^{\text{\'{e}t}}(\overline{X}; b, x)$ of the algebraic fundamental group,
should be an isomorphism. Here, the topological space
$H^1(\text{Gal}(\overline{\mathbb{Q}}/\mathbb{Q}), \pi_1^{\text{\'{e}t}}(\overline{X}, b))$
is the non-abelian continuous cohomology of the absolute Galois group
$\text{Gal}(\overline{\mathbb{Q}}/\mathbb{Q})$ of $\mathbb{Q}$ with
coefficients in the profinite \'{e}tale fundamental group of $X$. The
notation will suggest that a rational base-point $b \in X(\mathbb{Q})$
has been introduced. The context makes it clear that
$H^1(\text{Gal}(\overline{\mathbb{Q}}/\mathbb{Q}), \pi_1^{\text{\'{e}t}}(\overline{X}, b))$
can be understood as a non-abelian Jacobian in an \'{e}tale profinite
realization, where the analogy can be strengthened by the interpretation
of the $\text{Gal}(\overline{\mathbb{Q}}/\mathbb{Q})$-action as defining a
sheaf on $\text{Spec}(\mathbb{Q})$ and $H^1(\text{Gal}(\overline{\mathbb{Q}}/\mathbb{Q}), \pi_1^{\text{\'{e}t}}(\overline{X}, b))$ as the moduli space of torsors
for $\pi_1^{\text{\'{e}t}}(\overline{X}, b)$ in the \'{e}tale topos of
$\text{Spec}(\mathbb{Q})$. It is an analogue of the moduli space
$\text{Bun}_n(X)$ of rank $n$ vector bundles on $X$ for $n \geq 2$.
Their study was initiated by Weil (see \cite{Weil1938} as well as
\cite{Groth1958} and \cite{Se1999}) who regarded them also as non-abelian
Jacobians. In other words, Grothendieck's section conjecture states
that every torsor of the fundamental group should necessarily arise
from a rational point. This provides us with the possibility of
studying the set of such torsors in lieu of the set $X(\mathbb{Q})$.
However, the cohomology set that classifies these torsors does not
seem to have much structure with which we can work.

  The motivic fundamental group lies between the profinite fundamental
group and the homology group:
$$\begin{array}{cccccc}
  &\widehat{\pi}_1(\overline{X}, b) &\longleftrightarrow &\pi_1^{M}(\overline{X}, b)
  &\longleftrightarrow &H_1(\overline{X})\\
  &\text{anabelian geometry} &\longleftrightarrow &\text{intermediate theory}
  &\longleftrightarrow &\text{motives}\\
  &\text{Grothendieck} &\longleftrightarrow &\text{Deligne}
  &\longleftrightarrow &\text{Grothendieck}\\
\end{array}\eqno{(10.2)}$$
Instead of the above non-abelian Albanese map and as a replacement
of Grothendieck's section conjecture mentioned as above, Kim (see
\cite{Kim2005}) considered a motivic Albanese map
$$\kappa^{M}: X(\mathbb{Q}) \rightarrow H_M^1(\text{Gal}(\overline{\mathbb{Q}}/\mathbb{Q}),
  \pi_1^{M}(\overline{X}, b))\eqno{(10.3)}$$
that associates to points motivic torsors $\pi_1^M(\overline{X}; b, x)$
of paths. The fundamental group refers to the unipotent motivic fundamental
group in the sense of Deligne \cite{De1989}. It should be pointed out that
the technology of motives ends up contributing here as well because the
Diophantine aspect of this theory \cite{Kim2005} assigns an interesting
role to motivic fundamental groups \cite{De1989}, where Ext groups are
replaced by classifying spaces for non-abelian torsors. In fact, the
motivic $\pi_1$ is still quite close to homology, and hence, admits
homological techniques. But it is non-abelian. As an application, Kim gave
a motivic-$\pi_1$ proof of Siegel's theorem (see also \cite{Faltings2007})
on the finiteness of integral points for the thrice-punctured projective
line by using unipotent Galois representations associated to the unipotent
\'{e}tale fundamental group of projective line minus three points and
comparing the unipotent \'{e}tale and de Rham fundamental groups, where
$\pi_1$ denotes $\mathbb{Q}_p$-unipotent completions. The technical heart
of Kim's method in \cite{Kim2005} is to define \'{e}tale and de Rham period
maps and relate them using (non-abelian) $p$-adic Hodge theory. Moreover,
\'{e}tale and de Rham realizations of the unipotent fundamental group and
path torsors admit some extra structures also play key roles in the later
work (see \cite{Hadian} based on \cite{DeGon}). In particular, Kim pointed
out in \cite{Kim2005} that ``it is rather striking that the motivic theory
is capable of yielding Diophantine finiteness, even though the motivic
fundamental group could be viewed as a cruder invariant of a variety than
the profinite one. The $L$-functions, viewed as invariants of varieties,
seem to provide information in general only about linearized invariants,
e.g., Chow groups. It is then natural that results about nonlinear sets
should evoke nonlinear tools like the fundamental group. What is somewhat
surprising is that even a mild degree of nonlinearity (coming from the
subcategory of unipotent group objects in the category of motives) can
still provide substantial information.'' The motivic fundamental groups
are connected with the theory of mixed motives, which is related to
Diophantine geometry via $L$-functions. It can give information about
linearized Diophantine invariants, i.e., conjectures of Birch and
Swinnerton-Dyer, Deligne, Beilinson, Bloch-Kato, \ldots. Hence, it can
be regarded as a linearized section conjecture. However, the motivic
philosophy does not cover many nonlinear phenomena. So we still need a
fully nonlinear theory which is far away from the motivic philosophy.

  For modular curves $X=X(p)$, we have the following three distinct
kinds of geometric realizations by motives, motivic fundamental groups
and defining ideals (which can determine the fundamental groups),
respectively. In particular, they are abelian, nilpotent and utmost
nonabelian, respectively. At the same time, they are linear, mildly
nonlinear and fully nonlinear, respectively. In particular, as a
replacement of Grothendieck's section conjecture as well as the above
linearized section conjecture, we will use the following correspondence:
$$X(p)/\mathbb{Q} \longleftrightarrow \text{the defining ideal
  $I(\mathcal{L}(X(p)))/\mathbb{Q}$},\eqno{(10.4)}$$
which can be considered a fully nonlinear counterpart of Grothendieck's section
conjecture.

  Here, Chabauty's original method fits naturally into the setting of motivic
unipotent \'{e}tale fundamental groups and Selmer varieties, i.e., the
Chabauty-Kim method (rather than the Jacobians) as the technical foundation
of the analytic part now becomes non-abelian $p$-adic Hodge theory (see
\cite{Olsson}) and iterated integrals (see \cite{Kim2005}, \cite{Kim2009},
\cite{Kim2010}, \cite{Kim2012}, \cite{CoKim}, \cite{BDMTV} and \cite{BBBLMTV}).
Moreover, it relies on the Iwasawa theory (both commutative and non-commutative),
as Kim pointed out in the end of \cite{Kim2012}: ``Allowing ourselves a further
flight of fancy, the elusive function in general might eventually be the subject
of an Iwawasa theory rising out of a landscape radically more non-abelian and
non-linear than we have dared to dream''. This also leads to another proof
of the Mordell conjecture (see \cite{LV}, \cite{Poonen} and \cite{BBBLMTV}).
The abelian varieties are replaced by a more detailed analysis of the variation
of $p$-adic Galois representations in a family of algebraic varieties. The
key inputs into this analysis are the comparison theorems of $p$-adic Hodge
theory and explicit topological computations of monodromy.

  Now, we give some interpretation of the comparison (1.13). The defining
ideals of modular curves with a utmost degree of nonlinearity, which can
determine the fundamental groups of these curves, do provide substantial
information which is distinct from that the motivic fundamental group can
give. In particular, our theory (Theorem 1.1, Corollary 1.2 and Theorem
1.3) provides an arithmetic theory of $\pi_1$ as well as a nonlinear and
anabelian version of the Iwasawa theory, which can also be regarded as a
nonlinear version of Weil's construction. At the same time, our theory is
a fully nonlinear theory. In the case of motivic fundamental groups, they
give some Diophantine finiteness. In our context, three different perspectives
on non-abelian phenomena in number theory: the Langlands program and motives,
non-commutative Iwasawa theory and anabelian algebraic geometry fit together
in the following unified way: the anabelianization of the global Langlands
correspondence for $\text{GL}(2, \mathbb{Q})$ by the \'{e}tale cohomology
of modular curves gives rise to our correspondence for $\text{GL}(2, \mathbb{Q})$
by the defining ideals of modular curves; the non-abelianization of
Kummer-Herbrand-Ribet theorem for ideal class groups of cyclotomic fields
$\mathbb{Q}(\zeta_p)$ (i.e., cyclotomic Iwasawa theory) and anabelianization
of the theory of Kubert-Lang and Mazur-Wiles for the cuspidal divisor class
groups of modular curves $X(p)$ as well as Fukaya-Kato-Sharifi's theory
for homology relative to the cusps of the modular curve (i.e., cuspidal
Iwasawa theory associated with modular curves) give rise to Prasad-Shekhar's
theory for the eigenspace decomposition of the $p$-primary part of the ideal
class groups for division fields $\mathbb{Q}(E[p])$ of elliptic curves $E$
over $\mathbb{Q}$ (i.e., elliptic Iwasawa theory) and our theory for the
eigenspace of the space $V_{I(\mathcal{L}(X(p)))}$ associated with the
defining ideals $I(\mathcal{L}(X(p)))$ of modular curves $X(p)$ (i.e.,
anabelian Iwasawa theory associated with modular curves), respectively.
Finally, our defining ideals of modular curve gives a concrete and an
explicit realization of the non-abelian Jacobian appearing in Weil's
paper \cite{Weil1938} as well as a replacement of Grothendieck's
section conjecture (see \cite{GtoF} and \cite{Groth1958}).

\begin{center}
{\large\bf 11. An anabelian counterpart of the theory of Kubert-Lang and Mazur-Wiles
              on the cuspidal divisor class groups and the Eisenstein ideals}
\end{center}

  The Mazur-Wiles proof (\cite{MW}) of the Iwasawa main conjecture
for cyclotomic fields indicates that the geometry of modular curves
near the cusps has much to say about the arithmetic of cyclotomic
fields. That is, the Eisenstein ideal determining the congruences
between cusp forms and Eisenstein series is essentially the annihilator
of the cusp at $\infty$. So, when we look at residual reducible representations
attached to the cusp forms, we are essentially exploring the geometry
of a modular curve near $\infty$ (see \cite{Sharifi2018}). In the
present paper, we will show that the geometry of the full modular curves
has an analogue and a connection with the arithmetic of division fields
of elliptic curves over $\mathbb{Q}$. That is, the defining ideal is
essentially the annihilator of the locus of a modular curve. Thus, when
we look at Galois representations arising from the defining ideals, we
are essentially exploring the geometry of a full modular curve. Note that
the cusp forms come from $H_1$, i.e., the motives, on the other hand, a
full modular curve is closely related to $\pi_1$, i.e., the anabelian
algebraic geometry.

  For (1) in the correspondence (1.18), it is known that the cuspidal divisor
class group in the Jacobian variety of modular curves play an important role
in the theory of the arithmetic of modular curves. In a series of papers (see
\cite{Kubert-Lang}), Kubert and Lang gave a striking analogy between the
cuspidal divisor class groups of modular curves and the ideal class groups of
cyclotomic fields, which led into the deeper connections established subsequently
by Mazur (see \cite{Mazur1977}) on his proof of Ogg's conjecture as well as
Wiles (see \cite{Wiles1980}) and Mazur-Wiles (see \cite{MW}) on the proof of
the main conjecture of Iwasawa theory for cyclotomic fields. They found
that the Stickelberger element which annihilates ideal classes in cyclotomic
fields is related to and has analogies with certain Hecke operators (i.e.,
the Eisenstein ideals) which annihilate the divisor classes at infinity
on the modular curves. In particular, in the case of cyclotomic fields,
we have the class number formula which involves the first generalized
Bernoulli number $B_{1, \chi}$ (see Theorem 2.6 in \cite{Lang}):
$$h^{-}=w \prod_{\text{$\chi$ odd}} -\frac{1}{2} B_{1, \chi}.\eqno{(11.1)}$$
In the case of modular curves $X_1(p)$, we have the cuspidal divisor class
number formula which involves the second generalized Bernoulli number
$B_{2, \chi}$ (see Theorem 3.1 in \cite{Lang}):
$$h=p \prod_{\chi \neq 1} \pm \frac{1}{2} B_{2, \chi}.\eqno{(11.2)}$$
Similarly, in the case of modular curves $X(p)$, we have the cuspidal
divisor class number formula which also involves the second generalized
Bernoulli number $B_{2, \chi}$ (see p. 118, Theorem 3.1 in \cite{Kubert-Lang}):
$$h=\frac{6 p^3}{|G|} \prod_{\chi \neq 1} \frac{p}{2} B_{2, \chi},\eqno{(11.3)}$$
where $G=C(p)/ \pm 1$ for the Cartan group $C(p)$.

  Now, we find an analogue between the defining ideals of modular curves
and the division fields of elliptic curves, which can be regarded as both
a non-abelian and an anabelian counterpart of the above analogue. In particular,
the two generalized Bernoulli numbers appear as above are replaced by the
division equations and the transformation equations of elliptic functions,
respectively. At first, we will give a comparison between Eisenstein ideals
and defining ideals of modular curves (i.e., the second comparison between
motives and anabelian algebraic geometry). It is well-known that (see
\cite{Mazur1977}) an Eisenstein ideal annihilates the cuspidal divisor class
group which is the torsion subgroup of rational points on the Jacobian variety
of a modular curve (i.e, an object coming from $H_1$). On the other hand, a
defining ideal annihilates the locus of a modular curve (i.e., an object
coming from $\pi_1$).

  In the context of abelian side, the Herbrand-Ribet theorem establishes
a connection between (a) (the arithmetic side) the structure of the action
of $\text{Gal}(\mathbb{Q}(\zeta_p)/\mathbb{Q}) \cong (\mathbb{Z}/p \mathbb{Z})^{\times}$
on $\text{Cl}(\mathbb{Q}(\zeta_p))/p \text{Cl}(\mathbb{Q}(\zeta_p))$ and (b)
(the analytic side) the divisibility (or non-divisibility) by $p$ of the
numerator of certain Bernoulli numbers, by the $\chi$-eigen-spaces of
$(\mathbb{Z}/p \mathbb{Z})^{\times}$. In the context of non-abelian side,
our theory establishes a connection between (a) (the arithmetic side) the
structure of the action of
$\text{Gal}(K/\mathbb{Q}) \cong \text{PSL}(2, \mathbb{F}_p)$ on the defining
ideals of modular curves $X(p)$ and (b) (the analytic side) the transformation
equations of elliptic functions of order $p$, by the eigen-spaces associated
with $\text{PSL}(2, \mathbb{F}_p)$, i.e., some $\mathbb{Q}(\zeta_p)$-rational
irreducible representations of $\text{PSL}(2, \mathbb{F}_p)$, such as discrete
series, principal series and the Steinberg representations. One may summarize
Ribet's idea (see \cite{Ribet1976}) by saying that an Eisenstein congruence
relates the geometry (i.e. $H_1$) of $\text{GL}(2, \mathbb{Q})$ (e.g. cuspidal
eigen-forms or cuspidal divisor class groups) to the arithmetic of
$\text{GL}(1, \mathbb{Q})$ (e.g. ideal class groups of cyclotomic fields).
Correspondingly, our theory shows that there is an analogue or connection
between the anabelian algebraic geometry (i.e. $\pi_1$) of $\text{GL}(2, \mathbb{Q})$
(e.g. defining ideals of modular curves) and the arithmetic of $\text{GL}(2, \mathbb{Q})$
(e.g. ideal class groups of division fields of elliptic curves over $\mathbb{Q}$).
In particular, for the same modular curves, there are two distinct kinds of
realizations by arithmetic geometry (i.e., two distinct kinds of arithmetical
structures):

(1) realization by $H_1$ (motives): cuspidal divisor class groups in the Jacobian
   varieties of modular curves and the Eisenstein ideals;

(2) realization by $\pi_1$ (anabelian algebraic geometry): modular curves and
    their defining ideals.

  In Chapter 7 of \cite{Hida}, Hida pointed out that Iwasawa's theory for
cyclotomic fields (i.e., the so-called cyclotomic Iwasawa theory in \cite{Hida})
properly belongs to the automorphic theory for $\text{GL}(1)$, the Iwasawa
algebra (the base ring of the cyclicity) can be defined as a universal
deformation ring of $\text{GL}(1)$-representations and also as a Hecke
algebra for $\mathbf{G}_m=\text{GL}(1)$. The theory of Kubert-Lang (i.e.,
the so-called cuspidal Iwasawa theory in \cite{Hida}) is a bridge from the
$\text{GL}(1)$-theory to the $\text{GL}(2)$-theory, and the Hecke operators
acting on the cuspidal class group span the Eisenstein component of the
$\text{GL}(2)$ Hecke algebra and this fact is essential in the proof of
the main conjecture by Mazur-Wiles (see \cite{MW}). It is natural to
expect a cyclicity theory for $\text{GL}(2)$ to have $\mathbf{T}$ as its
base ring. Our theory provides such an answer in the context of anabelian
algebraic geometry as well as non-commutative Iwasawa theory. In particular,
the $\text{PSL}(2, \mathbb{F}_p)$-action on the defining ideal $I(\mathcal{L}(X(p)))$
of modular curves $X(p)$, which is an anabelian counterpart of the Hecke
operators acting on the cuspidal class group of modular curves, is a genuine
$\text{GL}(2)$-theory, For a number field $K$, we have the class group
$\text{Cl}(K)$ and its $p$-primary part $C_K$. According to \cite{Hida}, p. 237,
there are three incarnations of $C_K$: as the $p$-primary part of the class
group (field arithmetic), as the Galois group of the maximal abelian unramified
extension (Galois theory), and as the dual of a Selmer group (homology theory).
In its geometric counterpart, there is a cuspidal divisor class group (homology
theory). Our theory shows that there is a new interpretation as a defining
ideal of the modular curve (anabelian algebraic geometry), which goes beyond
the homology theory. In particular, the dimension of the defining ideal
$I(\mathcal{L}(X(p)))$ plays a similar role as the cuspidal divisor class
number.

  The proof of the Iwasawa Main Conjecture for cyclotomic fields (see \cite{MW},
\cite{Wiles1980} as well as \cite{Ribet1976}) arose from the use of abelian
varieties, and in particular of Jacobians. More precisely, the approach in
\cite{MW} was inspired by Ribet's proof of the converse of a theorem of
Kummer-Herbrand in that they used the structure of certain finite groups
of torsion points on abelian varieties arising as quotients of the Jacobian
varieties of some modular curves. An important role is played by the cuspidal
divisor class group whose structure is related to Stickelberger ideals,
and hence to Bernoulli numbers as well as $p$-adic $L$-functions, by
results of Kubert and Lang \cite{Kubert-Lang}.

  In the most classical situation, the Galois group of the cyclotomic
field $\mathbb{Q}(\mu_N)$ of roots of unity operates on ideal classes
and units modulo cyclotomic units. Classical problems of number theory
are concerned with the eigenspace decomposition of the $p$-primary part
of these groups when $N=p$, and of the structures as Galois modules in
general. Results include those of Kummer, Stickelberger, Herbrand, and
more recently Iwasawa, Leopoldt, Ribet, Mazur and Wiles.

  Kubert and Lang (see \cite{Kubert-Lang}) generalized Stickelberger's
theory of cyclotomic class groups to the setting of modular curves.
They observed that the Stickelberger ideal also arises naturally in
the study of the cusps on the modular curves. This theory gives a
base of the proof of the Iwasawa main conjecture by Mazur-Wiles (see
\cite{MW}). In \cite{Kubert-Lang}, Kubert and Lang were concerned
with the generic case of the modular function field of the modular
curve $X(N)$. The divisor class group generated by the cusps can be
represented as a quotient of the group ring of the Cartan group by a
Stickelberger ideal. This makes use of the characterization of units and
the quadratic relations, which determine which units have a given level
$N$. The connection with the geometry is made via the Fricke-Wohlfart
theorem. This allows Kubert and Lang to compute the order of the cuspidal
divisor class group in a way analogous to that of Iwasawa in the case
of cyclotomic fields. Whereas Iwasawa meets the first Bernoulli numbers
$B_{1, \chi}$, Kubert and Lang encounter the second Bernoulli numbers
$B_{2, \chi}$. Moreover, Kubert and Lang analyzed the eigen-space
decomposition at level $p$ on $X(p)$. Again they dealt with a purely
algebraic question in the group ring modulo the Stickelberger ideal,
and the eigen-spaces are generated by a product of the integralizing ideal
of quadratic relations times Bernoulli numbers, on the Cartan group.
Therefore, the cuspidal divisor class groups which is the torsion subgroup
of rational points on the Jacobian varieties of modular curves provide
geometric counterparts for the ideal class groups of cyclotomic fields in
algebraic number theory. The interplay of the algebraic number theory and
the theory of divisor classes for certain special number fields and special
curves is the topic of lectures given by \cite{Lang}. It was Wiles (see
\cite{Wiles1980}) who first showed precisely how this connection would come
about following work of Ribet (see \cite{Ribet1976}) in geometric language.
In particular, Wiles introduced Mazur's Eisenstein ideal (see \cite{Mazur1977})
which is closely related to the Stickelberger ideal into the picture relating
ideal class groups in the cyclotomic case with cuspidal divisor class groups.
Now, we will give an anabelian counterpart of this connection between cuspidal
divisor class groups in the modular Jacobians (motives) and cyclotomic fields
(abelian number fields), namely, there is a connection between defining ideals
of modular curves (anabelian algebraic geometry) and division fields of elliptic
curves (non-abelian number fields). In particular, the ideal of quadratic
relations is replaced by the ideal of quartic relations.

\begin{center}
{\large\bf 11.1. Cyclotomic fields and Stickelberger ideals: the Herbrand-Ribet
                theorem}
\end{center}

  An odd prime $p$ is called irregular if the class number of the field
$\mathbb{Q}(\mu_p)$ is divisible by $p$ ($\mu_p$ being, as usual, the
group of $p$-th roots of unity). According to Kummer's criterion, $p$ is
irregular if and only if there exists an even integer $k$ with $2 \leq k
\leq p-3$ such that $p$ divides (the numerator of) the $k$-th Bernoulli
number $B_k$, given by the expansion
$$\frac{t}{e^t-1}+\frac{t}{2}-1=\sum_{n \geq 2} \frac{B_n}{n!} t^n.$$
The purpose of \cite{Ribet1976} is to strengthen Kummer's criterion:
let $A$ be the ideal class group of $\mathbb{Q}(\mu_p)$, and let $C$
be the $\mathbb{F}_p$-vector space $A/A^p$. The Galois group
$\text{Gal}(\overline{\mathbb{Q}}/\mathbb{Q})$ acts on $C$ through its
quotient $\Delta=\text{Gal}(\mathbb{Q}(\mu_p)/\mathbb{Q})$. Since all
characters of $\Delta$ with values in $\overline{\mathbb{F}}_p^{*}$ are
powers of the standard character
$$\chi: \text{Gal}(\overline{\mathbb{Q}}/\mathbb{Q}) \rightarrow \Delta
  \stackrel{\thicksim}{\longrightarrow} \mathbb{F}_p^{*}$$
giving the action of $\text{Gal}(\overline{\mathbb{Q}}/\mathbb{Q})$ on
$\mu_p$, the vector space $C$ has a canonical decomposition
$$C=\bigoplus_{\text{$i$ mod $(p-1)$}} C(\chi^i),\eqno{(11.1.1)}$$
where the $\Delta$-eigen-spaces of $C$ is given by
$$C(\chi^i)=\{ c \in C | \sigma c=\chi^i(\sigma) c \quad \text{for all
               $\sigma \in \Delta$}\}.\eqno{(11.1.2)}$$

\textbf{Theorem 11.1.1. (Herbrand-Ribet theorem)}. {\it Let $k$ be even,
$2 \leq k \leq p-3$. Then $p | B_k$ if and only if $C(\chi^{1-k}) \neq 0$.}

  Furthermore, it is natural to ask for the precise order of $C(\chi^i)$
and it follows easily from Stickelberger's theorem that as a group
$C(\chi^i)$ is annihilated by $p^m$ where $m$ is defined as follows.
There is a canonical character $\theta: \Delta \rightarrow \mathbb{Z}_p^{*}$
which satisfies $\mu_p^g=\mu_p^{\theta(g)}$ for $g \in \Delta$. There
is also a natural isomorphism $\Delta \stackrel{\thicksim}{\longrightarrow}
(\mathbb{Z}/p \mathbb{Z})^{*}$ given by the action of $\Delta$ on the $p$-th
roots of unity, so using it we may define the first generalized Bernoulli
number (see \cite{Wiles1980})
$$B_{1, \theta^{-i}}=p^{-1} \sum_{a=1}^{p-1} \theta^{-i}(a) a \in \mathbb{Z}_p.$$

\textbf{Theorem 11.1.2.} (see \cite{Wiles1980} or Theorem 2.7 in \cite{Lang}). {\it
Let $i$ be an odd integer with $2<i \leq p-3$. If the eigenspace $C(\chi^{i})$
is cyclic, then its order is precisely $p^m$ where $m=v_p(B_{1, \theta^{-i}})$.}

  On the hand, the work of Kubert-Lang (\cite{Kubert-Lang}) shows that
there is a certain cyclic subgroup $C$ of $J_1(p)(\overline{\mathbb{Q}})$,
contained in the group generated by the cusps, and of order $p^n$ where
$n=v_p(B_{2, \theta^{-i-1}})$. The Hecke operators give rise to a ring of
endomorphisms $\mathbf{T}$ of $J_1(p)$, and we let $I$ be the ideal which
annihilates $C$. Then the representation considered in \cite{Wiles1980} is
essentially that of $\text{Gal}(\overline{\mathbb{Q}}/\mathbb{Q})$ on the
kernel of $I$ in $J_1(p)(\mathbb{Q})$. This leads to the following:

\begin{center}
{\large\bf 11.2. Cuspidal divisor class groups and Eisenstein ideals of
                modular curves: the theory of Kubert-Lang and Mazur-Wiles}
\end{center}

  According to \cite{Kamienny}, Let $p$ be a prime $\geq 13$, and let
$\Delta=\{ \langle a \rangle: a \in (\mathbb{Z}/p \mathbb{Z})^{*}/(\pm 1) \}$
denote the Galois group of the cover $\pi: X_1(p) \rightarrow X_0(p)$.
Let $J$ be the quotient of $J_1(p)$ by the image of $J_0(p)$ under the
map $\pi_{*}$ induced on Jacobians by $\pi$. We view the standard Hecke
operators $T_{\ell}$ (for $\ell \neq p$), $U_p$, $U_p^{\prime}$; and
$\Delta$ as endomorphisms of $J$ induced by Albanese functoriality, and
we define the Hecke algebra $\mathbf{T}$ to be the ring of endomorphisms
of $J$ generated over $\mathbb{Z}$ by these endomorphisms. We write $C_p$
for the projection to $J$ of the $p$-part of the subgroup of $J_1(p)$
generated by divisor classes supported only at the rational cusps. The
group $\Delta$ acts on $C_p$, and so we may decompose $C_p$ into eigen-spaces
for the action of $\Delta$:
$$C_p=\bigoplus C_p(\varepsilon)\eqno{(11.2.1)}$$
where the sum is taken over the even $\mathbb{Z}_p^{*}$-valued characters
$\varepsilon$ of $(\mathbb{Z}/p \mathbb{Z})^{*}$, and
$$C_p(\varepsilon)=\{ t \in C_p: \langle a \rangle t=\varepsilon(a) t \quad
  \text{for all $\langle a \rangle \in \Delta$} \}.\eqno{(11.2.2)}$$
The structure of $C_p(\varepsilon)$ has been computed by Kubert and Lang.
They proved the following:

\textbf{Theorem 11.2.1.} (see \cite{Kamienny} or Theorem 3.2 in \cite{Lang})
{\it If $\varepsilon \neq \omega^{-2}$ then
$C_p(\varepsilon) \cong \mathbb{Z}_p/p^{n_{\varepsilon}} \mathbb{Z}_p$
where $n_{\varepsilon}=v_p(B_{2, \varepsilon})$}.({\it Here $\omega$
is the Teichm\"{u}ller character, and $B_{2, \varepsilon}$ is the
second generalized Bernoulli number.})

  From now on we fix an even character $\varepsilon \neq \omega^{-2}$ of
$(\mathbb{Z}/p \mathbb{Z})^{*}$ with $n_{\varepsilon} \neq 0$, and write
$C$ for $C_p(\varepsilon)$. We define the Eisenstein ideal $I$ to be the
ideal of $\mathbf{T} \otimes \mathbb{Z}_p$ generated by $T_{\ell}-(1+\ell
\varepsilon(\ell))$ (for $\ell \neq p$), $U_p-1$, and
$\langle a \rangle-\varepsilon(a)$ for all $\langle a \rangle \in \Delta$.
Note that $I$ annihilates $C$. Let us consider the kernel $J_1(p)[I]$ on
$J_1(p)$ itself. This kernel, viewed as a Galois module, is an extension
of a twisted $\mu$-type group by the group $C$. It is shown that (see
\cite{Kamienny}) this extension is non-split. This shows that unramified
extensions of the type constructed by Ribet are seen already in the Galois
extension of $\mathbb{Q}$ cut out by $J_1(p)[I]$. This re-proves Ribet's
converse to Herbrand's criterion as follows:

\textbf{Corollary 11.2.2.} (see \cite{Kamienny}) {\it Let $i$ be an odd integer,
$2 \leq i \leq p-3$, and let $C(\chi^{-i})$ denote the $\chi^{-i}$-eigen-space
of the $p$-class group of $\mathbb{Q}(\mu_p)$}({\it where we have written
$\chi$ for the cyclotomic character.}) {\it If $v_p(B_{1, \chi^i}) \neq 0$,
then $C(\chi^{-i}) \neq 0$.}

  Thus, both Herbrand-Ribet theorem and Kubert-Lang theory are unified in the
context of generalized Bernoulli numbers. In the case of Herbrand-Ribet theorem,
it involves the first generalized Bernoulli number, while in the case of
Kubert-Lang theory, it involves the second generalized Bernoulli number.
We will show that in the non-abelian case, these analytic objects (i.e.,
two generalized Bernoulli numbers) are replaced by division equations
and transformation equations of elliptic functions, respectively.

  Let $N$ be a square-free integer and let $J_0(N)$ be the Jacobian of
the modular curve $X_0(N)$. In the case where $N$ is prime, Ogg (see
\cite{Ogg1972}) computed the order of the cuspidal subgroup (i.e., the
subgroup of $J_0(N)$ generated by linear equivalence classes of differences
of cusps) and conjectured that the cuspidal subgroup is the whole rational
torsion subgroup (see \cite{Ogg1975}). Mazur's proof of Ogg's conjecture
was one of the main arithmetic results of \cite{Mazur1977}:

\textbf{Theorem 11.2.3.} {\it Let $N \geq 5$ be a prime number,}
$n=\text{numerator}\left(\frac{N-1}{12}\right)$. {\it The torsion subgroup
of the Mordell-Weil group of $J_0(N)$ is a cyclic group of order $n$,
generated by the linear equivalence class of the difference of the two
cusps $(0)-(\infty)$.}

  In other words, the above theorem shows that the cuspidal divisor class
group of $X_0(N)$, the subgroup of $J_0(N)(\mathbb{Q})_{\text{tors}}$
generated by the difference of two cusps of $X_0(N)$, is of order
$(N-1)/(N-1, 12)$, coincides with the full rational subgroup
$J_0(N)(\mathbb{Q})_{\text{tors}}$.

  Ogg's computation of the order of the cuspidal subgroup has been generalized
to other modular curves by Kubert and Lang (see \cite{Kubert-Lang}). A number
of authors have considered the generalization of Ogg's conjecture to $J_0(N)$,
where $N$ is a positive integer that is not necessarily prime. Since the
cuspidal subgroup of $J_0(N)$ may not consist entirely of rational points,
the generalization states that the rational torsion subgroup of $J_0(N)$ is
contained in the cuspidal subgroup of this Jacobian. Roughly speaking, let
$C$ denote the $p$-primary part of the cuspidal subgroup and $T$ denote the
$p$-primary part of the rational torsion subgroup
$J_0(N)(\mathbb{Q})_{\text{tor}}$.

\textbf{Theorem 11.2.4.} (see \cite{Ribet2022} or \cite{Ohta2013}, \cite{Ohta2014})
{\it The rational torsion subgroup $T$ coincides with its subgroup $C$.}

  In particular, let $\mathbf{T}$ be the Hecke algebra acting on $J_0(N)$ and
let $I \subset \mathbf{T}$ be the Eisenstein ideal. Then it is proved that (see
\cite{Ribet2022}) the annihilators of $T$ and $C$ as $\mathbf{T}$-modules satisfy
$$\text{Ann}_{\mathbf{T}}(C)=I=\text{Ann}_{\mathbf{T}}(T).$$
Thus, the significance of the Eisenstein ideal comes from the fact that it is
the annihilator, in the endomorphism ring, of the group generated by the
cuspidal divisors in the Jacobian of the curve. In particular, there are five
principal structures in \cite{MW}, among which the cuspidal subgroups and
the Eisenstein ideal play an important role:

(a) {\it Good quotient abelian varieties} $B_{n/\mathbb{Q}}$ of the Jacobian of $X_1(p^n)$.

(b) {\it The Hecke algebra} $\mathbf{T}^{(n)}$ in the endomorphism ring of $B_n$.

(c) {\it The cuspidal subgroups} $C_{\chi}^{(n)} \subset B_n(\overline{\mathbb{Q}})$.
These are finite abelian $p$-groups, stable under the action of $\mathbf{T}^{(n)}$
and of $\text{Gal}(\overline{\mathbb{Q}}/\mathbb{Q})$.

(d) {\it The Eisenstein ideal} $I_{\chi}^{(n)}$, defined to be the annihilator in
$\mathbf{T}^{(n)}$ of $C_{\chi}^{(n)}$.

(e) {\it The kernel of the Eisenstein ideal}
$E_{\chi}^{(n)}=B_n(\overline{\mathbb{Q}})[I_{\chi}^{(n)}]$ in $B_n(\overline{\mathbb{Q}})$.
By construction we have
$$0 \longrightarrow C_{\chi}^{(n)} \longrightarrow B_n(\overline{\mathbb{Q}})[I_{\chi}^{(n)}]
  \longrightarrow M_{\chi}^{(n)} \longrightarrow 0$$
where $M_{\chi}^{(n)}$ is defined to make the sequence exact.

\begin{center}
{\large\bf 11.3. Modular curves and cyclotomic fields: Sharifi's conjecture}
\end{center}

  Recently, there is the other connection between modular curves and
cyclotomic fields which is different from the Herbrand-Ribet theorem.
This leads to the second approach (2) in the correspondence (1.18).

  Recall that by Serre duality theorem, we have
$$H^1(X(p), \mathcal{O}_{X(p)}) \cong H^0(X(p), \Omega_{X(p)}^1),$$
where $H^0(X(p), \Omega_{X(p)}^1)$ denotes the complex vector space of
global sections of the sheaf of holomorphic differentials on the Riemann
surface $X(p)$. Let us recall some facts from the theory of compact Riemann
surfaces: The first homology group $H_1(X(p), \mathbb{R})$ with real
coefficients is dual to $\Omega_{X(p)}^1$, as vector space over $\mathbb{R}$,
and we can identify the modular symbol $\{z_1, z_2 \}=\{ z_1, z_2 \}_{\Gamma(p)}$
as an element of $H_1(X(p), \mathbb{R})$ (see \cite{Lang1995}). In fact, let
$x$, $y$ be cusps. The Manin-Drinfeld theorem (see \cite{Lang1995}) shows
that $\{ x, y \}_{\Gamma(p)} \in H_1(X(p), \mathbb{Q})$ is in the homology
with rational coefficients. Moreover, we have a canonical isomorphism of
$\mathbb{C}$-vector spaces:
$$H_1(X(p), \partial_p, \mathbb{Z}) \otimes \mathbb{C} \rightarrow
  \text{Hom}_{\mathbb{C}}(H^0(X(p), \Omega_{X(p)}^1), \mathbb{C})
  \simeq H_1(X(p), \mathbb{R})$$
given by $c \rightarrow (\omega \rightarrow \int_{c} \omega)$, where
$\partial_p=\Gamma(p) \backslash \mathbb{H}$ is the set of cusps. It
defines a homomorphism of groups:
$$R: H_1(X(p), \partial_p, \mathbb{Z}) \rightarrow H_1(X(p), \mathbb{R}).$$
The fact that the image of $H_1(X(p), \partial_p, \mathbb{Z})$ under the
map $R$ is contained in $H_1(X(p), \mathbb{Q})$ is a reformulation of the
Manin-Drinfeld theorem (see \cite{Dr1973}, \cite{Manin1972}), which asserts
that the class of a divisor of degree zero supported on $\partial_p$ is
torsion in the Jacobian $J(p)$ of the modular curve $X(p)$. We call
Eisenstein classes the elements
$e \in H_1(X(p), \partial_p, \mathbb{Z}) \otimes \mathbb{R}
 \simeq H_1(X(p), \partial_p, \mathbb{R})$ such that $\int_e \omega=0$
for all $\omega \in H^0(X(p), \Omega_{X(p)}^1)$. The Eisenstein cycles
are paths in the Eisenstein classes. The Eisenstein classes constitute
a real vector $E(p)$. In other words, $e \in E(p)$ if and only if $e$
is in the kernel of the homomorphism $R$ obtained by extension of
scalars to $\mathbb{R}$. The Manin-Drinfeld theorem is reformulated
by saying that $E(p)$ admits a $\mathbb{Q}$-rational structure in
$H_1(X(p), \partial_p, \mathbb{Q})$. In \cite{BM}, Banerjee and Merel
determined explicitly a $\overline{\mathbb{Q}}$-basis for $E(p)$ in
terms of Manin symbol (see \cite{Manin1972}). Note that
$\text{PSL}(2, \mathbb{F}_p)=\pm \Gamma(p) \backslash \text{SL}(2, \mathbb{Z})$
and let
$$\xi: \text{PSL}(2, \mathbb{F}_p) \rightarrow H_1(X(p), \partial_p, \mathbb{Z})$$
be the map that takes the class of a matrix $g \in \text{SL}(2, \mathbb{Z})$
to the class in $H_1(X(p), \partial_p, \mathbb{Z})$ of the image of $X(p)$
of the geodesic in $\mathbb{H} \cup \mathbb{P}^!(\mathbb{Q})$ joining
$g(0)$ and $g(i \infty)$. We call $\xi(g)$ a Manin symbol associated to
$g$. The Manin symbols satisfy the Manin relations for all $g \in
\Gamma(p) \backslash \text{SL}(2, \mathbb{Z})$:
$$\xi(g)+\xi(gS)=0, \quad \xi(g)+\xi(gU)+\xi(gU^2)=0,$$
where $T=\left(\begin{matrix} 1 & 1\\ 0 & 1 \end{matrix}\right)$,
$S=\left(\begin{matrix} 0 & -1\\ 1 & 0 \end{matrix}\right)$ and
$U=ST=\left(\begin{matrix} 0 & -1\\ 1 & 0 \end{matrix}\right)$.
They expressed the Eisenstein classes in terms of Manin symbols with
coefficients in $\mathbb{Q}$. For $P=\left(\begin{matrix} \overline{x}\\
\overline{y} \end{matrix}\right) \in (\mathbb{Z}/p \mathbb{Z})^2$,
choose the representatives $x, y \in \mathbb{Z}$ of $\overline{x}$
and $\overline{y}$ respectively with $x-y$ odd. Define
$$F(P)=-\frac{1}{4} \left[\frac{\cos(\frac{\pi x}{p})+\cos (\frac{\pi y}{p})}
       {\cos(\frac{\pi x}{p})-\cos(\frac{\pi y}{p})}\right] \in \mathbb{R}.$$
They proved that $F(P) \in \overline{\mathbb{Q}}$ for all $P \in (\mathbb{Z}/
p \mathbb{Z})^2$. Let $\overline{F}$ be the function on $(\mathbb{Z}/p \mathbb{Z})^2
/\{ \pm 1 \}$ obtained from $F$ by passing to quotients. Define
$$\mathcal{E}_P=\sum_{\gamma \in \text{PSL}(2, \mathbb{F}_p)}
                \overline{F}(\gamma^{-1} P) \xi(\gamma)=\mathcal{E}_{-P}.$$

\textbf{Theorem 11.3.1.} (see \cite{BM}). {\it For $P \in (\mathbb{Z}/p \mathbb{Z})^2$,
the modular symbols $\mathcal{E}_P$ satisfy the following properties:

(1) Suppose $\ell$ is an odd prime number and} $\ell \equiv 1$ (mod $p$).
{\it Let $T_{\ell}$ be the Hecke operator at the prime $\ell$. The classes
$\mathcal{E}_P$ satisfy the equality
$$T_{\ell}(\mathcal{E}_P)=(\ell+1) \mathcal{E}_P.$$

(2) The classes of $\mathcal{E}_P$ lie in the kernel of $R$ and hence they are
    Eisenstein classes.

(3) The kernel of $R$ is spanned by the classes of $\mathcal{E}_P$ for $P \in
(\mathbb{Z}/p \mathbb{Z})^2$ with $P$ of order $p$.}

  This implies that the modular symbols $\mathcal{E}_P$ are annihilated
by the Eisenstein ideal. Hence, by the relative homology of modular curves,
Hecke's decomposition is related to the Eisenstein ideal. In general, the
role of modular symbols is determined by the fact that they constitute an
additional structure on the Betti realization of the cohomology of modular
varieties, i.e., they are modular motives. From the viewpoint of motives,
the modular symbols should be considered as a specific structure on the
$\mathbb{Q}$- or $\mathbb{Z}$-cohomology of the universal varieties (see
\cite{Manin1978}).

  There is a deeper relationship between the geometry of modular curves
and the arithmetic of cyclotomic fields in the form of Eisenstein maps
from the first homology of modular curves to the second $K$-groups of
cyclotomic integer rings, which are closely related to the class groups.
Sharifi viewed this as
$$\text{topology of $\text{GL}(2, \mathbb{Q})$ modulo Eisenstein}
  \Longleftrightarrow \text{arithmetic of $\text{GL}(1, \mathbb{Q})$}.$$
In his paper \cite{Sharifi2011} (see also \cite{Busuioc2008}), Sharifi
formulated a conjecture relating relative homology classes
$\{ \alpha \rightarrow \beta \}$ of paths between cusps on a modular
curve $X_1(p)$, and cup product $x \cup y$ of cyclotomic units in a
Galois cohomology group that agrees with $A^{-}$ modulo $p$, where
$A^{-}$ is the subgroup of $A$, the $p$-part of the ideal class group,
on which complex conjugation acts as $-1$. In fact, there is a simple
explicit map
$$\left\{ \frac{a}{c} \rightarrow \frac{b}{d} \right\} \mapsto
  (1-\zeta_p^c) \cup (1-\zeta_p^d)$$
taking one set of elements to the other (for $ad-bc=1$ and $p \nmid cd$).
On a certain Eisenstein quotients of the plus part of homology, Sharifi
conjectured this to provide an inverse to a canonical version of the
map that appeared in the proof of Ribet's theorem \cite{Ribet1976}.
More precisely, let $H_1(X_1(p), \text{cusps}, \mathbb{Z})$ be the singular
homology relative to the cusps of the modular curve $X_1(p)$. If $\alpha$
amd $\beta$ are in $\mathbb{P}^1(\mathbb{Q})$, let $\{ \alpha, \beta \}
\in H_1(X_1(p), \text{cusps}, \mathbb{Z})$ be the class of the closed
geodesics in the upper-half plane $\mathbb{H}$ between $\alpha$ and $\beta$,
which is called a modular symbol. If $(u, v) \in (\mathbb{Z}/p \mathbb{Z})^2$
is such that $\text{gcd}(u, v, p)=1$, let $\xi(u, v) \in H_1(X_1(p),
\text{cusps}, \mathbb{Z})$ be the modular symbol $\{ \frac{b}{d}, \frac{a}{c} \}$
where $\left(\begin{matrix} a & b\\ c & d \end{matrix}\right) \in
\text{SL}(2, \mathbb{Z})$ is such that $(c, d) \equiv (u, v)$ (mod $p$).
In fact, it can be proved that this does not depend on the choice of
$\left(\begin{matrix} a & b\\ c & d \end{matrix}\right)$ and the element
$\xi(u, v)$ (i.e. Manin symbols) generate $H_1(X_1(p), \text{cusps}, \mathbb{Z})$
(see \cite{Manin1972}). Using Manin's presentation of $H_1(X_1(p), \text{cusps},
\mathbb{Z})$, Sharifi defined a map
$$\varpi: H_1(X_1(p), C^{o}, \mathbb{Z}) \longrightarrow
  K_2\left(\mathbb{Z}\left[\zeta_p, \frac{1}{p}\right]\right)
  \otimes \mathbb{Z}\left[\frac{1}{2}\right]$$
where $C^{o}$ is the subset of cusps in $X_1(p)$ not lying over the cusp
$\infty$ of $X_0(p)$. Sharifi's conjecture gives an explicit map from
the homology of a modular curve to the second $K$-group of a cyclotomic
integer ring, sending Manin symbols to Steinberg symbols of cyclotomic
units (see \cite{FKS}):
$$\varpi(\xi(u, v))=\{1-\zeta_p^u, 1-\zeta_p^v\}$$
that relates the worlds of geometry/topology and arithmetic. Here,
$\xi(u, v)$ is a Manin symbol in the relative homology group $H_1(X_1(p),
\{ \text{cusps} \}, \mathbb{Z})$, and $\{ 1-\zeta_p^u, 1-\zeta_p^v \}$
is a Steinberg symbol in the algebraic $K$-group $K_2(\mathbb{Z}[\zeta_p,
\frac{1}{p}])$, where $u, v \in \mathbb{Z}/p \mathbb{Z}$ are nonzero
numbers with $(u, v)=1$. The above map connects two different worlds
in the following manner:
$$\text{geometric theory of $\text{GL}(2)$} \Rightarrow
  \text{arithmetic theory of $\text{GL}(1)$}$$
over the field $\mathbb{Q}$. Here, if we consider the geometry of the
modular curve $X_1(p)$ (i.e., relative homology) on the left, then we
consider the arithmetic of the cyclotomic field $\mathbb{Q}(\zeta_p)$
on the right. This map was conjectured to be, and has been largely proven
to be Eisenstein (see \cite{FK2024} and \cite{SV2024}), in the sense that
it factors through the quotient of homology by an ideal in the Hecke algebra
arising from coefficients of Eisenstein series:

\textbf{Conjecture 11.3.2. (Sharifi's Conjecture).} {\it The map $\varpi$ is
annihilated by the Eisenstein ideal $I$ generated by the Hecke operators
$T_{\ell}-\ell-\langle \ell \rangle$ for primes $\ell \neq p$, where
$\langle \ell \rangle$ is the $\ell$-th diamond operator.}

$$\begin{array}{ccc}
  \text{geometric theory of $\text{GL}(2)$} & \longleftrightarrow & \text{arithmetic theory}\\
  \text{(i.e. relative homology)} &  & \text{of $\text{GL}(1)$}\\
  \text{modulo the Eisenstein ideal} &  & \text{(i.e. cyclotomic fields)}\\
  \text{(Fukaya, Kato and Sharifi)} &  &
\end{array}\eqno{(11.3.1)}$$

  Let us recall some basic facts about the Beilinson elements in $K_2$ of
modular curves (see \cite{Kato2004}, \cite{FK2024} and \cite{Brunault}).
Let $p \geq 3$ be a prime and $Y(p)$ be the modular curve classifying
elliptic curves $E$ with a level $p$, that is, a basis $(e_1, e_2)$ of
$E[p]$ over $\mathbb{Z}/p \mathbb{Z}$. The curve $Y(p)$ is a smooth
projective curve defined over $\mathbb{Q}$, whose affine ring
$\mathcal{O}(Y(p))$ contains the cyclotomic field $\mathbb{Q}(\zeta_p)$.
The curve $Y(p)$ is not geometrically connected. Indeed, there is an isomorphism
$Y(p)(\mathbb{C}) \cong (\mathbb{Z}/p \mathbb{Z})^{*} \times (\Gamma(p) \backslash \mathbb{H})$.
The curve $Y(p)$ has a smooth compactification $X(p)$ over $\mathbb{Q}$
which is obtained by adding on the cusps. The function field of $X(p)$ will
be referred to by $\mathbb{Q}(X(p))$. Let us give the definition of Siegel
units (see \cite{Kubert-Lang}). The group of modular units of $X(p)$ will be
denoted by $\mathcal{O}^{*}(Y(p))$. The significance of modular units comes
from the fact that the trivial elements of the cuspidal divisor class group
are represented precisely by the divisors of modular units (see \cite{Kubert-Lang}).
In order to avoid torsion problems, Siegel units will always be considered
in the $\mathbb{Q}$-vector space
$\mathcal{O}^{*}(Y(p)) \otimes_{\mathbb{Z}} \mathbb{Q}$.
Let $B_2(X)=X^2-X+\frac{1}{6}$ be the second Bernoulli polynomial.

\textbf{Definition 11.3.3.} For any $(\alpha, \beta) \in (\mathbb{Z}/p \mathbb{Z})^2
-\{ (0, 0) \}$ let us define
$$g_{\alpha, \beta}(z)=q^{\frac{1}{2} B_2(\widetilde{\alpha}/p)} \prod_{n \geq 0}
  (1-q^n q^{\widetilde{\alpha}/p} \zeta_p^{\beta}) \prod_{n \geq 1}
  (1-q^n q^{-\widetilde{\alpha}/p} \zeta_p^{-\beta}),$$
where $\widetilde{\alpha} \in \mathbb{Z}$ is the unique representative of
$\alpha$ satisfying $0 \leq \widetilde{\alpha} <p$.

  It can be proved that $g_{\alpha, \beta}$ is well-defined as an element
of $\mathcal{O}^{*}(Y(p)) \otimes \mathbb{Q}$. The group $\text{GL}(2, \mathbb{F}_p)$
acts from the left on $Y(p)$ by the rule
$$\left(\begin{matrix}
  a & b\\ c & d \end{matrix}\right) \cdot (E, e_1, e_2)
 =(E, a e_1+b e_2, c e_1+d e_2) \quad \text{for} \quad
 \left(\begin{matrix}
  a & b\\ c & d \end{matrix}\right) \in \text{GL}(2, \mathbb{F}_p).$$
This induces on $\mathcal{O}^{*}(Y(p)) \otimes \mathbb{Q}$ a right action
of $\text{GL}(2, \mathbb{F}_p)$.  It turns out that $\text{GL}(2, \mathbb{F}_p)$
acts on the set of Siegel units, and the Siegel units of level $p$ generate
$\mathcal{O}^{*}(Y(p)) \otimes \mathbb{Q}$ (see \cite{Kubert-Lang}). Let us
consider the Quillen $K$-group $K_2(Y(p))$ which enjoys a right action of
$\text{GL}(2, \mathbb{F}_p)$ by functoriality. Beilinson constructed special
elements in it using cup-products of modular units. This leads to the map
$$\varrho: M_2(\mathbb{F}_p) \rightarrow K_2(Y(p)) \otimes_{\mathbb{Z}} \mathbb{Q},
  \quad \left(\begin{matrix} s & t\\ u & v \end{matrix}\right) \mapsto
  \{ g_{s, t}, g_{u, v} \}.$$
In \cite{Sharifi2011}, Sharifi expected that Beilinson elements would be
useful for the study of his conjectures. The most important of many keys
to the proofs of the results of Fukaya-Kato \cite{FK2024} is the fact that
Beilinson elements specialize to Steinberg symbols of cyclotomic units at
the cusps. Fukaya, Kato, and Sharifi (see \cite{FKS}) expected this to be
a special case of a general phenomenon of the geometry and topology of
locally symmetric spaces of higher dimension informing the arithmetic of
Galois representation attached to lower-dimensional automorphic forms. This
leads Sharifi, in the end of his survey paper on the Iwasawa theory
\cite{Sharifi2019}, to ask the following question: ``This begs the question
of its elliptic curve analogue, which is but one of a wealth of intriguing
possibilities for future directions in Iwasawa theory.''

  Now, our theory (Theorem 1.1, Corollary 1.2 and Theorem 1.3) gives such
an answer in the following:
$$\begin{array}{ccc}
  \text{anabelian geometric theory} & \longleftrightarrow & \text{arithmetic theory}\\
  \text{of $\text{GL}(2)$ (i.e. $\pi_1$) associated} &  & \text{of $\text{GL}(2)$}\\
  \text{with the defining ideal} &  & \text{(i.e. division fields of}\\
  \text{(ours)} &  & \text{elliptic curves over $\mathbb{Q}$)}
\end{array}\eqno{(11.3.2)}$$
It can be regarded as an anabelianization of the correspondence (11.3.1)
given by Fukaya, Kato and Sharifi. The right-hand side of (11.3.2) will
be given in the next section 11.4.

\begin{center}
{\large\bf 11.4. Division fields of elliptic curves: a non-abelian
                counterpart of the Herbrand-Ribet theorem}
\end{center}

  In \cite{Prasad} and \cite{PrasadS}, Prasad gave a non-abelian counterpart
of the Herbrand-Ribet theorem (see \cite{Ribet1976}) for division fields of
elliptic curves. More precisely, let $E$ be an elliptic curve over $\mathbb{Q}$,
$p$ an odd prime number, and $\rho: \text{Gal}(\overline{\mathbb{Q}}/\mathbb{Q})
\rightarrow \text{GL}(2, \mathbb{F}_p)$ the associated Galois representation on
elements of order $p$ on $E$. Assume that the image of the Galois representation
is all of $\text{GL}(2, \mathbb{F}_p)$. The representation $\rho$ then gives rise
to an extension $K$ of $\mathbb{Q}$ with Galois group $\text{GL}(2, \mathbb{F}_p)$.
Let $\text{Cl}_K$ denote the class group of $K$. The group $\text{GL}(2, \mathbb{F}_p)$
being the Galois group of $K$ over $\mathbb{Q}$, operates on $\text{Cl}_K$,
hence on the $\mathbb{F}_p$-vector space $\text{Cl}_K/p \text{Cl}_K$. Write the
semi-simplification of the representation of $\text{GL}(2, \mathbb{F}_p)$ on
$\text{Cl}_K/p \text{Cl}_K$ as $\sum V_{\alpha}$, where $V_{\alpha}$'s are the
various irreducible representations of $\text{GL}(2, \mathbb{F}_p)$ in
characteristic $p$. It is well-known that any irreducible representation
of $\text{GL}(2, \mathbb{F}_p)$ in characteristic $p$ is of the form
$V_{i, j}=\text{Sym}^i \otimes \det^j$, $0 \leq i \leq p-1$, $0 \leq j \leq p-2$
where $\text{Sym}^i$ refers to the $i$-th symmetric power of the standard
$2$-dimensional representation of $\text{GL}(2, \mathbb{F}_p)$ on
$\mathbb{F}_p+\mathbb{F}_p$, and $\det$ denotes the determinant character
of $\text{GL}(2, \mathbb{F}_p)$.

  It is a natural question to understand which $V_{i, j}$'s appear in
$\text{Cl}_K/p \text{Cl}_K$. In \cite{Prasad} and \cite{PrasadS}, Prasad
formulated some questions in this direction which can be viewed as a
$\text{GL}(2)$ analogue of the Herbrand-Ribet theorem (see \cite{Ribet1976}),
which can be considered to be a theorem for $\text{GL}(1)=\mathbb{G}_m$.
He hoped for a conjectural answer along the lines of Herbrand-Ribet to say
that representation $\text{Sym}^i(\mathbb{F}_p+\mathbb{F}_p) \otimes \det^j$
of $\text{GL}(2, \mathbb{F}_p)$ appears in $\text{Cl}_K/p \text{Cl}_K$
if and only if the algebraic part of the first nonzero derivative of
$L(s, \text{Sym}^i(E) \otimes \det^j)$ at $s=0$ is divisible by $p$.

  Prasad related \begin{otherlanguage}{russian} Sh \end{otherlanguage} of
an elliptic curve with the class group of $\mathbb{Q}(E[p])$ as follows:

{\bf Question 11.4.1.} Let $E$ be an elliptic curve over $\mathbb{Q}$ such
that $E(\mathbb{Q})=0$. Let $K=\mathbb{Q}(E[p])$ be the Galois extension of
$\mathbb{Q}$ obtained by attaching elements of order $p$ on $E$ where $p$
is an odd prime. We assume that $\text{Gal}(K/\mathbb{Q})=\text{GL}(2, \mathbb{F}_p)$,
and also that $p$ is coprime to $c_{\ell}=[E(\mathbb{Q}_{\ell}): E(\mathbb{Q}_{\ell})^0]$,
the Tamagawa numbers, for all finite primes $\ell$. Let $\text{Cl}_K$ denote
the class group of $K$ which comes equipped with a natural action of
$\text{Gal}(K/\mathbb{Q})=\text{GL}(2, \mathbb{F}_p)$. Then if $p$ divides
$|$\begin{otherlanguage}{russian} Sh \end{otherlanguage}$\!\!(E/\mathbb{Q})|$,
then is it true that the $\text{GL}(2, \mathbb{F}_p)$ representation
$\text{Cl}_K/p \text{Cl}_K$ contains the standard $2$-dimensional representation
of $\text{GL}(2, \mathbb{F}_p)$? What about the converse?

  Here, for the ideal class groups of division fields of elliptic curves,
the Tate-Shafarevich groups and the Selmer groups of elliptic curves are
central objects to study. Prasad and Shekhar (see \cite{PrasadS}) gave a
sufficient condition which implies that the semi-simplification of the
$p$-part of $\text{Cl}_K$ with $K=\mathbb{Q}(E[p])$ contains $E[p]$, in
terms of the Tate-Shafarevich group \begin{otherlanguage}{russian} Sh
\end{otherlanguage}$\!\!(E/\mathbb{Q})$ of $E$ over $\mathbb{Q}$ or the
Selmer group of $E$. Note that the semi-simplification is a kind of the
irreducible decomposition as a Galois module.

\textbf{Theorem 11.4.2.} {\it Let} $\rho_{E, p}: \text{Gal}(\overline{\mathbb{Q}}/\mathbb{Q})
\rightarrow \text{Aut}(E[p]) \cong \text{GL}(2, \mathbb{F}_p)$ {\it be the
$\mathbb{F}_p$-valued Galois representation associated to $E$. Suppose that
the following conditions on $E$ hold:}

(a) {\it $E$ has good reduction at $p$.}

(b) {\it In the case that $E$ has good ordinary reduction at $p$},
    $a_p(E) \equiv 1$ (mod $p$), {\it and $E$ has no complex multiplication
    over an extension of $\mathbb{Q}_p$, then $\rho_{E, p}$ is wildly
    ramified at $p$.}

(c) {\it For every prime number $\ell \neq p$, the Tamagawa number
    $c_{\ell}(E/\mathbb{Q}_{\ell})$ of $E/\mathbb{Q}_{\ell}$ is prime to $p$.}

(d) {\it $E[p]$ is an irreducible}
    $\text{Gal}(\overline{\mathbb{Q}}/\mathbb{Q})$-{\it module}.

\noindent {\it Then the condition on the Tate-Shafarevich group}
$\dim_{\mathbb{F}_p}($\begin{otherlanguage}{russian} Sh
\end{otherlanguage}$\!\!(E/\mathbb{Q})[p]) \geq 2$ {\it or the condition on
the Selmer group} $\dim_{\mathbb{F}_p}(\text{Sel}(\mathbb{Q}, E[p])) \geq 2$
{\it implies that the $\mathbb{F}_p$-representation} $\text{Cl}_K \otimes
\mathbb{F}_p$ {\it of} $\text{Gal}(\mathbb{Q}(E[p])/\mathbb{Q})$ {\it has
$E[p]$ as its quotient representation.}

  In his thesis (see \cite{Dainobu2022} and \cite{Dainobu2024}), Dainobu
gave a sufficient condition which implies that the semi-simplification of
the $p$-part of $\text{Cl}(\mathbb{Q}(E[p]))$ contains the $j$-th symmetric
product $\text{Sym}^j E[p]$ of $E[p]$, in terms of Bloch-Kato's Tate-Shafarevich
groups or Bloch-Kato Selmer group of $\text{sym}^j E[p]$. This result is
a partial generalization of the result of Prasad and Shekhar, which
considered the case $j=1$.

  Despite of this, Prasad pointed out the following (see \cite{Prasad}
and \cite{PrasadS}):

\textbf{Problem 11.4.3.} ``The authors have not seen any computation identifying
the $\text{GL}(2, \mathbb{F}_p)$ representation on $\text{Cl}_K/p \text{Cl}_K$
as a sum of irreducible pieces (after semi-simplification) where
$K=\mathbb{Q}(E[p])$ is a Galois extension of $\mathbb{Q}$ obtained by
attaching elements of order $p$ on an elliptic curve $E$ over $\mathbb{Q}$
with $\text{Gal}(K/\mathbb{Q})=\text{GL}(2, \mathbb{F}_p)$. Presumably
it is not beyond present computational powers to do such a computation
say for $p=5$; this would be very useful data to have for the problems
discussed in this paper. There are examples known due to K. Rubin and
A. Silverberg (see \cite{RubinS}) of families of elliptic curves with
the same field $\mathbb{Q}(E[5])$ with Galois group $\text{GL}(2, \mathbb{F}_5)$.''

\begin{center}
{\large\bf 11.5. Defining ideals of modular curves: an anabelian counterpart of the
                 eigenspace decomposition and the cuspidal divisor class formula}
\end{center}

  In order to go beyond the abelian setting (both abelian varieties in
the case of division fields of elliptic curves and abelian extensions
in the cases of cyclotomic fields as well as cuspidal divisor class
groups), we will give a counterpart of the above theories in the
anabelian setting. In fact, our decomposition theorem (Theorem 1.1 and
Corollary 1.2), which can be regarded as an anabelian counterpart of the
eigenspace decomposition, satisfies the semi-simplicity as well as the
multiplicity one property (at least for $p=7$, $11$ and $13$ comparing
with \cite{Ochiai1} and \cite{Ochiai2} for Iwasawa theory). In particular,
although the problem 8.4.3 remains open even for the simplest case $p=5$.
In the context of modular curves, Corollary 1.2 provides such an explicit
construction for the three cases $p=7$, $11$ and $13$. Problem 8.4.3 is
in the context of algebraic number theory, while Corollary 1.2 is in the
context of arithmetical algebraic geometry. Hence, Corollary 1.2, i.e.,
the dimension formula for defining ideals of modular curve can be regarded
as an anabelian counterpart of the cuspidal divisor class number formula in
(1.23). At the same time, Corollary 1.2 can also be considered as a modular
curve counterpart of the non-abelian Herbrand-Ribet theorem.

\begin{center}
{\large\bf 12. Other consequences and relations with various topics}
\end{center}

\begin{center}
{\large\bf 12.1. Relation with Langlands' work}
\end{center}

  Langlands' interest in Shimura varieties was in their relations with
automorphic representations. In particular, his first paper on Shimura
varieties is the fundamental article \cite{Langlands73} on the basic
one-dimensional Shimura varieties $S_K$ attached to the group $\text{GL}(2)$
(see \cite{Arthur2024}). In this case, the associated automorphic representations
are expected to be more manageable, in the sense that their archimedean
components $\pi_{\infty}$ should typically be square integrable (as in a
holomorphic modular form of weight $k \geq 2$), rather than an induced
representation (as in a Maass form). It is for this reason that the
reciprocity laws between $\ell$-adic representations and automorphic
forms tend to be more concrete. Langlands' goal was to establish
reciprocity laws between the two-dimensional $\lambda$-adic representations
that arise from the \'{e}tale cohomology of $S_K$ and the automorphic
representations of $\text{GL}(2)$. Before that, reciprocity laws between\
Frobenius classes and Hecke operators on modular curves were established
by Eichler and Shimura. It was an extension to higher weight of these
results that Deligne used to reduce the Ramanujan conjecture for
holomorphic modular forms to the last of the Weil conjectures, which he
later established in 1974. The congruence relation used by Eichler and
Shimura do not generalize easily beyond modular curves, whereas the trace
formulas extend in principle to arbitrary Shimura varieties. The problem
for Langlands was to extend the Lefschetz formula to the open Shimura
varieties $S_K$ for $\text{GL}(2)$, and to compare the results with
the Selberg trace formula. The $\lambda$-adic representations of
$\text{Gal}(\overline{\mathbb{Q}}/\mathbb{Q})$ are on the \'{e}tale
cohomology groups $H_{\text{\'{e}t}}^1(S_K)$ of $S_K$. What does this
have to do with automorphic representations? Langlands gave the precise
relationship between automorphic representations and automorphic forms.
The latter objects are closely related to differential forms on
$S_K(\mathbb{C})$, with values in a locally constant sheaf $\mathcal{F}$.
They can therefore be used to construct de Rham cohomology groups
$H_c^1(S_K(\mathbb{C}), \mathcal{F}(\mathbb{C}))$, which in turn
lead to an interpretation of these groups in terms of automorphic
representations. In fact, Hecke operators act on the analytic
cohomology $H_c^1(S_K(\mathbb{C}), \mathcal{F}(\mathbb{C}))$ because
they act on $S_K(\mathbb{C})$ as an analytic manifold. They act on the
arithmetic cohomology $H_{\text{\'{e}t}}^1(S_K, \mathcal{F})$ because
they act on $S_K$ as an algebraic variety over the reflex field
$E(G, h)=\mathbb{Q}$. The intertwining operators between these two
actions, combined with the theorem of strong multiplicity one for
$\text{GL}(2)$, led to a general bijective correspondence
$$\pi \rightarrow \sigma(\pi),\eqno{(12.1)}$$
between automorphic representations $\pi$ of $G$ and two-dimensional
$\lambda$-adic representations $\sigma$ of
$\text{Gal}(\overline{\mathbb{Q}}/\mathbb{Q})$. The reciprocity problem
was then to describe the correspondence explicitly. Note that the space
$S_K(\mathbb{C})$ is noncompact, so one account for the behavior of
differential forms at infinity. Langlands takes the image of the
cohomology of compact support $H_c^{*}( \cdot )$ in the full de Rham
cohomology. In general, it is better to work with $L^2$-cohomology
$H_{(2)}^{*}( \cdot )$, as has been the custom since the introduction
of intersection cohomology, with its role in Zucker's conjecture. One
is interested in the $L^2$-cohomology
$$H_{(2)}^{*}(S_K(\mathbb{C}), \mathcal{F})=\bigoplus_{d=0}^{2n}
  H_{(2)}^{d}(S_K(\mathbb{C}), \mathcal{F}), \quad n=\dim S_K(\mathbb{C}),$$
of $S_K(\mathbb{C})$ with coefficients in $\mathcal{F}$. We recall that
this is the cohomology of the complex of $\mathcal{F}$-valued, smooth
differentiable forms $\omega$ on $S_K(\mathbb{C})$ such that both $\omega$
and $d \omega$ are square integrable. This graded, complex vector space
has a spectral decomposition (see \cite{Arthur2024}, \cite{BW} and \cite{BC})
$$\bigoplus_{\pi} \left(m_2(\pi) \cdot H^{*}(\mathfrak{g}_{\mathbb{R}},
  K_{\mathbb{R}}; \pi_{\mathbb{R}} \otimes \xi) \otimes \pi_f^K\right),$$
where $\pi=\pi_{\mathbb{R}} \otimes \pi_f$ ranges over automorphic
representations of $G$ with archimedean and non-archimedean components
$\pi_{\mathbb{R}}=\pi_{\infty}$ and $\pi_f=\pi^{\infty}$, while $m_2(\pi)$
is the multiplicity with which $\pi$ occurs in the $L^2$-discrete spectrum
with appropriate central character determined by $\xi$, a fixed irreducible,
finite-dimensional representation of $G$, and $\pi_f^K$ is the finite
dimensional space of $K$-invariant vectors for $\pi_f$. For the special
case of the Shimura variety attached to the group $G=\text{GL}(2)$ above,
Langlands takes the representation $\xi$ of $G$ as $\xi=\xi_k \otimes (\det)^m$,
where $\xi_k=\text{Sym}^{k-1}(\text{St})$ is the $(k-1)$-symmetric power
of the standard representation of $G$ of dimension $k$, and $\det$ is the
one-dimensional determinant representation, for integers $k \in \mathbb{N}$
and $m$, with $0 \leq m <k$. It turns out that there is thus a decomposition
$$H_{(2)}^1(S_K(\mathbb{C}), \mathcal{F})=\bigoplus_{\{ \pi:
  \pi_{\mathbb{R}}=\pi_{\mathbb{R}}(\xi) \} }
  \left(m_2(\pi) (H^1(\pi_{\mathbb{R}}, \xi) \otimes \pi_f^K)\right)
  \eqno{(12.2)}$$
of the first cohomology group (see \cite{Arthur2024} for more details).
Langlands gave an interpretation of the multiplicities $m_2(\pi)$ by
automorphic representations. The existence of the above Langlands
correspondence then follows from the global properties of $\lambda$-adic
\'{e}tale cohomology. Here, the domain consists of the automorphic
representations $\pi$ that give nonzero summands in the above decomposition
formula, while the image is the corresponding set of two-dimensional
$\lambda$-adic representations $\sigma$. In particular, Langlands' work
\cite{Langlands73} is essentially complete and it has served as a model
for subsequent work on the cohomology of Shimura varieties, which
continues to this day.

  In general, our ultimate goal is founded on the two interpretations
of cohomology and the fundamental data they support. On the one hand,
we have the $L^2$-de Rham cohomology $H_{(2)}^{*}(S_K(\mathbb{C}), \mathcal{F})$,
and its decomposition in terms of automorphic representations. On the other
hand, we have the intersection cohomology
$IH^{*}(\overline{S}_K(\mathbb{C}), \mathcal{F})$ of the Baily-Borel-Satake
compactification $\overline{S}_K(\mathbb{C})$, equipped with the
correspondences defined by Hecke operators. In fact, these two complex
graded vector spaces are isomorphic. The ultimate goal of the theory
would be to deduce reciprocity laws for $S_K$ of the kind established
by Langlands for $G=\text{GL}(2)$ in \cite{Langlands73}.

  For $\text{GL}(2, \mathbb{Q})$, Langlands' interest in modular
curves was in their relations with automorphic representations for
$\text{GL}(2, \mathbb{A}_{\mathbb{Q}})$, i.e., the adelic version of
$\text{GL}(2, \mathbb{Q})$, by the various cohomology $H^1$ of modular
curves. On the other hand, our interest in modular curves is in their
relations with $\mathbb{Q}(\zeta_p)$-rational representations for
$\text{PSL}(2, \mathbb{F}_p)$, i.e., the mod $p$ version of
$\text{GL}(2, \mathbb{Q})$, by the defining ideals of modular curves.
In particular, for $\text{GL}(2, \mathbb{Q})$, (1.2) and (1.3) can be
regarded as an anabelian counterpart of (12.2), (1.8) and (1.9) can be
regarded as an anabelian counterpart of (12.1). Here, the (1.2), (1.3),
(1.8) and (1.9) are in the context of anabelian algebraic geometry, while
(12.2) and (12.1) are in the context of motives.

\begin{center}
{\large\bf 12.2. Relation with motives and automorphic representations}
\end{center}

  According to \cite{BCGP}, for algebraic curves defined over $\mathbb{Q}$
with genus $g$, there is a tetrachotomy corresponding to the cases $g=0$,
$g=1$, $g=2$ and $g>2$.

  The $g=0$ case goes back to the work of Riemann.

  The key point in the $g=1$ case (where the relevant objects are modular
forms of weight two) is that there are two very natural ways to study these
objects. The first (and more classical) way to think about a modular form
is as a holomorphic function on the upper half plane which satisfies specific
transformation properties under the action of a finite index subgroup of
$\text{SL}(2, \mathbb{Z})$. This gives a direct relationship between modular
forms and the coherent cohomology of modular curves; namely, cuspidal modular
forms of weight two and level $\Gamma_0(N)$ are exactly holomorphic differentials
on the modular curve $X_0(N)$. On the other hand, there is a second interpretation
of modular forms of weight two in terms of the Betti (or \'{e}tale or de Rham)
cohomology of the modular curve. A direct way to see this is that holomorphic
differentials can be thought of as smooth differentials, and these satisfy a
duality with the homology group $H_1(X_0(N), \mathbb{R})$ by integrating a
differential along a loop. And it is the second description (in terms of
\'{e}tale cohomology) which is vital for studying the arithmetic of modular
forms.

  When $g=2$, there is still a description of the relevant forms in terms
of coherent cohomology of Shimura varieties (now Siegel three-folds), but
there is no longer any direct link between these coherent cohomology groups
and \'{e}tale cohomology. As pointed out in \cite{BCGP}, the reason one can
not prove an analogue of Wiles' strategy (which proves that all ordinary
semistable elliptic curves over $\mathbb{Q}$ with $\overline{\rho}_{E, 3}$
absolutely irreducible are modular) is that there is no argument to reduce
the residual modularity of a surjective mod-$3$ representation $\overline{\rho}_3:
\text{Gal}(\overline{F}/F) \rightarrow \text{GSp}(4, \mathbb{F}_3)$ to special
cases of the Artin conjecture (proved by Langlands-Tunnell), where $F$ is a
totally real field. The difficulty is not simply that $\text{PSp}(4, \mathbb{F}_3)$
is not solvable, but also that Artin representations do not contribute to the
coherent cohomology of Shimura varieties in any setting other than holomorphic
Hilbert modular forms of weight one. On the automorphic side, weight two
Hilbert-Siegel modular cuspforms only occur in the coherent cohomology of the
Hilbert-Siegel modular variety. Since the Hecke eigenvalues associated to such
modular forms are not realized in the \'{e}tale cohomology of a Shimura variety,
we do not know how to associate a motive to a weight two Hilbert-Siegel modular
cuspidal eigenform, but only a compatible system of Galois representations which
should correspond to the system of $\ell$-adic realizations of this motive. In
fact, both strategies of \cite{W} and \cite{BCGP} use the division points on
elliptic curves or Abelian surfaces, and thus mod $p$ Galois representations.

  Finally, when $g>2$, even the relationship with coherent cohomology disappears,
the relevant automorphic objects have some description in terms of differential
equations on locally symmetric spaces, but there is no longer any way to get
a handle on these spaces. For those that know about Maass forms, the situation
for $g>2$ is at least as hard (probably much harder) than the notorious open
problem of constructing Galois representations associated to Maass forms of
eigenvalue $\frac{1}{4}$. In other words, it is probably very hard!

  Correspondingly, motives can be divided according to a tetrachotomy,
respectively (see \cite{Calegari}). The first form (corresponding to $g=0$)
are the Tate (and potentially Tate) motives, whose automorphy was known to
Riemann and Hecke. The second form (corresponding to $g=1$) are the motives
(conjecturally) associated to automorphic representations which are discrete
series at infinity and thus amenable to the Taylor-Wiles method. The third
form (corresponding to $g=2$) are the motives (conjecturally) associated
to automorphic representations which are at least seen by some flavour of
cohomology, either by the Betti cohomology of locally symmetric spaces or
the coherent cohomology of Shimura varieties (possibly in degrees greater
than zero) which are amenable in principle to the modified Taylor-Wiles
method. The fourth form (corresponding to $g>2$) consist of the rest, which
(besides a few that can be accessed by cyclic base change) are a complete
mystery.

  Moreover, the corresponding automorphic representations can also be divided
according to a tetrachotomy, respectively. It is well-known that an elliptic
curve (corresponding to $g=1$) over $\mathbb{Q}$ should arise from an
automorphic representation of $\text{GL}(2)$ whose archimedean component is
discrete series. (We know this to be the case by Wiles \cite{W}, Taylor-Wiles
\cite{TW} and Breuil-Conrad-Diamond-Taylor \cite{BCDT}.) An abelian surface
over $\mathbb{Q}$ (corresponding to $g=2$) should arise from an automorphic
representation of $\text{GSp}(4)$ whose archimedean component is a holomorphic
limit of discrete series. (We know that abelian surfaces over totally real
fields are potentially modular by Boxer-Calegari-Gee-Pillon \cite{BCGP}.) In
contrast with the above two cases, an abelian three-fold over $\mathbb{Q}$
(corresponding to $g=3$) should arise from an automorphic representation
$\pi$ of $\text{GSpin}(2, 5)$ whose archimedean component is a degenerate
limit of discrete series, but not from a $\pi$ whose archimedean component
is a non-degenerate limit of discrete series. Moreover, an abelian variety
of dimension $n \geq 4$ over $\mathbb{Q}$ (corresponding to $g>3$) should
arise from an automorphic representation $\pi$ of $\text{GSpin}(n, n+1)$
whose archimedean component is a degenerate limit of discrete series, but
neither from a $\pi$ whose archimedean component is a non-degenerate limit
of discrete series, nor from a $\pi$ of a group of Hermitian type (see
\cite{Goldring}).

  In fact, all generalizations of the Taylor-Wiles method to this point require
that the automorphic representations in question are associated to the Betti
cohomology groups of locally symmetric spaces, or the coherent cohomology
groups of Shimura varieties, which have integral structures and hence allow
one to talk about congruences between automorphic forms. The automorphic
representations contributing to the Betti cohomology groups of locally
symmetric spaces have regular infinitesimal characters, so can only be
used for $g=1$. The automorphic representations contributing to the coherent
cohomology of orthogonal Shimura varieties are representations of the inner
form $\text{SO}(2g-1, 2)$ of $SO(2g+1)$ (which is non-split if $g>1$), whose
infinity components $\pi_{\infty}$ are furthermore either discrete series,
or non-degenerate limits of discrete series. If $g=1$, the representations
are discrete series, and if $g=2$, they are non-degenerate limits of discrete
series, but if $g \geq 3$, then neither possibility occurs, so the automorphic
representations do not contribute to the cohomology (of any kind) of the
corresponding Shimura variety. Another way of seeing this is to compute the
possible infinitesimal characters of the automorphic representations
corresponding to automorphic vector bundles on the Shimura variety, or
equivalently the Hodge-Tate weights of the expected $2g$-dimensional symplectic
Galois representations; one finds that no Hodge-Tate weight can occur with
multiplicity bigger than $2$, while the symplectic Galois representations
coming from the \'{e}tale $H^1$ of a curve of genus $g$ have weights $0$,
$1$ each occurring with multiplicity $g$ (see \cite{BCGP}).

  The primary problem is that not all automorphic forms contribute to
cohomology (e.g. Maass forms), in particular, not all automorphic forms
appear in any context of cohomology of Shimura varieties at all. Since
cohomology of Shimura varieties is currently the only approach for passing
from automorphic forms to Galois representations, people are thinking a lot
about how to move from any given context to a Shimura variety context, by
applying functoriality or $p$-adic functoriality. Nevertheless, as a very
simple example, for $\text{GL}(2)$ over $\mathbb{Q}$, it is conjectured that
there is a canonical bijection between continuous even irreducible
two-dimensional representations $\text{Gal}(\overline{\mathbb{Q}}/\mathbb{Q})$
and normalized algebraic cuspidal Maass new eigenforms on the upper half plane.
This is a sort of non-holomorphic analogue of the Deligne-Serre theorem which
relates the odd irreducible Galois representations to holomorphic weight one
newforms. In the holomorphic modular form case we can use coherent cohomology
of the modular curve considered as an algebraic curve over $\mathbb{Q}$, or
(in weight two or more) singular cohomology of a typically non-trivial local
system on the curve. However, in the case of algebraic Maass forms, they are
not cohomological, so we can not expect to see them in singular cohomology
of a local system, and they are not holomorphic, so we can not expect to see
them in coherent cohomology either.

  Even in the context of motives, there are some very difficult problems
which can not be solved by cohomology. In particular, for curves with genus
greater than two, the automorphic representations do not contribute to the
cohomology of any kind of the corresponding Shimura variety. This leads us
to study these problems in the context of anabelian algebraic geometry. Instead
of various cohomology of Shimura varieties, we can use the defining ideals of
modular curves. Theorem 1.1 and Theorem 1.3 show that even for curves with
higher genus, i.e., $\mathcal{L}(X(p))$ $(p \geq 7$), we can still prove
some modularity about them by the connection between Galois representations
and the elliptic modular functions.

\begin{center}
{\large\bf 12.3. Relation with Kottwitz's conjecture and Arthur's conjecture}
\end{center}

  According to \cite{Arthur1988}, the ultimate goal is of course to prove
reciprocity laws between the arithmetic information conveyed by $\ell$-adic
representations of Galois groups, and the analytic information wrapped
up in the Hecke operators on $L^2$-cohomology. A comparison of the two
Lefschetz formulas would lead to generalizations of the results in
Langlands (see \cite{Langlands73}) for $\text{GL}(2)$. Here, we need
stable trace formulas and endoscopic groups. In fact, Arthur's global
conjectures begin with the understanding that automorphic representations
$\pi$ should also occur in packets. i.e., it asserts that any irreducible
representation in $\Pi(G(\mathbb{A}_F))$ which occurs in
$L^2(G^0(F) \backslash G^0(\mathbb{A}_F))$ must belong to one of the
packets $\Pi_{\psi}$. It also provides a multiplicity formula, which
requires some further description (see \cite{Arthur1989} and \cite{Arthur1984}).
In his paper \cite{Arthur1996}, Arthur studied a special case of the
Siegel moduli space $S(N)$ of level $N$.

  Let $G$ be a connected reductive group over $\mathbb{Q}$ and let
$X_{\infty}$ be a $G(\mathbb{R})$-conjugacy class of homomorphisms
$$h: R_{\mathbb{C}/\mathbb{R}} \mathbb{G}_m \rightarrow G_{\mathbb{R}}$$
satisfying conditions (2.1.1.1)--(2.1.1.3) of Deligne in \cite{De1979}.
We assume further that the derived group $G_{\text{der}}$ of $G$ is
simply connected and that the maximal $\mathbb{Q}$-split torus in the
center of $G$ coincides with the maximal $\mathbb{R}$-split torus in
the center of $G$. For any sufficiently small compact open subgroup
$K \subset G(\mathbb{A}_f)$ we get a smooth Shimura variety $S_K$ over
$\mathbb{C}$, whose complex points are given by
$$S_K(\mathbb{C})=G(\mathbb{Q}) \backslash (X_{\infty}
                  \times (G(\mathbb{A}_f)/K)).$$

  For $h \in X_{\infty}$ we write $\mu_h$ for the restriction of
$h_{\mathbb{C}}: (R_{\mathbb{C}/\mathbb{R}} \mathbb{G}_m)_{\mathbb{C}}
\rightarrow G_{\mathbb{C}}$ to the first factor of $(R_{\mathbb{C}/\mathbb{R}}
\mathbb{G}_m)_{\mathbb{C}}=\mathbb{G}_m \times \mathbb{G}_m$. The
$G(\mathbb{C})$-conjugacy class of $\mu_h$ is well-defined, and its
field of definition, a number field $E$ contained in $\mathbb{C}$, is
called the reflex field. It is known that there is a canonical model
of $S_K$ over $E$.

  Let $L$ be a number field and let $\xi$ be a finite dimensional
representation of $G$ on an $L$-vector space. Then $\xi$ gives rise to
a local system $\mathcal{F}$ of $L$-vector spaces on $S_K(\mathbb{C})$.
For any finite place $\lambda$ of $L$ the local system $\mathcal{F} \otimes
L_{\lambda}$ comes from a smooth $L_{\lambda}$-sheaf $\mathcal{F}_{\lambda}$
on $S_K$ over $E$.

  Let $\overline{S}_K$ be the Baily-Borel-Satake compactification of $S_K$;
it too has a canonical model over $E$. Consider the intersection cohomology
groups
$$W^i=IH^{i}(\overline{S}_K(\mathbb{C}), \mathcal{F}).$$
The Hecke algebra $\mathcal{H}_L$ of locally constant $K$-bi-invariant
$L$-valued functions on $G(\mathbb{A}_f)$ acts on $W^i$. For any finite
place $\lambda$ of $L$ the $\lambda$-adic vector space $W^i \otimes_L L_{\lambda}$
has algebraic meaning and carries an action of $\text{Gal}(\overline{\mathbb{Q}}/E)$
commuting with the action of $\mathcal{H}_L$. One would like to understand
these commuting actions.

  Let $p$ be a rational prime. Assume that $G$ and $K$ are unramified
at $p$. This means simply that $G$ is quasi-split over $\mathbb{Q}_p$
and split over an unramified extension of $\mathbb{Q}_p$, and that
$K$ is of the form $K^p \cdot K_p$, where $K^p$ is a compact open
subgroup of $G(\mathbb{A}_f^p)$ and $K_p=G(\mathbb{Z}_p)$ for some
extension of $G$ to a flat group scheme over $\mathbb{Z}_p$ with
connected reductive geometric fibers. Note that $E$ is necessarily
unramified at $p$. Let $\mathfrak{p}$ be a place of $E$ lying over
$p$. One expects that for every finite place $\lambda$ of $L$ not
lying over $p$, the representations $W^i \otimes L_{\lambda}$ of
$\text{Gal}(\overline{\mathbb{Q}}_p/E_{\mathfrak{p}})$ are unramified.
Granting this, a first step in understanding the commuting representations
of $\mathcal{H}_L$ and $\text{Gal}(\overline{\mathbb{Q}}/E)$ on $W^i \otimes
L_{\lambda}$ would be to find a suitable expression for the trace of
$f \times \Phi_{\mathfrak{p}}^j$ on the virtual representation
$$W_{\lambda}:=\bigoplus_{i=0}^{2 \dim S_K} (-1)^i W^i \otimes L_{\lambda},$$
where $f \in \mathcal{H}_L$, $\Phi_{\mathfrak{p}}$ denotes a geometric
Frobenius element of $\text{Gal}(\overline{\mathbb{Q}}_p/E_{\mathfrak{p}})$,
and $j$ is a positive integer. Complications arise unless $f$ is of the
form $f^p f_{K_p}$, where $f^p$ is a function on $G(\mathbb{A}_f^p)$ and
$f_{K_p}$ denotes the characteristic function of $K_p$ in $G(\mathbb{Q}_p)$,
divided by the measure of $K_p$; for simplicity, we assume that $f$ has
this special form.

  One expects to have an expression for $\text{tr}(f \times \Phi_{\mathfrak{p}}^j;
W_{\lambda})$. In order to do this, we need Arthur's conjectures. The
global version of Arthur's conjecture give a formula for the multiplicity
of an irreducible representation of an ad\`{e}le group in the discrete
spectrum. Let $A_G$ denote the maximal $\mathbb{Q}$-split torus in the
center of $G$ and let $A_G(\mathbb{R})^0$ denote the identity component
of the topological group $A_G(\mathbb{R})$. Let $\chi$ be a quasi-character
on $A_G(\mathbb{R})^0$. One expects that the cuspidal tempered part of
$L_{\chi}^2(G(\mathbb{Q}) \backslash G(\mathbb{A}))$ is isomorphic to
$$\bigoplus_{[\varphi]} \bigoplus_{\pi \in \Pi(\varphi)}
  m(\varphi, \pi) \cdot \pi,$$
where $m(\varphi, \pi)=|\mathfrak{S}_{\varphi}|^{-1}
\sum_{x \in \mathfrak{S}_{\varphi}} \langle x, \pi \rangle$.

  Here $[\varphi]$ stands for the equivalence class of an admissible
homomorphism $\varphi: \mathcal{L} \rightarrow {}^{L}G$ ($\mathcal{L}$
denotes the conjectural Langlands group of $\overline{\mathbb{Q}}/\mathbb{Q}$),
and the first direct sum is taken over all equivalence classes of elliptic
$\varphi$ such that the associated quasi-character on $A_G(\mathbb{R})^0$
equals $\chi$. We write $\Pi(\varphi)$ for the conjectural $L$-packet of
$\varphi$ and write $S_{\varphi}$ for the group of self-equivalences of
$\varphi$ and $\mathfrak{S}_{\varphi}$ for the quotient of $S_{\varphi}$
by $S_{\varphi}^0 \cdot Z(\widehat{G})$, where $S_{\varphi}^0$ denotes
the identity component of $S_{\varphi}$. Here, we choose a complex group
$\widehat{G}$ of dual reductive type to $G$ and write $Z(\widehat{G})$
for its center.

  Let $\psi: \mathcal{L} \times \text{SL}(2, \mathbb{C}) \rightarrow {}^{L}G$
be a homomorphism. We say that $\psi$ is an Arthur parameter if the restriction
of $\psi$ to $\mathcal{L}$ is an essentially tempered Langlands parameter,
and we say that $\psi$ is elliptic if $S_{\psi}^0$ is contained in
$Z(\widehat{G})$. Arthur conjectured (see \cite{Arthur1989}) that the
discrete part of $L_{\chi}^2(G(\mathbb{Q}) \backslash G(\mathbb{A}))$ is
isomorphic to
$$\bigoplus_{[\psi]} \bigoplus_{\pi \in \Pi(\psi)} m(\psi, \pi) \cdot \pi,$$
where $m(\psi, \pi)=|\mathfrak{S}_{\psi}|^{-1} \sum_{x \in \mathfrak{S}_{\psi}}
\epsilon_{\psi}(x) \langle x, \pi \rangle$. In particular, the multiplicities
are the deepest part of the problem.

  Here the first sum is taken over all equivalence classes of elliptic
Arthur parameters $\psi$ such that the associated quasi-character on
$A_G(\mathbb{R})^0$ equals $\chi$. We write $\Pi(\psi)$ for the conjectural
Arthur packet of $\psi$. The notion of equivalence and the groups $S_{\psi}$,
$\mathfrak{S}_{\psi}$ are defined in the same way as for $\varphi$. The
new ingredient is the character $\epsilon_{\psi}: \mathfrak{S}_{\psi}
\rightarrow \{ \pm 1 \}$.

  Now, we give Kottwitz's conjectural description (\cite{Kottwitz}) of
$IH^{*}(\overline{S}_K, \mathcal{F})$. Let $\mathcal{L}_E$ denote the
Langlands group of $\overline{\mathbb{Q}}/E$ and $\psi_E$ denote the
restriction of $\psi$ to $\mathcal{L}_E \times \text{SL}(2, \mathbb{C})$.
Let $V$ be the usual irreducible representation of ${}^{L}(G_E)$
determined by $X_{\infty}$. Via $\psi_E$ the group $\mathcal{L}_E \times
\text{SL}(2, \mathbb{C})$ acts on $V$. This action commutes with the
action of $S_{\psi}$ on $V$. Choose $h \in X_{\infty}$. As above we get
$\mu_h: \mathbb{G}_m \rightarrow G_{\mathbb{C}}$. Using $\overline{\mathbb{Q}}
\rightarrow \mathbb{C}$, we get a well-defined $G(\overline{\mathbb{Q}})$-conjugacy
class of homomorphisms $\mu_0: \mathbb{G}_m \rightarrow G_{\overline{\mathbb{Q}}}$.
Choose a maximal torus $\widehat{T}$ of $\widehat{G}$. Then dual to $\mu_0$
is a Weyl group orbit of elements $\mu^{*} \in X^{*}(\widehat{T})$. We define
$\mu_1$ to be the restriction of $\mu^{*}$ to $Z(\widehat{G})$. The subgroup
$Z(\widehat{G})$ of $S_{\psi}$ acts on $V$ by the inverse of the character
$\mu_1$. Let $\nu$ be a character on $S_{\psi}$ whose restriction to
$Z(\widehat{G})$ is $\mu_1$. We write $V_{\nu}$ for the largest subspace
of $V$ on which $S_{\psi}$ acts by the inverse of $\nu$, so that $V=\oplus_{\nu}
V_{\nu}$ as representations of $\mathcal{L}_E \times \text{SL}(2, \mathbb{C})$.
Let $\phi$ be the Langlands parameter associated to $\psi$, and let
$\phi_E$ be its restriction to $\mathcal{L}_E$. Then $\phi_E$ gives
another action of $\mathcal{L}_E$ on each $V_{\nu}$.

  Let $\pi_f$ belong to the Arthur packet of $G(\mathbb{A}_f)$ parameterized
by $\psi$. For any character $\nu$ as above we define a non-negative integer
$m(\pi_f, \nu)$ by choosing $\pi_{\infty} \in \Pi(\psi_{\infty})$ and taking
$m(\pi_f, \nu)$ to be multiplicity with which the one-dimensional character
$\nu-\lambda_{\pi_{\infty}}$ occurs in the character $\overline{x} \mapsto
\epsilon_{\psi}(\overline{x}) \cdot \langle \overline{x}, \pi \rangle$ of
$\mathfrak{S}_{\psi}$ (one expects that the function $\overline{x} \mapsto
\epsilon_{\psi}(\overline{x}) \cdot \langle \overline{x}, \pi \rangle$ is
the character of some finite dimensional representation of $\mathfrak{S}_{\psi}$).
Note that $m(\pi_f, \nu)$ should be independent of the choice of $\pi_{\infty}$.

  For $\psi$ and $\nu$ as above we write $V(\psi, \nu)$ for the twist
of the representation $V_{\nu}$ of $\mathcal{L}_E$ (acting via
$\phi_{\psi}$) by the unramified quasi-character $| \cdot |^{-(\dim S_K)/2}$
of $\mathcal{L}_E$. Kottwitz showed that (see \cite{Kottwitz})
$\text{tr}(f \times \Phi_{\mathfrak{p}}^j; W_{\lambda})$ should be
equal to
$$\sum_{[\psi]} \sum_{\pi_f} \text{tr} \pi_f(f_{\mathbb{C}}) \sum_{\nu}
  m(\pi_f, \nu) \cdot (-1)^{q(G)} \cdot \nu(s_{\psi}) \cdot
  \text{tr}(\Phi_{\mathfrak{p}}^{j}; V(\psi, \nu)).$$

  Let $\mathcal{M}_E$ denote the motivic Galois group of
$\overline{\mathbb{Q}}/E$ relative to the Betti fiber functor and the
inclusion of $E$ in $\mathbb{C}$. Then $V(\psi, \nu)$ should come from
a complex representation of $\mathcal{M}_E$, and, after $L$ is enlarged
suitably, this complex representation should come from a representation
of $\mathcal{M}_E$ on an $L$-vector space $V(\psi, \nu)_L$ (here we use
our choice of embedding $L \hookrightarrow \mathbb{C}$). Then, letting
$V(\psi, \nu)_{\lambda}$ denote the $\lambda$-adic completion of
$V(\psi, \nu)_L$, we should get a $\lambda$-adic representation of
$\text{Gal}(\overline{\mathbb{Q}}/E)$ on $V(\psi, \nu)_{\lambda}$, and the
above expression for $\text{tr}(f \times \Phi_{\mathfrak{p}}^j; W_{\lambda})$
suggests that the virtual representation
$$\bigoplus_{i=0}^{2 \dim S_K} (-1)^i IH^i(\overline{S}_K, \mathcal{F})
  \otimes L_{\lambda}$$
of $\mathcal{H}_L$ and $\text{Gal}(\overline{\mathbb{Q}}/E)$ is equal to
$$\bigoplus_{[\psi]} \bigoplus_{\pi_f} \bigoplus_{\nu} n(\pi_f, \nu)
  (\pi_f^K \otimes V(\psi, \nu)_{\lambda}),$$
where $[\psi]$ ranges over the equivalence classes of Arthur parameters
$\psi$ such that $\psi_{\infty}$ is cohomological for $\xi_{\mathbb{C}}$.
Here $\pi_f^K$ denotes the space of $K$-fixed vectors in $\pi_f$, and
$n(\pi_f, \nu)=(-1)^{q(G)} \cdot \nu(S_{\psi}) \cdot m(\pi_f, \nu)$. The
action of $\text{SL}(2)$ on $V(\psi, \nu)$ should give the Lefschetz
decomposition. The individual space $IH^i(\overline{S}_K, \mathcal{F})
\otimes L_{\lambda}$ can be reconstructed from the Euler characteristic
by taking $(-1)^i$ times the part of weight $i$; in terms of $\psi$ this
amounts to taking weight spaces for $\text{SL}(2)$ (up to a shift by
$\dim S_K$). This last decomposition could be regarded as the ultimate
goal.

  The simplest example of such a structure is the tower of modular curves
$M_n^0(\mathbb{C})=\mathbb{H}/\Gamma(n)$ (see \cite{Illusie}). It corresponds
to the Shimura datum $(G, X)$ where $G=\text{GL}(2)$,
$X=\mathbb{C}-\mathbb{R}=\mathbb{H} \cup -\mathbb{H}$ the conjugacy class
of the canonical inclusion $h_0: \mathbb{S} \hookrightarrow G$
$(x+iy \mapsto \left(\begin{matrix} x & y\\ -y & x \end{matrix}\right))$,
i.e., the homogeneous space
$\text{GL}(2, \mathbb{R})/\mathbb{R}^{*} \text{SO}(2, \mathbb{R})$, orbit
of $i$ in $\mathbb{C}$ under the natural action of $\text{GL}(2, \mathbb{R})$.
We have $M_n^{0}(\mathbb{C})={}_{K_n}M_{\mathbb{C}}(G, X)$ for
$K_n=\text{Ker}(G(\widehat{\mathbb{Z}}) \rightarrow G(\mathbb{Z}/n))$. The
projective limit $M_{\mathbb{C}}(G, X)=\varprojlim_{n} M_n^{0}(\mathbb{C})$
has an action of $G(\mathbb{A}_f)$.

  In contrast with Kottwitz's conjecture as well as Arthur's conjecture,
in particular Langlands's results for $\text{GL}(2)$ (see \cite{Langlands73}),
for modular curves $X(p)$ in the context of motives, Theorem 1.1 and Theorem
1.3 give reciprocity laws between the arithmetic information conveyed by
$\mathbb{Q}(\zeta_p)$-rational representations $\rho_p$ of Galois groups,
and the analytic information wrapped in the action of
$\text{PSL}(2, \mathbb{F}_p)$ on the defining ideals $I(\mathcal{L}(X(p)))$
of modular curves $X(p)$ in the context of anabelian algebraic geometry. In
particular, the analytic information leads to the elliptic modular functions.

\begin{center}
{\large\bf 12.4. Relation with the construction of interesting algebraic cycles}
\end{center}

  According to \cite{De1994} and \cite{DeHo}, little progress has been
made on this problem: the construction of interesting algebraic cycles.
However, in the context of anabelian algebraic geometry, the defining
ideals $I(\mathcal{L}(X(p)))$ and their decompositions under the action
of $\text{PSL}(2, \mathbb{F}_p)$ provide a class of explicit algebraic
cycles on which both the finite groups $\text{PSL}(2, \mathbb{F}_p)$
and the absolute Galois group $\text{Gal}(\overline{\mathbb{Q}}/\mathbb{Q})$
can act. In particular, the locus of modular curves $X(p)$ can be constructed
as an intersection of some $\text{PSL}(2, \mathbb{F}_p)$-invariant algebraic
cycles. For $p=7$, $11$ and $13$, we give such an explicit construction.

\begin{center}
{\large\bf 12.5. Relation with the realization of Shimura varieties}
\end{center}

  The Hodge and $\mathbb{A}_f$ realizations of a Shimura variety are a
fundamental part of its theory (see \cite{Arthur2024}, p. 157). Deligne
(see \cite{De1971} and \cite{De1979}) regarded Shimura varieties as
parameter spaces for certain Hodge structures. Following this viewpoint,
Langlands (see \cite{Langlands79}) treated an arbitrary Shimura variety
$S_K(G, X)$ as a moduli spaces of motives. This leads to the realization
of the motive. One is the Hodge realization of the Shimura variety $S_K$
on the de Rham cohomology $H_{(2)}^{*}(S_K(\mathbb{C}), \mathcal{F})$ of
$S_K(\mathbb{C})$, which can be regarded as the motive of the Shimura
variety $S_K$ over $E$. The most widely studied realization functor for
motives is defined by their \'{e}tale cohomology and its corresponding
compatible families of $\ell$-adic Galois representations. Known as the
$\mathbb{A}_f$-realization, it is of obvious arithmetic importance. In
this case, the Galois representations act on the $\ell$-adic \'{e}tale
version of the intersection cohomology $IH^{*}(\overline{S}_K, \mathcal{F})$
for the Baily-Borel-Satake compactification of $S_K$.

  The $\mathbb{A}_f$-realization of a motive is a compatible family of
$\ell$-adic representations of $\text{Gal}(\overline{\mathbb{Q}}/\mathbb{Q})$.
The Hodge realization is closely related to the period realization. Now,
by Theorem 1.1, Corollary 1.2 and Theorem 1.3, we give a completely
different realization by the defining ideals $I(S_K(G, X))$ of the
Shimura variety $S_K$, at least for $G=\text{GL}(2, \mathbb{Q})$.

\begin{center}
{\large\bf 12.6. Relation with the action of $\text{Gal}(\overline{\mathbb{Q}}/\mathbb{Q})$
                on Grothendieck's dessins d'enfants}
\end{center}

  In fact, a true understanding of the action of
$\text{Gal}(\overline{\mathbb{Q}}/\mathbb{Q})$ on Grothendieck's
dessins d'enfants will presumably require investigations of the relation
between the arithmetic features of $\text{Gal}(\overline{\mathbb{Q}}/\mathbb{Q})$
and the action of $\text{Gal}(\overline{\mathbb{Q}}/\mathbb{Q})$ on geometric
fundamental groups. Much progress has been made in this direction when the
fundamental group of $\mathbb{P}^1-\{ 0, 1, \infty \}$ is replaced by a
nilpotent or meta-abelian quotient (see \cite{De1989} and \cite{Ihara})
but the general case remains extremely unclear.

  In the spirit of \cite{Mumford} and \cite{SV}, there are five ways
of defining complex algebraic curves: (1) Writing an equation; (2)
Defining generators of the uniformizing Fuchsian groups; (3)
Specifying a point in the moduli space; (4) Introducing a metric;
(5) Defining Jacobian. Here, (1) is the most explicit way, and (2) is
closely related to the fundamental groups of the curves. However, (5)
is weaker than (1) and (2).

  According to \cite{Schneps}, dessins d'enfants are the drawings on
topological surfaces obtained by taking $\beta^{-1}([0, 1])$ where
$\beta: X \rightarrow \mathbb{P}^1(\mathbb{C})$ is a cover of Riemann
surfaces ramified over $0$, $1$ and $\infty$ only. They are in bijection
with the pairs $(\beta, X)$ up to isomorphism, and therefore with the
set of conjugacy classes of subgroups of finite index of
$\widehat{F}_2=\pi_1(\mathbb{P}^1-\{ 0, 1, \infty \})$. The link with
the genus zero mapping class groups is given by the fact that
$\mathbb{P}^1-\{ 0, 1, \infty \}$ is isomorphic to the first non-trivial
moduli space $\mathcal{M}_{0, 4}$. The two main problems of the theory
of dessins d'enfants are the following (see \cite{Schneps}): (i) given
a dessin, i.e. a purely combinatorial object, find the equations for
$\beta$ and $X$ explicitly; (ii) find a list of (combinatorial?
topological? algebraic) invariants of dessins which completely
identify their $\text{Gal}(\overline{\mathbb{Q}}/\mathbb{Q})$-orbits.
The second problem can be interestingly weakened from
$\text{Gal}(\overline{\mathbb{Q}}/\mathbb{Q})$ to the profinite version
of the Grothendieck-Teichm\"{u}ller group $\widehat{GT}$ (see \cite{Dr1990}),
but it remains absolutely non-trivial.

  In fact, Theorem 1.3 gives a connection between the Galois groups of the
$p$-th order transformation equations of the $j$-functions (i.e., the
arithmetic features of $\text{Gal}(\overline{\mathbb{Q}}/\mathbb{Q})$) and
the Galois representations $\rho_p: \text{Gal}(\overline{\mathbb{Q}}/\mathbb{Q})
\rightarrow \text{Aut}(\mathcal{L}(X(p)))$ (i.e., the action of
$\text{Gal}(\overline{\mathbb{Q}}/\mathbb{Q})$ on the defining ideals, which
can determine their geometric fundamental groups).

\vskip 2.0 cm

\noindent{Department of Mathematics, Peking University}

\noindent{Beijing 100871, P. R. China}

\noindent{\it E-mail address}: yanglei$@$math.pku.edu.cn

\end{document}